\documentclass[12pt,leqno]{article} %book me met 3 num\'eros % article % report
\usepackage{amsmath}
\usepackage{amssymb}
\usepackage{amsthm} % remplace :  \usepackage{theorem}
\usepackage{epsf}
\usepackage{graphicx}
\usepackage{esint} % added allow $\fint$
\usepackage{dsfont}
\usepackage{hyperref} % Jonas advise. and the next line
\hypersetup{pdfpagemode=FullScreen,  colorlinks=true} 
\numberwithin{equation}{section} 

\newtheorem{lem}{Lemma}[section] % section
\newtheorem{pro}[lem]{Proposition}
\newtheorem{defi}[lem]{Definition}
\newtheorem{thm}[lem]{Theorem}
\newtheorem{cor}[lem]{Corollary}
\theoremstyle{remark}
\newtheorem{rem}[lem]{\bf Remark}

\newcommand{\nn}{\nonumber}
\newcommand{\ms}{\medskip}
\newcommand{\msi}{\par\medskip\noindent}
\newcommand{\ssi}{\par\smallskip\noindent}
\newcommand{\pari}{\par\noindent}
\newcommand{\R}{\mathbb{R}}
\renewcommand{\H}{\mathcal H}
\newcommand{\B}{\mathbb{B}}

\newcommand{\bN}{\mathbb{N}}
\renewcommand{\S}{\mathbb {S}}
\newcommand{\bY}{\mathbb {Y}}
\newcommand{\bT}{\mathbb {T}}
\newcommand{\bV}{\mathbb {V}}
\newcommand{\bP}{\mathbb {P}}
\newcommand{\bH}{\mathbb {H}}
\newcommand{\bZ}{\mathbb {Z}}

\renewcommand{\d}{\partial}
\newcommand{\dist}{\,\mathrm{dist}}
\newcommand{\length}{\,\mathrm{length}}
\newcommand{\ddist}{\mathrm{dist}_\S} 
\newcommand{\Angle}{\,\mathrm{Angle}}
\newcommand{\dsp}{\displaystyle}
\newcommand{\sm}{\setminus}
\newcommand{\diam}{\mathrm{diam}}

\newcommand{\wt}{\widetilde}
\newcommand{\wh}{\widehat}
\newcommand{\ol}{\overline}
\newcommand{\ub}{\underbar}

\newcommand{\cE}{{\cal E}}
\newcommand{\cL}{{\cal L}}

\newcommand{\cX}{{\cal X}}
\newcommand{\cF}{{\cal F}}
\newcommand{\cI}{{\cal I}}
\newcommand{\cJ}{{\cal J}}

\newcommand{\cZ}{{\cal Z}}

\newcommand{\cR}{{\cal R}}
\newcommand{\cC}{{\mathfrak C}}
\newcommand{\1}{{\mathds 1}}
\newcommand{\2}{{\frac{2\pi}{3}}}

\begin{document}

\title{A local description of $2$-dimensional almost minimal sets bounded by a curve
% near some cones
}\author{
G. David\footnote{Univ Paris-Sud, Laboratoire de Math\'{e}matiques, 
UMR 8658 Orsay, F-91405}
\footnote{CNRS, Orsay, F-91405}
\footnote{
 G.\, David is/was partially supported by the Institut Universitaire de France
 the ANR, programme blanc GEOMETRYA ANR-12-BS01-0014, 
the European Community Marie Curie grant MANET 607643 and H2020 grant GHAIA 777822, 
and the Simons Collaborations in MPS grant 601941, GD.}
 }
\date{}
\maketitle

\ms\noindent{\bf Abstract.}
We study the local regularity of sliding almost minimal sets of dimension $2$ in $\R^n$, bounded by a
smooth curve $L$. These are a good way to model soap films bounded by a curve,
and their definition is similar to Almgren's. We aim for a local description, in particular near $L$ and
modulo $C^{1+\varepsilon}$ diffeomorphisms, of such sets $E$, but in the present
paper we only obtain a full description when $E$ is close enough to a half plane, 
a plane or a union of two half planes bounded by the same line, or a transverse minimal cone 
of type $\bY$ or $\bT$.
The main tools are adapted near monotonicity formulae for the density, including for balls that are not centered
on $L$, and the same sort of construction of competitors as for the generalization of J. Taylor's
regularity result far from the boundary.

\ms\noindent{\bf R\'esum\'e en Fran\c cais.}
On \'etudie la r\'egularit\'e locale des ensembles presque minimaux de dimension $2$ dans $\R^n$,
bord\'es par une courbe lisse $L$, et avec une condition glissante de bord semblable \`a celle d'Almgren.
Ces ensembles semblent le meilleur mod\`ele pour les films de savon bord\'es par une courbe.
Le but est d'obtenir une description locale de ces ensembles, en particulier pr\`es de $L$ et 
modulo un diff\'eomorphisme de classe $C^{1+\varepsilon}$. 
Dans ce papier on n'obtient une description compl\`ete que 
lorsque $E$ est assez proche d'un demi plan, un plan ou une union de deux demi plans bord\'es 
par la m\^eme droite, ou un c\^one minimal de type $\bY$ ou $\bT$ transverse \`a $L$. 
Les outils principaux sont des formules de presque monotonie adapt\'ees pour la densit\'e, 
y compris pour des boules qui ne sont pas centr\'ees sur $L$, et la construction du m\^eme genre de 
comp\'etiteurs que pour la g\'en\'eralisation du r\'esultat de J. Taylor sur la r\'egularit\'e loin du bord.

\ms\noindent{\bf Key words/Mots cl\'es.}
Almost minimal sets, sliding boundary condition, Plateau problem.

\ms\noindent

\vfill\eject
\tableofcontents

\listoffigures % affiche vraiment la liste

\part{Description of the results}
\section{Introduction}
\label{S1}

The goal of this paper is to start a study of the local behavior of two-dimensional soap
films near a smooth one-dimensional boundary. Our model for soap films, 
which will be discussed soon, is given by the notion of ``sliding almost minimal sets''. 
This is not so far from Almgren's notion of ``restricted set'' from \cite{AlmgrenMemoir}, 
and we would like to obtain along the boundary a description which is 
similar to Jean Taylor's regularity result \cite{Ta} far from the boundary.

\begin{figure}[!h]  
\centering
\includegraphics[width=6.cm]{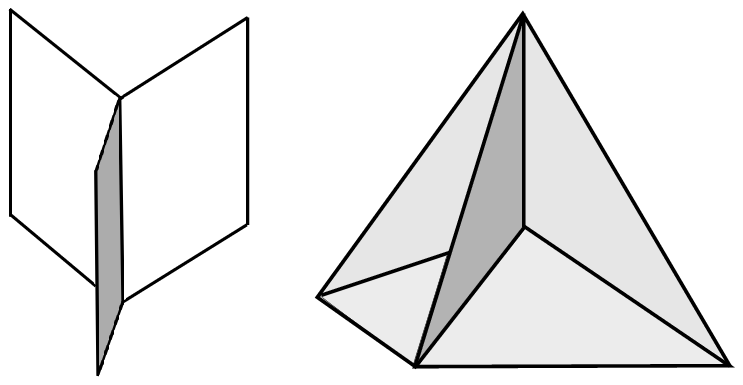}
\caption{A cone of type $\bY$ and a cone of type $\bT$ \label{f1.1}}
\end{figure}

Let us say a few words about the result of J. Taylor that we would like to imitate here.
There are actually two main steps to it, and the first one is a full description of the minimal cones 
(with the same definition of minimality as in \cite{AlmgrenMemoir} and roughly here) 
of dimension $2$ in $\R^3$. 
These are the planes, the cones of type $\bY$ composed of three half planes bounded by a same line 
and that make $\frac{2\pi}{3}$ angles along that line, and the cones of type $\bT$. 
A cone of type $\bT$ is the cone over the union of the $6$ edges of a regular tetrahedron 
centered at the origin; see Figure~\ref{f1.1}.  % p6
This first part is important because the blow-up limits of any almost minimal set $E$ at a Lebesgue 
density point of $E$ is a minimal cone (blow-up limits will be defined and commented a little more
near \eqref{17.7}). The second part consists in proving that under suitable assumptions,
all the blow-up limits of $E$ at such a point $x_0$ are equal, and that there is a small neighborhood 
of $x_0$ where $E$ is equivalent, through a $C^{1+\beta}$ diffeomorphism of $\R^3$, to this minimal cone. 
Thus J. Taylor's theorem gives a local classification of class $C^{1+\beta}$ of the almost minimal sets.

A partial generalization of this result was given in \cite{Holder} and \cite{C1}, that gives a local description of 
almost minimal sets of dimension $2$ in $\R^n$, but with two differences. First, the full list of minimal
cones of dimension $2$ in $\R^n$, $n\geq 4$, is not known; we just have a combinatoric description
in terms of faces. But also, if $X$ is a blow-up limit of the almost minimal set $E$ at $x_0$, 
we only prove the $C^{1+\beta}$ equivalence of $E$ to $X$ near $x_0$ when $X$ satisfies an 
additional property, the full length property. Otherwise, we only get a local biH\"older equivalence. 
The full length property, which is a metric property of the net of geodesics that compose $X\cap \d B(0,1)$,
is related to an epiperimetric inequality; it is satisfied by the planes and the cones of type $\bY$ and $\bT$,
but we do not know whether it is true in general, or whether $E$ is always $C^{1+\beta}$ equivalent
near $x_0$ to any of its blow-up limits at $x_0$, or even whether the blow-up limit is unique.

\subsection{Sliding almost minimal sets}

We would like to have similar theorems for almost minimal sets subject to a boundary 
constraint along some boundary set $L$ of dimension $1$, where we would describe $E$ 
near any point $x_0 \in E \cap L$, but we shall only be able to give such a description in some cases.
Before we do this, let us explain some of our definitions relative to \ub{almost minimal sets}, 
\ub{the sliding condition}, and the \ub{sliding Plateau problem}.

\ms
We give the definitions for arbitrary dimensions and boundaries, because this will not hurt.
In the discussion that follows, $L$ (our boundary) is a given closed subset of $\R^n$, and 
$d \in [1,n]$ is an integer, the dimension of our sets. Our putative almost minimizers will be closed
sets $E \subset \R^n$, with locally finite $d$-dimensional Hausdorff measure. 
That is, $\H^d(E \cap B(0,R)) < +\infty$ for $R > 0$.
We start with the notion of competitors.

\begin{defi}\label{t1.1} 
Let $E \subset \R^n$ be a closed set, and let $B = \overline B(x,r)$ be a closed ball.
A \underbar{deformation} of $E$ in $B$ (with sliding boundary $L$) is a one parameter family
$\{ \varphi_t \}, 0 \leq t \leq 1,$ of continuous mappings $\varphi_t : E \to \R^n$, such that
\begin{equation}\label{1.4}
\varphi(x,t) = \varphi_t(x) \ \text{ is a continuous function of $(x,t) \in E \times [0,1]$,}
\end{equation}
\begin{equation}\label{1.5}
\varphi_t(x) = x \text{ for $t = 0$ and for $x\in E \sm B$,}
\end{equation}
\begin{equation}\label{1.6}
\varphi_t(E \cap B) \subset B \ \text{ for } 0 \leq t \leq 1,
\end{equation}
\begin{equation}\label{1.7}
\varphi_t(x) \in L \ \text{ when } x \in E \cap L,
\end{equation}
and 
\begin{equation}\label{1.8}
\varphi_1 \ \text{ is Lipschitz on $E$.}
\end{equation}
A \underbar{sliding competitor} for $E$ in $B$ is a set $F = \varphi_1(E)$, where the family
$\{ \varphi_t \}$ is a deformation of $E$ in $B$.
\end{defi}

This is reasonably close to the initial definitions of Almgren in \cite{AlmgrenMemoir}; let us
comment on the differences.

Here we are really interested in what happens near the boundary $L$, but otherwise we could always take
$L = \emptyset$, forget about the condition \eqref{1.7}, and be in the same conditions as in
\cite{AlmgrenMemoir} or \cite{Ta}. 

We decided to keep the extra constraint \eqref{1.8}, because it was put forward by Almgren
and does not hurt. It makes it a tiny bit harder for $F$ to be a sliding competitor, hence a little
easier for $E$ to be an almost minimal set (as defined below), and our regularity theorems will
then be a tiny bit stronger.

In the analogous definition without sliding boundary condition, we took the habit of
defining the $\varphi_t$ on the whole $\R^n$, but this makes no difference when 
there is no condition \eqref{1.7}, as it would be easy to extend the $\varphi_t$
from $E$ to $\R^n$. The case when we work in a complicated domain $\Omega$, 
and we should require the $\varphi_t$ to take values in $\Omega$, will not arise in this paper.

For similar reasons, if we did not have \eqref{1.7}, we would not need to mention the whole homotopy
$\{ \varphi_t \}$, $0 \leq t \leq 1$, because given $\varphi_1$ we could simply complete the
homotopy by taking $\varphi_t(x) = t \varphi_1(x) + (1-t) x$. Because of \eqref{1.7}, we need to
be a little more careful. Yet, since most of the time in this paper $L$ will be a line, hence convex, 
it will often be enough to construct $\varphi_1$ and complete by convexity.

\begin{defi}\label{t1.2}
Let $U \subset \R^n$ be open, let $L\subset \R^n$ be closed, and let 
$h : (0,+\infty) \to [0,+\infty]$ be a \underbar{gauge function}. This just means that
$h$ is nondecreasing, and that
\begin{equation}\label{1.9}
\lim_{r \to 0} h(r) = 0.
\end{equation}
A sliding $(U,L,h)$-\underbar{almost minimal set} (of dimension $d$)
is a set $E \subset U$, which is closed in $U$, such that for every compact ball 
$B = \overline B(x,r) \subset U$,
\begin{equation}\label{1.10}
\H^d(E\cap B) < +\infty  
\end{equation}
and more importantly
\begin{equation}\label{1.11}
\H^d(E \cap B) \leq \H^d(F \cap B) + h(r) r^d
\end{equation}
for every sliding competitor $F$ of $E$ in $B$. When $h = 0$, 
we say that $E$ is $(U,L)$-minimal, or that $E$ is minimal in $U$, with sliding boundary $L$.
\end{defi}

See for instance \cite{Mattila} for the definition of the Hausdorff measure $\H^d$.
Some simple comments will be useful before we continue.
The open set $U$ may be useful to localize the notion, but $U = \R^n$ is already an interesting choice.
If $U$ is not convex, an equally good definition would only require $B$ in \eqref{1.11} to be 
a compact subset of $U$ that is contained in a ball of radius $r$; we shall not see the difference here
because all our results will be local.

The definition of almost minimal sets by Almgren \cite{AlmgrenMemoir} and \cite{Ta} is essentially 
the same as above, but with $L = \emptyset$ and hence no constraint \eqref{1.7}.
There is a slight difference, in the way we do the accounting in \eqref{1.11}, 
with the definition of restricted sets in \cite{AlmgrenMemoir},
or the definition of quasiminimal sets in \cite{DSQM} and further references, which is that here
we compare the measures of $E$ and $F$, but we could have compared the measures of $E \cap W$
and $\varphi_1(E\cap W)$, where $W = \big\{ x\in E \, ; \, \varphi_1(x) \neq x \big\}$.
We took what seems to be the simplest definition, but our results also work with the slightly
different way of accounting. We refer to \cite{Holder} for more detail about the alternate definitions
in the plain case without boundary, and why the basic properties that we use are also true with the 
other definitions. This is then generalized to the sliding case in \cite{Sliding}.

For the main results of this paper, we will take $d=2$, $L$ will be a smooth curve, and 
even $L$ will almost always be a line.
Most of the time, $U$, $L$, and $h$ will be given, and we shall just say that $E$ is 
a sliding almost minimal set, or even an almost minimal set, without further reference 
to $U$, $L$, and $h$.

\ms
As far as the author knows, the notion of sliding almost minimal set was only introduced 
(at least explicitly) in \cite{Montreal} and \cite{SteinConf}, even though similar notions 
existed in the past. This notion seems to give the best model for soap films attached to 
a set $L$. It comes with an associated \ub{Plateau problem}: suppose $L$ is compact 
for simplicity, start with a closed set $E_0$ such that $\H^d(E_0) < +\infty$, and minimize 
$\H^d(E)$, or a similar functional, among all the sliding competitors 
$E$ for $E_0$ (say, in a very large ball). This seems like a natural problem to consider, 
and it is nice that different initial sets $E_0$ will often yield different solutions (typically, 
with a different topology or combinatorics), as it happens in real life. The fact that the infimum 
may be $0$ if $E_0$ is not properly attached, or has lower-dimensional competitors, 
does not disturb us.

Unfortunately, we do not know whether this Plateau problem always has a solution, 
but if it does its minimizers $E$ will be sliding minimal sets, or almost minimal if we minimize 
a functional which is different, but not too much, from $\H^d$. 

This is one good reason for studying the regularity of sliding almost minimal sets, 
in particular near the boundary, but let us mention other ones. First of all, there are a few 
other ways to state a Plateau problem, and for which solutions (in some cases if they exist) 
are also given by, or associated to, sliding almost minimal sets.
For instance, the solutions of the Reifenberg homological Plateau problem, 
as in \cite{R}, \cite{A68}, \cite{Dp}, or \cite{Fang}, are sliding almost minimal sets. 
But also, the solutions of the similar problems posed in \cite{HaPu1}, \cite{HaPu2}, 
\cite{Ita1}, \cite{Ita2}, \cite{Ita3}, \cite{Ita4} or the supports of some size minimizing currents
with homological conditions at the boundary, for instance as in \cite{Dp}, are like this. 
We refer to \cite{SteinConf} for a little more information on these problems (leading to 
the interest of sliding almost minimal sets), and for the proof of sliding minimality for solutions of the 
Reifenberg problem or the support of size minimizing currents.

Next, the author believes that the best way to try to prove existence results for the sliding 
Plateau problem alluded to above, or its analogue with size minimizing currents, is by proving regularity 
for almost minimal sets.
This is not shocking; for the sliding Plateau problem, for instance, we are often able to produce 
a good candidate $E$, as a limit of some well chosen minimizing sequence, and limiting theorems
allow us to prove that $E$ is an almost minimal set; it then remains to show that $E$ 
is itself a competitor for $E_0$, and this will be easier if we have a good control on $E$.
For instance, proving that there is a Lipschitz retraction defined on a neighborhood of $E$
and which preserves $L$ would be very useful.

Let us also mention that almost minimality is a nicely flexible notion to study.
That is, we may want to minimize a minor variant of the functional $\H^d(E)$,
such as $\dsp\int_E f(x) d\H^d(x)$, with a H\"older function $f$ such that 
$C^{-1} \leq f \leq C$, or even $\dsp\int_E f(x,T_E(x)) d\H^d(x)$, with 
functions $f$ that depend also on the approximate tangent $d$-plane $T_E(x)$ to $E$ at $x$
(in a simple enough way); minimizers of such functionals are still almost minimal sets,
with a gauge function that depends on the (low) regularity of $f$, so we may apply
the results of this paper to them. 
By contrast, it is unlikely that the corresponding varifolds, for instance, have a locally finite
first variation. It is easy to believe that such functionals could be used
to model variants of the soap film problem, even though the author does not have specific
examples to show.

\ms
So we want to study the local regularity of sliding almost minimal sets. The following notion of 
\ub{coral} (or \ub{reduced}) set will help simplify the statements. 

\begin{defi}\label{t1.2c}
The \underbar{core} of the closed set $E$ (in a given open set $U \subset \R^n$) is 
the closed support of $\H^d_{\vert E}$, i.e.,
\begin{equation}\label{1.12}
{\rm core}(E) = \big\{ x\in  E \, ; \, \H^d(E \cap B(x,r)) > 0
\text{ for all } r > 0 \big\}.
\end{equation}
We say that $E$ is \underbar{coral}, or \underbar{reduced}, when ${\rm core}(E) = E$. 
\end{defi}

It is not so hard to see that when $E$ is almost minimal, its core is also almost minimal, 
and that hence it is enough to reduce our attention to coral almost minimal sets. 
See Proposition 3.3 in \cite{Sliding}. 
This makes the statements simpler, because we won't have to worry about additional thin sets
of vanishing measure that could have nearly dense tentacles, for instance. 
From now on we shall always assume that all our (sliding) almost minimal sets are coral, even though 
we do not always repeat this. Similarly, we exclude the empty set (and hence sets of vanishing measure)
from our discussions, even though it is minimal.

Notice however that we do not say that sets $E$ that minimize functionals, like in the Plateau problems 
discussed above, are coral. We just say that their $d$-dimensional part, or core, is still minimal,
so with some luck we can get a good description of those. For the rest of $E$, it is very hard to control
it, unless we ask for a more specific way to present $E$ in a clean way, so that for instance no proper
subset of $E$ is a competitor for $E$. This last way to see things was the initial way to proceed,
in the context of the Mumford-Shah functional (hence the name ``reduced''), but in some cases it is 
probably not so easy to pick a competitor $E' \subset E$ which is minimal for inclusion, especially if 
we don't want to deform it first; the reduction to coral sets, which we choose to do here, 
is simpler and does the job.

\ms
There does not seem to be too many regularity results on sliding almost minimal sets near the boundary,
especially for $d$ larger than $2$. The issue was taken rather brutally in \cite{Sliding} 
(see also a more digestible account in \cite{Montreal}), where some basic properties were proved 
(local Ahlfors regularity, rectifiability, and in some dimensions uniform rectifiability) 
under fairly general assumptions. What we will use most in the present paper is a nice collection 
of limiting theorems, that we will take from \cite{Sliding}. 
For instance, if the sets $E_k$, $k \geq 0$, are coral (see above) and almost minimal in $U$ 
with a given gauge function $h$, and if they converge (locally for the Hausdorff distance, 
as will  be explained below) to a limit $E$, then $E$ is almost minimal with the same gauge function $h$. 
In addition,
\begin{equation}\label{1.9a}
\H^d(E \cap V) \leq \liminf_{k \to +\infty} \H^d(E_k \cap V)
\end{equation}
for every open set $V \subset U$.
This is very useful because it allows many proofs by compactness.
We will recall all these results more precisely when we use them. 

It seems difficult to go much further in a general situation, and in particular for $d > 2$,
so we now turn to more precise regularity results in very specific situations, by which we mean when
$d=2$ and $L$ is a simple set. Because we are happy with J. Taylor's regularity result from
\cite{Ta} when there is no boundary (or equivalently, since this is a local result, far from the boundary),
and in higher ambient dimensions, with its generalization in \cite{Holder} and \cite{C1}, we shall
concentrate on regularity results near a point of $L$.

In \cite{Ta2}, J. Taylor gives a good description (similar to the regularity result of \cite{Ta}) 
for sets of finite perimeter in a bounded domain $U$ (say, with smooth boundary), that minimize a 
functional that looks like the perimeter, multiplied by a suitable constant $\alpha \in (0,1]$ on $\d U$. 
This is quite similar to what we want to do here, but she works in a different category. 

Much more recently X. Fang \cite{Fang2} started a similar study for sliding minimal
sets (thus, with no more direct constraint on the domains bounded by $E$) when $n=3$,
$L$ is a smooth surface of dimension $2$, and $E$ is required to stay on one side of $L$.
He was then able to prove that near any point $x_0 \in L$, $E$ is H\"older-equivalent to a minimal cone.
The minimal cones that show up for this problem are the tangent plane $P_0$ to $L$ at $x_0$,
or the union of $P_0$ with a half plane orthogonal to $P_0$, or the union of $P_0$ and a half 
set of type $\bY$ orthogonal to $P_0$. More recently, with methods similar to those of 
\cite{C1}, he even proved the more precise $C^{1+\varepsilon}$-equivalence \cite{Fang3}. 
Also some variants of this problem, for instance with mixed conditions as in \cite{Ta2},
are likely to be interesting and feasible.

\subsection{Towards a classification of singularities}
\label{S1.2}

In in the present paper we concentrate on the case of $2$-dimensional sets $E$ in (an open set of) 
$\R^n$, when the boundary set $L$ is a smooth (at least $C^{1+\varepsilon}$) curve. 
In fact, we will concentrate on the simpler case when $L$ is a line, 
and in Section \ref{S31} explain rapidly how to deal with the general case.
We would have liked a complete description of all the tangent objects (sliding minimal cones associated 
to a boundary which is a line), and then a precise local description, if possible, in the $C^{1+\varepsilon}$ 
category, in terms of the tangent cones. If we had all this, we would probably get a good existence result too, 
but as we shall see soon, we still have an important missing case.

Again there does not seem to be too much available information on this specific classification problem.
G. Lawlor and F. Morgan give in \cite{LM} a list of expected behaviors of minimal sets along a boundary
which is a curve, which the reader can also find in Figure~13.9.3  (Ten conjectured types...) of \cite{Mo}. 
Compared to the presentation below, there are of course common points, but also some 
small differences. Also, K. Brakke gives in \cite{B} a description of minimal surfaces bounded by a curve, which will be rapidly discussed below. 

\ms
We start the presentation of our result with a \ub{list of sliding minimal cones}.

We start with the case when $L = \emptyset$. Recall that when $n=3$, we have the full description 
completed by \cite{Ta} (but started by Plateau, Lamarle \cite{Lam}, and Heppes \cite{He}), 
which says that the (coral) minimal cones of dimension $2$ in $\R^3$ are the planes 
(also called cones of type $\bP$), and the cones of type $\bY$ or $\bT$ defined above.

In ambient dimensions $n \geq 4$, the cones of type $\bP$, $\bY$, or $\bT$ are still minimal,
but there are other ones. 
The union $P_1 \cup P_2$ of two planes through the origin that are orthogonal to each other is minimal, 
and in \cite{PUP}, X. Liang showed that this stays true when the two planes are nearly orthogonal.  
There is a conjecture of F. Morgan \cite{MoConj} % with the P1 \cup P_2 conjecture
on the precise condition on the angle of $P_1$ and $P_2$ under which $P_1 \cup P_2$ is minimal;
G. Lawlor \cite{Law} proved that this condition is necessary,
but we do not know whether it is sufficient.
In \cite{YxY} X. Liang showed that the product $Y \times Y$ of two
$Y$-sets of dimension $1$ contained in orthogonal $2$-planes is minimal. 
But there may be many other ones that we did not guess. 
Nonetheless \cite{Holder} gives a reasonable description of these cones $X$, that says that 
$X \cap \d B(0,1)$ is composed of a finite number of arcs of great circles with constraints 
on their length and how they meet.
We shall be more specific about this in Section \ref{S2}, 
because we need the description.

Now let $L$ be a line in $\R^n$, which we assume contains the origin.
We mentioned all the cones above, because they are still sliding minimal with the boundary
$L$ (there are more constraints on the competitors, hence the sliding minimality condition is weaker).
And so are their translations (the fact that they are not centered on $L$ does not matter).
But the set of planes that contain $L$, which we shall denote by $\bP(L)$, and the set 
$\bY(L)$ of cones of type $\bY$ whose spine is equal to $L$ (or equivalently, which are 
composed of three half planes bounded by $L$), will play a special role, so we give them a name.

In addition to all of these, we know of two more sliding minimal cones, and a possible third one.
We start with the sets of type $\bH$, which are just the half planes bounded by $L$.
That is, we take any $2$-plane $P$ that contains $L$, keep one of the two connected components
of $P \sm L$, and take its closure (i.e., add $L$ back). We will denote by $\bH(L)$
the collection of sets of type $\bH$ bounded by $L$. 

The sets of type $\bV$ (bounded by $L$) are the unions $V = H_1 \cup H_2$ of two half planes 
$H_1, H_2 \in \bH(L)$ which make an angle $\alpha \in [\2,\pi]$ with each other along $L$.
This last means that if $e_i$ is the unit vector in $H_i$ that is orthogonal to $L$, then
$\langle e_1, e_2 \rangle \leq -\frac12$. When $\alpha = \pi$, we get a plane of $\bP(L)$.
We denote by $\bV(L)$ the collection of sets of type $\bV$ bounded by $L$. 

Even though this is not needed for the main results of this paper, we decided
to include in Section~\ref{S32} a proof of the fact that sets of type $\bH$ and $\bV$ 
are sliding minimal, even with a possibly larger boundary set $L$ and the 
variant of Definition \ref{t1.1} where we do not require \eqref{1.8},
because this was apparently not written before. 

The following sets were suggested by X. Liang (in addition to those of \cite{LM} and \cite{Mo})
as possible sliding minimal cones.
Let $Q$ be a cube (in $\R^3$), and assume that one of the great diagonals of
$Q$ is contained in $L$. Then let $X$ denote the (positive) cone over the union of the edges of $Q$;
$X$ will be called a set of type $\mathbb Q$. Some experiments (including, with soap) suggest that 
the sets of type $\mathbb Q$ are sliding minimal, even though we know that they are not plain minimal 
(i.e., with no boundary constraint) because of the Plateau-Lamarle-Heppes-Taylor characterization above. 
But we do not have a proof of sliding minimality.

Again there may be lots of other sliding minimal cones that we do not know about, 
but at least we give in Section \ref{S2} a combinatorial description of these cones, 
similar to the one we have for plain minimal cones.

\ms
We now turn to our \ub{tentative classification of singularities}. By this we mean a local description of 
hopefully every sliding minimal set $E$, near any point $x_0 \in E \cap L$. Of course this description
will depend on the type of minimal cones $X$ that approximate $E$ on small balls centered at $x_0$.
This may mean, on the blow-up limits of $E$ at $x_0$, but we prefer to go directly to a quantitative 
statement with an approximation of $E$ by a minimal cone in a given ball 
$B(x_0,10r_0) \subset \R^n$ (we allow any ambient dimension $n$). 

We shall assume, for the following discussion, that
\begin{equation}\label{1.11a} 
E \ \text{ is a coral sliding $(B(0,10 r_0), L, h)$-almost minimal set,}
\end{equation}
with a gauge function $h$ such that
\begin{equation}\label{1.12a}
h(r) \leq C_h r^\beta \ \text{ for } 0 < r \leq 10r_0,
\end{equation}
for some $\beta\in (0,1]$ and some constant $C_h \geq 0$ such that $C_h r_0^\beta$ is small enough. 
Let us say,
\begin{equation}\label{1.13a}
C_h r_0^\beta \leq \varepsilon_0
\end{equation}
for some small $\varepsilon_0 > 0$ that we get to choose, depending in particular on $n$ and $\beta$.

We further assume that $L$ is a line through the origin; we shall explain in Section \ref{S31}
that similar statements hold when $L$ is a curve of class $C^{1+\varepsilon}$ 
which is flat enough in $B(0,10r_0)$, but let us try to keep things simple here.

Next we assume that $0 \in E$, and that we have a sliding minimal cone $X$, also
associated to the sliding boundary $L$, which is close enough to $E$ in $B(0,10 r_0)$. 
We shall systematically measure such things with the local normalized variants $d_{x,r}$ of the 
Hausdorff distance between sets, defined by
\begin{equation}\label{1.13}
d_{x,r}(E,F) = \frac{1}{r} \sup_{y\in E \cap B(x,r)} \dist(y,F)
+ \frac{1}{r} \sup_{z\in F \cap B(x,r)} \dist(y,E)
\end{equation}
when $E$ and $F$ are (nonempty, and most of the time closed) sets, $x\in \R^n$, and
$r > 0$. By convention $\sup_{y\in E \cap B(x,r)} \dist(y,F) = 0$ when $E \cap B(x,r)$ is
empty, and similarly for $\sup_{z\in F \cap B(x,r)} \dist(y,E)$. This distance does not exactly
satisfy the triangle inequality, but it localizes well and is very convenient to use.
So we assume that 
\begin{equation}\label{1.16b}
d_{0,10r_0}(E,X) \leq \varepsilon_0.
\end{equation}
We would love to prove that under these assumptions and if $\varepsilon_0$ is small enough, 
$E$ is $C^{1+\beta_1}$-equivalent to $X$ near $B(0,r_0)$, say. By this we mean that there is a 
constant $\beta_1 > 0$ (that depends only on $n$, $\beta$, and maybe on $X$) and a
$C^{1+\beta_1}$ diffeomorphism $\Phi : \R^n \to \R^n$, such that $\Phi(L) = L$, and
\begin{equation}\label{1.16a}
E \cap B(0,r_0) = \Phi(X) \cap B(0,r_0).
\end{equation}
Usually we also require a uniform control on the $C^{1+\beta_1}$ (uniform) norms of $\Phi$
and $\Phi^{-1}$ and that for some $\eta > 0$ that can be chosen as small as we want in advance
(and then $\varepsilon_0$ will depend on $\eta$),
\begin{equation}\label{1.17a}
|\Phi(x)-x| \leq \eta r_0 \text{ and } (1-\eta)|x-y| \leq |\Phi(x)-\Phi(y)| \leq (1-\eta)|x-y|
\ \text{ for } x, y \in \R^n,
\end{equation}
so that \eqref{1.16a} is just a little weaker than requiring that $E \cap B(0,r_0) \subset \Phi(X)$ and
$\Phi(X \cap B(0,2r_0)) \subset E$, which we could get with the same proof.

We shall see soon that the situation can be more complicated than this, depending on the 
approximating minimal cone $X$.

We start our discussion with the simplest case when $X$ is a \ub{half plane bounded by $L$}.
In this case we have the following perfect analogue of J. Taylor's theorem in \cite{Ta}.

\begin{thm}\label{t1.4a}
Let $E$ and $X \in \bH(L)$ satisfy the assumptions above, and in particular 
\eqref{1.11a}, \eqref{1.12a}, \eqref{1.13a}, and \eqref{1.16b}. If in addition $\varepsilon_0$
is small enough, depending on $n$, $C_h$, $\beta$ and $\eta$, then 
$E$ is $C^{1+\beta_1}$-equivalent to $X$ near $B(0,r_0)$, for some $\beta_1$ which 
depends only on $n$ and $\beta$.
\end{thm}

Let us comment on this statement before we go to more complicated cases.  
See Theorem~\ref{t29a.1} for a slightly more general statement and then the proof 
(given the rest of the paper).
Notice that the $C^{1+\beta_1}$-equivalence just means that near $B(0,r_0)$, 
$E$ is a $C^{1+\beta_1}$ surface bounded by $L$, and in fact a Lipschitz graph 
over the half plane $X$, with a small Lipschitz constant.
It may appear in our statements that we use the Reifenberg parameterization theorem 
from \cite{R}, or one of its later variants, but here we play in the $C^{1+\beta_1}$ category, 
where this theorem is much easier to prove than in its original, less regular setting. 
That is, here it essentially amounts to checking that there is a tangent plane $P(x)$ to $x$ 
at every point of $E \sm L$, and that the direction of $P(x)$ is H\"older continuous on $E$. 
In later statements, it would be H\"older continuous on each face of $E$.

Theorem \ref{t1.4a} is an extension of Corollary 1.7 on page 344 of \cite{Mono}, 
which essentially proves the same thing with a biH\"older equivalence only. 
The big difference is that we now prove an additional decay estimate on some quantities 
that measure closeness to planes or half planes. 
We shall discuss the proof ingredients in the next subsection. 

Notice also that Theorem \ref{t1.4a}, and already Corollary 1.7 in \cite{Mono}), say something on the 
topology of $E$ in $B(0,r_0)$: it has no holes or bubbles, and it stays attached to $L$ in the simple way
that one would expect. Plus we have some metric estimates on these properties.

Finally observe that Theorem \ref{t1.4a} implies a weaker statement with blow-up limits.
That is, if $E$ satisfies \eqref{1.11a} and \eqref{1.12a}, $L$ is a line through the origin 
(but a smooth curve would work too), and if one of the blow-up limits of $E$ at $0 \in E$ is 
a half plane $X \in \bH(L)$, then $E$ is $C^{1+\beta_1}$-equivalent to $X$ in some small ball 
centered at $0$, and in particular $E$ has a tangent cone (a unique blow-up limit) at $0$ equal to $X$. 
This is easy to check: just apply the theorem in a small enough ball, where \eqref{1.16b} holds.

\ms
Our next case is not really new, in the sense that it concerns the same minimal cones that were known
to work away from the boundary. Suppose that $X$ is a plain minimal cone that satisfies the 
\ub{full length property}. By plain, we mean with no boundary condition, or equivalently with $L = \emptyset$,
and the full length property is the sufficient condition given in \cite{C1} for the J. Taylor theorem to be 
satisfied for $X$; let us not give the definition for the moment, but only recall that the cones of
$\bP \cup \bY \cup \bT$ (i.e., the cones of type $\bP$, $\bY$, and $\bT$ as above) satisfy this. 
We say that a cone $X$ is \ub{fully transverse to $L$} when $X \cap L = \{ 0 \}$. 

\begin{thm}\label{t1.5a}
Let $X$ be a plain minimal cone that satisfies the full length property of \cite{C1} and is 
fully transverse to $L$, and suppose that $E$ and $X$ satisfy the assumptions above, 
and in particular \eqref{1.11a}, \eqref{1.12a}, \eqref{1.13a}, and \eqref{1.16b}. 
If in addition $\varepsilon_0$ is small enough, depending on $n$, $C_h$, $\beta$, $X$, 
and $\eta$, then $E$ is $C^{1+\beta_1}$-equivalent 
to $X$ near $B(0,r_0)$, for some $\beta_1$ which depends only on $n$, $X$, and $\beta$, but where
we no longer require that $\Phi(L) = L$ in the definition of $C^{1+\beta_1}$-equivalent.
\end{thm}

The dependence on $X$ is through some angles and the full length parameters, but we do not worry
too much about it because in practice we can discretize (i.e., use a finite number of cones $X$).
See Section \ref{S36.1} for a slight extension and the proof, which consists in reducing to
the case when there is no boundary. 
We cannot require that $\Phi(L) = L$ here, because $E$ could be a set of type $\bT$, for instance,
with a center very close to $0$, but not on $L$.
Other cases of this type are treated in Section \ref{S36}, but let us return to the main simple cases.

\ms
The last case that works almost perfectly is when $X$ is a \ub{generic cone of type $\bV$},
by which we mean that the two half planes $H_1$ and $H_2$ that compose $X$
(as in the definition above) make an angle $\alpha \in (\2,\pi)$, thus excluding planes
and what we shall call sharp $\bV$-sets.  

\begin{thm}\label{t1.6a}
Theorem \ref{t1.4a} is still true when $X$ is a generic cone of type $\bV$, but now 
$\varepsilon_0$ depends also on the angle of the two half planes $H_i \in \bH$ that compose $X$.
\end{thm}

\begin{figure}[!h]  
\centering
\includegraphics[width=6.cm]{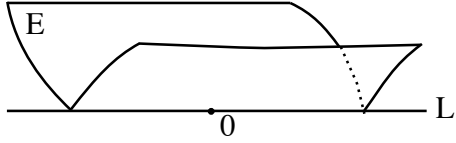}
\caption{$E$ near a generic cone of type $\bV$ \label{f1.2}}
\end{figure}

See Theorem \ref{t29b.1}.
Forgetting about the complicated mapping $\Phi$, the conclusion just means that near $B(0,r_0)$,
$E$ is composed of two faces $F_1$ and $F_2$ bounded by $L$, and that each $F_i$ is a
$C^{1+\beta_1}$ and Lipschitz graph over the corresponding half plane $H_i$ of $X$, with a Lipschitz
norm which is as small as we want, provided that we take $\varepsilon_0$ accordingly small.
In particular the two faces $F_i$ meet ``transversally'', with angles that are as close to $\alpha$
as we want. However, we do not say that $\Phi$ is conformal along $L$, or in simpler terms the angle 
that the two $F_i$ make at $x\in L$ is allowed to depend on $x$, although in a slow, H\"older way.
See Figure~\ref{f1.2} for a hint of what $E$ looks like in $B(0,r_0)$. 

% p15

As in the case of $\bH$, Theorem \ref{t1.6a} also contains topological information on $E$ 
that was not obvious a priori. In fact, the theorem excludes some behaviors that could have been 
considered possible, such as the behavior that is described for sharp $\bY$-cones below. 
Here the methods of \cite{Mono} are no longer enough, because the slow variation of the 
approximating minimal cones, which follows from the decay estimates in the present paper, 
seem to be needed to exclude these behaviors.

As before, the theorem implies that if $E$ satisfies \eqref{1.11a} and \eqref{1.12a}
and one of the blow-up limits of $E$ at $0 \in E$ is a generic cone $X$ of type $\bV$, then $E$ is 
$C^{1+\beta_1}$-equivalent to $X$ near $0$, and in particular $X$ is the tangent cone to $E$ at $0$.

\begin{figure}[!h]  
\centering
\includegraphics[width=12.cm]{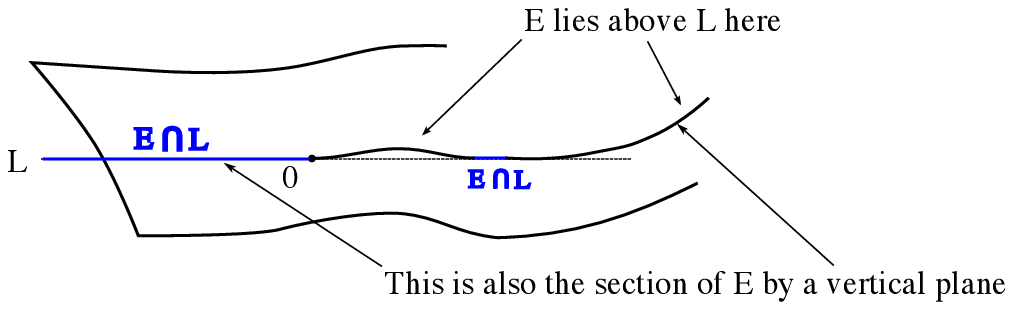}
\caption{ Behavior of $E$ near a plane through $L$ \label{f1.3}}
\end{figure}

\ms
The first case where we get into some (moderate) trouble is when $X$ is \ub{a plane that}
\ub{contains $L$} (and this is why we required $X$ to be transverse in Theorem \ref{t1.5a}). 
In the present case $E$ may be attached to $L$ along just about any closed subset of $L$, 
and not meet the rest of $L$. 
Along this set, $E$ may have a crease, i.e., have different tangent half planes, as depicted
by Figure~\ref{f1.3}. 

% mis page 15

That is, let $X$ be a plane that contains $L$, and assume otherwise that $E$ satisfies 
the same assumptions as for Theorems \ref{t1.4a}-\ref{t1.6a}. 
We claim that if $\varepsilon_0$ is small enough,
depending on $n$, $C_h$, $\beta$ and $\eta$, we have the following description of $E$ in
$B(0,r_0)$. First, there is an $\eta$-Lipschitz function $\psi : P \to P^\perp$ such that
$E \cap B(0,r_0) = \Gamma_\psi \cap B(0,r_0)$, where $\Gamma_\psi$ denotes the graph of $\psi$.
In addition, $\psi$ is of class $C^{1+\beta_1}$ on $(P \cap B(0,r_0)) \sm (E \cap L)$, with a uniform 
H\"older estimate for $\nabla \psi$ on $P \cap B(0,r_0) \sm (E \cap L)$ (with the geodesic distance).
Thus $\psi = 0$ on $L \cap E \cap B(0,r_0) \subset P$, and it has half derivatives from both 
accesses along $E \cap L$, but that may be slightly different from each other at interior points 
of $E \cap L$. And near interior points of $E \cap L$ where $\nabla\psi$ had two different limits, 
$E$ can be described by Theorem \ref{t1.6a}. 
See~\ref{t29c.1} for more details. 

\begin{figure}[!h]  
\centering
\includegraphics[width=12.cm]{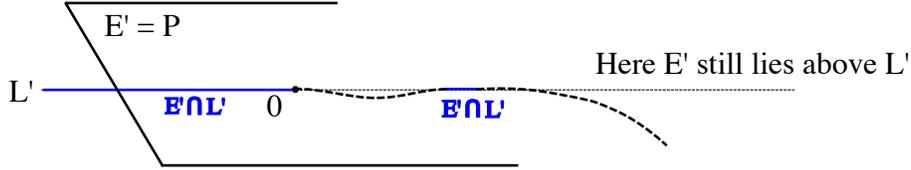}
\caption{ The set $E'=\Phi(P)$ is minimal because it is a plane \label{f1.4}}
\end{figure}

Notice that here $E$ is topologically the same as $P$, but not always in the $C^1$ category.
The description above is not shocking. Consider a nice deformation $\Phi$ of $\R^n$
that moves points downwards a little, sends the set $E$ depicted in Figure \ref{f1.3} 
to the plane $E' = P$, and $L$ to a new boundary $L'$ that coincides with $L$ on $E \cap L$; we 
know that $P$ is minimal, also with the sliding boundary $L'$, and we expect that $E = \Phi^{-1}(P)$
will stay almost minimal if $\Phi$ sufficiently flat. See Figure~\ref{f1.4}.

 % p16

We have the same sort of result when $X$ is cone of type $\bY$ or $\bT$, and one of the two half lines
(say $L_+$) that compose $L \sm \{ 0 \}$ is contained in the interior of a face of $X$.
In this case, we have the same description as in the previous case on a small open cone around $L_+$.
On the rest of $\R^n$, we can proceed as in the transverse case above. The argument also works
when $X$ is a sliding minimal cone that satisfies the full length property, where this time the full length 
property is as in Definition \ref{t3.1} below. See Section \ref{S36.2}.

\begin{figure}[!h]  
\centering
\includegraphics[width=12.cm]{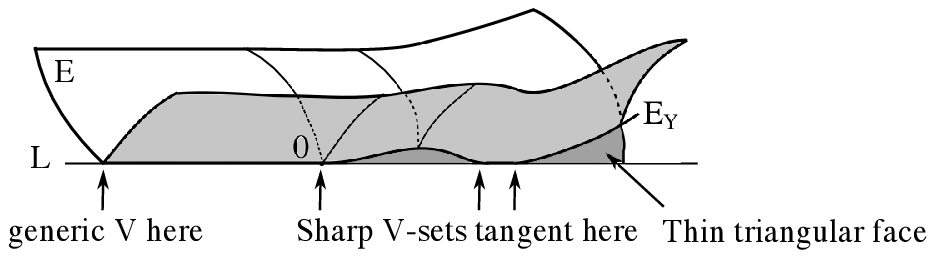}
\caption{ The set $E$ near a sharp $\bV$ set \label{f1.5}}
\end{figure}

\ms
Our next case is when $X$ is a \ub{sharp $\bV$-cone}, which means that $X = H_1 \cup H_2$
for some half planes $H_1, H_2 \in \bH$ that make a $\2$ angle. See Theorem \ref{t29d.1}.
The difference with the generic case is that now $E$ can partially detach itself from $L$, 
along a curve of $\bY$-points of $E \sm L$, as suggested by Figure~\ref{f1.5}. 

That is, assume now that $X$ is a sharp $\bV$-cone and that the other assumptions of 
Theorems~\ref{t1.4a}-\ref{t1.6a} are satisfied. We claim that in $B(0,r_0)$, we have the following description of~$E$. 

First, there is a curve $\gamma$, which is the graph of some function $g : L \to L^\perp$ that is
both $\eta$-Lipschitz and of class $C^{1+\beta_1}$, such that every point of 
$\gamma \cap B(0,2r_0) \sm L$ lies in the set $E_Y$ of points of $E \sm L$ that are type $\bY$.
This means, points $x\in E \sm L$ where $E$ is tangent to a $\bY$ set.
And at points $x\in \gamma \cap B(0,2r_0) \cap L$, $x$ has a tangent cone $V(x) \in \bV(L)$,
which may be generic (at interior points of $\gamma \cap L$), but always with an angle close to
$\2$. The curve $\gamma$ will play the role of a spine for $E$ that splits $E \cap B(0,r_0)$
into three faces that we try to describe now.

 % p17

Denote by $e_i$, $i = 1, 2$, the unit vector in $H_j$ that is orthogonal to $L$, set 
$e_3 = - (e_1+e_2)$, and denote by $H_3$ the half plane bounded by $L$ and pointing 
in the direction of $e_3$. 
Thus $H_1$, $H_2$, and $H_3$ would form a $\bY$-set.
Then denote by $P_i$, $1 \leq i \leq 3$, the plane that contains $H_i$, and by $\pi_i$ the orthogonal
projection on $P_i$. There are three sets $A_i$, $1 \leq i \leq 3$, with $A_i \subset H_i$, and
three functions $\Psi_i : A_i \to P_i^\perp$ which are both $\eta$-Lipschitz and of class 
$C^{1+\beta_1}$, so that if $F_i$ denotes the graph of $\psi_i$, then 
$E \cap B(0,r_0) = (F_1 \cup F_2 \cup F_3) \cap B(0,r_0)$.
The two faces $F_1$ and $F_2$ are $C^{1+\beta_1}$ surfaces bounded by $\gamma$,
which means in particular that for $i = 1, 2$, $A_i$ is the closure of the component of 
$H_i \sm \pi_i(\gamma)$ that leaves furthest from $L$. 
In the simplest case when $\gamma$ just leaves $L$ on one side, the face
$F_3$ looks like a thin triangular wall that connects $L$ to $E_Y$, but in general it may have 
infinitely many connected components. 
The face $A_3$ is bounded by $\pi_3(\gamma \sm L)$ on one side, and by $L \sm \gamma$ 
on the other side, and $F_3$ is bounded by $\gamma \sm L$ on one side, and $L \sm \gamma$ 
on the other side.

\begin{figure}[!h]  
\centering
\includegraphics[width=11.cm]{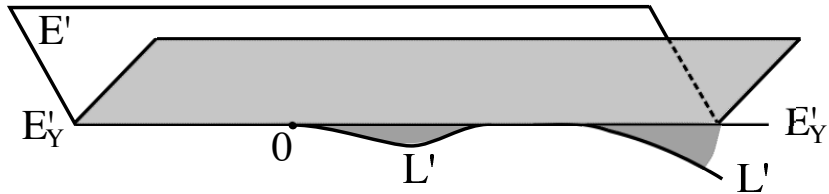}
\caption{ The image $E' = \Phi(E)$ is  probably minimal  because it is a truncated $\bY$ set \label{f1.6}}
\end{figure}

Hopefully this description (together with Figure~\ref{f1.5}) 
gives a good idea of what $E$ looks like near $0$. Another way to see it would be to say that
$E$ is $C^1$-equivalent to a set of type $\bY$, but truncated by the line $L$. This is also why
the description above looks logical: we could deform the set of Figure~\ref{f1.5}, 
a little as we suggested in Figure~\ref{f1.4}, 
by a nice mapping $\Phi$ that sends $E$ to a subset $E'$ of a cone of type $\bY$, 
but truncated by the curve $L' = \Phi(L)$. It is not too hard to believe that $\Phi(E)$ is 
minimal with the sliding boundary $L'$, and that if $\Phi$ is nice enough, $E$ is still sliding almost minimal. See Figure~\ref{f1.6}.

% p 18

A priori there was a possibility that the sort of behavior described here near sharp $\bV$ 
sets could also happen near generic $\bV$ sets, or even planes that contain $L$.
This was apparently suggested in \cite{B}, % leave this
but we claim that this does not happen.

This was the most interesting case that we can treat for the moment. 
Notice that this time $E$ does not even have the same topology as its model $X$. 
All this will be discussed a little more and proved with Theorem \ref{t29d.1}.

\begin{rem}\label{r1.6}
In the descriptions above, the fact that we consider general sliding almost minimal sets helped us 
claim that we probably have the right description, but this hides the fact that when $E$ is sliding minimal, 
there is probably some additional rigidity in the problem, that the author does not understand at all,
but that prevents the most complicated behaviors described above 
(near planes that contain $L$ and sharp $\bV$ sets) to occur. 
That is, planes could become $\bV$ sets and sharp $\bV$ sets could generate 
a curve $\gamma$ of points of type $\bY$ that leaves $L$, but no complicated 
hesitating limit sets would occur.
\end{rem}

Now we turn to the main case that we do not control, which is when $X \in \bY(L)$ is a cone of type
$\bY$ with a spine equal to $Y$. Assume that one of the blow-up limits of $E$ at $0$
is the cone $X \in \bY(L)$, and try to study $E$ near $0$.

The cone $X$ has the full length property (as in Definition \ref{t3.1}; see Theorem \ref{t30.1}), so we'll see in Theorem \ref{t22.2n} that $E$ has a tangent set of type $\bY$ at the origin, 
but nonetheless we do not have enough control on balls that are not centered on $L$,
to give a good description of $E$ near $0$. In this case, we expect that $E$ looks like (and in particular
has the topology of) a cone of type $\bY$, but with little creases (like in the case when $X$ was a plane
that contains $L$) along parts of $L$, but at this point we cannot exclude other, less beautiful options.

\begin{figure}[!h]  
\centering
\includegraphics[width=16.cm]{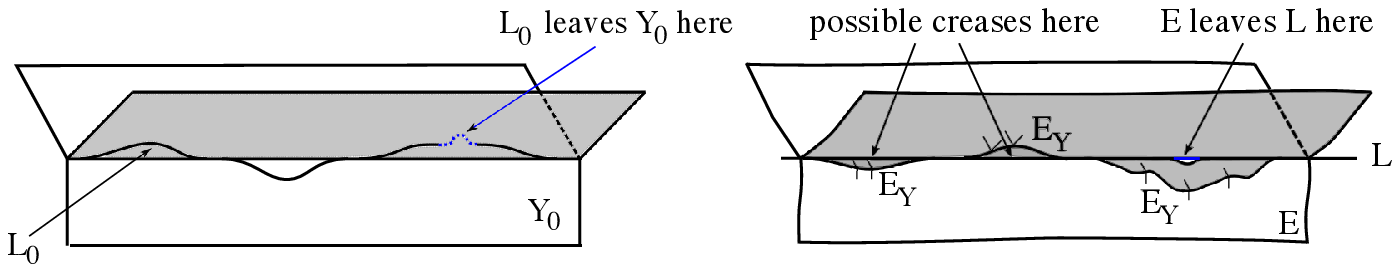}
\caption{Left: A minimal set $Y_0$ and a boundary curve.
Right: The sliding almost minimal set $E = \Phi(Y_0)$. \label{f1.7}}
\end{figure}

\begin{figure}[h]  
\centering
\includegraphics[width=7.cm]{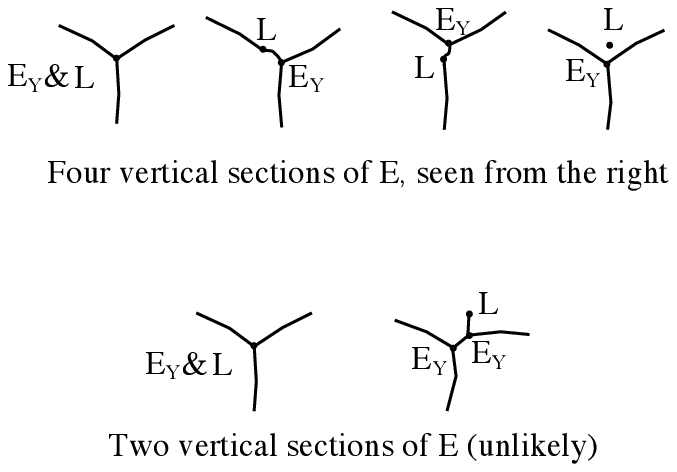}
\caption{Four sections of $E$ and two less probable sections \label{f1.8}}
\end{figure}

If the author had to guess the behavior of $E$ near $0$, he would start from a set $Y_0$
of type $\bY$, then draw a curve $L_0$ tangent to $E$ at $0$, not necessarily entirely drawn on
$Y_0$ (but it is more fun to travel on the various faces of $Y_0$). 
See the strongly exaggerated  Figure~\ref{f1.7} (left). 
Then he would send $Y_0$ and $L_0$ to $E = \Phi(Y_0)$ and hope that if $\Phi$ is gentle
enough, $E$ is still almost minimal with the sliding boundary $L = \Phi(L_0)$ depicted grossly
by Figure~\ref{f1.7} (right). 
The second hope is that nothing worse than that ever happens with sliding almost minimal sets. 
Figure~\ref{f1.8} shows four successive sections of our candidate $E$, 
and (below) two sections that could a priori exist, coming from a more complicated structure
of $E$ (but we hope not). We return to this in Section \ref{S35}.

 % p18

 % p19

Let us also say why the author believes that this is the main bad case.
Of course there are other cones $X$ for which we have the same problem. 
For instance, $X$ could be a set of type $\bT$, with a spine that contains a half of $L$. 
But such cases should be similar, in the sense that if we understand the case of $X \in \bY(L)$, 
we can probably deal with these other cases by restricting to cones around a half of $L$ first.

Now there are possibly many other cones $X$ that one should consider, but fortunately the 
points $x\in E \cap L$ where $x$ has an exotic blow-up limit $X$ like this are isolated, so
even though we may not have a very good control on $E$ precisely at those points 
(especially if $X$ does not satisfy the full length property), we can probably still get some control, 
by restricting to concentric annuli where there is no exotic point.
That is, the author believes that $X \in \bY(L)$ is the most complicated case because it may happen on a 
large set. See Section \ref{S36.3} for a slightly longer discussion.

\subsection{Decay for the density excess and approximation by cones}

Let us now describe elements of the proofs and estimates that lead to the results above.
In addition to the general regularity results of \cite{Sliding} that were mentioned
above, the key ingredient in the proofs will be related to the monotonicity of density, 
or a variant that will be discussed soon. This is not so different in spirit from what is was done 
far from the boundary, in \cite{Ta} and \cite{C1}.

Let $E$ be a coral sliding $(U, L, h)$-almost minimal set, and assume to simplify the discussion that
$L$ is a line through the origin. Define, for $x\in E$ and $r > 0$ such that $B(x,r) \subset U$
(we shall not need the other pairs) the density
\begin{equation}\label{1.1}
\theta(x,r) = r^{-2} \H^2(E \cap B(x,r)). 
\end{equation}
The local Ahlfors regularity of $E$ says that $C^{-1} \leq \theta(x,r) \leq Cr$ when 
$x\in E$ and $B(x,2r) \subset U$, and it is proved in \cite{Sliding} that $\theta(x,\cdot)$ is 
nondecreasing when $x\in L$ and $E$ is sliding minimal.
When $x\in L$ and $E$ is merely sliding almost minimal (but $h$ satisfies a Dini condition),
$\theta(x,\cdot)$ is still nearly nondecreasing; see Theorem 28.7 in \cite{Sliding}, quoted
as \eqref{18.10} below, for a precise estimate. The basic idea for the proof is the same 
as in the standard case, which is to compare $E \cap B(x,r)$ to the cone (centered at $x$) 
over $E \cap \d B(x,r)$; 
the fact that this cone is a limit of competitors for $E$ is still true here, because the deformations
that we generally use to prove this are radial and $L$ is a cone. 
This is the reason why we require $x$ to lie in $L$ for the near monotonicity property; 
we shall return to this issue below.

Because of the near monotonicity, the limit
\begin{equation}\label{1.2}
\theta(x) = \lim_{r \to 0} \theta(x,r)
\end{equation}
exists. Our main ingredient for the control of $E$ on balls that are centered on $L$ is a decay result for 
the density excess $f$ that we define now. Suppose that $0 \in E \cap L$, and set
\begin{equation}\label{1.3}
f(r) = \theta(x,r) - \theta(x) \ \text{ for $0 < r \leq \dist(0,\R^n \sm U)$.}
\end{equation}
Here we say that $f$ decays like a power, as soon as 
$h$ is small enough and $E$ is close enough to a good minimal cone.

\begin{thm}\label{t1.3}
Let $L$ be a line through the origin, $U$ an open set in $\R^n$, 
$r_1 > 0$ such that $B(0,r_1) \subset U$, and $E$ a coral sliding $(U, L, h)$-almost minimal set, 
with a gauge function $h$ such that $h(r) \leq C_h r^\beta$ for $0 < r \leq r_1$.
There exist constants $\varepsilon_0 > 0$ and $a \in (0,1)$, 
that depend only on $n$ and $\beta$, such that if in addition $C_h r_1^{\beta} \leq \varepsilon_0$ and 
there is a minimal cone $X$ (centered at $0$), of type $\bP$, $\bY$, $\bT$, $\bH(L)$, or $\bV(L)$, 
such that 
\begin{equation}\label{1.15}
\H^2(X \cap B(0,1)) = \theta(0) := \lim_{r \to 0} \theta(0,r)
\end{equation}
and
\begin{equation}\label{1.16}
d_{0,r_1}(E,X) \leq \varepsilon_0,
\end{equation}
then 
\begin{equation}\label{1.17}
f(r) \leq 10^{-10} (r/r_1)^{a} \ \text{ for } 0 \leq r \leq r_1/2.
\end{equation}
\end{thm}

Of course $10^{-10}$ could be replaced with any small constant, but $\varepsilon_0$ would have 
to be made even smaller.

In fact, there is a notion of (sliding) full length property for sliding minimal cones, that will be
explained in Section \ref{S3b} (see in particular Definition \ref{t3.1}),
and Theorem \ref{t1.3} remains valid for any minimal cone $X$ that satisfies this full length property
(and satisfies \eqref{1.15} and \eqref{1.16} as above). Then $\varepsilon_0$, $C_2$, 
and $a$ depend also on $X$ through its full length parameters. It just turns out that the 
standard cones mentioned above all satisfy the full length property
(see Theorem \ref{t30.1}), so that Theorem \ref{t1.3}
follows from its generalization, Theorem \ref{t22.2n} below.  
See Section \ref{S30} for the final steps of the verification of full length for the standard cones.

Here and below, we just found it easier to say that our constants depend on $n$, rather than trying
to check whether this is really true. 

We like the density excess $f$ because it decays and at the same time controls the geometry of $E$.
We give the basic consequence here, comment and explain some ideas about the proofs of both results, 
and refer to more specific statements later.

\begin{thm}\label{t1.4} 
Let $U$, $L$, $h$, $E$, and $r_0$ satisfy the assumptions of the previous statement.
Then there is a cone $X_0$ of type $\bP$, $\bY$, $\bT$, $\bH(L)$, or $\bV(L)$, centered 
at the origin, such that $\H^2(X_0 \cap B(0,1)) = \theta(0)$ and 
\begin{equation}\label{1.18}
d_{0,r}(E,X_0) \leq c(\varepsilon_0)  (r/r_1)^{a/4} \ \text{ for } 0 < r \leq r_1.
\end{equation}
Here $\varepsilon_0$ and $a \in (0,1]$ depend only on $n$ and $\beta$, 
and $c(\varepsilon_0)$ depends also on $\varepsilon_0$, but given $n$ and $\beta$, we can 
always choose $\varepsilon_0$ so small that $c(\varepsilon_0)$ is as small as we want.
\end{thm}

As before, there is a similar statement when $X$ is a full length minimal cone, and then $\beta_3$
and $c(\varepsilon_0)$ depend also on the full length parameters for $X$; 
see Theorem \ref{t22.2n}. Both Theorem~\ref{t22.2n} and Theorems \ref{t1.3} and \ref{t1.4}
will be proved in Section \ref{S21} (using the earlier sections).

The presentation of Theorems \ref{t1.3} and \ref{t1.4} as coming one after the other is slightly
misleading; for technical reasons we will need to prove the two of them together, event though
there are two main pieces, Proposition \ref{t16.2} that brings decay for $f$ and Theorem \ref{t18.1}
that gives a geometric control. We will return to this in detail in Section \ref{S21}. 
Only a simpler piece of Theorem \ref{t1.3}, Corollary \ref{t17.2}, will be proved directly from
the decay estimate in Proposition \ref{t16.2}.

Thus $E$ has a unique tangent cone (namely $X_0$) at $0$, of density $\theta(0)$, 
and we even have an estimate on how fast $r^{-1}E$ tends to $X_0$ in the unit ball. 
Of course $X_0$ may be slightly different from $X$, but not so much because they both 
approximate $E$ well in $B(0,r_0/2)$. In the specific case of Theorem \ref{t1.4}, it is even 
of the same type as $X$, because the types are determined by a finite number of densities.

In both statements we required the density of $X$ to match the density $\theta(0)$ of $E$;
if instead $\H^{2}(X\cap B(0,1)) > \theta(0)$, and even in the plain case (with no boundary), 
that fact that the density $\theta(0,\rho)$ may vary a lot between $0$ and $r_0$ seems 
to prevent us from proving any good quantitative estimate. % discuss this again some day?

The precise assumption that $h(r) \leq C_h r^{\beta_0}$ is not vital;
a slightly slower decay, like $h(r) \leq C \big[\ln\big(\frac{1+ r}{r}\big)\big]^{-B}$ for some large $B$,
would be enough to get a roughly similar decay for $f(r)$ and $d_{0,r}(E,X(r))$, but we shall skip the
computations and refer to a similar statement in \cite{C1}, where the computations were done 
that we may always copy.
Also see Section \ref{S31} for a discussion of what happens when $L$ is a smooth curve through the origin, rather than a line.

\ms
We now say a few words about how we intend to prove Theorems \ref{t1.3} and \ref{t1.4}.
Let $E$ and $r_0$ be as in the statements; the estimate \eqref{1.17} will follow from 
a differential inequality like 
\begin{equation}\label{1.25a}
r f'(r) \geq a f(r) - q(r),
\end{equation}
where $a> 0$ is a small constant that depends on the geometry (including the full length constants)
and $q(r)$ is a small error term that contains the contribution of the gauge function $h$.
This inequality will be proved for almost every $r \leq r_0$, say, and then integrated to get \eqref{1.17}. 
See Proposition~\ref{t16.2} for the statement, and Section \ref{S16} for how to derive 
\eqref{1.17} from that statement.

Notice that the near monotonicity comes from a similar statement with $a = 0$. 
This means that when $f(r) > 0$,
we have to improve on our proof of near monotonicity and save a quantity comparable to $f(r)$. 
Recall that for the near monotonicity we essentially compare $E$ with a cone; we will thus have to
find a better competitor than the cone. And indeed the main construction of the paper will be the construction, 
for almost every $r \leq r_0$, of a new competitor for $E$, which is at least as good as the cone over 
$E \cap \d B(0,r)$ and even significantly better if $E \cap \d B(0,r)$ is far from ``optimal''.

There is one basic case where we can do better than the cone, which is when $E \cap \d B(0,1)$
is composed of a simple net of Lipschitz curves with small constants (understand, small 
Lipschitz perturbations of geodesics), but which are not geodesics.
Then the cone is composed of small conic pieces that we can see as graphs of homogeneous functions
defined on triangular sectors, and in this case we can replace these homogeneous functions with harmonic
functions with the same boundary values, and save some surface measure if the Lipschitz curves are far
from geodesics. Here we use the fact that for small Lipschitz graphs, graphs of harmonic functions 
are almost as good as minimal surfaces. 

The next basic case where we can save some area is when there is a net of curves 
contained in $E \cap \d E$ that has some good separation properties, but is more complicated 
than a simple net of Lipschitz graphs with small constants, like the net of curves suggested above. 
In this case, we prove that we can replace $E \cap \d E$ with a simpler net, so that we can still 
use graphs to construct competitors, and moreover save some area when we compare to the cone 
(because the net of Lipschitz curves, even though not entirely contained in $E$, is also somewhat shorter).

We combine these two estimates with a third one, which is a little more surprising, and 
corresponds to the case when $E \cap \d B(0,r)$ is essentially a simple net of geodesics, 
but not necessarily arranged with the same angles and positions as $K = X \cap \d B(0,1)$. 
In this case the only way we found we could do better than the cone over $E \cap \d B(0,r)$ was 
to use competitors of deformations of $X$ and the definition of 
full length, which is the reason why we put it in the assumptions on $X$. 
Fortunately this combinatoric property, which is not unlike the existence of epiperimetric inequalities 
that can be found in the work of Reifenberg, Taylor, and others, is satisfied by the most familiar 
minimal cones.

The construction of a suitable net of curves, and then of competitors for $E$, 
is done in Sections \ref{S4}-\ref{S15},
which are then followed by estimates that lead to Theorems \ref{t1.3} and  \ref{t1.4}, done in 
Sections \ref{S16}-\ref{S21}. 

\ms
All this works well, in a way which is similar to what was done in \cite{C1} in particular, 
and we get good decay estimates and then approximation by minimal cones, but only for balls that
are centered on $L$. But for the classification and regularity results, it seems that we also need
a uniform control on balls that are centered a little off $L$. However, for $x\in E \sm L$, 
the density function $\theta(0,r)$ defined by \eqref{1.1} is no longer nondecreasing in general, 
even when $E$ is minimal. For instance, $E$ could be a half plane bounded by $L$ and that contains $x$,
in which case $\theta(x,r) = \pi$ for $r \leq \dist(x,L)$ and $\lim_{r \to +\infty}\theta(x,r) = \pi/2$.

Because of this, a variant of $\theta$ was introduced in \cite{Mono}, which at least is 
optimally monotone in some simple cases. Suppose that $0 \in E \sm L$, and denote by $S$ 
the shade of $L$, given by 
\begin{equation}\label{1.19}
S = \big\{ y\in \R^n \, ; \, \lambda y \in L \ \text{ for some } \lambda \in [0,1] \big\}.
\end{equation}
The substitute for $\theta(0,r)$ is the slightly larger function $F$ defined by
\begin{equation}\label{1.20}
F(r) = r^{-2} \Big[ \H^2(E \cap B(0,r)) + \H^2(S \cap B(0,r))\Big].
\end{equation}
One of the main points of \cite{Mono} is that when $E$ is a sliding minimal set
on $U \supset B(0,r_0)$, the function $F$ is nondecreasing on $[0,r_0)$; see Theorem 1.2 there.
Similarly, $F$ is nearly monotone when $E$ is a sliding almost minimal with a small enough 
gauge function $h$.

Thus even though $\theta(0,r)$ itself is not always monotone where $r \geq \dist(0,L)$, 
we add an increasing term $r^{-2}\H^2(S \cap B(0,r))$ that improves the situation. 
Of course this property is useful also because there are realistic situations where $F$ is constant,
so we may believe that we didn't add too much. 
Here are two instances of this.
The first one is when $E$ is a half plane bounded by $L$ (and that contains $0$ because 
we assumed that $0 \in E$). The second case is when $E$ a truncated cone 
of type $\bY$, i.e., when $E = \overline{Y \sm S}$, where $Y\in \bY(L)$ is a cone of type $\bY$
centered on $0$ and that contains $L$. In both cases, $F(r)$ is just the (constant) 
density of the completed set (a plane or the cone $Y \in \bY$).

We shall not use the more general dimension $d$ of $E$ that was allowed in \cite{Mono}, 
or the more general form of $L$, and this is rather good because this makes the proof somewhat easier.
But we shall use some of the variants or consequences of Theorem 1.2 in \cite{Mono}, because we
need to know that for sliding almost minimal sets, $F$ is nearly monotone, and that
$E$ is close to a half plane or a truncated cone of type $\bY$ through $L$ whenever $F$ 
is nearly constant. We shall be more specific later, during the proof.

We now state analogues of Theorems \ref{t1.3} and \ref{t1.4} for balls centered at $0 \in E \sm L$. 
We shall only worry here about two cases, when $E$ is close to a half plane or to a $\bV$-set in $B(0,r_0)$, 
and not more complicated sets for which the near monotonicity of $F$ does not really help.

\begin{thm}\label{t1.10a}
Let $L$ be a line that does not contain the origin, $U$ an open set in $\R^n$, 
$r_0 > \dist(0,L)$ such that $B(0,10r_0) \subset U$, and $E$ a coral sliding $(U, L, h)$-almost minimal set, 
with a gauge function $h$ such that $h(r) \leq C_h r^\beta$ for $0 < r \leq 10r_0$.
Also let $H$ denote the half plane bounded by $L$ that contains the origin.
There exist constants $\varepsilon_0 > 0$, $C_4 \geq 1$, and $\beta_4 \in (0,\beta]$, that depend only 
on $n$ and $\beta$, such that if in addition $C_h r^{\beta} \leq \varepsilon_0$ and 
\begin{equation}\label{1.28a}
F(3r_0) \leq \pi + \varepsilon_0
\end{equation}
or
\begin{equation}\label{1.29a}
d_{0,3r_0}(E,H) \leq \varepsilon_0,
\end{equation}
then 
\begin{equation}\label{1.30a}
F(r_1) -\pi \leq \Big(\frac{2 r_1}{r_2}\Big)^{\beta_4} [F(r_2) -\pi] 
+ C_4 C_h r_1^{\beta_4} r_2^{\beta-\beta_4}
\ \text{ for } 0 \leq r_1 \leq r_2 \leq r_0
\end{equation}
and in addition
\begin{equation}\label{1.31a}
d_{0,r}(E,H) \leq c(\varepsilon_0) \Big(\frac{r}{r_0}\Big)^{\beta_4/4} 
+ C_4 \big( C_h r^\beta \big)^{1/4}
\end{equation}
for $\dist(0,L) \leq r \leq r_0$, where $c(\varepsilon_0)$ can be made as small as we want by choosing
$\varepsilon_0$ above small enough (depending on $n$ and $\beta$).
\end{thm}

Here we do not try to control $E$ in $B(0,r)$ for $r < \dist(0,L)$, but this would follow easily from
the regularity far from the boundary, since \eqref{1.31a} for $r = \dist(0,L)$ shows that $E$ lies
close to a plane in $B(0, \dist(0,L))$. See Section \ref{S29a} for this type of argument.

This is a combination of Theorem \ref{t23.1} for \eqref{1.30a} and 
Theorem \ref{t23.5} for \eqref{1.31a}. In turn Theorem \ref{t23.1} comes from 
the differential inequality \eqref{23.15} in Proposition \ref{t23.3}, which will be obtained 
as before by constructing an appropriate competitor.

\ms
For the next result, recall that $\bV(L)$ is the set of unions $V = H_1 \cup H_2$ of two 
half planes bounded by $L$ and that make an angle at least $\2$ along $L$. 
Also, when $0 \notin L$ and $r > \dist(0,L)$, denote by $\bY(L,r)$ the set of cones $Y$ 
of type $\bY$ that are centered on $0$ and contain $L \cap B(0,r)$. 
Finally, for $Y \in \bY(L,r)$, denote by $Y^t$ the truncated
cone $\ol{Y \sm S}$; the truncation is not always perfect outside of $B(0,r)$, 
because the spine of $Y$ may be different from $L$, but all we shall care about is
the intersection with $B(0,r)$, where we neatly remove from $Y \cap B(0,r)$ a sector bounded by $L$ 
and contained in a face of $Y$.

\begin{thm}\label{t1.11a}
Let $L$ be a line that does not contain the origin, $U$ an open set in $\R^n$, 
$r_0 > \dist(0,L)$ such that $B(0,10r_0) \subset U$, and $E$ a coral sliding 
$(U, L, h)$-almost minimal set, with a gauge function $h$ such that $h(r) \leq C_h r^\beta$ 
for $0 < r \leq 10r_0$.
There exist constants $\varepsilon_0 > 0$, $C_5 \geq 1$, and $\beta_5 \in (0,\beta]$, 
that depend only on $n$ and $\beta$, such that if in addition $C_h r^{\beta} \leq \varepsilon_0$,
\begin{equation}\label{1.32a}
F(0) : = \theta(0,0) : = \lim_{r \to 0} \theta(0,r) = \frac{3\pi}{2}
\end{equation}
and 
\begin{equation}\label{1.33a}
d_{0,2r_0}(E,V) \leq \varepsilon_0
\end{equation}
for some set $V \in \bV(L)$, then 
\begin{equation}\label{1.34a}
F(r_1) -\2 \leq \Big(\frac{C_5 r_1}{r_2}\Big)^{\beta_5} [F(r_2) -\2] 
+ C_5 C_h r_1^{\beta_5} r_2^{\beta-\beta_5}
\ \text{ for } 0 \leq r_1 \leq r_2 \leq r_0
\end{equation}
and in addition, for $\dist(0,L) \leq r \leq r_0$ there is a set $Y = Y(r) \in \bY(Y,r)$, such that 
\begin{equation}\label{1.35a}
d_{0,r}(E,Y^t) \leq c(\varepsilon_0) \Big(\frac{r}{r_0}\Big)^{\beta_4/4} + C_4 \big( C_h r^\beta \big)^{1/4}.
\end{equation}
As before, the constant $c(\varepsilon_0)$ can be made as small as we want by choosing
$\varepsilon_0$ small enough (depending on $n$ and $\beta$).
\end{thm}

This time \eqref{1.34a} will come from Theorem \ref{t23.2} and \eqref{1.35a} 
from Theorem \ref{t28.2}, and the differential inequality that leads to Theorem \ref{t23.2} 
will be proved in Proposition \ref{t23.4}.

The statement looks a little strange because \eqref{1.33a} seems to authorize a set $V\in \bV(L)$
with an angle (much) larger than $\2$. But in effect, the fact that the density $\theta(0)$ of $E$ at the
origin is $\frac{3\pi}{2}$ forbids this, and indeed \eqref{1.35a} with $r=r_0$ implies 
that $E$ looks like a truncated $\bY$-set in $B(0,r_0)$. 
This last is not incompatible with \eqref{1.33a} (provided that $V$ is almost
sharp), and the reader should also keep in mind that the present situation is most interesting when 
$\dist(0,L)$ is much smaller than $\varepsilon_0$, so that \eqref{1.33a} only gives a rough idea of what
$E$ looks like in $B(0,r_0)$, while \eqref{1.35a}, at least when $\dist(0,L) \leq r << r_0$, 
is often much more precise. 

It seems that we find out that $V$ should be nearly sharp only after the proof, but we could 
also have guessed this earlier, by proving (as we will do for the proof of Theorem \ref{t1.6a}, 
for instance) that when $E$ is close enough to a generic $\bV$-set, or a plane, 
there is no point of type $\bY$ in $E \sm L$ near $L$. See Lemma \ref{t29b.2}.

The proof of Theorem \ref{t1.11a} is in the same spirit as for balls centered on $L$; 
some inequalities are harder to get because of the off-centered balls, and also we were forced 
to restrict to two simpler situations (in terms of combinatorics) because otherwise the near 
monotonicity of $F$ is too far from optimal. On the other hand the general construction is the same, 
and the combinatorics of the net of curves is simpler. In particular there is a notion of full length 
here too, which will be adressed in Sections \ref{S26} and \ref{S27}. 

We refer to the table of contents for more detail on the plan of the paper.

The author especially wishes to thank the Institut Universitaire de France for its invaluable help during 
the early stages of the preparation of this paper. The pictures were done with Inkscape.

\subsection{Notation that will be used extensively} 

As usual, $C$ is a generic notation for a constant,
often large, and whose value may change from line to line. Similarly, $c$ is a small positive constant;
\pari
$B(x,r)$ is the open Euclidean ball centered at $x$ with radius $r > 0$;
\pari
$\B = B(0,1)$ and $\S = \d B(0,1)$ are the unit ball and sphere;  
$\B_r = B(0,r)$ and $\S_r = \d B(0,r)$;
\pari
$L$ is our sliding boundary. Except in Section \ref{S31}, $L$ is a line, not always through the origin;
\pari
$E$ is our sliding $(U,L,h)$-almost minimal set, with sliding boundary $L$ and gauge function $h$;
\pari
$\H^2$ denotes the Hausdorff measure of dimension $2$;
\pari
$\theta(x,r) = r^{-2} \H^2(E \cap B(0,r))$ see \eqref{1.1}; then
$\theta(r) = \theta(0,r) = r^{-2} \H^2(E \cap B(0,r))$;
\pari
$F(r) = r^{-2} \Big[ \H^2(E \cap B(0,r)) + \H^2(S \cap B(0,r))\Big]$ where $S$ is the shade of $L$; 
see \eqref{1.19} and \eqref{1.20} or later \eqref{22.3};
\pari
$X$ is a sliding minimal cone (centered at $0$), often the one that approximates $E$ well, 
and $K = X \cap \d B(0,1)$;
\pari
$\bH = \bH(L)$, $\bP$, $\bP(L)$, $\bY$, $\bY(L)$, $\bT$ are special sets of minimal cones, 
see Subsection \ref{S1.2};
\pari
$d_{x,r}(E,F)$ is our normalized local Hausdorff distance between $E$ and $F$; see \eqref{1.13};
\pari 
MC(L) is the set of minimal cones with sliding boundary $L$; see above \eqref{2.1};
\pari
$V = V_0 \cup V_1 \cup V_2$ is the set of vertices of $K= X \cap \d B(0,1)$ 
(including artificial ones in $V_2$) in the standard decomposition of Section \ref{S3} (see \eqref{3.5});
\pari
the arcs $\cC_i$, $i\in \cI$, are the geodesics that compose $K$ in that standard decomposition above;
\pari
$\ddist$ is the geodesic distance on the sphere $\S$, 
\pari
$\rho(a,b)$ denotes the geodesic from $a$ to $b$ in $\S$ (see \eqref{3.4}); 
\pari
$v(a,b)$ is the unit vector that gives the direction of $\rho(a,b)$ at $a$; 
\pari
$\eta(X)$ controls the size of the smallest arcs of $X$, or its distance to points $\ell$, see \eqref{3.12};
\pari
$\Angle_a(x,y) = \Angle(v(a,x),v(a,y))$ is an angle of geodesics at $a$; see near \eqref{9.12}.
\pari
$\tau$ and $D_\pm(\tau)$ (small disks where we do surgery) appear in Section \ref{S5}
\pari
$\tau_4$, an extremely small number rather than a constant, appears in \eqref{13.2}
\pari
$\tau_1$ and $\lambda$ are rapidly discussed in Section \ref{S4}, but appears in Section \ref{S7}.
\pari
$\varepsilon$ appears in \eqref{4.3} to measure the distance to $X$, then is used all the time. 
It is chosen extremely small, at the end of the proof.

\section{Minimal cones bounded by a line}
\label{S2}

In this section we give a description of the sliding minimal cones of dimension $2$
in $\R^n$, associated to a sliding boundary $L$ which is a line through the origin.
Even when $n=3$, we do not know the exact list of these minimal cones, 
but the combinatoric description that follows will allow us to construct competitors in a 
fairly unified way. The description here is similar to the description of plain minimal cones
(that is, without a boundary condition) that was given in \cite{Holder}, Proposition~14.1,
and of course we will use its proof. 

So let $L \subset \R^n$ be a line through the origin.
We denote by $MC(L)$ the set of sliding minimal cones of dimension 2, with
sliding boundary $L$. That is, $X \in MC(L)$ if $X$ is a (reduced) sliding minimal set in $\R^n$, 
with sliding boundary $L$, and in addition $X$ is a cone. 

Fix $X \in MC(L)$ and set
\begin{equation} \label{2.1}
\B = B(0,1), \  \S = \d\B, \ \text{ and } K = X \cap \S ;
\end{equation}
we want a description of $K$. Let us give a statement now for future reference. 
If the reader is only interested in the small collection of known minimal cones of
dimension 2 in $\R^3$, he/she can just have a look at the statement, check that it fits
with the obvious decomposition of the minimal cones in question, and go to the next section.

\begin{pro} \label{t2.1}
There is a constant $\eta_0 > 0$, which depends only on the dimension $n$, such that
for each sliding minimal cone $X$ with sliding boundary $L$ (i.e., for $X \in MC(L)$), 
$K = \S \cap X$ is a finite union 
\begin{equation} \label{2.2}
K = \bigcup_{j \in \cJ} \cC_j,
\end{equation}
where the $\cC_j$, $j\in \cJ$, are either great circles or closed arcs of great circles. 
The great circles are disjoint from the rest of $K$, and even
\begin{equation} \label{2.3}
\dist(\cC_j, K \sm \cC_j) \geq \eta_0
\end{equation}
when $\cC_j$ is a great circle. 
The arcs of great circles have disjoint interiors, i.e., they can only meet at a common endpoint. 
No point of $L$ lies in the interior of one of our arcs of great circles (otherwise, we cut the arc in two).
We also have that
\begin{equation} \label{2.4}
\H^1(\cC_j) \geq \eta_0 \ \text{ for $j \in \cJ$, except perhaps when one of the endpoints
of $\cC_j$ lies in $L$.} 
\end{equation}
In addition, if $\ell \in L \cap K$ and $\H^1(\cC_j) < \eta_0$ for some  $\cC_j$ which admits 
$\ell$ as one of its endpoints, then there is at most another $\cC_i$ which admits $\ell$ as one 
of its endpoints, this $\cC_i$ (if it exists) makes an angle larger than $\frac{9 \pi}{10}$ with 
$\cC_j$ at $\ell$, and $\H^1(\cC_i) \geq \eta_0$.
\\
The arcs $\cC_j$ are also far from each other, i.e., 
\begin{equation} \label{2.5}
\dist(\cC_i, \cC_j) \geq \eta_0
\end{equation}
for $i, j \in \cJ$ such that $\cC_i \cap \cC_j = \emptyset$, i.e., 
when they do not share an endpoint, but again with the following possible exception:
if there is one one of the exceptional arcs $\cC_k$ for \eqref{2.4} such that the two endpoints 
of $\cC_k$ are also extremities of $\cC_i$ and $\cC_j$ respectively. 
Then instead we only get that $\dist(\cC_i, \cC_j) = \diam(\cC_k)$ in general, and 
$\dist(\cC_i, \cC_j) = \min(\diam(\cC_k),\diam(\cC'_k))$ if $\cC_i$ and $\cC_j$ are both 
almost half circles and happen to be also separated by an exceptional arc $\cC'_k$ near the antipodes.
\\
Finally, if $i\in \cJ$, $\cC_i$ is an arc of circle, and $a$ is one of the endpoints of $\cC_i$,
then one of the two following things happens:
\begin{equation} \label{2.6}
\begin{aligned}
& a \notin L,  \text{ there are exactly two other arcs of great circle } \cC_j \text{ and }\cC_k
\\ \hskip2cm &
\text{ that meet $\cC_i$ at $a$, and they make $\frac{2 \pi}{3}$ angles with $\cC_i$ at $a$;} 
\end{aligned}
\end{equation}
\begin{equation} \label{2.7}
\begin{aligned}
 a \in L &\text{ and all the other arcs of great circle that meet $\cC_i$ at $a$} 
\\ \hskip3cm &
\text{ make angles at least $\frac{2 \pi}{3}$ with $\cC_i$ at $a$.} 
\end{aligned}
\end{equation}
\end{pro}

\ms
We decided to require the arcs of geodesics not to contain a point of $L$ in their interior.
That is, we force the points of $K \cap L$ to be vertices of our description (that is,
when this is not the case, we just cut the arc at the point of $K\cap L$), unless
they lie on a full great circle. But even in this case, we shall later cut the great circles
into pieces, and we will cut at points of $L$ if we can.

The decomposition of Proposition \ref{t2.1} will some times be called
\underbar{the natural decomposition} of $K$. In the next section we will again cut
some of the arcs into smaller parts to get what we'll call 
\underbar{the standard decomposition}.

\ms
The rest of the section will be devoted to the proof of Proposition \ref{t2.1}, but we start 
with a few comments. By arcs of great circles, we mean geodesic arcs, but a priori 
they may be longer than $\pi$, although this will probably not happen. 
That is, the involved circles are centered at the origin. We write $\H^1(\cC_j)$, 
but we could equally have written $\length(\cC_j)$.
In \eqref{2.6}, $\cC_i$ may be the only arc that ends at $a \in L$, or there may be two,
or three, but no more. And when there are three, their directions at $a$ make 
$\frac{2 \pi}{3}$ angles and lie in a same $2$-plane orthogonal to $L$. 
 
We start the proof of the proposition with a first description of $K$ away from $L$. 
We claim that
\begin{equation} \label{2.8}
\begin{aligned}
&\text{Each $x\in K \sm L$ has a small neighborhood where $K$ coincides with a great circle}
\\
&\text{or a union of three arcs of great circles that start at $x$ and make 
$\frac{2 \pi}{3}$ angles there,}
\end{aligned}
\end{equation}
by the proof of Proposition~14.1 on page 88 of  \cite{Holder}.
Admittedly that proposition was announced when $K$ comes from a plain minimal cone
(with no sliding boundary condition), but the first part where we prove the conclusion 
of \eqref{2.8} only uses this information locally. Since the reader may not recall well how this goes,
let us sketch a rapid argument, which actually uses a little more information but is easier to believe.
First let $Z$ be any blow-up limit of $X$ at $x$; we know that it is a minimal cone
(with no sliding boundary), and since $X$ is a cone, a simple computation 
(that will be done soon in a slightly different context) shows that $Z$ is invariant by translations 
in the direction of $x$. When we look at the description of $Z \cap \S$ given in
Proposition~14.1 of  \cite{Holder}, we see that all the great circles involved in this
description are contained in $2$-planes that contain $x$, and it is easy to see that
$Z \in \bP \cup \bY$ (again, read the arguments below if you have a doubt). 

Suppose first that $Z$ is a plane; the local regularity result 
(of \cite{Ta} or \cite{C1}) says that near $x$, $X$ is a $C^{1+\varepsilon}$
surface, and its tangent plane at $x$ contains the radial direction. It follows
from the implicit function theorem that $K$ is a $C^1$ curve near $x$,
and then we can conclude, either as in Proposition~14.1 of  \cite{Holder}
(by constructing competitors by hand), or by saying that in fact (by the regularity
theory for elliptic PDE) $K$ is $C^2$ near $x$, then has vanishing curvature in the
direction of $\S$ (because the total mean curvature is zero, and $X$ has no
curvature in the radial direction). Thus $K$ is an arc of great circle in the
neighborhood of $x$, when $x\in K \sm L$ and $X$ has a blow-up limit at $x$ which is a plane.
The case when $Z \in \bY$ follows at once, because the regularity theorem says that near $x$, 
$K$ is composed three $C^1$ curves, and we just showed that they are arcs of great circles.
They make $\2$ angles because $Z \in \bY$.
This completes our sketch of \eqref{2.8}.

\ms
So we have a nice local description of $K$ away from $L$, and now we need to see what happens 
near a point of $K \cap L$; we start our study of $K$ near $L$ with a description of sliding
minimal sets of dimension $1$.

\begin{lem} \label{t2.2}
Let $Z$ be a (reduced) sliding minimal set of dimension $1$ in the whole $\R^n$, 
with sliding boundary $\{ 0 \}$. One possibility is that $Z$ is a line or a a set of type $Y$
(i.e., the union of three half lines that meet at a point with $\frac{2 \pi}{3}$ angles).
Otherwise, $0 \in Z$, and $Z$ is either a half line with its end at $0$, or a set of type $V$
(i.e., the union of two half lines with ends at $0$ and that make an angle at least $\frac{2 \pi}{3}$
at $0$), or a truncated $Y$ (i.e, a line segment $[0,a]$ with $a \neq 0$, plus two half lines leaving
from $a$, so that $[0,a]$ and the two half lines make $\frac{2 \pi}{3}$ angles at $a$.
\end{lem}

Let $Z$ be such a minimal set. Away from $0$, and for instance by
Chapter 10 of \cite{Holder}, $Z$ is composed of line segments, that can only meet 
by sets of three, with angles of $\2$, and at vertices that are isolated in 
$\R^n \sm \{ 0 \}$.

The argument that follows is obviously too heavy, as some parts could be replaced by constructions
of competitors with line segments, but hopefully it will convince the reader with less effort.

We may assume that the origin lies in $Z$, because otherwise $Z$ is a plain minimal
set of dimension $1$ (just check the definitions). Those were studied before, 
and they are lines or sets of type $Y$.
Then set $\theta(r) = r^{-1} \H^1(Z \cap B(0,r))$; we know, for instance from 
Section 28 of \cite{Sliding} (but again it is much easier in dimension $1$ because 
we just need to replace with cones over finite sets) that $\theta$ is a nondecreasing function. 
In addition, because of this and a theorem about limits (again \cite{Sliding} is a reference,
but in fact Golab's theorem does the job),
any blow-up limit $Z_0$ of $Z$ at $0$ is a sliding minimal cone
of constant density $\theta(0) = \lim_{r \to 0} \theta(r)$. Similarly, every 
blow-in limit $Z_\infty$ of $Z$ is a sliding minimal cone of constant density 
$\theta(\infty) = \lim_{r \to \infty} \theta(r)$. That is, $Z_0$ and $Z_\infty$ are finite unions 
of half lines emanating from $0$. In fact they can only be composed of $1$, $2$, or $3$
half lines, because a simple argument shows that the half lines make angles $\geq 2\pi/3$
with each other (otherwise, pinch a couple of them near the origin). 

Set $m = \theta(\infty)$; then every blow-in limit limit $Z_\infty$ is composed of $m$
half lines $\ell_i$, $1 \leq i \leq m$, and there are large radii $R$ such that 
$Z \cap \d B(0,R)$ is composed of exactly $m$ points that lie at distances larger than $R$ 
from each other. Indeed, notice that $Z \cap \d B(0,R)$ has at least $m$ points for $R$ 
large, one near each $\ell_i$, because otherwise we could contract a big piece of $Z$ near
$\ell_i$; in addition the presence of an additional point too often would make the density of $Z$
too large. Select such an $R$, and call $x_i$ the point of $Z \cap \d B(0,R)$ that lies close to
$\ell_i$. 

Set $Z_R = Z \cap \overline B(0,R)$, and 
let $C(0)$ denote the connected component of $0$ in $Z_R$.
We want to show that 
\begin{equation} \label{2.9}
C(0) \text{ contains $x_i$ for } 1 \leq i \leq m.
\end{equation}
Let us first assume that $C(0)$ contains none of the $x_i$. For $\varepsilon > 0$ small,
denote by $C_\varepsilon(0)$ the set of points $z\in Z_R$ that can be connected to $0$
by an $\varepsilon$-chain in $Z_R$, i.e., a finite chain of points $\zeta_j \in Z_R$ such that
$\zeta_0 = 0$, $|\zeta_j-\zeta_{j-1}| \leq \varepsilon$ for $j \geq 1$, and $z$ is the last $\zeta_j$.
Since $Z_R$ is a compact set of finite length, it is easy to see that if for every $\varepsilon > 0$
the point $z\in Z_R$ can be connected to $0$ by an $\varepsilon$-chain in $Z_r$, then
there is a path of finite length in $Z_R$ that goes from $0$ to $z$. 
See for instance \cite{Falconer}, or Chapter 30 of \cite{MSBook}. 
In other words, $C(0)$ is the intersection of the $C_\varepsilon(r)$. Since the $C_\varepsilon(0)$
are open in $Z_R$ and $C(0)$ is closed (for instance because the $C_\varepsilon(0)$ 
are also closed), we get that $C(0) = C_\varepsilon(0)$ for some $\varepsilon > 0$, 
and our assumption implies that $C_\varepsilon(0)$ does not contain any $x_i$,
and hence does not meet $Z \sm B(0,R)$ 
(recall that $C_\varepsilon(0) \subset Z_R = Z \cap \overline B(0,R)$).
By compactness, $\dist(C_\varepsilon(0), Z \sm B(0,R)) > 0$. We shall now check that 
this is impossible because it implies the existence of a competitor $Z' = \varphi(Z)$ 
which is strictly better than $Z$.

First observe that if $\varphi : Z \to \R^n$ is Lipschitz, $\varphi(x) = x$ for $x\in Z \sm B$, 
and $\varphi(0)=0$, then $Z' = \varphi(Z)$ is automatically a competitor for $Z$,
because we can interpolate linearly between the identity and $\varphi$ to get a one parameter
family $\{ \varphi_t \}$, and all the mapings $\varphi_t$ satisfy the sliding condition $\varphi_t(0)=0$.
See Definition \ref{t1.1}.

Now we define $\varphi$ on $Z$ by $\varphi(x) = 0$ for $x\in C_\varepsilon(0)$ and
$\varphi_t(x)=x$ on the rest of $Z$. Notice that the rest of $Z$ lies at positive distance
from $C_\varepsilon(0)$, so $\varphi$ is Lipschitz. It is easy to see that $Z' = \varphi(Z)$
does better than $Z$, because we simply removed the measure of $C_\varepsilon(0)$
which contains $Z \cap B(0,\varepsilon)$, whose measure is positive because $0 \in Z$ and $Z$ is Ahlfors-regular.
So $C(0)$ contains at least one $x_i$.

Now suppose that (for instance) the connected component $C(1)$ of $x_1$
does not contain $0$ or any other $x_i$. Let $C_\varepsilon(1)$ denote
the set of points $z\in Z_R$ that can be connected to $x_1$ by an $\varepsilon$-chain
in $Z_R$. As before, $C_\varepsilon(1)$ is both open and closed in $Z_R$,
$C(1)$ is the intersection of the $C_\varepsilon(1)$, and for $\varepsilon$ small enough
$C_\varepsilon(1)$ does not contain $0$ or any other $x_i$ and stays at positive distance
from the rest of $Z_R$.

This time we define a competitor $Z' = \varphi(Z)$ with a function $\varphi$ such that
$\varphi(x)=x$ on $Z \sm B(0,R)$ and on $Z_R \sm C_\varepsilon(1)$, and 
$\varphi(x) = x_1$ on $C_\varepsilon(1)$. For the verification, first observe that
$\varphi(0)=0$ because $0 \in Z_R \sm C_\varepsilon(1)$, 
and that it is enough to check the Lipschitz property
of $\varphi$ near $x_1$ (because $C_\varepsilon(1)$ is far from the rest of $Z_R$).
This is easier if we observe that we could choose $R$ such that near each $x_i$,
$Z$ is in fact a line segment that crosses $\d B(0,R)$ transversally. Indeed almost every
$R$ is like this, because the set of vertices for $Z$ is at most countable, and by Sard's theorem
(to exclude segments that are tangent to $\d B(0,R)$).
Now the Lipschitz property is easy, and it is also clear that $Z' = \varphi(Z)$ is better than $Z$,
because we contract at least a segment to $x_1$.

If \eqref{2.9} fails, we are in one of the following situations. 
Suppose for the sake of definiteness that $C(0)$ does not contain $x_1$.
Recall that it contains some $x_i \,$; without loss of generality we can assume that
$x_2 \in C(0)$. But $x_1$ must be connected to some point, and the only choice left
is $x_3$ (recall that there are at most three $x_i$).
Let us now say why this is impossible. As before, if $\varepsilon$ is small enough,
the set $C_\varepsilon(0)$ defined above coincides with $C(0)$, and thus contains
$0$ and $x_2$, and the set $C_\varepsilon(1)$ coincides with $C(1)$ and therefore
contains $x_1$ and $x_3$, but lies at distance at least $\varepsilon$ from $C_\varepsilon(0)$.

We now define $\varphi$ and the competitor $Z' = \varphi(Z)$ as follows.
As usual, we take $\varphi(x) = x$ on $Z \sm B(0,R)$. On $Z_R \sm C_\varepsilon(1)$,
we let $\varphi$ coincide with a Lipschitz retraction from $\overline B(0,R)$ onto the
line segment $[0,x_2]$. Finally, on the rest of $Z_R$, that is, on $C_\varepsilon(1)$,
we let $\varphi$ coincide with a Lipschitz retraction from $\overline B(0,R)$ onto the
line segment $[x_1,x_3]$. Notice that $\varphi(0) = 0$ because $0\in Z_R \sm C_\varepsilon(1)$.
Again the Lipschitz property of $\varphi$ only needs to be checked near the $x_i$, where we 
just need to know that $[0,x_2]$ and $[x_1,x_3]$ are transverse to $\d B(0,R)$.
Finally, $\H^1(Z' \cap B(0,R)) \leq \H^1([0,x_2] \cup [x_1,x_3]) \leq R + |x_1 - x_3|
< 29 R/10$ if $R$ is large enough, because the $x_i$ lie close to the minimal cone $\Sigma_\infty$
and thus almost make angles of $2\pi/3$. On the other hand, 
$\H^1(Z \cap B(0,R))$ tends to $\theta(\infty) R = 3R$ when $R$ tends to $+\infty$,
so $\H^1(Z' \cap B(0,R)) < \H^1(Z \cap B(0,R))$, $Z'$ is a better competitor, and
this contradiction proves \eqref{2.9}.

Now the set $Z_R$ contains a connected set that connects $0$ and the $x_i$.
This implies (because $\H^1(Z_R) < +\infty$; see again \cite{Falconer} 
or Chapter 30 of \cite{MSBook}) 
that there is a simple arc $\xi_1$ in $Z_R$
that goes from $0$ to $x_1$. If $m\geq 2$, there is also an arc in $Z_R$ that goes
from $x_2$ to $0$; we call $\zeta_2$ the first point of this arc (leaving from $x_2$)
that lies in $\xi_1$. We call $\xi_2$ the portion of this arc between $x_2$ and $\zeta_2$.
Thus $\xi_2$ is essentially disjoint from $\xi_1$, and their union connects $0$, $x_1$, and $x_2$.
If $m=3$, we also find an arc from $x_3$ to $0$, stop it at the first point $\zeta_3$ 
of $\xi_1 \cup \xi_2$, and thus get a third arc $\xi_3$.

Let us assume that $m=3$ (the other cases are simpler). We see that
\begin{equation} \label{2.10}
\H^1(Z_R) \geq \H^1(\xi_1 \cup \xi_2 \cup \xi_3) 
=  \H^1(\xi_1)+\H^1(\xi_2)+\H^1(\xi_3).
\end{equation}
Notice that $\xi_1 \cup \xi_2 \cup \xi_3$ is composed of (at most) five essentially disjoint
curves that connect $0$ and the $x_i$ (in the worse case we cut 
$\xi_1$ in two at $\zeta_2$ and $\xi_2$ or one of the two pieces of $\xi_1$ at $\zeta_3$); 
if we replace each of these arcs with a line segment with the same endpoints,
we get a connected set $F$ such that $\H^1(F) \leq  \H^1(Z_R)$, with a strict inequality
if $Z_R \neq F$.

Denote by $\cal F$ the class of connected unions of at most five line segments contained in
$\overline B(0,R)$, and that contain $0$ and the three $x_i$. Thus $F \in \cal F$.
Let $F_0 \in \cal F$ be such that
\begin{equation} \label{2.11}
\H^1(F_0) = \inf_{G \in \cal F} \H^1(G) \leq \H^1(F) \leq \H^1(Z_R).
\end{equation}
Existence is not an issue, because there are finitely many combinations of intervals,
with endpoints that lie in the compact set $\ol B(0,R)$.

First suppose that $F_0$ has no vertex in $\overline B(0,R) \sm \{ 0 \}$, which means
that $F_0$ is the intersection of $\overline B(0,R)$ with an array of $1$, $2$, or $3$
half lines emanating from $0$. These segments make angles at least $2\pi/3$ at $0$,
because otherwise we may pinch two of them near $0$ and make $F_0$ shorter. We consider 
this good and go to the next case.

Suppose next that $F_0$ has exactly one vertex in $B(0,R) \sm \{ 0 \}$. 
Call this vertex $v$, and observe that the three segments of $F_0$ that leave from $v$
make $2\pi/3$ angles with each other (otherwise, move $v$ a little and this gives a shorter $F_0$).
They either end at points $x_i \in \d B(0,R)$, or at the origin. Call $V_0$ the union of these three
segments; this is a piece of $Y$-set.

Let us first assume that $V_0$ ends at the three $x_i$. One possibility is that $0$ lies 
in $V_0$. Then $F_0 = V_0$ (no need to add anything), we shall consider that $0$ is a vertex
$F_0$ is in fact composed of four segments (three that make a smaller piece of $Y$-set centered
at $v$, and a segment $[0,x_i]$ opposite to it, and this will be a good enough description.
Otherwise, $0$ is also connected to one of the $x_i$ (there is no other inside vertex, and
$v$ already has three segments leaving from it). Notice that $v$ lies very close to the origin,
because the branches of $V_0$ make $\2$ angles, and the $x_i$ are seen from $0$ with
angles that are arbitrarily close to $2\pi/3$.
This is impossible, because we could easily make 
$F_0$ shorter by replacing the long segment $[0,x_i]$ with $[0,v]$, for instance.

Now assume that $V_0$ ends at $0$ and, say, $x_1$ and $x_2$. Again, $v$ lies very close to $0$.
If $m=2$, then $F_0 = V_0$ and we declare ourselves happy. Otherwise, $m=3$, there is
another segment that goes from $x_3$ to either $0$ or $x_1$ or $x_2$ (the other vertex $v$
is already full), and no more, because $F_0$ is minimal. But $x_1$ (for instance) is impossible,
because we would make $F_0$ shorter by replacing $[x_3,x_1]$ with the shorter $[x_3,0]$.
So $F_0 = V_0 \cup [x_3,0]$, and we like this case too. Notice that $[0,v]$ and $[0,x_3]$
make an angle at least $2\pi/3$ at the origin, because otherwise we could pinch.

We end our discussion with the case when $F_0$ has (at least) two vertices $v_j$ 
in $B(0,R) \sm \{ 0 \}$. 

Let us count vertices and edges to reduce to one possibility.
First observe that $F_0$ contains no cycle since it is minimal. That is, $F_0$ is a tree. 
It has some vertices of valence $3$ (the $v_j$, and maybe some other), 
maybe some vertices of valence $2$, and some vertices of valence $1$ that we call extremities.
There are at most $4$ extremities, the origin and the $x_i$, because the other vertices 
have valence $3$.
It is easy to see that such a tree has $2$ extremities if it has no vertex of valence $3$,
$3$ extremities if it has one vertex of valence $3$, $4$ extremities if it has two vertices 
of valence $3$, and more otherwise (you may remove the vertices of valence $2$ to do this computation). Here we have at least two vertices $v_i$ and at most $4$ extremities, so in fact
we have exactly $2$ vertices $v_i$ and $4$ extremities, which are $0$ and three points $x_i$.
That is, $F_0$ is a simple graph with $5$ segments, and after renaming the $x_i$ and the $v_j$ 
we may assume that
$F_0 = [0,v_1] \cup [x_1,v_1] \cup [v_1,v_2] \cup [v_2,x_2] \cup [v_2,x_3]$, with segments that do not meet except at the $v_j$, and with $\2$ angles as usual.
We are a little less happy with this last case, but keep it anyway.

In all our cases, we claim that set $F_0$ gives a competitor for $Z$ in $\ol B(0,R)$.
That is, due to the simple shape of $F_0$, we can find a Lipschitz mapping $\varphi : \ol B(0,R) \to F_0$,
such that $\varphi(z) = z$ for $z\in F_0$, and in particular $\varphi(0)=0$. We extend $\varphi$ to 
$Z \sm B(0,R)$ by setting $\varphi(z) = z$ there. Notice that because near the points $x_i$, 
$Z$ is composed of a $C^1$ curve which is transverse to $\d B(0,R)$, this makes $\varphi$ Lipschitz on
$Z$. We do not care about the Lipschitz constant, and $\varphi$ is the endpoint of the family 
$\{ \varphi_t \}$, $0 \leq t \leq 1$, obtained by linear interpolation with the identity. 
Thus $\varphi(Z)$ is a sliding competitor for $Z$ in $\ol B(0,R)$ and, since $Z$ is minimal, 
$\H^1(Z \cap \ol B(0,R)) \leq \H^1(\varphi(Z \cap \ol B(0,R)) \leq \H^1(F_0)$. 
Recall that $Z \cap \ol B(0,R) = Z_R$, so \eqref{2.11} says that in fact $\H^1(Z_R) = \H^1(F_0) = \H^1(F)$,
and by its proof $Z_R$ is actually equal to $F$ (every curve in the decomposition is a line segment).
In addition, $F$ is minimal, so the discussion above, with $F_0 = F = Z_R$, gives a description of
$F_0 = Z_R = Z \cap \ol B(0,R)$.

\ms
Notice that all this happens for radii $R$ that we can take as large as we want. 
Suppose that we ever encounter the bad case when $F_0$ has five pieces. 
Then for all the radii $R'$ larger than $R$ (and for which the argument works), our description of 
$Z \cap B(0,R')$ coincides in $B(0,R)$ with the description of $Z \cap B(0,R)$, 
which means that the two vertices $v_j$ are always the same, 
and $F_0 = F_0(R')$ is just obtained from $F_0(R)$ by extending the three branches by straight lines, 
past the three $x_i$. Since we can take $R'$ as large as we want,
we see that $Z$ is the union of the two segments $[0,v_1]$ and $[v_1,v_2]$, plus three half lines,
namely the half line $L_1$ that starts from $v_1$ and goes in the direction of $x_1$,
and the two half lines $L_2$ and $L_3$ that leave from $v_2$ and go in the directions of $x_2$ 
and $x_3$ respectively. 

Denote by $e_i$ the direction of $L_i$. Since the blow-in limits of $Z$ are $Y$-sets, we see that the three
$e_i$ make $2\pi/3$ angles with each other. In particular, they lie in a same plane. Now $[v_2,v_1]$
makes $2\pi/3$ angles with $L_2$ and $L_3$ at the point $v_2$, so it lies in the same plane $P'$
(parallel to $P$) that contains $L_2 \cup L_3$. This plane contains $L_1$ too (because it contains $v_1$
and its direction contains $e_1$), and since $L_1$, $[v_1,0]$, and $[v_1,v_2]$ also make 
$2\pi/3$ angles at $v_1$, we see that $0 \in P'$ as well. It is good to know that the picture is done
in $P'$, because now $e_2$ and $e_3$ are easily seen to make angles of $\2 \pm \frac{2\pi}{6}$ with
$e_1$, a contradiction. So we may assume that our last bad case never happens for $R$ large.

Our next case is when for some $R > 0$, $F_0$ is of the form $V_0 \cup [0,x_3]$, i.e., a
truncated $Y$-set, plus a segment that goes roughly in the opposite direction. As before, for every
$R' > R$ for which we can make the description above, the set $F_0(R')$ extends $F_0(R)$. This implies
that $Z$ is a set of type $Y$, truncated at the origin, plus a half line emanating from $0$.
The blow-in limits of $Z$ are unions of three half lines leaving from $0$, and since these blow-in limits 
are minimal, the three three half lines make $2\pi/3$ angles. That is, $Z$ is a cone of type $Y$.

Now assume that this never happens, and that there is an $R$ for which $F_0 = F_0(R)$ is a truncated
$Y$-set. Then as before we can extend, and $Z$ itself is a truncated $Y$-set. Similarly, if $F_0(R)$
is composed of radii starting from the origin, and the descriptions above never occur for any $R$,
we see that $Z$ is a union of $1$, $2$, or $3$ half lines emanating from $0$ with the usual condition that
they make angles at least $2\pi/3$ at the origin.

Thus we have a description of $Z$ which fits what was announced in the statement;
Lemma~\ref{t2.2} follows.
\qed

\ms
We deduce from this a description of translation invariant sliding minimal sets of
dimension $2$.

\begin{lem} \label{t2.3}
Let $T$ be a (reduced) sliding minimal set of dimension $2$ in the whole $\R^n$, with sliding boundary $L$,
and suppose that $T$ is invariant by translations parallel to the line $L$. 
Then $T$ is either a plane, a set of type $\bV$ (two half planes bounded by $L$ and that make an angle at least $\2$ along $L$), or a set of type $\bY$, parallel to $L$ but not necessarily
containing $L$, or else a half plane bounded by $L$ (i.e., $T \in \bH(L)$) or a truncated set of type
$\bY$ (i.e., a set of the form $(Y \sm H) \cup L$, where $Y \in \bY$ has a spine parallel to
$L$, and $H \in \bH(L)$ is a half plane contained in $Y$).
\end{lem}

To prove this, write $T = L \times Z$, where $Z$ is a subset of the vector hyperplane $P$ 
perpendicular to $L$. We want to show that $Z$ is a one-dimensional minimal set in $P$,
with a sliding boundary reduced to the origin $0$, and then we'll use Lemma \ref{t2.2}. 
This is fairly a standard argument, so we  just sketch the proof and refer to Lemma 2.1 of \cite{Lu1}
for a more detailed argument.

Let $\varphi : Z \to P$ be a Lipschitz mapping such that
$\varphi(x) = x$ for $|x|$ large, and $\varphi(0) = 0$ if $0 \in Z$.
This last is enough to take care of the sliding boundary condition. That is,
in principle our competitors are of the form $\varphi_1(Z)$, where $\{ \varphi_t \}$
is a one-parameter family of continuous functions that satisfy the sliding condition
that $\varphi_t(x) \in \{ 0 \}$ when $x\in \{0 \}$ (our sliding boundary is $\{ 0 \}$). 
But we'll take $\varphi_t(x) = t \varphi(x) + (1-t) x$, and our condition that $\varphi(0)=0$
is enough for the sliding condition. 
 
Let $B$ be a ball such that $\varphi(x) = x$ for $x\in Z \sm B$
and $\varphi(Z\cap B) \subset B$. 
Suppose that, in contradiction with our claim, we can choose $\varphi$
so that
\begin{equation} \label{2.12}
\Delta : = \H^1(\varphi(Z) \cap B)-\H^1(Z \cap B) 
= \H^1(\varphi(Z \cap B))-\H^1(Z \cap B) < 0.
\end{equation}

Let $I \subset L$ denote a very long interval and let $\psi : I \to [0,1]$ be a nice cut-off function
on $I$. For the sake of definiteness, we can identify $L$ with $\R$, take $I = [-N-1, N+1]$ for
some large $N$ and choose $\psi(y) = \max(0,\min(1, N+1-|y|))$ for $y\in \R$.
Denote by $(x,y)$ the generic point of $\R^n$, with $x\in P \simeq \R^{n-1}$
and $y\in L \simeq \R$. 
A good competitor for $T$ is $f(T)$, where $f : T \to \R^n$ is defined by
$f(x,y) = (\psi(y) \varphi(x) + (1-\psi(y)) x, y)$. 
It is easy to see that $f(T)$ is a sliding competitor
for $T$ in the rectangular shaped set $R = B \times I$, in particular because $f(x,y) = (x,y)$ 
when $x=0$ and because it is easy to interpolate between the identity and $f$.

The minimality of $T$ says that $\H^d(T \cap R) \leq \H^d(f(T) \cap R)$.
Set $R' = B \times [-N,N]$, and observe that 
\begin{equation} \label{2.13}
\H^2(T\cap R') = \H^2((Z \cap B) \times [-N,N]) = 2N \H^1(Z \cap B)
\end{equation}
not completely trivially, but because $Z$ is rectifiable. 
See for instance the computations of pages 530-531 in \cite{MSBook}, although in a slightly different context.
The rectifiability of $Z$ itself comes from the rectifiability of $T = Z \times L$; we leave the details. Similarly, the 
$2$-rectifiability of $f(T \cap R') = (\varphi(Z\cap B)) \times [-N,N]$ (recall that $\psi(y) = 1$
on $[-N,N]$) yields
\begin{equation} \label{2.14}
\H^2(f(T \cap R')) = \H^2((\varphi(Z \cap B)) \times [-N,N])
= 2N \H^1(\varphi(Z \cap B)),
\end{equation}
so we win $2N \Delta$ from the contribution of $R'$. 
We still need to estimate the contribution of $R \sm R'$.
Since $\varphi$ is Lipschitz, $\H^2(f(T\cap (R\sm R'))\leq C$, where $C$ depends on the Lipschitz
constant for $\varphi$. It is possibly huge, but it does not depend on $N$. We take $N$ large,
add this to the estimates from inside of $R'$, and get the desired contradiction with the minimality 
of $T$. 

So $Z$ is a sliding minimal set, Lemma \ref{t2.2} gives a good description of $Z$,
and the description of $T = Z \times L$ needed for Lemma \ref{t2.3} follows.
\qed

\ms 
When we restrict to cones, Lemma \ref{t2.3} yields that 
(with the notation of Subsection \ref{S1.2}, and still assuming that $L$ is a line through the origin)
\begin{equation} \label{2.15}
\begin{aligned}
&\text{if $T$ is a sliding minimal cone with sliding boundary $L$, and $T$ is invariant }
\\ &\hskip1cm
\text{by translations parallel to $L$, then } T \in \bH(L) \cup \bV(L) \cup \bY(L).
\end{aligned}
\end{equation}

\ms
We will often use Lemma \ref{t2.3} and \eqref{2.15} to control limits of minimal cones,
and then obtain information in the direction of Proposition \ref{t2.1}. The standard notation for this
is the following. We have a sequence $\{ X_k \}$ of sliding minimal cones associated to the
boundary $L$ (a line through the origin). We select points $a_k \in K_k = X_k \cap \S$ and radii
$r_k > 0$, with $\lim_{k \to +\infty} r_k = 0$, and consider
\begin{equation}\label{2.16a}
Y_k = r_k^{-1}(X_k-a_k).
\end{equation}
Notice that $0 \in Y_k$; this allows us to take a subsequence, which we shall still denote the same way,
so that $\{ Y_k \}$ converges to a closed set $Y$, and $\{ a_k \}$ converges to a limit $a \in \S$.
We will need to know that $Y$ is invariant by translations in the direction of $a$, i.e., that
\begin{equation} \label{2.17}
\xi+ta \in Y \ \text{ for $\xi \in Y$ and $t \in \R$.} 
\end{equation}
Indeed, we can find $\xi_k \in Y_k$, so that $\xi_k$ tends to $\xi$. Set
$\zeta_k = a_k + r_k \xi_k$; then $\zeta_k \in X_k$, and since $X_k$ is a cone,
$s \zeta_k \in X_k$ for $s > 0$. Then $r_k^{-1}(s \zeta_k - a_k) \in Y_k$ for $s > 0$.
But $r_k^{-1}(s \zeta_k - a_k) = r_k^{-1} (s a_k + s r_k \xi_k - a_k)
=  s \xi_k + (s-1)r_k^{-1} a_k$. We apply this with $s= 1 + r_k t$, get that
$(1 + r_k t) \xi_k + t a_k \in Y_k$ for $k$ large, take a limit, and get \eqref{2.17}.

Let $z_k \in L$ minimize the distance to $a_k$, and notice that $Y_k$ is a sliding minimal
set, with respect to the boundary $L_k = L - r_k^{-1} a_k = L + r_k^{-1}(z_k - a_k)$.

There will be two main cases. The first one is when $\lim_{k \to +\infty} r_k^{-1} \dist(a_k,L) = +\infty$,
or equivalently $\lim_{k \to +\infty}  \dist(0,L_k) = +\infty$.
In this case, since $Y_k$ is a plain minimal set in $B(0,\dist(0,L_k))$, 
then by Theorem 4.1 (and Definition 2.4) in \cite{limits}, 
$Y$ is a minimal set in $\R^n$, with no sliding boundary condition. Since by \eqref{2.17} it is also
invariant by translations the direction of $a$, the simpler variant of Lemma \ref{t2.3} where there is 
no boundary constraints implies that 
\begin{equation} \label{2.18a}
\text{$Y$ is a plane or a set of type $\bY$ (possibly not centered at $0$).} 
\end{equation}
The other possibility is that $ \dist(0,L_k) = r_k^{-1} \dist(a_k,L)$ stays bounded; 
then, modulo a new sequence extraction, we may assume that $\{ L_k \}$ converges to 
a line $L_\infty$, which is parallel to $L$ (and the $L_k$). 
Theorem 10.8 or 21.3 in \cite{Sliding} says that $Y$ is a sliding minimal set, 
with boundary $L_\infty$, and since $Y$ is still invariant by translations in the direction of $a$
(which happens to be the direction of $L$, since $a = \lim_{k \to +\infty} a_k$ and $\dist(a_k,L)$
tends to $0$), Lemma~\ref{t2.3} says that 
\begin{equation} \label{2.19a}
\begin{aligned}
&\text{$Y$ is a plane, a set of type $\bV$ (bounded by $L$), a set of type $\bY$
 (with a spine }
\cr&\hskip1.1cm
\text {parallel to $L$), or a half plane or a truncated set of type $\bY$.} 
\end{aligned}
\end{equation}

We return to the proof of Proposition \ref{t2.1}.
The following lemmas will help with the relative position and length of the arcs $\cC_j$
that compose $K = X \cap \S$. We start with a description of $K$ far from $L$, which is
more precise than what we did near \eqref{2.8} because we give a lower bound
for the radius of the good balls.

\begin{lem} \label{t2.4}
If $\eta_1$ is small enough, depending only on $n$, not on $X$,
then if $a$ is a vertex of $K$ in the description near 
\eqref{2.8}, $K \cap B(a, \eta_1 \dist(a,L))$ is the union of three geodesics that
leave from $a$ with equal angles of $\frac{2 \pi}{3}$.
\end{lem}

We shall prove this with a contradiction and compactness argument. 
Suppose that the lemma fails, and let $X_k$, $L_k$, $K_k = X_k \cap \S$, 
and $a_k \in K_k \sm L_k$ provide a counterexample, with $\eta_1(k) = 2^{-k}$. 
By rotation invariance, we may assume that $L_k = L$ stays the same.
By \eqref{2.8}, there is a neighborhood of $a_k$ where 
$K_k$ is composed of three arcs of geodesic. That is, for each $k$ we can find $r > 0$ such that 
\begin{equation} \label{2.16}
K_k \cap B(a_k,r) = (\gamma_1 \cup \gamma_2 \cup \gamma_3) \cap B(a_k,r)
\end{equation}
for some choice of three geodesics $\gamma_j$, $1 \leq j \leq 3$, that leave from
$a_k$, make $\2$ angles at $a_k$, and go at least to $\d B(a_k,r)$. Let $r_k$ denotes the
largest $r > 0$ such that the representation \eqref{2.16} holds.
Since the description of the lemma fails for $r = 2^{-k} \dist(a_k,L)$, we see
that $r_k \leq 2^{-k} \dist(a_k,L) \leq 2^{-k}$.

Consider $Y_k = r_k^{-1}(X_k-a_k)$ as above, and replace our sequence with a subsequence for 
which $Y_k$ tends to a limit $Y$. Since $r_k^{-1} \dist(a_k,L) \geq 2^k$ tends to $+\infty$,
we see that $Y$ is a plane or a set of type $\bY$, as in \eqref{2.18a}.

Since by definition of $r_k$ \eqref{2.16} holds for $r = r_k/2$, we see that
$X_k$ has a beautiful description as a set of type $\bY$ in $B(a_k, r_k/3)$,
$Y_k$ has a similar description in $B(0,1/3)$, and when we take a limit $Y$
we get a cone of type $\bY$ (this time centered at $0$).

Return to $Y_k$ and $X_k$. Since $Y_k$ tends to $Y$, we get that
$d_{0,10}(Y_k,Y)$ tends to $0$, or equivalently
$d_{0,10r_k}(X_k,a_k + Y) = d_{0,10}(Y_k,Y)$ tends to $0$
(see the definition \eqref{1.13}). Let us again be  slightly brutal and apply the regularity
theorem from \cite{C1}; for $k$ large enough, we get that in $B(a_k, 3r_k)$, 
$X_k$ is a smooth version of $a_k + Y$, with small $C^1$ constants, to the point that 
$K \cap B(a_k, 2r_k)$ is, by the implicit function theorem, composed of exactly $3$
smooth curves that meet at $a_k$. These smooth curves are arcs of geodesics
(by \eqref{2.8}) and this contradicts the definition of $r_k$. Lemma \ref{t2.4} follows. 
\qed

\ms
With almost the same proof, we can also get a uniform control of $X$ near 
$\ell \in L \cap \S$, provided that $\ell \notin K$.

\begin{lem} \label{t2.5}
If $\eta_1$ is small enough, depending only on $n$ but not on $X$,
and if $\ell \in L \cap \S \sm K$, then for each $a \in K \cap B(\ell, 10^{-1})$
which is a vertex of $K$, $K \cap B(a, \eta_1)$ is the union of three geodesics that
leave from $a$ with equal $\frac{2 \pi}{3}$ angles.
\end{lem}

\ms
The difference is that the size of the ball no longer depends on $\dist(a,L)$.
We start the proof the same way. 
By rotation invariance, it is enough to prove this for a fixed $L$ and $\ell$.
Then we proceed by contradiction and suppose that for $k \geq 1$,
$X_k$ and $a_k$ define a counterexample with $\eta_1 = 2^{-k}$. We define
$Y_k$ as before, i.e., let $r_k$ be the smallest radius $r$ such that \eqref{2.16} fails,
and set $Y_k = r_k^{-1} (X_k - a_k)$. Notice that $r_k \leq 2^{-k}$ because $X_k$
is a counterexample.

Switching to a subsequence if needed, we can assume that $Y_k$ 
converges to a limit $Y$. Now we claim that $Y$ is a minimal set in $\R^n$, 
with no sliding boundary condition, but for a different reason as before.

For our proof of \eqref{2.18a}, we used the fact that $L$ was too far. 
Here $X_k$ is sliding minimal in $B(a_k, 1/2)$,  with a sliding boundary $L$ that could be very
close to $a_k$. But $X_k$ does not meet $L$ in that ball (because we assumed that $K_k$ does not 
contain $\ell$), so we easily deduce from the definitions that $X_k$ is a (plain) minimal set in $B(a_k, 1/2)$. 
This is because \eqref{1.7} is void here. Then $Y_k$ is a plain minimal set in $B(0,r_k^{-1}/2)$,
and by Theorem 4.1 in \cite{limits}, $Y$ is a minimal set in $\R^n$. Again we have \eqref{2.18a},
i.e., $Y$ is a plane or a set of type $\bY$.

The rest of the proof is as above: $Y$ is actually a set of type $\bY$ with a spine that contains $a$,
then for $k$ large $X_k$ is so close to $r_k Y + a_k$ in $B(a_k,10r_k)$ that it coincides
with a $C^1$ version of that set in $B(a_k,2r_k)$. By \eqref{2.8}, $K$ is composed of three
geodesics inside $B(a_k,2r_k)$, this contradicts the definition of $r_k$, and Lemma \ref{t2.5} follows.
\qed

\ms
\begin{lem} \label{t2.6}
There is a small $\eta_2 >0$, depending on $n$ but not on $X$, such that
if $K$ contains the point $\ell \in L \cap \S$, $a\in K \sm L$
is one of the vertices of $K \sm L$, and $|a-\ell| \leq \eta_2$,
then there is no other vertex of $K$ in $B(\ell,10|a-\ell|) \sm B(\ell,10^{-1}|a-\ell|)$.
\end{lem}

\ms
Once more we prove this by contradiction and compactness. 
Suppose the lemma fails, and let $X_k$, $K_k$, $\ell_k \in L_k\cap K_k$, and
$a_k \in K_k$ provide a counterexample with $\eta_2 = 2^{-k}$. By rotation
invariance we may assume that $L_k = L$ and $\ell_k =\ell$ are always the same.
Set $r_k = |a_k - \ell| \leq 2^{-k}$; thus $r_k$ tends to $0$
and $r_k^{-1}\dist(a_k,L)$ tends to $1$ (because $a_k$ tends to $\ell$).
We are in the situation of \eqref{2.19a}, where modulo a sequence extraction
$Y_k = r_k^{-1} (X_k - a_k)$ tends to a sliding minimal set $Y$, which is a plane, 
a set of type $\bY$, a cone of type $\bV$, a half plane, or a truncated $\bY$-set,
each time bounded by a half line $L_\infty$ parallel to $L$.

Recall that $a_k$ is a vertex of $K_k \sm L$; Lemma \ref{t2.4} says that in 
$B(a_k, \eta_1 r_k)$, $K_k$ is composed of three geodesics $g_1$, $g_2$, $g_3$ that 
meet at $a_k$ with $120^\circ$ angles. In the same ball, $X_k$ coincides with the cone 
$H_k$ over  $g_1 \cup g_2 \cup g_3$. Or equivalently,
$Y_k$ coincides with $r_k^{-1}(H_k - a_k)$ in $B(0, \eta_1)$.
Thus $Y$ has a singularity of type $\bY$ at the origin, and is a $\bY$-set, possibly truncated,
with a spine parallel to $L$ (because it is invariant by translations in the direction of 
$\ell = \lim_{k \to +\infty} a_k$).

But the contradiction assumption says that $K_k$ has another vertex 
$b_k \in B(\ell,10 r_k) \sm B(\ell,10^{-1} r_k)$, and Lemma \ref{t2.4} says that
in $B(b_k, \eta_1 r_k/10)$, $K_k$ is composed of three geodesics that meet at $b_k$.
In particular, $a_k$ lies outside of this ball, hence $|b_k-a_k| \geq \eta_1 r_k/10$.
We may extract a new subsequence so that $\wt b_k =  r_k^{-1}(b_k - a_k)$
converges to a limit $\wt b$, and the same argument as above says that $\wt b$ also lies on the
spine (the singular set) of $Y$, just like $0$. But $\wt b$ lies on $L^\perp$
(the hyperplane orthogonal to $L$), because both $a_k$ and $b_k$ lie in $\S$ and tend to $\ell$,
and in addition $|\wt b - \wt a| \geq \eta_1/10$ (because $|b_k-a_k| \geq \eta_1 r_k/10$).
This is impossible; Lemma~\ref{t2.6} follows.
\qed

\ms
The same argument says a little more. 
Let $X$, $K$, and $a$ be as in Lemma~\ref{t2.6}.
We claim that not only $K \cap B(a,\frac{9 |a-\ell|}{10})$ is composed of three
geodesics (with no other vertex of $K$), but also that (again if $\eta_2$ is small enough), 
one of these geodesics makes an angle less than $\pi/100$ with the geodesic $\rho(a,\ell)$ 
from $a$ to $\ell$.

Indeed, otherwise we proceed as in the proof of Lemma~\ref{t2.6}, with a sequence
$\{ X_k \}$ for which the three geodesics that compose $K_k \cap B(a_k, \eta_1 r_k)$
make angles at least $\pi/100$ with $\rho(a_k,\ell)$. 
As before, we can extract a subsequence for which $Y_k = r_k^{-1} (X_k - a_k)$ tends to a sliding minimal 
set $Y$, which is either a $\bY$-cone with a spine parallel to $L$, or such a $\bY$-cone, 
truncated by a line $L_\infty$ parallel to $L$. Either way, $\ell \in X_k$ by assumption, 
so $z_k = r_k^{-1} (\ell - a_k) \in Y_k$, and we can extract a subsequence so that 
$z = \lim_{k \to +\infty} z_k \in Y$. Notice that $|z|= 1$ because
$|z_k| = 1$ since $r_k = |a_k - \ell|$. 

It is easy to see that the direction of $\rho(a_k,\ell)$ at $a_k$ tends
to $z$. But on the other hand the directions at $a_k$ of the three geodesics that compose 
$K_k \cap B(a,\frac{9 |a_k-\ell|}{10})$, or equivalently the directions at $0$ of the three geodesics that
compose $\wt K_k \cap B(0,\frac{9}{10})$, with $\wt K_k = r_k^{-1} (K_k - a_k)$,
tend to the unit directions of the faces of $Y$ (intersected by the orthogonal of $L$).
Thus one of these directions tends to $z$ (because $z\in Y$), a contradiction with our assumption that 
they all make large angles with the direction of $\rho(a_k,\ell)$ (that also tends to $z$).
This proves our claim.

\ms
Let us continue with our assumption that $K \cap L$ contains a point $\ell$. 
We now claim that 
\begin{equation} \label{2.18}
K \cap B(\ell,\eta_2/10)\sm L \ \text{ contains at most one vertex.} 
\end{equation}
Indeed, suppose that $K \cap B(\ell,\eta_2/10) \sm L$ has two vertices $a$ and $b$.
We may assume that $|b-\ell| \leq |a-\ell|$. Notice that Lemma \ref{t2.6} says that the
values of $|x-\ell|$, where $x \in K \cap B(\ell,\eta_2) \sm L$ is a vertex of $K$, are
lacunary, so we may assume that $b$ was chosen so that $|b-\ell|$ is maximal once $a$
is chosen. Also, $|b-\ell| \leq 10^{-1}|a-\ell|$ by lemma \ref{t2.6}.

\begin{figure}[!h]  
\centering
\includegraphics[width=10cm]{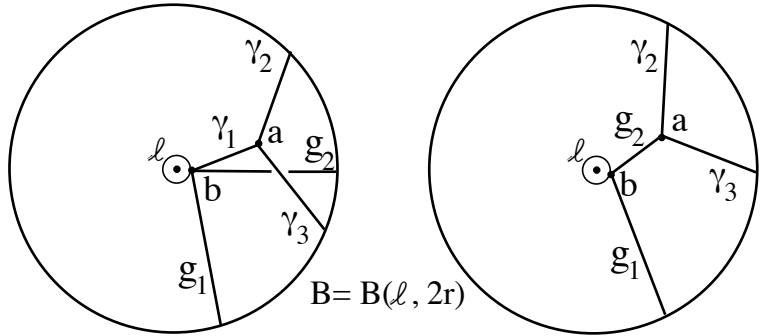}
\caption{These arcs of geodesics are contained in $K$ (two cases). On the left, $\gamma_3$
does not need to cross $g_2$ because $n \geq 4$ is allowed
\label{f2.0}}
\end{figure}

Set $r = |a-\ell|$, $B = B(\ell, 2r)$, and $A = B \sm B(\ell, |b-\ell|)$; 
we now give a description of $K \cap A$. 
First we have two geodesics $g_1$ and $g_2$ that leave from $b$,
making $120^\circ$ angles with each other and also roughly with $\rho(b,\ell)$.
Because of this, they go away from $B(\ell, |b-\ell|)$, i.e., they start in $A$.
They stay in $K$ as long as they stay in $B$ and they don't meet a vertex of $K$;
since there is only one vertex in $A$ (namely, $a$), we only have two options
(see Figure~\ref{f2.0}).
Either $g_1$ and $g_2$ both miss $a$, and then $K$ contains $A \cap (g_1 \cup g_2)$.
Or else one of them, say $g_2$, contains $a$, and then we only know that
$K$ contains $A \cap g_1$ and $\rho(b,a)$.

We also know that $B \sm \frac{1}{5} B$ contains three arcs
of geodesics $\gamma_i$, that leave from $a$ with $120^\circ$ angles and
go all the way to the boundary of $B \sm \frac{1}{5} B$ (because they don't meet
a vertex).

Altogether, we found a collection of reasonably long geodesics that are contained
in $K \cap B$, either $4$ of them (three that make a $Y$ centered at $a$, 
plus one leg that leaves from $b$ with a $120^\circ$ angle), 
or $5$ (three that make a $Y$ centered at $a$, and two other ones that make 
a disjoint $V$). It is important for the present argument that long means, of diameter
at least $|\ell - a|/10$, say.
We claim that this is impossible. We proceed by contradiction
and compactness as in the previous lemma, and get a description of the limit $Y$
of a convergent subsequence of normalizations $Y_k$ of counterexamples $X_k$.
As above, $Y$ is a possibly truncated cone of type $\bY$ that is centered at $0$
(a normalized limit of the $a_k$) and contains a line $L_\infty$ parallel to $L$. 
But then $X_k$ looks a lot in $B$ like the image of $Y$ by a translation and a dilation, 
and this does not fit the fact that $K_k = X_k \cap \S$ contains the four or five long 
geodesics above. This proves \eqref{2.18}.

\ms
Let us now check that 
\begin{equation} \label{2.19}
|a-b| \geq \frac{\eta_1\eta_2}{20} \ \text{ when $a$, $b$ are different vertices of } K \sm L.
\end{equation}
Indeed if $\dist(a,L) \geq \eta_2/20$, Lemma \ref{t2.4} says that 
$|b-a| \geq \eta_1\dist(a,L) \geq \eta_1 \eta_2/20$. If $|a-\ell| \leq \eta_2/15$
for some $\ell \in K \cap L$, \eqref{2.18} says that $|b-\ell| \geq \eta_2/10$ and
hence $|a-b| \geq \eta_2/60$. Finally if $|a-\ell| \leq \eta_2/15$ for some 
$\ell \in L \cap \S \sm K$, Lemma \ref{t2.5} says that $|b-a| \geq \eta_1$.

\ms
Because of this, $K$ has a finite number of vertices, hence it is composed of a finite 
number of geodesic arcs, plus some full great circles. Recall that when an arc $\cC$ meets $L$,
we consider the points of $K \cap \cC$ as vertices, i.e., we cut $\cC$ at these points.
This gives our decomposition of $K$ into the $\cC_j$, $j\in J$. 

Before we start the verification of the various properties stated in Proposition \ref{t2.1}, let us say 
two last words about the minimizing properties of $K$ itself. It will be good to know that
\begin{equation} \label{2.20}
K \ \text{ is a weak almost minimal set in $\S$, with sliding boundary } K \cap L,
\end{equation}
even though we shall also try to provide proofs that do not use this fact.
Let us first say what \eqref{2.20} means. The vocabulary comes from 
Definition 9.1 of \cite{Holder}, where a similar notion (without sliding boundary)
was used to record some easy properties of $K = X \cap \S$ when $X$ is a minimal
cone, in order to get the description in terms of geodesics that we used for \eqref{2.8}.
By \eqref{2.20} we mean that if $f : \S \to \S$ is an $M$-Lipschitz mapping 
and $B(x,r)$ is a ball centered on $\S$ such that
\begin{equation} \label{2.21}
f(y) = y \ \text{ for } y\in \S \sm B(x,r) \ \text{ and } \ 
f(\S \cap B(x,r)) \subset \S \cap B(x,r), 
\end{equation}
and also $f(\ell) = \ell$ for every $\ell \in K \cap L$,
then 
\begin{equation} \label{2.22}
\H^1(K \cap B(x,r)) \leq \H^1(f(K) \cap B(x,r)) + C (1+M) r^2.
\end{equation}
The present definition is a little less demanding than Definition 9.1 of \cite{Holder},
where we also required \eqref{2.22} when $f$ is piecewise $M$-Lipschitz, but this will be
enough for our purposes. We also use a specific gauge function (namely, $Cr$, with a $C$
that depends only on the dimension $n$) in \eqref{2.22}, again because this is what we get 
from the proof. On the other hand, we added the requirement that $f(\ell) = \ell$ for 
$\ell \in K \cap L$, to account for the sliding condition for $X$.

Now the proof of Proposition 9.4 in \cite{Holder} applies to the present
situation and shows that \eqref{2.20} holds for every sliding minimal cone $X$
of dimension $2$, with boundary $L$. Our extra condition is of course used to ensure that the
competitors build in \cite{Holder} come from one-parameter families $\{ \varphi_t \}$
such that $\varphi_t(y) \in L$ when $y \in X \cap L$.
The details of Proposition~9.4 in \cite{Holder} are easy; we just use any map $f$ as above
to construct a competitor for $X$; however we find it easier to refer to \cite{Holder} rather
than doing the verification here. Of course we'll use \eqref{2.20} a few times to derive a 
contradiction when needed.

\ms
Return to the properties of Proposition \ref{t2.1}, and let us first check what happens at the vertices. 
The fact \eqref{2.6} that near each vertex $a\in K \sm L$,
$K$ consists of three geodesics that leave from $a$ with $120^\circ$ angles comes from
\eqref{2.8}, and was already used many times. Similarly, let us check \eqref{2.7}, 
which says that when $\ell \in K \cap L$, there is a small neighborhood of $\ell$ where 
$K$ is composed of one, two, or three geodesic arcs that leave from $\ell$ and
make angles of at least $120^\circ$.

First, we may restrict to the geodesics that contain $\ell$, because there is a finite number
of geodesics, and the ones that don't contain $\ell$ don't meet some small ball centered at $\ell$.
The remaining geodesics all start from $\ell$ (because for the other ones we added $\ell$ as a vertex),
and is easy to check that they make angles of at least $120^\circ$ at $\ell$, because otherwise 
we may pinch two of them in a small ball $B(\ell,r)$, make $K$ sorter by at least $r/C$, and contradict 
\eqref{2.20}.
Alternatively (if you don't like weak almost minimality), we could say that a blow-up limit of $X$ 
at $\ell$ is composed of as many half planes bounded by $L$ as there are geodesics in $K$ near $\ell$, 
and that make the same angles along $L$ as the geodesics at $\ell$; then Lemma \ref{t2.3}
gives the desired result.

\ms
Next we show that each full great circle in the list of $\cC_j$ is far from the rest of $K$,
i.e., that there is $\eta_3 > 0$, that depends only on $n$, such that
\begin{equation} \label{2.23}
\dist(\cC_j, K\sm \cC_j) \geq \eta_3 \ \text{ when $\cC_j$ is a full great circle of } K. 
\end{equation}

Of course \eqref{2.3} will follow from this (we'll take $\eta_0$ very small at the end). 
As usual, we proceed by contradiction and compactness. 
So suppose that for each $k \geq 0$, \eqref{2.23} fails for $\eta_3 = 2^{-k}$, 
and let $K_k$ provide a counterexample. 
By rotation invariance, we may assume that $L$ is the same for each $K$,
and our assumption is that we can find a great circle $\cC_k$ and a point $a_k \in K_k \sm \cC_k$ 
such that $\dist(a_k, \cC_k) \leq 2^{-k}$. 

First observe that $K_k \sm \cC_k$ is closed. Indeed, otherwise we can find some $\xi \in \cC_k$
which is the limit of a sequence in $K_k \sm \cC_k$. Then $\xi \in L$ (because otherwise this 
contradicts \eqref{2.6}); even that way the two arcs of $\cC_k$ near $\xi$ make a $180^\circ$ 
angle at $\xi$, which by \eqref{2.7} excludes the possibility that other geodesics of $K$ end at $\xi$.
Recall also that $K$ has a finite number of vertices, a contradiction that shows that $K_k \sm \cC_k$ is closed.

By compactness of $K_k \sm \cC_k$, we may assume that $a_k$ minimizes the
distance to $\cC_k$ (in $K \sm \cC_k$). Set $r_k = \dist(a_k, \cC_k)$; thus
$0 < r_k \leq 2^{-k}$.

We start with the more interesting case when $\{ r_k^{-1} \dist(a_k,L) \}$ is a bounded sequence.
We may assume that (out of the two possibilities) there is a fixed $\ell \in L$ such that 
$r_k^{-1} \dist(a_k, \ell) \leq C$. Set $Y_k = r_k^{-1}(X_k-a_k)$ as usual, take
a converging subsequence, denote by $Y$ the limit, and notice that $Y$ satisfies \eqref{2.19a}.
But $X_k$ contains the plane $P_k$ that contains $\cC_k$, and which lies at distance a little smaller
than $r_k$ from $a_k$ (a little smaller because the closest point of $P_k$ lies a little inside of $\S$);
This means that $Y_k$ contains $\wt P_k = r_k^{-1}(P_k-a_k)$, which almost lies at distance $1$
from $0$; at the limit, $Y$ contains a plane $P$ at distance $1$ from the origin. By \eqref{2.19a},
$Y = P$; this is impossible because $a_k \in K_k$ and hence $0 \in Y$.

We are left with the case when, modulo a sequence extraction, $r_k^{-1}\dist(a_k,L)$ 
tends to $+\infty$. This time, modulo extraction, $Y_k = r_k^{-1}(X_k-a_k)$, and
tends to a set $Y$ which is a plane or set of type $\bY$, as in \eqref{2.18a}.
As before, $\dist(P,0) = 1$, which contradicts the fact that $0 \in Y$ because $a_k \in Y_k$.
This last contradiction completes our proof of \eqref{2.23}, and again, \eqref{2.3} follows.

\ms 
Next we want to check \eqref{2.4}. The following will be useful.

\begin{lem} \label{t2.7}
If $\eta_4$ is chosen small enough, the following happens. 
Suppose that $\ell \in K\cap L$ and $a\in K \sm L$ are such that
$|a-\ell| \leq \eta_4$. Then the geodesic $\rho(\ell,a)$ is contained in $K$,
there is at most one other arc $\gamma$ of $K$ that leaves from $\ell$, 
and (if $\gamma$ exists) $\H^1(\gamma) \geq \eta_2/10$ and $\gamma$ makes an 
angle at least $\frac{9 \pi}{10}$ with $\rho(\ell,a)$ at $\ell$.
\end{lem}

Let $\ell \in K \cap L$ and $a\in K \sm L$ be as in the statement.
We know from \eqref{2.18} that there is no other vertex of $K$ in
$B(\ell,\eta_2/10)$, so, except perhaps for the geodesic $\rho(\ell,a)$
if it lies in $K$, all the arcs of $K$ that leave from $a$ or $\ell$ are at least
$\eta_2/20$ long.

Let us first assume that $K$ does not contain $\rho(\ell,a)$. 
Then $K$ contains three arcs $\gamma_i = \rho(a,b_i)$, $1 \leq i \leq 3$, of length 
at least $\eta_2/20$, that leave from $a$ with $120^\circ$ angles, and also an arc 
$\gamma = \rho(\ell, b)$, of length least $\eta_2/20$ too, and that leaves from $\ell$.
All these arcs are disjoint, except perhaps for their endpoints, by \eqref{2.8} and the definition 
of our decomposition of $K$.

\begin{figure}[!h]  
\centering
\includegraphics[width=7cm]{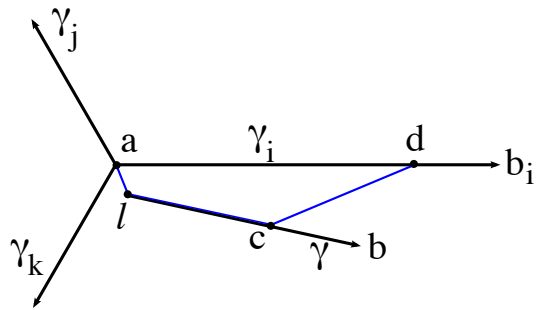}
\caption{A replacement fro $\gamma_i$. Scales are not respected
\label{f2.1}}
\end{figure}

Notice that $|a-\ell| \leq \eta_4$, so for $\eta_4$ small enough, it is not hard to imagine that
we could construct a competitor for $K$ that contradicts the weak almost minimality 
property \eqref{2.20}. See Figure \ref{f2.1} for a hint, 
but don't forget that even though $\gamma$ does not meet the $\gamma_j$, the point
$\ell$ could be more or less anywhere on the sphere $\S(a,|\ell-a|)$ (and not just in a triangular sector 
as the picture suggests) because we allow subsets of $\R^n$, $n > 3$. 
Also, we do not exclude the case when other pieces of $K$ pass by, but this is not a 
real problem because they stay at positive distance from the rest of the picture, 
so our Lipschitz deformation $f$ can be chosen so that $f(x)=x$ on them, and
as alluded to above, \eqref{2.20} is also valid for piecewise $M$-Lipschitz functions $f$.

But we shall avoid using \eqref{2.20}, and instead we will use compactness again. 
So suppose that for $k \geq 0$, we can find an example $X_k$, with $\ell \in K_k \cap L
= X_k \cap \S \cap L$, a vertex $a_k$ of $K_k$ such that $d_k = \dist(\ell,a_k) \leq 2^{-k} \eta_2$, 
three geodesics arcs $\gamma_{i,k} \subset K_k$ of length $\eta_2/20$ that leave from $a_k$,
and a fourth one, $\gamma_k \subset K_k$, of length $\eta_2/20$ too, that leaves from $\ell$ 
and is disjoint from the others.

Consider $Y_k = d_k^{-1}(X_k - \ell)$, which contains the origin, and as usual take a subsequence
for which $Y_k$ converges to a limit $Y$. Then, by the proof of \eqref{2.19a}, 
$Y$ is a sliding minimal cone with boundary $L$, which is invariant by translations in the direction 
of $L$, so it is a plane or a $\bY$ set (through the origin, since $\ell \in K_k$), 
or else a half plane, a $\bV$ set, or a $\bY$-set truncated by $L$.

The geodesics $\wt \gamma_{i,k} = d_k^{-1}(\gamma_{i,k} - \ell)$ and 
$\wt \gamma_{k} = d_k^{-1}(\gamma_{k} - \ell)$ (in the spheres $d_k^{-1}(\S - \ell)$)
tend to four half lines (maybe at the price of extracting a new subsequence, if you prefer),
and these half lines are contained in $Y \cap L^\perp$. Because of this, there is only one possibility:
$Y$ is a $\bY$-set, that contains $L$ but is centered at the limit $\wt a$ of the 
$\wt a_k = d_k^{-1}(a_{k} - \ell)$; $\wt a$ lies at distance $1$ from $L$, and the limit
$\wt \gamma$ of the $\wt \gamma_{k} = d_k^{-1}(\gamma_{k} - \ell)$ is contained in the branch
of $Y \cap L^\perp$ that contains $0$. 

Select $z\in \wt \gamma$, at distance $1$ from $0$ (and hence $z=-\wt a$). 
Then set $B = B(z, 10^{-1})$, and notice that $Y_k$ converges to $Y$ (or equivalently the plane 
$P$ that contains the face of $Y$ that contains $z$) in $9B$. Then for $k$ large enough, we may apply the 
standard regularity theorem to $X_k$ in $8B$, and find that in $2B$, $X_k$ is a $C^1$ surface, 
and at the same time a small Lipschitz graph over $P$. Set $S(k) = d_k^{-1}(\S - \ell)$; this is a very
large sphere, with tangent directions that tend to $L^\perp$; by the implicit function theorem,
$X_k \cap S(k) = d_k^{-1}(K_k - \ell)$ is a $C^1$ curve in $2B$, and also a small 
Lipschitz graph over the line that contains $\wt \gamma$.

But at the same time $d_k^{-1}(K_k - \ell)$ contains the two geodesics 
$\wt \gamma_k$ and $\wt \gamma_{i,k}$ (for some $i$), which are disjoint and both converge
to $\wt \gamma$ in $2B$. This is impossible; we are left with the other case when $K$ contains 
$\rho(\ell,a)$.

In this case, in addition to $\rho(\ell,a)$, $K$ contains two geodesics $\rho_1$ and $\rho_2$
that leave from $a$ with $\2$ angles, and maybe one or two geodesics $\rho'_j$ that leave from
$0$ (again with angles at least $\2$). All these geodesics have lengths at least $\eta_2/20$
(because there is no vertex of $K$ nearby where they could stop). A priori there may also be other 
pieces of $K$ that pass near $\ell$, but all we have to do now is prove that there is at most one $\rho'_j$,
and that it goes in a direction almost opposite to $\rho(\ell,a)$ at $\ell$.

Again it is simpler to prove this by compactness. Suppose not, let $X_k$ provide an example, with
$d_k = \dist(\ell,a_k) \leq 2^{-k} \eta_2$. This means that in addition to $\rho(\ell,a_k)$
and the two geodesics $\gamma_{i,k}$ that leave from $a_k$ with $\2$ angles and are at least
$\eta_2/20$-long, we have at least one more geodesic $\rho_k \subset K$ that leaves from $0$, 
is at least $\eta_2/20$-long, and makes an angle smaller than $\frac{9\pi}{10}$ with the direction of 
$\rho(\ell,a_k)$ at $\ell$. Indeed, if we have two, they make $\2$ angles with the direction of 
$\rho(\ell,a_k)$, which is even worse. 

As usual, set $Y_k = d_k^{-1}(X_k - \ell)$, and extract a subsequence
for which $Y_k$ converges to a limit $Y$. The same argument as above shows that $Y$ is a sliding minimal
set which is invariant by translations in the direction of $L$, then is one of the examples allowed by 
\eqref{2.19a}, and because of $\rho(\ell,a_k)$, the $\gamma_{i,k}$, and $\gamma_k$, is a set of
type $\bY$ that contains $L$ and is centered at $\wt a$, which lies at distance $1$ from $L$. 
But the geodesics $\wt\gamma_k = d_k^{-1}(\gamma_{k} - \ell)$ converge (modulo extraction 
if the reader wishes) to a half line that makes an angle at most $\frac{9\pi}{10}$ with the direction of 
the half line $[0,\wt a)$, and is contained in $Y$; 
this contradiction completes the proof of Lemma \ref{t2.7}.
\qed

\ms
It looks like we forgot some additional information that we could prove, the fact that in the situation of 
Lemma \ref{t2.7}, there in no other piece of $K$ in $B(\ell, \eta_2/40)$, but we shall return to this soon. 

\ms
We easily deduce \eqref{2.4} from the lemma. Let $\cC_j$ be one of the arcs that
compose $K$, and suppose that $\H^1(\cC_j) \leq \frac{\eta_1 \eta_2}{20}$. 
By \eqref{2.19}, at least one of its endpoints lies in $K \cap L$ (call it $\ell$), 
and by \eqref{2.18} the other one (call it $a$) is the the only point of
$K \cap B(\ell, \eta_2/10) \sm L$. In addition, if $\H^1(\cC_j) \leq \eta_4$,
Lemma \ref{t2.7} says that there is at most one other arc that leaves from $\ell$, 
and that it makes an an angle at least $\frac{9 \pi}{10}$ with $\rho(\ell,a)$ at $\ell$.
This proves \eqref{2.4} and the description of the exceptions, and for this we can choose any
constant $\eta_0 \leq \eta_4$.

\ms
Let us now prove \eqref{2.5}, first modulo its exception.
Suppose that for some $\eta_0 < \eta_4/4$, we can find two arcs $\cC_i$
and $\cC_j$, that do not share an endpoint and furthermore are not connected
by one of the exceptions of \eqref{2.4}, but for which
\begin{equation} \label{2.26aa}
\dist(\cC_i, \cC_j) \leq \eta_0.
\end{equation}
We want to derive a contradiction. Let $x_i \in \cC_i$ and $x_j \in \cC_j$
be such that $|x_i - x_j| \leq \eta_0$. First assume that we can find endpoints
$a_i$ of $\cC_i$ and $a_j$ of $\cC_j$ such that $|a_i-x_i| \leq \eta_4/4$
and $|a_j-x_j| \leq \eta_4/4$. Then $|a_i-a_j| \leq 6 \eta_4/4$, 
\eqref{2.19} says that one of them (say, $a_i$ for definiteness) lies in $L$.
By \eqref{2.18}, $a_j$ is the the only point of $K \cap B(a_i, \eta_2/10) \sm L$. 
By Lemma~\ref{t2.7}, $\rho(a_i,a_j) \subset K$. This contradicts our assumption
that $\cC_i$ and $\cC_j$ were not connected by one of the exceptional arcs of
length $\leq \eta_4$.

So $x_i$, for instance, lies at distance at least $\eta_4/4$ from both endpoints 
of $\cC_i$. 

We are now going to follow the proof of \eqref{2.23}, and in particular 
proceed by contradiction.
Suppose that, for all $k \geq 0$, we can find a minimal cone $X_k$,
arcs $\cC_{i,k}$ and $\cC_{j,k}$ for which \eqref{2.26aa} holds with $\eta_0 = 2^{-k}$, 
and also points $x_{i,k} \in \cC_{i,k}$ and $y_{j,k} \in \cC_{j,k}$ 
such that $|x_{i,k}-y_{j,k}| = \dist(\cC_i, \cC_j)$, and yet $x_{i,k}$ lies at distance
at least $\eta_4/4$ from both endpoints of $\cC_i$.
By rotation invariance, we can work with a fixed $L$.
Set $r_k = |x_{i,k}-x_{j,k}| = \dist(\cC_{i,k},\cC_{j,k}) \leq 2^{-k}$.

We first assume that $r_k^{-1} \dist(x_{i,k}, L) \leq C$, and that there is a fixed 
$\ell \in \S \cap L$ such that $|x_{i,k} - \ell| \leq 2C r_k$. Set
$Y_k = r_k^{-1} (X_k - \ell)$, and extract a subsequence for which $Y_k$ tends to some $Y$. 
Then as usual $Y$ is a minimal set with sliding boundary $L$, $Y$ is invariant by translations 
parallel to $L$, and Lemma \ref{t2.3} gives a description of $Y$. 

Notice that $K_k$ contains an arc of geodesic $\rho_k$ of length $\eta_4/2$
centered at $x_{i,k}$. Let $P_k$ denote the plane that contains $\rho_k$ and the origin,
and let $D_k = P_k \cap B(x_{i,k},\eta_4/4)$. We know that $D_k \subset X_k$,
and hence $Y_k$ contains  $D'_k = r_k^{-1} (D_k - \ell)$, which is a 
planar disk parallel to $P_k$, but centered at $x'_{i,k} = r_k^{-1} (x_{i,k} - \ell)$,
and with a large radius $r_k^{-1} \eta_4/4$.

Notice that $|x'_{i,k}| = r_k^{-1} |x_{i,k}-\ell| \leq 2C$, hence we may assume 
that $x'_{i,k}$ has a limit $x'$, and that the direction of $P_k$ admits a limit too.
Then $D'_k$ has a limit, which is a whole plane $P'$ centered at $a'$ (because of the large radius), 
and in addition $P' \subset Y$ because $D'_k \subset Y_k$. 
By the description of Lemma \ref{t2.3}, $Y = P'$. 

Now consider the point $x_{j,k} \in K_k$. By definitions,
$\dist(x_{j,k}, \cC_{i,k}) = r_k$, and (since $\cC_{i,k}$ is orthogonal to $x_{i,k}-x_{j,k}$ 
at $x_{i,k}$ and by elementary geometry)
\begin{equation} \label{2.27}
\dist(x_{j,k}, D_k) \geq r_k/2.
\end{equation}
Set $x'_{j,k} = r_k^{-1} (x_{j,k} - \ell)$; then 
\begin{equation} \label{2.28}
\dist(x'_{j,k}, D'_k) = r_k^{-1} \dist(x_{j,k}, D_k) \geq 1/2.
\end{equation}
As before $|x'_{j,k}| = r_k^{-1} |x_{j,k}-\ell| \leq 2C+1$, so we may assume
that $x'_{j,k}$ tends to a limit $b$. Then $b\in Y$ because $x'_{j,k} \in Y_k$
for all $k$, and yet $\dist(b, P') \geq 1/2$ by \eqref{2.28}. 
This contradicts the fact that $Y = P'$, and we are left with the case where 
$r_k^{-1} \dist(x_{i,k}, L)$ is unbounded.

In this case, we set $Y_k = r_k^{-1} (X_k - x_{i,k})$, extract a subsequence for which
$Y_k$ converges to a limit $Y$, and notice that $Y$ is minimal, without a sliding boundary
condition. We may also assume that $x_{i,k}$ has a limit $x$, and then $Y$ is invariant
by translations in the direction of $x$. Thus $Y$ is a plane or a $\bY$-set.

We proceed as before, find disks $D_k \subset X_k$, then big disks $D'_k \subset Y_k$, 
and obtain that $Y$ contains a plane $P'$ (the limit of the $D'_k$). 
But at the same time, $x_{j,k}$ is far from $\cC_{i,k}$,
which leads to \eqref{2.27} and \eqref{2.28}. Again this is impossible, because
$Y$ contains a limit of the $x'_{j,k}$, which is at positive distance from $P'$.

This completes our proof of \eqref{2.5} with its exception: whenever
\eqref{2.26aa} fails, $\cC_i$ and $\cC_j$ are merely separated by a short arc
$\cC_k$ of $K$. So far we said that $\H^1(\cC_k) \leq \eta_4$,
but we also want to compute $\dist(\cC_i, \cC_j)$. 
Near $\cC_k$, the situation is the following.
We have the short arc $\cC_k$, with endpoints
$\ell \in L$ and $a \in K \sm L$; then $\cC_i$ and $\cC_j$ are two geodesic arcs
of length at least $\eta_2/20$, one leaving from $a$ with a $120^\circ$ angle
with $\cC_k$, and the other one leaving from $\ell$ with an angle of at least a $120^\circ$
with $\cC_k$ (we just applied \eqref{2.7}, but we could even get more by applying 
Lemma \ref{t2.7}). 

The standard case is when $\cC_i$ and $\cC_j$ are not too long,
and merely get away from each other when they leave $\cC_k$;
then $\dist(\cC_i, \cC_j) = \diam(\cC_k)$, and in particular $\cC_k$ was
also an exception of \eqref{2.5}, even with the smaller constant $\eta_0$.

The second case is when they get together again, near the antipode, so as to get
within $\eta_0$ from each other. Then they are merely separated by another
(in fact the only other) exceptional arc $\cC'_k$, and of course
$\dist(\cC_i, \cC_j) = \min(\diam(\cC_k),\diam(\cC'_k))$.

This completes our discussion of \eqref{2.5} with its exceptional case, 
and we are also finished with Proposition \ref{t2.1}.
\qed

\ms
We end this section with a short remark. Although we proved all our estimates by compactness,
this was mostly out of laziness. It is quite probable that we could get an explicit bound for $\eta_0$,
but we shall not try to do this here and doubt that it would be interestingly large.

\section{The standard decomposition of $K$}
\label{S3}

In this section we define the standard decomposition of a minimal set of dimension $2$
in $\R^n$, with sliding boundary $L$. This decomposition will be used to construct our main
competitors for the almost minimal set $E$ (in the next sections).
The full length property defined in the next section will use this decomposition as well.

Let $L \subset \R^n$ be a line, and $X$ be a sliding minimal cone with boundary $L$.
Recall that in Proposition \ref{t2.1} we defined a natural decomposition of $K = X \cap \S$,
into a finite and almost disjoint collection of sets $\cC_j$, $j\in J$, which are either 
arcs of geodesics, or full great circles, drawn on $\S$.

We modify this decomposition slightly, to get what we'll call a 
\underbar{standard decomposition} of $K$. In fact, we just take some
pieces $\cC_j$ and cut them into $2$, $3$, $4$, or $5$ pieces, so as to get
arcs of geodesics of length at most $\pi/2$. For the full great circles $\cC_j$,
we just cut them in $4$ equal parts, by adding four vertices. If
$\dist(\cC_j, L) < 1/4$, say, let us choose the two points of $\cC_j$ that lie closest to
$L$ as (two of the) cutting points. 

When $\cC_j$ is just an arc of geodesic and its length
is more than $\pi/2$, we cut it into sub-arcs of length between $\pi/4$ and $\pi/2$.
We may use the latitude that we have to choose the additional vertices as close to the points
$\ell \in L\cap S$ as possible (when $K$ does not already contain $\ell$), but the author does
not recall ever using this possibility after all.
For the moment, we do not care much if $\cC_j$ is not cut into equal pieces.

Once we cut all the long arcs as we just explained, we get a standard decomposition of $K$.
It is usually not unique, but this does not matter. It is again a decomposition of $K$
into arcs of geodesics $\cC_i$, $i \in \cI$, and we now review some of its properties. 
First, each $\cC_i$ is an arc of geodesic, 
\begin{equation} \label{3.1}
K = \bigcup_{i\in \cI} \cC_i \, ,
\end{equation}
and 
\begin{equation} \label{3.2}
\text{the interiors of the $\cC_i$, $i \in \cI$, are disjoint and do not meet $L$.} 
\end{equation}
Denote by $a_i$ and $b_i$ the endpoints of $\cC_i$. 
We do not pay attention to which is which, i.e., we may exchange the names at any time 
for the convenience of notation, but anyway
\begin{equation} \label{3.3}
\cC_i = \rho(a_i,b_i),
\end{equation}
where $\rho(a,b)$ systematically denotes the closed (shortest) geodesic from $a$ to $b$ in $\S$
(when $b \neq -a$). We will never take antipodal points, hence $\rho(a,b)$ will be uniquely
determined. As a general rule, $\rho$ will denote a geodesic or a union of geodesics.

Let us use this opportunity to introduce the geodesic distance in $\S$ defined by:
\begin{equation} \label{3.4}
\ddist(a,b) = \H^1(\rho(a,b)) \in [0,\pi];
\end{equation}
when $a=-b$ we set $\ddist(a,b) =\pi$.

We denote by $V$ the set of vertices of the standard decomposition, i.e., the 
collection of endpoints $a_i$ and $b_i$. We write
\begin{equation} \label{3.5}
V = V_0 \cup V_1 \cup V_2,
\end{equation}
where $V_0 = K \cap L$, $V_1$ is the set of vertices of the natural decomposition 
that do not lie in $V_0$, and $V_2$ is the set of vertices that we added to 
cut some of the initial arcs to make them shorter, and (for the case of full great circles)
that are not points of $V_0$. Thus the three $V_i$ are disjoint.

We said that the arcs $\cC_i$, $i \in \cI$, only meet at their endpoints, and there are
rules about how they can meet. For $\ell \in V_0$, there can only be one, two, 
or three $\cC_i$ that start from $\ell$, and always with angles $\geq \2$.
This comes from \eqref{2.7}. For $a\in V_1$, there are exactly three $\cC_i$ that start from $a$,
and they make angles of $\2$ at $a$ (see \eqref{2.6}). Finally, at $a\in V_2$, 
there are exactly two $\cC_i$ that start from $a$, and they make angles of $\pi$ at $a$
(their tangent half lines lie in opposite directions); this is clear, we just cut a geodesic at $a$.

We also control the length of the $\cC_i$. The general rule is that
\begin{equation} \label{3.6}
\frac{\pi}{4} \leq \H^1(\cC_i) \leq \frac{\pi}{2}
\ \text{ when at least one of the endpoints of $\cC_i$ lies in $V_2$;}  
\end{equation}
\begin{equation} \label{3.7}
\eta_0 \leq \H^1(\cC_i) \leq \frac{\pi}{2}  
\end{equation}
when both endpoints of $\cC_i$ lies in $V_1$;
and for $\ell \in V_0$, \eqref{3.7} holds for all the $\cC_i$ that end
at $\ell$, except perhaps one. When this exception happens, there is at most
one other arc $\cC_i$ that leaves from $\ell$, and it makes an angle at least
$\frac{9\pi}{10}$ with $\cC_i$ at $\ell$. See below \eqref{2.4}.

We also have the following consequence of \eqref{2.5} and the discussion that
follows it. When $i, j \in \cI$, are such that $\cC_i$ and $\cC_j$ don't have 
a common endpoint, the general rule is that
\begin{equation} \label{3.8}
\dist(\cC_i, \cC_j) \geq \eta_0,
\end{equation}
and this may only fail when there is an arc $\cC_k$ such that $\diam(\cC_k) \leq \eta_0$,
with one common endpoint with $\cC_i$ and one common endpoint with $\cC_j$.
Since $\cC_i$ and $\cC_j$ are now short (as in \eqref{3.7}), the strange special
case when $\cC_i$ and $\cC_j$ are close at both ends does not happen any more,
and we get the simpler formula $\dist(\cC_i, \cC_j) = \diam(\cC_k)$.
It will also good to know that for $i\in \cI$,
\begin{equation} \label{3.9}
\dist(x,K \sm \cC_i) \geq \min(\eta_0, |x-a_i|, |x-b_i|) \ \text{ for $x \in \cC_i$,} 
\end{equation}
where $a_i$ and $b_i$ still denote the endpoints of $\cC_i$. Indeed the distance 
to the direct neighbors of $\cC_i$ is controlled by our angle conditions, and the distance
to the other arcs $\cC_j$ is controlled by \eqref{3.8}, except when $\cC_i$ and
$\cC_j$ are separated by a short arc $\cC_k$. But even in this case, \eqref{3.9} follows
from the fact that $\cC_i$ and $\cC_j$ leave from $\cC_k$ in directions that make
an angle larger than $\frac{2 \pi}{3} - \frac{\pi}{5}$. See below \eqref{2.5}.

\section{The full length condition}
\label{S3b}

Now we define the full length condition, which will be our replacement for the 
epiperimetric inequality of Reifenberg. 
This will be a relatively simple condition on the position of the geodesics
that compose $K = X \cap \S$, which will be sufficient for our proof to give
good density decay, and then some regularity, at points of an almost minimal set where 
$X$ is a blow-up limit. See Definition \ref{t3.1}.

Let $X$ be a minimal cone, and choose a standard decomposition of $K = X \cap \S$, 
as in the previous section. We first discuss how to construct perturbations of $X$ by moving
the vertices $x\in V$.
We do not want to move them too much, because we want to modify the structure of $K$ 
as little as possible, and in order to measure how far we will be allowed to go we set
\begin{equation} \label{3.10}
\eta_L(X) = \min_{\ell \in L \cap \S \sm K} \dist(\ell, K) > 0,
\end{equation}
\begin{equation} \label{3.11}
\eta_V(X) = \min_{\ell \in  K \cap L} \dist(\ell, V_1 \cup V_2) > 0,
\end{equation}
and
\begin{equation} \label{3.12}
\eta(X) = 10^{-1} \min(\eta_0, \eta_L(X), \eta_V(X)),
\end{equation}
where the absolute constant $\eta_0$ comes from Proposition \ref{t2.1}.
Notice that when $\eta_V(X) < \eta_0$, it is the diameter of the smallest of the exceptional
arcs $\cC_i$ for \eqref{3.7}.

\ms
The basic tool to generate perturbations of $X$ is the set $\Phi_X(\eta)$
of mappings $\varphi : V \to \S$ such that
\begin{equation} \label{3.13}
|\varphi(x) - x | < \eta \ \text{ for } x\in V,
\end{equation}
and in fact we will restrict to $\eta < \eta(X)$.
We want to use the mappings $\varphi \in \Phi_X(\eta)$ to modify the geodesics $\cC_i$,
$i \in \cI$, and this will be easier when both endpoints of $\cC_i$ lie in $V_1 \cup V_2$.
Denote by $a_i$ and $b_i$ the endpoints of $\cC_i$, and set
\begin{equation}\label{3.14}
\cI_1 = \big\{ i\in \cI \, ; \, a_i \text{ and $b_i$ lie in } V_1 \cup V_2 \big\}.
\end{equation}
When $\varphi \in \Phi_X(\eta)$ and $i\in \cI_1$, we simply set
\begin{equation}\label{3.15}
\varphi_\ast(\cC_i) = \rho(\varphi(a_i),\varphi(b_i)).
\end{equation}
Things are a little more complicated when $i\in \cI_0 = \cI \sm \cI_1$.
When $i\in \cI_0$, we use the convention that $b_i \in V_0$, and $a_i \notin V_0$.
We want to leave more options, so we will need to append to $\varphi$ some additional
information. 

Let $\ell \in V_0$ be given, denote by $\cI(\ell) \subset \cI_0$ the set of indices $i\in \cI$
such that $\ell$ is an endpoint of $\cC_i$. Also call $m(\ell) \in \{ 1, 2, 3 \}$ the number
of elements of $\cI(\ell)$. If $m(\ell) = 1$, we don't need more information, and
in fact we can even forget about $\varphi(\ell)$, because we set
\begin{equation}\label{3.16}
\varphi_\ast(\cC_i) = \rho(\varphi(a_i),\ell)
\ \text{ when $i \in \cI(\ell)$, $m(\ell) = 1$, and $\cC_i = \rho(a_i,\ell)$.}
\end{equation}
When $m(\ell) = 2$, we add to $\varphi$ a component $\varphi_\ell \in \{ -1, 1\}$,
and we set
\begin{equation}\label{3.17}
\varphi_\ast(\cC_i) = \rho(\varphi(a_i),\varphi(\ell))
\ \text{ when $i \in \cI(\ell)$, $m(\ell) = 2$, $\varphi_\ell = -1$,
and $\cC_i = \rho(a_i,\ell)$}
\end{equation}
(the free option), and
\begin{equation}\label{3.18}
\varphi_\ast(\cC_i) = \rho(\varphi(a_i),\varphi(\ell)) \cup \rho(\varphi(\ell),\ell)
\ \text{ when $i \in \cI(\ell)$, $m(\ell) = 2$, $\varphi_\ell = 1$,
and $\cC_i = \rho(a_i,\ell)$}
\end{equation}
(the attached option). In this case, we added the same connecting arc $\rho(\varphi(\ell),\ell)$
to the two $\varphi_\ast(\cC_i)$, $i\in \cI(\ell)$, but this is just to avoid more complicated
notation, and we will never count this arc with multiplicity.

When $m(\ell) = 3$, we add a component $\varphi_\ell \in \{ -1 \} \cup \cI(\ell)$
(i.e., choose the free option or the attached option, and in this last case choose one
of the three $\cC_i$ that end at $\ell$, and set
\begin{equation}\label{3.19}
\varphi_\ast(\cC_i) = \rho(\varphi(a_i),\varphi(\ell))
\ \text{ when $i \in \cI(\ell)$, $m(\ell) = 3$, $\varphi_\ell = -1$,
and $\cC_i = \rho(a_i,\ell)$}
\end{equation}
(as before, the free option where we just move the center and let the $\cC_i$ follow).
In the last case when $m(\ell) = 3$, and $\varphi_\ell = j\in \cI(\ell)$,
we set
\begin{equation}\label{3.20}
\varphi_\ast(\cC_i) = \rho(\varphi(a_i),\ell) \ \text{ if } i = j
\end{equation}
and, for the two other indices $i \in \cI(\ell) \sm \{ j \}$
\begin{equation}\label{3.21}
\varphi_\ast(\cC_i) = \rho(\varphi(a_i),\varphi(\ell)) \cup \rho(\varphi(\ell),\ell).
\end{equation}
Again we put the same arc $\rho(\varphi(\ell),\ell)$ twice when once would have been enough,
and this time we transformed the union of the three arcs $\cC_i$ that looks like a $Y$
into a truncated $Y$ plus an arc, both leaving from $\ell$.
Notice that when $\varphi(\ell) = \ell$, some of the pictures above get simpler, and
we don't even need the free option.

Since we may have added coordinates to $\varphi$, let us denote
by $\Phi^+_X(\eta)$ the set of enlarged mappings $\varphi$. We do not give a different name
to $\varphi$ and its extension, so as not to exaggerate with notation; when we really want
to know which one we consider, we will say that $\varphi \in \Phi_X(\eta)$ or 
$\varphi \in \Phi_X^+(\eta)$,
but the truth is that we shall work with $\Phi_X^+(\eta)$.

For $\varphi \in \Phi_X^+(\eta)$, we define a perturbation of $K$ by
\begin{equation}\label{3.22}
\varphi_\ast(K) = \bigcup_{i \in \cI} \varphi_\ast(\cC_i) \subset \S,
\end{equation}
and then a modified cone
\begin{equation}\label{3.23}
\varphi_\ast(X) = \big\{ tx \, ; \, x\in \varphi_\ast(K) \text{ and } t\in [0,+\infty) \big\}.
\end{equation}

With the present definition, it may happen that even after we remove the arcs 
$\rho(\ell,\varphi(\ell))$ that we counted twice, some of our arcs $\varphi_\ast(\cC_i)$
still cross (i.e., meet somewhere else than their common endpoints). 
Since we took $\eta < \eta(X)$, the arcs of \eqref{3.16} don't do that, 
so all the crossing happens near the points of $V_0$.

Let us decide to forbid this, and restrict to the subset $\Phi_X^{+,i}(\eta) \subset \Phi_X^+(\eta)$
of mappings $\varphi$ for which this does not happen, i.e., for which the arcs $\varphi_\ast(\cC_i)$ 
are disjoint, except for common endpoints and for the double occurrence of some
$\rho(\ell,\varphi(\ell))$. We will call the mappings $\varphi \in \Phi_X^{+,i}(\eta)$ injective.

The main reason for restricting to injective mappings $\varphi$ is that later in the proof,
the only cones $\varphi_\ast(X)$ for which we need to use the full length condition
below are injective anyway, so we save ourselves a little bit of unpleasant verification
for the full length condition, at essentially no expense.

We are finally ready to say what we mean by full length.

\begin{defi}\label{t3.1} 
We say that the minimal cone $X$ (with sliding boundary $L$) satisfies the 
\ub{full length} condition when there is a standard decomposition of $K = X \cap \S$
and small numbers $\eta \in (0, \eta(X))$ and $c > 0$ such that for all 
injective mappings $\varphi \in \Phi_X^{+,i}(\eta)$ such that
\begin{equation}\label{3.24}
\Delta(\varphi) := \H^1(\varphi_\ast(K)) - \H^1(K) > 0,
\end{equation}
there is a sliding competitor $\wt X$ for $\varphi_\ast(X)$ in $\ol B(0,1)$ 
(see Definition \ref{t1.1}) such that
\begin{equation}\label{3.25}
\H^2(\wt X \cap B(0,1)) \leq \H^2(\varphi_\ast(X)\cap B(0,1)) - c \Delta(\varphi).
\end{equation}
\end{defi}

\ms
This looks complicated (just as the initial definition of full length for plain minimal cones,
see \cite{C1}), and the only justification for it is that it makes the machine work.
At least we shall be able to check the full length condition on the simple examples that 
we mostly care about, and we can always think that the definition is simpler than
some notions of epiperimetric inequalities.

One paradox that may be worth mentioning is that the full length condition only makes sense
once we know that $X$ is minimal. Otherwise, finding $\wt X$ so that \eqref{3.25} is just too
easy: use a better competitor of $X$ and deform it a little to fit $\varphi_\ast(X)$.

Probably choosing a standard decomposition by force, instead of allowing some flexibility
as we did, would not change the notion of full length. But having the choice will allow us to
simplify some computations when we check the full length.

When we choose the free option in the description of 
$\varphi_\ast(K)$ near a point $\ell \in K$, we seem to save some area when we 
omit to add $\rho(\ell,\varphi(\ell))$, but at the same time we allow more competitors 
for \eqref{3.25}, because we don't necessarily need to check \eqref{1.7} near $\ell$
any more. Requiring the condition of Definition \ref{t3.1} also for the deformations with the
free option seems to be more stable. Think about the case when $X$ is a plane, 
which may contain $L$ or just pass very close to $L$.

We will have to return to the notion of full length later, and play a little more with the definitions.
First, there may be circumstances where we will need to check the existence of $\wt X$
only for perturbation that are free (we will also say detached) at one or two of the points 
$\ell \in V_0$, typically because anyway the set $E$ that we study is detached from $L$. 

Also, we will need a similar notion of full length when we prove decay estimates, in some cases, for the
functional $F$ of \eqref{1.20} adapted to balls that are not centered exactly on $L$. This will happen at the end of Section \ref{S25}, with proofs in Sections \ref{S26} and \ref{S27}.

Finally, instead of requiring the existence of a sliding competitor $\wt X$ such that \eqref{3.25} holds,
we will some times be able to manage with a simpler version of this. See the end of Section \ref{S25}, and in particular Lemme \ref{t25.1}.

The cones of type $\bP$, $\bY$, $\bT$, $\bH$, and $\bV$, all satisfy the full length property;
a good part of this will be checked in Sections \ref{S26} and \ref{S27}, 
when we need full length estimates related to $F$ in the more general case of balls that are not
necessarily centered on $L$, and the remaining verifications will be done in Section \ref{S30}.
See Theorem \ref{t30.1}.

\vfill\eject
\part{Density decay for balls centered on $L$}

\section{The initial setup and two words about the constants}
\label{S4}

A very large portion of this paper consists in the construction of some competitors
for an almost minimal set $E$. In this section we give some notation and our basic assumptions
for the sections that follow.

We work in $\R^n$, and with a reduced (or coral) sliding almost minimal set $E$ of dimension $2$, 
in an open set $U$ that contains $B(0,2r)$, and with a sliding condition coming from the line 
$L$ through the origin. See Definitions \ref{t1.2} and \ref{t1.2c}.
We shall assume that the gauge function $h$ is small enough, and
more precisely that
\begin{equation}\label{4.1}
h(s) \leq C_h s^{\beta} \ \text{ for } 0 < s \leq 2r,
\end{equation} 
where $\beta > 0$ and $C_h \geq 0$ are constants such that
\begin{equation} \label{4.2}
C_h r^{\beta} \leq \varepsilon,
\end{equation}
where $\varepsilon > 0$ is a very small constant that will be chosen much later,
to make the argument work. 
This is a little more restrictive than what we actually need. 
For instance taking $h(s) = \varepsilon \big(\ln(2r/s)\big)^{-b}$ for some sufficiently large $b$ 
that depends on $n$ would be enough, by looking at similar statements in \cite{C1} and checking
that they adapt.
In what follows, we could decouple \eqref{4.1} from \eqref{4.2}; in effect, we shall use 
\eqref{4.1} because it implies some nice general properties for $E$, like the fact that
it has a $C^1$ description far from $L$, or the existence of the density $\theta(x)$ of \eqref{1.2}. 
Then the estimates of the next sections will use (these general properties, plus) the size of $h(2r)$.

We also assume that
\begin{equation}\label{4.3}
d_{0,2r}(E,X) \leq \varepsilon,
\end{equation}
where $d_{0,2r}$ is the local Hausdorff distance of \eqref{1.13} and $X$ is a sliding 
minimal cone (centered at $0$), also with a sliding boundary condition coming from $L$.

Our main task will be to construct, under various assumptions on $E$ and $r$,
some good competitors for $E$ in the closure of $B = B(0,r)$; then 
we will use this to get differential inequalities on a density excess function $f(r)$, 
in principle associated to the standard monotonicity formula from \cite{Sliding},
although later on we also want to use a slightly different monotonicity formula
from \cite{Mono}, with balls that are not centered on $L$.

We assume some additional properties of $E$, which will be used in
the proof. First assume that 
\begin{equation}\label{4.4}
\H^{1}(E \cap \d B) < +\infty.
\end{equation}
This is true for almost every $r > 0$ such that $B(0,r) \subset U$ (for instance by 
the co-area formula), so it costs virtually nothing, and will be useful in some proofs.
Next we require some maximal function to be small at $r$. Define a measure $\mu_\pi$
on $[0,2r)$ by 
\begin{equation}\label{4.5}
\mu_\pi(A) = \H^d(E \cap \pi^{-1}(A))
\end{equation}
for Borel sets $A \subset [0,2r)$, and where $\pi$ is the radial projection defined by
\begin{equation}\label{4.6}
\pi(z) = |z| \ \text{ for } z\in \R^n.
\end{equation}
We require that there is a number $C \geq 0$, that may depend wildly on $r$, such that
\begin{equation}\label{4.7} 
\H^2(E \cap A_\xi) \leq C \xi 
\ \text{ for } 0 < \xi < r, \ \text{ where we set } A_\xi = B(0,r) \sm B(0,r-\xi).
\end{equation}
That is, we just require the one-sided variant of the Hardy-Littlewood maximal 
function of $\mu_\pi$ to be finite at the point $r$. This is like (4.5) in \cite{C1},
and we shall see in the proof of Proposition \ref{t16.2}
that \eqref{4.7} holds for almost every $r$.

We also require that for every continuous nonnegative function $f$ on $\R^n$,
\begin{equation}\label{4.8}
\lim_{\rho \to 0} \rho^{-1}\int_{t \in (r-\rho,r)} \int_{E \cap \d B(0,t)} f(z) d\H^1(z) dt
= \int_{E \cap \d B(0,r)} f(z) d\H^1(z).
\end{equation}
This is the same thing as (4.3) in \cite{C1}, and it turns out that this is also satisfied for 
almost every $r > 0$. The proof is given in Lemma 4.12 of \cite{C1}, 
and works here as well. So assuming \eqref{4.8} for our standard $r$ 
costs us nothing.

\ms
Starting with the next section, we shall fix $E$, $X$, and $r$ as above,
and even assume that $r=1$ to simplify the notation. Our general goal is to
modify the set $E$ inside of $\B = B(0,1)$ to get a better competitor if we can.

The construction of the good competitor will keep us busy for 
Sections \ref{S5}-\ref{S15}. It will use a few different constants, and 
maybe it is the right time to announce in which order they will be chosen.

We already have a constant $\eta(X)$, which may be very small
(depending on $X$, and in particular on the distance between $L$
and some faces or edges of $X$), but we see it as a geometric constant.

In Section \ref{S5}, we introduce a small constant $\tau$, which gives
the size of the disks $D$ near the points of $\S \cap L$ where most of the 
action will take place. 

Then there is a small Lipschitz constant $\lambda$, which we use to construct
Lipschitz graphs in Section \ref{S7}.
We will need $\lambda$ to be small enough (we often use it as a small parameter
to control some angles), and in particular so that the estimates of Section \ref{S8} apply.

For a long time, the only constraint on $\tau$ will be to be small enough,
in particular compared to $\eta(X)$ (see \eqref{5.2bis}), but for the estimates
in Section \ref{S13} to give small enough errors compared to what we win with
the construction of Section \ref{S13}, we will need $\tau$ to be small enough,
depending on $\lambda$.

There will be a short occasion or two, in Section \ref{S10}, where we briefly
use a smaller value of $\tau$, but this will be explained then. See Lemma \ref{t10.1}.
At the same time, we will use a small $\alpha > 0$, which will be chosen in
Section \ref{S10} and will depend on $\lambda$ and $\tau$.

There is a small constant $\tau_1$ in \eqref{7.3}, which we may
as well take very small, compared to both $\lambda$ and $\tau$. 
But it is not exactly of the same nature as $\tau$.

Our last real constant is $\varepsilon$ in \eqref{4.3}, which will be chosen at the end, 
extremely small, and depending on all the constants above.

We mention $\tau_3$ (in Proposition \ref{t12.1}) and $\tau_4$ 
(in Lemma \ref{t13.1}) for completeness, but they are extremely small numbers,
not constants, as they may depend on the radius $r$ above.

\section{A local description of $E \cap \S$ far from $L$}
\label{S5}

From now on we fix the line $L$, the minimal cone $X$,
the reduced almost minimal set $E$,  and the radius $r=1$, as in Section \ref{S4}.

In this section we first record simple properties of the $\cC_i$
concerning their distances (see Lemma \ref{t5.1}), 
and then use known local regularity results for plain almost minimal sets 
(i.e., with no sliding boundary condition) to give, at least far from $L$, 
a local description of $E \cap \S$ as a finite union of $C^1$ curves $\cL_i$ that follow 
the curves $\cC_i$ of the standard description of $K$. The reader may want to check
Proposition~\ref{t5.4} below to convince herself that no real surprise will come out of 
this section. 

The description of  $E \cap \S$ that we'll give in this section relies on regularity results
for (plain) almost minimal sets. We shall quote \cite{Holder} and \cite{C1} for convenience,
but when $n=3$ we could as well use \cite{Ta}. Also, we proceed this way because
we prefer to insist right away on the places where new difficulties appear (i.e., close to $L$);
with a more complicated version of the present paper, we would quite probably be able to
prove the local regularity of $E$ far from $L$ at the same time. 
But this would not really make things simpler: the proof of the present paper essentially 
contains the proof in \cite{C1} anyway.

\ms
We start with more notation. 
Recall from Section \ref{S3} that the standard decomposition of $K$ is composed of arcs 
$\cC_i$, $i\in \cI$, of geodesics, and that $a_i$, $b_i$ denote the endpoints of $\cC_i$,
so that $\cC_i = \rho(a_i, b_i)$. Set
\begin{equation} \label{5.1}
\cI_0 = \big\{ i \in \cI \, ; \, \cC_i \text{ meets } L \big\}
= \big\{ i \in \cI \, ; \, \text{ $a_i \in L$ or $b_i \in L$} \big\}.
\end{equation}
Also denote by $\cI_1 = \cI \sm \cI_0$ the set of $i\in \cI$ such that $a_i, b_i \in V_1 \cup V_2$.

In this section will stay at some distance from $L$. More precisely, 
denote by $\ell_+$ and $\ell_-$ the two points of $L \cap \S$. Then set 
\begin{equation} \label{5.2}
D_\pm(\tau) = \S \cap B(\ell_\pm,\tau)
\ \text{ and } \ 
\d D_\pm(\tau) = \S \cap \d B(\ell_\pm,\tau)
\end{equation}
for $\tau \leq \eta(X)$. We will use various small numbers $\tau > 0$, to be chosen later, 
but always such that
\begin{equation} \label{5.2bis}
\tau \leq 10^{-3} \eta(X).
\end{equation}

\subsection{More about distances between vertices and curves $\cC_i$}

Before we come to $E \cap \S \sm (D_+(\tau) \cup D_-(\tau))$, we need some 
additional information on $V$ and the $\cC_i$.

\begin{lem} \label{t5.1}
Let $\eta(X)$ be as in \eqref{3.12}. Then
\begin{equation} \label{5.3}
\H^1(\cC_i) \geq 10 \eta(X) \ \text{ for } i \in \cI,
\end{equation}
\begin{equation} \label{5.4}
\dist(\cC_i, \cC_j) \geq 9 \eta(X) \ \text{ for $i, j \in \cI$ such that } 
\cC_i \cap \cC_j = \emptyset,
\end{equation}
and 
\begin{equation} \label{5.5}
\dist(x, \cC_i) \geq 9 \eta(X) \ \text{ when $x \in V$ and the arc
$\cC_i$ does not contain $x$.}
\end{equation}
\end{lem}

We start with a proof of \eqref{5.3}.
When both endpoints of $\cC_i$ lie on $V_1 \cup V_2$, this
follows from \eqref{3.7} and \eqref{3.12}. Otherwise, 
$\diam(\cC_i) \geq \eta_V(X) \geq 10\eta(X)$ by \eqref{3.11} and \eqref{3.12};
\eqref{5.3} follows.

Next we check \eqref{5.4}. Suppose that $\cC_i \cap \cC_j = \emptyset$. In general,
\eqref{5.4} follows from \eqref{3.8} and \eqref{3.12}, and the only exception, as explained
below \eqref{3.8}, is when $\cC_i$ and $\cC_j$ are only separated by a short arc $\cC_k$. 
Even so, we said below \eqref{3.8} that
$\dist(\cC_i, \cC_j) = \diam(\cC_k)$. But \eqref{5.3} says that $\H^1(\cC_k) \geq 10 \eta(X)$,
and since we may safely assume that $\eta(X) \leq \eta_0$ is very small, we get that
$\diam(\cC_k) \geq 9\eta(X)$, as needed.

Finally let $x \in V$ and $\cC_i$ be as in \eqref{5.5}. 
First assume that there is an arc $\gamma$ of $K$ that goes from $x$ to $\cC_i$;
then $\H^1(\gamma) \geq 10 \eta(X)$ by \eqref{5.3}, hence
$\diam(\gamma) \geq 9\eta(X)$. In this case
$\dist(x, \cC_j) \geq \diam(\gamma) \geq 9\eta(X)$, because
$\gamma$ and $\cC_j$ make an angle of at least $120^\circ$ at their
common endpoint.

Assume now that $x$ is not directly connected to $\cC_i$, and let $\cC_j$ be any arc that 
contains $x$; then $\cC_i \cap \cC_j = \emptyset$, and 
$\dist(x, \cC_j) \geq \dist(\cC_i, \cC_j) \geq 9 \eta(X)$ by \eqref{5.4};
\eqref{5.5} and Lemma~\ref{t5.1} follow.
\qed

\subsection{A description of $E \cap \S$ near a vertex $x\in V_1$}

We start our description of $E \cap \S$ with what happens near the vertices of $V_1$.
We fix $x\in V_1$, and denote by $\gamma_1$, $\gamma_2$, $\gamma_3$ 
the three arcs $\cC_i$ that leave from $x$.
Recall that the $\gamma_i$ meet with $120^\circ$ angles, so there is a cone $Y(x) \in \bY$,
whose spine (understand, singularity set) contains the line through $x$, and
that contains $\gamma = \gamma_1 \cup \gamma_2 \cup \gamma_3$. 
Let us check that
\begin{equation} \label{5.6}
K \cap B(x,9\eta(X)) = \gamma \cap B(x,9\eta(X)) = Y(x) \cap \S \cap B(x,9\eta(X)).
\end{equation}
Since $\gamma \subset K$, for the first identity we just need to check that the only $\cC_i$
that meet $B(x,9\eta(X))$ are the $\gamma_i$. This follows from \eqref{5.5},
because the $\gamma_i$ are the only arcs that contain $x$.
For the second identity, observe that $\gamma \subset Y(x)$ by definition of $Y(x)$;
the other inclusion holds because $\H^1(\gamma_i) \geq 10 \eta(X)$ by \eqref{5.3}.

The next description will come from Corollary 12.25 in \cite{C1}.

\ms
\begin{lem} \label{t5.2} Suppose that $\tau \leq 10^{-3} \eta(X)$.
If $\varepsilon$ is small enough, depending on $n$ and $\tau$,
there is a $C^1$ diffeomorphism $\Phi : B(x,20\tau) \to \Phi(B(x,20\tau))$ with 
the following properties. First 
\begin{equation} \label{5.7}
|\Phi(y) - y| \leq 10^{-10} \tau \ \text{ for } y\in B(x,20\tau);
\end{equation}
\begin{equation} \label{5.8}
|D\Phi(y) - Id| \leq 10^{-2} \ \text{ for } y\in B(x,20\tau);
\end{equation}
\begin{equation} \label{5.9}
E \cap B(x,9\tau) = \Phi(Y(x)\cap B(x,20\tau)) \cap B(x,9\tau).
\end{equation}
\end{lem}

\ms
Some preparation will be needed before we apply Corollary 12.25 in \cite{C1} to get this.
First observe that because of \eqref{3.10}-\eqref{3.12},
$\dist(x, \S \cap L) \geq 10 \eta(X)$; thus 
\begin{equation} \label{5.10}
\text{$E$ is a plain almost minimal set in $B(x,9 \eta(X))$, with gauge function $h(s) = C_h s^\beta$}
\end{equation} 
(see \eqref{4.1} too). Next we check that $E$ is close to $Y(x)$ near $x$. First
observe that
\begin{equation} \label{5.11}
X \cap B(x,8\eta(X)) = Y(x) \cap B(x,8\eta(X)).
\end{equation}
Indeed $B(x,8\eta(X))$ is contained in the cone over $\S \cap B(x,9\eta(X))$.
By \eqref{5.6}, the two cones $X$ and $Y(x)$ coincide on $\S \cap B(x,9\eta(X))$;
then they also coincides in $B(x,8\eta(X))$, as needed. It now follows from
\eqref{4.3} (and if $\varepsilon$ is small enough) that
\begin{equation} \label{5.12}
d_{x,8\eta(X)}(E,Y(x)) = \frac{2}{8 \eta(X)}\, d_{0,2}(E,X)
\leq \frac{\varepsilon }{4 \eta(X)}
\end{equation}
(recall that $r=1$ here).

Next we find a point of type $Y$ near $x$: we claim that there exists
$x_0 \in E$ such that
\begin{equation} \label{5.13}
\theta(x_0) = 3\pi \ \text{ and } \ |x_0 - x| \leq  C 
\varepsilon.
\end{equation}
Here the density $\theta(x_0) = \lim_{r \to 0}\theta(x_0,r)$
is defined by \eqref{1.1} and \eqref{1.2}, and the first condition is another way to
say that all the blow-up limits of $E$ at $x_0$ lie in $\bY$.

To find $x_0$ we apply Proposition 16.24 in \cite{Holder} to $E$ and the small ball 
$B(x,r)$, where $r$ will be chosen in a second. 
Let $\varepsilon_2$ denote the small constant in that proposition;
we choose $r = 10 \varepsilon_2^{-1} \varepsilon$; most assumptions
are satisfied readily (for instance, \eqref{4.2} takes care of the size of the gauge function);
the main one is that $d_{x,r}(E, Y(x)) \leq \varepsilon_2$, and it follows from 
the second part of \eqref{5.12} and our choice of $r$. 
The conclusion of Proposition 16.24 in \cite{Holder}
is that $E \cap B(x, 10^{-2}r)$ contains a point $x_0$ of type $Y$.
This point satisfies \eqref{5.13}, with $C = (10\varepsilon_2)^{-1}$.

We also claim that $E$ is close to $Y(x)$ in measure; actually we shall just need to know that
for each $\varepsilon_1 > 0$, 
\begin{equation} \label{5.14}
\H^2(E \cap B(x_0,1100\tau))
\leq \H^d(Y(x) \cap B(x_0,1100(1+\varepsilon_1)\tau)) + \varepsilon_1 (1100\tau)^2.
\end{equation}
if $\varepsilon$ is small enough. To see this, apply Lemma 16.43 in \cite{Holder}
to the almost minimal sets $E$ and $Y(x)$, in the ball $B(x_0,1100\tau)$,
and with $\delta = \varepsilon_1$. Since $\tau \leq 10^{-3} \eta(X)$, \eqref{5.10}
gives ample room to do this. The main assumption, that
$d_{x_0, 11000\tau/9}(E, Y(x))$ be small, follows from \eqref{5.12} if $\varepsilon$
is small enough (depending on $\tau$ and $\varepsilon_1$).

Notice that $\rho^{-2} \H^d(Y(x) \cap B(x_0, \rho)) \leq 3\pi$ for all $\rho > 0$,
for instance because the left-hand side is a nondecreasing function of $\rho$
(recall that $Y(x)$ is a minimal set), and tends to $3\pi$ at $\infty$. 
Thus \eqref{5.14} yields
\begin{equation} \label{5.15}
\theta(x_0, 1100\tau) = (1100\tau)^{-2} \H^2(E \cap B(x_0,1100\tau)) 
\leq 3\pi + C \varepsilon_1.
\end{equation}

We are now ready to apply Corollary 12.25 in \cite{C1}, to the set $E - x_0$
(because the corollary applies to a ball centered on the set), 
and with the radius $r_0 = 10\tau$. 

There are a few assumptions to check. First, \eqref{5.10} says that $E-x_0$
is almost minimal (with no sliding condition) in $B(0,110r_0)$, because 
$\tau \leq 10^{-3} \eta(X)$ and by \eqref{5.13}.
Next there are assumptions on the size of the gauge function, in particular evaluated
at $r_0$; these are satisfied if $\varepsilon$ in \eqref{4.2} is small enough.
Then there is the assumption (12.27) on the distance from $E-x_0$ to a minimal cone.
We now that $E-x_0$ is $2\varepsilon$-close to $Y(x)-x_0$ (by \eqref{5.12}), but since 
Corollary 12.25 in \cite{C1} requires a minimal cone centered at the origin
and $Y(x)-x_0$ is centered at $x-x_0$, we translate it by $x_0-x$ and get a
minimal cone $Y = Y(x)-x$ centered at $0$. Fortunately $|x_0-x|$ is as small as we wish 
(use \eqref{5.13} and take $\varepsilon$ small), 
so $Y$ is as close to $E-x_0$ as we wish in $B(x_0,110r_0)$.

Finally we need to check (12.27); half of it concerns the size of the gauge function 
and follows from \eqref{4.1}--\eqref{4.2}, and the other half requires that the density
excess $f(110r_0)$ be sufficiently small. 
With the definitions of \cite{C1} (see (3.5) et (3.2) there),
\begin{equation} \label{5.16}
f(110 r_0) = \theta(x, 1100 \tau) - \theta(x_0) 
= (1100\tau)^{-2 }H^2(E \cap B(x_0, 1100\tau)) - 3\pi,
\end{equation}
by definitions and \eqref{5.13}. This is as small as we want, by \eqref{5.15}
(we choose $\varepsilon_1$ small, depending on the constants in Corollary 12.25,
then take $\varepsilon$ small). So we can apply the corollary.

The conclusion (just applied to $r = 10\tau \leq r_0$) is that there is a 
$C^{1+\beta_1}$ diffeomorphism $\Psi : B(0,20\tau) \to \Psi(B(0,20\tau))$
(for some small $\beta_1 > 0$ that depends on the $\beta > 0$
of \eqref{4.1} and \eqref{4.2}) such that 
\begin{equation} \label{5.17}
\Psi(0) = 0, \, |\Psi(y)-y| \leq 10^{-2}\tau \ \text{ for } y \in B(0,20\tau),
\end{equation}
and 
\begin{equation} \label{5.18}
(E-x_0) \cap B(0,10\tau) = \Psi(Y\cap B(0,20\tau)) \cap B(0,10\tau).
\end{equation}
We take 
\begin{equation} \label{5.19}
\Phi(z) = x_0 + \Psi(z-x) \ \text{ for } z\in B(x,20\tau)
\end{equation}
to translate things back; notice that this way 
\begin{equation} \label{6.21n}
\Phi(x) = x_0 + \Psi(0) = x_0.
\end{equation}
Let us see what we get. We start with the good news: $\Phi$ is well defined on
$B(x,10\tau)$, and is a $C^{1+\beta_1}$ diffeomorphism whose image is
$x_0 + \Psi(B(0,20\tau))$, which is almost the same as $B(x_0,20\tau)$. 
Next \eqref{5.9} holds, because
\begin{equation} \label{5.20}
\begin{aligned}
E \cap B(x,9\tau) &= x_0 + [(E-x_0) \cap B(x-x_0,9\tau)]
\cr&= x_0 + [\Psi(Y\cap B(0,20\tau)) \cap B(x-x_0,9\tau)]
\cr&= [x_0 + \Psi(Y\cap B(0,20\tau))] \cap B(x,9\tau)
\cr&= \Phi(x + Y\cap B(0,20\tau))] \cap B(x,9\tau)
\end{aligned}
\end{equation}
by \eqref{5.18}. But $x+Y = Y(x)$, so $x + Y\cap B(0,20\tau) =
Y(x) \cap B(x,20\tau)$; \eqref{5.9} follows.

Now \eqref{5.17} only yields $|\Phi(z)-z| \leq |(x_0 +\Psi(z-x)-z| 
\leq |x_0 - x| + |\Psi(z-x)-(z-x)| \leq  10^{-2}\tau + C\varepsilon$
for $z \in B(x,20\tau)$, while we announced $|\Phi(z)-z| \leq 10^{-10} \tau$
in \eqref{5.7}. We could try to reduce the difference
and keep same proof by taking $r_0 = 10^{-8}\tau$, but it is more honest to say that 
the constant $10^{-2}$ in Corollary 12.25 can be replaced by any small number we wish, 
at the only expense of taking $\varepsilon$ much smaller. 
In the construction of a Reifenberg parameterization,
this amounts to starting to move points only after a certain number of generations; the
price to pay is precisely to force the initial set $Y$ to be close enough to $E$ so that
we still get a good enough approximation at the scale where we really start things.

The second difference is that we announced $|D\Phi - Id| \leq 10^{-2}$, and 
Corollary 12.25 only says that $\Psi \in C^{1+\beta_1}$. 
But in fact it is a uniform $C^{1+\beta_1}$ estimate, which means that we even get a 
uniform control on $(|y-z|/r)^{-\beta_1} |D\Phi(y)-D\Phi(z)|$ for $y,z \in B(x,20\tau)$. 
With this, a very tight uniform control on $|\Phi(y)-y|$ (take it even better than \eqref{5.13}, 
since it costs nothing), and some interpolation, we rather easily get \eqref{5.14}.
Another way to put this is to notice that the proof of existence for the Reifenberg
parameterization also gives a derivative which is as close to the identity as we want, 
again if $\varepsilon$ is taken small enough.
Finally, what really matters to us is the fact that if $y, z$ lie in $B(x,20\tau)$
and on the same face of $E$, then the distance between the directions of the tangent planes
to $E$ at $y$ and $z$ is at most $10^{-3}$, say. This is what we prove in \cite{C1},
in estimates like Lemma 12.35 and 12.50 (where we can take $\varepsilon_0$ as small as we want),
which show that approximate minimal cones vary slowly.
Hopefully the reader will trust one of these arguments;
this completes our proof of Lemma \ref{t5.2}.
\qed

\begin{rem}\label{r6.3n}
We can replace the constants $10^{-10}$ in \eqref{5.7} and $10^{-2}$ in \eqref{5.8}
by any small number $a_0>0 $ that we wish, but then $\varepsilon$ has to be taken small enough,
depending on $n$ and $\tau$ as above, but also on $a_0$. The proof is the same; as explained above,
we just need to know that in Corollary 12.25 of \cite{C1}, the mapping $\Psi$ that we get can be required
to be close enough to the identity in $C^1$ norm. Using this remark, we will be given the opportunity of
simplifying our construction slightly at the beginning of Section \ref{S12}. 
\end{rem}

\ms
Let us now say why Lemma \ref{t5.2} also gives a good control on 
$E \cap \S \cap B(x,8\tau)$. Some more notation will be useful.
Denote by $F_1$, $F_2$, and $F_3$ denote the three faces of $Y(x)$,
and choose the labels so that $F_j$ is the half plane that contains 
$\gamma_j$ (and is bounded by the spine, which is the vector line through $x$). 
Then \eqref{5.9} gives a decomposition of $E \cap B(x,9\tau)$ into three
faces $F'_j = B(x,9\tau) \cap \Phi(F_j \cap B(x,20\tau))$, which are
intersections with $B(x,9\tau)$ of three pieces of $C^1$ surfaces.
Each of these surfaces $F'_j$ is also a small Lipschitz graph (over 
a part of the plane that contain $F_j$), by \eqref{5.8}. Finally we know that
the three $F'_j$ meet along a $C^1$ curve (a piece of the image by $\Phi$ of the spine
of $Y(x)$), which is also a small Lipschitz curve, and they meet with $120^\circ$
angles. This is really the description of $E$ near $x$ that we shall use; the fact
that we have a parameterization by $Y(x)$ is of lesser importance.

We shall use the description above in a later section, 
but for the moment we only care about its consequence on $E \cap \S$.
For this we apply the implicit function theorem to each face $F'_j$, to get
a description of $F'_j \cap \S$ (the zero set of $f(z) = |z|^2-1$). The relevant
derivative is the derivative of the distance to the origin, which stays large
because $F'_j$ is a small Lipschitz graph over a plane that is orthogonal to $\S$.
Or equivalenty, we could use \eqref{5.8} to estimate the partial derivative
of $f \circ \Phi$ in the direction of $x$. We get that
\begin{equation} \label{5.21}
F'_j \cap \S \cap B(x,9\tau) \ \text{ is a $C^{1+\beta_1}$ curve.} 
\end{equation}
Call this curve $\cL_j$. We have a little more information on the $\cL_j$.
First,
\begin{equation} \label{5.22}
\text{ the three  $\cL_j$ start from a same point $x^\ast$,}
\end{equation}
which is also the unique point of $F'_1 \cap F'_2 \cap F'_3 \cap \S$
(or the only point of $\Phi(F_1\cap F_2 \cap F_3)$ that lies on $\S$,
apply the implicit function theorem to that curve). 

Next the $\cL_j$ are small Lipschitz graphs. Let us state this in terms of
the oscillation of their unit tangent direction. For $z\in \cL_j$, denote
by $v_j(z)$ a unit tangent vector to $\cL_j$ at $z$. We define $v$ so that
it is continuous, and $v(x^\ast)$ points in the direction of $\cL_j$.
Also denote by $v_j$ the direction of $\gamma_j$ (or $K$) at $x$,
again going away from $x$; we claim that
\begin{equation} \label{5.23}
|v_j(z) - v_j| \leq 30^{-1} \ \text{ for } z\in \cL_j.
\end{equation}
Indeed $v_j(z)$ lies in the intersection of the hyperplane $H_z$ tangent to $z$
at $\S$, and the tangent plane $P_z$ to $E$ at $z$. If $y\in B(x,20\tau)$ is the point
of $Y(x)$ such that $\Phi(y) = y$, $P_z$ is the image by $D\Phi(y)$ of the plane $P$ that
contains $F_j$. But $F_j$ is orthogonal to $H_x$ at $x$, and contains the tangent vector $v_j$,
and \eqref{5.8} says that $|D\Phi(y) - Id| \leq 10^{-2}$. Thus $P_z$ is quite close to $P$,
$P_z \cap H_z$ is quite close to $P \cap H_x$, and \eqref{5.23} follows.

Notice that, by Remark \ref{r6.3n} (and if $\tau$ is small enough compared to $a_1$), 
we can also make sure that
\begin{equation} \label{6.26n}
|v_j(z) - v_j| \leq a_1 \ \text{ for } z\in \cL_j,
\end{equation}
where $a_1>0$ is any small number given in advance. That is, $\cL_j$ is a Lipschitz graph, with
a Lipschitz constant that is as small as we want.

Of course it follows from \eqref{5.23} that
\begin{equation} \label{5.24}
\text{the three $\cL_j$ make angles of at least $100^\circ$ at  $x^\ast$,}
\end{equation}
because the $v_j$ make $120^\circ$ angles. 
We also want to show that for $1 \leq j \leq 3$,
\begin{equation} \label{5.25}
d_{x,8\tau}(\cL_j,\gamma_j) \leq 10^{-10} \tau.
\end{equation}
First let $z\in \gamma_j \cap B(x,8\tau)$ be given; by \eqref{5.9} we can find
$y\in Y(x) \cap B(x,20\tau)$ such that $\Phi(y)=z$. By \eqref{5.7}, $|y-z| \leq 10^{-10}\tau$.
Then $\dist(y,\S) \leq 10^{-10}\tau$, and $y' = y/|y|$ lies close to $y$ and $z$. Also,
$y$ lies on the face $F_j$ (because $z\in F'_j$), and $y'\in F_j$ too. 
But $F_j \cap \S \cap B(x,9\tau) \subset \gamma_j$ (see \eqref{5.11} and 
the definition of $\gamma_j$), so $\dist(x,\gamma_j) \leq 2 \cdot 10^{-10}\tau$.

Conversely let $y\in \gamma_j \cap B(x,8\tau)$ be given. 
Consider the radial line segment $I$ centered at $y$ and with length $3 \cdot 10^{-10}\tau$.
By \eqref{5.7} its image by $\Phi$ crosses $\S$ (one extremity in $B(0,1)$, the other one outside),
so there is a point $y'\in I$ such that $\Phi(y') \in \S$. The point $z' = \Phi(y')$ lies in
$F'_j$ (because $y'\in F_j$), and of course $|z'-z| \leq 3 \cdot 10^{-10}\tau$.
Hence $z' \in \cL_i = \S \cap F'_j \cap B(x,9\tau)$ and \eqref{5.25} follows.

It follows directly from \eqref{5.7} and the fact that $x^\ast = \Phi(y)$ for some 
$y \in F_1 \cup F_2 \cup F_3$ that
\begin{equation} \label{5.26}
|x^\ast - x| \leq 2 \cdot 10^{-10} \tau.
\end{equation}
Finally let us record the fact that \eqref{5.25} also implies that
\begin{equation} \label{5.27}
d_{x,8\tau}(K, \cL_1 \cup \cL_2 \cup \cL_3) = 
d_{x,8\tau}(\gamma, \cL_1 \cup \cL_2 \cup \cL_3) \leq 10^{-10} \tau,
\end{equation}
where the first part comes from \eqref{5.6} (recall that 
$\gamma = \gamma_1\cup\gamma_2\cup\gamma_3$).
This completes our rather precise description of $E$ and $E\cap \S$
near the real vertices $x\in V_1$.

\subsection{A description of $E \cap \S$ near a flat point $x\in E \cap \S$ far from $L \cup V_1$}

Next we want to do something similar near the regular points of $K$, i.e., the points
near which $K$ is a geodesic of $\S$. This includes the vertices of $V_2$, since
we only cut $K$ artificially there. We take a such a point $x$ and also assume that
$x$ lies far enough from $L$ or $V_1$. More precisely, we take $x$ such that
\begin{equation} \label{5.28}
x \in K \ \text{ and }\   \dist(x, L \cup V_1) \geq \tau.
\end{equation}
The constraint on the distance to $V_1$ will not cost us much, since we already have a 
good description of $E$ near $V_1$. What will be missing is what happens near $L$, 
but of course this is the main point of the paper.

We fix $x\in K \sm L$ such that \eqref{5.28} holds, and do exactly as in the previous
subsection. Near $x$, $K$ coincides with an arc of geodesic $\gamma$, which we may as well
choose maximal. It could be that $\gamma$ is composed of two or more successive arcs
$\cC_i$, $i\in \cI$, because we cut some arcs artificially with vertices of $V_2$, one of
which may even be close to $x$. 

Denote by $P(x)$ the plane that contains $\gamma$. Let us check that
if we take $\tau \leq \eta(X)$,
\begin{equation} \label{5.29}
K \cap B(x,\tau) = \gamma \cap B(x,\tau) = P(x) \cap \S \cap B(x,\tau)
\end{equation}
(as in \eqref{5.6}). For the first part, we just need to check that every 
arc $\cC_i$ that meet $B(x,\tau)$ is contained in $\gamma$. Let $\cC_j$
be the arc of the standard decomposition that contains $x$. 
If $\cC_i \cap \cC_j = \emptyset$, \eqref{5.4} says that 
$\dist(\cC_i, \cC_j) \geq 9 \eta(X) \geq 9\tau$; this is impossible 
because $\dist(x,\cC_i) \leq \tau$. So $\cC_j$ meets $\cC_i$, and
since there is no true vertex of $K$ near $x$, this means that $\cC_i$
is part of $\gamma$ too; the first part of \eqref{5.29} follows.
For the second part we just need to observe that none of the two branches
of $\gamma$ (when we leave from $x$) stops before we reach a point of $V_1 \cup L$.
By \eqref{5.28}, this does not happen as long as we stay in $B(x,\tau)$.

From \eqref{5.28} we also deduce that
\begin{equation} \label{5.30}
\text{$E$ is a plain almost minimal set in $B(x,\tau)$, with gauge function $h(s) = C_h s^\beta$.}
\end{equation}
The radius is somewhat smaller than in \eqref{5.10}, so this will force us to apply
the lemmas from \cite{Holder} and \cite{C1} with slightly smaller radii, but otherwise
things will be as easy as before.
The analogue of Lemma \ref{t5.2} for this case is the following.

\ms
\begin{lem} \label{t5.3}
If $\varepsilon$ is small enough, depending on $n$ and $\tau$, then for $x$ 
as in \eqref{5.28} there is a $C^1$ diffeomorphism 
$\Phi : B(x,2 \cdot 10^{-3}\tau) \to \Phi(B(x,2\cdot10^{-3}\tau)$ 
with the following properties: 
\begin{equation} \label{5.31}
|\Phi(y) - y| \leq 10^{-10} \tau \ \text{ for } y\in B(x,2 \cdot 10^{-3}\tau);
\end{equation}
\begin{equation} \label{5.32}
|D\Phi(y) - Id| \leq 10^{-2} \ \text{ for } y\in B(x, 2 \cdot 10^{-3}\tau);
\end{equation}
\begin{equation} \label{5.33}
E \cap B(x,10^{-3}\tau) = \Phi(P(x)\cup B(x, 2 \cdot 10^{-3}\tau)) \cap B(x,10^{-3}\tau).
\end{equation}
\end{lem}

\ms
We skip the proof, which is just the same as for Lemma \ref{t5.2}.
Again notice that there is nothing special about points $x\in V_2$, we do not
have singularities of $K$ near these points, they were just added to simplify some 
estimates in later sections.

As in Remark \ref{r6.3n}, we can even replace $10^{-10}$ in \eqref{5.31} and $10^{-2}$ in \eqref{5.32}
with any small constant $a_0$ decided in advance, but then $\varepsilon$ has to depend on $a_0$ too.

By the same discussion as for $x\in V_1$ (but simpler because we have no branching), 
$E \cap B(x,10^{-3}\tau)$ is also a small Lipschitz graph over $P(x)$.
Then the implicit function theorem allows us to say (as in \eqref{5.21}) that
\begin{equation} \label{5.34}
E \cap \S \cap B(x,9 \cdot 10^{-4} \tau) \ \text{ is a $C^{1+\beta_1}$ curve $\cL_x$.} 
\end{equation}
Moreover $\cL_x$ is a small Lipschitz graph, in the sense that if $v_x$ denotes
a unit tangent vector to $K$ at $x$, and if similarly $v(z)$ denotes a continuous choice
of unit tangent vector to $\cL_x$ at $z\in \cL_x$, then the proof of \eqref{5.23} also yields
\begin{equation} \label{5.35}
|v(z) - v_x| \leq 30^{-1} \ \text{ for } z\in \cL_x
\ \text{ or } |v(z) + v_x| \leq 30^{-1} \ \text{ for } z\in \cL_x
\end{equation}
(if we choose the opposite orientations by mistake). And as before, by choosing $\varepsilon$ even
smaller, we can even arrange that
\begin{equation} \label{6.39n}
|v(z) - v_x| \leq a_1 \ \text{ for } z\in \cL_x
\ \text{ or } |v(z) + v_x| \leq a_1 \ \text{ for } z\in \cL_x
\end{equation}
for any given small constant $a_1 > 0$.
The analogue of \eqref{5.25}, namely the fact that
\begin{equation} \label{5.36}
d_{x,8 \cdot 10^{-4}\tau}(\cL_x,\gamma) 
= d_{x,8 \cdot 10^{-4}\tau}(\cL_x,K) 
\leq 10^{-6},
\end{equation}
is proved the same way (the easy first part comes from \eqref{5.29}).

In the special case when $x\in V_2$, we will need to define a vertex $x^\ast$
where we cut $\cL_x$ in two. We simply choose $x^\ast \in \cL_x$ so that
$|x^\ast-x|$ is minimal, and by the proof of \eqref{5.25}, we get that 
\begin{equation} \label{6.41n}
|x^\ast-x| \leq 10^{-10} \tau.
\end{equation}

\subsection{The desired description of $E\cap \S \sm [D_+(\tau) \cup D_-(\tau)]$}

We now have a nice description of $E \cap \S$ near all the points of $K \sm L$
which lie far from $L$, which we put together to get a relatively simple statement.
Recall the definition of $D_\pm(\tau)$ (two small spherical balls centered at the points $\ell_\pm$
of $L\cap \S$).

\begin{pro}\label{t5.4} 
For each choice of $\tau  \leq 10^{-3}\eta(X)$, there exists $\varepsilon > 0$ such that
if $X$ and $E$ satisfy the assumptions of Section \ref{S4} with $r=1$, then
we can find $C^1$ curves $\cL_i \subset \S$, $i\in \cI$, such that
\begin{equation}\label{5.37}
E \cap \S \sm (D_+(\tau) \cup D_-(\tau)) = \bigcup_{i\in \cI} \cL_i \, ,
\end{equation}
the curves $\cL_i$, $i\in I$, are disjoint, except perhaps for their endpoints,
and they are related to the $\cC_i$ in the following way.
For each vertex $x\in V_1 \cup V_2$, we can find a point $x^\ast \in E \cap \S$,
such that 
\begin{equation} \label{5.38}
|x^\ast - x| \leq 10^{-9} \tau
\end{equation}
and, for $i \in \cI_1$, $\cL_i$ is a simple $C^1$ curve in $\S$, with
endpoints $a_i^\ast$ and $b_i^\ast$ (the two points $x^\ast$ associated to 
$x=a_i$ and $x=b_i$ respectively), such that
\begin{equation} \label{5.39}
\dist(z,\cC_i) \leq 10^{-8} \tau \text{ for } z\in \cL_i \ \text{ and } \ 
\dist(z,\cL_i) \leq 10^{-8} \tau \text{ for } z\in \cC_i.
\end{equation}
For $i \in \cI_0$, $\cC_i$ meets $\d D_+(\tau) \cup \d D_-(\tau)$
at a single point $c_i$, $\cL_i$ ends at a point $c_i^\ast$ such that
\begin{equation} \label{5.40}
|c_i^\ast - c_i| \leq 10^{-8} \tau
\end{equation}
and, if $a_i$ denotes the endpoint of $\cC_i$ that does not lie on $L$,
$\cL_i$ goes from $a_i^\ast$ to $c_i^\ast$, and 
\begin{equation} \label{5.41}
\dist(z,\rho(a_i,c_i)) \leq 10^{-8} \tau \text{ for } z\in \cL_i \ \text{ and } \ 
\dist(z,\cL_i) \leq 10^{-8} \tau \text{ for } z\in \rho(a_i,c_i).
\end{equation}
\end{pro}

\ms
We already have most of the needed information, but need to make a few remarks
to put things together. Also, the statement above misses some information that we
obtained in the last subsections; we will refer to them concerning $E$ itself and 
Lipschitz graph properties, for instance.

Before we put our local arcs together, let us say how we intend to end
our curves near the two points of $\S \cap L$.

Let $\ell \in L\cap \S$ be given, and first suppose that $\ell \notin K$. 
Then $\dist(\ell, K) \geq \eta_L(X) \geq 10\eta(X)$
by \eqref{3.10} and \eqref{3.12}, and then
\begin{equation} \label{5.42}
B(\ell,6\eta(X)) \cap (X \cup E) = \emptyset,
\end{equation}
by \eqref{4.3} (and if $\varepsilon$ is small enough). In this case we'll not need
to do anything to cut the arcs of $E \cap \S$ near $\ell$.

Next assume that $\ell \in K$. Then the set $\cI_0(\ell)$ of indices $i$ such that
$\cC_i$ (starts or) ends at $\ell$ is not empty. The arcs $\cC_i$ are not too short,
because \eqref{3.11} and \eqref{3.12} say that $\dist(\ell,V_1 \cup V_2) \geq 10 \eta(X)$.
Hence each of them cuts each sphere $\d B(\ell, t)$, $t \leq 5 \eta(X)$, at a point $x_i(t)$.
Call $m \in \{ 1, 2, 3\}$ the number of elements in $\cI_0(\ell)$.

We shall restrict to $t\in [\tau, 2\tau]$ for some $\tau < 10^{-3} \eta(X)$. We claim that
if $\varepsilon > 0$ is small enough (depending on $\tau$ and $n$), then for each choice 
of  $t \in (\tau,3\eta(X))$ and $i\in \cI(\ell)$, 
\begin{equation} \label{5.43}
\begin{aligned}
\d B(\ell, t) &\ \text{ intersects $E \cap \S$ exactly $m$ times, transversally,}
\cr&\hskip1cm 
\text{ and at points $x_i^\ast(t)$, $1 \leq i \leq m$, such that } 
| x_i^\ast(t) - x_i(t) | \leq 10^{-8} \tau.
\end{aligned}
\end{equation}
This is in fact easy. Lemma \ref{t5.3}, applied to $x=x_i(t) \in K \cap \d B(\ell,t)$,
shows that $E \cap \S \cap B(x,\tau)$ (or equivalently, $\cL_x$) is a small Lipschitz 
curve with (by \eqref{5.15}) a tangent direction which is almost the same as the
direction of $K$ at $x$. This curve crosses $\d B(\ell, t)$ transversally (in fact, almost
perpendicularly). Thus, near the $x_i(t)$, we get a unique point 
$x_i^\ast(t) \in E \cap \S \cap \d B(\ell, t)$, and $| x_i^\ast(t) - x_i(t) | \leq 10^{-8} \tau$.
But there is no other point, because all points of $E \cap B(0,2) \cap \d B(\ell, t)$ 
lie close to $X$ (by \eqref{4.3}), hence close to one of the $x_i(t)$.

We are now ready to say how we organize the local description of $E \cap \S \sm B$
to make curves $\cL_i$, $i\in \cI$. Fix $\tau \in (0,10^{-3}\eta(X)]$.
Also set $B = B(\ell_+,\tau) \cup B(\ell_-,\tau)$; we have a nice local description,
in balls centered on $K$, of $\Gamma = E \cap \S \sm B$, and by \eqref{4.3} the balls with
the same centers and half the radii cover $\Gamma$. We cut $\Gamma$ at the points
$x^\ast$, $x\in V_1 \cup V_2$, and we get a collection of connected components $\cL$
which we can also describe: the two endpoints of $\cL$ are points $x^\ast$,
we can use \eqref{5.25} or \eqref{5.36} to follow $K$ along $\cL$, and find out
that $\cL$ stays close to some $\cC_i$, $i\in \cI$. In addition,
if  $\cC_i = \rho(a_i, b_i)$, $\cL$ has the endpoints $a_i^\ast$ and $b_i^\ast$.
Or, when $i\in \cI_0$, one of the endpoint of $\cL$ is the point $x_i^\ast(\tau)$ associated
to $t=\tau$ as in \eqref{5.43}.
This completes our proof of Proposition \ref{t5.4}.
\qed

\ms
Notice that by \eqref{5.43}, the spheres $\d B(\ell, t)$, $t\in [\tau,2\tau]$, also cross the $\cL_i$,
$i \in \cI_0(\ell)$ once and transversally, at least if $\ell \in V_0 = K \cap L$.
Otherwise, they don't touch, by \eqref{5.42}. 

Notice also that the $\cL_i$ also satisfy a version of Lemma \ref{t5.1},
since they are very close in distance to the $\cC_i$.

Finally observe that we may apply Remark \ref{r6.3n} and the ensuing comments and get that 
each $\cL_i$ is a Lipschitz graph over the geodesic $\cC_i$ (maybe made a little longer or shorter
to accommodate the endpoints), with a Lipchitz constant that can be taken as small as we want. 
The only price to pay is that we have to choose smaller constants $\tau$ and $\varepsilon$. 
We will have the option to use this to simplify the construction of our 
competitors, at the beginning of Section \ref{S12}. 

\section{Connectedness configurations near $\ell \in L \cap \S$, and a first net of curves} 
\label{S6} 

In this section and the next ones we fix a point $\ell \in L \cap \S$ and restrict our 
attention to the small spherical disk
\begin{equation} \label{6.1}
D = \S \cap \overline B(\ell,\tau),
\end{equation}
where $\tau$ is a small constant, to be chosen later. We shall assume that
$\tau \leq 10^{-3} \eta(X)$, so as to be able to apply the results of Section \ref{S5}.

We shall distinguish between a few different cases, depending on the number of
points of $E \cap \d D$ and the way $E \cap D$ connects them to each other
and to $\ell$, and then we shall construct a first net of simple curves 
$\gamma \subset E \cap D$, with the same basic connecting properties.
We will do this independently for each of the two points of $L \cap \S$.

The first number that we care about is the number $m = m(\ell)$ of arcs of
$K$ that start from $\ell$. Thus $m$ is the cardinality of the set $\cI_0(\ell)$
of indices $i$ such that $\cC_i$ starts from (or ends at) $\ell$. 
When $i \in \cI_0(\ell)$, we shall systematically denote by $a_i$ the other endpoint of $\cC_i$.
As a general rule, we shall say that we are in {\bf \underbar{Configuration m}} 
when $\cI_0(\ell)$ has $m$ elements. We know from \eqref{2.7}
that $0 \leq m \leq 3$. But there will be numerous subcases, 
depending on connectivity properties.

We start with the easy {\bf \underbar{Configuration 0}} when $m=0$.
In this case, $\ell$ is connected to no vertex of $V_1 \cup V_2$,
\eqref{3.10} and \eqref{3.12} say that $\dist(\ell,K) \geq \eta_L(X) \geq 10 \eta(X)
\geq 10^4 \tau$, and then
\begin{equation} \label{6.2}
X \cap B(\ell, 2\tau) = \emptyset \ \text{ and } \ 
E \cap B(\ell, 2\tau) = \emptyset
\end{equation}
(by \eqref{4.3}). In this case we shall do nothing in $D$ in this section or the next ones.

\ms
In the other cases, let us renumber the curves $\cC_i$ that touch $0$
(or equivalently, the set $\cI_0(\ell)$), and just call them $\cC_i$, $1 \leq i \leq m$.

For $1 \leq i \leq m$, Proposition \ref{t5.4} gives a curve $\cL_i \subset E$, that leaves from
a point $a_i^\ast$ that lies very close to $a_i$ (the other endpoint of $\cC_i$), and
ends at a point of $\d D = \S \cap \d B(\ell,\tau)$ that we called $c_i^\ast$.
Recall from \eqref{5.40} that $c_i^\ast$ lies quite close to $c_i$, the point of 
$\cC_i \cap \d D$. Also, by \eqref{5.43},
\begin{equation} \label{6.3}
E \cap \d D = \big\{ c_i^\ast \, ; \, 1 \leq i \leq m \big\}.
\end{equation}
For each $x\in E \cap D$, we denote by $H(x)$ the connected component of 
$x$ in $E \cap D$. Also set $H_i = H(c_i^\ast)$ for $1 \leq i \leq m$ and
$H_\ell = H(\ell)$. 

One case that we particularly like is the case of a \underbar{hanging curve}.
We define it as the case when for some $i \leq m$, $c_i^\ast$ is not connected
to any of the other special points, i.e., when
\begin{equation} \label{6.4}
c_i^\ast \notin H_\ell \cup \bigcup_{j \neq i} H_j.
\end{equation}
We refer to this as {\bf \underbar{Configuration H}}. 
In this case we are happy, because we will be able to contract a large piece of $E$
along the piece of $E$ near $\cL_i$, and this will give estimates that are quite favorable.
Nonetheless we define a set $\gamma$ by (short) induction.

We define a first set $\gamma_i = \{ c_i^\ast \}$, remove the point $c_i^\ast$
(or equivalently the whole $H_i$) from the game, and get a new configuration with $m-1$
points. If $m=1$, just set $\gamma = \gamma_i$. Otherwise, define the net $\gamma'$ 
associated to this smaller configuration (as will be explained in the next cases), 
and take $\gamma = \gamma_i \cup \gamma'$. As we shall see, $\gamma'$ is
contained in the union of the $H_j$, $j\neq i$, so we get a disjoint union. For instance,
if we had three hanging curves, $\gamma$ would be composed of the three $c_i^\ast$.

The next simple case is {\bf \underbar{Configuration 1}}, where $m=1$ and $H_1 = H_\ell$.
In this case we choose $\gamma$ so that
\begin{equation} \label{6.5}
\gamma \text{ is a simple curve in $E \cap D$ that goes from $c_1^\ast$ to $\ell$}.
\end{equation}
The existence of such a curve follows rather easily from the fact that $H_1$ is connected
and contains $c_1^\ast$ and $\ell$, plus \eqref{4.4} which says that 
\begin{equation} \label{6.6}
\H^1(H_1) \leq \H^1(E \cap \S) < +\infty;
\end{equation}
see for instance \cite{Falconer}, or Chapter 30 of \cite{MSBook}.
Later on we shall replace this curve with a small Lipschitz curve, and then we shall
retract onto it, but for the moment we continue with our list of cases.
Notice that only Configurations H and 1 are possible when $m=1$.

Now assume that $m=2$, and that we are not in Configuration H. A first option is that
$H_1 = H_2 \neq H_\ell$, i.e., $c_1^\ast$ and $c_2^\ast$ are connected to each
other (in $E \cap D$), but not to $\ell$. We shall call this {\bf \underbar{Configuration 2-}}.
By \eqref{6.6} and the same argument as above, we can find $\gamma$ such that
\begin{equation} \label{6.7}
\gamma \text{ is a simple curve in $E \cap D$ that goes from $c_1^\ast$ to $c_2^\ast$}.
\end{equation}
Notice that $\ell \notin \gamma$, because $\gamma \subset H_1$.

The next case, called {\bf \underbar{Configuration 2+}}, is when $H_1 = H_2 = H_\ell$,
i.e., our three special points are connected. First select (as above) a simple curve 
$\gamma_{1,2} \subset E \cap D$, that goes from $c_1^\ast$ to $c_2^\ast$.
Also choose a simple curve $\gamma_0 \subset \d D$, that goes from $\ell$
to $c_1^\ast$. Finally denote by $\gamma_\ell$ the part of $\gamma_0$
that lies between $\ell$ and the first point of $\gamma_{1,2}$ that we hit when
we start from $\ell$ and run in the direction of $c_1^\ast$. We include both endpoints.
It could happen that $\ell$ already lies on $\gamma_{1,2}$, and then $\gamma_\ell$
is reduced to the point $\ell$; this is all right too. We set
\begin{equation} \label{6.8}
\gamma = \gamma_{1,2} \cup \gamma_\ell.
\end{equation}
This was our last case when $m=2$, since each $c_i^\ast$ is connected to another
special point when we are not in case $H$. Recall that in case $H$, either the two curves $\cL_i$
are hanging, and then we set $\gamma = \{ c_1^\ast, c_2^\ast \}$ or $\cL_2$, say, is hanging
and $H_1 = H_\ell$, and then we select $\gamma_1 \subset H_1$ as in Configuration 1, and set
$\gamma = \gamma_1 \cup \{ c_2^\ast \}$.

Now suppose that $m = 3$ and we are not in Configuration H. A first possibility, that we shall call
{\bf \underbar{Configuration 3 = 2+1}}, is that two of the $c_i^\ast$ are connected to each
other, but not to $\ell$, and the third one is connected to $\ell$ (hence not to the others).
Let us relabel the $c_i^\ast$, if needed, so that in fact
\begin{equation} \label{6.9}
H_1 = H_2  \text{ and  $H_3 = H_\ell$, but } H_1 \neq H_\ell.
\end{equation}
In this case we take
\begin{equation} \label{6.10}
\gamma = \gamma_{1,2} \cup \gamma_\ell,
\end{equation}
where $\gamma_{1,2}$ is a simple curve in $E \cap D$ that goes from $c_1^\ast$
to $c_2^\ast$, and $\gamma_\ell$ is a simple curve in $E \cap D$ that goes from $\ell$
to $c_3^\ast$. Notice that these two curves are disjoint, since $H_1 \neq H_3$.

Next we turn to {\bf \underbar{Configuration 3-}}, where $H_1 = H_2 = H_3 \neq H_\ell$.
In this case we select a simple arc $\gamma_{1,2}$ in $E \cap D$ that goes from $c_1^\ast$
to $c_2^\ast$ and a simple arc $\gamma_{3,1}$ in $E \cap D$ that goes from 
$c_3^\ast$ to $c_1^\ast$. Then we let $\gamma_3$ be, as in Configuration 2+, 
the closed sub-arc of $\gamma_{3,1}$ between $c_3^\ast$ and the first time 
we hit $\gamma_{1,2}$. We set
\begin{equation} \label{6.11}
\gamma = \gamma_{1,2} \cup \gamma_3;
\end{equation}
thus $\gamma$ is a (possibly degenerate) three-legged spider that connects the
$c_i^\ast$; it does not contain $\ell$. 

We are left with only one case, which we call {\bf \underbar{Configuration 3+}}, where
all the points are connected, i.e., $H_1 = H_2 = H_3 = H_\ell$. We define $\gamma_{1,2}$
as before, then $\gamma_{3,1}$ and $\gamma_3 \subset \gamma_{3,1}$,
but now also pick a simple arc $\gamma_{\ell,1}$ in $E \cap D$ that goes from
$\ell$ to $c_1^\ast$, but only keep the arc $\gamma_\ell$ that goes from $\ell$
to the first time we hit $\gamma_{1,2} \cup \gamma_3$. Finally we take
\begin{equation} \label{6.12}
\gamma = \gamma_{1,2} \cup \gamma_3 \cup \gamma_\ell
\end{equation}
(a three-legged spider with a short leash to $\ell$, which again can also be degenerate in 
different ways). This is our most complicated case; it will turn out that, later in the proof, 
we shall replace $\gamma$ by a connected set with a slightly simpler shape, but for the moment 
$\gamma$ is good enough.

\ms
At this stage, we constructed a net $\gamma \subset E \cap D$, with the same
connecting properties (regarding our special points) as $E \cap D$.
We will have to modify this net, though, for a few concurrent % checked
reasons. 

The first one is that we may find a significantly shorter net in $D$, still with the same
connecting properties, but which may not be entirely contained in $E$. 
The interest is that the cone over this new net will have a smaller surface, 
and then we should be able to use this new cone to define an interesting competitor. 
This looks like a good idea, but we need to show some restraint, and only do this 
transform when we win enough length (and then surface for the cone) to compensate 
for the extra cost that we will need to pay when we glue together pieces of different surfaces.
The ideal thing would be to find a shorter net in $E \cap \S$, but this is probably not
going to be possible in general, so we'll have to add small connecting pieces.

The second reason is that our construction of competitors (in particular related to the
new curves) will use Lipschitz retractions from a neighborhood of the net to the net.
They will be easier to find if the net has some regularity: we do not want to have to
retract on a curve that almost makes a closed loop. Thus we will like better nets that
are composed of a small number of small Lipschitz curves that meet with large angles
and make no loops.

Fortunately, the two reasons go together: if $\gamma$ makes unnecessarily long 
connections (think about a long arc that connects two close points), we may replace
some parts of it with shortcuts, save a nontrivial amount of length, and at the same
time increase the chances for a nice retraction.

There is a third reason for which we like small Lipschitz graphs.
One of the main engines of our proof is the comparison between cones and
graphs of harmonic functions. Suppose for instance that $C$ is the unit circle
in a plane $P$, and that $\Gamma$ is the graph over $C$ of a Lipschitz function
$A : C \to P^\perp$. Let $A_1$ be the homogeneous extension of $A$ (of degree $1$),
so that the graph of $A_1$ (call it $\wh \Gamma$) is also the cone over $\Gamma$. 
Then let $A_h$ denote the harmonic extension of $A$ to $P \cap \overline B(0,1)$,
and denote by $\Sigma$ its graph.
We expect that $\H^2(\Sigma) < \H^2(\wh\Gamma \cap \overline B(0,1))$, especially
if $A$ is far from an affine function, but this is much more pleasant to prove when 
$A$ has a small Lipschitz constant, because then we can use expansions of order $2$
to estimate the surface measure. (If $\Gamma$ is not even a graph, it is much harder 
to imagine an analogue of the harmonic graph, and we are a little desperate.)
This comparison argument will also work with the circle replaced by small sectors
that we glue to each other, and this is the reason why we wish to replace our $\gamma$ with
shorter nets of small Lipschitz curves.

\section{The standard replacement by a Lipschitz graph}
\label{S7}

As we just explained, we intend to modify the simple nets $\gamma$ of the previous section,
and our main tool for this is a construction, which we will take from \cite{C1}
but which was not especially original anyway, that takes a simple curve $\gamma$
(typically, one piece of the nets above) and creates a Lipschitz curve $\Gamma$
with the same endpoints.

The constraints of the game are that we don't want $\Gamma$ to be longer than 
$\gamma$, and we only want to introduce parts of $\Gamma \sm \gamma$
when we are sure that this will make the curve significantly shorter, or that we will
win something proportional in the next section. We don't need to be specific yet, let us
just remember that we do not want to change the curve for no reason.

We do the construction in this section, and apply it later. 
So let $\gamma$ be a simple curve. To complicate matters (and in particular parameterizations), 
$\gamma \subset \S$ and we want $\Gamma \subset \S$ too. But we will not be disturbed 
by antipodal problems, because we will very fast be able to assume that 
$\H^1(\gamma) \leq \frac{4 \pi}{5}$.

Let $a$ and $b \neq a$ denote the endpoints of $\gamma$. We assume that
\begin{equation}\label{7.1}
\ddist(a,b) := \H^1(\rho(a,b)) < \frac{3 \pi}{4},
\end{equation}
where $\ddist$ denotes the geodesic distance in the sphere, and
$\rho = \rho(a,b)$ is the geodesic with endpoints $a$ and $b$. 

So we want to construct a new curve $\Gamma$ in $\S$, which is a Lipschitz graph 
over $\rho$ (we shall explain what this means soon), has the same endpoints $a$ and $b$, 
and coincides with $g$ on a set which is as large as possible.

The main parameter in this construction is a small constant $\lambda < 1$,
which is essentially the Lipschitz constant for desired graph $\Gamma$. 
It used to be called $\eta$ (in \cite{C1}), but we want to avoid notation conflicts.
It is required to be small, depending on the dimension, and the main reason for this
is that we want to be able to apply the results of Section \ref{S8}. Recall also that
we want to choose $\tau$ small, depending on $\lambda$, so $\lambda$ should not depend
on $\tau$.

Before we start for real, let us eliminate a simple case. 
Let $\tau_1 > 0$ be small (to be chosen later, depending on $\lambda$). 
We take  
\begin{equation}\label{7.2}
\Gamma = \rho(a,b) \ \text{ when } \ 
\H^1(\gamma) \geq (1+\tau_1) \ddist(a,b),
\end{equation}
and feel  happy because although we added a big set, we also saved at least
$\tau_1 \ddist(a,b)$ in length. This works even if $\H^1(\gamma) = +\infty$,
but we do not need the information.

From now on we assume that the condition in \eqref{7.2} fails, i.e., we assume that
\begin{equation}\label{7.3}
\length(\gamma) = \H^1(\gamma) < (1+\tau_1)\, \ddist(a,b),
\end{equation}
where the first equality is just a change of notation because $\gamma$ is simple
(see for instance \cite{Falconer}, but the truth is that we could use $\H^1$ all along). 

Set $\rho = \rho(a,b)$ and let $P$ denote the $2$-plane that contains $\rho$ 
(or equivalently $a$, $b$, and $0$). We shall often use the fact that because of \eqref{7.3},
\begin{equation}\label{7.4}
\dist(z,\rho) \leq \tau_2 \, \ddist(a,b) \ \text{ for } z\in \gamma,
\end{equation}
where $\tau_2 = 10\sqrt{\tau_1} > \tau_1$ (because $\tau_1 < 1$).

It should be said now that (contrary to what we may have implied so far)
we do not always try to make $\H^1(\Gamma)$ significantly shorter than $\gamma$, 
but sometimes we want to control $\H^1(\Gamma \sm \gamma)$ in terms of something else, 
the $L^2$-norm of the derivative of a function whose graph describes $\Gamma$.
This is because we shall see in the next section that a harmonic replacement of the 
cone over $\Gamma$ will make us save a comparable amount of area.

The argument will be essentially imported from Section 7 of \cite{C1}, and we use similar notation
(except that $\eta_1$ is now called $\lambda$). 
There is a small difference with what was done in \cite{C1}, where for convenience we assumed 
(in (7.1) there) that $\length(\gamma) \geq 9 \eta_0$ for some small geometric constant $\eta_0$. 
Here we do not want to assume this, and this will force us to be some times a little more careful 
with the normalization; for instance, in \cite{C1} which just required that 
$\length(\gamma) \leq  \ddist(a,b) + \tau_1$.
Here we assume the stronger \eqref{7.3}, with an error term of at most $\tau_1 \ddist(a,b)$.
The necessary modifications, to adapt the construction of \cite{C1}, will all be of that type.

With this in mind, the main assumptions (7.1)-(7.3) in Section 7 of \cite{C1} are satisfied, by
\eqref{7.1}, \eqref{7.3} and \eqref{7.4}. An important quantity, that we want to use to estimate 
various terms, is the length excess
\begin{equation}\label{7.5}
\Delta L = \length(\gamma) - \ddist(a,b).
\end{equation}

We now describe the highlights of the construction of \cite{C1}.
We let $z : I \to \S$ denote a parameterization of $\gamma$ by arclength, 
so that $|I| = \length(\gamma)$, and we write $z(t) = (z_1(t), z_2(t), v(t))$, 
where the first two coordinates are in $P$, and $v(t) \in P^\perp$.
A simple estimate with Fourier series shows that
\begin{equation}\label{7.6}
\int_I |v'(t)|^2 dt \leq 14 \Delta L;
\end{equation}
see Lemma 7.8 in \cite{C1}; the reader should not worry about normalization here, 
as even in \cite{C1} the constant 14 does not depend on $\eta_0$.

Next write $(z_1(t), z_2(t)) = (w(t) \cos\theta(t), w(t) \sin\theta(t)$,
as in (7.5) of \cite{C1} (notice that by \eqref{7.4}, $(z_1(t), z_2(t))$ stays close to $\rho$);
we also need to know that (if $a$ and $b$ are chosen in trigonometric order) 
$\theta'(t)$ is rather large on average. To measure this, we define $f$ on $I$ by
\begin{equation}\label{7.7}
f(t) = 1+2|v(t)|^2 - \theta'(t)
\end{equation}
(as in (7.20) there), observe that $f(t) \geq 0$ almost everywhere 
(see the line below (7.20) there, which uses the fact that $|z'(t)| = 1$ almost everywhere), 
and use \eqref{7.6} and \eqref{7.5} to show that
\begin{equation}\label{7.8}
\int_I f(t) dt \leq 30 \Delta L;
\end{equation}
see (7.21) in \cite{C1}. We now use a maximal function argument, 
based on the two estimates \eqref{7.6} and \eqref{7.8}, to find an open set $Z$ in $I$ 
(in fact, $Z$ is the set where one of the two maximal functions $f^\ast(t)$ or $(v')^\ast$ is large) 
with the following two properties. First,
\begin{equation}\label{7.9}
|Z| \leq C \lambda^{-2} \Delta L = C \lambda^{-2} (\length(\gamma) - \ddist(a,b))
\leq C \lambda^{-2} \tau_1  \ddist(a,b)
\end{equation}
as in (7.26) of \cite{C1} (and by \eqref{7.5} and \eqref{7.3}).
Here $C$ is an absolute constant that comes from the Hardy-Littlewood maximal 
theorem on $I$. But also, $z$ has good Lipschitz properties away from $Z$, 
that we shall explain soon.

Write $Z$ as a countable disjoint union of open intervals $I_j = (a_j,b_j)$, with 
possibly two exceptions: if the initial endpoint of $I$ lies in $Z$, 
the corresponding interval is of the form $[a_j,b_j)$, and if the final endpoint 
of $I$ lies in $Z$, then $I_j = (a_j,b_j]$.
Both things do not happen at the same time, for the following reason:
we shall choose $\tau_1$ small, depending on $\lambda$, and in 
particular, we can make sure that $C \lambda^{-2} \tau_1 < 1$ in \eqref{7.9},
so that $|Z| < \ddist(a,b) \leq |I|$.

We come to the good Lipschitz properties. The definition of $Z$ in terms of 
maximal functions yields (see (7.33) in \cite{C1})
\begin{equation}\label{7.10}
|v(b_j)-v(a_j)| \leq \frac{\lambda (b_j-a_j)}{4}
\ \text{ and }\ 
\theta(b_j)-\theta(a_j) \geq \frac{b_j-a_j}{2}.
\end{equation}

We now construct $\Gamma$. We directly define a parameterization $\wt z : I \to \S$
of $\Gamma$. On $I \sm Z$, we simply keep $\wt z(t) = z(t)$, and on each interval $I_j$, 
we let $\wt z$ be a parameterization with constant speed of the arc of geodesic 
$\rho_j = \rho(z(a_j),z(b_j))$.
Then $\wt z$ is continuous; it is even $1$-Lipschitz, because 
\begin{equation}\label{7.11}
\length(\rho_j) = \ddist(z(a_j),z(b_j)) \leq b_j-a_j
\end{equation}
because $z$  is $1$-Lipschitz. Notice also that (as in (7.30) of \cite{C1})
\begin{equation}\label{7.12}
\Gamma  \text{ has the same endpoints $a$ and $b$ as $\gamma$ and $\rho$.}
\end{equation}
Next we can use \eqref{7.10} to show that
\begin{equation}\label{7.13}
\Gamma \ \text{ is a Lipschitz graph with constant } \leq \lambda.
\end{equation}
See (7.32) in \cite{C1}, and (7.42) or (7.44) there for definitions in terms of parameterizations, 
but for here the simplest is to notice (and take as a definition of Lipschitz curve)
that we have the simpler-to-state property that
\begin{equation}\label{7.14}
\pi^\perp(\wt z(t)) \text{ is a $\lambda$-Lipschitz function of } \pi(\wt z(t)),
\end{equation}
where $\pi$ and $\pi^\perp$ denote the orthogonal projections on $P$ and its orthogonal complement;
see (7.45) in \cite{C1}. Notice that
\begin{equation}\label{7.15}
\H^1(\Gamma) = \length(\Gamma) \leq \length(\gamma) = \H^1(\gamma)
\end{equation}
because both curves are simple, and by \eqref{7.11}; hence
\begin{equation}\label{7.16}
\H^1(\Gamma \sm \gamma) \leq \H^1(\gamma \sm \Gamma) 
\leq \sum_j (b_j-a_j) \leq C \lambda^{-2} \Delta L,
\end{equation}
by \eqref{7.9}, and as in (7.31) of \cite{C1}.

\section{Harmonic graphs usually do better than cones}
\label{S8}

Let $\Gamma$ be a small Lipschitz graph over a reasonably short geodesic 
$\rho(a,b)$ (say, so that \eqref{7.1} holds), and denote by $P$ the vector plane that 
contains $\rho(a,b)$.
Our main example will be the curve that we constructed in the previous section, 
starting from $\gamma$, but we could use slightly different $\Gamma$.
In this section we use Section 8 of \cite{C1} to construct a small Lipschitz graph over 
a sector of $P$, whose area is often significantly smaller than the area of the cone 
over $\Gamma$.
The estimates below will work as soon as a bound $\lambda$ on the Lipschitz constant 
is small enough (depending on $n$ only).

When $\Gamma$ is a geodesic, the cone over $\Gamma$ is a plane sector, and we shall not 
modify anything in this section, but one could also follow the construction below and find out
at the end that we did nothing. This means that when $\Gamma$ comes from a curve $\gamma$
as in the last section, we can assume that \eqref{7.3} holds (because otherwise $\gamma$
is a geodesic). The amount of area that we are able to save will be essentially proportional
to the quantity
\begin{equation} \label{8.1}
\Delta_\Gamma = \length(\Gamma) - \ddist(a,b) 
\leq \tau_1 \ddist(a,b).
\end{equation}

First we define a homogeneous function $F$, defined on a sector $D_T$ of $P$, 
and whose graph coincides in $B(0,1)$ with the cone over $\Gamma$. 
We start with the sector $D_T$. Choose coordinates on $P$ so that $a = (1,0)$ and 
$b= (\cos T, \sin T)$, where $T = \ddist(a,b) < \frac{3\pi}{ 4}$ by \eqref{7.1}. 
Then set
\begin{equation}\label{8.2}
D_T = \big\{ (r \cos t, r \sin t) \, ; \, r\in (0,1) \text{ and } t\in (0,T) \big\}.
\end{equation}
We assume that we can parameterize $\Gamma$ in the following way:
we can find an $\lambda$-Lipschitz function $v : [0,T] \to P^\perp$, 
with $v(0) = v(T) = 0$, such that if we set
\begin{equation}\label{8.3}
w(t) = \big(1-|v(t)|^2 \big)^{1/2}
\end{equation}
and then
\begin{equation}\label{8.4}
h(t) = (w(t)\cos t , w(t)\sin t , v(t)) \in P \times P^\perp 
\end{equation}
for $t\in [0,T]$, then $h$ is a parameterization of $\Gamma$.
Thus $t\in [0,T]$ is the angle with the direction of $a$ of the orthogonal projection 
on $P$ of the running point.

In the special case when $\Gamma$ comes from $\gamma$ as in Section \ref{S7},
the existence of $v$ is checked in Remark 8.3 of \cite{C1}. Let us also check that 
we can find $v$ as above, except maybe only $2\lambda$-Lipschitz,
when $\Gamma$ is a $\lambda$-Lipschitz graph over the geodesic $\rho(a,b)$.
By definition, this last means that $\Gamma$ is a curve in $\S$, from $a$ to $b$, and that
\begin{equation} \label{8.5}
|\pi^\perp(z)-\pi^\perp(z')| \leq \lambda |\pi(z)-\pi(z')|
\end{equation}
for $z,z'\in \Gamma$. 
Here we denote by $\pi$ and $\pi^\perp$ the orthogonal projections on $P$
and $P^\perp$ respectively. We also write $\pi(z) = w e^{it}$ and 
$\pi(z') = w' e^{it'}$, with $w, w' \geq 0$. Notice that since $\pi^\perp(a) = 0$, we see that
$|\pi^\perp(z)| \leq 2\lambda$ is very small, and (since $z \in \S$)
\begin{equation} \label{8.6}
w = (1-|\pi^\perp(z)|^2)^{1/2} \geq (1-4\lambda^2)^{1/2} =: \rho_0,
\end{equation}
with a $\rho_0 \in (0,1)$ that is as close to $1$ as we want.
We deduce from \eqref{8.5} and the first part of \eqref{8.6} that 
\begin{equation} \label{8.7}
|w-w'| \leq \rho_0^{-1} \lambda |\pi(z)-\pi(z')|. 
\end{equation}
Then since
\begin{equation} \label{8.8}
\pi(z)-\pi(z') = w e^{it} - w' e^{it'} 
= w [e^{it}-e^{it'}] + [w-w']e^{it'},
\end{equation}
\begin{equation} \label{8.9}
|e^{it}-e^{it'}| \geq w |e^{it}-e^{it'}|
\geq |\pi(z)-\pi(z')| - |w-w'| \geq 
(1- \rho_0^{-1} \lambda) |\pi(z)-\pi(z')|
\end{equation}
and
\begin{equation} \label{8.8a}
|e^{it}-e^{it'}| \leq \rho_0^{-1} \big( |\pi(z)-\pi(z')| + |w-w'| \big)
= \rho_0^{-1}(1+\rho_0^{-1} \lambda) |\pi(z)-\pi(z')|.
\end{equation}
By \eqref{8.1}, \eqref{8.8a}, and since $\H^1(\rho(a,b)) < \frac{3\pi}{4}$ by \eqref{7.1},
$e^{it}$ stays in an arc of circle of length at most $\frac{4\pi}{5}$.
There is a unique continuous determination of $t$ that comes
from inverting $e^{it}$ on that arc (in a $\frac{\pi}{2}$ Lipschitz way),
and now \eqref{8.9} implies that $\pi(z)$ is a $2$-Lipschitz function of $t$.
We write this $\pi(z) = \varphi(t)$.
Also $\pi^\perp(z)$ is a Lipschitz function of $\pi(z)$; we write this as 
$\pi^\perp(z) = \psi(\pi(z))$.
Then set $v= \psi \circ \varphi$; we see that $v$ is $2 \lambda$-Lipschitz.
With our notation, $\pi^\perp(z) = \psi \circ \varphi(t) = v(t)$,
\eqref{8.3} is the same as \eqref{8.6}, and the fact that \eqref{8.4} parameterizes 
$\Gamma$ comes from the fact that 
\begin{equation} \label{8.10}
z = \pi(z) + \pi^\perp(z) = w e^{it} + \psi \circ \varphi(t) = (w \cos t, w\sin t, v(t)).
\end{equation}

Return to the construction of $\Gamma$.
We define $F : D_T \to P^\perp$ by  
\begin{equation}\label{8.11}
F(r\cos t, r \sin t) = \frac{r v(t)}{w(t)}
\ \text{ for $r \geq 0$ and } t\in [0,T];
\end{equation}
notice that $w(t) \neq 0$, and even $w(t)-1$ stays small, because $|v(t)| \leq \lambda T$ is small.

Denote by $\Sigma'_F$ the graph of $F$ over $\overline D_T$; it easily follows from
the definitions that $\Gamma \subset \Sigma'_F$, because $v(t) = F(w(t)\cos t , w(t)\sin t)$.

The function $G$ that we construct has the following properties
(see (8.6), (8.7), (8.9) and (8.10) in \cite{C1}).
It is defined on $\overline D_T$, and it coincides with $F$ on the outer ring, i.e.,
\begin{equation}\label{8.12}
G(r\cos t, r \sin t) = F(r\cos t, r \sin t)
\ \text{ for $\frac{9}{10} \leq r \leq 1$ and } t\in [0,T].
\end{equation}
We also preserve a small region near the origin, where we may further modify the 
resulting surface: there is a small absolute constant $\kappa >0$, for which we take
\begin{equation}\label{8.13}
G(r\cos t, r \sin t) = 0
\ \text{ for $0 \leq r \leq 2\kappa$ and } t\in [0,T]
\end{equation}
(see (8.7) in \cite{C1}). Next,
\begin{equation}\label{8.14}
G \ \text{ is $C\lambda$-Lipschitz on }\overline D_T,
\end{equation}
and satisfies the Dirichlet condition
\begin{equation}\label{8.15}
G(r\cos t, r \sin t) = 0
\ \text{ when $0 \leq r \leq 1$ and } t\in \{ 0,T \}.
\end{equation}
Finally, the graph $\Sigma'_G$ of $G$ over $\overline D_T$ has a significantly smaller measure
\begin{equation}\label{8.16}
\begin{aligned}
\H^2(\Sigma'_G) &\leq \H^2(\Sigma'_F) - 10^{-4} \int_0^T |v'(t)|^2 dt 
\cr&\leq \H^2(\Sigma'_F) - 10^{-4} [\length(\Gamma) - T]
= \H^2(\Sigma'_F) - 10^{-4} \Delta_\Gamma.
\end{aligned}
\end{equation} 
For this one too, the fact that we no longer assume that $T \geq \eta_0$ does not interfere
(and indeed our bound in (8.10) of \cite{C1} does not depend on $\eta_0$).
But this is the main place where we need $\lambda$ to be small enough, so that the approximation
of the area functional by the Dirichlet energy is precise enough. 

We are interested in the intersections with the unit ball $\B$, which we denote by
\begin{equation}\label{8.17}
\Sigma_F = \Sigma'_F \cap \B \ \text{ and }\ \Sigma_G = \Sigma'_G \cap \B,
\end{equation}
and \eqref{8.16} immediately yields
\begin{equation} \label{8.18}
\H^2(\Sigma_G) \leq \H^2(\Sigma_F) - 10^{-4} \Delta_\Gamma
= \H^2(\Sigma_F) - 10^{-4} [\length(\Gamma) - \ddist(a,b)]
\end{equation}
because \eqref{8.12} says that $\Sigma'_G = \Sigma'_F$ outside of $\B$
(recall that $||v||_\infty \leq \lambda T$ is small), and by \eqref{8.1}.

In the special case when $\Gamma$ comes from $\gamma$ as in Section \ref{S7},
we can also compare with the cone 
\begin{equation} \label{8.19}
X(\gamma) = \big\{ tx \, ; \, x\in \gamma \text{ and } 0 \leq t \leq 1 \big\}
\end{equation}
over $\gamma$. Notice that
\begin{equation} \label{8.20}
\H^2(\Sigma_G) = \frac{1}{2} \length(\Gamma) \ \text{ and } \ 
\H^2(X(\gamma)) = \frac{1}{2} \length(\gamma),
\end{equation}
for instance by the area (or co-area) formula. So 
\begin{eqnarray} \label{8.21} 
\H^2(X(\gamma)) - \H^2(\Sigma_G)
\nn &=& \frac{1}{2} [\length(\gamma) - \length(\Gamma)] + 
\H^2(\Sigma_F) -\H^2(\Sigma_G)
\\ \nn
&\geq& \frac{1}{2} [\length(\gamma) - \length(\Gamma)] 
+ 10^{-4} [\length(\Gamma) - \ddist(a,b)]
\\ 
&\geq&   10^{-4} [\length(\gamma) - \ddist(a,b)]
=  10^{-4} \Delta L 
\\ \nn&\geq&  C(\lambda)^{-1} [\H^1(\Gamma \sm \gamma) + \H^1(\gamma \sm \Gamma)],
\end{eqnarray}
where we used \eqref{7.15} and \eqref{7.16}, and 
$C(\lambda)$ is a constant that depends on $\lambda$.
As we will see later, this will often mean that it is worth replacing $\gamma$ with $\Gamma$,
because the cost of gluing is often much smaller than
$\H^1(\Gamma \sm \gamma) + \H^1(\gamma \sm \Gamma)$.

\section{ Our Lipschitz net $\Gamma$: basic rules and easy cases}
\label{S9}

Recall that in Section \ref{S6} we fixed a point $\ell \in L \cap K$
(when $\ell \notin K$, $D$ stays far from $E$ and we decided to do nothing),
and we constructed a first net $\gamma$ of curves in $E \cap D$.

There are a few different configurations, but each time $\gamma$ is a union
of simple curves in $E\cap D$ (between one and five of them), 
and $\gamma$ contains the points $c_i^\ast$, $1 \leq i \leq m$, of $E \cap \d D$.
We extend $\gamma$ by adding to it the arcs $\cL_i$. Recall that 
$\cL_i$ is the arc of $E \cap \S$ that goes from $c_i^\ast$ to $a_i^\ast$,
where $a_i^\ast$ is associated to $a_i$ as in Proposition \ref{t5.4}, and
$a_i$ is the other endpoint of the arc $\cC_i$ of $K$ that passes near $c_i^\ast$.
We denote by $\overline \gamma$ the extended set, i.e., set
\begin{equation} \label{9.1}
\overline \gamma = \gamma \cup \bigcup_{1 \leq i \leq m} \cL_i.
\end{equation}

We intend to replace $\overline\gamma$ with a possibly shorter net
$\Gamma$ of Lipschitz graphs, typically constructed with the help of Section \ref{S7},
but before we start with the long list of configurations and subcases, let us explain the
main properties that we want our Lipschitz net to have.

First, $\Gamma$ should be composed of a small number (in fact, at most $4$) 
of Lipschitz curves $\Gamma_j$, disjoint except perhaps for their endpoints, and such that 
\begin{equation}\label{9.2}
\begin{aligned}
&\text{no more than $3$ curves $\Gamma_j$ ever meet at a common endpoint $z$,} 
\\&
\text{and when they do they always make angles larger than $\frac{\pi}{2}$ at $z$.}
\end{aligned}
\end{equation} 
The statement about angles contains a small abuse of notation, but we shall
fix this and say more precisely what it means near \eqref{9.8}, when we prove a similar
statement for the first time.
This condition will be very useful because later on we want to construct local Lipschitz retractions 
near $\Gamma$. When there are many curves $\Gamma_j$, we'll have additional properties
that make this possible, but which would be awkward to state here.

Also, we should say that later on, we shall consider the natural decomposition of $\gamma$
and $\Gamma$ into connected components. That is, we shall consider the component $H_i$
of $c_i^\ast$ in (the corresponding) $D\cap E$, and for each one we will rename it as 
$c \in CC$ (because more than one $i$ could give the same component $c$), consider the
piece $\gamma_c$ of $\gamma$ which is attached to $H_i$ (it is connected by construction),
and there will naturally be one connected piece $\Gamma_c$ of $\Gamma$ which corresponds.
What we will really use in later sections is more the collection of pieces $\Gamma_c$ than the 
nets $\Gamma$ themselves.

We also want a good junction with the rest of $E \cap \S$.
Recall that we started with $m$ points $c_i^\ast \in E \cap \d D$, corresponding to
$m$ curves $\cL_i$ that go from $c_i^\ast$ to $a_i^\ast$. We will make sure that
for each $i$, there is a unique index $j(i)$ such that $\Gamma_{j(i)}$ ends at $a_i^\ast$,
and moreover (with the same abuse of notation as above)
\begin{equation} \label{9.3}
\Gamma_{j(i)} \text{ makes at $a_i^\ast$ an angle larger than $\frac{\pi}{2}$
with all the $\cL_{i'}$, $i' \neq i$, that end at } a_i^\ast.
\end{equation} 

Our second demand concerns the size of the modification. We want a good estimate on the
symmetric difference
\begin{equation} \label{9.4}
\Delta(\overline\gamma ,\Gamma) 
= (\overline\gamma \sm \Gamma) \cup (\Gamma \sm \overline\gamma)
\end{equation}
in terms of the amount of surface measure that we can win. To measure this,
denote by $\rho_j$ the geodesic in $\S$ with the same endpoints as $\Gamma_j$,
and then set
\begin{equation} \label{9.5}
\rho = \cup_j \rho_j.
\end{equation}
We require that 
\begin{equation} \label{9.6}
\H^1(\Gamma) \leq \H^1(\overline \gamma)
\end{equation}
and, for some constant $C$ that depends on $\lambda$ (but not on $\tau$),
\begin{equation} \label{9.7}
\H^1(\Delta(\overline\gamma ,\Gamma)) \leq C [\H^1(\overline\gamma)- \H^1(\Gamma)]
+ C [\H^1(\Gamma)-\H^1(\rho)].
\end{equation}
We write \eqref{9.7} in this strange way because the two terms on the right-hand side 
are nonnegative (for the last one, because the $\Gamma_j$ are essentially disjoint), so
majorising by any nonnegative combination of the two pieces will be enough. 
As the reader may have guessed, we intend to win an area comparable to 
$\H^1(\Gamma)-\H^1(\rho)$ because we will apply the construction
of Section \ref{S8} to all the $\Gamma_j$ and by \eqref{8.18}.

As a last comment before we start, notice that for Configuration H and Configuration 3 = 2+1,
$\gamma$ is composed of disjoint pieces. In this case we shall construct $\Gamma$
piece by piece (i.e., independently), and take the union (it will be disjoint too).

\ms
Let us now do the construction of $\Gamma$ in the simplest cases; this will also
help us understand better as it goes.  The most interesting case will be Configuration 3+, 
which will take some time and is kept for later.

\ms
In \underbar{Configuration 0}, we have no $\gamma$ and we do nothing.

\ms
In \underbar{Configuration 1}, \eqref{6.5} says that $\gamma$ is a simple curve from $\ell$
to the unique point $x_1^\ast$ of $E \cap \d D$, $\overline\gamma$ is the simple
curve obtained by concatenating $\gamma$ and $\cL_1$, and it goes from $\ell$
to $a_1^\ast$. We apply the construction of Section \ref{S7} to $\overline \gamma$,
and get a curve $\Gamma$ with the same endpoints $\ell$ and $a_1^\ast$.
In this case \eqref{9.2} is true but pointless (there is only one $\Gamma_j$),
\eqref{9.6} comes from \eqref{7.15}, and \eqref{9.7} from the end of \eqref{8.21}.

So we are only left with \eqref{9.3} to check, and at the same time we should say
what this means. Indeed \eqref{9.3} seems to assume that $\Gamma_{j(i)}$ has a tangent
at $a_1^\ast$, but we only know that it is Lipschitz. This is easy to fix, but we need additional
definitions and notation.

When $v_1$ and $v_2$ are two unit vectors of $\R^n$, we define
the angle $\Angle(v_1,v_2)$ by
\begin{equation} \label{9.8}
\Angle(v_1,v_2) \in [0,\pi] \ \text{ and } \cos(\Angle(v_1,v_2)) = \langle v_1, v_2 \rangle.
\end{equation}
If the simple curve $\Gamma$ ends at $a$, we call \underbar{direction} of $\Gamma$
at $a$ any unit vector $v$ obtained as
\begin{equation} \label{9.9}
v = \lim_{k \to +\infty} \frac{x_k - a}{|x_k-a|},
\end{equation}
where $\{ x_k \}$ is a sequence in $\Gamma \sm \{ a \}$ that tends to $a$.
Finally, if $\Gamma_1$ and $\Gamma_2$ are two simple curves that share the endpoint $a$,
we say that
\begin{equation} \label{9.10}
\Gamma_1 \text{ and $\Gamma_2$ make an angle at least $\alpha$ at } a
\end{equation}
when 
\begin{equation} \label{9.11}
\begin{aligned}
&\Angle(v_1,v_2) \geq \alpha \text{ for every tangent direction $v_1$ of $\Gamma_1$
at $a$}
\cr& \hskip1.5cm\text{and every tangent direction $v_2$ of $\Gamma_2$ at $a$.} 
\end{aligned}
\end{equation}

We shall also use the following notation concerning geodesic directions and angles.
When $a \in \S$ and $x \in \S \sm \{-a\}$, we denote by $v(a,x)$ the direction
at $a$ of the geodesic $\rho(a,x)$ from $a$ to $x$. With this we can also compute
angles between points: we set
\begin{equation} \label{9.12}
\Angle_a(x,y) = \Angle(v(a,x),v(a,y)) \in [0,\pi];
\end{equation}
this is the angle that the geodesics $\rho(a,x)$ and $\rho(a,y)$ make at $a$.

\ms
Now we return to \eqref{9.3}. 
The proof below will actually work in many other configurations, with a minor modification that
will be explained at the end of the proof.

We shall actually prove it with an angle larger than
$\frac{7 \pi }{ 12} > \frac{\pi}{2}$. That is, we shall check that if
$j \in \cI$ is such that $j \neq i$ but $a_i^\ast$ is also an endpoint of $\cL_j$,
and $v_2$ is the tangent direction of $\cL_j$ at $a_i^\ast$ (we know that there is only one, 
since that curve is $C^1$), then 
\begin{equation} \label{9.13}
\Angle(v_1,v_2) > \frac{7 \pi}{12}
\end{equation}
whenever $v_1$ is a tangent direction of $\Gamma$ at $a_i^\ast$.
We first check that
\begin{equation} \label{9.14}
|v_1- v(a_i^\ast,\ell)| \leq 3 \lambda.
\end{equation} 
Set $a = a_i^\ast$ to simplify.
Recall from \eqref{7.13} that $\Gamma$ is a $\lambda$-Lipschitz
graph over the geodesic $\rho(a,\ell)$ with the same endpoints.
Recall from \eqref{7.14} that this means that if $\pi$ and $\pi^\perp$ denote the orthogonal
projections on the plane $P$ that contains $\rho(a,\ell)$, and on $P^\perp$ respectively,
then $\pi^\perp(z)$ is a $\lambda$-Lipschitz function of $\pi(z)$ on $\Gamma$. 
Then write our tangent direction of $\Gamma$ at $a$ as
$v_1 = \lim_{k \to +\infty} w_k$, where $w_k = (z_k-a)/|z_k-a|)$ for some
$z_k \in \Gamma$ that tends to $a$. We know that
$|\pi^\perp(z_k)| = |\pi^\perp(z_k)-\pi^\perp(a)|
\leq \lambda |\pi(z_k)-\pi(a)|$ and hence, since 
$|z_k-a| \geq |\pi(z_k)-\pi(a)| - |\pi^\perp(z_k)-\pi^\perp(a)|
\geq (1 - \lambda)|\pi(z_k)-\pi(a)|$, that
\begin{equation} \label{9.15}
|\pi^\perp (w_k)| = \frac{|\pi^\perp (z_k-a)|}{|z_k-a|} \leq \frac{\lambda}{1 - \lambda} \leq 2\lambda 
\end{equation}
if $\lambda$ is chosen small enough. Then $|\pi^\perp (v_1)| \leq 2\lambda$ too.
Since $v_1$ lies on the tangent hyperplane to $\S$ at $a$, we get that its projection
$\pi(v_1)$ lies in the direction of $\rho(a,\ell)$ at $a$ and finally \eqref{9.14} follows.
In many other configurations, we will still know that $\Gamma_{j(i)}$ is a small Lipschitz graph,
but sometimes over a slightly different geodesic $\rho(a_i^\ast,b)$, where $b$ is quite close
to $\ell$. Then \eqref{9.14} will follow as above, but maybe with $C\lambda$ instead
of $\lambda$ (coming from $|v(a_i^\ast,b)-v(a_i^\ast,\ell)|$, or a slightly larger Lipschitz constant).
The rest of the proof will work unchanged.

Notice that $|a_i^\ast-a_i| \leq 2\cdot 10^{-10} \tau$ by \eqref{5.26},
and $\ddist(a_i,\ell) = \H^1(\cC_i) \geq 10\eta(X) > 10^4 \tau$ by \eqref{5.3}
and \eqref{5.2bis}, so we also get that
\begin{equation} \label{9.16}
|v(a_i^\ast,\ell)-v(a_i,\ell)| \leq 10^{-10}.
\end{equation}

Now we consider $v_2$. First assume that $a_i$ is a true vertex, i.e., that $a_i \in V_1$.
One of the arcs of $K$ leaving from $a_i$ is $\cC_i = \rho(a_i,\ell)$, and $\cC_j$,
the arc of $K$ that lies close to $\cL_j$, is another one. Write
$\cC_j = \rho(a_i,b)$; then
\begin{equation} \label{9.17}
\Angle(v(a_i,b), v(a_i,\ell)) = \frac{2 \pi}{3} = \frac{8 \pi}{12},
\end{equation}
because it is the angle of $\cC_i$ and $\cC_j$ at $a_i$.
Now we apply \eqref{5.23}, to the vertex $x = a_i$ and the point $z=a \in \cL_j$. 
With the notation below \eqref{5.22}, $v_j(z) = v_2$ (the direction of $\cL_j$ at $a$)
and $v_j = v(a_i,b)$ (the direction at $x=a_i$ of the tangent to $\cC_j = \rho(a_i,b)$).
Thus \eqref{5.23} says that
\begin{equation} \label{9.18}
|v_2 - v(a_i,b)| = |v_j(z) - v_j| \leq 1/30.
\end{equation}
In the present case, our goal \eqref{9.13} follows easily from \eqref{9.14}-\eqref{9.18},
but we still need to treat the case when $a_i \in V_2$.

We still have \eqref{9.14} (for the same reasons) and \eqref{9.16} (see below \eqref{5.36}).
Now there is only one other $\cC_j$ leaving from $a_j$, and 
\begin{equation} \label{9.19}
\Angle(v(a_i,b), v(a_i,\ell)) = \pi,
\end{equation}
i.e., $\cC_j$ is a continuation of $\cC_i$. We still have \eqref{9.18}, but this time
we apply \eqref{5.35}, and then we conclude as above.

This completes our verification of \eqref{9.2}-\eqref{9.7} in the case of Configuration 1.

Notice that in our proof of \eqref{9.3}, if instead of ending at $\ell$ the curve $\Gamma$ 
ends at some other point $x_0 \in D$, then we just need to replace $v(a_i^\ast,\ell)$ 
with $v(a_i^\ast,x_0)$ in \eqref{9.14}, and add $|v(a_i^\ast,x_0) - v(a_i^\ast,\ell)|$ 
to the error term in \eqref{9.16}; but 
\begin{equation} \label{9.19aa}
|v(a_i^\ast,x_0) - v(a_i^\ast,\ell)| \leq 2 |x_0-\ell| |a_i^\ast-\ell|^{-1}
\leq \frac{2 \tau }{ 8 \cdot 10^{3} \tau} \leq \frac{1}{4000}
\end{equation}
(see the estimate just above \eqref{9.16}); this still gives \eqref{9.13} and \eqref{9.3}.

\ms
Our next case is \underbar{Configuration $2-$}.
In this case $\gamma$ is a simple curve in $E \cap D$ that goes from $c_1^\ast$
to $c_2^\ast$, and does not contain $\ell$. 
Select a point $x_0 \in \gamma \cap D$, for instance a point that  
minimizes the distance to $\ell$, and cut $\gamma$ into two essentially disjoint
simple curves $\gamma_1$ and $\gamma_2$, where $\gamma_i$ goes from 
$x_0$ to $c_i^\ast$.
Then extend $\gamma_i$, by adding to it the corresponding arc $\cL_i$; 
we assume that our notation is such that $\cL_i$ is the curve in $E \cap \S$
that contains $c_i^\ast$. This gives a curve $\ol\gamma_i$, that goes from 
$x_0$ to $a_i^\ast$.

We apply the construction of Section \ref{S7} to $\ol\gamma_i$
and get a curve $\Gamma_i$ with the same endpoints $x_0$ and $a_i^\ast$.
Then we take $\Gamma = \Gamma_1 \cup \Gamma_2$.

Let us check that
\begin{equation}\label{9.19a}
\text{the two $\Gamma_i$ make an angle larger than $110^\circ$ at $x_0$;}
\end{equation}
of course  \eqref{9.2} will follow (the only interior vertex of $\Gamma$ is $x_0$).
We shall merely use the fact that $|x_0 - \ell| \leq \tau$, even though we could rather 
easily deduce from \eqref{4.3} that $\gamma$ comes within $10\varepsilon$ of $\ell$.
We first control the direction of the geodesic $\rho(x_0,a_i^\ast)$ over which
$\Gamma_i$ is a small Lipschitz graph. Recall from \eqref{5.38} that
$|a_i^\ast-a_i| \leq 10^{-9}\tau$, while \eqref{5.3} and \eqref{5.4}
imply that $|a_i-\ell| \geq 9 \cdot 10^3 \tau$. 
Thus $|x_0-a_i^\ast| \geq 8 \cdot 10^3 \tau$ and
\begin{equation} \label{9.19b}
\begin{aligned}
|v(x_0,a_i^\ast) - v(\ell,a_i)|
&\leq |v(x_0,a_i^\ast) - v(\ell,a_i^\ast)|+|v(\ell,a_i^\ast) - v(\ell,a_i)|
\cr&
\leq 2|x_0-\ell| |x_0-a_i^\ast|^{-1} + 2 |a_i^\ast-a_i| |\ell - a_i|^{-1} \leq 10^{-3}.
\end{aligned}
\end{equation}
We know from \eqref{2.7} that $\Angle(v(\ell,a_1),v(\ell,a_2)) \geq \frac{2\pi}{3}$
so \eqref{9.19b} gives a good control on $\Angle(v(x_0,a_1^\ast),v(x_0,a_2^\ast))$,
and \eqref{9.19a} will follow as soon as we check that for $i=1,2$,
\begin{equation} \label{9.19c}
\Angle(v_i,v(x_0,a_i^\ast)) \leq 3\lambda
\ \text{ when $v_i$ is any tangent direction to $\Gamma_i$ at $x_0$.} 
\end{equation}
But this is true, and the proof is the same as for \eqref{9.14}.

Next we are supposed to check \eqref{9.3}, but the proof of \eqref{9.13} still
works in the present case, as explained near \eqref{9.19aa}.
Then \eqref{9.6} follows from \eqref{7.15} (we apply it to each piece,
and then sum), and \eqref{9.7} from the end of \eqref{8.21} (again sum the 
two pieces). This completes our verifications for Configuration $2-$.

\ms
Now we switch to \underbar{Configuration $3 = 2+1$}. This will just be a combination
of Configurations $1$ and $2-$.
Here, possibly after relabelling, $\gamma$ is composed of an arc $\gamma_{1,2}$
that goes from $c_1^\ast$ to $c_2^\ast$, and an arc $\gamma_3$ from
$\ell$ to $c_3^\ast$. We select an origin $x_0 \in \gamma_{1,2} \cap D$,
and in fact the simplest is to take $x_0 = c_1^\ast$. This way we have
three arcs, $\gamma_1 = \{ c_1^\ast \}$,  $\gamma_2 = \gamma_{1,2}$,
and $\gamma_3$, which we extend as before. This gives three arcs,
$\overline\gamma_{i}$, two that leave from $x_0 = c_1^\ast$ and one from $\ell$. 
Then we apply the construction of Section \ref{S7} independently to the 
three $\overline\gamma_{i}$ and get small Lipschitz graphs $\Gamma_i$. 
Finally we set $\Gamma = \Gamma_1 \cup \Gamma_2 \cup \Gamma_3$.

The curves $\Gamma_1$ and $\Gamma_2$ have a common endpoint $x_0$,
and by the proof of \eqref{9.19a} they make an angle larger than $110^\circ$
at $x_0$. We claim that 
\begin{equation} \label{9.19cc}
\text{$\Gamma_3$ does not meet $\Gamma_1 \cup \Gamma_2$.}
\end{equation}
Let us first check that for $i=1,2$,
\begin{equation} \label{9.19d}
\dist(z,\rho(x_0,a_i^\ast)) \leq 2\lambda |z-x_0| 
\ \text{ for } z\in \Gamma_i.
\end{equation}
Let $P$ be that plane that contains $\rho_i = \rho(x_0,a_i^\ast)$, and denote
by $\pi$ and $\pi^\perp$ the orthogonal projections on $P$ and its
orthogonal complement; by \eqref{7.14}, $\pi^\perp(z)$ is a $\lambda$-Lipschitz
function of $\pi(z)$ (hence also of $z$) on $\Gamma_i$. 
This implies that 
\begin{equation} \label{9.19e}
|\pi^\perp(z)| = |\pi^\perp(z)-\pi^\perp(x_0)| \leq \lambda |z-x_0|
\ \text{ for } z\in \Gamma_i.
\end{equation}
Next $|\pi(z)| = (1-|\pi^\perp(z)|^2)^{1/2}$ is a $2\lambda$-Lipschitz
function of $\pi(z)$ (differentiate $f(x) = (1-x^2)^{1/2}$ near $1$).
Now set $\xi(z) = \pi(z)/|\pi(z)|$ (a projection on the circle that contains $\rho_i$); then
\begin{equation} \label{9.19f}
\begin{aligned}
|\xi(z) -\xi(z')| &\geq \frac{|\pi(z)-\pi(z')|}{|\pi(z)|} 
- |\pi(z')| \Big|\frac{1}{|\pi(z)|}-\frac{1}{|\pi(z')|} \Big|
\cr&\geq \frac{|\pi(z)-\pi(z')| }{ 1+\lambda|z-x_0|} - \Big|\frac{|\pi(z)|-|\pi(z')|}{|\pi(z)|}\Big|
\cr&\geq \frac{|\pi(z)-\pi(z')| }{ 1+2\lambda} - 3\lambda |\pi(z)-\pi(z')| 
\geq (1-5\lambda) |\pi(z)-\pi(z')| 
\end{aligned}
\end{equation}
by \eqref{9.19e} and because $|z-x_0| \leq 2$. This shows that $\pi(z)$
is a Lipschitz function of $\xi(z)$, and (since $\pi^\perp(z)$ and hence $z$
are Lipschitz functions of $\pi(z)$), we see that $z$ is a Lipschitz function of 
$\xi(z)$. In particular, $\xi(z) \neq \xi(z')$ when $z\neq z'$, and this implies
that $\xi(z)$ stays on the geodesic $\rho_i$ (instead of wandering
somewhere else on the circle). Also (we don't need this now), we can use
$\xi(z) \in \rho_i$ to parameterize $\Gamma_i$ in a Lipschitz way.
But since $\xi(z)\in\rho_i$,
\begin{equation} \label{9.19g}
\begin{aligned}
\dist(z,\rho_i) &\leq |z-\xi(z)| \leq |z-\pi(z)| + |\pi(z)-\xi(z)|
\cr&= |\pi^\perp(z)| + (1-|\pi(z)|) \leq 2|\pi^\perp(z)| \leq 2\lambda |z-x_0|
\end{aligned}
\end{equation}
by \eqref{9.19e}. This proves \eqref{9.19d}; the same argument shows that
\begin{equation} \label{9.19h}
\dist(z,\rho(\ell,a_3^\ast)) \leq 2\lambda |z-\ell| \ \text{ for } z\in \Gamma_3.
\end{equation}
This is where our choice of $x_0=c_1^\ast$ makes our life more comfortable.
Recall from \eqref{2.7} that $\cC_1$ and $\cC_3$ make an angle of at least $\frac{2\pi}{3}$
at $\ell$ (in the present case, we have $3$ curves $\cC_i$, so the angle is in fact $\frac{2\pi}{3}$).
This implies that $\Angle(v(\ell,x_0), v(\ell,a_i^\ast)) \geq \frac{2\pi}{3}-10^{-2}$, say,
and then \eqref{9.19cc} follows rather easily from \eqref{9.19d} and \eqref{9.19h} 
(but we skip the details and instead encourage the reader to draw a picture).

So our set $\Gamma$ is composed of two connected pieces, $\Gamma_3$
and $\Gamma_1 \cup \Gamma_2$, which are disjoint (one could even check
that their distance is at least $\tau/2$). They both satisfy \eqref{9.2}: for
$\Gamma_3$ this is trivial, and for $\Gamma_1 \cup \Gamma_2$ the proof is
the same as for Configuration $2-$. They also satisfy \eqref{9.3}, by the
proof of \eqref{9.13} and the remark near \eqref{9.19aa}.
Finally \eqref{9.6} and \eqref{9.7} are proved piece by piece, and follow from 
\eqref{7.15} and the end of \eqref{8.21}, as before.

This completes our verification for Configuration $3=2+1$.
Notice however that the net $\Gamma$ that we construct is far from
optimal: in the present situation, since $\gamma_1$ and $\gamma_2$ make
an angle of nearly $120^\circ$ near $\ell$, we could easily organize a much
more brutal shortcut, and save a lot of length. But we choose a way which is
easier to handle with the same estimates as in the other cases.
The fact that our competitor is not so good will show up later, when
we will see that if our competitor looks like a cone over $\Gamma$ in a small ball,
we can easily improve on it.

\ms
We are almost ready for \underbar{Configuration H}. 
For each of the hanging curves $\cL_i$ (those for which $c_i^\ast$ is 
not connected to $\ell$ or any other $c_i^\ast$), we kept the curve
$\gamma_i = \cL_i$, and the simplest is to take $\Gamma_i = \cL_i$ too.
This is, if we are ready to use the fact that if we took $\tau$ and $\varepsilon$ 
small enough, depending on $\lambda$, the curve $\cL_i$ is automatically a $\lambda$-Lipschitz
graph. Otherwise, we apply the construction of Section \ref{S7} to $\cL_i$, as we did 
in the previous cases, to get a Lipschitz graph $\Gamma_i$.

Of course this does not look glorious: we should rather have cut off the whole $\cL_i$
and saved a lot of length, but this is a way for us to make our construction more uniform.
Later on, we will notice with apparent surprise that we can still cut off the geodesic 
$\rho(c_i^\ast,a_i^\ast)$ from a net of geodesics, and save some length, and this will
compensate the present laziness.

There still may be one or two $c_j^\ast$ left, that are connected to something. 
If they are connected as in Configuration $1$, i.e., if there is only one $c_j^\ast$ 
left and it is connected to $\ell$, let $\gamma_j$ be the arc of $E\cap D$
that was selected above, extend it to get an arc $\ol\gamma_j$ that goes from
$\ell$ to $a_j^\ast$, and let $\Gamma_j$ be obtained by applying
the construction of Section \ref{S7} to $\ol\gamma_j$. Then let 
$\Gamma$ be the union of $\Gamma_j$ with the hanging graphs $\Gamma_i$
that we already selected. The proof of \eqref{9.19cc} still works here
and shows that $\Gamma_j$ is disjoint from these curves.

When the remaining $c_j^\ast$ are connected as in Configuration $2-$, 
we have one index $i$ and two indices $j$, which we label so that $i=3$.
We construct $\Gamma_1$ and $\Gamma_2$ exactly as we did in
Configuration $2-$, and set $\Gamma = \Gamma_1 \cup \Gamma_2 \cup \Gamma_3$.
Again $\Gamma_1 \cup \Gamma_2$ does not meet $\Gamma_3 = \cL_3$,
by the proof of \eqref{9.19cc} (and you may find it more convenient to choose 
$x_0 = c^\ast_1$ as the center of $\gamma$ where you cut the curve).

We are left with the case when there is only one hanging $c_i^\ast$, which
we call $c_3^\ast$, and $c_1^\ast$ and $c_2^\ast$ are connected as in
Configuration $2+$. We did not treat the case of Configuration $2+$ yet,
but we shall do it later, and there will be no loophole. The construction described
below, performed with the connected set that connects $c_1^\ast$, $c_2^\ast$, 
and $\ell$, will give a net of curves $\Gamma'$; then we take 
$\Gamma = \Gamma_3 \cup \Gamma'$, the local description (with \eqref{9.2} and
\eqref{9.3}) can be proved independently for the two pieces, and the fact
that $\Gamma' \cap \Gamma_3 = \emptyset$ will be true, as in \eqref{9.19cc}.
See the remarks below \eqref{11.4} and above \eqref{11.14}.

In all these subcases, we get a disjoint union of curves or nets that satisfy
the conditions \eqref{9.2} and \eqref{9.3}, as in the single configurations and 
for the same reasons. 

Notice that for the first time we get curves that end at a point 
other than $\ell \in L$. This is not bad in itself; it means that our future competitor 
is rather poor, but this is all right. In fact it means that Configuration H will not happen.

Finally, \eqref{9.6} and \eqref{9.7} are checked piece by piece, with the same estimates as for
the other configurations.
This completes our discussion in Configuration H.

\ms
The last simple case is \underbar{Configuration 3-}. In this case 
$\gamma = \gamma_1 \cup \gamma_2 \cup \gamma_3$, three almost disjoint curves
that start from the same origin $x_0$. We add the corresponding $\cL_i$ and get 
curves $\ol\gamma_i$ from $x_0$ to the $a_i^\ast$. Finally 
$\ol\gamma = \bigcup_i \ol\gamma_i$. We apply the construction of Section \ref{S7}
and get three curves $\Gamma_i$, with the same endpoints as the $\ol\gamma_i$.
Finally we take $\Gamma = \cup_i \Gamma_i$.

The fact that \eqref{9.2} holds, and in fact
\begin{equation}\label{9.31}
\text{the three $\Gamma_i$ make angles larger than $110^\circ$ at $x_0$,}
\end{equation}
is proved just like \eqref{9.19a} above; we could also use the argument that will be given
for Configurations $3+$, above \eqref{10.23}.  

As usual \eqref{9.3} holds for the same reason as in Configuration 1;
see the proof of \eqref{9.13} and the comment near \eqref{9.19aa}.

Finally \eqref{9.6} follows from \eqref{7.15} and \eqref{9.7} from the end of \eqref{8.21};
as before we just have to add the three estimates for the three $\Gamma_i$.
This completes our verification for Configuration~$3-$.

\ms
We are left with two more complicated cases, Configurations $2+$
and $3+$, which we deal with in the next two sections.

\section{Our net $\Gamma$ in Configuration 3+} 
\label{S10} 

In the two remaining cases, there is a small additional difficulty, due to the fact
that the construction of Section \ref{S7} was meant to cut curves and get shorter
Lipschitz curves, and we do not seem to have a corresponding simple construction for 
$3$-legged spiders. Instead we will distinguish between many cases, and construct
different acceptable nets of Lipschitz curves.
Again we want to be prudent, because we do not want to replace large portions of our spiders 
if we do not save a comparable amount of surface later. As before, this saving will also come
from comparing cones with harmonic graphs,  but often we shall first try to make 
$\H^1(\ol\gamma) - \H^1(\Gamma)$ large.

Thus, rather than trying to make a nice general construction for spiders, we shall use our
construction for curves and try to fix by hand the obvious problems near the center. 

In this section we study the case of Configuration 3+, 
which appears to be the most complicated. 
Configuration 2+ will be slightly easier, and will be treated in Section \ref{S11}.

\subsection{Preparation}

We start with some notation.
Recall that we constructed in Section \ref{S6} a net $\gamma$, which is a possibly
degenerate spider with three long legs and a short tail $\gamma_\ell$. 
The short tail ends at $\ell$, and the three legs end at points $c_i^\ast$, 
$1 \leq i \leq 3$. Denote by $x_0$ the center of the spider, i.e., the point 
where $\gamma_3$ meets $\gamma_{1,2}$. Also denote by $\gamma_1$
and $\gamma_2$, respectively, the arc of $\gamma_{1,2}$ between $x_0$
and $c_1^\ast$ and $c_2^\ast$. Thus the three $\gamma_i$ are essentially disjoint, and
$\gamma = \gamma_\ell \cup \big(\bigcup_{i=1}^3 \gamma_i \big)$.

As usual we extend the three legs $\gamma_i$ by adding the corresponding 
curves $\cL_i \subset E \cap \S$ that go from the $c_i^\ast$ to the $a_i^\ast$; 
this gives three essentially disjoints simple curves $\ol\gamma_i \subset E\cap \S$. 
We set 
\begin{equation} \label{10.1}
\ol \gamma = \gamma_\ell \cup \big(\bigcup_{i=1}^3 \ol\gamma_i \big).
\end{equation}

Let us apply the construction of Section \ref{S7} to each of the curves $\ol\gamma_i$; 
we get a Lipschitz graph $\Gamma_i$ with a small constant $\lambda$, with the same endpoints 
$x_0$ and $a_i^\ast$. Then we set 
\begin{equation}\label{10.2}
\Gamma^\ast = \bigcup_{i=1}^3 \Gamma_i \, .
\end{equation}
We like $\Gamma^\ast$ because, as we shall see, it is a nice looking spiral.
In Configuration 3-, we decided to take $\Gamma = \Gamma^\ast$; here things will not be
so simple, because we have to take care of the special point $\ell$. In the mean time,
let us derive some simple properties of $\Gamma^\ast$.
The next lemma is also valid in Configuration 3-.

\begin{lem}\label{t10.1} 
For each small constant $\alpha < 1$, we can find $\varepsilon(\alpha) > 0$ 
such that if we take $\varepsilon < \varepsilon(\alpha)$ in \eqref{4.3}, then 
\begin{equation} \label{10.3}
|x_0 - \ell| + \sum_{i=1}^3 |a_i^\ast-a_i| \leq 2\alpha^2 \tau.
\end{equation}
\end{lem}

We state this with quantifiers to avoid any suspicion % checked
of loopholes. In practice, we will apply this with a small constant $\alpha > 0$,
that will be chosen later in this section, depending on various geometric constants
and our choice of $\lambda$.
And we shall make sure that $\varepsilon$ is so small that \eqref{10.3} holds. 

Let us apply Proposition \ref{t5.4}, but with the smaller constant $\alpha^2 \tau$; 
this forces us to take $\varepsilon$ even smaller than before, but this is all right.
We get a description of $E \cap \S \sm (D_+(\alpha^2 \tau) \cup D_-(\alpha^2 \tau))$ 
as a union of simple curves $\cL'_i$, $i\in \cI$. Of course this description matches
the description that we used for $\tau$ (i.e., with the $\cL_i$); in particular,
the vertices $x^\ast$ that show up in \eqref{5.38} are the same for $\alpha^2\tau$
as for $\tau$, even when $x \in V_2$, because of the way we chose them
(below \eqref{5.36}, so that $|x^\ast-x|$ is minimal). Thus the part of \eqref{10.3} 
that comes from the $a_i$ follows from \eqref{5.38} with $\alpha^2\tau$.

Now we concentrate on what happens in the spherical annulus
$A = \S \cap B(\ell,2\tau) \sm B(\ell,\alpha^2\tau)$. 
Here the curves $\cL'_i$ lie at distances at least $\alpha^2\tau/10$ 
from each other (by \eqref{5.39} and because the $\cC_i$ are far from 
each other in $A$); then $E \cap \S$ has no triple point in $A$, i.e.,
points like $x_0$ near which $E\cap \S$ is composed of three short simple curves 
leaving from $x_0$, that are disjoint except for $x_0$. 
This proves that $x_0 \in B(\ell,\alpha^2\tau)$, as needed.
\qed

\ms
We will do lots of little computations with small Lipschitz graphs over geodesics, and
the definition \eqref{7.14} that we gave in Section \ref{S7} is not so pleasant. 
Next we observe that when we restrict to a small enough spherical disks, \eqref{7.14}
yields a definition of small Lipschitz graphs that looks a lot like the usual one.
Some notation will be useful. Set
\begin{equation} \label{10.4}
B_1 = \S \cap B(x_0,\lambda) 
\end{equation}
($\lambda$ is the scale at which our approximation will start being less good)
and, for $1 \leq i \leq 3$, 
\begin{equation} \label{10.5}
e_i = v(x_0,a_i^\ast).
\end{equation}
Also denote by $P_i$ the vector plane that contains $\rho(x_0,a_i^\ast)$,
by $P_i^\perp$ its orthogonal complement, by $\pi_i$ and $\pi_i^\perp$
the orthogonal projections on $P_i$ and $P_i^\perp$, 
and by $p_i$ and $p_i^\perp$ the orthogonal projections on the vector lines
through $e_i$ and $x_0$ respectively. Notice that
\begin{equation} \label{10.6}
I = \pi_i + \pi_i^\perp = p_i + p_i^\perp + \pi_i^\perp.
\end{equation}
Recall that 
\begin{equation} \label{10.7}
\text{$\Gamma_i$ is a $\lambda$-Lipschitz graph over $\rho(x_0,a_i^\ast)$.}
\end{equation}
By \eqref{7.14}, this means that ($\Gamma_i$ is a curve with the given endpoints and that)
on $\Gamma_i$, $\pi_i^\perp(z)$ is a $\lambda$-Lipschitz function of $\pi_i(z)$. 
Since $\pi_i$ is $1$-Lipschitz, we immediately get that
\begin{equation} \label{10.8}
\pi_i^\perp \ \text{ is $\lambda$-Lipschitz on } \Gamma_i. 
\end{equation}
In addition, we claim that 
\begin{equation} \label{10.9}
p_i^\perp \ \text{ is $\frac{10\lambda}{9}$-Lipschitz on } B_1. 
\end{equation}
This is easy, but we prove it anyway. Let $z_1, z_2 \in B_1$ be given; 
for $j = 1, 2$, write $z_j = p_i^\perp(z_j) + w_j$,
with $w_j \perp x_0$. Then $|w_j| \leq \lambda$ ($w_j$ is a $1$-Lipschitz function of $z_j$, 
null when $z=x_0$), $|p_i^\perp(z_j)|^2 = 1 - |w_j|^2$, and hence
$\langle x_0, z_j \rangle = (1 - |w_j|^2)^{1/2}$ (it is obviously positive, since $z_j$
is close to $x_0$). Now
\begin{equation} \label{10.10}
\begin{aligned}
|p_i^\perp(z_1) - p_i^\perp(z_2)| &= |\langle x_0, z_1-z_2 \rangle|
= \big|(1 - |w_1|^2)^{1/2} - (1 - |w_2|^2)^{1/2} \big|
\cr&
\leq \frac{10}{9}\, \lambda |w_1-w_2| \leq  \frac{10}{9}\, \lambda |z_1-z_2| 
\end{aligned}
\end{equation}
(just notice that the derivative of $(1 - x^2)^{1/2}$ is $x(1 - x^2)^{-1/2}$ and
estimate). So \eqref{10.9} holds.

We deduce from \eqref{10.6}-\eqref{10.9} that for $z, z' \in \Gamma_i \cap B_1$,
\begin{equation} \label{10.11}
|p_i(z)-p_i(z')| \geq |z-z'| - |p_i^\perp(z)-p_i^\perp(z')| - |\pi_i^\perp(z)-\pi_i^\perp(z')| 
\geq \Big(1-\frac{19 \lambda}{9} \Big) |z-z'|,
\end{equation}
and so (using \eqref{10.6}-\eqref{10.9} again)
\begin{equation} \label{10.12}
\Gamma_i \cap B_1 \ \text{ is a $3\lambda$-Lipschitz graph over } 
p_i(\Gamma_i \cap B_1) \subset {\rm Vect}(e_i);
\end{equation}
the fact that
\begin{equation} \label{10.13}
p_i(\Gamma_i \cap B_1) \supset [0,1-4\lambda] e_i
\end{equation}
easily follows from \eqref{10.12}, the fact that $\Gamma_i$ starts from $x_0$
in the direction of $e_i$, and a continuity argument.

This description of $\Gamma_1 \cap B_1$ will be easier to use than the initial
definition with \eqref{7.14}. 
There is also a converse that we want to record.

\begin{lem} \label{t10.2}
Let $\Gamma'$ be a curve that goes from $x_0$ to $a_i^\ast$ and
coincides with $\Gamma_i$ on $\S \sm B(x_0, \lambda/10)$. Suppose in
addition that for some $A \in [1,100]$,
\begin{equation} \label{10.14}
\Gamma'_i \cap B_1 \ \text{ is a $A\lambda$-Lipschitz graph over } 
p_i(\Gamma_i \cap B_1).
\end{equation}
Then $\Gamma'$ is a $2A\lambda$-Lipschitz graph over $\rho(x_0,a_i^\ast)$.
\end{lem}

We just need to check that on $\Gamma'$, $\pi_i^\perp(z)$
is a $2A \lambda$-Lipschitz function of $\pi_i(z)$. 
This is true on $\Gamma' \cap B_1$, because
the orthogonal projection on the direction perpendicular to $e_i$ 
(call it $p = I - p_i$) dominates the orthogonal projection $\pi_i^\perp$, so that
\begin{equation} \label{10.15n} 
|\pi_i^\perp(z)-\pi_i^\perp(z')| \leq |p(z)-p(z')| 
\leq A\lambda |p_i(z)-p_i(z')| 
\leq A\lambda |\pi_i(z) - \pi_i(z')|  
\end{equation}
for $z,z' \in \Gamma' \cap B_1$. We also have this on $\Gamma' \sm B(x_0,\lambda/10)$,
by definition, so we just need to show that
\begin{equation} \label{10.15}
|\pi_i^\perp(z) - \pi_i^\perp(z')| \leq 2A\lambda |\pi_i(z) - \pi_i(z')|
\end{equation}
when $z'\in \Gamma' \cap B(x_0, \lambda/10)$ and 
$z\in \Gamma' \sm B_1 = \Gamma_i \sm B_1$. By \eqref{10.8} 
\begin{equation} \label{10.16}
|\pi_i^\perp(z)-\pi_i^\perp(x_0)| \leq \lambda |z-x_0| 
\end{equation}
so
\begin{equation} \label{10.17}
|\pi_i(z)-\pi_i(x_0)| \geq |z-x_0|-|\pi_i^\perp(z)-\pi_i^\perp(x_0)| 
\geq (1-\lambda) |z-x_0|.
\end{equation}
Similarly, \eqref{10.15n} 
implies that
\begin{equation} \label{10.18}
|\pi_i^\perp(z')-\pi_i^\perp(x_0)| 
\leq A\lambda |z'-x_0| \leq A \lambda^2/10
\end{equation}
and, since 
\begin{equation} \label{10.19}
|\pi_i(z')-\pi_i(x_0)| \leq |z'-x_0| \leq \lambda/10,
\end{equation}
we get that
\begin{equation} \label{10.20}
\begin{aligned}
|\pi_i(z)-\pi_i(z')| &\geq |\pi_i(z)-\pi_i(x_0)| - \lambda/10
\geq (1-\lambda) |z-x_0| - \lambda/10 
\cr& \geq (1-\lambda- \frac{1}{10}) |z-x_0| \geq \frac{8 |z-x_0|}{10}
\end{aligned}
\end{equation}
because $|z-x_0| \geq \lambda$. 
In addition
\begin{equation} \label{10.21}
|\pi_i^\perp(z')-\pi_i^\perp(z)| \leq \lambda |z-x_0| + A \lambda^2/10
\leq (\lambda+A\lambda/10) |z-x_0|, 
\end{equation}
and \eqref{10.15} follows, because $\frac{10}{8} \, (\lambda+A\lambda/10) \leq 2A\lambda$
when $A \geq 1$.
\qed

\ms
Let us also record that for $1 \leq i \leq 3$, $e_i$ is quite close to
the direction $v(\ell,a_i)$ of $\cC_i$ at $\ell$:
\begin{eqnarray} \label{10.22}
|e_i - v(\ell,a_i)| &=&
|v(x_0,a_i^\ast)-v(\ell,a_i)| 
\leq |v(x_0,a_i^\ast)-v(\ell,a_i^\ast)|+|v(\ell,a_i^\ast)- v(\ell,a_i)|
\nn\\
&\leq& 2 |x_0-\ell| |\ell-a_i^\ast|^{-1} + 2|a_i^\ast-a_i| |\ell-a_i^\ast|^{-1} 
\leq 4 \alpha^2 \tau |\ell-a_i^\ast|^{-1} 
\nn\\
&\leq& 4 \alpha^2 \tau [5 \eta(X)]^{-1} \leq 10^{-3} \alpha^2
\end{eqnarray}
by \eqref{10.3}, then \eqref{3.11}, \eqref{3.12}, and \eqref{5.2bis}.

It follows from \eqref{10.22} and the fact that $\Gamma_i$ is a small
Lipschitz graph over $\rho(x_0,a_i^\ast)$ (or more directly \eqref{10.12})
that
\begin{equation}\label{10.23}
\text{the three $\Gamma_i$ make angles larger than $110^\circ$ at $x_0$.}
\end{equation}
Notice that we only used \eqref{10.7} here, so \eqref{10.23} is also valid 
in the case of Configuration 3-, therefore proving \eqref{9.31} and completing
the discussion for this case.

\subsection{Case A, where we force $\Gamma^\ast$ to be centered at $\ell$}

We return to Configuration 3+.
Even though $\Gamma^\ast = \cup_i \Gamma_i$ is nice, we shall need to modify it
because we want $\Gamma$ to contain $\ell$ too, and the success of the construction 
will depend on various parameters such as the relative position of $\ell$ and the $\Gamma_i$.

In this subsection we try to modify $\Gamma^\ast$ in the following simple way:
we shall select points $z_i \in \Gamma_i$, rather far from the center, and replace
the three arcs of $\Gamma_i$ between $x_0$ and the $z_i$ with a spider $Y$ 
centered at $\ell$ and composed of geodesic arcs. This will turn out to work well when
\begin{equation}\label{10.24} 
\H^1(\gamma_\ell) + \sum_{i=1}^3 [\H^1(\ol\gamma_i) - \H^1(\Gamma_i)] 
\geq 32 \lambda |x_0 - \ell|.  
\end{equation}
We call this \underbar{Case A}. Incidentally, the constant 32 is computed backwards to make the proof
work; a mistake in the computations would probably force us to  make it larger, but this would not be bad.
Let $\alpha > 0$ be small (compared to $\lambda$), decide to choose 
$\varepsilon$ smaller than $\varepsilon(\alpha)$ from Lemma \ref{t10.1},  and set
\begin{equation} \label{10.25}
r = \alpha^{-1} |x_0 - \ell|,  
D = \S \cap B(x_0,r), \text{ and } \d D = \S \cap \d B(x_0,r).
\end{equation}
The notation is the same as with the disks $D_\pm$ centered at $\ell_\pm$ above,
but this will be a different spherical disk and sphere. We promise no conflict of notation.
Notice that since $|x_0 - \ell| \leq 2 \alpha^2 \tau$, we get that
$r = \alpha^{-1} |x_0 - \ell| \leq 2 \alpha \tau$ and 
\begin{equation} \label{10.26}
D \subset \S \cap B(x_0,2\alpha \tau) \subset B(x_0, 10^{-3}\lambda),
\end{equation}
if $\alpha$ is small enough (we could also have relied on $\tau$ being small), 
so we can use the Lipschitz description \eqref{10.12}
of $\Gamma_i \cap D$. In particular, each $\Gamma_i$
meets $\d D$ exactly once, at a point which we call $z_i$. Set
\begin{equation} \label{10.27}
\Gamma''_i = \Gamma_i \cap D \ \text{ and } 
\Gamma'' = \bigcup_{i=1}^3 \Gamma''_i = \Gamma^\ast \cap D.
\end{equation}
Thus $\Gamma''_i$ is the arc of $\Gamma_i$ that goes from $x_0$ to $z_i$.
Also set
\begin{equation} \label{10.28}
\Gamma'_i = \Gamma_i \sm D \ \text{ and } 
\Gamma' = \bigcup_{i=1}^3 \Gamma'_i = \Gamma^\ast \sm D
\end{equation}
(the exterior part); we want to replace $\Gamma''$ with the spider
\begin{equation} \label{10.29}
Y = \bigcup_{i=1}^3 \rho(\ell,z_i),
\end{equation}
which has the advantage of containing $\ell$. So we set
\begin{equation} \label{10.30}
\wt\Gamma_i = \rho(\ell,z_i) \cup \Gamma'_i
\ \text{ and } \  
\Gamma = Y \cup \Gamma' = \bigcup_{i=1}^3 \wt\Gamma_i.
\end{equation}
We see $\Gamma$ as a three-legged spider centered at $\ell$,
whose legs are the arcs $\wt \Gamma_i$. We want to be able to apply
the results of Section \ref{S8} to the $\wt\Gamma$, so let us check
that they are small Lipschitz graphs.

\begin{lem} \label{t10.3}
For $1 \leq i \leq 3$, $\wt\Gamma_i$ is a $8\lambda$-Lipschitz 
graph over $\rho(\ell,a_i^\ast)$.
\end{lem}

Of course the difference between $\lambda$ and $8\lambda$ will not prevent
us from applying Section \ref{S8}.
This looks like Lemma \ref{t10.2}, but we will need to worry a little
because we slightly change one endpoint and the orientation. 
Fix $i$; in addition to $P_i$ (the plane that contains $\rho(x_0,a_i^\ast)$),
$\pi_i$, and $\pi_i^\perp$ (see below \eqref{10.5}), we introduce the plane $\ol P_i$ 
that contains $\rho(\ell,a_i^\ast)$ and the corresponding projections $\ol\pi_i$
and $\ol\pi_i^\perp$. Notice that
\begin{equation} \label{10.31}
|v(x_0,a_i^\ast)-v(\ell,a_i^\ast)| \leq 10^{-3} \alpha^2
\end{equation}
by the proof of \eqref{10.22}; then
\begin{equation} \label{10.32}
|| \ol\pi_i^\perp - \ol\pi_i^\perp || = || \ol\pi_i - \ol\pi_i || \leq 10^{-2} \alpha^2.
\end{equation}
We first look outside of the disk $D_1 = \S \cap B(x_0,10r)$. There
$\wt\Gamma_i = \Gamma_i$, and the definition \eqref{7.14} says that
\begin{equation} \label{10.33}
|\pi^\perp(z) - \pi^\perp(z')| \leq \lambda |\pi(z) - \pi(z')| \leq \lambda |z-z'|
\end{equation}
for $z, z' \in \wt\Gamma_i \sm D_1$ and, 
(if $\alpha$ is small enough compared to $\lambda$), \eqref{10.32} yields
\begin{equation} \label{10.34}
|\ol\pi^\perp(z) - \ol\pi^\perp(z')| \leq 2\lambda |\pi(z) - \pi(z')|.
\end{equation}
Next we look inside $D_2 = \S \cap B(x_0,500r)$. 
Recall from \eqref{10.26} and \eqref{10.25} that 
\begin{equation} \label{10.35}
D_1 \subset B(x_0,\lambda/2) \subset B(\ell, \lambda).
\end{equation}
On $D_2$, we can use \eqref{10.12}, which says that
$\Gamma_i \cap D_2$ is a $3\lambda$-Lipschitz graph over (a part of) the
line through $v(x_0,a_i^\ast)$. By \eqref{10.31}, it is also a 
$4\lambda$-Lipschitz graph over the line through $v(\ell,a_i^\ast)$. 
But we modified it, and replaced the arc between $x_0$ and $z_i$
with the arc $\rho(\ell,z_i)$. Let $z$ be any point of $\rho(\ell,z_i)$ and
$v$ denote a tangent vector to $\rho(\ell,z_i)$ at $z$, oriented in the direction of $z_i$. 
Then
\begin{equation} \label{10.36}
|v-v(\ell,z_i)| \leq 2|\ell-z_i| \leq 2|\ell-x_0| + 2r \leq 2(1+\alpha^{-1}) |\ell-x_0|
\leq 4(1+\alpha^{-1}) \alpha^2 \tau < 10^{-1} \lambda
\end{equation}
because $z_i \in \d D$, by \eqref{10.25} and \eqref{10.3}, and if $\alpha$ is small enough.
Next
\begin{equation} \label{10.37}
|v(\ell,z_i)-v(x_0,z_i)| \leq 2 |\ell-x_0| |x_0-z_i|^{-1} = 2 |\ell-x_0| r^{-1} \leq 2\alpha
< 10^{-1} \lambda
\end{equation}
because $r = \alpha^{-1}|\ell-x_0|$ (by \eqref{10.25}), 
\begin{equation} \label{10.38}
|v(x_0,z_i)-v(x_0,a_i^\ast)| \leq 3 \lambda
\end{equation}
by the Lipschitz description \eqref{10.12}, and
\begin{equation} \label{10.39}
|v(x_0,a_i^\ast)-v(\ell,a_i^\ast)| \leq 2 |x_0- \ell| |\ell-a_i^\ast|^{-1}
\leq 2 |x_0- \ell | (5\eta(X))^{-1} \leq \alpha^2 \leq 10^{-1} \lambda
\end{equation}
because $|x_0- \ell | + |a_i - a_i^\ast| \leq 2 \alpha^2 \tau$ by \eqref{10.3},
and $|a_i-\ell| \geq 10 \eta(X) \geq 10^4 \tau$ by \eqref{3.11}, \eqref{3.12},
and our choice of $\tau$.

Altogether $|v-v(\ell,a_i^\ast)| \leq 4 \lambda$, $\rho(\ell,z_i)$ is a
$4\lambda$-Lipschitz graph over the line ${\rm Vect}(v(\ell,a_i^\ast))$,
and since we already know this about $\Gamma_i \cap D_2$, we 
also get that $\wt\Gamma_i \cap D_2$ is a $4\lambda$-Lipschitz graph over that line.
Now we can apply Lemma \ref{t10.2}, transposed for curves that start from $\ell$
(and, if we want to be precise, with a radius a little smaller than $\lambda$ 
for the analogue of $B_1$), and we get that $\wt\Gamma_i$ is a 
$8\lambda$-Lipschitz graph over $\rho(\ell,a_i^\ast)$.
\qed

\ms
Next we check \eqref{9.2}-\eqref{9.7} for $\Gamma$ (which we see as a union of
three curves $\wt\Gamma_i$). For \eqref{9.2}, we just need to know that
the three branches of $Y$ make large angles at $\ell$.
But if $v_i$ denotes the tangent direction of $\rho(\ell,z_i)$ at $\ell$, 
we know from \eqref{10.36}-\eqref{10.39} that $|v_i - v(\ell, a_i^\ast)| \leq 4\lambda$.
Since the $v(\ell,a_i)$ make $120^\circ$ angles, we see that
\begin{equation}\label{10.40}
\text{the three legs of $Y$ make angles larger than $110^\circ$ at $\ell$.}
\end{equation}

Next we need to check \eqref{9.3}, i.e. that each $\Gamma_i$ makes a large angle with the
other curves $\cL_k$ that arrive at $a_i^\ast$. The verification is the same as 
what we did below \eqref{9.12}.

\ms
Now we turn to the length estimates. 
First we want to compare $\H^1(Y)$ with $\H^1(\Gamma'')$, and to this effect
we shall differentiate 
\begin{equation} \label{10.41}
f(z) = \sum_{i=1}^3 \ddist(z,z_i)
\end{equation}
in the interior of $D$. First notice that for $1 \leq i \leq 3$, 
$\ddist(z,z_i)$ is differentiable on $\S \sm \{ z_i \}$, with
\begin{equation} \label{10.42}
\nabla_z \ddist(z,z_i) = -v(z,z_i).
\end{equation}
Thus
\begin{equation} \label{10.43}
|\nabla f(x_0)| = \big| \sum_{i=1}^3 v(x_0,z_i) \big|
= \big| \sum_{i=1}^3 [v(x_0,z_i)-v(\ell,a_i)] \big|
\end{equation}
because $\sum_i v(\ell,a_i) = 0$ (the three $\cC_i$ make $120^\circ$ angles). But
\begin{equation} \label{10.44}
|v(x_0,z_i)-v(\ell,a_i)| \leq |v(x_0,z_i)-v(x_0,a_i^\ast)| + |v(x_0,a_i^\ast)-v(\ell,a_i)| 
\leq 3 \lambda + 10^{-3} \alpha^2 \leq 4\lambda
\end{equation}
by \eqref{10.38} and \eqref{10.22}, so
\begin{equation} \label{10.45}
|\nabla f(x_0)| \leq 12\lambda.
\end{equation}
 Also, $v(z,z_i)$ is differentiable, with
$|\nabla_z v(z,z_i) | \leq |z-z_i|^{-1}$. For $z \in \rho(x_0,\ell)$,
\begin{equation} \label{10.46}
|z-x_0| \leq |\ell - x_0| = \alpha r \leq r/2
\end{equation}
by \eqref{10.25}, so $|\nabla_z v(z,z_i)| \leq 2 r^{-1}$ (because $|z_i-x_0| = r$). 
We sum over $i$, integrate on a part of $\rho(x_0,\ell)$, and get that
\begin{equation} \label{10.47}
|\nabla f(x_0) - \nabla f(z)| \leq 6 r^{-1} \ddist(x_0,\ell).
\end{equation}
Then we integrate again on $\rho(x_0,\ell)$ and get that
\begin{equation} \label{10.48}
\begin{aligned}
\H^1(Y) = f(\ell) 
&\leq f(x_0) + \ddist(x_0,\ell) |\nabla f(x_0)| + 6r^{-1} \ddist(x_0,\ell)^2
\cr&\leq f(x_0) + [12\lambda+ 6r^{-1} \ddist(x_0,\ell)] \ddist(x_0,\ell) 
\cr&\leq f(x_0) + [12\lambda+ 9\alpha] \ddist(x_0,\ell)
\leq f(x_0) + 13\lambda |x_0 -\ell|
\end{aligned}
\end{equation}
by \eqref{10.45} \eqref{10.47}\eqref{10.46}, and if $\alpha$ is small enough.

Notice that $f(x_0) \leq \H^1(\Gamma'')$, because $\Gamma''$ is composed of 
three essentially disjoint curves that go from $x_0$ to the $z_j$ 
(see near \eqref{10.27}). Hence
\begin{equation} \label{10.49}
\H^1(Y) \leq f(x_0) + 13\lambda |x_0 -\ell| \leq \H^1(\Gamma'') +  13\lambda |x_0 -\ell|.
\end{equation}
We add the missing piece $\Gamma'$, and get
\begin{equation} \label{10.50}
\H^1(\Gamma) = \H^1(\Gamma') + \H^1(Y) 
\leq \H^1(\Gamma') + \H^1(\Gamma'') +  13\lambda |x_0 -\ell|
= \H^1(\Gamma^\ast) + 13\lambda |x_0 -\ell|
\end{equation}
by \eqref{10.30}, \eqref{10.27}, \eqref{10.28}, and \eqref{10.49}.
So
\begin{equation} \label{10.51}
\H^1(\ol\gamma) - \H^1(\Gamma) 
\geq \H^1(\ol\gamma) - \H^1(\Gamma^\ast) - 13\lambda |x_0 -\ell|.
\end{equation}
Recall from \eqref{10.1} and \eqref{10.2} that
\begin{equation} \label{10.52}
\H^1(\ol\gamma) - \H^1(\Gamma^\ast) 
= \H^1(\gamma_\ell) + \sum_{i=1}^3 [\H^1(\ol \gamma_i)-\H^1(\Gamma_i)]
\geq 32\lambda |x_0 -\ell|  
\end{equation}
because $\gamma_\ell$ and the $\ol \gamma_i$ are disjoint, and
the $\Gamma_j$ are disjoint, and then by the defining condition \eqref{10.24}.
We also have that
\begin{equation} \label{10.53}
\H^1(\ol \gamma_i)-\H^1(\Gamma_i) \geq 0
\end{equation}
by \eqref{7.16}, so 
\begin{equation} \label{10.54}
\H^1(\ol\gamma) - \H^1(\Gamma^\ast) \geq \H^1(\gamma_\ell),
\end{equation}
and now \eqref{10.51}, \eqref{10.53}, and \eqref{10.54} yield
\begin{eqnarray} \label{10.55}
\H^1(\ol\gamma) - \H^1(\Gamma) 
&\geq& \frac{1}{4}[\H^1(\ol\gamma) - \H^1(\Gamma^\ast)] 
+ \frac{3}{4}[\H^1(\ol\gamma) - \H^1(\Gamma^\ast)] -13\lambda |x_0 -\ell| 
\nn\\
&\geq& \frac{1}{4}[\H^1(\ol\gamma) - \H^1(\Gamma^\ast)]
+ \frac{1}{8} \H^1(\gamma_\ell) + \frac{5}{8}[\H^1(\ol\gamma) - \H^1(\Gamma^\ast)]
- 13\lambda |x_0 -\ell| \big)
\nn\\
&\geq& \frac{1}{4}[\H^1(\ol\gamma) - \H^1(\Gamma^\ast)] 
+ \frac{1}{8} \H^1(\gamma_\ell) + 7\lambda |x_0-\ell|. 
\end{eqnarray}
The three terms are nonnegative, so \eqref{10.55} is stronger than \eqref{9.6}.

\ms
We are left with \eqref{9.7} to check, i.e.,
\begin{equation} \label{10.9.7}
\H^1(\Delta(\overline\gamma ,\Gamma)) \leq C [\H^1(\overline\gamma)- \H^1(\Gamma)]
+ C [\H^1(\Gamma)-\H^1(\rho)],
\end{equation}
where $\ol \gamma$ is the extended version of our initial net $\gamma$ (see \eqref{9.1}), 
which is an extended spider with a short tail, and $\rho$ will be discussed soon. We write
\begin{equation} \label{10.56}
\Delta(\ol\gamma,\Gamma) 
\subset \Delta(\ol\gamma,\wt\gamma) 
\cup \Delta(\wt\gamma,\Gamma^\ast) 
\cup \Delta(\Gamma^\ast,\Gamma)
\end{equation}
(all symmetric differences), where $\wt\gamma = \ol\gamma \sm \gamma_\ell$
is the spider without its tail (see \eqref{10.1}), $\Gamma^\ast$ is our initial Lipschitz
spider with the same ends, and $\Gamma$ is our final pick (centered at $\ell$).

We start with $\Delta(\ol\gamma,\wt\gamma)
= \ol\gamma \sm\wt\gamma = \gamma_\ell$. Recall from \eqref{10.55} 
that $\H^1(\gamma_\ell) \leq 8[\H^1(\ol\gamma) - \H^1(\Gamma)]$,
which is dominated by the right-hand side of \eqref{10.9.7}, so this term
is all right.

Next we consider $\Delta(\Gamma^\ast,\Gamma)$.
Since $\Gamma^\ast$ and $\Gamma$ have the three exterior curves
$\Gamma'_i$ in common, we are left with $Y$ for $\Gamma$,
and $\Gamma''$ for $\Gamma^\ast$ (see \eqref{10.27}-\eqref{10.30}); 
thus
\begin{equation} \label{10.57}
\H^1(\Delta(\Gamma^\ast,\Gamma)) 
\leq \H^1(Y) + \H^1(\Gamma'') 
\leq 7 r + \H^1(\Gamma'').
\end{equation}
The simplest case is when $\H^1(\Gamma'') \leq 14r$, say. Then
\begin{equation} \label{10.58}
\H^1(\Delta(\Gamma^\ast,\Gamma)) \leq 21r 
= 21 \alpha^{-1} |x_0-\ell|
\leq 3 \lambda^{-1}\alpha^{-1} [\H^1(\ol\gamma) - \H^1(\Gamma)]
\end{equation}
by \eqref{10.25} and \eqref{10.55}, which again is enough for \eqref{9.7}.

If instead $\H^1(\Gamma'') > 14r$, we can revise some of our earlier
pessimistic estimates, because 
\begin{equation} \label{10.59}
\H^1(Y) \leq 7r \leq \frac{1}{2} \H^1(\Gamma'')
\end{equation}
and then, after adding the exterior part $\H^1(\Gamma')$ to both sides,
\begin{equation} \label{10.60}
\H^1(\Gamma) \leq \H^1(\Gamma^\ast) 
- \frac{1}{2} \H^1(\Gamma'').
\end{equation}
This is better than what we had before (see \eqref{10.50}); it implies that
\begin{equation} \label{10.61}
\H^1(\ol\gamma) - \H^1(\Gamma) \geq \H^1(\ol\gamma) - \H^1(\Gamma^\ast) 
+ \frac{1}{2} \H^1(\Gamma'') .
\end{equation}
We forget the term $\H^1(\ol\gamma) - \H^1(\Gamma^\ast)$, which
is nonnegative by \eqref{10.52}, and get that
$\H^1(\Gamma'') \leq 2[\H^1(\ol\gamma) - \H^1(\Gamma)]$; finally
\begin{equation} \label{10.62}
\H^1(\Delta(\Gamma^\ast,\Gamma)) 
\leq 7 r + \H^1(\Gamma'') \leq 2\lambda^{-1}\alpha^{-1} [\H^1(\ol\gamma) - \H^1(\Gamma)]
\end{equation}
by \eqref{10.57} and the second part of \eqref{10.58}.

We are left with the middle term $\Delta(\wt\gamma,\Gamma^\ast)$ from \eqref{10.56}.
Recall that $\Gamma^\ast$ was obtained by applying the construction of Section \ref{S7}
to the three curves $\ol \gamma_i$ that compose 
$\wt\gamma = \ol\gamma \sm \gamma_\ell$. 
Thus by \eqref{7.16} and \eqref{7.5}
\begin{equation} \label{10.63}
\H^1(\Gamma^\ast \sm \wt\gamma) \leq \H^1(\wt\gamma \sm \Gamma^\ast)
\leq C(\lambda)[\H^1(\Gamma^\ast) - \H^1(\wt \rho)],
\end{equation}
where $\wt\rho$ is the union of the three geodesics $\rho(x_0,a_i^\ast)$
that we used to construct $\Gamma^\ast$.
Hence 
\begin{equation} \label{10.63aa}
\H^1(\Delta(\wt\gamma,\Gamma^\ast)) \leq C(\lambda)[\H^1(\Gamma^\ast) - \H^1(\wt \rho)],
\end{equation}
which we just need to bound by the right-hand side of \eqref{9.7} and \eqref{10.9.7}.
We already did this for
\begin{equation} \label{10.63ab}
\H^1(\Gamma^\ast) - \H^1(\Gamma) \leq \H^1(\Delta(\Gamma^\ast,\Gamma))
\end{equation}
(by \eqref{10.62}), and since $\H^1(\Gamma) - \H^1(\rho)$ is a part of 
the right-hand side of \eqref{9.7} and \eqref{10.9.7} (and the other one is nonnegative
by \eqref{10.51} and \eqref{10.52}), we just need to control
$\H^1(\rho)-\H^1(\wt \rho)$. Recall that $\wt \rho$ is the union of the geodesics 
$\rho(x_0,a_i^\ast)$ with the same endpoints as the arcs of the spider $\Gamma^\ast$,
while $\rho$ is the union of the geodesics $\rho(\ell,a_i^\ast)$ that correspond to the 
decomposition of the spider $\Gamma$ that we want to use. Since
$\H^1(\rho(\ell,a_i^\ast)) \leq \H^1(\rho(x_0,a_i^\ast))+\ddist(\ell,x_0)$
(a brutal estimate), we see that
\begin{equation} \label{10.63ac}
\H^1(\rho)-\H^1(\wt \rho) \leq 3 \ddist(\ell,x_0) \leq 4 |x_0-\ell| \leq 4\alpha r
\leq 12 \lambda^{-1} [\H^1(\ol\gamma) - \H^1(\Gamma)]
\end{equation}
by \eqref{10.25} and \eqref{10.58} or \eqref{10.62}. This completes our proof of 
\eqref{10.9.7} and \eqref{9.7} and the verification of \eqref{9.2}-\eqref{9.7} in Case A;
we may now turn to the next case.

\subsection{Case B : consequences of the definition on the geometry of $\ol\gamma$}

\begin{figure}[!h]  
\centering
\includegraphics[width=7cm]{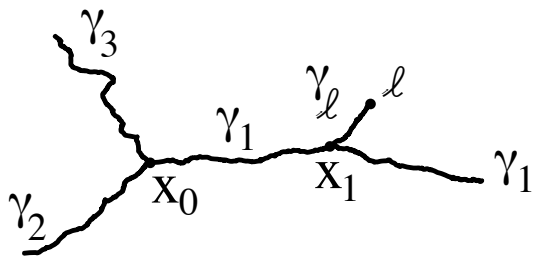}
\caption{The initial setting for Case B
\label{f10.a}}
\end{figure}

Since we are happy in Case A, we shall now assume that its defining condition \eqref{10.24}
fails, i.e. that
\begin{equation}\label{10.65} 
\H^1(\gamma_\ell) + \sum_{i=1}^3 [\H^1(\ol\gamma_i) - \H^1(\Gamma_i)]  
< 32 \lambda |x_0 - \ell|. 
\end{equation}
We shall call this \ub{Case B}; see Figure \ref{f10.a} % 10a
First notice that \eqref{10.65} implies that
\begin{equation} \label{10.66}
|x_1 - \ell| \leq \H^1(\gamma_\ell) \leq 32\lambda |x_0 - \ell|, 
\end{equation}
where $x_1$ is the point where $\gamma_\ell$ is attached to $\gamma$.
Recall that $\lambda$ is small, so $x_1$ lies relatively far from $x_0$ (compared to $\ell$). 
Without loss of generality, we can assume that
$x_1 \in \gamma_1$. The next lemma says that $x_1$ lies in the expected direction
(seen from $x_0$).

 %p89

\begin{lem} \label{t10.4} The direction of $\rho(x_0, x_1)$ at $x_0$ is such that
\begin{equation} \label{10.67}
|v(x_0,x_1) - v(x_0,a_1^\ast)| \leq 30 \sqrt\lambda.
\end{equation}
\end{lem}

\begin{figure}[!h]  
\centering
\includegraphics[width=9cm]{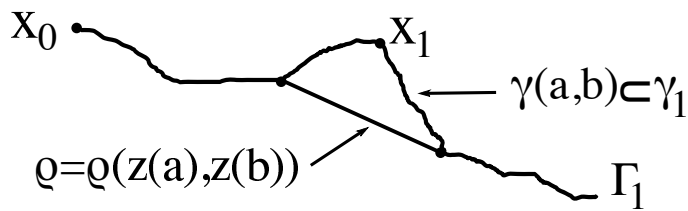} % Figure10-b.eps marche aussi
\caption{The point $x_1$ lies close to $\rho$ because of (11.72) 
\label{f10.b}} %%% attention reference a la main en code c'est \eqref{10.68}
\end{figure}

\ms
Here again, 30 is what we get from the proof, but a larger number would still be fine.
Suppose not. Since $\Gamma_1$ is a $\lambda$-Lipschitz graph over $\rho(x_0,a_1^\ast)$
(and by \eqref{10.12} to make things simpler), $x_1 \in \gamma_1 \sm \Gamma_1$.
Recall how $\Gamma_1$ was constructed. We started from a parameterization $z : I \to \S$
of $\gamma_1$, selected a certain number of intervals $I_j$, and replaced $\gamma_1$
on $I_j$ by the constant speed parameterization of the geodesic $\rho_j$ with the same endpoints.
See below \eqref{7.10}. Here $x_1 \in \gamma_1 \sm \Gamma_1$, so the parameter $t$
such that $x_1 = z(t)$ lies in some $I_j$. Write $I_j = I$, and denote by $a$ and $b$ its endpoints.
We choose the names so that $z(a)$ lies between $x_0$ and $x_1$ on $\gamma_1$,
and hence $z(b)$ lies between $x_1$ and $a_i^\ast$. Also call
$\gamma(a,b)$ the portion of $\gamma_1$ between $z(a)$ and $z(b)$.
See Figure \ref{f10.b} already.  

Thus we replaced $\gamma(a,b)$ with the geodesic $\rho = \rho(z(a),z(b))$ 
in the construction of $\Gamma_1$. There was a similar replacement of other arcs of 
$\gamma_1$ on other intervals, and of course each time the length of the geodesic 
was no longer than the length of the arc of $\gamma$ it replaced. 
See near \eqref{7.15}. Because of this
\begin{equation} \label{10.68}
\H^1(\gamma(a,b)) - \H^1(\rho) 
\leq \H^1(\gamma_1) - \H^1(\Gamma_1)
\leq \sum_i [\H^1(\gamma_i) - \H^1(\Gamma_i)]
\leq 32\lambda |x_0 - \ell| 
\end{equation}
because the three numbers $\H^1(\gamma_i) - \H^1(\Gamma_i)$ are nonnegative 
(by \eqref{7.15}) and by \eqref{10.65}. We shall now complete this with a lower bound for 
$\H^1(\gamma(a,b)) - \H^1(\rho)$ which yields the desired contradiction. The computations
that follow seem shockingly long to the author, who feels compelled to do them but
hopes the pictures will be convincing enough.

%p89

Set $d = |x_0 - x_1|$. Also let $p$ denote the point of $\rho$ that minimizes the 
distance to $x_1$; notice that since $p \in \rho \subset \Gamma_1$ and $\Gamma_1$
is a small Lipschitz graph, 
\begin{equation} \label{10.69}
|v(x_0,p) - v(x_0,a_1^\ast)| \leq 2\lambda.
\end{equation}
Then $|v(x_0,p) - v(x_0,x_1)| \geq 30 \sqrt\lambda- 2\lambda \geq 29\sqrt\lambda$
because we assumed \eqref{10.67} to fail and if $\lambda$ is small enough; hence
\begin{equation} \label{10.70}
|p-x_1| \geq 29 \sqrt\lambda d.
\end{equation}
We shall first assume that
\begin{equation} \label{10.71}
\text{$p$ lies in the interior of $\rho$;} 
\end{equation}
this is supposed to be the main case, and pictures will be easier to draw.
In particular observe that the geodesic $\rho(x_1,p)$ is perpendicular to $\rho$ at $p$.

We will feel better if we know that all the geometric arguments that follow happen
in a tiny ball, so let us check that
\begin{equation} \label{10.72}
x_1, \, z(a),  \text{ and } p \ \text{ all lie in } \ol B(x_0,4d) \subset B(x_0,12\alpha^2 \tau)
\end{equation} 
(that is, a very small ball where we can expect curvature to play almost no role). 

\begin{figure}[!h]  
\centering
\includegraphics[width=9.cm]{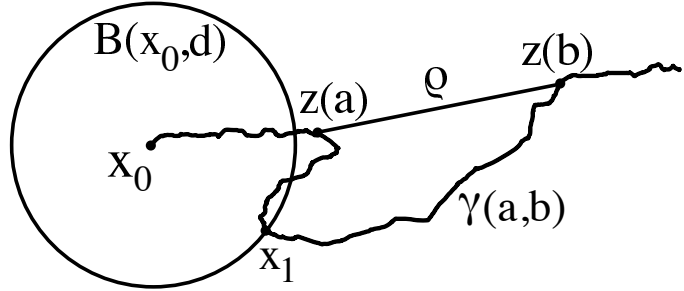}
\caption{ In this case already (and more if $|z(a)-x_0| \geq 2d$), $p=z(a)$.
\label{f10.1}}
\end{figure}

First we claim that $|z(a)-x_0| \leq 2d$.
Suppose not, and recall that $\rho \subset \Gamma_1$, 
$\Gamma_1$ is a small Lipschitz graph over $\rho(x_0,a_1^\ast)$, and
$\rho$ leaves from $z(a)$ in the direction opposite to $x_0$. Then, as suggested by
Figure \ref{f10.1} (and $z(a)$ should even lie further on the right), 
$z(a)$ should be the point of $\rho$ that lies closest to $x_1$, a 
contradiction with \eqref{10.71}. So $|z(a)-x_0| \leq 2d$.

Then $|p-x_0| \leq 4d$, because otherwise 
$|p-x_1| \geq |p-x_0|-|x_0-x_1| > 3d \geq |z(a)-x_0|+|x_0-x_1| \geq |z(a)-x_1|$. 

 %p90

Finally the second part of \eqref{10.72} holds because 
$d \leq |x_0-\ell| + |x_1-\ell| \leq (1+32\lambda) |x_0-\ell| < 3\alpha^2 \tau$ 
by \eqref{10.66} and Lemma \ref{t10.1}; so \eqref{10.72} holds.

We will need lower bounds for 
\begin{equation} \label{10.73}
\delta_a = \ddist(x_1,z(a)) - \ddist(p, z(a)) 
\ \text{ and  } \delta_b = \ddist(x_1,z(b)) - \ddist(p, z(b)).
\end{equation}
Let us first consider $\delta_b$. 
Recall that $\rho(x_1,p)$ is perpendicular to $\rho$ at $p$;
we can choose orthonormal coordinates of $\R^n$ where $p$, $z(b)$, and $x_1$ 
lie in $\R^3$ (so that we won't even need to write all the other coordinates), and
\begin{equation} \label{10.74}
p=(1,0,0), z(b) = (\cos s_b, \sin s_b, 0), \text{ and } x_1 = (\cos t, 0, \sin t),
\end{equation}
where in fact we can take $s_b = \ddist(z(b),p)$ and $t = \ddist(x_1,p)$,
maybe at the price of changing the orientation. Then
\begin{equation} \label{10.75}
|x_1 - z(b)|^2 = |\cos s_b - \cos t|^2 + \sin^2 s_b + \sin^2 t
= 2 - 2 \cos s_b \cos t.
\end{equation}
Here we shall not try to win much, because $z(b)$ may be quite far and then 
$\delta_b$ small. Let us just observe that 
\begin{equation} \label{10.76}
\begin{aligned}
s_b &= \ddist(z(b),p) \leq \ddist(a_1^\ast,p) 
\cr&\leq \H^1(\cC_i) + \ddist(a_1,a_1^\ast) 
+\ddist(\ell, x_0)+\ddist(x_0,p)
\leq \frac{\pi}{2}+25\alpha^2 \tau
\end{aligned}
\end{equation}
because $a_1^\ast$ is clearly the furthest point of $\Gamma$ from $p$,
because $\H^1(\cC_i) \leq \frac{\pi}{2}$ (see the construction above \eqref{3.1}; 
if we had forced the lengths of the $\cC_i$ to be a little shorter, we would have a 
slightly better estimate now, but $\frac{\pi}{2}$ looks natural), 
then by Lemma \ref{t10.1} and \eqref{10.72}.
Notice that setting $t=0$ in \eqref{10.75} corresponds to $p=x_1$ and a computation
of $|p - z(b)|^2$. If $s_b \leq \frac{\pi}{2}$, then $\cos s_b \geq 0$ and 
$|x_1 - z(b)|^2 \geq |p - z(b)|^2$ by \eqref{10.75}
(recall that $t = \ddist(x_1,p) \leq 20\alpha^2 \tau$ by \eqref{10.72}, so $\cos t > 0$). 
Then $\delta_b \geq 0$. Otherwise, even though $\cos s_b < 0$,
$-2\cos s_b \leq 50\alpha^2 \tau$ by \eqref{10.76} so
\begin{equation} \label{10.77}
 |p - z(b)|^2 - |x_1 - z(b)|^2 = (-2 \cos s_b)(1-\cos t) \in [0, 50\alpha^2 \tau t^2]
\end{equation}
by \eqref{10.76} (again recall that $t \leq 20\alpha^2 \tau$).
In this region where the distances $\ddist(p,z(b))$ and $\ddist(x_1,z(b))$
are very close to $\pi/2$, we can recover them in a $3$-Lipschitz way from 
$|p - z(b)|^2$ and $|x_1 - z(b)|^2$; we then deduce from \eqref{10.77} that 
\begin{equation} \label{10.78}
\delta_b \geq - 150 \alpha^2 \tau t^2 \geq - 10^4 \alpha^2 \tau d^2
\end{equation}
because $t = \ddist(x_1,p) \leq 9d$ by \eqref{10.72}.

Next we estimate $\delta_a$. The same computation as for \eqref{10.75} yields
\begin{equation} \label{10.79}
|x_1 - z(a)|^2 = |\cos s_a - \cos t|^2 + \sin^2 s_a + \sin^2 t
= 2 - 2 \cos s_a \cos t,
\end{equation}
with $s_a = \ddist(p, z(a))$. Now both $s_a$ and $t$ are small, 
by \eqref{10.72}, so there is no sign issue, and comparing with $t=0$ yields
\begin{equation} \label{10.80}
|x_1 - z(a)|^2 - |p - z(a)|^2 = 2 \cos s_a (1-\cos t) \geq t^2 \cos s_a
\geq \frac{28 t^2}{29}. 
\end{equation}
But \eqref{10.70} implies that $t = \ddist(x_1,p) \geq 29 \sqrt\lambda d$,
so \eqref{10.80} says that
\begin{equation} \label{10.81}
|x_1 - z(a)|^2 - |p - z(a)|^2 \geq 28 \cdot 29 \lambda d^2 = 812 \lambda d^2.
\end{equation}
Set $\alpha_0 = \frac{1}{2} \ddist(z(a),p)$ and $\alpha_1 = \frac{1}{2} \ddist(x_1,p)$;
thus $\alpha_0 \leq \alpha_1$ by \eqref{10.80},
$|p - z(a)| = 2\sin \alpha_0$ and $|x_1 - z(a)| = 2\sin \alpha_1$, and, by the fundamental
theorem of calculus,
\begin{equation} \label{10.82}
|x_1 - z(a)|^2 - |p - z(a)|^2 = 2[\sin^2 \alpha_1-\sin^2 \alpha_0] 
= 4 (\alpha_1-\alpha_0) \sin\alpha \cos\alpha 
\end{equation}
for some $\alpha\in [\alpha_0,\alpha_1]$. 
In addition $\alpha_1 \leq 5d$ because \eqref{10.72} says that $|x_1-z(a)| \leq 8d$, so 
$2\sin \alpha \cos \alpha = \sin 2\alpha \leq 2\alpha_1 \leq 10d$, and
\begin{equation} \label{10.83}
\begin{aligned}
\delta_a &= 2 (\alpha_1-\alpha_0) 
= (2 \sin \alpha \cos \alpha)^{-1} \Big(|x_1 - z(a)|^2 - |p - z(a)|^2\Big)
\cr&
\geq (10d)^{-1} \big(|x_1 - z(a)|^2 - |p - z(a)|^2\big)
\geq \frac{812 \lambda d^2}{10d} \geq 81 \lambda d 
\end{aligned}
\end{equation}
by \eqref{10.73} and \eqref{10.81}. On the other hand,
\begin{equation} \label{10.84}
\begin{aligned}
\delta_a+\delta_b &= \ddist(x_1,z(a)) + \ddist(x_1,z(b)) - \ddist(p, z(a)) - \ddist(p, z(b))
\cr&= \ddist(x_1,z(a)) + \ddist(x_1, z(b)) - \ddist(z(a),z(b))
\cr& 
\leq \H^1(\gamma(a,b)) - \H^1(\rho) 
\leq 32\lambda |x_0 - \ell| \leq 32 (1-32\lambda)^{-1} \lambda d 
\leq 33 \lambda d 
\end{aligned}
\end{equation}
because $p$ lies between $z(a)$ and $z(b)$ on the geodesic $\rho = \rho(z(a),z(b))$,
then because $\gamma(a,b)$ goes from $z(a)$ to $x_1$ to $z(b)$, and finally by
\eqref{10.68} and because \eqref{10.66} says that $d = |x_0-x_1| \geq (1-32\lambda)|x_0-\ell|$. 
We get the desired contradiction by comparing this to \eqref{10.83} and \eqref{10.78}.

\ms
We are not quite finished yet, because we still need to deal with the case when 
\eqref{10.71} fails, i.e, when $p = z(a)$ or $z(b)$.
Suppose first that $p = z(a)$, i.e., $\dist(x_1,\rho) = |x_1-z(a)|$
(as in Figure \ref{f10.1}). %
We claim that
\begin{equation} \label{10.85}
\dist(x_1, z(b)) \geq \dist(z(a), z(b)) - 10^4 \alpha^2 \tau d^2.
\end{equation}
Indeed, let $H$ be the vector hyperplane through $z(a)$ and perpendicular to $\rho$.
Since $\dist(x_1,\rho) = |x_1-z(a)|$, $x_1$ lies on $H$, or on the other side of
$H$ as $z(b)$. Call $\xi$ the intersection of $H$ with the geodesic $\rho(x_1, z(b))$.
Also denote by $P$ the plane that contains $\rho$. It is easy to see that 
$\dist(\xi,P) \leq \dist(x_1,P)$ (the geodesic $\rho(x_1, z(b))$ goes through $\xi$ and ends on $\rho$);
then $\dist(\xi,P) \leq |x_1 - z(a)|$ (because $z(a) \in P$). Then, since $\xi \in H$
and the geodesic distance is a monotone function of the Euclidean distance,
$\ddist(\xi,z(a)) = \ddist(\xi,P) \leq \ddist(x_1,z(a))$. The proof of \eqref{10.78} 
applies to $\xi$ too, because of the orthogonality that comes from the fact that $\xi \in H$,
and we get that $\dist(\xi, z(b)) \geq \dist(z(a), z(b)) - 10^4 \alpha^2 \tau d^2$, and 
\eqref{10.85} follows because $\dist(x_1, z(b))$ is at least as large.
Then
\begin{eqnarray} \label{10.86}
\H^1(\gamma(a,b)) &\geq& \ddist(z(a),x_1) + \ddist(x_1,z(b))
\geq 29\sqrt\lambda d + \ddist(x_1,z(b))
\nn\\
&\geq& 29\sqrt\lambda d + \ddist(z(a), z(b)) - 10^5 \alpha^2 \tau d^2
= 29\sqrt\lambda d + \H^1(\rho) - 10^5 \alpha^2 \tau d^2
\end{eqnarray}
because $\gamma(a,b)$ goes from $z(a)$ to $x_1$ to $z(b)$, by \eqref{10.70}
and because $p=z(a)$ and $\rho = \rho(z(a),z(b))$, then by \eqref{10.85} and because 
$\ddist$ is a $10$-Lipschitz function of $\dist$ in the current range. 
This is not compatible with \eqref{10.68} (with the same sort of verification as above).

We are left with the case when $p = z(b)$. We claim that 
\begin{equation} \label{10.87}
\dist(x_1, z(a)) \geq \dist(z(a), z(b)) - 10^4 \alpha^2 \tau d^2.
\end{equation}
We prove this as with \eqref{10.85}, but with $z(a)$ and $z(b)$ exchanged.
We arrive to the information that $\ddist(\xi,z(b)) = \ddist(\xi,P) \leq \ddist(x_1,z(b))$.
We may now follow the same argument above, with the proof of \eqref{10.78}, and get \eqref{10.87},
but we may also observe that if $p = z(b)$, then $z(b)$ is not far from
$x_0$ and $x_1$, hence the proof of \eqref{10.83} also allow us to get rid of the ugly term 
$-10^4 \alpha^2 \tau d^2$. Anyway, may conclude as in \eqref{10.86}, and get the desired contradiction.
This completes our proof of Lemma \ref{t10.4}.
\qed

\subsection{Construction of $\Gamma$ in Case B}
\begin{figure}[!h]  %% p94
\centering 
\includegraphics[width=10.cm]{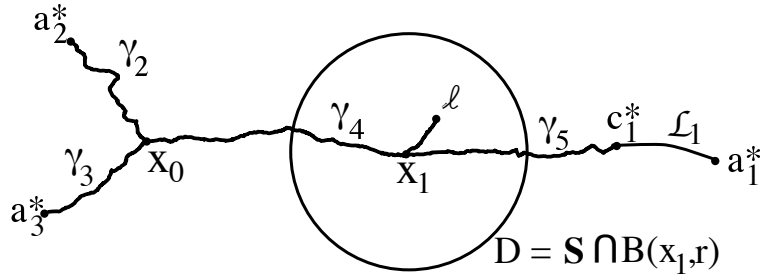}
\caption{ The $\gamma_j$ near $D$.
\label{f10.2}}
\end{figure}

\ms
We stay in \ub{Case B} (defined by \eqref{10.65}), and now we build the net $\Gamma$.
The general principle will be the same as in Case A, where we forced $\Gamma^\ast$
to make a small detour through $\ell$, but now everything will happen near $x_1$ and 
(beause of \eqref{10.66}) relatively far from $x_0$.
Set
\begin{equation} \label{10.88}
r = (100\lambda)^{-1} |x_1 - \ell|,  
D = \S \cap B(x_1,r), \text{ and } \d D = \S \cap \d B(x_1,r).
\end{equation}
We choose this radius because this way,
\begin{equation} \label{10.88b}
100 \lambda r =  |x_1 - \ell| \leq 32 \lambda |x_0 - \ell|
\end{equation}
by \eqref{10.66}, hence
\begin{equation} \label{10.88c}
|x_0 - x_1| \geq |x_0 - \ell|-|\ell - x_1| \geq (1-32\lambda) |x_0 - \ell|
\geq \frac{(1-32\lambda)}{32\lambda} |x_1 - \ell| \geq 3r
\end{equation}
and our construction will not involve $\Gamma_2$ or $\Gamma_3$. 
Indeed, not only does $|x_0-x_1| > 2r$, but 
Lemma \ref{t10.4}, the fact that $\Gamma_i$ is a small Lipschitz graph
over $\rho(x_0,a_i^\ast)$, and \eqref{10.23}, imply that
\begin{equation} \label{10.89}
\dist(\Gamma_2\cup\Gamma_3, D) > 2r.
\end{equation}
So we leave $\Gamma_2$ and $\Gamma_3$ alone, but we change $\Gamma_1$.
Denote by $\gamma_4$ the arc of $\gamma_1$ between $x_0$ and $x_1$,
and by $\gamma_5$ the rest of $\gamma_1$, i.e., between $x_1$ and $a_1^\ast$.
See Figure \ref{f10.2}. %
Then set $\ol\gamma_4 = \gamma_4$ and $\ol\gamma_5 = \gamma_5 \cup \cL_1$. 
Apply the construction of Section \ref{S7} to $\gamma_4$ and $\ol\gamma_5$; 
this gives two small Lipschitz graphs $\Gamma_4$ and $\Gamma_5$. 
Now set $\ol\gamma = \cup_{i=2}^5 \ol\gamma_i$ and the analogue of $\Gamma^\ast$
is $\cup_{i=2}^5 \Gamma_i$.

\begin{figure}[!h]  %% p94
\centering 
\includegraphics[width=10.cm]{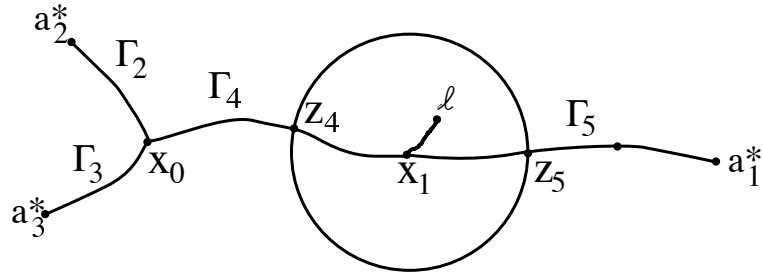}
\caption{ The corresponding $\Gamma_j$.
\label{f11.3}}
\end{figure} 

\begin{figure}[!h]  %% p94
\centering 
\includegraphics[width=10.cm]{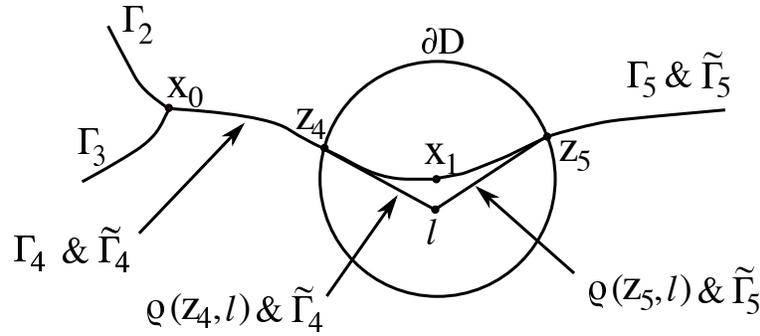}
\caption{ Our choice of $\Gamma$ (this time, with $\ell$ below).
\label{f11.4}}
\end{figure}

Since $\Gamma_4$ and $\Gamma_5$ are small Lipschitz graphs starting from
$x_1$, they meet $\d D$ exactly once, at points which we call $z_4$ and $z_5$
(see Figure \ref{f11.3}). %
Denote by $\Gamma'_4$ the arc of $\Gamma_4$ between $x_0$ and $z_4$,
and set $\wt\Gamma_4 = \Gamma'_4 \cup \rho(z_4,\ell)$. Similarly,
denote by $\Gamma'_5$ the arc of $\Gamma_5$ between $z_5$ and $a_1^\ast$,
and set $\wt\Gamma_5 = \rho(\ell,z_5) \cup \Gamma'_5$. 
We shall take
\begin{equation} \label{10.90}
\Gamma = \Gamma_2 \cup \Gamma_3 \cup \wt\Gamma_4 \cup \wt\Gamma_5,
\end{equation}
which is now composed of $6$ pieces, but looks a lot like a $3$-legged spider with
a small detour organized along one of its legs. See Figure \ref{f11.4}. %

Our next task consists in checking that this description is right, and that the angles are
large enough.

\begin{lem} \label{t10.5}
The curve $\wt\Gamma_4$ is a $10^3\lambda$-Lipschitz graph
over $\rho(x_0,\ell)$, and $\wt\Gamma_5$ is a $10^3\lambda$-Lipschitz graph
over $\rho(\ell,a_1^\ast)$. In addition,
\begin{equation} \label{10.91}
\text{$\wt\Gamma_4$ and $\wt\Gamma_5$ make an angle larger than 
$\pi - 40\sqrt\lambda$ at $\ell$.} 
\end{equation}
\end{lem}

\ms
Let us first consider $\Gamma_4$. A small advantage of the situation is that
by \eqref{10.3} and \eqref{10.66}, 
\begin{equation} \label{10.92}
\text{$x_1$ and $\ell$ both lie in the very small ball $B(x_0,3\alpha^2\tau)$,}
\end{equation}
so the Lipschitz geometry will be simpler.
Set $e = v(x_0,x_1)$; by the proof of \eqref{10.12}, 
\begin{equation} \label{10.93}
\Gamma_4 \text{ is a $3\lambda$-Lipschitz graph over } {\rm Vect}(e).
\end{equation}
As for the geodesic piece, notice that
\begin{equation} \label{10.94}
|v(z_4, x_1) - v(z_4,\ell)| \leq 2  |z_4-x_1|^{-1} |x_1-\ell| 
= 2 r^{-1} |x_1 - \ell| \leq 200 \lambda
\end{equation}
by \eqref{10.88}. Since $z_4$ and $x_1$ lie on $\Gamma_4$, \eqref{10.93} 
also says that $|v(z_4, x_1) - e| \leq 4\lambda$, and so $|v(z_4,\ell)-e| \leq 204 \lambda$.
In addition, the geodesic $\rho(z_4,\ell)$ is too short to turn much: if $v$ denotes
a tangent direction to $\rho(z_4,\ell)$ (oriented in the direction of $\ell$), then
$|v-v(z_4,\ell)| \leq \ddist(z_4,\ell) \leq 2 |z_4-\ell|
\leq 2r + 2|x_1-\ell| \leq 3r \leq 2 |x_0-\ell| \leq 4 \alpha^2 \tau$
by \eqref{10.88b}, \eqref{10.88c}, and \eqref{10.3}.

Altogether $|v-e| \leq 208\lambda$, and $\wt\Gamma_4$ is a 
$208 \lambda$-Lipschitz graph over ${\rm Vect}(e)$.
Now the easy part of the proof of Lemma \ref{t10.2} says that $\wt\Gamma_4$
is also a $10^3 \lambda$-Lipschitz graph over $\rho(x_0,\ell)$.

Next consider $\wt\Gamma_5$. 
Here the proof is very similar to what we did for Lemma \ref{t10.3}, 
so we shall skip some details. We first control $\wt\Gamma_5$ outside
of $D_1 = \S \cap B(x_1,10r)$, and for this we just copy the proof of 
Lemma \ref{t10.3}, with $x_0$ replaced with $x_1$, up to \eqref{10.34} included.
Then we look inside $D_2 = \S \cap B(x_1,500r)$. We continue as before,
but modify slightly the angle estimates. We start with \eqref{10.36}; instead we say that 
when $v$ is a tangent direction to $\rho(\ell,z_5)$, 
\begin{equation} \label{10.95}
|v-v(\ell,z_5)| \leq 2|\ell-z_5| \leq 2|\ell - x_1| +2r \leq (40\lambda)^{-1} |\ell - x_1|
\leq \lambda^{-1} \alpha^2 \tau \leq \lambda
\end{equation}
again because the geodesic is too short to turn, and by \eqref{10.88} and \eqref{10.92}.
Then, instead of \eqref{10.37} and as in \eqref{10.94},
\begin{equation} \label{10.96}
|v(\ell,z_5)-v(x_1,z_5)| \leq 2 |\ell-x_1| |x_1-z_5|^{-1}
= 2 r^{-1} |x_1 - \ell|  \leq 200 \lambda,
\end{equation}
and (as in \eqref{10.38})
\begin{equation} \label{10.97}
|v(x_1,z_5)-v(x_1,a_1^\ast)| \leq 3 \lambda
\end{equation}
because $z_5$ lies in the small Lipschitz graph $\Gamma_5$
over $\rho(x_1,a_1^\ast)$ (and by the analogue of \eqref{10.12}).
Finally, 
\begin{equation} \label{10.98}
|v(x_1,a_1^\ast)-v(\ell,a_1^\ast)| \leq 2 |x_1- \ell| |\ell-a_1^\ast|^{-1}
\leq 2 |x_1- \ell | (5\eta(X))^{-1} \leq \alpha^2 \leq \lambda
\end{equation}
as in \eqref{10.39}, and because 
$|x_1- \ell| \leq 6\alpha^2 \tau < 6 \alpha^2\cdot 10^{-3} \eta(X)$ 
by \eqref{10.92} and \eqref{5.2bis}. Recall also that $\alpha$ can be chosen small, 
depending on $\lambda$ (see below Lemma \ref{t10.1}).
Altogether $|v-v(\ell,a_1^\ast)| \leq 205 \lambda$, $\rho(\ell,z_1)$ is a
$205\lambda$-Lipschitz graph over the line ${\rm Vect}(v(\ell,a_1
^\ast))$.
Since we already know something like this about $\Gamma_1 \cap D_2$, 
we also get that $\wt\Gamma_1 \cap D_2$ is a $208\lambda$-Lipschitz graph over that line.
Then we apply the same version of Lemma \ref{t10.2}, transposed for curves that start 
from $\ell$, and conclude as in Lemma \ref{t10.3}.

Now we prove \eqref{10.91}. Because we already know that $\Gamma_4$ and $\Gamma_5$
are $10^3\lambda$-Lipschitz graphs starting from $\ell$, we just need to show that 
the two corresponding geodesics leaving from $\ell$ are almost opposed. That is, it is enough
to prove that
\begin{equation} \label{10.99}
|v(x_0,\ell) - v(\ell,a_1^\ast)| \leq 35 \sqrt\lambda,
\end{equation}
say. But $|v(x_0, \ell) - v(x_0,x_1)| \leq 2 |x_0-\ell|^{-1} |x_1-\ell| \leq 64\lambda$
by \eqref{10.66}, $|v(\ell,a_1^\ast)-v(x_0,a_1^\ast)| \leq \lambda$ 
by \eqref{10.39}, and $|v(x_0,x_1) - v(x_0,a_1^\ast)| \leq 30 \sqrt\lambda$
by \eqref{10.67}; \eqref{10.99}, \eqref{10.91}, and Lemma \ref{t10.5} follow. 
\qed

\ms
We are ready to check \eqref{9.2}. For the fact that
\begin{equation} \label{10.100}
\text{$\Gamma_2$, $\Gamma_3$ and $\wt\Gamma_4$ make angles larger than 
$100^\circ$ at $x_0$,} 
\end{equation}
we use \eqref{10.23} (for the angle of $\Gamma_2$ and $\Gamma_3$),
and (for the two other angles) the fact that $\wt\Gamma_4$ coincides with 
$\Gamma_4$ near $x_0$, the proof of \eqref{10.23}, and the fact that by 
Lemma \ref{t10.4} the general direction $v(x_0,x_1)$ of $\Gamma_4$ is 
almost the same as the general direction $v(x_0,a_1^\ast)$ of $\Gamma_1$ 
in \eqref{10.23}.

Then the angle of $\wt\Gamma_4$ and $\wt\Gamma_5$ is controlled by
\eqref{10.91}, so \eqref{9.2} holds. The verification of \eqref{9.3} is the same as
usual; the fact that $\Gamma_5$ comes from $x_1$ rather than $x_0$ does not matter.

\ms
Now we prove the length estimates \eqref{9.6} and \eqref{9.7}. First we want to estimate
the extra length for the detour through $\ell$, and use the function
\begin{equation} \label{10.101}
f(z) = \ddist(z,z_4)+\ddist(z,z_5)
\end{equation}
defined on $\S$. We can still use \eqref{10.42} to differentiate 
(away from $\pm z_4$ and $\pm z_5$), and get that 
\begin{equation} \label{10.102}
- \nabla f(z) = v(z,z_4)+v(z,z_5).
\end{equation}
Then
\begin{eqnarray} \label{10.103}
|\nabla f(x_1)| &=& | v(x_1,z_4)+v(x_1,z_5)|
\nn\\
&\leq& | v(x_1,x_0)+v(x_1,a_1^\ast)| + | v(x_1,x_0)-v(x_1,z_4)|
+ | v(x_1,a_1^\ast)-v(x_1,z_5)|
\nn\\
&\leq& | v(x_1,x_0)+v(x_1,a_1^\ast)| + 6\lambda
\\
&\leq& | v(x_0,x_1)-v(x_0,a_1^\ast)| + |v(x_0,a_1^\ast) - v(x_1,a_1^\ast)|
+ 6\lambda
\leq 30\sqrt\lambda+ 7\lambda \leq 35\sqrt\lambda
\nn
\end{eqnarray}
because $z_4$ lies between $x_0$ and $x_1$ on a $3\lambda$-Lipschitz graph,
for the same reason for $z_5$, and by \eqref{10.67} and \eqref{10.92} 
(as in \eqref{10.39}, say). For $i=1,2$, $v(z,z_i)$ is differentiable, with
$|\nabla_z v(z,z_i) | \leq |z-z_i|^{-1}$. We shall use this for $z\in \rho(\ell,x_1)$.
Then 
\begin{equation} \label{10.104}
|z-x_1| \leq |\ell - x_1| < r/2
\end{equation}
by \eqref{10.88}, so $|z-z_i| \geq r/2$ (recall that $|z_i-x_1| = r$) and
$|\nabla_z v(z,z_i)| \leq 2 r^{-1}$. We sum over $i$ and integrate on $\rho(z,x_1)$;
this yields
\begin{equation} \label{10.105}
|\nabla f(x_1) - \nabla f(z)| \leq 2 r^{-1} \ddist(z,x_1) \leq 2 r^{-1} \ddist(\ell,x_1).
\end{equation}
We integrate on $\rho(x_1,\ell)$ and get that
\begin{equation} \label{10.106}
\begin{aligned}
f(\ell) &\leq f(x_1) + \ddist(\ell,x_1) |\nabla f(x_1)| + 2r^{-1} \ddist(\ell,x_1)^2
\cr&\leq f(x_1) + [35\sqrt\lambda+ 2r^{-1} \ddist(\ell,x_1)] \ddist(\ell,x_1) 
\cr&\leq f(x_1) + [35\sqrt\lambda+ 300\lambda] \ddist(\ell,x_1)
\leq f(x_1) + 36\sqrt\lambda |\ell-x_1|
\end{aligned}
\end{equation}
by \eqref{10.103} and \eqref{10.88}. But
\begin{equation} \label{10.107}
\begin{aligned}
\H^1(\wt\Gamma_4 \cup \wt\Gamma_5)
&= \H^1(\Gamma'_4)+\H^1(\Gamma'_5)+ \H^1(\rho(z_4,\ell))+\H^1(\rho(\ell,z_5))
\cr&= \H^1(\Gamma'_4)+\H^1(\Gamma'_5)+f(\ell)
\end{aligned}
\end{equation}
by definition of the $\wt\Gamma_i$, and
\begin{equation} \label{10.108}
\H^1(\Gamma_4 \cup \Gamma_5) \geq \H^1(\Gamma'_4)+\H^1(\Gamma'_5)+ f(x_1)
\end{equation}
because $\Gamma_4$ is composed of $\Gamma'_4$ and an arc from $z_4$ to $x_1$,
and similarly for $\Gamma_5$, so
\begin{equation} \label{10.109}
\H^1(\wt\Gamma_4 \cup \wt\Gamma_5) 
\leq \H^1(\Gamma_4 \cup \Gamma_5)+ f(\ell) - f(x_1)
\leq \H^1(\Gamma_4 \cup \Gamma_5) + 36\sqrt\lambda |\ell-x_1|.
\end{equation}
We add the contributions of $\Gamma_2$ and $\Gamma_3$ and get that
\begin{equation} \label{10.110}
\H^1(\Gamma) \leq \H^1\Big(\bigcup_{i=2}^5 \Gamma_i \Big) + 36\sqrt\lambda |\ell-x_1|
\leq \H^1\Big(\bigcup_{i=2}^5 \ol\gamma_i \Big) + 36\sqrt\lambda |\ell-x_1|
\end{equation}
by \eqref{7.16}.
But $\ol\gamma = \gamma_\ell \cup \big(\cup_{i=2}^5 \ol\gamma_i\big)$
(an essentially disjoint union), and $\H^1(\gamma_\ell) \geq |\ell-x_1|$ 
(because $\gamma_\ell$ goes from $\ell$ to $x_1$), so
\begin{equation} \label{10.111}
\H^1(\Gamma) \leq \H^1(\ol\gamma) - \frac{1}{2} \H^1(\gamma_\ell).
\end{equation}
This proves \eqref{9.6}; we are left with \eqref{9.7} to check.
Before we start, let us record the fact that
\begin{equation} \label{10.112}
r = (100\lambda)^{-1} |x_1 - \ell| \leq (100\lambda)^{-1} \H^1(\gamma_\ell)
\end{equation} 
by \eqref{10.88}. Next observe that for the symmetric difference $\Delta(\ol\gamma,\Gamma)$
of \eqref{9.4},
\begin{equation} \label{10.113}
\Delta(\ol\gamma,\Gamma) 
\subset \gamma_\ell \cup \Big(\bigcup_{i=2}^5 \Delta(\ol\gamma_i,\Gamma_i)\Big)
\cup \rho(z_4,\ell) \cup \rho(z_5,\ell) \cup \Gamma''_4 \cup \Gamma''_5,
\end{equation}
where $\Gamma''_4$ is the arc of $\Gamma_4$ between $z_4$ and $x_1$,
and similarly for $\Gamma''_5$. Since
\begin{equation} \label{10.114}
\begin{aligned}
\H^1(\gamma_\ell) + \H^1&(\rho(z_4,\ell)) +  \H^1(\rho(z_5,\ell)) + 
 \H^1(\Gamma''_4) + \H^1(\Gamma''_5)
 \cr&\leq \H^1(\gamma_\ell) + 10r \leq C(\lambda) \H^1(\gamma_\ell)
 \leq 2C(\lambda) [\H^1(\ol\gamma) -\H^1(\Gamma)]
 \end{aligned}
\end{equation}
by \eqref{10.112} and \eqref{10.111} and the right-hand side of \eqref{10.114}
is controlled by the right-hand side of \eqref{9.7},
we are left with the four $\Delta(\ol\gamma_i,\Gamma_i)$.
By \eqref{7.16} and \eqref{7.5},
\begin{equation} \label{10.115}
\H^1(\Gamma_i \sm \ol\gamma_i) \leq \H^1(\ol\gamma_i \sm \Gamma_i)
\leq C(\lambda)[\H^1(\Gamma_i) - \H^1(\rho_i)], 
\end{equation}
where $\rho_i$ is the geodesic arc between the endpoints of $\Gamma_i$, $2 \leq i \leq 5$.

What we want for \eqref{9.7} is $\sum_{i=2}^5 \H^1(\wt\rho_i)$, where we may keep 
$\wt \rho_i = \rho_i$ for $i=2,3$, but for $i=4,5$, 
$\wt\rho_i$ is the geodesic with the same endpoints as the corresponding arc $\wt\Gamma$.

We sum \eqref{10.115} over $i$ and get that
\begin{equation} \label{10.116}
\begin{aligned}
\sum_i \H^1(\Delta(\ol\gamma_i,\Gamma_i))
&\leq 2 C(\lambda) \sum_i [\H^1(\Gamma_i) - \H^1(\rho_i)] 
\cr&\leq 2 C(\lambda) \Big[\H^1(\Gamma) - \sum_i \H^1(\wt\rho_i) \Big] 
+ 2 C(\lambda) (\Sigma_1 + \Sigma_2),
\end{aligned}
\end{equation}
with
\begin{equation} \label{10.117}
\Sigma_1 = \sum_{i=2}^5 \H^1(\Gamma_i) - \H^1(\Gamma) 
\ \text{ and } \ 
\Sigma_2 = \sum_{i=2}^5 [\H^1(\wt\rho_i) - \H^1(\rho_i)].
\end{equation}
Notice that 
\begin{equation} \label{10.118}
\Sigma_1 = \sum_{i=4}^5 [\H^1(\Gamma_i) - \H^1(\wt\Gamma_i)] 
\leq \H^1(\Gamma''_4) + \H^1(\Gamma''_5)
\leq 2C(\lambda) [\H^1(\ol\gamma) -\H^1(\Gamma)]
\end{equation}
because $\H^1(\Gamma) = \H^1(\Gamma_2)+\H^1(\Gamma_3)
+\H^1(\wt\Gamma_4)+\H^1(\wt\Gamma_4)$
(since the union in \eqref{10.90} is essentially disjoint), then because 
the other part of $\Gamma_i$, namely $\Gamma'_i$, is contained in $\wt\Gamma_i$,
and finally by \eqref{10.114}. This part is dominated by the right-hand side of \eqref{9.7}.

As for $\Sigma_2$, first observe that $\rho_i = \wt \rho_i$ when $i \in \{ 2, 3 \}$.
When $i \in \{ 4, 5 \}$, the difference is that one endpoint is $\ell$ instead of $x_1$.
That is,
\begin{equation} \label{10.119}
\H^1(\wt \rho_4) - \H^1(\rho_4) = \H^1(\rho(x_0,\ell) - \H^1(\rho(x_0,x_1)) 
\leq \ddist(x_1,\ell),
\end{equation}
we have a similar estimate for $\H^1(\wt\rho_5) - \H^1(\rho_5)$
(just replace $x_0$ with $a_1^\ast$), and by \eqref{10.111}
\begin{equation} \label{10.120}
\Sigma_2 \leq 2 \ddist(x_1,\ell) \leq 3 |x_1-\ell| \leq 3 \H^1(\gamma_\ell) 
\leq 6[\H^1(\ol\gamma) -\H^1(\Gamma)],
\end{equation}
which is also dominated by the right-hand side of \eqref{9.7}.
Since the main part of the right-hand side of \eqref{10.116}
shows up as $\H^1(\Gamma)-\H^1(\rho)$ in \eqref{9.7}, we get the desired 
estimate for our last term $\sum_i \H^1(\Delta(\ol\gamma_i,\Gamma_i))$.
This completes our proof of \eqref{9.7}, and the verification of \eqref{9.2}-\eqref{9.7}
in Case B, the last case for Configuration 3+.

\section{The net $\Gamma$ for Configuration 2+}
\label{S11}

As was observed at the end of Section \ref{S9}, we still need to construct $\Gamma$
in the case of Configuration 2+. As for Configuration 3+, we will have two cases, 
one where we keep the same center $x_0$, and one where we go directly to $\ell$.

Recall that in the present case $E \cap \d D$ has only two points $c_i^\ast$
(or else there is a hanging curve, to be discussed later, but which we ignore for the moment), 
and $\gamma$ is composed of a simple curve $\gamma_{1,2}$ that goes from
$c_1^\ast$ to $c_2^\ast$, plus a simple curve $\gamma_\ell$, possibly reduced
to the point $\ell$, that goes from $\ell$ to a point of $\gamma_{1,2}$. We
call this point $x_0$, and denote by $\gamma_i$, $i\in \{ 1, 2 \}$, the arc of $\gamma_{1,2}$
that goes from $x_0$ to $c_i^\ast$. We also denote by $\ol\gamma_i$ the union of
$\gamma_i$ and the arc $\cL_i \subset E$ that goes from $c_i^\ast$ to $a_i^\ast$.
Finally set 
\begin{equation} \label{11.1}
\ol\gamma = \gamma_\ell \cup \ol\gamma_1 \cup \ol\gamma_2.
\end{equation}

We start as in Section \ref{S10}, apply the construction of Section \ref{S7}
to the three curves $\ol\gamma_1$, $\ol\gamma_2$, and $\gamma_\ell$,
and this gives three Lipschitz curves $\Gamma_1$, $\Gamma_2$, and
$\Gamma_3$. The simplest case, which we shall call \ub{Case A}, is when
\begin{equation} \label{11.2}
\text{$v(x_0,a_1^\ast)$, $v(x_0,a_2^\ast)$ and $v(x_0,\ell)$ make angles larger
than $\frac{2\pi}{3} - \frac{\pi}{10}$ with each other.} 
\end{equation}
In this case, we set
\begin{equation} \label{11.3}
\Gamma = \Gamma_1 \cup \Gamma_2 \cup \Gamma_3,
\end{equation}
and we can check \eqref{9.2}-\eqref{9.7} right away. There is only one angle condition
to check for \eqref{9.2}, at $x_0$, and it is satisfied because we claim that
\begin{equation} \label{11.4}
\text{the three $\Gamma_i$ make angles larger than $\frac{\pi}{2}$ at $x_0$.}
\end{equation} 
The claim follows readily from \eqref{11.2} and the Lipschitz graph description
of the $\Gamma_i$; see for instance the proof of \eqref{10.23}.

If there was also a hanging curve in the configuration, then there are three $\cC_i$ 
that leave from $\ell$, they make $120^\circ$ angles there, and the hanging curve 
is attached to the third point $c_3^\ast$ of $\d D$.
As promised in the description of Configuration H, we add the corresponding 
$\cL_i$ (just relabel if needed, and call it $\cL_3$) to $\Gamma$. 
The graphs $\Gamma_1$ and $\Gamma_2$ go essentially straight in the direction of
$a_1^\ast$ and $a_2^\ast$, and the same argument as for \eqref{9.19cc} shows that 
$\cL_3$ does not meet $\Gamma_1 \cup \Gamma_2$ outside of $D$. 
It does not meet $\Gamma_1$ or $\Gamma_2$ inside $D$ either, because it does not 
get inside $D$. So $\cL_3$ does not meet $\Gamma_1 \cup \Gamma_2 \cup \Gamma_3$, 
and we feel better.

Next \eqref{9.3} holds for the same reason as before (see the proof below \eqref{9.12}),
and we are left with the length estimates. First of all, \eqref{9.6} holds, simply by adding the
three estimates \eqref{7.16} coming from the three curves.
The symmetric difference $\Gamma \Delta \ol\gamma$ of \eqref{9.6}
is contained in the union of the symmetric differences
$\Gamma_1\Delta \ol\gamma_1$, $\Gamma_2 \Delta \ol\gamma_2$,
and $\Gamma_3 \Delta \gamma_\ell$, so \eqref{9.7} 
follows by adding the three estimates from the end of \eqref{8.21}.

\ms
We may now switch to \ub{Case B}, which is when \eqref{11.2} fails. 
We shall try a set $\Gamma$ that goes more directly through $\ell$, without passing 
through $x_0$; the construction will look like what we did for Case B of Configuration 3+.

Let $\alpha > 0$ be a small constant, which is allowed to depend on $\lambda$
and will be chosen near the end of the section. 
We need an analogue of Lemma \ref{t10.1}, which says that if $\varepsilon$ is chosen small
enough in \eqref{4.3} (depending on $\alpha$ and $\tau$),
\begin{equation} \label{11.5}
|x_0 - \ell| + |a_1-a_1^\ast| + |a_2-a_2^\ast| \leq 2\alpha^2 \tau.
\end{equation}
The proof is the same as for Lemma \ref{t10.1}, with only two branches coming from
the $c_i^\ast$. The reader may be worried about the special case when there is a third
point $c_3^\ast \in D_\pm(\tau)$, that leads to a hanging curve. But, as long as we 
stay in the spherical annulus $A = \S \cap B(\ell,2\tau) \sm B(\ell,\alpha^2\tau)$,
this third curve $\cL_3'$ stays far from the other two and does not interfere with the proof
(which, as the reader recalls, consists in saying that we don't meet a triple point like $x_0$). 

We set (more or less as usual)
\begin{equation} \label{11.6}
r = \alpha^{-1} |x_0 - \ell|,  
D = \S \cap B(x_0,r), \text{ and } \d D = \S \cap \d B(x_0,r)
\end{equation}
(compare with \eqref{10.25}); notice that $r \leq 2\alpha \tau$, by \eqref{11.5},
so it is still very small. Since $\Gamma_1$ and $\Gamma_2$ are small Lipschitz graphs, 
they meet $\d D$ exactly once, at points that we call $z_1$ and $z_2$. 
Set, for $i=1,2$,
\begin{equation} \label{11.7}
\Gamma'_i = \Gamma_i \sm D, \  \Gamma''_i = \Gamma_i \cap D, 
\ \text{ and } \wt\Gamma_i = \rho(\ell,z_i) \cup \Gamma'_i.
\end{equation}
Thus $\wt\Gamma_i$ is a curve that goes from $\ell$ to $a_i^\ast$,
while $\Gamma_i$ goes from $x_0$ to $a_i^\ast$; both curves go through $z_i$.
Finally we set
\begin{equation} \label{11.8}
\Gamma = \wt\Gamma_1 \cup \wt\Gamma_2,
\end{equation}
which we really see as a collection of two graphs. So we want to claim that for
$i= 1,2$,
\begin{equation} \label{11.9}
\wt\Gamma_i  \text{ is a $8\lambda$-Lipschitz graph over } \rho(\ell,a_i^\ast). 
\end{equation}
Fortunately, the proof is the same as for Lemma \ref{t10.3}, so we can skip it.
The point is that since $z_i$ and $x_0$ both lie on the small Lipschitz graph 
$\Gamma_i$, the geodesic $\rho(z_i, x_0)$ is almost aligned with $\Gamma'_i$,
and then $\rho(z_i,\ell)$ makes a small angle because \eqref{11.6} says that
$|x_0-\ell| << r$.

Next we claim that
\begin{equation} \label{11.10}
\text{$\wt\Gamma_1$ and $\wt\Gamma_2$ make an angle larger
than $\frac{\pi }{ 2}$ at $\ell$.}
\end{equation}
We start from the description of Proposition \ref{t2.1}, which says that
$\cC_1$ and $\cC_2$ make an angle at least $\frac{2\pi}{3}$ at $\ell$
(see \eqref{2.7}). That is,
\begin{equation} \label{11.11}
\Angle(v(\ell,a_1),v(\ell,a_2)) \geq \frac{2\pi}{3}.
\end{equation}
Let $v_i$ be a tangent direction to $\wt\Gamma_i$ at $\ell$.
From \eqref{11.9} (and \eqref{10.12}) we deduce that 
\begin{equation} \label{11.12}
|v_i - v(\ell,a_i^\ast)| \leq 32\lambda.
\end{equation}
In addition,
\begin{equation} \label{11.13}
|v(\ell,a_i^\ast) - v(\ell,a_i)| 
\leq 2 |a_i-a_i^\ast| |\ell - a_i|^{-1} 
\leq 4 \alpha^2\tau (5\eta(X))^{-1} \leq \alpha^2 < \lambda
\end{equation}
by \eqref{11.5}, \eqref{3.11}, \eqref{3.12}, and \eqref{5.2bis},
and if $\alpha$ is small enough. Now \eqref{11.10} follows from \eqref{11.11},
\eqref{11.12}, and \eqref{11.13}.

Notice that the description of $\Gamma$ as a union of small Lipschitz 
graphs follows from \eqref{11.9}, \eqref{9.2} (the control on the inside angles) 
follows from \eqref{11.10}, 
and \eqref{9.3} holds for the usual reason (see below \eqref{9.12}).

There may also be a hanging curve in the configuration. Then there are three
$\cC_i$ that leave from $\ell$, they make $120^\circ$ angles, and the hanging
curve is attached to the third point of $E\cap \d D$, i.e., $c_3^\ast$.
We add to $\Gamma$ the corresponding curve $\cL_3$, and it is good to know that $\cL_3$
does not meet $\Gamma$. This is the case, because the $\wt\Gamma_i$ are
small Lipschitz graphs over the $\rho(\ell,a_i^\ast)$ (by \eqref{11.9}), which
go away from $c_3^\ast$ and $\cL_3$; the proof goes as for \eqref{9.19cc}.

So we just need to prove the two usual length estimates.
For this we introduce
\begin{equation} \label{11.14}
f(z) = \ddist(z,z_1) + \ddist(z,z_2) 
\end{equation}
and estimate its gradient
\begin{equation} \label{11.15}
\nabla f(z) = -[v(z,z_1)+v(z,z_2)].
\end{equation}
We need to estimate some angles. Set 
\begin{equation} \label{11.16}
e_1 = v(x_0,a_1^\ast), \ e_2 = v(x_0,a_2^\ast), \text{ and } e_3 = v(x_0,\ell).
\end{equation}
Observe that for $i=1,2$ and $z\in \rho(x_0,\ell)$
\begin{equation} \label{11.17}
|v(z,z_i)-v(x_0,z_i)| \leq 2 |z-x_0| |x_0-z_i|^{-1}
\leq 2|\ell-x_0| r^{-1} = 2\alpha
\end{equation}
by \eqref{11.6}, and
\begin{equation} \label{11.18}
|v(x_0,z_i)-e_i| = |v(x_0,z_i)-v(x_0,a_i^\ast)| \leq 4\lambda
\end{equation}
because $z_i \in \Gamma_i$, which is a $\lambda$-Lipschitz graph
over $\rho(x_0,a_i^\ast)$ and by \eqref{10.12}; then
\begin{equation} \label{11.19}
|\nabla f(z) + e_1 + e_2| \leq 12\lambda
\end{equation}
if $\alpha$ is small enough. Next we want to check that
\begin{equation} \label{11.20}
\langle e_1+e_2, - e_3 \rangle \leq 1 - 10^{-2}.
\end{equation}
First observe that for $i=1,2$,
\begin{eqnarray} \label{11.21}
|e_i- v(\ell,a_i)| 
&=& |v(x_0,a_i^\ast) - v(\ell,a_i)| 
\leq |v(x_0,a_i^\ast)-v(\ell,a_i^\ast)| + |v(\ell,a_i^\ast) - v(\ell,a_i)|
\nn\\
&\leq&   2 |\ell-x_0| |\ell - a_i^\ast|^{-1} + \alpha^2
\leq 4 \alpha^2\tau (5\eta(X))^{-1} + \alpha^2 \leq 2\alpha^2
\end{eqnarray}
by \eqref{11.13}, \eqref{11.5} \eqref{3.11}, \eqref{3.12}, and \eqref{5.2bis}.
Then by \eqref{11.11}
\begin{equation} \label{11.22}
\Angle(e_1,e_2) > \frac{2\pi}{3} - 5\alpha^2.
\end{equation}
Since \eqref{11.2} fails (by definition of Case B), the three vectors $e_1$, $e_2$,
$e_3$ do not all make angles larger than $\frac{2\pi}{3} - \frac{\pi}{10}$.
Since this is the case for $e_1$ and $e_2$, we may assume, without loss of generality,
that
\begin{equation} \label{11.23}
\Angle(e_1, e_3) \leq \frac{2\pi}{3} - \frac{\pi}{10}.
\end{equation}
If \eqref{11.20} fails, $|e_1+e_2| \geq 1-10^{-2}$. Set $\theta = \frac{1}{2}\Angle(e_1,e_2)$;
then $|e_1+e_2| = 2\cos\theta$, hence 
\begin{equation} \label{11.24}
\theta \leq \arccos 
\Big(\frac{1}{2}-\frac{10^{-2}}{2}\Big) \leq \frac{\pi}{3} + 10^{-2}.
\end{equation}
Then 
\begin{eqnarray} \label{11.25}
\Angle(e_3,e_1+e_2) 
&\leq& \Angle(e_3,e_1) + \Angle(e_1,e_1+e_2) 
\nn\\
&\leq& \Angle(e_3,e_1) + \theta \leq \pi - \frac{\pi}{10} + 10^{-2}
\leq \pi - \frac{2}{10},
\end{eqnarray}
$\Angle(- e_3,e_1+e_2) \geq \frac{2}{10}$, 
and $\cos(\Angle(-e_3,e_1+e_2)) \leq 1 - 2 \cdot 10^{-2}$. Since by \eqref{11.22}
\begin{equation} \label{11.26}
|e_1+e_2| =  2\cos\theta \leq 2\cos\big(\frac{\pi}{3} - \frac{5\alpha^2}{2}\big) \leq 1+10^{-4},
\end{equation}
\begin{equation} \label{11.27}
\langle e_1+e_2, - e_3 \rangle = |e_1+e_2| \cos(\Angle(-e_3,e_1+e_2))
\leq (1+10^{-4})(1 - 2 \cdot 10^{-2}) 
\leq 1 - 10^{-2}
\end{equation}
and \eqref{11.20} holds after all. 

We may now return to the computation of $\nabla f(z)$. We integrate on $\rho(x_0,\ell)$
and get that
\begin{equation} \label{11.28}
\begin{aligned}
f(\ell) &= f(x_0) + \int_{\rho(x_0,\ell)} \langle \nabla f(z), v(z,\ell) \rangle \, d\H^1(z)
\cr& \leq f(x_0) + \int_{\rho(x_0,\ell)}
\big[ \langle e_1+e_2, - v(z,\ell) \rangle + 10\lambda \big] \, d\H^1(z) 
\cr& \leq f(x_0) 
+ \int_{\rho(x_0,\ell)} \big[ \langle e_1+e_2, - e_3 \rangle + 11\lambda \big] \, d\H^1(z) 
\end{aligned}
\end{equation}
by \eqref{11.19} and because 
\begin{equation} \label{11.29}
|v(z,\ell)-e_3| = |v(z,\ell)-v(x_0,\ell)| \leq |x_0-\ell | \leq 2 \alpha^2 \tau < \frac{\lambda}{2}
\end{equation}
for $z\in \rho(x_0,\ell)$, because geodesics do not turn too fast, by \eqref{11.5}, 
and if $\alpha$ is small enough.
By \eqref{11.20} and the definition \eqref{11.14}, this yields
\begin{eqnarray} \label{11.30}
\ddist(\ell,z_1) + \ddist(\ell,z_2) &=& f(\ell)  
\leq f(x_0) + (1 - 10^{-2} + 11\lambda) \ddist(x_0,\ell)
\nn\\
&\leq& \ddist(x_0,z_1) + \ddist(x_0,z_2)  + (1 - 10^{-3}) \ddist(x_0,\ell)
\nn\\
&\leq& \H^1(\Gamma''_1) + \H^1(\Gamma''_2) + (1 - 10^{-3}) \ddist(x_0,\ell)
\end{eqnarray}
because $\Gamma''_i$ precisely goes from $x_0$ to $z_i$ (see near \eqref{11.7}).

We are about ready for \eqref{9.6} and \eqref{9.7}. Because of \eqref{11.1},
\begin{equation} \label{11.31}
\H^1(\ol\gamma) = \H^1(\gamma_\ell) + \H^1(\ol\gamma_1)+\H^1(\ol\gamma_2)
\geq \H^1(\gamma_\ell) + \H^1(\Gamma_1)+\H^1(\Gamma_2)
\end{equation}
by \eqref{7.16}. Besides, by \eqref{11.8} and \eqref{11.7},
\begin{equation} \label{11.32}
\H^1(\Gamma) = \sum_{i=1}^2 \H^1(\wt \Gamma_i)
= \sum_{i=1}^2 [\H^1(\Gamma'_i) + \ddist(x_0,z_i)]
\end{equation}
and, since $\H^1(\Gamma_i) = \H^1(\Gamma'_i)+\H^1(\Gamma''_i)$,
we get that
\begin{equation} \label{11.33}
\begin{aligned}
\H^1(\ol\gamma) - \H^1(\Gamma)
&\geq \H^1(\gamma_\ell) + 
\sum_{i=1}^2 \big[\H^1(\Gamma_i)-\H^1(\Gamma'_i) - \ddist(x_0,z_i)\big]
\cr& = \H^1(\gamma_\ell) + 
\sum_{i=1}^2 \big[\H^1(\Gamma''_i) - \ddist(x_0,z_i)\big]
\cr&\geq \H^1(\gamma_\ell) - (1 - 10^{-3}) \ddist(x_0,\ell)
\geq 10^{-3} \H^1(\gamma_\ell) \geq 10^{-3} |x_0-\ell |
\end{aligned}
\end{equation}
by \eqref{11.8} and \eqref{11.7}, and because 
$\H^1(\gamma_\ell) \geq \ddist(x_0,\ell) \geq |x_0-\ell|$.

This is better than \eqref{9.6}. For \eqref{9.7}, first observe that
\begin{equation} \label{11.34}
r = \alpha^{-1} |x_0-\ell| \leq 10^3 \alpha^{-1} [\H^1(\ol\gamma) - \H^1(\Gamma)],
\end{equation}
which is therefore controlled by the right-hand side of \eqref{9.7}. We write
\begin{equation} \label{11.35}
\ol\gamma\Delta\Gamma \subset \gamma_\ell 
\cup \big(\bigcup_{i=1}^2 \ol\gamma_i \Delta\Gamma_i \big)
\cup \big(\bigcup_{i=1}^2 \Gamma_i \Delta \wt\Gamma_i \big).
\end{equation}
The first part is in order, since $\H^1(\gamma_\ell) 
\leq 10^3[\H^1(\ol\gamma) - \H^1(\Gamma)]$
by \eqref{11.33}. The second one as well, because 
\begin{equation} \label{11.36}
\H^1(\Gamma_i \Delta \wt\Gamma_i) \leq \H^1(\Gamma''_i) + \H^1(\rho(z_i,\ell)) \leq 5r,
\end{equation}
and by \eqref{11.34}. We are left with the
\begin{equation} \label{11.37}
\H^1(\ol\gamma_i \Delta\Gamma_i) 
\leq \H^1(\ol\gamma_i \sm \Gamma_i)+\H^1(\Gamma_i \sm \ol\gamma_i)
\leq C \lambda^{-2} [\H^1(\Gamma_i)-\ddist(x_0,a_i^\ast)]
\end{equation}
by \eqref{7.16} and \eqref{7.5} (recall that $\Gamma_i$ goes from $a_i^\ast$ to $x_0$).
In turn $\H^1(\Gamma_i) \leq \H^1(\wt\Gamma_i) + 2r$ and 
$\ddist(\ell,a_i^\ast) \leq \ddist(x_0,a_i^\ast) + 2r$, so
\begin{equation} \label{11.38}
\H^1(\Gamma_i)-\ddist(x_0,a_i^\ast) 
\leq \H^1(\wt\Gamma_i) - \ddist(\ell,a_i^\ast) + 4r,
\end{equation}
which is also controlled by the right-hand side of \eqref{9.7}.
This completes the verifications in Case B of Configuration 2+.
We finally constructed the Lipschitz net $\Gamma$ in all cases.

\section{Lipschitz projections near $E \cap \S$}
\label{S12}

In Sections \ref{S9}-\ref{S11} we started from a net $\ol\gamma$ 
of curves in $E\cap \S$, in fact near a point $\ell \in K\cap L$,
and constructed a corresponding net $\Gamma$ of Lipschitz curves on $\S$.

We did this for each $\ell$, and in some cases (Configurations H and $3= 2+1$)
independently for the two or three configurations present near $\ell$.
For the other curves $\cL_i$, the ones for which both endpoints of $\cC_i$
lie in $V_1\cup V_2 = V \sm V_0$, the simplest is to set $\Gamma = \ol\gamma = \cL_i$.

For this to work well, it will be better to know that the $\cL_i$ that do not get close to
$L$ are $\lambda$-Lipschitz graphs over the geodesics with the same endpoints.
This is why we made Remark \ref{r6.3n} and the similar later ones. If we did not do this, we still 
would be all right, but we would need to replace each of these $\cL_i$ with the small Lipschitz
curve $\Gamma_i$ obtained from $\cL_i$ by the construction of Section \ref{S7}.
We would also need to check that this replacement does not alter much the angle conditions
\eqref{9.3} with the other curves or nets $\Gamma$, but nothing dramatic.

We now let $\gamma^\ast$ denote the union of all the curves $\ol\gamma$ that
we have here, and $\Gamma^\ast$ the union of all the $\Gamma$ that we constructed.
In the cases (as Configuration $3-$) where some points of $K\cap L$ do not lie in any
constructed $\Gamma$, we just add them to $\Gamma^\ast$ as isolated points.
Thus $\Gamma^\ast$ can be decomposed into nets of one to four small Lipschitz
curves, plus maybe one or two points of $K\cap L$. 

In this section we shall build a Lipschitz projection on $\Gamma^\ast$. 
In fact, the term is a little inappropriate, because what we are interested in is a collection
of Lipschitz mappings, defined in small neighborhoods of the main connected components of 
$E \cap \S$ and with values in $\Gamma^\ast$. Let us explain what we want.

\begin{pro}\label{t12.1} 
We can find a small number $\tau_3 > 0$ and a Lipschitz mapping $p$,
defined on
\begin{equation} \label{12.1}
E_+ = E_+(\tau_3) = \big\{ x \in \S \, ; \, \dist(x,E\cap \S) \leq \tau_3 \big\},
\end{equation}
with values in $\Gamma^\ast$, such that
\begin{equation} \label{12.2}
|p(x)-x| \leq 60\tau 
\ \text{ for } x\in E_+  \, ,
\end{equation}
\begin{equation} \label{12.3}
p \ \text{ is $30$-Lipschitz on } E_+ \cap B(x,2\tau_3)
\end{equation}
for each $x\in E_+$, and
\begin{equation} \label{12.4}
\text{$p(\ell) = \ell$ for $\ell \in K \cap L$.}
\end{equation}
\end{pro}

\ms
Here $\tau$ in \eqref{12.2} is as above, but 
$\tau_3$ and the Lipschitz constant for $p$ may depend very badly on
the set $E$ and the initial radius (here normalized to $1$) that we took in Section \ref{S4}.
So we will need to be careful when we apply the proposition; what will save us is the
local Lipschitz bound \eqref{12.3}, which will be used to control the measure of the image.

We cannot hope to get a continuous projection which is defined on the sphere,
because even if $n=3$ and $\Gamma^\ast$ is a great circle, there is a topological issue
(where do we send the two poles?). Also, $\Gamma^\ast$ may have more than one connected
component, different domains of the sphere will be sent to different components, and so 
we count on small gaps in $E+$, coming from the fact that $E\cap \S$ also is not connected,
to patch the various local Lipschitz mappings that we take.

The distance estimate \eqref{12.2} is rather poor (we can expect much better in many
cases), but it will be enough. 

Anyway, we shall start our proof with the construction of local Lipschitz projections 
defined on relatively large pieces of $\S$, that we will then need to patch together.

We shall use the description of $E\cap S$ near the curves $\cC_i$
that was given in Proposition~\ref{t5.4}; in particular $\tau$ in \eqref{12.2} and 
below is still coming from
this proposition, we assume that $\tau \leq 10^{-3}\eta(X)$, and the
$\cL_i$ are the curves in $E \cap \S$ provided by the proposition.

Recall that we split $\cI$ into $\cI_0$ (the indices for which $\cC_i$
has an endpoint on $K \cap L$) and $\cI_1 = \cI \sm \cI_0$.
When $i\in \cI_0$, we shall denote by $\ell(i)$ the point of $K \cap L$
where $\cC_i$ ends, and by $D_i$ the spherical disk
$D_i = \S \cap B(\ell(i),\tau)$ associated to $\ell(i)$.

Recall also that when $i\in \cI_1$, the curve $\cL_i$ given by Proposition \ref{t5.4}
connects two vertices $a_i^\ast$ and $b_i^\ast$ (that lie close to the vertices $a_i$
and $b_i$ of $\cC_i$), while for $i \in \cI_0$, $\cL_i$ start at a vertex $a_i^\ast$ but 
end at a point $c_i^\ast$ of $\d D_i$.

For $i\in \cI$ we define a region of influence $R_i$ by
\begin{equation}\label{12.5}
R_i = \big\{ z\in \S \, ; \, \dist(z,\cL_i) \leq 10^{-1}\tau
\text{ and } \dist(z,\cL_i) \leq \dist(z,\cL_j) \text{ for } j \in \cI \sm \{ i \} \big\}.
\end{equation}
Then, when $i\in \cI_1$, we define a projection $p_i$ such that
\begin{equation} \label{12.6}
p_i : R_i \to \cL_i \ \text{ is $3$-Lipschitz}, 
\end{equation}
and
\begin{equation} \label{12.7}
|p_i(x) - x| \leq 3 \dist(x,\cL_i) < \tau \ \text{ for } x\in R_i,
\end{equation}
where the second inequality follows from \eqref{12.5}. For the moment, this is easy
to arrange because $\cL_i$ is such a nice curve.

Now it could be that $\cL_i$ shares an endpoint $a^\ast$ with one or two other $\cL_j$,
and to avoid conflicts, we require that when this happens we take
\begin{equation} \label{12.8}
p_i(z) = a^\ast \ \text{ for } 
z\in R_i \cap R_j = \big\{ x\in R_i \, ; \, \dist(z,\cL_i) = \dist(z,\cL_j) \big\}.
\end{equation}
This is easy to arrange: recall that at $a^\ast$, the worse that can happen is that 
two other $\cL_j$ end at $a^\ast$, in a nice $C^1$ way and with large angles:
see \eqref{5.24} (for endpoints in $V_1$) and \eqref{5.34} (for endpoints in $V_2$).

Also, we claim that this is enough to guarantee that if $j\in \cI_1\sm \{ i \}$ and if
we set $p(x) = p_i(x)$ for $x\in R_i \cap B(a^\ast,\tau)$ and
$p(x) = p_j(x)$ for $x\in R_j \cap B(a^\ast,\tau)$, then
\begin{equation} \label{12.9}
\text{$p$ is $10$-Lipschitz on } (R_i \cup R_j) \cap B(a^\ast,\tau).
\end{equation}
The point is that if $x\in R_i \cap B(a^\ast,\tau)$ and $y \in R_j \cap B(a^\ast,\tau)$,
then there is a path from $x$ to $y$ in $(R_i \cup R_j) \cap B(a^\ast,\tau)$, of length
at most $10|x-y|$, and that goes through some point $z\in R_i \cap R_j$.
Then $p_i(z) = p_j(z) = a^\ast$, and the fact that 
$|p(x)-p(y)| \leq |p_i(x)-p_i(z)| + |p_j(z)-p_j(y)|$ gives the right estimate.

Notice also that when $i, j \in \cI_1$ are such that $\cL_i$ and $\cL_j$ do not share an endpoint,
then $\dist(R_i,R_j) \geq \frac{\tau}{10}$, by \eqref{5.39}, \eqref{5.4}, and \eqref{5.2bis}.
This means that if we set $p(z) = p_i(z)$ for $i\in \cI_1$ and $z\in R_i$, not only the definitions are compatible, but we get a Lipschitz mapping on $\cup_{i\in \cI_1} R_i$.

This will take care of most of the sphere, but the most interesting part will be what we
do near the two points of $\S\cap L$. In fact, only the points of $K \cap L = \S \cap X \cap L$
matter, because if $\ell \in \S\cap L \sm K$, \eqref{3.10}, \eqref{3.12}, \eqref{5.2bis}, 
and \eqref{4.3} say that neither $X$ nor $E$ gets within $2\tau$ of $\ell$.

Let us review a little what we did in Section \ref{S6} and add some notation.
For each $\ell \in K \cap L$, we introduced a small disk $D = D_\ell$, then we wrote 
the curves $\cC_i$ that end at $\ell$ as $\cC_1, \cdots \cC_m$, introduced the
components $H_i$ of $c_i^\ast$ (the endpoint of $\cL_i$) in $E\cap D$, and then
grouped the $\cC_i$ by components. Let us denote by $CC(\ell)$ the set of connected
components $H_i$ (we need a different name, because some different indices $i$ may
give the same component). In Configuration $3 = 2+1$, for instance, $CC(\ell)$ has two
elements; in Configuration $3+$ or $3-$, $CC(\ell)$ has just one element.

For each $c\in CC(\ell)$, we have a connected set $\gamma = \gamma_c$, 
which we eventually completed into the larger $\ol\gamma = \ol\gamma_c$, 
and modified to get a net $\Gamma=\Gamma_c$. We may need to use
the set $\cI(c)$ of indices $i\in \cI_0$ such that $a_i^\ast \in c$
(or equivalently $H_i \subset c$).

We also complete $CC(\ell)$: if $\ell$ lies in one the components $c\in CC(\ell)$
we keep $CC_+(\ell) = CC(\ell)$. Otherwise, we add the special component 
$c_\ell = H_\ell$ (the component of $\ell$ in $E\cap D$, 
which is disjoint from the other ones), and associate to it the degenerate curves
$\gamma_{c_\ell} = \{ \ell \}$ and $\Gamma_{c_\ell} = \{ \ell \}$.
We do the same thing (i.e., add $\gamma_{c_\ell} = \{ \ell \}$ and $\Gamma_{c_\ell} = \{ \ell \}$)
if $\ell$ does not even lie in $E$. 
Then we set $CC_+(\ell) = CC(\ell) \cup \{c_\ell\}$.

Also denote by $CC$ the union of the $CC(\ell)$, $\ell \in K\cap L$,
and $CC_+$ the union of the $CC_+(\ell)$, $\ell \in K\cap L$.
Finally, if $c\in CC(\ell)$, we set $\ell(c) = \ell$ and $D_c = \S \cap B(\ell,\tau)$;
this is unambiguous, because a single curve $\cC_i$ never has both points of $\S \cap L$
as endpoints.

Our next step is the construction of mappings $p_c$, $c\in CC_+$.
When $c$ is one of the special components $c_\ell$, $\ell \in K\cap L$, 
we have set $\Gamma_c = \{ \ell \}$ and now we take
\begin{equation} \label{12.10}
p_{c_\ell}(z) = \ell \ \text{ for } z\in \S.
\end{equation}
The more interesting case of $c\in CC$ is treated in the next lemma.

\begin{lem} \label{t12.2}
For $c\in CC(\ell)$, set $D(c) = D_c \cup \bigcup_{i \in \cI(c)} R_i$. There
is a mapping $p_c$ such that
\begin{equation} \label{12.11}
p_c : D(c) \to \Gamma_c \ \text{ is $10$-Lipschitz,}
\end{equation}
\begin{equation} \label{12.12}
|p_c(x) - x| \leq 10 \dist(x,\Gamma_c) \ \text{ for } x\in D(c), 
\end{equation}
and, for each $i \in \cI(c)$ and each index $j\in \cI$ such that $a_i^\ast$ is also
an endpoint of $\cL_j$, \eqref{12.8} holds with $a^\ast = a_i^\ast$.
\end{lem}

\ms
Recall that $a_i^\ast$ is the endpoint of $\cL_i$ that lies far from $D_c$ (i.e., which is not
$c_i^\ast$). Also, all the indices $j\in \cI$ such that $a_i^\ast$ is also an endpoint of $\cL_j$
lie in $\cI_1$ and we already defined $p$ on the corresponding $R_j$.

The domain $D(c) = D_c \cup \bigcup_{i \in \cI(c)} R_i$ 
is composed of a central disk, which is so small that it is bilipschitz equivalent to a ball in 
$\R^{n-1}$ with a constant close to $1$, plus a small number (between one and three) 
of appendices that are thin tubes $R_i$ around $C^1$ curves $\cL_i$, 
and leave from $D$ in directions that make large angles. 
The set $\Gamma_c$ itself is a net of at most four small Lipschitz curves
(that make large angles when they meet), and $\Gamma_c$ reaches the same 
$a_i^\ast$, $i\in \cI(c)$.
In each case, the construction of $p_c$ is rather easy, but may be painful to write explicitly. 
This is why we shall simply review the different cases that we encounter, and hopefully 
the reader will agree that $p_c$ is not hard to find. 

In the case of Configuration $0$ (when there is no curve near $\ell$), 
there is no $\Gamma$ and we still do nothing.

When $c$ comes from a configuration of type $1$, $\Gamma$ is a small 
Lipschitz curve that goes from $\ell$ to $a_i^\ast$ (where $i$ is the 
only index in $\cI(c)$), and projecting on $\Gamma$ is easy. The additional condition 
\eqref{12.8} is not hard to get either, and we could easily get a $3$-Lipschitz function.

When $c$ is of type H, and we consider one of the hanging curves, 
recall that we started from $\ol\gamma_c = \cL_i$, where $i\in \cI_0$ is 
the index such that the hanging curve contains $c_i^\ast$, and we kept
$\Gamma_c = \ol\gamma_c = \cL_i$. In this case too $p_c$ is easy to find.

When $c$ is of type $2-$, $\gamma_c$ is the union of two simple curves $\gamma_i$
that leave from a same center $x_0$, and $\Gamma_c$ is the union of two Lipschitz
curves $\Gamma_i$ with the same endpoints $x_0$ and $a_i^\ast$, and that make
a large angle at $x_0$. Here too $p_c$ is easy to construct.

When $c$ comes from a configuration of type $3=2+1$, 
we combine the types $1$ and $2-$ above. We don't even need to know
that the two corresponding sets $\Gamma_c$ are disjoint, because we 
build two independent projections $p_c$ on different sets $\Gamma_c$.
The fact that the $p_c$ share a piece of their domains of definitions will be compensated
by the fact that we will later restrict the $p_c$ to disjoint domains at positive distances
from each other.

When $c$ is of type $2+$, $\Gamma_c$ is either a truncated $Y$ that connects
$\ell$ to $a_i^\ast$ and $a_j^\ast$ (where $\cI(c) = \{ i, j \}$), or composed
of two small Lipschitz graphs from $\ell$ to $a_i^\ast$ and $a_j^\ast$ (and thus
make a large angle at $\ell$. This case and the next one are just a little harder to
treat than the previous ones, but we shall only comment on the last one because it looks uglier.

When $c$ is of type $3-$, $\Gamma_c$ is a small Lipschitz spider that goes from a 
center $x_0$ to the three relevant $a_i^\ast$, and is not hard to project on.

Finally, when $c$ comes from a configuration of type $3+$, 
$\Gamma_c$ is either a small Lipschitz spider that goes from $\ell$ to the 
three relevant $a_i^\ast$ (as in Case A), or a slightly more complicated union
of $4$ small Lipschitz graphs, coming from case $B$.
As in the previous cases, all the angles between the curves are larger than $\pi/2$.

\begin{figure}[!h]  
\centering
\includegraphics[width=8cm]{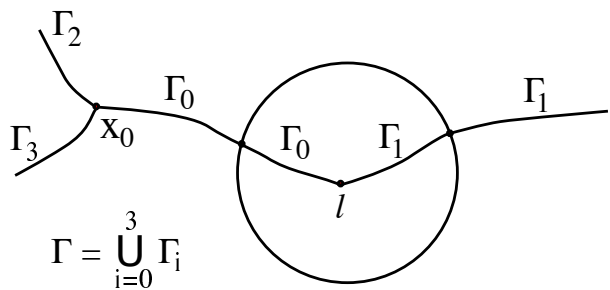}
\caption{ The curve $\Gamma = \Gamma_c$.
\label{f13-1b}}
\end{figure}

Let us only explain how we find $p_c$ in the apparently most complicated 
case B of type $3+$. Here (see Figure \ref{f13-1b}) %
$\Gamma_c$ is composed of two long curves $\Gamma_2$ and $\Gamma_3$, 
that connect a center $x_0$ to exterior points $a_2^\ast$ and $a_3^\ast$, 
a short curve $\Gamma_0$ (previously composed of a Lipschitz curve and a piece of 
geodesic, but we put them together) that goes from $x_0$ to $\ell$, 
and a third long curve (again originally composed of a geodesic and a piece of curve)
$\Gamma_1$ from $\ell$ to $a_1^\ast$. 
As in the previous cases, all these curves $\Gamma_i$ are small Lipschitz
graphs over the geodesics $\rho_i$ with the same endpoints, and they make large angles 
where they meet. 

 % p108

\begin{figure}[!h]  
\centering
\includegraphics[width=6cm]{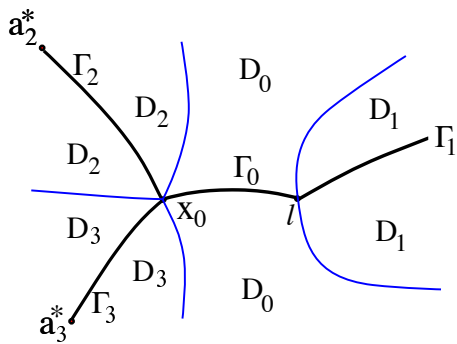}
\caption{ The domains $D_i = D_i(c)$.
\label{f13-1c}}
\end{figure}

We will cut $D(c)$ into a few simple regions $D_i(c)$, and then take a simple definition
for $p_c$ on each piece. Set $\rho = \bigcup_{i=0}^3 \rho_i$, and then
\begin{equation} \label{12.13}
D_1(c) = \big\{ x\in D(c) \, ; \,  \dist(x,\rho_1) 
\leq \frac{1}{3} \dist(x,\rho \sm \rho_1)\big\}.
\end{equation}
For $i \in \{2,3\}$, choose
\begin{equation} \label{12.14}
D_i(c) = \big\{ x\in D(c) \, ; \,  \dist(x,\rho_i) = \dist(x,\rho) \text{ and }
\dist(x,\rho_i) \leq  \frac{1}{3} \dist(x,\rho_0 \cup \rho_1)\big\},
\end{equation}
and finally set $D_0(c) = D(c) \sm \bigcup_{i=1}^3 D_i(c)$. 
See Figure \ref{f13-1c} %
for a sketch of our four domains (sitting in the two-dimensional sphere $\S$) when $n=3$. 
The case when $n > 3$ is not different; the common boundaries just have a larger dimension,
and the three domains $D_i(c)$, $0 \leq i \leq 2$, now have a $(n-3)$-dimensional common
boundary that goes through $x_0$ (when $n=4$, think about a curve through
$x_0$ that crosses the plane of the  picture).

 % p108

\begin{figure}[!h]  
\centering
\includegraphics[width=6cm]{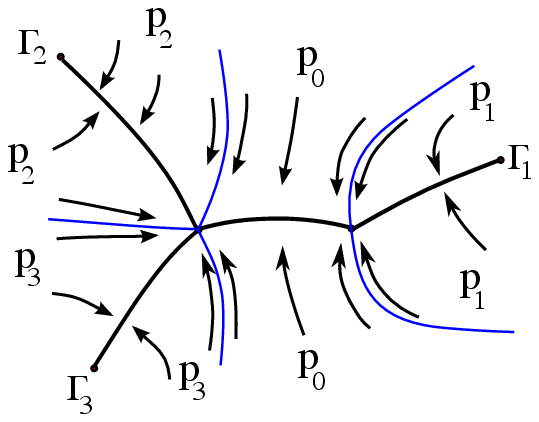}
\caption{ How the mappings $p_i$ (and hence $p_c$) act.
\label{f13-1d}}
\end{figure}

\begin{figure}[!h]  
\centering
\includegraphics[width=6cm]{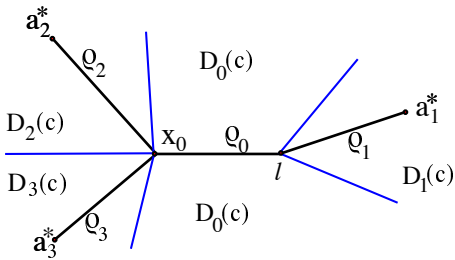}
\caption{ The same picture after a small change of variable.
\label{f13-1e}}
\end{figure}

Notice that the $D_i(c)$ cover $D(c)$, and also that for $i\in \{ 1, 2, 3 \}$,
$R_i \sm D \subset D_i(c)$.
The strange choice of a constant $1/3$ is to make sure (as in the picture,
and because the $\rho_i$ make large angles, as well as the directions of $\rho_1$, $\rho_2$,
and $\rho_3$) that $D_1(c)$ does not get close to $D_2(c) \cup D_3(c)$.
Because of this, we can find $5$-Lipschitz projections $p_i : D_i(c) \to \Gamma_i$
such that $|p_i(x)-x| \leq 2 \dist(x, \Gamma_i)$ for $x\in D_i(c)$, 
\begin{equation} \label{12.15}
p_i(x) = x_0 \ \text{ when $i \in \{ 2, 3, 0 \}$ and $x\in D_i(c) \cap D_j(c)$
for some other $j\in \{ 2, 3, 0 \}$},
\end{equation}
and 
\begin{equation} \label{12.16}
p_i(x) = \ell \ \text{ when $i \in \{ 0, 1 \}$ and $x\in D_0(c) \cap D_1(c)$.}
\end{equation}
Of course this would have been hard to arrange if $D_1(c) \cap D_0(c)$ had been
too close to $D_0(c) \cap (D_2(c) \cup D_3(c))$, but otherwise it is easy.
See Figure \ref{f13-1d} %
for a hint of what the desired projections should do, and Figure \ref{f13-1e}
for an equivalent model where the $p_i$ could be defined explicitly. 
We can also make sure that for $1 \leq i \leq 3$, $p_i(x) = a_i^\ast$
on $R_i \cap R_j$, where $j$ is any index $j\in \cI \sm \{ i \}$ such that
$\cL_i$ and $\cL_j$ share the endpoint $a_i^\ast$.

 % p109

 % p109

The mapping $p_c$ defined by $p_c(x) = p_i(x)$ for $x\in D_i(c)$ does the 
job. In particular, \eqref{12.12} can be arranged on each piece $D_i(c)$ separately,
and \eqref{12.11} is true because when $x\in D_i(c)$ and $y\in D_j(c)$ for some
$j \neq i$, the shortest path from $x$ to $y$ in $D(c)$ passes through boundaries
where the definitions coincides. That is, if this path $\gamma$ goes for instance
from $D_i(c)$ to some $D_k(c)$ to $D_j(c)$, it goes through points 
$z \in D_i(c) \cap D_k(c)$ and then $z'\in D_k(c) \cap D_i(c)$, and
\begin{eqnarray} \label{12.17}
|p(x)-p(y)| &=& |p_i(x) - p_j(y)| \leq |p_i(x) - p_i(z)|+ |p_k(z) - p_k(z')|+|p_j(z') - p_j(y)|
\nn\\&\leq& 5 \length(\gamma) \leq 10|x-y|
\end{eqnarray}
because $p_i(z)=p_k(z)$, $p_k(z') = p_j(z')$, and because the geometry of $D(c)$
is not that complicated.
As was announced earlier, the other cases are simpler; Lemma \ref{t12.2} follows.
\qed

\ms
At this point we have defined local projections $p_i$, $i\in \cI_1$
and $p_c$, $c\in CC_+$, and now we should glue them to make the mapping $p$ 
of Proposition \ref{t12.1}. The interesting part for the gluing will be near the points 
$\ell \in K \cap L$, where we want to attribute the points of $E \cap D_\ell$ 
to the various $c\in CC_+(\ell)$ (see below \eqref{12.9}).
For this we need a separation lemma.

\begin{lem} \label{t12.3}
We can find a small number $\tau_3 > 0$ and closed disjoint sets 
$T_c$, $c\in CC_+(\ell)$, such that
\begin{equation} \label{12.18}
c \subset T_c \subset E \cap D \ \text{ for } c\in CC_+(\ell), 
\end{equation}
\begin{equation} \label{12.19}
E \cap D \subset \bigcup_{c\in CC_+(\ell)} T_c \, ,
\end{equation}
and
\begin{equation} \label{12.20}
\dist(T_c, T_{c'}) \geq 10\tau_3 \ \text{ for } c, c' \in CC_+(\ell), c \neq c'.
\end{equation}
\end{lem}

\ms
Again, we can make no claim on the size of $\tau_3$. It may be extremely small, and 
it depends on $E$ and our earlier choice of radius for $\S$ (now normalized to be $1$).

This is mostly question of connectedness (say that each component $c$
is the intersection of the open and closed sets in $E\cap D$ that contain $c$),
but we shall cheat a little and use the fact that by \eqref{4.4}, 
\begin{equation} \label{12.21}
\H^1(E \cap D) \leq \H^1(E \cap \S) < +\infty.
\end{equation}
We try an argument by hands, with strings of small balls. 
Set $E_0 = E \cap D$.
For each integer $m>0$, we select a set $X(m) \subset E_0$ which is maximal under the constraint
that $|x-y| \geq 2^{-m}$ for $x, y \in X(m)$, $x\neq y$. Thus the balls
$\overline B(x,2^{-m})$, $x\in X(m)$, cover $E_0$. Set
\begin{equation} \label{12.22}
X_0(m) = \big\{ x\in X(m) \, ; \, B(x,2^{-m+1}) \text{ meets some } c \in CC_+(\ell) \big\},
\end{equation}
and $X_1(m) = X(m) \sm X_0(m)$.

Declare a point $x\in X_1(m)$ bad, or useless, if there is a radius $r \in (2^{-m}, 2^{-m+1})$
such that $E_0 \cap \d B(x,r) = \emptyset$. Denote by $X_b(m) \subset X_1(m)$ 
the set of bad points, and set $X_g(m) = X_1(m) \sm X_b(m)$.
Notice that if $x\in X_g(m)$, the radial projection $\pi_x$ defined by $\pi_x(z) = |z-x|$
maps $B(x,2^{-m+1})$ onto a set that contains $(2^{-m}, 2^{-m+1})$ 
(otherwise, $x\in X_b(m)$); hence
\begin{equation} \label{12.23}
2^{-m} \leq \H^1(\pi_x(E_0 \cap B(x,2^{-m+1}))) \leq \H^1(E_0 \cap B(x,2^{-m+1})).
\end{equation}
Because of this, the cardinality of $X_g(m)$ is
\begin{equation} \label{12.24}
\sharp(X_g(m)) \leq 2^m \sum_{x\in X_g(m)} \H^1(E_0 \cap B(x,2^{-m+1}))
\leq C 2^m \H^1(E_0)
\end{equation}
because the $B(x,2^{-m+1})$, $x\in X(m)$, have bounded covering. 

By ``$m$-string'', we shall mean a finite sequence of points $x_k \in X_0(m) \cup X_g(m)$, 
$0 \leq k \leq k_{max}$, such that, if we set 
$B_k = \ol B(x_k,2^{-m+2})$ for $0 \leq k \leq k_{max}$, we have that
\begin{equation} \label{12.25}
B_k \cap B_{k+1} \cap E_0 \neq \emptyset \ \text{ for } 0 \leq k < k_{max}.
\end{equation}

First suppose that for some choice of $\ell$ and $m$, there is no $m$-string as
above such that $B_0$ meets some component $c\in CC_+(\ell)$ and $B_{m_{max}}$ 
meets some other $c' \in CC_+(\ell)$. In this case, we can define the $T_c$ as follows.
For each $x\in X_b(m)$, choose $r \in (2^{-m}, 2^{-m+1})$ such that 
$E_0 \cap \d B(x,r) = \emptyset$ and set $B_x = B(x,r)$; then set
\begin{equation} \label{12.26}
E_b = E_0 \cap \bigcup_{x\in X_b(m)} B_x.
\end{equation}
This set is both open and closed in $E_0$, because each $B_x$ is. In fact,
for each $x\in X_b(m)$ there is a minuscule $\tau(x) > 0$ such that 
\begin{equation} \label{12.27}
\dist(E_0\cap B_x, E_0 \sm B_x) \geq \tau(x).
\end{equation}
Next, for $c\in CC_+(\ell)$, we denote by $T(c)$ the set of points $y\in E_0 \sm E_b$
that can be connected to $c$ by an $m$-string. 
This last means that we can find an $m$-string as above, such that $B_0$ 
meets $c$ and $B_{k_{max}}$ contains $x$. 
The sets $T(c)$, $c\in CC_+(\ell)$, are disjoint, because if $T(c)$ meets $T(c')$,
then there is an $m$-string that connects some point of $c$ to some point of $c'$.
They are also closed, because each $T(c)$ is in fact a finite union of sets 
$(E_0\sm E_b) \cap B_k$, and we made sure to take closed balls $B_k$.
Similarly, if we denote by $T_\infty$ the set of points of $E_0 \sm E_b$
that cannot be connected to any $c\in CC_+(\ell)$ by an $m$-string,
this set is also the union of the $\ol B(x,2^{m+2})$, $x\in X_g(m)$ that
meets it, and it is closed and disjoint of the others.

Finally, each $T(c)$ contains the corresponding $c$. Indeed, let $y\in c$ be given.
We know that the balls $\overline B(x,2^m)$, $x\in X(m)$, cover $E_0$, so we can find 
$x\in X(m)$ such that $y\in \overline B(x,2^m)$. Then $x\in X_0(m)$, by \eqref{12.22},
so it not bad. For the same reason, $y \notin E_b$, because no ball $B_x$ meets $c$.
We use the single $\ol B(x,2^{m+2})$ to connect $y$ to itself, and this shows that $y\in T(c)$.

Now set $T_c = T(c)$ for every component $c$ except one, and 
$T_c = T(c) \cup E_b \cup T_\infty$ for the last one. It is a little nicer to choose 
the special component $c_\ell = H_\ell$ as the last one, if $c_\ell \in CC_+(\ell)$,
because this way $p_c$ sends the whole $T_c$ to $\ell$. But really it does not matter.

The $T_c$ are disjoint and closed by construction. They cover $E_0$ by construction too,
and of course they lie at positive distances from each other, so \eqref{12.20} holds
for some $\tau_3 >0$. Thus the lemma holds in this case.

We are left with the case when for some $\ell$ and all choices of $m$,
we can find two different components $c, c' \in CC_+(\ell)$ that can be connected
by an $m$-string. Since $CC_+(\ell)$ has at most $4$ points, we may assume that
for a sequence of $m$ that goes to $+\infty$, the components $c$ and $c'$ are
the same. Choose an $m$-string that connects $c$ to $c'$, with a minimal number
of elements. Then the same ball $B_k$ does not appear twice in the sequence
(otherwise, drop all the intermediate balls), and similarly $B_0$ is the
only ball that meets $c$ and $B_{k_{max}}$ is the only $B_k$ that meets
$c'$. All the other $B_k$ are thus centered at points $x_k \in X_g(m)$,
and so $k_{max} \leq C 2^{m}+1$, by \eqref{12.24}.

For $0 \leq k \leq k_{max}-1$, connect $x_k$ to $x_{k+1}$ by a line segment.
This gives a curve $\Gamma_m$, that goes from $c$ to $c'$. Since
$|x_{k+1} - x_k| \leq 2^{-m+3}$ because $B_{k+1}$ meets $B_k$,
we get that $\length(\Gamma_m) \leq 8(C+2^{-m})$. Also, every point of
$\Gamma_m$ lies within $2^{-m+3}$ of $E_0$, because $B_k$ meets $E_0$.

We can parameterize $\Gamma_m$ with a mapping $z_m :[0,1] \to \Gamma_m$, 
in such a way that $z_m$ is $9C$-Lipschitz; then we can extract a sequence for which
the $z_m$ converge to a limit $z$, and $z([0,1])$ is a connected set in $E_0$
that goes from $c$ to $c'$. This contradiction with the fact that $c$ and $c'$
are different components proves that our second case does not happen, and 
Lemma~\ref{t12.3} follows.
\qed

\ms
We may now return to the construction of a global projection $p$
from the various $p_i$, $i\in \cI_1$ and $p_c$, $c\in CC_+$.
We now give a zone of influence to each $c\in CC_+$; for $i\in \cI_1$,
this was already done in \eqref{12.5}.

So fix $\ell \in K \cap L$ and $c\in CC_+(\ell)$. Set
\begin{equation} \label{12.28}
T_c^+ = \big\{ x\in D_c \, ; \, \dist(x,T_c) \leq 3\tau_3 \big\}
\ \text{ and } \ 
R_c = T_c^+ \cup \Big(\bigcup_{i\in \cI(c)} R_i \Big).
\end{equation}
Here $R_i$ is still defined by \eqref{12.5}, and $D_c$ is the disk $D$ associated to
the $\ell \in K \cap L$ such that $c\in CC_+(\ell)$. We should be able to avoid confusion 
between $R_i$, with $i\in \cI$ and the larger $R_c$, $c\in CC_+$.
The definition may look a little strange, but away from $D$, we are happy to keep 
$\bigcup_{i\in \cI(c)} R_i$, and not more (to avoid complications with the gluing), 
and in $D$ it is better to add a small neighborhood of $T_c$, because we want to cover 
a small neighborhood of $E\cap \S$.
Our domain of definition will be 
\begin{equation} \label{12.29}
R_+ = \Big(\bigcup_{i\in \cI_1} R_i \Big) \cup \Big(\bigcup_{c\in CC_+} R_c \Big)
= \Big(\bigcup_{i\in \cI} R_i \Big) \cup \Big(\bigcup_{c\in CC_+} T_c^+ \Big)
\end{equation}
and we want to set
\begin{equation} \label{12.30}
p(x) = p_i(x) \text{ for $i\in \cI_1$ and $x\in R_i$,} 
\end{equation}
and
\begin{equation} \label{12.31}
p(x) = p_c(x) \text{ for $c\in CC_+$ and $x\in R_c \subset D(c)$,}
\end{equation}
where the inclusion is easy (compare \eqref{12.28} with the first line of Lemma \ref{t12.2})
and implies that $p_c(x)$ is defined.
We need to check that all this is compatible, and produces a Lipschitz
function. We will cut $R_+$ into three regions that overlap, and first check things
on each one.

We start with $R_+(1) = \bigcup_{i\in \cI} R_i$ (where we also include $\cI_0$).
Let us first check that
\begin{equation} \label{12.32}
\text{if $i, j \in \cI$ are such that $\dist(R_i,R_j) \leq \frac{\tau}{10}$, then 
$\cL_i$ and $\cL_j$ have a common endpoint $a^\ast$.} 
\end{equation}
Suppose $i \neq j$ and $\dist(R_i,R_j) \leq 10^{-1}\tau$. 
Then $\dist(\cL_i,\cL_j) \leq \tau$ by \eqref{12.5}, $\dist(\cC_i,\cC_j) \leq 2\tau$
by \eqref{5.39}, and by \eqref{5.4} and \eqref{5.2bis} $\cC_i$ meets $\cC_j$.
This means that they have a common endpoint, which we call $a$.

We can say a bit more. 
If $i$ or $j$ lies in $\cI_1$, then $\dist(a,K\cap L) \geq 10\eta(X)$
by \eqref{3.10} and \eqref{3.12}, and Proposition \ref{t5.4} says that
$\cL_i$ and $\cL_j$ have a common endpoint $a^\ast$, with
$|a^\ast-a| \leq 10^{-9}\tau$ (see \eqref{5.38}).
Notice that the case when $i, j$ lie in $\cI_0$ and end at the same $\ell$
does not arise, because in this case  
\begin{equation} \label{12.33}
\dist(R_i,R_j) \geq \dist(\cL_i,\cL_j) - 2 \cdot 10^{-1} \tau
\geq 8 \cdot 10^{-1} \tau
\end{equation}
because the $\cL_j$ do not get inside $D$, start from points $c_i^\ast$ and
$c_j^\ast$ such that $|c_i^\ast-c_j^\ast| \geq \tau$, and go away in the direction
opposite to $\ell$. For a proof, use \eqref{5.41} and the fact that $\cC_i$ and $\cC_j$
make angles of at least $120^\circ$.
The case when $i,j \in \cI_0$ but come from different $\ell \in K\cap L$ goes like
when $i$ or $j$ lies in $\cI_1$; so \eqref{12.32} holds.

Because of \eqref{12.32} and our precautions \eqref{12.8} and below \eqref{12.12},
not only $p_i(x)=p_j(x)$ when $x\in R_i \cap R_j$, but the proof of \eqref{12.9} shows that
the mapping $p$ on $R_+(1)$ that we construct in this way is locally Lipschitz, in the sense that
\begin{equation} \label{12.34}
\text{$p$ is $30$-Lipschitz on $R_+(1) \cap B$ for every ball $B$ of radius $10^{-2}\tau$.} 
\end{equation}
Next we pick $\ell \in K\cap L$ and consider 
\begin{equation} \label{12.35}
R_+(2,\ell) = R_+ \cap A, \ \text{ with } A_0 = B(\ell,2\tau) \sm B(\ell,2\tau/3).
\end{equation}
Let us apply Proposition \ref{t5.4}, but with the smaller constant $\tau' = \tau/3$.
We get a nice description of $E \cap \S$ in the complement of 
$B(\ell,\tau') \cup B(-\ell,\tau')$, but we only care about what happens on the 
annulus $A = B(\ell,3\tau) \sm B(\ell,\tau/3)$. We get that
\begin{equation} \label{12.36}
E \cap \S \cap A = \bigcup_{i\in \cI(\ell)} \cL'_i \cap A,
\end{equation}
where $\cI(\ell)$ is the set of concerned indices $i$, i.e., those for which $\ell \in \cC_i$,
and the $\cL_i'$ are nice $C^1$ curves that go from $\d B(\ell,3\tau)$ to $\d B(\ell,\tau/3)$.

On $A \sm B(\ell,\tau)$, we have two representations of the same set 
$E \cap A\sm B(\ell,\tau)$, given by applications of Proposition \ref{t5.4} 
with different values of $\tau$, but which must coincide anyway. Thus the
$\cL'_i$ coincide with the $\cL_i$ on $A \sm B(\ell,\tau)$.
We may assume that we chose the labels correctly, so that in fact 
\begin{equation} \label{12.37}
\cL_i \cap A \sm B(\ell,\tau) = \cL'_i \cap A \sm B(\ell,\tau)
\ \text{ for } i\in \cI(\ell).
\end{equation}
Set $\gamma_i = \cL'_i \cap A \cap D$. By Proposition \ref{t5.4},
$\gamma_i$ is a $C^1$ curve that starts at $c_i^\ast$
(the only point of $\cL'_i \cap \d B(\ell, \tau) = \cL_i \cap \d B(\ell, \tau)$) 
and goes to $\d B(\ell,\tau/3)$.
Since each $\cL'_i$ stays close to the corresponding $\cC_i$, we also have that
\begin{equation} \label{12.38}
\dist(\gamma_i,\gamma_j) \geq \tau/3 \ \text{ for } i, j \in \cI(\ell), i\neq j.
\end{equation}

For each $i\in \cI(\ell)$, $c_i^\ast$ lies in some component $c \in CC(\ell)$, 
which we call $c(i)$. Then
\begin{equation} \label{12.39}
\gamma_i \cap D \subset c(i) \subset T_{c(i)}
\end{equation}
because $\cL'_i \cap A \cap D$ is a connected subset of $E \cap D$ 
that contains $c_i^\ast$, and by \eqref{12.18}. Since the $T_c$ are disjoint
and contained in $E$, and the $\gamma_i$ already cover $E \cap A \cap D$
(by \eqref{12.36}), we see that the only $T_c$ that meet $A\cap D$ are
the $T_{c(i)}$, and in addition, for each $i$
\begin{equation} \label{12.40}
T_{c(i)} \cap A\cap D = \bigcup_{j \in \cI(\ell) ; c(j)=c(i)} \gamma_j.
\end{equation}
Recall from \eqref{12.28} that $T_{c}^+$ is just a $3\tau_3$-neighborhood
of $T_c$ in $D$. Then the only $T_c^+$ that meet $A_0\cap D$ are the $T_{c(i)}^+$,
and each $T_{c(i)}^+$ is just the $3\tau_3$-neighborhood in $A_0\cap D$ of the 
$\bigcup_{j \in \cI(\ell) ; c(j)=c(i)} \gamma_j$.

Now $R_+(2,\ell)$ is the union of at most three pieces $R_+(2,\ell,i)$, $i\in \cI(\ell)$,
where each $R_+(2,\ell,i)$ is composed of $R_i$, plus the 
$3\tau_3$-neighborhood in $A_0\cap D$ of $\gamma_i$. 
Each $R_+(2,\ell,i)$ is contained in a single $R_c$ (and meets no other),
hence $p$ is well defined on $R_+(2,\ell,i)$ by \eqref{12.31}, and $10$-Lipchitz 
by Lemma \ref{t12.2}. In addition 
\begin{equation} \label{12.41}
\dist(R_+(2,\ell,i),R_+(2,\ell,j)) \geq \tau/4
\end{equation}
by \eqref{12.38}, so $p$ is also Lipschitz on their union $R_+(2,\ell)$.

We turn to our last sets 
\begin{equation} \label{12.42}
R_+(3,\ell) = R_+ \cap B(\ell, 4\tau/5) \subset \bigcup_{c\in CC_+} T_c^+
\end{equation}
by the second part of \eqref{12.29} 
and because the $R_i$, $i\in \cI$, never go that far inside $D$ (since the
$\cL_i$ don't meet $B(\ell,\tau)$). By \eqref{12.20} and \eqref{12.28},
\begin{equation} \label{12.43}
\dist(T_c^+,T_{c'}^+) \geq 4\tau_3  \ \text{ when } c \neq c' ,
\end{equation}
and $p$ is well defined and $10$-Lipchitz on each $T_c^+$
(by Lemma \ref{t12.2}), so $p$ is well defined and Lipschitz on $R_+(3,\ell)$,
and $10$-Lipschitz on each open ball of radius $2\tau_3$.

\ms
At this point we have a coherent definition of $p$ on $R_+$,
and proved Lipschitz bounds for $p$ on the various pieces that compose $R_+$.
These pieces have sufficient overlap, so we get the local Lipschitz property \eqref{12.3}
required for Proposition \ref{t12.1}. Then $p$ is automatically Lipschitz on $R_+$,
although perhaps only with the very bad norm $\tau_3^{-1}$: if $x,y \in R_+$,
either $|x-y| < 2\tau_1$ and then $|p(x)-p(y)|\leq 20|x-y|$, or else
$|p(x)-p(y)| \leq 2 \leq \tau_1^{-1}|x-y|$.

Next we check that the domain $E_+$ promised in Proposition \ref{t12.1}
is contained in $R_+$. Let $z\in E_+$ be given, and pick $x\in E \cap \S$
such that $|x-z| \leq \tau_3$. If $x$ lies in a disk $D_\ell = \S \cap B(\ell, \tau)$,
then \eqref{12.19} says that it lies in some $T_c$, and $z\in T_c^+$ by \eqref{12.28},
unless by bad luck $z$ falls outside of $D_\ell$. But if this 
happens, 
$\dist(x, \d B(\ell,\tau)) \leq \tau_3$, \eqref{12.36} and \eqref{12.37} say that
$\dist(x,\cL_i) \leq 2\tau_3$ for some $i\in \cI(\ell)$, and then $z \in R_i$.
In both cases, $z\in R_+$ (see \eqref{12.29}). The other case is when
$x$ lies in no $D_\ell$. Then it lies very close to some $\cL_i$
(by Proposition~\ref{t5.4}), and then $z\in R_i \subset R_+$
by \eqref{12.5} and \eqref{12.29}. So $E_+ \subset R_+$.

Next we check \eqref{12.2}. Let $z\in R_+$ be given.
When $z\in R_i$ for some $i\in \cI_1$,
\eqref{12.2} follows from \eqref{12.30} and \eqref{12.7}.
Otherwise, $z\in R_c$ for some $c\in CC_+$, and 
\begin{equation} \label{12.44}
|p(z)-z| = |p_c(z)-z| \leq 10 \dist(z,\Gamma_c)
\end{equation}
by \eqref{12.31} and \eqref{12.12}. 
Let $\ell$ be such that $c\in CC_+(\ell)$. If $z\in T_c^+$, then
$z\in D_\ell$ (see \eqref{12.28}), and $\dist(z,\Gamma_c) \leq 2\tau$
because every $\Gamma_c$ contains at least a point in $D_\ell$.
In this case \eqref{12.2} follows from \eqref{12.44}.

By \eqref{12.28}, we are left with the case when $z\in R_i$ for some $i\in \cI(c)$,
and we still want to evaluate $\dist(z,\Gamma_c)$.
Let $\ell$ be such that $c \in CC_+(\ell)$; notice that in fact $c\in CC(\ell)$;
the special components $c_\ell$ that were artificially added don't come with 
a set $\cI(c)$.

We shall now use the fact that for all the components $c$ such that $i\in \cI(c)$,
$\Gamma_c$ contains a small Lipschitz graph $\Gamma$, over some geodesic
$\rho = \rho(a_i^\ast,x)$, and where the other endpoint $x$ lies in $D$.
This is why we did not remove $\cL_i$ in Configuration H, for instance.

We want to see where $\Gamma$ is localized.
Recall that $\Gamma$ was constructed by applying Section \ref{S7}
to a curve $\gamma$ with the same endpoints. There are a few ways in which $\gamma$
was chosen, depending on the configuration, but in all the cases $\gamma$ was
contained in $\cL_i \cup D$. In the algorithm of Section \ref{S7}, $\Gamma$ is obtained
from $\gamma$ by replacing some of its sub-arcs with the geodesics with the same endpoints;
because of this,
\begin{equation}\label{12.45}
\Gamma \subset {\rm Hull}(\gamma) \subset {\rm Hull}(\cL_i \cup D),
\end{equation}
where the convex hulls ${\rm Hull}(\gamma)$ and ${\rm Hull}(\cL_i \cup D)$ are defined in terms of
geodesics in $\S$. There is no ambiguity about geodesics, because we shall see that
$\cL_i \cup D$ stays quite close to $\cC_i$, which is a geodesic of length at most $\pi/2$.
More precisely, \eqref{5.41} says that $\dist(x,\cC_i) \leq 10^{-8}\tau$ for $x\in \cL_i$,
and since $\ell$ is an endpoint of $\cC_i$, we deduce from \eqref{12.45} that
\begin{equation}\label{12.46}
\dist(x,\cC_i) \leq 2\tau \ \text{ for } x\in \Gamma.
\end{equation}
Let us check that this implies that
\begin{equation}\label{12.47}
\dist(x,\Gamma) \leq 5\tau \ \text{ for } x\in \cC_i.
\end{equation}
Since $\Gamma$ starts from $a_i^\ast$ (very close to the endpoint $a_i$ of $\cC_i$,
and ends in $D$, we can assume that $\tau < |x-\ell| < |a_i^\ast-\ell|$. Since
$\Gamma$  is connected, we can find $y\in \Gamma$ such that $|y-\ell| =  |x-\ell|$,
and by \eqref{12.46} a point $z\in \cC_i$ such that $|z-y| \leq 2\tau$.
Then $\big| |z-\ell| - |x-\ell| \big| = \big| |z-\ell| - |y-\ell| \big| \leq 2\tau$, and hence
$|z-x| \leq 3\tau$ because $x$ and $z$ both lie on the geodesic $\cC_i$ that starts from $\ell$. 
Finally $|y-x| \leq |y-z|+|z-x| \leq 5\tau$,
as needed for \eqref{12.47}.

We may now return to $z\in R_i$ and chase points.
By \eqref{12.5}, there is a point $z_1\in \cL_i$ such that $|z_1-z| \leq 10^{-1}\tau$. 
By \eqref{5.41}, we can find $z_2 \in \cC_i$ such that $|z_2-z_1| \leq 10^{-8}\tau$. 
By \eqref{12.47}, we can find $z_3 \in \Gamma$ such that $|z_3-z_2| \leq 5\tau$.
Since $\Gamma \subset \Gamma_c$, we get that
\begin{equation}\label{12.48}
\dist(z,\Gamma_c) \leq \dist(z,\Gamma) \leq |z-z_3| \leq 6\tau,
\end{equation}
and \eqref{12.2} follows from \eqref{12.44} in this last case as well.

\ms
Finally we need to check \eqref{12.4}. When $\ell \in c$ for some $c\in CC(\ell)$,
we made sure to keep $\ell \in \Gamma_c$. Then $p_c(\ell) = \ell$ by \eqref{12.12},
and \eqref{12.4} follows from \eqref{12.31}. Otherwise, we added a special component 
$c_\ell = \H_\ell$ to $CC_+(\ell)$, and took $p_{c_\ell}(z) = \ell$ for all $z$,
in particular $z=\ell$.
This completes our proof of Proposition \ref{t12.1}.
\qed

\section{Our first competitor and the contribution from the thin gluing annulus}
\label{S13}

We now have a net $\gamma^\ast$ of curves in $E\cap \S$, another net $\Gamma^\ast$
of Lipschitz graphs, and (by Proposition \ref{t12.1}) a projection $p$ from a neighborhood 
$E_+$ of $E\cap \S$ to the net $\Gamma^\ast$. 
We want to use these to construct a first competitor for $E$.
We use the following lemma to choose another very small number $\tau_4 > 0$.

\begin{lem} \label{t13.1}
Set 
\begin{equation}\label{13.1}
A(t) = \overline B(0,1) \sm B(0, 1-t) \ \text{ for } 0 < t < 10^{-1}.
\end{equation}
We can find $\tau_4 > 0$ such that
\begin{equation}\label{13.2}
\frac{x}{|x|} \in E_+(\tau_3) \ \text{ for } x \in E \cap A(2\tau_4).
\end{equation}
\end{lem}

\ms
Here $E_+(\tau_3)$ is defined by \eqref{12.1}, and $\tau_3$ was chosen in 
Proposition \ref{t12.1}.
Of course we don't get any uniform control on $\tau_4$; we did not even
get a uniform control on $\tau_3$.

The proof is easy. If we could not find $\tau_4$, we would be able to find a sequence of 
points $x_k \in E \cap A(2^{-k})$, that tends to a limit $x_\infty$, but so that 
$\frac{x_k}{|x_k|}$ stays at distance at least $\tau_3$ from $E \cap \S$.
This is impossible because $x_\infty \in E \cap S$. 
\qed

\ms
Let $\tau_4$ satisfy the conclusion of the lemma, and take any 
$\sigma \in (0,\tau_4)$; we see $\sigma$ as a small parameter that we may
chose later.

Extend $p$ so that it is homogeneous of degree $0$.  That is, set
\begin{equation}\label{13.3}
p(x) = p(x/|x|) \ \text{ when } \frac{x}{|x|} \in E_+.
\end{equation}
Then, by \eqref{13.2}, $p$ is well defined (and Lipschitz with a bad norm)
on $E \cap A(2\tau_4)$.

We are ready to define a new competitor for $E$, which we write as 
\begin{equation} \label{13.4}
F^0 = \varphi^0(E),
\end{equation}
for some $\varphi^0$ that will be defined soon. The main part of $F^0$
will be a subset of the cone over $\Gamma^\ast$. We will not be finished yet,
$F^0$ will need to be further improved. First we set
\begin{equation}\label{13.5}
\varphi^0(x) = x \ \text{ for } x\in E \sm B(0,1).
\end{equation}
On the exterior part of $B(0,1)$, we use $p$ to contract reasonably slowly on $\Gamma^\ast$. 
That is, we set
\begin{equation}\label{13.6}
\varphi^0(x) = \frac{|x| + \sigma - 1 }{ \sigma}\, x + \frac{1 - |x| }{ \sigma}\, p(x)
\ \text{ for } x \in E \cap A(\sigma).
\end{equation}
Notice that this makes sense because $p(x)$ is well defined there, and also that
the two definitions yield $\varphi^0(x) = x$ on $\S$. 
On the other sphere, 
\begin{equation} \label{13.7}
\varphi^0(x) = p(x) \in \Gamma^\ast \subset \S \ \text{ for } x\in E \cap \d B(0,1-\sigma).
\end{equation}

Now we contract very brutally along the cone over $\Gamma^\ast$. Set
\begin{equation}\label{13.8}
\varphi^0(x) = \frac{|x| + 2\sigma -1 }{ \sigma}\, p(x)
\ \text{ for } x \in E \cap A(2\sigma) \sm A(\sigma).
\end{equation}
Again this is continuous across $\d B(0,1-\sigma)$, and 
$\varphi^0(x) = 0$ on $\d B(0,1-2\sigma)$. Thus we can safely take
\begin{equation}\label{13.9}
\varphi^0(x) = 0 \ \text{ for } x\in E \cap B(0,1-2\sigma).
\end{equation}
This gives a Lipschitz mapping $\varphi^0$ defined on $E$;
its Lipschitz constant depends on $\sigma$, $\tau_4$, and $\tau_3$ and may be really huge, 
so we will be careful not to use this directly in the estimates. 
Since we like to define competitors in terms of deformations, we are also led to set
\begin{equation}\label{13.10}
\varphi^0_t(x) = (1-t) x + t \varphi^0(x) \ \text{ for $x\in E$ and } t\in [0,1],
\end{equation}
and check that 
\begin{equation} \label{13.11}
\text{the $\varphi_t^0$ define an acceptable deformation for $E$ in $\ol B(0,1)$,} 
\end{equation}
as in Definition \ref{t1.1}. As often, \eqref{1.4} and \eqref{1.5}
are trivial, \eqref{1.6} holds because $|p(x)| \leq 1$ and the unit ball is convex,
\eqref{1.8} holds because $\varphi^0$ is Lipschitz, and the only interesting piece
is the boundary condition \eqref{1.7}.

Let $x\in E \cap L$ be given. 
If $|x| \geq 1$, $\varphi_t(x) = x$ by \eqref{13.5} and \eqref{13.10},
and $\varphi_t^0(x) \in L$ trivially. If $|x| < 1-2\sigma$, $\varphi_0(x) = 0$
and hence $\varphi_t^0(x) = (1-t)x \in L$. So let us assume that
$|x| \geq 1-2\sigma \geq 1/2$. 

Set $\ell = x/|x| \in L$. We claim that $\ell \in K$ too. Indeed, otherwise
\eqref{3.10} and \eqref{3.12} say that $\dist(\ell, K) \geq \eta_L(X) \geq 10 \eta(X)$,
but yet $\dist(x,X) \leq 2\varepsilon$ by \eqref{4.3}, hence 
$\dist(\ell,X) \leq 2|x|^{-1} \varepsilon \leq 5\varepsilon$, a contradiction.
So $\ell \in K\cap L$, $p(\ell) = \ell$ by \eqref{12.4}, and the various formulae
yield $\varphi_t^0(x) \in L$; \eqref{13.11} follows.

\ms
It is amusing that the very brutal part \eqref{13.8} works so well.
We like it because it allows us to concentrate on the set $E \cap A(2\sigma)$,
and essentially disregard any bad behavior that $E$ may have in a smaller
ball. Of course we will still need to know that $E$ is nice on the thin annulus $A(\sigma)$,
and we shall get part of this with a maximal function argument.

The main part of $F^0 \cap B(0,1)$ is contained in 
\begin{equation} \label{13.12}
\Sigma_F(\Gamma^\ast) = \big\{ tx \, ; \, t\in [0,1] \text{ and } x\in \Gamma^\ast \big\},
\end{equation}
the (truncated) cone over $\Gamma^\ast$. Indeed,
\begin{equation} \label{13.13}
\varphi^0(E \cap B(0,1-\sigma)) 
= \{ 0 \} \cup \varphi^0(E \cap B(0,1-\sigma) \sm B(0,1-2\sigma))
\subset \Sigma_F(\Gamma^\ast)
\end{equation}
by \eqref{13.9}, because $p(E\cap A(2\sigma)) \subset p(E_+) \subset \Gamma^\ast$ 
(by Proposition \ref{t12.1}), and then by \eqref{13.8}.

\ms
In the rest of this section, we control the remaining piece of $\varphi^0(E \cap \ol\B)$, 
which is the set
\begin{equation} \label{13.14}
F(\sigma) = \varphi^0(E \cap A(\sigma)).
\end{equation}
Once this is done, we shall still want to improve on the cone $\Sigma_F(\Gamma^\ast)$, 
and construct other competitors. But we shall be able to use the next estimates 
on $F(\sigma)$ for those too.

We shall leave the dependence on $\sigma$ explicit in estimates, because we shall need to check
that some of our estimates do not depend on $\sigma$,
but we set $A = A(\sigma)$ to save some space. 
We want to estimate
\begin{equation} \label{13.15}
M(\sigma) = \H^2(F(\sigma)).
\end{equation}

In next lemma we use some of the additional properties of our radius $r=1$ that 
we required in Section \ref{S4}.

\begin{lem}\label{t13.2} 
If $\H^1(E \cap \S) < +\infty$ and the assumption \eqref{4.8} holds for $r=1$, then
\begin{equation}\label{13.16}
\limsup_{\sigma \to 0} M(\sigma) \leq C \int_{E \cap \S} \dist(x, p(x)) d\H^1(x).
\end{equation}
\end{lem}

\ms
This lemma is essentially measure-theoretic; 
then we shall estimate the right-hand side of \eqref{13.16}, 
and this will use the construction of $p$.

Before we prove this, let us explain roughly why it may be true.
We shall use the area formula to estimate $M(\sigma)$,
but the point is that $\varphi^0(E \cap A)$ is like a curtain, composed 
from all the segments $[x,p(x)]$; thus \eqref{13.16} looks a little like Fubini's theorem.

The proof is not very complicated, but since it is also done in [C1] (see (9.46) there and its proof), 
we only give the great lines. First, we use the rectifiability of $E$ and the area theorem
(Corollary 3.2.20 in [Fe]) %%
to write
\begin{equation}\label{13.17}
M(\sigma) \leq \int_{E \cap A} J_{\varphi^0}(x) d\H^2(x),
\end{equation}
where $J_{\varphi^0}(x)$ is the Jacobian of the approximate differential $D_{\varphi^0}(x)$ 
of $\varphi^0$ along $E$, which is defined for $\H^2$-almost every $x\in E$. 
Then we estimate the size of $D_{\varphi^0}(x)$ on an orthonormal basis of $(v,w)$ 
of the approximate tangent plane to $E$ at $x$. We choose $(v,w)$
so that $v$ is orthogonal to the radial direction $[0,x]$; then \eqref{13.6} and our local
Lipschitz estimate \eqref{12.3} for $p$ yield $|D_{\varphi^0}(x) \cdot v| \leq C$ 
(we may assume that $|x| \geq 1/2$ so that \eqref{13.3} is tame, and in the direction
of $v$ the differential of the radial cut-off function in \eqref{13.6} vanishes). 
In the direction $w$, we get the estimate
\begin{equation}\label{13.18}
|D_{\varphi^0}(x) \cdot w| \leq C + \sigma^{-1} \cos\theta(x) |p(x)-x|,
\end{equation}
where $\theta(x)$ is the angle of $w$ with the radial direction, or equivalently
$\cos\theta(x) = |\langle w, \frac{x}{|x|} \rangle|$.
Then
\begin{equation}\label{13.19}
J_{\varphi^0}(x) \leq |D_{\varphi^0}(x) \cdot v| \, |D_{\varphi^0}(x) \cdot w| 
\leq C + C \sigma^{-1} \cos\theta(x) |p(x)-x|.
\end{equation}
Then we apply the coarea theorem (3.2.22 in [Fe]), to the mapping $h$ 
defined by $h(x) = |x|$, integrated against the continuous function $x \to |p(x)-x|$, 
and get that
\begin{equation}\label{13.20}
\int_{E \cap A} |p(x)-x| J_h(x) d \H^2(x) 
= \int_{t \in (1-\sigma,1)}\int_{E \cap \d B(0,t)} |p(x)-x| d\H^1(x) dt,
\end{equation}
with the one-dimensional jacobian $J_h = |D_h(x)\cdot w| = \cos\theta(x)$. 
Thus by \eqref{13.17}, \eqref{13.19} and \eqref{13.20},
\begin{equation}\label{13.21}
\begin{aligned}
M(\sigma) &\leq C \int_{E \cap A} \big[1 + \sigma^{-1} \cos\theta(x) |p(x)-x| \big] d\H^2(x)
\\&\leq C \H^2(E \cap A) 
+ C \sigma^{-1} \int_{t \in (1-\sigma,1)}\int_{E \cap \d B(0,t)} |p(x)-x| d\H^1(x) dt.
\end{aligned}
\end{equation}
Now we let $\sigma$ tend to $0$. Notice that $\H^2(E \cap \S) = 0$ because
$\H^1(E \cap \S)<+\infty$. Next $\H^2(E \cap B(0,1) \sm B(0,1-\sigma))$
tends to $0$, because $\H^2(E \cap B(0,1)) < +\infty$ and by the monotone convergence
theorem (or the definition of a measure). 
Thus $\H^2(E \cap A(\sigma))$ tends to $0$. 
The other term in \eqref{13.21} tends to $\int_{E \cap \d B(0,1)} |p(x)-x| d\H^1(x)$, 
by our special assumption \eqref{4.8}; 
\eqref{13.16} and Lemma \ref{t13.2} follow.
\qed

\ms
We now estimate the right-hand of \eqref{13.16}, and proceed as in the end
of Section 9 in [C1]. % near 9.57 p60
Notice that we are very happy that we have to estimate an integral on $E \cap \S$ 
(as opposed to an annulus), because this is precisely the place that we control. 
We know from \eqref{12.2} that
\begin{equation}\label{13.22}
\dist(x, p(x)) \leq 60\tau \ \text{ for } x\in E \cap \S,
\end{equation}
but there are lots of points of $E\cap \S$ for which $p(x)=x$, and which we can 
take out of the estimates. First we check that
\begin{equation}\label{13.23}
p(x)=x \text{ for } x \in \bigcup_{i\in \cI_1} \cL_i.
\end{equation}
Indeed, $x\in R_i$ (see the definition \eqref{12.5}), and 
$p(x) = p_i(x) = x$ by \eqref{12.30} and \eqref{12.7}.
Next we claim that
\begin{equation}\label{13.24}
p(x)=x \text{ for } x \in \bigcup_{c\in CC} \big(\ol\gamma_c \cap \Gamma_c \big).
\end{equation}
Recall that $\ol\gamma_c = \gamma_c \cup \big(\bigcup_{i\in \cI(c)} \cL_i \big)$
(we add the curves $\cL_i$ that touch the points of $c\cap \d D$).
The arcs $\cL_i$ above are contained in $R_c$, by \eqref{12.28}, 
so $p(x) = p_c(x) = x$ when $x\in \cL_i \cap \Gamma_c$,
by \eqref{12.31} and \eqref{12.12}. 

We are left with $x\in \gamma_c \cap \Gamma_c$.
Let us recall why $\gamma_c \subset c$. In Sections \ref{S9}-\ref{S11}, 
when we constructed the nets $\Gamma = \Gamma_c$, we always started from a set 
$\gamma \subset E \cap D$. This set was connected and contained at least one point 
of $\d D$; this is how we defined the different configurations. Then $\gamma$
is contained in the component $c$ (often called $H_i$) that contains any point of 
$\gamma \cap \d D$. So $\gamma_c \subset c$. 
In addition, $c\subset T_c \subset T_c^+ \subset R_c$
by \eqref{12.18} and \eqref{12.28}. Then $p(x) = p_c(x)$ by \eqref{12.31},
and since $x\in \Gamma_c$, \eqref{12.44} says that $p(x) = x$.
So \eqref{13.24} holds. Thus, in the set
\begin{equation} \label{13.25}
\gamma^\ast = \Big(\bigcup_{i\in \cI_1} \cL_i \Big)\cup \Big(\bigcup_{c\in CC} \ol\gamma_c\Big),
\end{equation}
the only part that remains is $\bigcup_{c\in CC} [\ol\gamma_c \sm \Gamma_c]$.
For the rest of $E \cap \S$, we don't try anything, and just keep the set 
$E \cap \S \sm \gamma^\ast$.

Finally, it follows from \eqref{13.16} and \eqref{13.22} that
\begin{equation}\label{13.26}
\limsup_{\sigma \to 0} M(\sigma) \leq C \tau \sum_{c\in CC} \H^1(\ol\gamma_c \sm \Gamma_c)
+C\tau \H^1(E \cap \S \sm \gamma^\ast).
\end{equation}
We are reasonably happy about this. We consider $M(\sigma)$ as a loss in the estimates,
and we expect it to be compensated by larger wins. In this respect,
$E \cap \S \sm \gamma^\ast$ is a part of $E\cap \S$ that we just dropped 
to get $\gamma^\ast$, so we will save much more than 
$C\tau \H^1(E \cap \S \sm \gamma^\ast)$ by removing it from the picture
(when we construct cones), and similarly $\H^1(\ol\gamma_c \sm \Gamma_c)$
is controlled by \eqref{9.7}, and the corresponding term of \eqref{13.26}
will be compensated by a win of area when we replace the cone $\Sigma_F$ over
$\Gamma^\ast$ with a bunch of harmonic graphs.

\section{A second competitor build with harmonic graphs}
\label{S14}

The competitor $F^0$ that we constructed in the previous section was just
a first attempt, which still looks a little like the cone over $E\cap \S$. 
It is better in some sense, because we got rid of $E \cap \S \sm \gamma^\ast$,
but the advantage of replacing $\gamma^\ast$ with $\Gamma^\ast$ is not
clear yet.

In this section we construct our first serious competitor, obtained from $F^0$
by replacing parts of $\Sigma_F(\Gamma^\ast)$ by better surfaces constructed in Section \ref{S8}.

Recall that $\Gamma^\ast$ is the union of the Lipschitz nets $\Gamma$ that were constructed
in Sections~\ref{S9}-\ref{S11}. Since the construction was done by connected components in
configurations, our best description so far is that 
\begin{equation}\label{14.1}
\Gamma^\ast 
= \Big(\bigcup_{i\in \cI_1} \Gamma_i \Big)\cup \Big(\bigcup_{c\in CC_+} \Gamma_c\Big),
\end{equation}
where the notation is the same as in Section \ref{S12} (see below \eqref{12.9}),
and where, if we use Remark \ref{r6.3n}, we managed to take $\Gamma_i = \cL_i$ for $i\in \cI_1$.
But we could also have kept things the way they were at the beginning, but replaced $\cL_i$
by a small Lipschitz $\Gamma_i$ with the same endpoints, obtained as in Section \ref{S7}.

To save some space, we condense \eqref{14.1} into
\begin{equation}\label{14.2}
\Gamma^\ast = \bigcup_{c\in \cI_1 \cup CC_+} \Gamma_c.
\end{equation}
Then the main part of $F^0 \cap \B$ is the cone
\begin{equation}\label{14.3}
\Sigma_F(\Gamma^\ast) = \bigcup_{c\in \cI_1 \cup CC_+} \Sigma_F(\Gamma_c),
\end{equation}
where $\Sigma_F$ is our standard notation for a cone. Thus, as in \eqref{13.12},
\begin{equation} \label{14.4}
\Sigma_F(\Gamma_c) = \big\{ tx \, ; \, t\in [0,1] \text{ and } x\in \Gamma_c \big\}.
\end{equation}
Recall that the sets $\Gamma_c$ are disjoint, except for common endpoints at vertices $a_i^\ast$.

Now we want a finer decomposition of $\Gamma^\ast$ into single Lipschitz curves, 
which we shall write as 
\begin{equation}\label{14.5}
\Gamma^\ast = \bigcup_{j\in J^\ast} \Gamma_j
\end{equation}
for some new set of indices $J^\ast$. Let us say how we do it, so as not to create
too much confusion. When $i \in \cI_1$, we have a single curve $\Gamma_i = \cL_i$; 
we keep it as it is, just put the set of indices $\cI_1$ in $J^\ast$
and keep the same curves with the same names. We could also do this when $c\in CC$
and $\Gamma_c$ is composed of a single curve, but let us not bother. Instead, for $c\in CC$,
we observe that $\Gamma_c$ us composed of at most $4$ small Lipschitz graphs $\Gamma_j$, 
as in the description of Section \ref{S9}, and write this as
\begin{equation} \label{14.6}
\Gamma_c = \bigcup_{j\in J(c)} \Gamma_j.
\end{equation}
We also include (the elements of) $J(c)$ in the index set $J^\ast$. 
Finally, there is the case of the special components $c_\ell \in CC_+ \sm CC$. 
If $c_\ell$ is such a component, we took $\Gamma_c = \{ \ell \}$, it is a single 
degenerate curve, and we also put it in our bag $J^\ast$ with the same name. 
Thus our new set of indices is
\begin{equation}\label{14.7}
J^\ast = \cI_1 \cup \Big(\bigcup_{c\in CC} J(c) \Big) \cup (CC_+\sm CC) .
\end{equation}
But again, we just took all the nets we had, decomposed them into single curves
(some times, just points), and got a bunch of curves $\Gamma_j$.
With our new notation, \eqref{14.3} becomes
\begin{equation}\label{14.8}
\Sigma_F(\Gamma^\ast) = \bigcup_{j\in J^\ast} \Sigma_F(\Gamma_j).
\end{equation}

Now we want to replace each $\Sigma_F(\Gamma_j)$, $j\in \cI_1 \cup \bigcup_{c\in CC} J(c)$,
by a better surface $\Sigma_G(\Gamma_j)$, and this is the place where we shall use Section \ref{S8}.

We start with the case when we change nothing. When $j \in CC_+\sm CC$, i.e., when 
$j = c_\ell$ comes from one of our special components, we just keep
\begin{equation}\label{14.9}
\Sigma_G(\Gamma_j) = \Sigma_F(\Gamma_j) = [0,\ell],
\end{equation}
where the second part comes from the fact that $\Gamma_{c_\ell} = \{ \ell \}$.

But we do the modification for all the other $\Gamma_j$, including the
$\Gamma_j$, $j\in \cI_1$. In this  last case, $\Gamma_j$ does not come
from another curve through the construction of Section \ref{S7}, but (if $\varepsilon$
in \eqref{4.3} is small enough), the proof of Proposition \ref{t5.4} shows that it still 
satisfies the assumptions of Section \ref{S8}; see Remark \ref{r6.3n}.

The same remark applies to $\Gamma_j$ when it comes from Configuration H; 
in this case we decided in Section \ref{S9} to keep $\Gamma_j = \cL_j$
(i.e., without applying the construction of Section \ref{S7} to it), precisely
because Proposition \ref{t5.4} tells us that it is essentially useless.

So fix $j\in J^\ast \sm (CC_+\sm CC)$; recall that $\Gamma_j$ is a small Lipschitz graph 
over some geodesic $\rho_j$, and by changing coordinates in $\R^n$ we can assume that
\begin{equation}\label{14.10}
\rho_j = 
\big\{(\cos\theta, \sin\theta, 0) \, ; \, \theta\in [0,T_j] \big\} \subset \R^2 \times \{ 0 \}.
\end{equation}
We do not have a lower bound on $T_j = \length(\rho_j)$ as in \cite{C1}, 
but the construction yields $T_j \leq \frac{\pi }{ 2} + 2\tau$, 
since $\length(\cC_i) \leq \frac{\pi}{2}$ for $i\in \cI$, and this will be enough.

It may be that $\Gamma_j$ is only $10^3\lambda$-Lipschitz
(if it comes from Lemma \ref{t10.5}), but we shall assume that $\lambda$ is so small
that the results of Section \ref{S8} apply to $\Gamma_j$ anyway. Then we get a new surface
$\Sigma_G(\Gamma_j)$, with the same boundary as $\Sigma_F(\Gamma_j)$, and a few
additional properties. Let us say more and recall a little bit of Section \ref{S8} at the same time.

We started from the (infinite positive) cone over $\Gamma_j$, which is the graph of some
function $F$ which is defined on a sector of $\R^2$ and naturally homogeneous of degree $1$. 
In fact, we restricted $F$ to the domain $D_T$ of \eqref{8.2}, and obtained a set
$\Sigma'_F$, which is the graph of $F$ over $D_T$ and contains $\Sigma_F = \Sigma_F(\Gamma_j)$.

Then we constructed a new function $G$ on the same domain $D_T$, 
which is also null on the segments $[0,a]$ and $[0,b]$, where $a$ and $b$ 
denote the endpoints of $\rho_j$. 
In addition, $G = F$ on $D_T \sm B(0,9/10)$ (see \eqref{8.12}),
so that the graph $\Sigma'_G$ coincides with $\Sigma'_F$ in a small neighborhood of $\d \B$.
So, in some way, $\Sigma'_G$ and $\Sigma'_F$ have the same boundary on $\B$. 

Then we set $\Sigma_G(\Gamma_j) = \Sigma_G = \Sigma'_G \cap \B$. 
Just like $\Sigma_F$, $\Sigma_G$ is bounded by $\Gamma_j$ on $\d \B$, and by
the two segments $[0,a]$ and $[0,b]$.
Let us say why it is a little better than $\Sigma_F$.

There is an additional condition \eqref{8.13}, which implies that for some small 
constant $\kappa >0$,
\begin{equation}\label{14.11}
\Sigma_G(\Gamma_j) \cap \ol B(0,\kappa) = \Sigma_F(\rho_j) \cap \ol B(0,\kappa),
\end{equation}
where as usual $\Sigma_F(\rho_j)$ is the cone over $\rho_j$; we will use this later, 
but for this section we do not care.

Also, $G$ is $C\lambda$-Lipschitz (by \eqref{8.14}), so $\Sigma_G$ is not too
wild, and the important new information is that $\Sigma_G$ has less area than $\Sigma_F$,
since \eqref{8.18} says that
\begin{equation}\label{14.12}
\H^2(\Sigma_G(\Gamma_j)) \leq \H^2(\Sigma_F(\Gamma_j))
- 10^{-4}[\length(\Gamma_j) - \length(\rho_j)].
\end{equation}

Now we glue the $\Sigma_G(\Gamma_j)$ together, and get the set
\begin{equation}\label{14.13}
\Sigma_G(\Gamma^\ast) = \bigcup_{j\in J^\ast} \Sigma_G(\Gamma_j).
\end{equation}

Our next task is to construct a competitor $F^1$ whose main piece is contained
in $\Sigma_G(\Gamma^\ast)$ rather than $\Sigma_F(\Gamma^\ast)$ 
(compare to \eqref{13.13}).
For this the simplest is to continue the deformation that led to $F^0$, i.e., 
deform $\Sigma_F(\Gamma^\ast)$ into $\Sigma_G(\Gamma^\ast)$. 
We proceed piece by piece, and for each $j\in J^\ast$ find a mapping 
$\psi_j : \Sigma_F(\Gamma_j) \to \Sigma_G(\Gamma_j)$.

Fix $j\in J^\ast$.
If $\Gamma_j$ comes from one of our special components $c_\ell$, then 
$\Sigma_G(\Gamma_j) = \Sigma_F(\Gamma_j) = [0,\ell]$ by \eqref{14.9}, 
and we just take $\psi_j(z)=z$.
So let us assume that $\Gamma_j$ is a real curve, either an $\cL_j$,
$j\in \cI_1$, or coming from the construction of Section \ref{S7}.

Choose as above coordinates in $\R^n$ so that $\rho_j$ is, as in \eqref{14.10}, 
an arc of circle inside $P = \R^2$. Call $\pi$ and $\pi^\perp$ the orthogonal projections on
$P$ and its orthogonal complement $P^\perp$, and otherwise keep the same notation as above. 
Recall that $\Sigma'_F$ is the graph of $F: P \to P^\perp$ over $D_T$, and similarly with
$\Sigma'_G$ and $G$. We define $\psi_j : \Sigma'_F \to \Sigma'_G$ by
\begin{equation}\label{14.14}
\psi_j(z) = \pi(z) + G(\pi(z));
\end{equation}
in other words, we project along the direction of $P^\perp$. We are only interested in the restriction
of $\psi_i$ to $\Sigma_F(\Gamma_j) = \Sigma_F = \Sigma'_F \cap \B$. Let us check that 
\begin{equation}\label{14.15}
\psi_j(z) = z \ \text{ for } z\in \d\Sigma_F(\Gamma_j),
\end{equation}
where the boundary $\d\Sigma_F(\Gamma_j)$ is composed of the two line segments
$[0,a]$ and $[0,b]$ that go from $0$ to the endpoints of $\rho_j$, and
$\Sigma'_F \cap \d \B$. For $[0,a]$ and $[0,b]$, this is just because $F=G=0$
on these two segments. For $\Sigma'_F \cap \d \B$, we can even see that
\begin{equation}\label{14.16}
\psi_j(z) = z \ \text{ for } z\in \Sigma_F(\Gamma_j) \sm B\big(0,\frac{99}{100}\big),
\end{equation}
because $F=G$ outside of $\frac{9}{10} \B$ (by \eqref{8.12}); of course we also use the fact that
$||F||_\infty$ is small if $\lambda$ is small enough, to make sure that the two graphs
$\Sigma'_F$ and $\Sigma'_G$ coincide outside of $\frac{99}{100} \B$.

So we have mappings $\psi_j : \Sigma_F(\Gamma_j) \to \Sigma_G(\Gamma_j)$, and we put them 
together to get a mapping $\psi : \Sigma_F(\Gamma^\ast) \to \Sigma_G(\Gamma^\ast)$.
Here we use the fact that the curves $\Gamma_j$ only meet at their endpoints, and with large
angles; the result is that the $\Sigma_F(\Gamma_j)$ only meet along the segments that
go from $0$ to these endpoints, and with the same angles. Since $\psi_j(z)=z$ along
these segments, we get that $\psi$ is well defined, and even Lipschitz (because each piece
is Lipschitz). We do not care if the Lipschitz constant is large (typically, if two of the
$\Gamma_j$ get close to each other somewhere else than the common endpoints),
so we shall not try to check that this does not happen.
Similarly, we shall not try to show that the sets $\Sigma_G(\Gamma_j)$ are disjoint; 
they probably are, but our argument does not need this.
Finally let us observe that because of \eqref{14.16},
\begin{equation}\label{14.17}
\psi(z) = z \ \text{ for } z\in \Sigma_F(\Gamma^\ast) \sm B\big(0,\frac{99}{100}\big).
\end{equation}
Our second competitor is $F^1 = \varphi^1(E)$, where $\varphi^1$ is defined by
\begin{equation}\label{14.18}
\varphi^1(x) = \varphi^0(x) \ \text{ for } x \in E \sm B(0,1-\sigma)
\end{equation}
and
\begin{equation}\label{14.19}
\varphi^1(x) = \psi \circ \varphi^0(x) \ \text{ for } x \in E \cap B(0,1-\sigma).
\end{equation}
This last part makes sense, because $\varphi^0(x) \in \Sigma_F(\Gamma^\ast)$
for $x \in E \cap B(0,1-\sigma)$ (by \eqref{13.13}).

Let us check that $\varphi^1$ is Lipschitz on $E$. The only potential problem
is across $\d B(0,1-\sigma)$. Recall from \eqref{13.7} that
$\varphi^0(x) \in \Gamma^\ast \subset \S$ for $x\in E \cap \d B(0,1-\sigma)$,
so we can find $\sigma'>\sigma$ such that $|\varphi^0(x)| > \frac{99}{100}$
for $\in E \cap \sm B(0,1-\sigma')$. Then \eqref{14.19} actually yields
\begin{equation}\label{14.20}
\varphi^1(x) = \varphi^0(x) \ \text{ for } x \in E \cap B(0,1-\sigma)\sm B(0,1-\sigma'),
\end{equation}
so there is an annulus where the two definitions coincide, and $\varphi^1$ is Lipschitz.
We can of course define a one parameter family $\{ \varphi^1_t \}$, $0 \leq t \leq 1$,
by linear interpolation, as we did in \eqref{13.10}, and, as before, the fact that
\begin{equation} \label{14.21}
\text{the $\varphi_t^1$, $0 \leq t \leq 1$, define an acceptable deformation for $E$ in $\ol B(0,1)$} 
\end{equation}
will follow as soon as we check the boundary constraint, i.e., that
\begin{equation}\label{14.22}
\varphi^1_t(x) \in L \ \text{ for } x\in E \cap L.
\end{equation}
When $|x| \geq 1-\sigma$, $\varphi^1_t(x) = \varphi^0_t(x)$, and we already checked
this. When $|x| < 1-2\sigma$, $\varphi^0(x) = 0$ by \eqref{13.9}, and then
$\varphi^1(x) = 0$ and $\varphi^1_t(x) = (1-t)x \in L$.
We are left with the case when $x\in A(2\sigma) \sm A(\sigma)$.
We already checked below \eqref{13.11} that $\ell = x/|x|$ lies in $K\cap L$. 
By \eqref{12.4}, $p(\ell) = \ell$.
By \eqref{13.8}, $\varphi^0(x) = \alpha \ell$ for some $\alpha \in [0,1]$.

Now there are two cases. If $\ell$ lies in one of our special components $c =c_\ell$,
then $\Sigma_F(\Gamma_c) = \Sigma_G(\Gamma_c) = [0,\ell]$, and we took 
$\psi_j(z) = z$ on $\Sigma_F(\Gamma_c)$ (where $j$ is the element of $J^\ast$
that comes from $c_\ell$). Otherwise, $\ell$ lies in one of the regular components
$c \in CC(\ell)$. In this case, which comes from Configuration $1$, $2+$, $3=2+1$,
or $3+$, we made sure to include $\ell$ in the Lipschitz net $\Gamma_c$, not only
as a point, but as a vertex. This means that $\ell$ is actually an endpoint of one of our
curves $\Gamma_j$, $j\in J^\ast$, and by \eqref{14.15} $\psi(z) = z$ on $[0,\ell]$.
Thus $\varphi^1(x) = \varphi^0(x) = \alpha \ell$ by \eqref{14.19}, and
$\varphi^1_t(x)\in L$ by the analogue of \eqref{13.10} for the $\varphi^1_t$;
\eqref{14.22} and \eqref{14.21} follow.

With the terminology of Definition \ref{t1.1}, $F^1$ is a sliding competitor for 
$E$ in $\overline \B$ (in fact our first interesting competitor), and since $E$
is a sliding almost minimal set, Definition \ref{t1.2} says that \eqref{1.11} holds for
$F^1$. That is,
\begin{equation}\label{14.23}
\H^2(E \cap \ol\B) \leq \H^d(F^1 \cap \overline \B) + h(1) 
\leq \H^d(F^1 \cap \overline \B) + \varepsilon
\end{equation}
by \eqref{4.1} and \eqref{4.2}. Now we cut $F^1 \cap \ol\B$ into pieces.
First observe that
\begin{equation}\label{14.24}
F^1 \cap \overline \B = \varphi^1(E\cap \ol \B)
\end{equation}
because $\varphi^1(x) = \varphi^0(x) = x$ for $x\in E \sm \ol\B$
(by \eqref{14.18} and \eqref{13.5}). 
We start with an exterior part which is the same as before; that is,
\begin{equation}\label{14.25}
\varphi^1(E \cap \ol\B \sm B(0,1-\sigma))
= \varphi^0(E \cap \ol\B \sm B(0,1-\sigma))
= \varphi^0(E \cap A(\sigma)) = F(\sigma)
\end{equation}
by \eqref{14.18}, \eqref{13.1} and \eqref{13.14}. The size of this part will be estimated
by \eqref{13.15} and \eqref{13.26}. We are left with
\begin{equation}\label{14.26}
\varphi^1(E \cap B(0,1-\sigma)) \subset \psi(\varphi^0(E \cap B(0,1-\sigma)))
\subset \psi(\Sigma_F(\Gamma^\ast)) \subset \Sigma_G(\Gamma^\ast),
\end{equation}
by \eqref{14.19}, \eqref{13.13}, and the definition of $\psi$.
We shall thus need to estimate $\H^2(\Sigma_G(\Gamma^\ast))$.
We start with an easier estimate for the cone $\Sigma_F(\Gamma^\ast)$.
We claim that
\begin{equation}\label{14.27}
\H^2(\Sigma_F(\Gamma^\ast))
= \sum_{j\in J^\ast} \H^2(\Sigma_F(\Gamma_j))
= \frac{1}{2} \sum_{j\in J^\ast} \H^1(\Gamma_j)
= \frac{1}{2} \H^1(\Gamma^\ast)
= \frac{1}{2} \sum_{c\in \cI_1 \cup CC} \H^1(\Gamma_c).
\end{equation}
The first equality is true because the union is disjoint, except for segments that
come from the endpoints of the $\Gamma_j$.
For the second part, the simplest is to use the area formula.
Let $z : I \to \Gamma_j$ denote a parameterization of $\Gamma_j$
by arclength. Then we have a parameterization of $\Sigma_F(\Gamma_j)$
by $(t,x) \in [0,1] \times I\to t z(x) \in \Sigma_F(\Gamma_j)$
(compare with the definition \eqref{14.2} if needed). The area formula says that
\begin{equation}\label{14.28}
\H^2(\Sigma_F(\Gamma_j)) = \int_{[0,1] \times I} J(t,x) dx dt,
\end{equation}
where $J$ is the appropriate Jacobian. Since $z$ takes values in the sphere, a simple 
computation says that $J(t,x) = t$; then
\begin{equation}\label{14.29}
\int_{[0,1] \times I} J(t,x) dx dt = \int_{[0,1] \times I} t dx dt = \frac{| I | }{ 2}
= \frac{\length(\Gamma_j)}{2} ,
\end{equation}
as needed. The third identity comes from \eqref{14.5} (the $\Gamma_j$ are essentially
disjoint), and in the last one we used \eqref{14.7}, regrouped indices $j \in J(c)$ for $c\in CC$, and
simply dropped the exceptional sets $\Gamma_{c_\ell}$ coming from $CC_+ \sm CC$, 
because they are singletons $\{ \ell \}$ with no $H^1$-measure.

With the present notation, \eqref{9.6} says that $\H^1(\Gamma_c) \leq \H^1(\ol\gamma_c)$ 
for $c\in CC$; since we do not want to forget too fast what we win, set
\begin{equation}\label{14.30}
\Delta_1 = \sum_{i\in \cI_1} [\H^1(\cL_i) - \H^1(\Gamma_i)]
+\sum_{c \in CC} [\H^1(\ol\gamma_c) - \H^1(\Gamma_c)],
\end{equation}
where we observe that all the terms are nonnegative. With the presentation we chose
(using Remark \ref{r6.3n} and then taking $\Gamma_i = \cL_i$ for $i\in \cI_1$), the
first terms disappear; if we had chosen the other option where $\Gamma_i$
is obtained from $\cL_i$ by the method of Section \ref{S7}, they would exist but 
would not harm (by \eqref{7.16}). Then
\begin{equation}\label{14.31}
\begin{aligned}
\sum_{c\in \cI_1 \cup CC} \H^1(\Gamma_c) 
&= \sum_{i\in \cI_1} \H^1(\Gamma_i) + \sum_{c\in CC} \H^1(\Gamma_c)  
\cr&
\leq \sum_{i\in \cI_1} \H^1(\cL_i) + \sum_{c\in CC} \H^1(\ol\gamma_c) - \Delta_1
= \H^1(\gamma^\ast) - \Delta_1
\end{aligned}
\end{equation}
because $\Gamma_i = \cL_i$ for $i\in \cI_1$, by \eqref{13.25}, and because the union 
in \eqref{13.25} is essentially disjoint. 
Notice that $\gamma^\ast \subset E\cap \S$ by construction, 
so we will really save $\H^1(E \cap \S \sm \gamma^\ast)$ here. 
In the mean time, we return to \eqref{14.27} and get that
\begin{equation}\label{14.32}
\H^2(\Sigma_F(\Gamma^\ast)) = \frac{1}{2} \sum_{c\in \cI_1 \cup CC} \H^1(\Gamma_c)
\leq \frac{1}{2} \H^1(\gamma^\ast) - \frac{\Delta_1}{2}.
\end{equation}
Next we record what we win in \eqref{14.12}. Set
\begin{equation}\label{14.33}
\Delta_2 = \sum_{j\in J^\ast} [\length(\Gamma_j - \length(\rho_j)] \geq 0;
\end{equation}
then by \eqref{14.12}
\begin{equation}\label{14.34}
\H^2(\Sigma_G(\Gamma^\ast)) -\H^2(\Sigma_F(\Gamma^\ast)) 
\leq \sum_{j\in J^\ast} \H^2(\Sigma_G(\Gamma_j)) -\H^2(\Sigma_F(\Gamma_j))
\leq - 10^{-4} \Delta_2,
\end{equation}
and
\begin{eqnarray}\label{14.35}
\H^2(E \cap \ol\B) 
&\leq& \H^2(F^1 \cap \ol\B) + h(1)
\leq \H^2(F(\sigma))+ \H^2(\Sigma_G(\Gamma^\ast)) + h(1)
\nn\\
&\leq&  M(\sigma) + [\H^2(\Sigma_G(\Gamma^\ast))-\H^2(\Sigma_F(\Gamma^\ast))] 
+ \frac{1}{2} \H^1(\gamma^\ast) - \frac{\Delta_1}{2} + h(1)
\nn\\
&\leq& M(\sigma)+ \frac{1}{2} \H^1(\gamma^\ast) 
- \frac{\Delta_1}{2} - 10^{-4} \Delta_2 +  h(1)
\end{eqnarray}
by the first part of \eqref{14.23}, \eqref{14.24}, \eqref{14.25}, and \eqref{14.26},
then \eqref{13.15}, \eqref{14.32}, and \eqref{14.34}. 

Next we estimate $M(\sigma)$. Let $\varepsilon'$ be a very small number, to be chosen
later. We may now choose $\sigma$ very small, such that
\begin{equation} \label{14.36}
M(\sigma) \leq
\varepsilon' + \limsup_{\sigma \to 0} M(\sigma) 
\leq \varepsilon' + C \tau \sum_{c\in CC} \H^1(\ol\gamma_c \sm \Gamma_c)
+C\tau \H^1(E \cap \S \sm \gamma^\ast),
\end{equation}
where the second estimate comes from \eqref{13.26}.
For each $c\in CC$, $\ol\gamma_c \sm \Gamma_c$ is contained in the
symmetric difference $\Delta(\overline\gamma_c ,\Gamma_c)$ that shows up in \eqref{9.4}.
Then by \eqref{9.7}
\begin{equation}\label{14.37}
\begin{aligned}
\sum_{c\in CC} \H^1(\ol\gamma_c \sm \Gamma_c)
&\leq \sum_{c\in CC} \H^1(\Delta(\overline\gamma_c ,\Gamma_c))
\cr&\leq C(\lambda) \sum_{c\in CC} [\H^1(\overline\gamma_c)- \H^1(\Gamma_c)]
+ [\H^1(\Gamma_c)-\H^1(\rho_c)]
\end{aligned}
\end{equation}
where the notation has been adapted, and $\rho_c = \cup_{j\in J(c)} \rho_j$, 
by \eqref{9.5} and the notation of \eqref{14.6};
also see the definition of $\rho_j$ above \eqref{14.10}.
In the first sum we recognize $\Delta_1$ from \eqref{14.30}, and for the second sum
we notice that
\begin{equation}\label{14.38}
\H^1(\Gamma_c)-\H^1(\rho_c) =  \sum_{j\in J(c)} [\H^1(\Gamma_j)-\H^1(\rho_j)]
\end{equation}
because $\Gamma_c$ is the disjoint union of the $\Gamma_j$, $j\in J(c)$
(see \eqref{14.6}), and the $\rho_j$ also are disjoint (again by construction, by the same proof).
We recognize a partial sum of $\Delta_2$; thus \eqref{14.36} and \eqref{14.37} yield
\begin{equation}\label{14.39}
M(\sigma) \leq C'(\lambda) \tau (\Delta_1+\Delta_2) + C \tau \H^1(E \cap \S \sm \gamma^\ast)
+\varepsilon'.
\end{equation}
For the second term, we we just observe
that since $\gamma^\ast \subset E \cap \S$,
\begin{equation}\label{14.40}
\H^1(E \cap \S) = \H^1(\gamma^\ast) + \H^1(E \cap \S \sm \gamma^\ast).
\end{equation}
Now comes the main relation between $\tau$ and $\lambda$: we require $\tau$ to be so small,
depending on $\lambda$, that $C'(\lambda) \tau \leq 10^{-5}$, and $C\tau < 1/4$ 
for the second term; this way both terms of \eqref{14.39} are eaten and \eqref{14.35} yields
\begin{equation}\label{14.41}
\begin{aligned}
\H^2(E \cap \ol\B) 
&\leq \frac{1}{2} \H^1(\gamma^\ast) + C \tau \H^1(E \cap \S \sm \gamma^\ast)
- \frac{\Delta_1}{2} - 10^{-5} \Delta_2 +  h(1) +\varepsilon'
\cr&\leq \frac{1}{2} \H^1(E \cap \S)  - \frac{1}{4} \H^1(E \cap \S \sm \gamma^\ast) 
- \frac{\Delta_1}{2} - 10^{-5} \Delta_2 +  h(1) +\varepsilon'.
\end{aligned}
\end{equation}
We kept $\varepsilon'$ obediently, but since it can be taken arbitrarily small, 
we may now drop it from \eqref{14.41}.
We want a more concise version of this, so let us reorganize some of the terms.
We start with 
\begin{equation}\label{14.42}
\begin{aligned}
\H^1(\gamma^\ast) - \H^1(\Gamma^\ast)
&= \sum_{i\in \cI_1} \H^1(\cL_i) + \sum_{c\in CC} \H^1(\ol\gamma_c) 
- \sum_{c\in \cI_1 \cup CC} \H^1(\Gamma_c)
\cr& = \sum_{i\in \cI_1} [\H^1(\cL_i)-\H^1(\Gamma_i)]
+ \sum_{c\in CC} [\H^1(\ol\gamma_c) - \H^1(\Gamma_c)]
= \Delta_1
\end{aligned}
\end{equation}
by \eqref{13.25}, the end of \eqref{14.27}, and \eqref{14.30}.
Next we introduce the geodesic net
\begin{equation}\label{14.43}
\rho^\ast = \bigcup_{j\in J^\ast} \rho_j \, ;
\end{equation}
these curves are essentially disjoint (by the same proof as for the $\Gamma_j$);
we possibly included one or two degenerate curves $\{ \ell \}$.
In terms of estimates, these degenerate curves will not count, so we may also have 
dropped them too. But for the moment we keep them. Anyway,
\begin{equation}\label{14.44}
\H^1(\Gamma^\ast) - \H^1(\rho^\ast) 
= \sum_{j\in J^\ast} [\H^1(\Gamma_j)-\H^1(\rho_j)] = \Delta_2
\end{equation}
by \eqref{14.5}, because the $\Gamma_j$ are also essentially disjoint,  and by \eqref{14.33}.
By \eqref{14.40} and this,
\begin{equation}\label{14.45}
\H^1(E\cap \S) - \H^1(\rho^\ast) 
= \H^1(E \cap \S \sm \gamma^\ast) + \H^1(\gamma^\ast) - \H^1(\rho^\ast)
= \H^1(E \cap \S \sm \gamma^\ast) + \Delta_1 + \Delta_2
\end{equation}
and \eqref{14.41} implies that
\begin{equation}\label{14.46}
\H^2(E \cap \ol\B) \leq \frac{1}{2} \H^1(E \cap \S) 
- 10^{-5} [\H^1(E\cap \S) - \H^1(\rho^\ast)] +h(1).
\end{equation}
We can be confident that this will lead to reasonable differential inequalities
in some cases, because it looks a lot like (9.69) in \cite{C1}. But we also expect, because
this is what happens in \cite{C1}, that it will not be great in some other cases, and this is the
reason why we introduce a last competitor $F^2$ in the next section, which uses the
full length property. Modulo computations that will be done below,
the estimate above shows that the most delicate case
is probably when $\H^1(E\cap \S) - \H^1(\rho^\ast)$ is very small,
i.e., when $E\cap \S$ actually looks a lot like a collection of geodesics $\rho_j$.
The point of the full length property is somehow to take care of this situation, at
least at the level of definitions.

\section{A third competitor that uses the full length}
\label{S15}

As was discussed near (9.70) in \cite{C1} (and we propose to trust this for the moment), 
\eqref{14.46} will give good differential inequalities when 
\begin{equation}\label{15.1}
\H^1(\rho^\ast) \leq 2 \H^2(X \cap \B).
\end{equation}
Otherwise, we shall need to improve \eqref{14.46} a little bit, by a quantity which is roughly 
proportional to
\begin{equation}\label{15.2}
\Delta_L = \big[\H^1(\rho^\ast) - 2 \H^2(X \cap \B) \big]_+ 
= \max\big(0,\H^1(\rho^\ast) - 2 \H^2(X \cap \B)\big).
\end{equation}
In this section we assume that 
\begin{equation}\label{15.3}
X \ \text{ has the full length property,}
\end{equation}
improve our second competitor in the small tip near the origin, and use this show that
under the assumptions of Section \ref{S4} and if $\varepsilon$ in \eqref{4.3} 
is chosen small enough,
\begin{equation}\label{15.4}
\H^2(E \cap \ol\B) \leq \frac{1}{2} \H^1(E \cap \S) - 10^{-5} [\H^1(E\cap \S) - \H^1(\rho^\ast)]
- C^{-1} \Delta_L+h(1).
\end{equation}
This will be our main comparison estimate, the one that leads to nice differential inequalities.
The constant $C$ in \eqref{15.4} depends on $X$, in particular through $\eta(X)$ in \eqref{3.12}
and the small constants $\eta$ and $c$ in the full length property. Also, $\varepsilon$ will need to be small,
depending on our usual constants $\lambda$ and $\tau$, but also the small $\eta$ in the full 
length property.

Notice that when \eqref{15.1} holds, \eqref{15.4} is just the same as \eqref{14.46};
thus we can assume that \eqref{15.1} fails.

We start from \eqref{14.11}, which says that for $j\in J^\ast$,
$\Sigma_G(\Gamma_j)$ coincides with the cone $\Sigma_F(\rho_j)$ in a small ball
$\ol B(0,\kappa)$. We use \eqref{14.13} to take the union, and get that
\begin{equation}\label{15.5}
\Sigma_G(\Gamma^\ast) \cap \ol B(0,\kappa) = \Sigma_F(\rho^\ast) \cap \ol B(0,\kappa);
\end{equation}
we shall see later that \eqref{14.24}-\eqref{14.26} give the same description for 
$F^1\cap \ol B(0,\kappa)$.
Here $\kappa$ is the small absolute constant of \eqref{8.13}.
% which does not even depend on the dimension $n$.

We shall start our discussion by assuming that
\begin{equation}\label{15.6}
\rho^\ast = \varphi_\ast(K) \ \text{ for some $\varphi \in \Phi_X^{+}(\eta)$,}
\end{equation}
for some $\eta > 0$ which is small enough for us to apply the full length property \eqref{15.3}.
This is not always the case, but we shall first assume \eqref{15.6}, prove our main estimate, 
and then return and take care of the differences. Set $X_1 = \varphi_\ast(X)$
(the full cone over $\varphi_\ast(K)$). Notice that since \eqref{15.6} implies that
\begin{equation}\label{15.7}
\Sigma_F(\rho^\ast) = \varphi_\ast(X) \cap \ol\B = X_1 \cap \ol\B, 
\end{equation}
we can try to use competitors for $X_1$ to improve $\Sigma_F(\rho^\ast)$,
$\Sigma_G(\Gamma^\ast)$, and then $F^1$.
Observe also that $\varphi$ is ``injective'', i.e., 
$\varphi \in \Phi_X^{+,i}(\eta)$, because the arcs $\rho_j$ that compose it
(as in \eqref{14.43}) only meet at their common endpoints. The quantity
$\Delta(\varphi)$ of \eqref{3.24} is
\begin{equation}\label{15.8}
\Delta(\varphi) = \H^1(\varphi_\ast(K)) - \H^1(K) = \H^1(\rho^\ast) - \H^1(K)
= \H^1(\rho^\ast) - 2\H^1(X \cap \B) = \Delta_L > 0
\end{equation}
by the proof of the second part of \eqref{14.27}, the definition \eqref{15.2},
and because \eqref{15.1} fails.
So we are in postion to apply our assumption \eqref{15.3}. 
Definition \ref{t3.1} says that there is a sliding competitor $\wt X$ for $X_1$ in $\ol B(0,1)$ 
such that \eqref{3.25} holds. By Definition \ref{t1.1}, there 
is a one parameter family of functions $g_t : X_1 \to \R^n$
(we change the name because $\varphi$ is already used), such that \eqref{1.4}-\eqref{1.8}
hold with $B = \ol B(0,1)$, and for which $\wt X = g_1(X_1)$.

We want to use the $g_t$ to construct a competitor $F^2$ for $F^1$, and by the same
token for $E$. Initially, the mappings $g_t$ are only defined on $X_1$, but we can extend $g_t$
to $X_1 \cup (\R^n \sm B(0,2))$ by setting
\begin{equation}\label{15.9}
g_t(x) = x \ \text{ for } x \in \R^n \sm B(0,2).
\end{equation}
This gives mappings $g_t$ that are still continuous, by \eqref{1.5}, and such that
$g_t(x) \in L$ when $x\in L$ (by \eqref{1.7}). Also, $g_1$ is Lipschitz (by \eqref{1.8}).
We set 
\begin{equation}\label{15.10}
\varphi^2(x) = \frac{\kappa}{2}\, g_1(2\kappa^{-1}\varphi^1(x)) \ \text{ for } x\in E.
\end{equation}
Let us check that $\varphi^2$ is well defined. Notice that we can use \eqref{15.9}
as soon as $|\varphi^1(x)| \geq \kappa$.

For $x\in E \sm \B$, $\varphi^1(x) = \varphi^0(x) = x$ by \eqref{14.18} and \eqref{13.5},
we can use \eqref{15.9}, and we get that $\varphi^2(x)=x$.

For $x \in A(\sigma) = \ol\B \sm B(0,1-\sigma)$, $\varphi^1(x) = \varphi^0(x)$ by \eqref{14.18},
and $\varphi^0(x) \in [x, p(x)]$ by \eqref{13.6}. Since $p(x) = p(x/|x|)$ by \eqref{13.3}
and $p(x/|x|) \in \S \cap B(x/|x| ,60\tau)$ by \eqref{12.2}, we see that
$|\varphi^1(x)| > 1/2 > \kappa$, and we can apply \eqref{15.9} again. In this case
\begin{equation}\label{15.11}
\varphi^2(x) = \varphi^0(x) \in F(\sigma),
\end{equation}
by \eqref{13.14}. We are left with the case when $x\in E \cap B(0,1-\sigma)$.
In this case, \eqref{14.26} says that $\varphi^1(x) \in \Sigma_G(\Gamma^\ast)$.
If $|\varphi^1(x)| \geq \kappa$, we use \eqref{15.9} and we get that
\begin{equation}\label{15.12}
\varphi^2(x) = \varphi^1(x) \in \Sigma_G(\Gamma^\ast) \sm B(0,\kappa).
\end{equation}
Otherwise, if $|\varphi^1(x)| < \kappa$,
\begin{equation}\label{15.13}
\varphi^1(x) \in \Sigma_G(\Gamma^\ast) \cap B(0,\kappa) = 
\Sigma_F(\rho^\ast) \cap  B(0,\kappa) 
= X_1 \cap B(0,\kappa)
\end{equation}
by \eqref{15.5} and \eqref{15.7},
\begin{equation}\label{15.14}
2\kappa^{-1}\varphi^1(x) \in X_1 \cap B(0,2)
\end{equation}
because $X_1$ is a cone, and this allows us to
use the initial definition of $g_1$ in \eqref{15.10}. We get that
\begin{equation}\label{15.15}
\varphi^2(x) = \frac{\kappa}{2}\, g_1(2\kappa^{-1}\varphi^1(x))
\in \frac{\kappa}{2}\, g_1(X_1 \cap B(0,2)).
\end{equation}
So $\varphi^2$ is well defined.

The fact that $\varphi^2$ is Lipschitz comes directly from the definitions (in particular, the
fact that the extended $g_1$ is Lipschitz); it is easy to find a one-parameter family
$\{ \varphi_t^2 \}$ that has $\varphi^2$ as its endpoint, and as before the simplest
is to use a formula like \eqref{13.10} and the convexity of $\B$, the boundary property 
\eqref{1.7} holds because $g_1(L) \subset L$ (and $L$ is convex), as usual, so
$F^2 = \varphi^2(E)$ is a sliding competitor for $E$ in $\ol\B$. Thus Definition \ref{t1.2}
yields
\begin{equation}\label{15.16}
\H^2(E \cap \ol \B) \leq \H^2(F^2\cap \ol \B) + h(1) 
\end{equation}
(as in \eqref{14.23}). Now
\begin{equation}\label{15.17}
F^2\cap \ol \B \subset F(\sigma) \cup \big[\Sigma_G(\Gamma^\ast) \sm B(0,\kappa)\big]
\cup \big[\frac{\kappa}{2}\, g_1(X_1 \cap B(0,2)) \big]
\end{equation}
by the discussion above and \eqref{15.11}-\eqref{15.13}.
We cut the last set in two. If $z\in X_1 \cap B(0,2)\sm \B$, then $g_1(z) = z$
by \eqref{1.5} for $g_1$ (i.e., the fact that $\wt X$ is a competitor for $X_1$ in $\ol\B$);
then
\begin{equation}\label{15.18}
\frac{\kappa}{2}\, g_1(z) = \frac{\kappa z}{2} \in X_1 \cap B(0,\kappa)\sm B(0,\kappa/2)
= \Sigma_G(\Gamma^\ast) \cap B(0,\kappa)\sm B(0,\kappa/2)
\end{equation}
by \eqref{15.13}. If instead $z\in X_1 \cap \B$, then
\begin{equation}\label{15.19}
\frac{\kappa}{2}\, g_1(z) \in \frac{\kappa}{2}\, g_1(X_1 \cap \B) 
= \frac{\kappa}{2}\, \wt X \cap \B
\end{equation}
because $g_1(z) = z$ for $z \in X_1 \sm \B$ and $g_1(X_1 \cap \B) \subset \B$
by \eqref{1.5} and \eqref{1.6}.
Thus \eqref{15.17} yields
\begin{equation}\label{15.20}
F^2\cap \ol B \subset F(\sigma) \cup \big[\Sigma_G(\Gamma^\ast) \sm B(0,\kappa/2)\big]
\cup \frac{\kappa}{2}\, \wt X \cap \B
\end{equation}
and 
\begin{eqnarray}\label{15.21}
\H^2(F^2\cap \ol B) &\leq& \H^2(F(\sigma)) + \H^2(\Sigma_G(\Gamma^\ast))
- \H^2(\Sigma_G(\Gamma^\ast) \cap B(0,\kappa/2))
+ \H^2(\frac{\kappa}{2}\, \wt X \cap \B)
\nn\\
& \leq& \H^2(F(\sigma)) + \H^2(\Sigma_G(\Gamma^\ast)) 
- \H^2(X_1 \cap B(0,\kappa/2))
+ \H^2(\frac{\kappa}{2}\, \wt X \cap \B)
\nn\\
&\leq& \H^2(F(\sigma)) + \H^2(\Sigma_G(\Gamma^\ast)) 
+ \frac{\kappa^2}{4} [-\H^2(X_1 \cap \B) + \H^2(\wt X \cap \B)]
\\
&\leq& \H^2(F(\sigma)) + \H^2(\Sigma_G(\Gamma^\ast))
- \frac{c \kappa^2 \Delta(\varphi)}{4} 
\nn
\end{eqnarray}
by \eqref{15.13} and \eqref{3.25}. This is the same estimate as we had for $F^1$,
at the end of the first line of \eqref{14.35}, except that we saved an extra 
$c \kappa^2 \Delta(\varphi) / 4$. Then we continue the computations exactly as in
Section \ref{S14}, and get that
\begin{equation}\label{15.22}
\H^2(E \cap \ol\B) 
\leq \frac{1}{2} \H^1(E \cap \S) - 10^{-5} [\H^1(E\cap \S) - \H^1(\rho^\ast)]
- \frac{c \kappa^2 \Delta_L}{4} + h(1) 
\end{equation}
instead of \eqref{14.46}, and where $\Delta_L$ is given by \eqref{15.2}. 
This proves the desired estimate \eqref{15.4}, with $C^{-1} = c\kappa^2/4$,
but only in the case when when \eqref{15.6} holds.

\ms
Let us now discuss the reasons why \eqref{15.6} may fail, and what to do then.
The problem is with some of the configurations of Section \ref{S6}, which
may not always produce nets $\rho_j$ that follow the description of Section \ref{S3}.

First assume that for some $\ell \in K \cap L$, Configuration H shows up in our 
construction of $\Gamma^\ast$ near $\ell$. 
We intend to show that we do not even need the full length condition to find
better competitors for $X_1$ and $E$, because we can contract a hanging curve 
in $\rho^\ast$.

Recall that when $E \cap \S$ contains a hanging curve that starts from $c_i^\ast$,
we kept the corresponding curve $\cL_i$ both in $\gamma_c$ and $\Gamma_c$,
where $c\in CC$ is the component that contains $c_i^\ast$.
The geodesic $\rho_i = \rho(a_i^\ast,c_i^\ast)$ with the same endpoints as $\Gamma_c$ 
is contained in $\rho^\ast$. Let us identify $i$ with the only index $j \in J(c)$, so that
$\rho_i$ shows up with the same name in the union of \eqref{14.43}.
Notice that its endpoint $c_i^\ast$ is still hanging in $\rho^\ast$, 
which means that it does not lie in any other $\rho_j$, $j\in J^\ast \sm \{ i \}$. 
Set $\rho'_i = \rho_i \sm \{ a_i \}$; then $\rho'_i$ does not meet
any other $\rho_j$, and this means that the mapping 
$f : \rho^\ast \to \rho^\ast \sm \rho'_j$ defined by $f(z)=z$ for $z\in \rho^\ast \sm \rho'_j$
and $f(z) = a_i^\ast$ for $z\in \rho_i^\ast$ is Lipschitz (recall that the $\rho_j$ that meet 
$\rho_i$ at $a_i^\ast$ make large angles with $\rho_i$ there).

Let us use $f$ to define a nice competitor for the cone over $\rho^\ast$.
Set $X_2 = \big\{ tz \, ; \, z\in \rho^\ast \text{ and } t \geq 0 \big\}$
and define $g : X_2 \cup (\R^n \sm \B) \to \R^n$ by
\begin{equation}\label{15.23}
g(x) = x \ \text{ for } x \in \R^n \sm \B,
\end{equation}
\begin{equation}\label{15.24}
g(tz) = 2(1-t) tf(z) + (2t-1) t z \ \text{ for $z\in \rho^\ast$ and } \frac{1}{2} \leq t \leq 1,
\end{equation}
where we choose the coefficients so that $g(tz) = tz = z$ when $t=1$, 
and $g(tz) = t f(z)$ when $t=1/2$, and finally
\begin{equation}\label{15.25}
g(tz) = t f(z) \ \text{ for $z\in \rho^\ast$ and } t \leq \frac{1}{2}.
\end{equation}

Let $X_1$ be the (full positive) cone over $\rho^\ast$ (just as in \eqref{15.7}), and
set $\wt X =  g(X_1)$. Let us check that $\wt X$ is a sliding competitor for $X_1$ in $\ol\B$.
Of course we use $g$, and the one parameter family naturally associated with it, defined by 
$g_s(x) = s g(x) + (1-s)x$; the usual simple estimates \eqref{1.4}, \eqref{1.5}, and \eqref{1.9}
hold because $g$ is Lipschitz, \eqref{1.6} holds because $\ol\B$ is convex, and \eqref{1.7}
holds because the only place where $g(x) \neq x$ is the cone over $\rho_j'$, which does not
meet $L$ because $\rho_j$ starts at $c_j^\ast \in \d D$, and goes in the direction of 
$a_i^\ast$ which is away from $\ell$. Next $\wt X =  g(X_1)$ is contained in $X_1$,
but inside $B(0,1/2)$ the cone over $\rho_j'$ is missing. Thus
\begin{equation}\label{15.26}
\begin{aligned}
\H^2(\wt X \cap \ol \B) &\leq \H^2(X_1 \cap \ol \B) - \H^2(\Sigma_F(\rho_j') \cap B(0,1/2)) 
\cr&=  \H^2(X_1 \cap \ol \B) - \frac{1}{8} \H^1(\rho_j')
\leq \H^2(X_1 \cap \ol \B) - \eta(X)
\end{aligned}
\end{equation}
because $\ddist(c_j^\ast,a_j^\ast) \geq 8 \eta(X)$ by 
\eqref{3.11}, \eqref{3.12}, \eqref{5.2bis}, and \eqref{5.38}.
This is even better than the information we obtained from \eqref{15.6} and \eqref{15.3}: 
the proof of \eqref{15.22} yields
\begin{equation}\label{15.27}
\H^2(E \cap \ol\B) 
\leq \frac{1}{2} \H^1(E \cap \S) - 10^{-5} [\H^1(E\cap \S) - \H^1(\rho^\ast)]
- \eta(X) + h(1),
\end{equation}
without even having to assume that \eqref{15.1} fails. Notice that the constant
$\Delta_L$ in \eqref{15.2} is bounded by $\H^1(\rho^\ast) \leq 2 \H^1(K)$
(if $\eta$ is small enough in the definition of full length), so \eqref{15.27} is stronger
than \eqref{15.4}, and we are happy in this case.

A second case when \eqref{15.6} fails is when we encounter Configuration $3 = 2+1$
in the construction of Sections \ref{S9}-\ref{S11}. Recall that in this case we 
chose a center $x_0$, in fact $x_0 = c_1^\ast$ because this was simpler,
then the corresponding $\Gamma$ was composed of three Lipschitz curves, 
one leaving from $\ell$ and two leaving from $x_0 = c_1^\ast$.
At the end of the game, near $\ell$, $\rho^\ast$ is composed of three geodesics
$\rho(c_1^\ast,a_1^\ast)$, $\rho(c_1^\ast,a_2^\ast)$, and $\rho(\ell,a_3^\ast)$. 
More precisely, we claim that
\begin{equation}\label{15.28}
\rho^\ast \cap B(\ell, 9\eta(X)) 
= [\rho(c_1^\ast,a_1^\ast)\cup\rho(c_1^\ast,a_2^\ast)\cup\rho(\ell,a_3^\ast)] 
\cap B(\ell, 9\eta(X)).
\end{equation}
Indeed, all the other $\rho_j$, $j\in J^\ast$ such that meet $B(\ell, 9\eta(X))$ have to come
from curves $\cL_j$, $j\in \cI_1$ (the other option, that they would come from
curves that come from $-D$, is impossible because our curves are not too long).
But in this case \eqref{5.38} says that the two endpoints of $\cL_j$ lie quite close to $\cC_j$,
so does the geodesic $\rho_j$ with the same endpoints, and \eqref{15.23} follows from 
the fact that $\dist(\cC_i, \ell) \geq 10\eta(X)$. This last fact is true,
by \eqref{2.5} and \eqref{3.12}, or the description of the counterexamples that follows \eqref{2.5}, 
plus the fact that the diameter of any exceptional arc  $\cC_k$ for \eqref{2.4} 
(so, that $\cC_k$ ends at $\ell$ or $-\ell$) is controlled by \eqref{3.11}.

In this case, we can find a sliding competitor $\wt X$ for $X_1$ in $\ol\B$, a little bit like 
the one given by $g$ in \eqref{15.23}-\eqref{15.24}, except that (instead of just removing 
it progressively as above) we deform the union $\rho(c_1^\ast,a_1^\ast)\cup\rho(c_1^\ast,a_2^\ast)$ 
into a shorter arc with the same endpoints, such as the union 
$\rho(x_1,a_1^\ast)\cup\rho(x_1,a_2^\ast)$, for some $x_1$ that is a little
closer to $a_1^\ast$ and $a_2^\ast$. 
The reason why we can easily find $x_1$ is that, since the three $\cC_i$, $1 \leq i \leq 3$,
make $120^\circ$ angles at $\ell$, the geodesics $\rho_1$ and $\rho_2$ make an angle smaller 
than $130^\circ$ at $c_1^\ast$. Notice also that $\rho(x_1,a_1^\ast)\cup\rho(x_1,a_2^\ast)$
does not meet $\rho(\ell,a_3^\ast)$ either, which is comforting even though it is not needed.

Now we claim that because of this we can find a sliding competitor $\wt X$ for $X_1$ in $\ol\B$,
such that
\begin{equation}\label{15.29}
\H^2(\wt X \cap \ol B) \leq \H^2(X_1 \cap \ol B) - C^{-1}\eta(X),
\end{equation}
where $C$ is a geometric constant; the verification is rather easy (but a little long),
and we skip it. The interested reader may find more or less the same argument in \cite{C1},
and slightly more elaborate versions, with three branches instead of two, in Section \ref{S26},
starting below \eqref{26.4}.

Now \eqref{15.29} is nearly as good as \eqref{15.26}; so, when Configuration $3 = 2+1$
shows up in the construction, we can still prove that
\begin{equation}\label{15.30}
\H^2(E \cap \ol\B) 
\leq \frac{1}{2} \H^1(E \cap \S) - 10^{-5} [\H^1(E\cap \S) - \H^1(\rho^\ast)]
- C^{-1}\eta(X) + h(1),
\end{equation}
still regardless of whether \eqref{15.1} holds or not. As before, this estimate is better
than \eqref{15.4} because $\Delta_L \leq C$.

Now let us assume that Configurations $H$ and $3 = 2+1$ do not occur.
We have a last case where \eqref{15.6} may fail. Recall that when some
$\ell \in V_0$ does not lie in the net of curves that we constructed, we added 
an element $c_\ell$ to $CC(\ell)$, to get the extended $CC_+(\ell)$, and we also added the 
point $\ell$ to $\rho^\ast$. Denote by $V'_0$ the set of (at most two) points
$\ell$ that we added this way, by $L'$ the (full) positive cone over $V'_0$,
and also set $\rho' = \rho^\ast \sm L'$. Finally denote by $X_1'$
the (full) positive cone over $\rho'$.

First observe that $\rho'$ satisfies \eqref{15.6} 
(if $\varepsilon$ is small enough, as before); this is the reason why we added 
the free option in the definition of $\Phi_X^+(\eta)$ in Section \ref{S3}.
So we can apply the full length condition, and we get a sliding competitor 
$\wt X'$ for $X_1'$ in the ball $\ol\B$. Let $\{ g'_t \}$ denote the associated 
one parameter family of mappings. The $g'_t$ are defined on $X_1'$, and we want an 
extension of $g'_1$ to the full $X_1 = X'_1 \cup L'$. [As was noticed before, 
we only need $g_1$ here, we always compute the one parameter extensions at the end.]

Set $a = g_1(0)$, it is not clear that $a = 0$, but at least $a\in L \cap \ol \B$,
by \eqref{1.6} and \eqref{1.7}. Extend $g_1$ to $L'$, so that it is Lipschitz on 
$L'$, with $g_1(0) = a$, $g_1(x) = x$ for $x\in L' \sm \B$, and 
$g_1(L'\cap \ol\B) \subset L' \cap\ol\B$. This gives a mapping $g_1$, 
now defined on $X_1$, and we want to check that it is Lipschitz.

Clearly it is enough to control $|g'_1(x)-g_1(y)|$ when $x\in X'_1$ and $y\in L'$.
Write $x=tz$, with $t \geq 0$ and $z\in \rho'$; by construction 
$\dist(z,L') > \alpha$ for some (possibly very small) $\alpha > 0$ that depends on 
$\rho^\ast$ and $L'$; then $|x-y| \geq \alpha t$. 
If $|y| \leq 2t$, we say that $|g'_1(x)-g_1(y)| \leq |g'_1(x)-g_1(0)| + |g_1(0)-g_1(y)|
\leq C_1 |x| + C_2 |y| \leq (C_1+2C_2) t \leq (C_1+2C_2) \alpha^{-1}|x-y|$,
which may be very bad but is enough.
Otherwise, $|x-y| \geq |y|/2$ and we just need to change the end of the estimate.

So $g_1$ is Lipschitz, $g_1(L \cap X_1) \subset L$ by construction, and we can 
use the same linear interpolation trick as in \eqref{13.10} to construct a one parameter family 
of mappings that shows that $\wt X = g_1(X_1)$ is a sliding competitor for $X_1$ in $\ol\B$.
Now $\H^2(L') = 0$, so $\H^2(X_1 \cap \ol\B) = \H^2(X'_1\cap \ol\B)$ and 
(when we take take the Lipschitz images by $g_1$)
$\H^2(\wt X\cap \ol\B) = \H^2(\wt X' \cap \ol\B)$. In other words, we still have
\eqref{3.25} for $X_1$ and $\wt X'$, and we may conclude as in the main case.

This completes our verification of \eqref{15.4}.

\ms 
We end this section with a small cosmetic modification of \eqref{15.4}.
Set 
\begin{equation}\label{15.31}
\alpha = \alpha(X) = \min(10^{-5},C^{-1}),
\end{equation}
where $C$ is as in \eqref{15.4}, and observe that in \eqref{15.4} the two main correction terms 
are nonpositive. That is, $\Delta_L \geq 0$ by \eqref{15.2}, and 
$\H^1(E\cap \S) \geq H^1(\rho^\ast)$ by \eqref{14.45} and earlier parts of the proof.
Then \eqref{15.4} implies that
\begin{eqnarray}\label{15.32}
\H^2(E \cap \ol\B) 
&\leq& \frac{1}{2} \H^1(E \cap \S) - \alpha [\H^1(E\cap \S) - \H^1(\rho^\ast)]
- \alpha \Delta_L+h(1) 
\nn\\
&\leq& \frac{1}{2} \H^1(E \cap \S) - \alpha [\H^1(E\cap \S) - \H^1(\rho^\ast)]
-\alpha [\H^1(\rho^\ast) - 2 \H^2(X \cap \B)] + h(1) 
\nn\\
&=& \frac{1}{2} \H^1(E \cap \S) - \alpha [\H^1(E\cap \S) - 2 \H^2(X \cap \B)] + h(1)
\end{eqnarray}
by \eqref{15.2}.

\section{Density excess and a differential inequality}
\label{S16}

Our next goal is to transform our main estimate \eqref{15.4} into a differential inequality,
and then we will integrate it on intervals to get decay estimates for a density excess $f(r)$.

In this section we fix an open set $U$ of $\R^n$ that contains the origin
and a line $L$ through the origin, and we consider a reduced sliding almost minimal 
set $E$ of dimension $2$ in $U$, with sliding boundary $L$.
We shall restrict to radii $r \in (0,r_0)$, where $r_0$ is so small that
\begin{equation}\label{16.1}
B(0,2r_0) \subset U.
\end{equation} 
As in Section \ref{S4}, we shall assume that the gauge function $h$ 
(in the definition of sliding minimal sets) is such that
\begin{equation}\label{16.2}
h(s) \leq C_h s^{\beta} \ \text{ for } 0 < s < 2r_0
\end{equation}
for some constants $C_h \geq 0$ and $\beta > 0$. 

We shall also give ourselves a fixed number $\theta_0 > 0$ and consider 
the density $\theta$ and the density excess $f$ defined on $(0, 2r_0)$ by
\begin{equation}\label{16.3}
\theta(r) = r^{-2} \H^2(E \cap B(0,r)) \ \text{ and } \ 
f(r) = \theta(r) - \theta_0.
\end{equation}
In practice, we will take for $\theta_0$ the density of $E$ at the origin, i.e., 
\begin{equation}\label{16.4}
\theta_0 = \lim_{r \to 0} \theta(r)
\end{equation}
(which exists, as mentioned near \eqref{1.2}), but let us not require this for the moment.
We start with differentiability properties that don't use much.

\begin{lem}\label{t16.1} 
Let $E$ satisfy the assumptions above. Set
\begin{equation}\label{16.5}
v(r) = \H^2(E \cap B(0,r)) \ \text{ for } 0 < r < 2r_0.
\end{equation}
Then $v$ is differentiable almost everywhere on $(0,2r_0)$.
Also, if $b$ is a $C^1$ function on $(0,2r_0)$, the product $bv$ is also
differentiable almost everywhere on $(0,2r_0)$, with $(bv)' = bv'+b'v$ almost everywhere.
In addition, 
\begin{equation}\label{16.6}
(bv)(r_2) \geq (bv)(r_1) + \int_{r_1}^{r_2} (bv)'(r) dr
\ \text{ for } 0 < r_1 \leq r_2 < r_0.
\end{equation}
\end{lem}

The simplest is to refer to Lemma 5.1 in \cite{C1}, but anyway this is not hard: $v$
is nondecreasing, so it is differentiable almost everywhere; it also has a distribution 
derivative $\mu$, and $v' dx \leq d\mu$. This proves \eqref{16.6} for $b=1$.
For general $b$, the differentiability of the product is easy to prove by hand, 
and \eqref{16.6} is proved with a soft integration by parts (i.e., apply Fubini's 
theorem to the right integral).
\qed

\ms
Lemma \ref{t16.1} shows that $\theta$ and $f$ in \eqref{16.3} are differentiable
almost-everywhere on $(0,2r_0)$; next we want to use the previous sections
to derive differential inequalities for $f$, and after this we'll get some decay for $f$.

Before we state our main differential inequality, we introduce some notation concerning
minimal cones and the full length condition. We work with $n$ and $L$ fixed, as above.
Denote by $MC(L)$ the set of minimal cones in $\R^n$ with sliding boundary $L$
(as above \eqref{2.1}). To each cone $X\in MC(L)$, we associate a standard decomposition
as in Section \ref{S3}, and then a geometric constant $\eta(X)$ as in \eqref{3.12}. 
It is fairly easy to see that when $\eta(X)$ is fairly small, it does not depend on our 
choice of standard decomposition, but it would not matter if it did.

Next denote by $FL(L)$ the set of cones $X\in MC(L)$ that also satisfy the full length property.
To $X \in FL(L)$ we also associate as in Definition \ref{t3.1} two small constants 
$\eta > 0$ and $c > 0$, which we call the full length constants for $X$. 
For each choice of positive constants
$c_{fl}$, $\eta_{fl}$, and $\eta_g$, with $\eta_{fl} < \eta_g < 10^{-2}$, say, 
we shall denote by $FL(L,c_{fl}, \eta_{fl}, \eta_g)$ the set of cones
$X\in FL(L)$, which admit a geometric constant $\eta(X) \geq \eta_g$ and full length
constants $c \geq c_{fl}$ and $\eta \geq \eta_{fl}$.
We also associate to this choice a small number $\varepsilon(c_{fl}, \eta_{fl}, \eta_g)$,
which we choose so that the construction and results of Sections \ref{S4}-\ref{S15}
are valid as soon as \eqref{4.1}-\eqref{4.8} are satisfied with 
$\varepsilon \leq \varepsilon(c_{fl}, \eta_{fl}, \eta_g)$, and the constant
$\alpha(c_{fl}, \eta_{fl}, \eta_g)$ that we get in \eqref{15.32} when this happens.

The new assumptions for the next proposition are that for almost each radius $r \in (0,r_0)$,
we can find some constants $c_{fl}(r)$, $\eta_{fl}(r)$, and $\eta_g(r)$, and a 
minimal cone $X(r) \in FL(L,c_{fl}(r), \eta_{fl}(r), \eta_g(r))$, with the following properties.
First
\begin{equation}\label{16.7}
d_{0,2r}(E,X(r)) \leq \varepsilon(c_{fl}(r), \eta_{fl}(r), \eta_g(r))
\end{equation}
(a local Hausdorff distance, as in \eqref{1.13}), where 
$\varepsilon(c_{fl}(r), \eta_{fl}(r), \eta_g(r))$ is the small constant that we get from
the previous sections. We also require that
\begin{equation}\label{16.8}
C_h r_0^{\beta} \leq \varepsilon(c_{fl}(r), \eta_{fl}(r), \eta_g(r)).
\end{equation}

As the reader may have noticed, we are just copying the assumptions of Section \ref{S4}.
Our result will be better if we have a good control on the density 
\begin{equation}\label{16.9}
\theta(X(r)) = \H^2(X(r) \cap B(0,1))
\end{equation}
of the minimal cone $X(r)$; for the moment let us just assume that we have a 
number $q(r) \geq 0$ such that  
\begin{equation}\label{16.10}
\theta(X(r)) \leq \theta_0 + q(r) \ \text{ for } 0 < r < r_0.
\end{equation}
In fact, for our simple applications, we will simply have $\theta(X(r)) = \theta_0$
and $q(r)=0$. We do not need to assume that $q$ is small, but Proposition \ref{t16.2}
below will be hard to apply otherwise.

It is important to let the minimal cone $X(r)$ depend on $r$, even in the good cases where
we take $\theta_0 = \lim_{r \to 0} \theta(r)$ as in \eqref{16.4} and require
$\theta(X(r)) = \theta_0$. The point is that the $X(r)$ could be various blow-up
limits of $E$ at $0$; we do not want to assume that they are all the same, we want
to get this as a conclusion.

Similarly, it would be tempting to require that all the $X(r)$ lie in a same 
$FL(L,c_{fl}, \eta_{fl}, \eta_g)$, but we may have more trouble finding the cones $X(r)$.
We find it more flexible to allow some cones $X(r)$ to have different 
geometric constants $\eta(X(r))$, for instance. We shall see in the next section how to 
choose the $X(r)$ in some simple cases.

\begin{pro} \label{t16.2} 
Let $E$ satisfy the assumptions above (that is, \eqref{16.1}, \eqref{16.2},
\eqref{16.7}, \eqref{16.8}, and \eqref{16.10}. Then
\begin{equation}\label{16.11}
r f'(r) \geq \frac{4\alpha}{(1-2\alpha)} f(r) - 3(h(r)+2\alpha q(r))
\ \text{ for almost every } r \in (0,r_0),
\end{equation} 
where $\alpha = \alpha(c_{fl}, \eta_{fl}, \eta_g)$ 
is the small constant that is associated to $c_{fl}(r)$, $\eta_{fl}(r)$, and $\eta_g(r)$
as in \eqref{15.31}.
\end{pro}

Let us prove the proposition. It turns out that we already did the hard work; the proof will 
be derived softly from the previous sections.
The first thing we have to do is check that the assumptions of Section \ref{S4}
are satisfied (now, without the renormalization $r=1$) for almost every $r\in r_0$.
The three first assumptions \eqref{4.1}-\eqref{4.3} were just copied above.
Next, \eqref{4.4}, the fact that $\H^1(E \cap \d B(0,r)) < +\infty$ is true almost-everywhere,
holds because $\H^2(E \cap B(0,s)) < +\infty$ for $0 < s < r_0$. Since $E$ is
rectifiable, we can deduce this from the coarea formula, but in fact the estimate that we need
here is just is the easy part, which can be obtained directly with a covering lemma. 
% say where, for the poor man's coarea formula

We said earlier that \eqref{4.7} is just requiring that the one-sided Hardy-Littlewood
maximal function of the measure $\mu_\pi$ of \eqref{4.5} is finite at the point $r$,
and since $\mu_\pi$ is a finite measure, the fact that this is true almost everywhere
(and even with weak integral estimates on $C$) is a direct consequence of the 
weak $L^1$ Hardy-Littlewood estimate; see the first pages of \cite{Stein}.

We are left with \eqref{4.8}, which requires maximal function estimates like \eqref{4.7},
but also some manipulation and a density argument; fortunately the proof is done in
Lemma 4.12 of \cite{C1}, and applies here.

So we can use the estimates of the previous sections, and \eqref{15.32} holds
for almost every $r\in (0,r_0)$. Written with the variable $r$, the correct homogeneity,
the notation $\S_r = \d B(0,r)$, and with $\alpha = \alpha(c_{fl}, \eta_{fl}, \eta_g)$,
it says that
\begin{eqnarray}\label{16.12}
r \theta(r) &:=& r^{-1} \H^2(E \cap B(0,r)) \leq r^{-1} \H^2(E \cap \ol B(0,r))
\nn\\ 
&\leq& \frac{1}{2} \H^1(E \cap \S_r) 
- \alpha [\H^1(E\cap \S_r) - 2 r \H^2(X(r) \cap B(0,1))] + r h(r)
\nn\\
&=& \frac{1}{2} \H^1(E \cap \S_r) - \alpha [\H^1(E\cap \S_r) - 2r \theta(X(r))] + rh(r)
\nn\\
&= & \frac{1}{2} \H^1(E \cap \S_r) - \alpha [\H^1(E\cap \S_r) - 2r\theta_0] 
+ 2\alpha r [\theta(X(r)) - \theta_0]+ r h(r) 
\\
&\leq& \frac{1}{2} \H^1(E \cap \S_r) - \alpha [\H^1(E\cap \S_r) - 2r\theta_0] 
+ r (2\alpha q(r)+h(r)) 
\nn
\end{eqnarray}
by \eqref{16.9} and \eqref{16.10}. Next write $\H^1(E\cap \S_r) = 2r x(r)$
for the duration of the computation. We claim that (with $v$ as in \eqref{16.5})
\begin{equation}\label{16.13}
v'(r) \geq \H^1(E\cap \S_r) = 2r x(r)
\end{equation}
almost everywhere on $(0,r_0)$; for a rapid proof with heavy material, apply the
co-area formula to $E$ and the function $x \to |x|$; for a slow one, see (5.8) in \cite{C1}.
Recall that since $f(r) = r^{-2} v(r) - \theta_0$ by \eqref{16.3}, 
Lemma \ref{t16.1} says that $f$ also is differentiable almost everywhere, with 
\begin{equation} \label{16.14}
r f'(r) =  r^{-1} v'(r) -  2 r^{-2} v(r) \geq 2x(r) - 2r^{-2} v(r)
\end{equation}
by \eqref{16.13}. Recall that by \eqref{16.12},
\begin{equation}\label{16.15}
r^{-1} v(r) = r \theta(r)
\leq r x(r) - \alpha [2rx(r)-2r\theta_0] + r (2\alpha q(r)+h(r)) .
\end{equation}
That is,
\begin{equation}\label{16.16}
r x(r) (1-2\alpha) \geq r^{-1} v(r) - 2 \alpha r \theta_0 - r (2\alpha q(r)+h(r))  
\end{equation}
or equivalently
\begin{equation}\label{16.17}
x(r) \geq \frac{v(r)}{(1-2\alpha) r^2} - \frac{2 \alpha \theta_0}{1-2\alpha} 
- \frac{2\alpha q(r)+h(r)}{1-2\alpha}.
\end{equation}
Then we return to \eqref{16.14}, replace, and get that
\begin{equation}\label{16.18}
\begin{aligned}
r f'(r) &\geq 2 x(r) - 2 r^{-2} v(r) 
\cr&\geq - 2 r^{-2} v(r) + \frac{2v(r)}{(1-2\alpha) r^2} 
- \frac{4 \alpha \theta_0}{1-2\alpha} - \frac{2(2\alpha q(r)+h(r))}{1-2\alpha}
\cr&
\geq \frac{4\alpha v(r)}{(1-2\alpha) r^2} - \frac{4 \alpha \theta_0}{1-2\alpha} 
- \frac{2(2\alpha q(r)+h(r))}{1-2\alpha}
\cr&
= \frac{4\alpha \theta(r)}{(1-2\alpha) } - \frac{4 \alpha \theta_0}{1-2\alpha} 
- \frac{2(2\alpha q(r)+h(r))}{1-2\alpha}
\cr&
\geq \frac{4\alpha}{(1-2\alpha) }\, f(r) - 3(2\alpha q(r)+h(r)) 
\end{aligned}
\end{equation} 
by \eqref{16.3} and because $\alpha$ is small (see \eqref{15.31}) and $q(r) \geq 0$. 
This proves \eqref{16.11}; the proposition follows.
\qed

\ms
We now make some additional comments on Proposition \ref{t16.2} and
then show how it may imply decay estimates; the true examples are in the next sections.

We decided not to require that $\theta_0$ is given by \eqref{16.4}, or that the
cones $X(r)$ have a density equal to $\theta_0$, but this will be our main example.

The proposition is also valid on an interval. That is, if $E$ is a reduced sliding almost minimal 
set (relative to $L$) in a domain $U$ that contains $B(0,2r_0)$, if \eqref{16.2} holds,
and if the assumptions \eqref{16.7}-\eqref{16.10} hold on an interval $(r_{00}, r_0)$, then 
\eqref{16.11} holds on $(r_{00}, r_0)$ too. The proof is the same.

\ms
The differential inequality \eqref{16.11} is not hard to integrate. 
Let $E$ and $r_0$ be as in Proposition \ref{t16.2}, and suppose in addition that the 
constants $c_{fl}(r)$, $\eta_{fl}(r)$, and $\eta_g(r)$ are such that
\begin{equation}\label{16.19}
\alpha(c_{fl}(r), \eta_{fl}(r), \eta_g(r)) \geq \alpha \ \text{ for almost every } 0 < r < r_0
\end{equation}
for some $\alpha > 0$ that does not depend on $r$. Then set 
\begin{equation}\label{16.20}
a = \frac{4\alpha}{1-2\alpha}
\end{equation}
and consider the auxiliary function $g(r) = r^{-a} f(r)$; \eqref{16.11} says that
\begin{equation}\label{16.21}
g'(r) = -a r^{-a-1} f(r) + r^{-a} f'(r) \geq -3 r^{-a-1} (h(r)+2\alpha q(r)),
\end{equation}
which we interpret as saying that $g$ is nearly nondecreasing.
And indeed, Lemma \ref{t16.1}  
says that for $0 < r_1 \leq r_2 < r_0$, 
\begin{equation}\label{16.22}
g(r_2) \geq g(r_1) - 3 \int_{r_1}^{r_2} (h(r)+2\alpha q(r)) \frac{dr}{r^{a+1}}
\end{equation}
or equivalently (since we are more often interested in letting $r_1$ tend to $0$),
\begin{equation}\label{16.23}
f(r_1) = r_1^a g(r_1) \leq \Big(\frac{r_1}{ r_2}\Big)^a f(r_2)
+ 3 r_1^a \int_{r_1}^{r_2} (h(r)+2\alpha q(r)) \frac{dr }{ r^{a+1}}.
\end{equation}
If the right-hand side cooperates, this says that $f(r_1)$ decays at some speed
when $r_1$ tends to $0$. For instance, if 
\begin{equation}\label{16.24}
h(r)+q(r) \leq C r^b \ \text{ for some } b < a
\end{equation}
(to simplify the computation), we get that near $0$,
\begin{equation}\label{16.25}
f(r_1) \leq \Big(\frac{r_1}{r_2}\Big)^a f(r_2) + C r_1^b
\end{equation}
with a constant $C$ that depends on $a$ and $b$, but not on $r_2$.

This will be good when we get it, and we will see examples in the next section. 
Then we will not be finished, because it will be much better to show that
some significant geometric quantities, rather than $f$ alone, decay near the origin. 
This will be the object of Part C (i.e., Sections \ref{S18}-\ref{S21}).

\section{Compactness, blow-up limits, and decay for $f$}
\label{S17}

In this section we fix the dimension $n$, the line $L$ through the origin, 
a sliding almost minimal set $E$ that contains $0$, and we use the compactness 
of the set $MC(L)$ of sliding minimal cones (with respect to $L$) to prove that
if in addition to the usual assumptions, all the blow-up limits of $E$ at $0$ 
satisfy the full length property, then the assumptions of Section \ref{S16} 
are satisfied for $r_0$ small.
See Proposition \ref{t17.1} for the ensuing statement. 

So we fix $n$, $L$, a radius $r_1 > 0$, and a closed set $E$ in $B(0,r_1)$, 
and assume that 
\begin{equation} \label{17.1}
\begin{aligned}
E \ &\text{ is a reduced sliding almost minimal set in } B(0,r_1), 
\cr&\hskip1cm \text{with boundary condition coming from } L, 
\end{aligned}
\end{equation}
with a gauge function $h$ such that
\begin{equation} \label{17.2}
h(r) \leq C_h r^\beta \ \text{ for } 0 < r \leq r_1
\end{equation}
for some constants $C_h \geq 0$, $\beta > 0$. 
Also we assume that
\begin{equation} \label{17.3}
0 \in E \cap L.
\end{equation}

Let us say a few words about $MC(L)$ before we discuss the blow-up limits
of $E$ at $0$, and then state the main result of this section. 
So far we have a definition of local Hausdorff convergence on 
closed subsets of $\R^n$, which is defined with the local Hausdorff distances
$d_{x,r}$ of \eqref{1.13}, and for which $\{ X_k \}$ converges to $X$
if $\lim_{k \to +\infty} d_{x,r}(X_k,X) = 0$ for every choice of $x\in \R^n$
and $r >0$, or equivalently for $x=0$ and every $r > 0$.
But since $MC(L)$ is composed of cones, 
$d_{0,1}(X,X') = d_{0,r}(X,X')$ for $X,X' \in MC(L)$ and $r > 0$, and it is enough
to use the ``distance function'' 
\begin{equation} \label{17.4}
\begin{aligned}
d_{0,1}(X,X') &= \sup\big\{ \dist(x,X') \, ; \, x \in X \cap B(0,1) \}
\cr&\hskip2cm
+ \sup\big\{ \dist(x',X) \, ; \, x' \in X' \cap B(0,1) \},
\end{aligned}
\end{equation}
for $X, X' \in MC(L)$. It is easy to see that 
$d_{0,1}(X,X')$ is also equivalent to the usual Hausdorff distance between
$K = X \cap \S$ and $K' = X' \cap \S$, defined by
\begin{equation} \label{17.5}
d_\H^c(X,X') := d_{\H}(K,K') = \sup\big\{ \dist(x,K') \, ; \, x \in K \}
+ \sup\big\{ \dist(x',K') \, ; \, x' \in K' \}.
\end{equation}
The small advantage of this is that it is well known that $d^c_{\H}$ is a distance
(i.e., in particular the triangular inequality holds with the constant $1$) on the set
of closed cones, and that some simple facts are very well known in this context. 
We shall use the fact that, with either of these distances, 
\begin{equation} \label{17.6}
MC(L) \ \text{ is a compact set.} 
\end{equation}
Given the fact that the set of closed subset of $\S$, with the Hausdorff distance $d_\H$,
is compact, this simply amounts to checking that if $X$ is the limit of the sequence
$\{ X_k \}$ in $MC(L)$, then $X \in MC(L)$ too. This is not trivial, but follows at once 
from Theorem 10.8 in \cite{Sliding}. 

\ms
Let us also recall some simple facts about blow-up limits.
Let $E$ be as above, and denote by $\cX$ the set of  blow-up limits 
of $E$ at $0$.  Recall that $\cX$ is the collection of sets $X$ such that
\begin{equation} \label{17.7}
X = \lim_{k \to +\infty} r_k^{-1} E
\end{equation}
for some sequence $\{ r_k \}$ of positive numbers that tends to $0$.
Recall also that this means that $\lim_{k \to +\infty} d_{0,R}(X,r_k^{-1} E) = 0$
for every $R > 0$, with $d_{0,R}$ as in \eqref{1.13}. Let us say why
\begin{equation} \label{17.8}
\cX \text{ is a closed subset of } MC(L). 
\end{equation}
The fact that if $X \in \cX$, then $X$ is a sliding minimal cone is Corollary 29.53
in \cite{Sliding}; we even get that the density of $X$ is 
\begin{equation} \label{17.9}
\H^2(X\cap B(0,1)) = \lim_{r \to 0} r^{-2} \H^2(E \cap B(0,r))
\end{equation}
(where the limit exist by near monotonicity of $r^{-2} \H^2(E \cap B(0,r))$,
as in Theorem 28.7 of \cite{Sliding}). So we just need to show that $\cX$
is closed.

Suppose $X$ is the limit in $MC(L)$ of the sequence $\{ X_j \}$ in $\cX$, and write
$X_j = \lim_{k \to +\infty} r_{j,k}^{-1} E$ for some sequence $\{ r_{j,k} \}$,
$k \geq 0$, that tends to $0$. By standard manipulations of sequence extraction,
we can find a sequence $\{ k(j) \}$, $j \geq 0$, such that $r_{j, k(j)}$ tends to $0$
and $X = \lim_{j \to +\infty} r_{j,k(j)}^{-1} E$. That is, $X\in \cX$; \eqref{17.8} follows.

We are ready to state the main result of this section.

\begin{pro} \label{t17.1}
Let the sliding minimal set $E$ satisfy the assumptions \eqref{17.1}, \eqref{17.2}, 
and \eqref{17.3}. Suppose in addition that
\begin{equation} \label{17.10}
\text{every blow-up limit of $E$ at $0$ satisfies the full length condition.}
\end{equation}
Then we can find $a \in (0,1)$ and a radius $r_0 \in (0, r_1]$ such that 
\begin{equation} \label{17.11}
r f'(r) \geq a f(r) - 3 h(r)
\ \text{ for } 0 < r < r_0,
\end{equation}
where $f(r)$ is still defined by \eqref{16.3}, but with 
\begin{equation} \label{17.12}
\theta_0 = \lim_{r \to 0} r^{-2} \H^2(E \cap B(0,r))
\end{equation}
as in \eqref{16.4}.
\end{pro}

Recall from Lemma \ref{t16.1} that we already knew that $f$ is differentiable almost
everywhere on $(0,r_0)$, and that we can partially recover the variations
of $f$ from $f'$. We will see how to use this after the proof. 

The proof will use Proposition \ref{t16.2} and a small compactness argument on $MC(L)$.
For each $X\in \cX$, the definition of full length gives a small constant $c = c(X)$
and a small radius $\eta = \eta_{fl}(X) \in (0,\eta(X))$ such that \eqref{3.25} holds
for every injective deformation parameter $\varphi \in \Phi^{+,i}_X(\eta)$ that satisfies
\eqref{3.24} (see Definition \ref{t3.1}). Then the construction of Sections \ref{S4}-\ref{S15}
also gives a small constant $\varepsilon(X)$, such that if the assumptions of Section \ref{S4},
and in particular \eqref{4.2} and \eqref{4.3}, are satisfied with $\varepsilon = \varepsilon(X)$,
we get the main conclusion \eqref{15.32}, with some small constant $\alpha = \alpha(X)$.
With the notation of the previous section, 
$\varepsilon(X) = \varepsilon(c_{fl},\eta_{fl},\eta_g)$ and $\alpha(X) = \alpha(c_{fl},\eta_{fl},\eta_g)$,
where $\eta_g = \eta(X)$, $\eta_{fl} = \eta_{fl}(X)$, and $c_{fl} = c(X)$; this notation has the advantage 
of making it plain that $\varepsilon(X)$ and $\alpha(X)$ depend only on the constants above.
We define a small ball $V_X$ in $MC(L)$ by
\begin{equation} \label{17.13}
V_X = \big\{ Y \in MC(L) \, ; \, d_\H^c(X,Y) < 10^{-1} \varepsilon(X)) \big\}.
\end{equation}
Now $\cX$ is a closed set in the compact $MC(L)$, so there is a finite set $\cX_0 \subset \cX$,
such that the $V_X$, $X\in \cX_0$, cover $\cX$. In other words,
\begin{equation} \label{17.14}
\text{ for $Y\in \cX$ we can find $X \in \cX_0$ such that } 
d_\H^c(X,Y) < 10^{-1} \varepsilon(X).
\end{equation}
We also need to know that 
\begin{equation}\label{17.15}
\lim_{r \to 0} \Big\{\inf_{X \in \cX} d_{0,3r}(E,X) \Big\} = 0.
\end{equation}
Suppose not. Then there is an $\varepsilon > 0$ and a sequence $\{ r_k \}$, that tends to $0$, such that
$d_{0,3r_k}(E,X) \geq \varepsilon$ for all $k$. We may replace $\{ r_k \}$ with a subsequence, for which
the sets $E_k = r_k^{-1} E$ converge, locally for the Hausdorff distance, but on the whole $\R^n$, to a 
closed set $X$. By definitions, $X \in \cX$, and by \eqref{17.8}, $X$ is a sliding minimal cone.
But the local convergence says that $d_{0,3}(E_k,X)$ tends to $0$, which contradicts the definition of $r_k$.

Set $\varepsilon_0 = \inf_{X \in \cX_0} \varepsilon(X)$, and let $r_0$ be such that 
\begin{equation}\label{17.16}
0 < r_0 < \frac{r_1}{3}  \ \text{ and }\ C_h r_0 \leq 10^{-1}\varepsilon_0,
\end{equation}
where $C_h$ is the same constant as in \eqref{17.2}, and
\begin{equation}\label{17.17}
\inf_{X \in \cX} d_{0,3r}(E,X) \leq 10^{-1}\varepsilon_0
\ \text{ for } 0 < r < r_0.
\end{equation}
We want to show that $r_0$ satisfies all the assumptions of Proposition \ref{t16.2}
with $U = B(0,r_1)$. 
This is clear for \eqref{16.1} and \eqref{16.2}. For the other assumptions, we fix $r \in (0,r_0)$
and we want to find a cone $X(r)$. But \eqref{17.17} gives a cone $X\in \cX$ such that
$d_{0,2r}(E,X) \leq 10^{-1}\varepsilon_0$, and then \eqref{17.14} yields an $X(r) \in \cX_0$
such that $d_\H^c(X(r),X) < 10^{-1} \varepsilon(X(r))$.
We take $\eta_g(r) = \eta(X(r))$, $\eta_{fl}(r) = \eta_{fl}(X(r))$, and $c_{fl}(r) = c(X(r))$,
and then $\varepsilon(c_{fl}X(r),\eta_{fl}X(r),\eta_gX(r)) =  \varepsilon(X(r))$ 
and $\alpha(c_{fl}X(r),\eta_{fl}X(r),\eta_gX(r)) =  \alpha(X(r))$ by the definitions above.

Then we need to check \eqref{16.7}, i.e., that $d_{0,2r}(E,X(r)) \leq \varepsilon(X(r))$,
and this easily follows from the definitions of $X$ and $X(r)$, plus the fact that 
$\varepsilon_0 \leq \varepsilon(X(r))$ since $X(r) \in \cX_0$. 
Next \eqref{16.8} holds by \eqref{17.16}, and \eqref{16.10} holds, with $q(r) = 0$,
because $X(r) \in \cX$ and by \eqref{17.9}. We apply Proposition \ref{t16.2} and get \eqref{16.11},
with $q(r) = 0$ and $\alpha = \alpha(X)$. But the  same proof would also yield \eqref{16.11} 
with the constant
\begin{equation}\label{17.18}
\alpha = \inf_{X\in \cX_0} \alpha(X) > 0.
\end{equation}
So we get \eqref{17.11} with $a = \frac{4\alpha}{1-2\alpha}$. 
We prefer to say things like this, because formally we do not know that \eqref{16.11} with some 
$\alpha$ implies \eqref{16.11} with any smaller $\alpha$;
even here with our special choice of $\theta_0 = \lim_{r \to 0} \theta(r)$, we do not know for sure 
that $f(r) \geq 0$,  because we only know that $\theta(r)$ is almost monotone.
On the other hand, if $f(r) \leq 0$, we should be happy anyway, because the goal of all the story 
is to show that $f(r)$ is small (but don't worry, we don't need this remark). 
This completes the proof of Proposition \ref{t17.1}.
\qed

\ms
Let us comment on Proposition \ref{t17.1}.
We had to be slightly careful, because with the proof above the constants $\varepsilon$ 
in \eqref{4.2} and \eqref{4.3} depend not only on the full length constant $c$ for $X$, 
but also on the more geometric constants $\eta(X)$ and $\eta_{fl}(X)$;
so we don't want to use cones $X(r)$ that come extremely close to $L \cap \S$ without actually meeting it, for instance.

Our proof of Proposition \ref{t17.1} relies on compactness, but in concrete cases, 
the covering of $\cX$ by balls $V_X$ can be obtained explicitly (and then we get a better control 
on the constant $C$). Suppose for instance that $\theta_0 = \lim_{r \to 0} \theta(r)$ is equal 
to $3\pi/2$; then $\cX$ is contained in the set of cones $X\in MC(L)$ such that 
$\theta(X) = 3\pi/2$. If in addition $n = 3$, say, we know that $\cX$ is exclusively 
composed of cones of type $\bY$. Now some of them contain $L$ in their spine, 
others don't but contain half of $L$ in one of their faces, and some meet $L$ only at $0$. 
Given $r > 0$ as above, and if $X \in \cX$ approximates $E$ well in $B(0,3r)$, 
we choose to take $X(r) = X$ if $X$ is of the first type, but otherwise we will replace $X$ with
an $X(r)$ of the first type if its spine is very close to $L$, and an $X(r)$ of the second type if 
half of $L$ is very close to some face of $X$, but the other one is reasonably far from $X$. 
In this case, the full length property is proved in Section \ref{S30} below. % thm30.1 en fait
Of course this concrete way of proving Proposition \ref{t17.1} is harder to do when we don't 
know well the list of minimal cones of density $\theta_0$, not to mention the fact that 
we cannot be sure that they satisfy the full length property.

Once we have \eqref{17.11}, with a constant $a>0$ that does not depend on $r$,
but only on $\cX$, we can use Lemma \ref{t16.1} to integrate it and get the decay 
estimate \eqref{16.22}. Here $q(r) = 0$, so we get that for $0 < r < s < r_0$,
\begin{equation} \label{17.19}
r^{-a} f(r) \leq s^{-a} f(s) + 3 \int_{r}^{s} \frac{h(r) dr}{r^{a+1}}
\leq s^{-a} f(s) + 3 C_h \int_{r}^{s} r^{\beta - a - 1} dr.
\end{equation}
We may as well assume that $a < \beta$ (in fact, $a$ depends on $\beta$, 
we expect it to be very small, and anyway we can always make it smaller); 
then \eqref{17.19} yields
\begin{equation} \label{17.20}
r^{-a} f(r) \leq s^{-a} f(s) + \frac{3 C_h}{\beta-a} \, s^{\beta - a} dr
\ \text{ for } 0 < r < s < r_0.
\end{equation}
This is good: for $s \in (0,r_0)$ fixed, this means that $f$ decays like $r^a$ near $0$.

\ms
We end this section with a corollary of the discussion above.

\begin{cor} \label{t17.2} 
Let $E$ satisfy \eqref{17.1}-\eqref{17.3}, let $\theta_0$ be as in \eqref{17.12},
and define $f$ by \eqref{16.3}. If
\begin{equation} \label{17.21}
\text{every blow-up limit of $E$ at $0$ satisfies the full length condition,}
\end{equation}
there exist a constant $a > 0$ and numbers $r_0 > 0$ and $C \geq 0$ such that
\begin{equation} \label{17.22}
f(r) \leq C r^{a} \ \text{ for } 0 < r < r_0.
\end{equation}
The constant $a$ depends only on $n$ and a full length constant coming from the family of 
blow-up limits of $E$ at $0$. But $r_0$ and $C$ depend on the specific situation (and in particular $E$).
\end{cor}

\ms
Indeed, the assumptions of Proposition \ref{t17.1} are satisfied, so we can find 
$r_0 > 0$ and $a \in (0,1)$ (that depends on the class $\cX$ of blow-up limits of $E$ at $0$,
in particular through the full length constants of a finite number of cones used for a covering) 
such that \eqref{17.11} holds. The estimate \eqref{17.22} now follows from \eqref{17.20} 
(we just dropped the more explicit computation of constants) and the discussion that leads to it.
\qed

\ms
We would like to say that the assumption \eqref{17.21} holds automatically when 
$\theta_0 \leq \frac{3\pi}{2}$, but for this we would need to know that 
\begin{equation} \label{17.23}
\text{if $X \in MC(L)$ and $\theta(X) \leq \frac{3\pi}{2}$, then }
X \in \bH \cup \bV \cup \bY,
\end{equation}
where $\bH$, $\bV$, and $\bY$ are as below Subsection \ref{S1.2} and define the 
same cones as in Theorem~\ref{t1.3}.

This looks reasonable, but the author did not find a simple proof. 
But see Lemma \ref{t22.2} for the simpler special case when $\theta(X) \leq \pi + \varepsilon_n$.
As soon as we can prove \eqref{17.23}, we observe that if $\theta_0 \leq \frac{3\pi}{2}$, 
\eqref{17.23} says that every blow-up limit of $E$ at $0$ lies in $\bH \cup \bV \cup \bY$, 
hence satisfies the full length property by Theorem \ref{t30.1}, and we can apply Corollary \ref{t17.2}.

\begin{rem} \label{r17.2}
There will be a better statement, Corollary \ref{t22.1n}, 
where we only assume that some blow-up limit of $E$ at $0$ satisfies the full length condition, 
but it will be harder to prove, because we need to be able to find good approximating cones $Z(r)$ 
at all the scales $r < r_0$, so that we can apply Proposition \ref{t16.2}. 
For this we will need the approximation results of the next part.
\end{rem}

\vfill\eject
\part{Approximation by cones for balls centered on $L$}

\section{The density excess controls the distance to a cone}
\label{S18}

In this part we still consider balls centered at $0 \in E \cap L$,
assume that the density excess $f(r)$ is small, and use this to obtain geometric information
on $E$, and in particular its Hausdorff distance to minimal cones, first on most spheres,
and then on thicker annuli. 
The rough idea is that if $f(r)$ is small, it cannot vary much between $r/2$ and $r$, hence
$f'(\rho)$ is often small on $[r/2,r]$, and the proof of the differential inequality \eqref{17.11}
will allow us to control various quantities when $f'(\rho)$ is small.

We are given a reduced sliding almost minimal set $E$ of dimension $2$, in an open set
$U \subset \R^n$ which contains the origin, with a sliding condition that comes from the line $L$ 
through $0$ and a small enough gauge function $h$. We suppose that $0 \in E$, set
\begin{equation}\label{18.1}
\theta_0 = \lim_{t \to 0} t^{-2} \H^2(E\cap B(0,t))
\end{equation}
as in \eqref{16.4} (we shall soon remind the reader of why it exists when $h$ is small enough), and
\begin{equation}\label{18.2}
f(r) = \theta(r) - \theta_0 = r^{-2} \H^2(E\cap B(0,r)) - \theta_0
\end{equation}
for $r < \dist(0,\R^n \sm U)$ (as in \eqref{16.3}). 
We want to show that $f(r)$ controls the local Hausdorff distance from $E$ 
to small modifications of minimal cones. We will roughly proceed as in 
Section 11 of \cite{C1}, where we established this for two-dimensional almost minimal
sets with no sliding boundary condition.

We start with a discussion of the list of modifications of minimal cones that we allow,
and how we measure the distances.

Let us first consider a fixed minimal cone $X$, and use the (in fact, any) standard decomposition of  
$K = X\cap \partial B(0,1)$ into arcs of circles $\cC_i$, $i\in \cI$, that was described in  
Section \ref{S3}. We consider deformations of $K$, which we construct as for the definition
of the full length condition near Definition \ref{t3.1}. That is, we select a small constant $\eta > 0$
(for instance choose any number smaller than $\eta(X)$ in \eqref{3.12}; 
the actual choice won't matter), and we define the extended class $\Phi_X^+(\eta)$ of 
enlarged mappings $\varphi$, as near \eqref{3.21}.
Most of the information of $\varphi$ is a mapping defined on the set of endpoints of the $\cC_i$,
which says where we send each one, but $\varphi$ also contains some information relative to the way we
glue the pieces near vertices of $L \cap \S$. For each $\varphi \in \Phi_X^+(\eta)$ we define a set
$\varphi_\ast(K)$, which is the deformation of $K$ associated to $\varphi$, 
and the cone $\varphi_\ast(X)$ over $\varphi_\ast(K)$. 
Recall that modulo some small modifications of the protocole near the points of 
$L \cap \S$, $\varphi_\ast(K)$ is obtained by replacing each arc $\cC_i = \rho(a_i,b_i)$ of $K$ by
the arc $\rho(\varphi(a_i),\varphi(b_i))$.

Let us denote by $\cZ(X,\eta)$ the set of cones $Z = \varphi_\ast(X)$, where
$\varphi \in \Phi_X^+(\eta)$. 
These are not  exactly minimal cones, because the angles 
between the $\cC_i$, for instance, may have changed a bit, but they are close to $X$ if 
$\eta$ is small enough (which we can assume). 
For some estimates, it may be interesting to measure how far they are from
being minimal, so we introduce a number $\alpha(Z)$ which records this. In Section 11 of \cite{C1},
we have used a partial measurement of minimality based on the angles made by the geodesics that compose 
$Z$; here we find it more pleasant not to describe these angles (the distances of the edges to $L$
should also be taken into account), and measure the lack of minimality more directly (but less
geometrically). For $Z \in \cZ(X,\eta)$,  set 
\begin{equation}\label{18.3}
\alpha(Z) = \inf\big\{ \H^2(Z \cap \ol B) - \H^2(\wt Z \cap \ol B)\, ; \, 
\wt Z \text{ is a sliding competitor for $Z$ in $\ol B$ }\big\},
\end{equation}
where $\ol B = \ol B(0,1)$ and the notion of sliding competitor in $\ol B(0,1)$ is explained in 
Definition~\ref{t1.1}.

We also want to allow $X$ to vary, so we let $\cX$ denote a class of sliding minimal cones
centered at the origin; for instance, we shall use 
\begin{equation}\label{18.4}
\cX(\theta_0) = \big\{ X \, ; \, X \text{ is a reduced sliding minimal cone centered at $0$, with }
\H^2(X\cap B) = \theta_0 \big\}.
\end{equation}

Then we fix a small number $\eta > 0$ and set
\begin{equation}\label{18.5}
\beta_{\cX,\eta}(E,r) = \inf\big\{ d_{0,r}(E,Z) + \alpha(Z)^{1/4}  \, ; \, 
Z \in \cZ(X,\eta) \text{ for some } X \in \cX \big\},
\end{equation}
where the local distance $d_{0,r}$ is still as in \eqref{1.13}, and we put a power $1/4$
to simplify the statement of the next theorem without losing too much information
(notice that with this definition, $\beta_{\cX,\eta}(E,r)$ tends to be larger).

The next result summarize what we want to do in the next sections. We state it in a 
normalized ball to simplify some expressions (such as the precise form of \eqref{18.6}).

\begin{thm}\label{t18.1} 
Let $E$ be a sliding almost minimal set in an open set $U \subset \R^n$ which contains
$B(0,400)$, with sliding conditions coming from the line $L$ through the origin, 
and with a gauge function $h$ such that
\begin{equation}\label{18.6}
h(t) \leq C_0 t^{\beta_0} \ \text{ for } 0 < t < 400
\end{equation}
for some constants $\beta_0 \in (0,1]$ and $C_0 \geq 0$.
Suppose $\theta_0$ is as in \eqref{18.1}, let $\eta > 0$ be given, and let 
$\cX = \cX(\theta_0)$ be as in \eqref{18.4}.
Assume that $C_0$ is small enough, depending on $n$, $\eta$, and $\theta_0$.
Then
\begin{equation}\label{18.7}
\beta_{\cX,\eta}(E,1) \leq C \Big[ f(200) + \int_0^{400} \frac{h(t) dt}{t} \Big]^{1/4},
\end{equation}
where $f$ is as in \eqref{18.2}, and the constant $C$ depends only on $\beta_0$, $C_0$,  
$\theta_0$, $\eta$, and $n$. 
\end{thm}

\ms 
See Theorem 11.4 in \cite{C1} for the analogous statement away from the boundary.
The power $1/4$ is certainly not optimal, and the same sort of contortion as in
\cite{C1} should probably lead to the power $1/3$. See the end of the proof of Lemma \ref{t20.2}
for this.
But $1/3$ is probably not optimal either, and $1/2$ would sound more right; 
we know that in the proof below (and without the possible improvement on
Lemma \ref{t20.2}), we will probably lose something significant when we go from 
isolated estimates on spheres to a global estimate on the ball. 
Similarly, we noted a dependence on $\theta_0$, because this
is the way that we shall prove things, but probably a more clever argument would allow us to get
rid of this.

Here we decided to assume a geometric decay in \eqref{18.6}, but a weaker condition (probably
a Dini condition) would be enough. We decided for a complicated way to state \eqref{18.6},
where $C$ has some dependence on the constants $C_0$ and $\beta_0$, but where we allow 
the possibility that $\int_0^{400} \frac{h(t) dt}{t}$ is smaller than suggested by \eqref{18.5} 
and then we get a better estimate.

Notice that since we proved in the earlier sections that $f(r)$ often decays like a power of $r$,
the theorem will imply a similar decay of $\beta_{\cX,\eta}(E,r)$. For the moment, we allow
the reference minimal cone $X$ in the computation of $\beta_{\cX,\eta}(E,r)$ to vary with
$r$, but once we get a power decay, we will know that we can take for $X$ any blow-up
limit of $E$ at $0$, and this will imply the uniqueness of the blow-up limit in question.
But in the mean time it is better to allow $X$ to vary. On the opposite side, we could have allowed 
$\cX$ to be the whole class of minimal cones, but then \eqref{18.7} would have been
less precise.

\ms
Before we start the proof (which will take some time), let us record that 
it is enough to prove \eqref{18.7} when
\begin{equation}\label{18.8}
f(200) + \int_0^{400} \frac{h(t) dt}{t} \leq \varepsilon_1,
\end{equation}
where the very small constant $\varepsilon_1 > 0$ will be chosen later 
(depending on $n$, $\theta_0$ and $\eta$). Indeed, if \eqref{18.8} fails, then \eqref{18.7} holds with
$C = 2\varepsilon_1^{-1}$ (try for $X$ a bow-up limit of $E$ at the origin, and observe that
$d_{0,r}(E,X) \leq 2$). So let us assume that \eqref{18.8} holds.

The following lemma will allow us to use the same construction of competitors as in the previous section.

\begin{lem} \label{t18.2} Let $\tau_1 > 0$ be small. If \eqref{18.8} hold for a small enough
$\varepsilon_1 > 0$ (that depends also on $\theta_0$ and $n$), we can find a minimal cone
$X\in \cX(\theta_0)$ such that
\begin{equation}\label{18.9}
d_{0,180}(E,X) \leq \tau_1.
\end{equation}
\end{lem}

So let $E$ be as in the theorem, and assume \eqref{18.8}.
We want to show that because of \eqref{18.8}, the density $\theta(r) = r^{-2}\H^2(E\cap B(0,r))$ 
is nearly constant near $0$, and then use this to show that $E$ looks a lot like a minimal cone in $B(0,180)$. 
 
We start with the near monotonicity of $\theta$. Recall from Theorem 28.7 in \cite{Sliding} that
there is a constant $\alpha_n$, which here depends only on $n$ (because the geometry of $L$ is simple),
such that 
\begin{equation}\label{18.10}
\theta(r) {\rm exp}\Big(\alpha_n \int_0^r \frac{h(2t) dt}{t} \Big) 
\ \text{ is nondecreasing on } (0,200).
\end{equation}
Then \eqref{18.1} and \eqref{18.8} imply that
\begin{equation}\label{18.11}
\theta(r) \geq \theta_0 \, {\rm exp}\Big(-\alpha_n \int_0^r \frac{h(2t) dt}{t}\Big) 
\geq \theta_0  \, e^{-\alpha_n \varepsilon_1} 
\ \text{ for } 0 < r \leq  200
\end{equation}
(where the endpoint $r=200$ is obtained by taking a limit), so that
\begin{equation}\label{18.12}
\theta_0 e^{-\alpha_n \varepsilon_1} \leq \theta(200) = \theta_0 + f(200)
\leq \theta_0 + \varepsilon_1
\end{equation}
by \eqref{18.2} and \eqref{18.8}. We deduce from this and \eqref{18.10} again that
\begin{equation}\label{18.13}
\theta(r) \leq e^{\alpha_n \varepsilon_1} \theta(200)
\leq e^{\alpha_n \varepsilon_1}[\theta_0 + \varepsilon_1]
\ \text{ for } 0 < r < 200.
\end{equation}

Let us now apply an almost constant density result from \cite{Sliding} to say that
$E$ looks like a minimal cone. Let $\tau > 0$ be small, to be chosen soon,
and let $\varepsilon > 0$ be associated to $\tau$ as in Proposition 30.19 in \cite{Sliding}.
We want to apply that proposition to $E$ and the radii $r_2 = r_0 = 200$. 
The bilipschitz assumption on the boundary $L$ (up to (30.20) in \cite{Sliding}) is trivially satisfied,
(30.21) holds if $\varepsilon_1$ is small enough and because \eqref{18.8} controls $h(300)$,
and the more important near constant density assumption (30.22) holds because 
$\theta(r_2) = \theta(200) \leq \theta_0 + \varepsilon_1$ while
$\theta(r) \geq e^{-\alpha_n \varepsilon_1} \theta_0$ for $0 < r < 10^{-3} r_0$.

By Proposition 30.19 in \cite{Sliding}, there is a sliding minimal cone $T$ such that 
\begin{equation} \label{18.14}
d_{0,190}(E,T) \leq 2\tau,
\end{equation}
(see (30.24) and (30.24) there), and
\begin{equation} \label{18.15}
|\H^2(E\cap B(0,r))- \H^2(T\cap B(0,r))| \leq 200^2 \tau 
\ \text{ for } 0 \leq r \leq 190
\end{equation}
(see (30.26)). We apply this to $r=190$, then use \eqref{18.11} and \eqref{18.13}
to estimate $\theta(r)$, and get that 
\begin{equation}\label{18.16}
\begin{aligned}
|\H^2(T\cap B(0,1)) - \theta_0| &\leq |\theta(190)-\theta_0| + (200/190)^2 \tau
\cr&
\leq [e^{\alpha_n \varepsilon_1}-1] \theta_0 + 2\varepsilon_1 + (200/190)^2 \tau \leq 2\tau
\end{aligned}
\end{equation}
if $\varepsilon_1$ is small enough.

We cannot use $X=T$ in the lemma, because the density $\H^2(T\cap B(0,1))$
may be a little different from $\theta_0$. So we shall use \eqref{18.16} and a compactness
argument to find $X \in \cX(\theta_0)$ that is very close to $T$, and then deduce
\eqref{18.9} from \eqref{18.14}.

We claim that for all small $\tau_1 > 0$, we can choose $\tau > 0$ such that if 
$T$ is a sliding minimal cone such that \eqref{18.16} holds, then there is a minimal
cone $X\in \cX(\theta_0)$ such that 
\begin{equation}\label{18.17}
d_{0,1}(X,T) \leq \frac{\tau_1}{4}.
\end{equation}
Indeed, otherwise there is a sequence $\{ T_k \}$ of sliding minimal cones such that $T_k$ 
satisfies \eqref{18.16} with $\tau = 2^{-k}$, and yet $d_{0,1}(X,T_k) > \frac{\tau_1}{4}$
for $X\in \cX(\theta_0)$. We can extract a subsequence
(which we still denote the same way), for which $T_k$ converges to a limit $X$ in local Hausdorff
distance (say, for $d_{0,1}$, and see the discussion above \eqref{17.6} if you are worried about which 
notion of convergence). By the compactness property \eqref{17.6}, or directly by 
Theorem 10.8 in \cite{Sliding}, $X$ is a sliding minimal cone. 

By Theorem 10.97 in \cite{Sliding} (the lower semicontinuity of $\H^d$ along sequences of 
quasiminimal sets) and Theorem 22.1 in \cite{Sliding} (the upper semicontinuity along 
sequences of almost minimal sets), plus the fact that $\H^2(X \cap \S) = 0$, we get that 
\begin{equation}\label{18.18}
\theta_0 = \lim_{k \to +\infty} \H^2(T_k \cap B) = \H^2(X\cap B),
\end{equation}
so $X \in \cX(\theta_0)$ and this contradicts the fact that the $T_k$ were chosen far from 
$\cX(\theta_0)$. This proves our claim.

We choose $\tau$ with this property, also also such that $\tau \leq 10^{-1}\tau_1$,
apply this to the cone $T$ of \eqref{18.16}, and find $X \in \cX(\theta_0)$ that satisfies
\eqref{18.17}. It is easy to see that $X$ satisfies \eqref{18.9}; the lemma follows.
\qed

\ms
Now we want to apply the construction of Sections \ref{S4}-\ref{S15} to find cones
$Z(r)$, $0 < r < 180$, that approximate $E$ well on the circles $\S_r = \d B(0,r)$.
Since the dependance on $X$ of the constants $\tau$, and then $\varepsilon$ in \eqref{4.3} 
seems to be a little shady at first sight, let us spend some time discussing these constants.

\begin{rem} \label{R18.3}
We claim that we can apply the construction of Sections \ref{S4}-\ref{S15}, 
with a value of the various constants $\varepsilon$ and $\tau$ that depends only on $\beta_0$,
$C_0$, $\theta_0$, $\eta$, and $n$.

To see this we shall use the same compactness argument as in Section \ref{S17}, below \eqref{17.12}.
To each $X\in \cX(\theta_0)$ we can associate a standard decomposition as in Section \ref{S3} and 
a small number $\eta(X) > 0$, that satisfies the requirements of Section \ref{S3}. 
Let us even choose $\eta(X)$ somewhat smaller than the constant $\eta$ of the statement of
Theorem \ref{t18.1}. 
This is a brutal way to make sure that the deformations of $X$ that we construct later
will come from $\varphi \in \Phi_X^+(\eta)$.

Then there is a small number $\varepsilon(X)$, that depends on $\eta(X)$, such that 
if \eqref{4.2} and \eqref{4.3} are satisfied with $\varepsilon(X)$, then we can apply the 
construction of Sections \ref{S4}-\ref{S15}, except that we do not intend to assume, or use,
the full length property.

Recall from \eqref{17.6} that the class $MC(L)$ of sliding minimal cones is compact.
The proof of existence for \eqref{18.17} shows that the extra condition 
$\H^2(X\cap B) = \theta_0$ that defines $\cX(\theta_0)$ is closed, so  
$\cX(\theta_0) \subset MC(L)$ is compact too.

For $X \in \cX(\theta_0)$ define the small ball $V_X$ centered at $X$ and with radius
$10^{-5}\varepsilon(X)$, essentially as we did in \eqref{17.13}, except that now we are 
only interested in $V_X \cap \cX(\theta_0)$. By compactness, 
we can find a finite family $\cX_0 \subset \cX(\theta_0)$
such that the sets $V_X$, $X \in \cX_0$, cover $\cX(\theta_0)$.

Then set $\varepsilon_0 = 10^{-3} \inf_{X \in \cX_0} \varepsilon(X)$,
and apply Lemma \ref{t18.2} with $\tau_1 = \varepsilon_0$. 
This gives a small constant $\varepsilon_1$, and if \eqref{18.8} holds with this $\varepsilon_1$,
we know that we can find $X \in \cX(\theta_0)$ such that 
$d_{0,180}(E,X') \leq \varepsilon_0$. Then $X'$ lies in a ball $V_X$, $X \in \cX_0$,
and this implies that $d_{0,180}(E,X) \leq 2\varepsilon_0$ (compare with \eqref{17.5}, 
but also note that we can modify the constant $10^{-5}$ above as we want).
We took $\varepsilon_0$ so small because this way
\begin{equation}\label{18.19}
d_{0,2r}(E,X) \leq \frac{180}{2r} d_{0,180}(E,X) \leq \varepsilon(X)
\ \text{ for } 10^{-2} \leq r \leq 90,
\end{equation}
so \eqref{4.3} holds for these $r$. Because of our assumption that $C_0$ is small enough,
depending on $n$, $\eta$, and $\theta_0$, we also have \eqref{4.2}, and so if we want to apply the construction of Sections~\ref{S4}-\ref{S15} (excluding the full length property) to $r$, 
we shall just need to check the conditions \eqref{4.4}-\eqref{4.8}. 
We shall see now that they hold for almost every $r\in [10^{-2}, 90]$.
\end{rem}

\ms
Return to the proof of Theorem \ref{t18.1}. Pick $\varepsilon_1$, and then 
a cone $X \in \cX_0 = \cX_0(\theta_0)$, as in the remark above. Set 
\begin{equation}\label{18.20}
v(r) = {\H}^2(E \cap B(0,r)) = r^2 \theta(r)
\end{equation}
for $0 < r < 200$ and denote by ${\cal R}$ the set of radii $r \in (10^{-2},90)$
such that $\theta$ and $v$ are differentiable at $r$, 
\begin{equation} \label{18.21}
\theta'(r) = r^{-2} v'(r) - 2 r^{-3} v(r),
\end{equation}
\begin{equation}\label{18.22}
v'(r) \geq {\H}^1(E\cap \partial B(0,r))
\end{equation}
(which incidentally implies \eqref{4.4}), and in addition \eqref{4.7} and \eqref{4.8} hold. Then 
\begin{equation} \label{18.23}
\H^1((10^{-2},90) \setminus {\cal R}) = 0
\end{equation}
by Lemma \ref{t16.1}, \eqref{16.13},
the Hardy-Littlewood maximal theorem (see below \eqref{16.11}),
and the proof of Lemma 4.12 in \cite{C1} (see below \eqref{4.8}).
Now each $r\in \cR$ satisfies the constraints \eqref{4.3}-\eqref{4.8}, and 
we can apply to it the results of Sections \ref{S4}-\ref{S15}, excluding those that use
the full length property.

Observe that we shall use the same cone $X$ for all the radii $r\in \cR$; for smaller values of $r$, 
we could still use the same proof, but we would need to apply Lemma \ref{t18.2} to a different
radius, get another cone $X'$, and possibly a different type of deformation $Z(r)$ in the argument below.
But we do not need to do this for the moment.

We intend to use the results of the previous sections to get information on 
$E \cap \S_r$, where we set $\S_r = \d B(0,r)$, for $r\in \cR$. 
We shall be able to get better estimates when $j(r)$ is small, where
\begin{eqnarray}\label{18.24} 
j(r) &=& rf'(r) + f(r) + (1+2 \theta_0\alpha_n) h(2r) 
+ (1+\theta_0\alpha_n) \int_0^r \frac{h(2t) dt}{t}
\nn\\
&=& r\theta'(r) + f(r) + (1+2 \theta_0\alpha_n) h(2r) 
+ (1+\theta_0\alpha_n) \int_0^r \frac{h(2t) dt}{t}. 
\end{eqnarray}
We added the constant $\alpha_n$ (from \eqref{18.11})
because we want to make sure that $j(r) \geq 0$.
Indeed, it could be that $\theta'(r)$ and $f(r)$ are slightly negative. Nonetheless, it follows from
\eqref{18.11} that
\begin{equation}\label{19.25n}
f(r) = \theta(r) - \theta_0 
\geq \theta_0 \, \Big({\rm exp}\Big(-\alpha_n \int_0^r \frac{h(2t) dt}{t}\Big) - 1\Big)
\geq - \theta_0 \alpha_n \int_0^r \frac{h(2t) dt}{t}
\end{equation}
(because \eqref{18.8} says that $\int_0^r \frac{h(2t) dt}{t} \leq \varepsilon_1$ is small). Similarly,
it follows from \eqref{18.10} (and in fact this is the way \eqref{18.10} is proved) that for $r\in \cR$;
\begin{equation}\label{19.26n}
r f'(r) = r \theta'(r) \geq - \alpha_n \theta(r) h(2r) \geq - 2 \theta_0 \alpha_n h(2r)
\end{equation}
by \eqref{18.13} and if $\varepsilon_1$ is chosen small enough (recall that $\theta_0 \geq \pi$). Thus
\begin{equation}\label{19.27n}
j(r) \geq (r \theta'(r))_+ + f(r)_+ + h(2r) + \int_0^r \frac{h(2t) dt}{t} \ \text{ for } r\in \cR,
\end{equation}
with positive parts, which will be simpler to use than \eqref{18.24} for some estimates.

\begin{lem} \label{t18.4}
For $r \in \cR$, we can find a compact set $\gamma^\ast(r) \subset E \cap \S_r$
and a cone $Z(r) \in \cZ(X,\eta)$ such that
\begin{equation} \label{18.25}
\H^1(E \cap \S_r \sm \gamma^\ast(r)) \leq C j(r)
\end{equation}
and
\begin{equation} \label{18.26}
\sup_{z\in Z(r) \cap \S_r} \dist(z,\gamma^\ast(r)) 
+ \sup_{z\in \gamma^\ast(r)} \dist(z, Z(r) \cap \S_r)
\leq C r^{1/2} j(r)^{1/2}.
\end{equation}
\end{lem}

We will see other properties of $Z(r)$ along the way. 
The constant $C$ in \eqref{18.25} depends on $\eta$ too, through the choice of $\eta(X)$
in Remark \ref{R18.3}. 
We added $r^{1/2}$ in \eqref{18.26} to show the homogeneity,
but this was not needed because $r \in \cR \subset (10^{-2},90)$.

Let $j_0 >0$ be small, to be chosen later. We shall keep in mind that it is enough to
prove the conclusion of the lemma when
\begin{equation}\label{18.27}
r \in \cR \ \text{ and } j(r) < j_0,
\end{equation}
where the small constant $j_0 > 0$ will be chosen later, depending also
on $\eta$. Indeed otherwise we just take $Z(r) = X$ and $\gamma^\ast(r) = E \cap \S_r$, 
and the conclusion holds because the Hausdorff distance between 
$E \cap \S_r$ and $X \cap \S_r$ is bounded.

So let $r\in \cal R$ be given. First notice that
\begin{equation} \label{18.28}
\begin{aligned}
\H^1(E \cap \S_r) &\leq v'(r) = r^2 \theta'(r) +2r^{-1} v(r)
= r^2 \theta'(r) +2r \theta(r)  
\cr&
= r^2 f'(r) +2r f(r) +2 r \theta_0 \leq 2r j(r) + 2 r \theta_0
\end{aligned}
\end{equation}
by \eqref{18.22}, \eqref{18.21}, \eqref{18.20}, because $f(r) = \theta(r)-\theta_0$, and finally
by \eqref{18.24} (or rather \eqref{19.27n}, because of the strange case when $f'(r) < 0$).

Recall that for $r\in {\cal R}$, we can apply Sections \ref{S5}-\ref{S15}, where we constructed a 
few competitors and used them to prove estimates on $\H^2(E \cap B(x,r))$. 
In particular we have \eqref{14.46}, which says that if we normalized everything so that $r=1$,
\begin{equation} \label{18.29}
\H^2(E \cap \ol B(0,1)) \leq \frac{1}{2} \H^1(E \cap \S) 
- 10^{-5} [\H^1(E\cap \S) - \H^1(\rho^\ast)] +h(1) .
\end{equation}
Here $\rho^\ast = \rho^\ast(r)$ is the net of geodesic that was chosen during the proof; 
see \eqref{14.43}. 

Let us observe that this estimate was obtained without using the full length property.
It will be all right to use this if we do not want to include estimates on $\alpha(Z)$ in \eqref{18.3}; 
otherwise we will need to correct the estimate as we did in Section \ref{S15}. This will be done later,
but for the moment we ignore this and work with \eqref{18.29}. Let us rewrite it 
without the normalization by $r=1$. We get that for $r\in \cR$ 
\begin{equation} \label{18.30}
\H^2(E \cap \ol B(0,r)) \leq \frac{r}{2} \H^1(E \cap \S_r) 
- 10^{-5} r [\H^1(E\cap \S_r) - \H^1(\rho^\ast(r))] +r^2 h(r).
\end{equation}
The author feels a little better with the powers of $r$ because they give the homogeneity,
but here $r \in \cR \subset [10^{-2},90]$ so we should not need to worry much. Next we write
\begin{equation} \label{18.31}
\H^1(E\cap \S_r) - \H^1(\rho^\ast(r)) = \Delta_0(r) + \Delta_1(r) + \Delta_2(r),
\end{equation}
where the $\Delta_i(r)$ are defined as follows. First
\begin{equation}\label{18.32}
\Delta_0(r) = \H^1(E\cap \S_r) - \H^1(\gamma^\ast(r))
= \H^1(E \cap \S_r\sm \gamma^\ast(r)),
\end{equation}
where $\gamma^\ast(r)$ is the union of the the curves $\gamma$ that were
constructed in Section \ref{S6}. The identity in \eqref{18.32} holds because the 
curves $\gamma$ are contained in $E$; also, it could be observed that the only
contribution comes from the two small disks $D$ near the points of $L\cap \S_r$,
because outside of these disks, $E\cap \S_r$ is composed of nice curves.
Next
\begin{equation} \label{18.33}
\Delta_1(r) = \H^1(\gamma^\ast(r)) - \H^1(\Gamma^\ast(r)) \geq 0,
\end{equation}
where $\Gamma^\ast(r)$ is the union of the Lipschitz graphs $\Gamma_j$ that we construct in 
Sections \ref{S7}-\ref{S11}; the fact that $\Delta_1(r) \geq 0$ comes by adding up its analogue for 
each configuration; see the comment below \eqref{9.7}. Notice also that $\Delta_1(r)$
is like $\Delta_1$ in \eqref{14.42}. Finally
\begin{equation}\label{18.34}
\Delta_2(r) = \H^1(\Gamma^\ast(r)) - \H^1(\rho^\ast(r)) \geq 0,
\end{equation}
because $\rho^\ast(r)$ simply consists in replacing each arc $\Gamma$ of $\Gamma^\ast$
with the geodesic $\rho$ with the same endpoints; this is the same as $\Delta_2$ in \eqref{14.44} 
(and $\rho^\ast$ is defined by \eqref{14.43}). 
Thus \eqref{18.31} is essentially the same as \eqref{14.45}.
Also \eqref{18.30} can be rewritten as
\begin{equation}\label{19.38n}
\Delta_0(r) + \Delta_1(r) + \Delta_2(r) 
\leq 10^5 \Big[\frac{1}{2}\H^1(E \cap \S_r) - r^{-1}\H^2(E\cap \ol B(0,r))\Big] + 10^5 r h(r)
\end{equation}
and since
\begin{eqnarray}\label{19.39n}
\H^1(E \cap \S_r) \leq v'(r) &=& r^2 \theta'(r) + 2 r^{-1} v(r) = r^2 \theta'(r) + 2 r \theta(r)
\nn\\
&=& r^2 \theta'(r) + 2 r f(r) + 2 r \theta_0 \leq 2 r j(r) + 2 r \theta_0
\end{eqnarray}
by \eqref{18.22}, \eqref{18.21}, and \eqref{19.27n}, \eqref{19.38n} yields 
\begin{equation}\label{18.35}
\Delta_0(r) + \Delta_1(r) + \Delta_2(r) 
\leq 10^5 \big[r j(r) + r \theta_0 - r^{-1}\H^2(E\cap \ol B(0,r))\big] + 10^5 r h(r).
\end{equation}
But
\begin{equation}\label{18.36}
r^{-1}\H^2(E\cap \ol B(0,r)) \geq r \theta(r) 
\geq r \theta_0 \Big(1 -\alpha_n \int_0^r  \frac{h(2t) dt}{t}\Big)
\end{equation}
by \eqref{19.25n}, so
\begin{equation}\label{18.37}
\Delta_0(r) + \Delta_1(r) + \Delta_2(r) 
\leq 10^5 r \, \Big[  j(r) + \theta_0 \alpha_n \int_0^r  \frac{h(2t) dt}{t} + h(r)\Big] 
\leq 10^6 r j(r)
\end{equation}
by \eqref{18.24}.

This was our basic estimate, but we can try to improve this in the same way as 
in Section~\ref{S15}, when we modified the tip of our second competitor to get a third one. 
There is a special case when things will be easier, which we want to mention first; this is when
(after the normalization that makes $r=1$)
\begin{equation}\label{18.38}
\rho^\ast = \varphi_\ast(K) \ \text{ for some } \varphi \in \Phi_X^+(\eta).
\end{equation}
This is the same statement as in \eqref{15.6}, but here $\eta$ comes from
the definition \eqref{18.5} and the statement of Theorem \ref{t18.1}. Apart from
this, we can still discuss as in Section \ref{S15}. 

Assume first that \eqref{18.38} holds, and let $Z(r)$ denote the cone over $\rho^\ast(r)$;
that is, 
\begin{equation} \label{18.39}
Z(r) = \big\{ \lambda \xi \, ; \, \xi \in \rho^\ast(r) \text{ and } \lambda\in [0,+\infty) \big\}
= \varphi_\ast(X)
\end{equation}
by definition of $\varphi_\ast(X)$, and  $Z(r) \in \cZ(X,\eta)$ by definition of $\cZ(X,\eta)$
(see above \eqref{18.3}). In addition, we can modify our first competitor near its tip, exactly as
we did below \eqref{15.6}, to construct an improved competitor and derive a better estimate
than \eqref{18.29} (or, if we renormalize back, \eqref{18.30}).
Indeed the competitor that we used so far coincides with $Z$ in a small ball
$\ol B(0,\kappa r)$, and we can further replace this tip with the intersection with $\ol B(0,\kappa r)$
of a competitor for $Z$ in $\ol B(0,\kappa r)$. In Section \ref{S15} we used the full length property
of $X$ to find this competitor; here we just use the definition \eqref{18.3} of $\alpha(Z)$,
which says that we can find a competitor $\wt Z$ for $Z$ in $\ol B(0,\kappa r)$, such that
\begin{equation} \label{18.40}
\H^2(\wt Z \cap \ol B(0,\kappa r)) \leq \H^2(Z \cap \ol B(0,\kappa r))
- \frac{\alpha(Z) \kappa^2 r^2}{2}.
\end{equation}
The construction of the new competitor for $E$, and in particular the gluing argument, is the 
same as in Section \ref{S15} (all the way up to \eqref{15.22}).
Thus we can save an extra $\frac{\alpha(Z) \kappa^2 r^2}{2}$ in the estimate \eqref{18.30},
and the proof of \eqref{18.37}, with this extra negative term, also yields
\begin{equation}\label{18.41}
\alpha(Z(r)) \leq C j(r),
\end{equation}
where the dependence of $C$ on $\kappa$ does not matter, because $\kappa$ is an absolute
constant.

We shall continue with the argument later, but let us now return to the case when 
\eqref{18.38} fails. As was explained below \eqref{15.22}, there may be a few different
reasons why this may happen. The first one is when Configuration H shows up in our construction.
In this case, we showed that, without using the full length condition, we can construct a 
new competitor (essentially obtained by contracting a hanging curve), and improve our estimate
\eqref{14.46} (the one that was used above to prove \eqref{18.29}) by an amount of $\eta(X)$; 
see \eqref{15.27}.
This means that we can subtract $\eta(X)$ from the right-hand side of \eqref{18.29},
or subtract $r^2 \eta(X)$ from the right-hand side of \eqref{18.30}. We claim that this is too much
to win if $j(r)$ is small enough. Indeed recall from \eqref{18.31} and the discussion below that
$\H^1(E\cap \S_r) - \H^1(\rho^\ast(r)) = \Delta_0(r) + \Delta_1(r) + \Delta_2(r) \geq 0$,
hence the improved \eqref{18.30} says that
\begin{eqnarray} \label{18.42}
\H^2(E \cap \ol B(0,r)) - \frac{r}{2} \H^1(E \cap \S_r) 
&\leq& - 10^{-5} r [\H^1(E\cap \S_r) - \H^1(\rho^\ast(r))] - r^2 \eta(X) +r^2 h(r)
\nn\\
&\leq&  - r^2 \eta(X) +r^2 h(r) 
\leq - r^2 \eta(X) + \varepsilon_1
\leq - \frac{1}{2} r^2 \eta(X)
\end{eqnarray}
by \eqref{18.8} and if $\varepsilon_1$ is small enough. Here again we feel good because
Remark \ref{R18.3} allows us to use a constant $\eta(X)$ that depends only on $n$, $\theta_0$
(through a covering of $\cX(\theta_0)$), and $\eta$ (because we forced $\eta(X) \leq \eta$).
At the same time 
\begin{equation}\label{18.43}
\H^2(E\cap \ol B(0,r)) \geq  r^2 \theta_0 \Big(1 -\alpha_n \int_0^r \frac{h(2t) dt}{t}\Big)
\geq r^2 \theta_0 \Big(1 -\alpha_n \varepsilon_1\Big)
\end{equation}
by \eqref{18.36} and \eqref{18.8}, and 
\begin{equation}\label{18.44}
\frac{r}{2} \H^1(E \cap \S_r) \leq r^2 j(r) + r^2 \theta_0
\end{equation}
by \eqref{18.28}; so \eqref{18.42} implies that
$\theta_0 \Big(1 -\alpha_n \varepsilon_1\Big) \leq \theta_0 + j(r) - \frac{1}{2} \eta(X)$,
which is impossible when $j(r) \leq j_0$, if $j_0$ and $\varepsilon_1$ are chosen small enough
(depending on a lower bound for $\eta(X)$, which itself depends on $\eta$). 
This proves that this first case when \eqref{18.38} fails does not happen when \eqref{18.27} holds.

The second reason why \eqref{18.38} may fail is explained below \eqref{15.27}; it corresponds to
the occurence of Configuration 3=2+1. In this case too, we constructed (without the help of the full length) 
a modification of our second competitor for $E$, that allows us to save $C^{-1} \eta(X)$ in the
estimate; see \eqref{15.30}. This case does not happen either, for the same reasons as for the previous case.

We are left with the case, described below \eqref{15.30}, where for some $\ell \in K \cap L$,
we added an element $c_\ell$ to $CC(\ell)$, to get the extended $CC_+(\ell)$, and we also added the 
point $\ell$ to $\rho^\ast$. If $\ell$ was already present in some $\rho_j$, $j\in J^\ast$,
we do not even need to worry; otherwise it is an isolated point of $\rho^\ast$ and we remove it.
That is, denote by $V'_0$ the set of (at most two) points $\ell$ that we added this way,
or equivalently the set of isolated points of $\rho^\ast$, and set $\rho' = \rho^\ast \sm V'_0$
(with this notation, we still normalize so that $r=1$). In this case we change a little the definition
of $Z(r)$, and set 
\begin{equation} \label{18.45}
Z(r) = \big\{ \lambda \xi \, ; \, \xi \in \rho' \text{ and } \lambda\in [0,+\infty) \big\}.
\end{equation}
Notice that when we have \eqref{18.38}, this new definition is the same as \eqref{18.39}, because
$\rho^\ast = \varphi_\ast(K)$ does not have isolated points.

We checked in Section \ref{S15} that $\rho'$ satisfies \eqref{15.6}, or in other words that
$\rho' = \varphi_\ast(K)$  for some $\varphi \in \Phi_X^+(\eta)$. 
Thus $Z(r) \in \cZ(X,\eta)$ (see the definition above \eqref{18.3}). 
Also, we can use any competitor for $Z(r)$ in the unit ball to modify the tip of our 
second candidate, essentially as in the case when \eqref{18.38} holds. The point is to extend
the deformation (originally defined on $Z(r)$) to the cone over $\rho^\ast$, get a competitor
for the cone over $\rho^\ast$, and glue it to the tip of our second competitor. The verifications
are done below \eqref{15.30}. This way we also get the additional estimate \eqref{18.41} in 
this remaining case.

We will still need to check \eqref{18.25} and \eqref{18.26} with this choice of $Z(r)$;
we will do this at the same time as we do it for the other case.
In the mean time, observe that there is yet another way to improve on our main estimate
\eqref{18.32}, which is to notice that the first inequality in \eqref{18.28} may be strict.
That is, set
\begin{equation}\label{18.46}
\Delta_3(r) = v'(r) - \H^1(E \cap \S_r) \geq 0.
\end{equation}
Notice that \eqref{18.30} implies that 
\begin{equation}\label{18.47}
\H^2(E \cap \ol B(0,r)) \leq \frac{r}{2} \H^1(E \cap \S_r) + r^2 h(r)
\end{equation}
by \eqref{18.31} and because $\Delta_i(r) \geq 0$ for $0 \leq i \leq 2$. Then
\begin{equation}\label{18.48}
\H^1(E \cap \S_r) =  v'(r) -  \Delta_3(r)
\leq 2r j(r) + 2 r \theta_0 -  \Delta_3(r)
\end{equation}
by \eqref{18.28}, so
\begin{equation}\label{18.49}
\Delta_3(r) \leq 2r j(r) + 2 r \theta_0 - \H^1(E \cap \S_r)
\leq 2r j(r) + 2 r \theta_0 - \frac{2}{r}\, \H^2(E \cap \ol B(0,r)) + 2r h(r).
\end{equation}
Since by \eqref{18.36}
\begin{equation}\label{18.50}
\frac{1}{r} \, \H^2(E\cap \ol B(0,r)) 
\geq r \theta_0 \Big(1 -\alpha_n \int_0^r \frac{h(2t) dt}{t}\Big),
\end{equation}
we see that
\begin{equation}\label{18.51}
\Delta_3(r) 
\leq 2r j(r) + 2 r \alpha_n \int_0^r \frac{h(2t) dt}{t} + 2r h(r)
\leq 4 r j(r)
\end{equation}
by \eqref{19.27n}. This completes the list of our basic estimates on the $\Delta_i(r)$.

\ms
We return to the proof of \eqref{18.25} and \eqref{18.26} for $Z(r)$.
Here we take for $\gamma^\ast(r)$ the set $r \gamma^\ast$, where $\gamma^\ast$
is the set of \eqref{13.25}, with the normalization by $r=1$. Thus
$\H^1(E \cap \S_r \sm \gamma^\ast(r)) = \Delta_0(r) \leq 10^6 r j(r)$ by
\eqref{18.32} and \eqref{18.37}, and \eqref{18.25} holds.

For \eqref{18.26} we first estimate 
\begin{equation}\label{18.52}
\Delta_4(r) 
= \H^1(\Gamma^\ast(r) \sm \gamma^\ast(r)))+\H^1(\gamma^\ast(r) \sm \Gamma^\ast(r))),
\end{equation}
where $\Gamma^\ast(r) = r \Gamma^\ast$ and $\Gamma^\ast$
is the net of Lipschitz graphs that shows up in \eqref{14.1}, for instance. We claim that
\begin{equation}\label{18.53}
\Delta_4(r) \leq C [\Delta_1(r) + \Delta_2(r)] \leq C r j(r).
\end{equation}
The last part comes from \eqref{18.37}. For the first part recall the decompositions
\begin{equation} \label{18.54}
\gamma^\ast = \Big(\bigcup_{i\in \cI_1} \cL_i \Big)\cup \Big(\bigcup_{c\in CC} \ol\gamma_c\Big),
\end{equation}
from \eqref{13.25} and 
\begin{equation} \label{18.55}
\Gamma^\ast 
= \Big(\bigcup_{i\in \cI_1} \Gamma_i \Big)\cup \Big(\bigcup_{c\in CC_+} \Gamma_c\Big),
\end{equation}
from \eqref{14.1}. The difference between $CC$ and $CC_+$ is just that maybe for some
$\ell \in K \cap L$ we added the trivial component $\{ c_\ell\}$ to $CC(\ell)$ to get
$CC_+(\ell)$; see above \eqref{12.10}. When $c\in CC_+\sm CC$,
we took $\Gamma_c = \{ \ell \}$, and this does not contribute
to the measure of the symmetric difference. Thus, returning to the normalization with $r=1$ and
using $\Delta$ to denote symmetric differences, 
\begin{equation} \label{18.56}
r^{-1}\Delta_4(r) \leq \sum_{i\in \cI_1} \H^1(\Delta(\cL_i,\Gamma_i))
+ \sum_{c\in CC} \H^1(\Delta(\ol\gamma_c,\Gamma_c)).
\end{equation}

For $i\in \cI_1$, we use the fact that $\Gamma_i = \cL_i$ when we dare to apply Remark \ref{r6.3n},
so we get no contribution, but even if we did not dare, $\Gamma_i$ would be obtained from $\cL_i$
by the construction of Section \ref{S7}, so \eqref{7.16} and \eqref{7.5} would yield
$\H^1(\Delta(\cL_i,\Gamma_i)) \leq C [\H^1(\cL_i)-\H^1(\rho_i)]$,
where $\rho_i$ is the geodesic with the same endpoints as $\cL_i$.
When we sum we would get a contribution which is dominated by $r^{-1}\Delta_2(r)$;
see \eqref{14.44}.

For $c\in CC$, we already observed in \eqref{14.37} that when we apply \eqref{9.7}
to each configuration $c\in CC$ and then sum, we get that
\begin{equation}\label{18.57}
\sum_{c\in CC} \H^1(\Delta(\ol\gamma_c,\Gamma_c)) 
\leq C(\lambda) \sum_{c\in CC} [\H^1(\overline\gamma_c)- \H^1(\Gamma_c)]
+ [\H^1(\Gamma_c)-\H^1(\rho_c)]
\end{equation}
which is then dominated by $r^{-1}\Delta_1(r) + r^{-1}\Delta_2(r)$; 
see the argument below \eqref{14.37}, and compare our definitions with \eqref{14.30} and \eqref{14.33}.
This completes our proof of \eqref{18.53}.

Next we use $\Delta_4(r)$ to control some distances. Set 
\begin{equation} \label{18.58}
\Gamma' 
= \Big(\bigcup_{i\in \cI_1} \Gamma_i \Big)\cup \Big(\bigcup_{c\in CC} \Gamma_c\Big),
\end{equation}
where the only difference with $\Gamma^\ast$ is that we removed $CC_+ \sm CC$ from the
indices, which means that we may have removed one or two points $\ell \in K \cap L$ 
from $\Gamma^\ast$. Let us first check that
\begin{equation} \label{18.59}
\dist(z,\Gamma') \leq r^{-1}\Delta_4(r) \leq C j(r)
\ \text{ for } z\in \gamma^\ast.
\end{equation}
First assume that $z\in \ol\gamma_c$ for some $c\in CC$; notice that 
$\H^1(\ol\gamma_c) \geq \eta(X)$ because $\ol\gamma_c$ meets 
$\d D = \S \cap \d B(\ell,\tau)$ for some $\ell \in K\cap L$ 
and it reaches out to some other endpoint $a_i^\ast$ (see \eqref{6.1}, \eqref{9.1},
\eqref{3.10} and \eqref{3.11}, and \eqref{5.38}). Then, since \eqref{18.27}
implies that $j(r) \leq j_0 < \eta(X)/C$ (if $j_0$ is chosen small enough; recall that
by Remark \ref{R18.3} we have a lower bound on $\eta(X)$ that depends on
$\eta$ and the other usual constants, but not on $X$), 
\eqref{18.53} implies that $\ol\gamma_c$ meets $\Gamma_c$,
and $\dist(z,\Gamma^\ast) \leq \dist(z,\ol\gamma_c \cap \Gamma^\ast)
\leq \H^1(\ol\gamma_c \sm \Gamma^\ast) \leq \Delta_4(r)$ because $\ol\gamma_c$
is connected. The case when $z\in \cL_i$ for some $i\in \cI_1$ is treated the same way,
because $\H^1(\cL_i) \geq \eta(X)$ too; \eqref{18.59} follows.
Conversely,
\begin{equation} \label{18.60}
\dist(z,\gamma^\ast) \leq r^{-1} \Delta_4(r) \leq C j(r)
\ \text{ for } z\in \Gamma',
\end{equation}
by the same argument as above, but this time using the fact that $\Gamma_c$
and $\Gamma_i$ are connected too. So we control the Hausdorff distance between
$\gamma^\ast$ and $\Gamma'$. We still need to compare $\Gamma'$ and 
$\rho' = r^{-1} (Z(r) \cap \S_r)$ (see \eqref{18.45}).

For each of the Lipschitz curves $\Gamma_j$ that compose $\Gamma^\ast$ 
(this time, with the condensed notation of \eqref{14.5}, but avoiding the trivial curves
$\{ \ell \}$ that come from $CC_+ \sm CC$), $\rho_j$ is the geodesic with 
the same endpoints, and by Pythagorus (and a tiny bit of spherical geometry, but recall that
the diameter of $\rho_j$ is less than $11/10$, say),
\begin{eqnarray} \label{18.61}
\sup_{z\in \rho_j} \dist(z,\Gamma_j) + \sup_{z\in \Gamma_j} \dist(z,\rho_j) 
&\leq& 10 [\length(\Gamma_i) - \length(\rho_j)]^{1/2} \length(\rho_j)^{1/2}
\nn\\&\leq& C j(r)^{1/2} \length(\rho_j)^{1/2} \leq C j(r)^{1/2}.
\end{eqnarray}
We take a supremum and get that
\begin{equation} \label{18.62}
\sup_{z\in \rho^\ast} \dist(z,\Gamma') + \sup_{z\in \Gamma'} \dist(z,\rho^\ast) 
\leq C j(r)^{1/2}.
\end{equation}
Now \eqref{18.26} follows from \eqref{18.62}, \eqref{18.59}, and \eqref{18.60}.
This completes our proof of Lemma~\ref{t18.4}.
\qed

\section{Where we control the variations of $Z(r)$}
\label{S19}

At this stage, we found for most $r\in \cR$ a nice cone $Z(r)$, which 
approximates $E$ well on $\S_r = \partial B(0,r)$.
The next stage is to show that $Z(r)$ varies slowly with $r$, and for this we start
with a study of some almost radial curves drawn on $E$, and that cross the annulus
\begin{equation}\label{19.1}
A = \ol B(0,90) \sm B(0,10^{-2}).
\end{equation}

Let $X \in \cX(\theta_0)$ be as in the last section, and recall from
the discussion over \eqref{18.19} that
\begin{equation}\label{19.2}
d_{0,180}(E,X) \leq 2\varepsilon_0 \leq 2 \cdot 10^{-3} \varepsilon(X).
\end{equation}
We shall use again the standard decomposition 
of $X$ into arcs $\cC_j$, $j\in J$, that is given by Section \ref{S3}. For each $j\in J$, 
denote by $\d\cC_j$ the boundary of $\cC_j$ (composed of its two endpoints), and let
\begin{equation}\label{19.3}
\cC'_j = \big\{ z\in \cC_j \, ; \, \dist(z, \d\cC_j) \geq 10^{-1} \eta(X)\big\}
\end{equation}
denote a slightly smaller arc where we remove a bit of $\cC_j$ at each end.

Denote by $P(j)$ the $2$-plane that contains $\cC_j$; we shall think of $P(j)$ as being horizontal.
For each $z\in \cC'_j$, denote by $P_z = P_{j,z}$ the vector hyperplane that contains
$z$ and is orthogonal to $\cC_j$ at $z$; we think of $P_z$ as the vertical hyperplane through $z$. 
Also denote by $L(z)$ the half line through $z$, and set 
\begin{equation} \label{19.4}
T(z) = \big\{ \xi \in A \, ; \, \dist(\xi,L_z) \leq 10^4\varepsilon_0 \big\}
\end{equation} % pas un cone, mais c'est OK
(a thin tube around $L(z)$), 
\begin{equation}\label{19.5}
T = \bigcup_{z\in \cC'_j} T(z)
\end{equation}
and 
\begin{equation}\label{19.6}
G_z = T(z) \cap P_z \cap E.
\end{equation}
\begin{lem} \label{t19.1}
For $j\in J$ and $z \in \cC_j'$, the set $G_z$ is a $C^1$ and $\frac{1}{10}$-Lipschitz graph 
over some segment of $L(z)$, and it crosses $A$.
\end{lem}

We shall even prove that $T(z) \cap E$ is a $C^1$ and $\frac{1}{20}$-Lipschitz graph,
over a piece of $P(j)$, and then the Lipschitz description of $G_z$ will follow from the Implicit Function Theorem.

We shall use the interior $C^1$ regularity theorem, that we pick from \cite{C1}. 
Set $\tau = 10^{-5}\eta(X)$, and let $w\in L(z) \cap A$ be given. By 
\eqref{19.2} $E$ is $360 \varepsilon_0$-close to $X$ in $B(w,10\tau)$, and 
we can pick $x_0 \in E$ such that $|x_0-w| \leq 360 \varepsilon_0$.

We want to apply Corollary 12.25 of \cite{C1} to $E$, in a small ball centered at $x_0$,
but there will be a few assumptions to check. We first take care of the distance to a plane.
Observe that
\begin{equation}\label{19.7}
d_{w,10\tau}(E,X) \leq \frac{18}{\tau} d_{0,180}(E,X) \leq \frac{36 \varepsilon_0}{\tau}  
\end{equation}
by \eqref{19.2}, and because we may safely assume that $\varepsilon_0$ is much smaller
than $\eta(X)$ and $\tau$ (so that $B(w,10\tau) \subset B(0,180)$).

Let us check that $X$ coincides with $P(j)$ near $w$. Recall from \eqref{3.9} that 
\begin{equation}\label{19.8}
\dist(z,K\sm \cC_j) \geq \min(\eta_0, \dist(z,\d\cC_j)) \geq 10^{-1} \eta(X)
\end{equation}
by \eqref{19.3} and \eqref{3.12}. Then $\dist(z, X \sm P(j)) \geq 10^{-1} \eta(X)/2$
(because $X$ is a cone and $P(j)$ contains the cone over $\cC_j$), and by homogeneity
\begin{equation}\label{19.9}
\dist(w, X \sm P(j)) \geq 10^{-3} \eta(X)/2 \geq 50 \tau.
\end{equation}
Conversely, $\dist(w, P(j) \sm X) \geq 50 \tau$ even more easily, because $\cC_j$ contains a 
$10^{-1}\eta(X)$-neighborhood of $\cC'_j$ in $P(j) \cap \S$ by \eqref{19.3}.
Thus \eqref{19.7} implies that 
\begin{equation}\label{20.10n}
d_{w,10\tau}(E,P(j)) \leq \frac{36 \varepsilon_0}{\tau}.
\end{equation}
Let $P'$ denotes the $2$-plane parallel to $P(j)$ and that contains $x_0$; notice that
$P'$ lies within $360 \varepsilon_0$ of $P(j)$, so we get that
\begin{equation}\label{19.10}
d_{x_0,9\tau}(E,P') \leq 100 \tau^{-1}\varepsilon_0,
\end{equation}
again with $\tau^{-1}\varepsilon_0$ as small as we want, and because 
$B(x_0,9\tau) \subset B(w,10\tau)$, with some extra space to take care of the difference
between $P(j)$ and $P'$.

This will take care of the distance assumption in Corollary 12.25 of \cite{C1}.
But we also have a density requirement, which will be fulfilled because if $E$ is so close to
$P'$ in $B(x_0,9\tau)$, then its density in $B(x_0,8\tau)$ cannot be too large.

More precisely, we want to apply Lemma~16.43 in \cite{Holder} to the ball $B(x_0,\rho)$,
with $\rho = 8\tau$, and  with a small constant $\delta$ that will be chosen soon.
For this there are only three things to check. First, that $E$ is almost minimal in $B(x_0,10\rho/9)$
(without a sliding condition). It is clear that $B(x_0,10\rho/9) \subset B(x_0,10\tau) \subset B(0,180)$, 
so we just need to check that $x_0$ is far from $L$. But
\begin{eqnarray}\label{19.11}
\dist(z, L \cap \S) &\geq& \min(\dist(z,L \cap \cC_j), \dist(z,K\sm \cC_j), \dist(z, L\cap \S \sm K))
\nn\\
&\geq& \min(\dist(z,\d\cC_j), \dist(z,K\sm \cC_j), \dist(K,L\cap \S \sm K))
\geq 10^{-1} \eta(X)
\end{eqnarray}
because the interior of $\cC_j$ does not meet $L$ (by \eqref{3.2}), and by 
\eqref{3.10}, \eqref{3.12}, and \eqref{19.8}. Then
$\dist(z,L) \geq \frac{2}{30}\eta(X)$, and 
\begin{equation}\label{19.12}
\dist(x_0,L) \geq \dist(w,L) - \frac{36 \varepsilon_0}{\tau} 
\geq  \frac{2}{30} 10^{-2} \eta(X) - \frac{36 \varepsilon_0}{\tau} 
\geq \frac{\eta(X)}{2000} \geq 50\tau,
\end{equation}
and $E$ is (plain) almost minimal even in $B(x_0,50\tau)$. 

Next $h(20\rho/9)$ should be small enough (again as in (16.44) of \cite{Holder}), 
but this follows from  \eqref{18.8} if $\varepsilon_1$ is small enough, depending also on $\delta$. 
Finally, $d_{x_0,10\rho/9}(E,P')$ should be small enough (depending on $\delta$), and this follows from 
\eqref{19.10} if $\varepsilon_0$ is chosen small enough (depending on both $\tau$ and $\delta$). 
So we may apply Lemma~16.43 in \cite{Holder}, and we get that
\begin{equation}\label{19.13}
\H^2(E \cap B(x_0,\rho)) \leq  \H^2(P' \cap B(x_0,(1+\delta)\rho)) + \delta \rho^2 
\leq (1+3\delta) \pi \rho^2. 
\end{equation}
Because of the near monotonicity of the density $t \to t^{-2}\H^2(E \cap B(x_0,t))$,
we easily deduce from this that the density of $E$ at $x_0$ is $\pi$ (because the density is always 
$\geq \pi$ at a point of $E$, and the next smallest density is $3\pi/2$, that corresponds to points 
of type $\bY$). So the analogue of the density excess for $E$ at $x_0$ is
\begin{equation}\label{19.14}
\wt f(\rho) = \rho^{-2}\H^2(E \cap B(x_0,\rho)) - \pi  \leq 3 \delta \pi
\end{equation}
by \eqref{19.13}. This is the first part of the requirement (12.26) of Corollary 12.25 in \cite{C1},
for the radius $r_0 = \rho/110 = 8\tau/110$. We just need to be sure that $\delta$ is small enough,
depending on the constant $\varepsilon_1$ from \cite{C1}.

The second requirement, about the size of $h$ (i.e., in the present case, of $C_0$),
follows from the assumptions of Theorem \ref{t18.1}. The final requirement is that 
$d_{x_0,100r_0}(E,P')$ be small enough, and follows from \eqref{19.10} if $\varepsilon_0$
is small enough. Then Corollary 12.25 in \cite{C1} says that $E$ is $C^{1+\beta}$-equivalent
to a plane in $B(x_0,r_0)$, with some additional precisions on the way it is equivalent, and
an exponent $\beta > 0$ that could be computed in terms of our various constants.

In addition to this, and as described at the beginning of Section \ref{S5}, $E \cap B(x_0,r_0)$
is also a Lipschitz graph with small constant (as small as we want, if the constants $C_0$, 
$\varepsilon_1$, and $\varepsilon_0$ are chosen small enough) over a subset of $P'$ that 
contains $P' \cap B(x_0,r_0/2)$.
See the discussion below \eqref{5.20}. Thus there is a neighborhood of $T$ where $E$ is a 
$C^1$, and small Lipschitz, graph over its projection on $P(j)$; recall that the width of $T$ is smaller
than the radius of the balls where we get a $C^1$ and flat description above, so that we neither
get a hole in the projection, or two layers (we skip some of the details here).

Now we can apply the implicit function theorem and find that $E \cap P_z \cap B(x_0,r_0/4)$
is a Lipschitz graph over a segment of $L(z)$ that contains $L(z) \cap B(w,r_0/8)$.
Recall also that $r_0 = 8\tau/110$ is much larger than the width $10^4\varepsilon_0$
of $T(z)$. Lemma 19.1 then follows from our local Lipschitz description of $E \cap P_z$ 
near $L(z) \cap A$.
\qed

\ms
We want to relate average flatness estimates for the graphs $G_z$ to the 
variations of the density excess  $f(r) = \theta(r)-\theta_0$.
The connection will be through the coarea theorem, the computation of a Jacobian, 
and the following angle $\alpha(x)$.

For almost every $x\in E \cap A$, $E$ has a tangent plane $T_E(x)$ at $x$
(the rectifiability gives an approximate tangent plane, which would be enough here, but the
local Ahlfors regularity, or more brutally the fact that $E$ is $C^1$ in a neighborhood of almost every
point of $E \sm L$, give a true tangent plane). For these $x$, we denote by $\alpha(x) \in [0,\pi/2]$
the (smallest) angle between the radial direction $(0,x)$ and a unit vector in $T_E(x)$.

% Ceci, de l'ancienne version, est faux!! : 
%Notice that when $x\in G_z$ for some $z\in \cC_j$, our $C^1$ description of $E$ near
%$x$ implies that $\alpha(x)$ is also the angle between $(0,x)$ and the tangent direction of $G_z$ at $x$. 
We want to show that $\alpha(x)$ is small on average, and this will mean something
about the average regularity of $G_z$

\begin{lem}\label{t19.2} 
There is  a constant $C \geq 0$, that depends only on $n$, such that 
\begin{equation}\label{19.15}
\int_{x\in E \cap A \,;\, \cos\alpha(x) > 0}  [1-\cos\alpha(x)] \, d\H^2(x) \leq C \cE,
% added the constraint otherwise hard to control. Done however later below \eqref{20.14}
\end{equation}
where we set $\cE = f(90) + \int_{0}^{180} h(r) \frac{dr}{r}$. 
\end{lem}

\ms
Let us apply the coarea formula (i.e., Theorem 3.2.22 in \cite{Federer})
to some nonnegative measurable function $g$, on the rectifiable set $E \cap B(0,90)$,
and with the level sets of the function $x \to |x|$; this yields
\begin{equation}\label{19.16}
\int_{E \cap B(0,90)} g(x) J(x) d\H^2(x) 
= \int_{0}^{90} \Big\{ \int_{E \cap \S_r} g(x) d\H^1(x) \Big\} dr,
\end{equation}
where $J(r)$ is the appropriate Jacobian. In the present context, a simple computation shows that
$J(x) = \cos\alpha(x)$. See (4.13) and (4.14) in\cite{C1}  (but this is not so hard to check anyway). 

Let us take $g(x) = (\cos\alpha(x))^{-1}$ when $\cos\alpha(x) > 0$, and $g(x)=0$ otherwise, 
but first restrict to $B(0,r') \sm B(0,r)$, with $0 < r < r' < 90$;
notice that $g$ is integrable on $E$ against $J(x) d\H^2(x)$, and by \eqref{19.16}
\begin{equation}\label{19.17}
v(r')-v(r) = \int_{E \cap B(0,r') \sm B(0,r)} d\H^2(x) 
\geq \int_{r}^{r'} \Big\{ \int_{E \cap \S_r} g(x) d\H^1(x) \Big\} dr
\end{equation}
since $g(x) J(x) = g(x) \cos\alpha(x) \leq 1$ everywhere. The measurability of the inside
integral is part of the coarea formula. Also, when we divide by $r'-r$ and 
let $r'$ tend to $r$ in the formula above, we get that
\begin{equation}\label{19.18}
v'(r) \geq \int_{E \cap \S_r} g(x) d\H^1(x)
\end{equation}
for almost every $r \in (0,90)$ (both sides exist almost everywhere, since 
both sides of \eqref{19.17} are monotone functions of $r$ and $r'$).
Next, for almost every $r\in \cR$,
\begin{equation}\label{19.19}
\int_{E \cap \S_r} [g(x)-1] d\H^1(x) \leq v'(r) - \H^1(E \cap \S_r) = \Delta_3(r)
\leq 4 r j(r)
\end{equation}
by \eqref{19.18}, \eqref{18.46}, and \eqref{18.51}.
We apply the coarea formula in the other direction and get that
\begin{equation}\label{19.20}
\int_{r\in (10^{-2},90)}\int_{E \cap \S_r} [g(x)-1] d\H^1(x) dr
= \int_{E \cap A} [g(x)-1] \cos\alpha(x) d\H^2(x),
\end{equation}
which is the left-hand side of \eqref{19.15} by definition of $g$ (because the set where
$\cos\alpha(x)=0$ does not contribute).
To complete the proof, we just need to show that 
\begin{equation}\label{19.21}
\int_{r\in (10^{-2},90)} r j(r) \leq C\cE, 
\ \text{ with }  \cE=  f(90) + \int_{0}^{180} h(r) \frac{dr}{r}.
\end{equation}
Recall from \eqref{18.24} that $j(r) = r\theta'(r) + f(r) + (1+2 \theta_0\alpha_n) h(2r) 
+ (1+\theta_0\alpha_n) \int_0^r \frac{h(2t) dt}{t}$. Since here $r \leq 90$, the last two terms
are clearly dominated by the second half of $\cE$.
For $f(r)$, we observe that for $0 < r \leq 90$,
\begin{equation}\label{19.22}
\begin{aligned}
f(r) &= \theta(r)-\theta_0 \leq \theta(90) \exp\Big(\alpha_n \int_0^{90} \frac{h(2t) dt}{t}\Big) - \theta_0
\cr&= f(90) + \theta(90) \Big[\exp\Big(\alpha_n \int_0^{90} \frac{h(2t) dt}{t}\Big)-1 \Big]
\end{aligned}
\end{equation}
by the almost monotonicity formula \eqref{18.10}. We multiply by $r \leq 90$, integrate, 
and get less than $C \cE$. 
We are left with $\theta'$. But
\begin{eqnarray}\label{19.23}
\int_0^{90} r^2 \theta'(r)dr 
\leq 90^2 \int_{r \in (0,90); \theta'(r) \geq 0} \theta'(r) dr
&\leq& 90^2 \Big[\int_0^{90} \theta'(r)dr - \int_{r \in (0,90); \theta'(r) < 0} \theta'(r) dr \Big]
\nn\\
&\leq& 90^2 \Big[\theta(90) - \theta_0 + \int_{r \in (0,90); \theta'(r) < 0} 
\alpha_n  h(2r) \frac{dr}{r} \Big]
\nn\\
&\leq& 90^2 \Big[ f(90) + \alpha_n \int_0^{90} h(2r) \frac{dr}{r} \Big]
\end{eqnarray}
by Lemma \ref{t16.1} and 
because we know from \eqref{18.10} that $\theta'(r) \geq - \alpha_n r^{-1} h(2r)$ 
almost everywhere. This proves \eqref{19.21} and Lemma \ref{t19.2}.
\qed

\ms
We shall now use Lemma \ref{t19.2} to control the variations of 
the cone $Z(r)$ from the previous section. 
Let $j\in J$ and $z\in \cC'_j$ be given, and let $G_z$ be as in \eqref{19.6}.
Lemma \ref{t19.1} says that for $r \in (10^{-2},90)$, there is a   
unique point of $G_z \cap \S_r$, which we denote by $w_{z}(r)$. Set 
$\xi_{z}(r) = w_{z}(r)/|w_{z}(r)|$. Then set
\begin{equation}\label{19.24}
\delta_j(z) = \sup\big\{ |\xi_z(r)-\xi_z(r')| \, ; \, 10^{-2} < r, r' < 90 \big\};
\end{equation}
let us check that
\begin{equation}\label{19.25}
\int_{z\in \cC_j'} \delta_j(z) d\H^1(z) \leq C \cE^{1/2}.
\end{equation}
First we fix $z \in \cC_j'$ and study the variations of $\xi_z(r)$. By Lemma \ref{t19.1},
$\xi_z$ is $C^1$. We want to show that 
\begin{equation}\label{19.26}
|\xi_{z}'(r)| \leq C \sin\alpha(w_z(r)),
\end{equation}
but let us try not to get confused by the various angles. 
Set $x= w_z(r)$, $\alpha = \alpha(w_z(r))$, $e = \xi_{z}(r) = x/|x|$,
denote by $T$ the direction of the tangent plane to $E$ at $x$, and let $v\in T$ be a unit vector
that minimizes the angle with $e$. Thus $\langle v, e \rangle = \cos\alpha$.

Recall from the proof of Lemma \ref{t19.1} that near $x$, $E$ is a Lipschitz graph over
$P(j)$ (the plane that contains $\cC_j$) with a constant as small as we want. This means that
$T$ is as close to $P(j)$ as we want. 
In particular $T$ is not contained in $P_z$, and we can find a unit vector $a \in T$, which is 
orthogonal to $P_z$. Notice that $x \in G_z \subset P_z$, so $e\in P_z$ and 
$\dist(v,P_z) \leq \dist(v, \R e) = \sin\alpha$. Since \eqref{19.26} is trivial when $\alpha \geq 10^{-1}$,
we may assume that $\alpha \leq 10^{-1}$ (we could also have proved this too); 
then, denoting by $\pi_z$ the orthogonal projection on $P_z$, 
\begin{equation}\label{20.28n}
|\langle v, a \rangle| = |\langle v-\pi_z(v), a \rangle| \leq \dist(v,P_z) \leq \sin\alpha \leq 10^{-1},
\end{equation}
the basis $(v,a)$ is nearly orthogonal, and the norm (in $T$) of the projection on the direction of $v$
parallel to $a$ is less than $2$.

Denote by $w$ a unit tangent vector to $G_z$ at $x$; of course $w \in T$, and we can 
write $w = \lambda v + \mu a$, with $|\lambda| \leq 2$. Recall that we are interested in the angle
between $w$ and the radial direction $e$. Denote by $\pi_\perp$ the orthogonal projection
on the direction orthogonal to $e$; then
\begin{equation}\label{20.29n}
|\pi_\perp(w)| \leq |\lambda| |\pi_\perp(v)| + |\mu||\pi_\perp(a)|
= |\lambda| |\pi_\perp(v)| \leq 2 |\pi_\perp(v)| = 2 \sin\alpha
\end{equation}
because $a$ is orthogonal to $P_z$, hence to $e$, and by definition of $\alpha$.

Now we compute $\xi_{z}'(r)$ brutally. Since the derivative of $|w_z(r)|$ is 
$\langle w_z(r), w'_z(r) \rangle |w_z(r)|^{-1}$, 
\begin{equation}\label{20.30n}
\xi_{z}'(r) = \frac{w'_z(r)}{| w_z(r)|} - \frac{w_z(r) \langle w_z(r), w'_z(r)\rangle}{| w_z(r)|^3}.
\end{equation}
That is,
\begin{equation}\label{20.31n}
| w_z(r)| \xi_{z}'(r) = w'_z(r)  - w_z(r) \langle w_z(r), w'_z(r)\rangle | w_z(r)|^{-2}.
\end{equation}
As could be expected, the total contribution of $w'_z(r)$ in the direction of $w_z(r)$ 
(or equivalently, with the notation above, of $e$) vanishes. Also, $| w_z(r)| = r$;  we are left with 
\begin{equation}\label{20.32n}
|r \xi_{z}'(r)| = |\pi_\perp(w'_z(r))| \leq 2\sin\alpha |w'_z(r)|  \leq 3 \sin\alpha 
\end{equation}
by \eqref{20.29n} and because $G_z$ is a small Lipschitz graph; \eqref{19.26} follows.

We integrate \eqref{19.26} on a subinterval of $(10^{-2},90)$ and find that 
\begin{equation}\label{19.27}
\delta_j(z) \leq \int_{10^{-2}}^{90} |\xi_{z}(r)'| dr
\leq C \int_{10^{-2}}^{90} \sin\alpha(w_z(r)) dr.
\end{equation}
Then we integrate on $\cC_j'$ and get 
\begin{equation} \label{19.28}
\int_{z\in \cC_j'} \delta_j(z) d\H^1(z)
\leq C \int_{z\in \cC_j'} \int_{10^{-2}}^{90} \sin\alpha(w_z(r)) dr d\H^1(z).
\end{equation}
Now the double integral looks like an integral on a piece of $E \cap A$.
Indeed, denote by $Gr(j)$ the union of the graphs $G_z$, $z\in \cC_j'$; 
that is, 
\begin{equation}\label{19.29}
Gr(j) = E \cap A \cap \bigcup_{z\in \cC_j'} (T(z) \cap P_z)
=  E \cap A \cap T \cap \bigcup_{z\in \cC_j'} P_z.
\end{equation}
By the proof of Lemma \ref{t19.1}, $Gr(j)$ is a $\frac{1}{20}$-Lipschitz graph over 
(a subset of) $P(j)$. In addition, \eqref{20.10n} says that it stays as close as we want to $P(j)$,
and therefore $\cos\alpha(w) > 0$ on $Gr(j)$.
Now \eqref{19.28} yields
\begin{eqnarray}\label{19.30}
\int_{z\in \cC_j'} \delta_j(z) d\H^1(z) 
&\leq& C \int_{w\in Gr(j)} \sin\alpha(w) d\H^2(w)
\leq C \Big\{\int_{Gr(j)} \sin^2\alpha(w) d\H^2(w) \Big\}^{1/2} 
\nn\\
&\leq& C \Big\{\int_{w \in E \cap A \, ; \, \cos\alpha(w) > 0} 
\big[1-\cos\alpha(w)\big] d\H^2(w) \Big\}^{1/2} 
\leq C \cE^{1/2}
\end{eqnarray}
by Cauchy-Schwarz, because 
$Gr(j) \subset \big\{ w \in E \cap A \, ; \, \cos\alpha(w) > 0\big\}$,
then $\sin^2\alpha(w) \leq 2(1-\cos\alpha(w))$, and by \eqref{19.15}.
This proves \eqref{19.25}.

\begin{lem}\label{t19.3} 
Let $Z(r)$, $r\in \cR$, denote the cone of Section \ref{S18}. Then 
\begin{equation}\label{19.31}
d_{0,1}(Z(r),Z(s)) \leq C j(r)^{1/2} + C j(s)^{1/2} + C \cE^{1/2}
\end{equation}
for $r, s \in \cR$.
\end{lem}

Here and below, $C$ is allowed to depend on constants like $\eta(X)$.
It will be enough to prove \eqref{19.31} when
\begin{equation}\label{19.32}
\text{$j(r)$ and $j(s)$ are small enough}
\end{equation}
(depending on $\eta(X)$ in particular), because otherwise it is trivial. This will allow us 
to avoid some unpleasant cases.

First we construct some points.
Fix $r, s \in \cR$, with \eqref{19.32}, and let an index $j\in J$ be given. 
By Chebyshev, we can find $z_1 = z_1(j)$ and $z_2=z_2(j) \in \cC'_j$, with the following properties:
\begin{equation}\label{19.33}
\dist(z_1,z_2) \geq C^{-1} \eta(X),
\end{equation}
(we just use the fact that $\length(\cC_i) \geq 10^{-1} \eta(X)$ here)
\begin{equation}\label{19.34}
\delta_j(z_1) + \delta_j(z_2) \leq C \cE^{1/2}
\end{equation}
(by \eqref{19.25}, and if we choose $C$ in \eqref{19.34} large enough), and, 
for $i=1,2$ and with the notation of Lemma \ref{t18.4},
\begin{equation}\label{19.35}
w_{z_i}(r) \in \gamma^\ast(r) \ \text{ and } \ w_{z_i}(s) \in \gamma^\ast(s).
\end{equation}
Let us check that \eqref{18.25} allows us to arrange this last condition as well.
Recall that $w_{z_i}(r)$ and $w_{z_i}(s)$ lie in $E \cap A \cap T$,
where $T$ is the thin region near the cone over $\cC_j'$ that was defined in \eqref{19.5}.
On this region, the projection which to a point $w$ associates the point $z\in \cC_j'$
such that $w\in P_z$ is $C$-Lipschitz, and now the exceptional set of $z\in \cC_j'$
for which \eqref{19.35} fails is contained in the projection of the union of the bad sets
for \eqref{18.25}; we assume that $j(r)$ is so small that we have a lot of choices left,
and use Chebyshev to get \eqref{19.33} and \eqref{19.34}.

\ms
Now we want to use these points to control $Z(r)$, so let us first remind the reader of how we chose
$Z(r)$ and at the same time introduce more notation.
Recall from the discussion below \eqref{18.41} that since we may assume that $j(r)$ is small enough, 
as in \eqref{19.32}, we may assume that \eqref{18.27} holds. Then there are only two options. 
The first one is when \eqref{18.38} holds, i.e., when $\rho^\ast(r) = \varphi_\ast(K)$
is an acceptable small deformation of $K = X \cap \S$, and then we took 
$Z(r) = \varphi_\ast(X)$ (the cone over $\rho^\ast(r)$, or equivalently the corresponding
deformation of $X$), as in \eqref{18.39}. The other option is described below \eqref{18.44}), 
where $\rho^\ast(r)$ has one or two isolated additional points (vertices of $K \cap \S$), 
which we remove from $\rho^\ast$ to get $\rho'$, and then we take for $Z(r)$ the cone over $\rho'$, 
as in \eqref{18.45}.
Let us set $\rho' = \rho^\ast(r)$ in the first case (when \eqref{18.38} holds), so that 
$Z(r)$ is the cone over $\rho' = \rho'(r)$ in both cases. Of course it will be enough to
control $\rho'$.

By construction, $\rho'$ is composed of a collection of geodesics. 
Most of them are obtained from an arc $\cC_j$, $j\in J$, by moving a tiny bit 
one or two of its endpoints. Let us write $\rho'_j$ the
arc of geodesic that comes like this from $\cC_j$. When $\ell \in K \cap \S$ is one of the two
endpoints of $\cC_j$, it may be that the corresponding endpoint of $\rho'_j$ is of the form
$\varphi(\ell)$, or just $\ell$ itself, depending on the configurations. Also, it can happen that 
in addition to the $\rho'_j$, $\rho'$ contains one or two very short additional arcs, 
that go from some $\ell \in K \cap \S$ to $\varphi(\ell)$. 
For each $\ell$, there is only (at most) one such arc, which we call
$\rho_\ell = \rho(\ell, \varphi(\ell))$. Since we do this construction with both radii $r$ and $s$,
we shall often add this in the notation. Thus $Z(r)$ is the cone over
\begin{equation}\label{19.36}
\rho'(r) = \bigcup_{j\in J} \rho'_j(r) \cup \bigcup_{\ell} \rho_\ell(r),
\end{equation}
where the last union may be empty and concerns at most two vertices $\ell \in K\cap L$,
and there is a similar description for $\rho'(s)$.

We want to place each $\rho'_j(r)$ by finding two points in it; we will take care of 
the $\rho_\ell(r)$ later. 
We start from the two points $w_{z_i}(r)$, $i=1,2$. 
Since $w_{z_i}(r) \in \gamma^\ast(r)$,
\eqref{18.26} tells us that we can find $y_i(r) = y_i(r, j) \in Z(r) \cap \S_r$
such that $\dist(y_i(r,j), w_{z_i}(r)) \leq C j(r)^{1/2}$. 
We claim that
\begin{equation}\label{19.37}
y_i(r,j) \in \rho'_j(r) \text{ for } i=1,2 \  \text{ and } 
\dist(y_1(r,j),y_2(r,j)) \geq C^{-1}\eta(X).
\end{equation}
Recall that $w_{z_i}(r) \in E \cap T \cap \S_r$, so $w_{z_i}(r)$  
lies within $10^4\varepsilon_0 << \eta(X)$ of the cone over $\cC_j'$,
and (since $j(r)$ is small by \eqref{19.32}) $y_i(r,j)$ lies very close too.
Since $\cC_j'$ lies at distance at least $\eta(X)/10$ from $L \cap \S$ and $K \sm \cC_j$,
and at the same time $\varphi$ does not move points much, we see that $y_i(r,j)$
cannot lie in any other $\rho'_k(r)$, $k\in J \sm \{ j \}$, nor any $\rho'_\ell(r)$.
That is, $y_i(r,j) \in \rho'_j(r)$. 
In addition, the two $z_i$ are far from each other (by \eqref{19.33}), hence also the $w_{z_i}(r)$ 
and the $y_i(r,j)$. This proves \eqref{19.37}.

Let $\wh\rho_j(r)$ denote the great circle in $\S_r$ that contains the geodesic $\rho'_j$.
Then $r^{-1}\wh\rho_j(r)$ is the great circle in $\S$ that contains the two points
$r^{-1} y_i(r,j)$. 
Similarly define the great circles $\wh\rho_j(s)$
(starting from $Z(s)$), and points $y_i(s,j) \in Z(s) \cap \S_{s}$, and notice that
$s^{-1}\wh\rho_j(s)$ is the great circle in $\S$ that contains the two points
$s^{-1} y_i(s,j)$. In addition, for $i=1,2$,
\begin{eqnarray}\label{19.38}
&\,&\hskip-0.4cm |s^{-1}y_i(s,j)-r^{-1} y_i(r,j)| 
\nn\\
&\,& \hskip0.7cm\leq s^{-1} |y_i(s,j)- w_{z_i}(s)|
+ |s^{-1} w_{z_i}(s) - r^{-1} w_{z_i}(r)|
+r^{-1} |w_{z_i}(r)- y_i(r,j)|
\nn\\
&\,& \hskip0.7cm \leq C j(r)^{1/2} + |s^{-1} w_{z_i}(s) - r^{-1} w_{z_i}(r)| + C j(s)^{1/2}
\leq  C j(r)^{1/2} + C\cE^{1/2} + C j(s)^{1/2}.
\end{eqnarray}
We deduce from this and \eqref{19.37} that 
\begin{equation}\label{19.39}
d_{\H}(r^{-1}\wh\rho_j(r),s^{-1}\wh\rho_j(s)) \leq C j(r)^{1/2} + C j(s)^{1/2} + C\cE^{1/2} = C\cE',
\end{equation}
where the Hausdorff distance $d_\H$ is defined as in \eqref{17.5}, and we set
$\cE' = j(r)^{1/2} + j(s)^{1/2} + \cE^{1/2}$ to save some space. 
This is good, but we also want to control the position of the endpoints of the $\rho'_j(r)$ and the
$\rho'_j(s)$, because we want to show that
\begin{equation}\label{19.40}
d_{\H}(r^{-1}\rho'(r),s^{-1}\rho'(s)) \leq C \cE'.
\end{equation}
Indeed, \eqref{19.31} will follow from \eqref{19.40}, since $Z(r)$ is the cone over $r^{-1}\rho'(r)$,
and similarly for $Z(s)$.

We intend to prove this locally, in balls of radius roughly equal to $C^{-1}\eta(X)$
and centered on $K$. We start away from $K \cap L$, and first consider balls centered on the 
vertices of $V_1 \cup V_2$ of our standard decomposition (see the definitions near \eqref{3.5}).

Let $a_0 \in V_1 \cup V_2$ be given, call $\cC_j$, $\cC_k$, and maybe $\cC_l$
(that is, if $a_0 \in V_1$) the two or three arcs of $K$ that end at $a_0$.
By the various definitions, $\rho'_j(r)$, $\rho'_k(r)$, and maybe $\rho'_l(r)$ are arcs of geodesics 
that end at some point $r a(r)$, with $a(r) \in \S$, and $a(r)$ lies very close to $a_0$ (because it is of the
form $\varphi(a_0)$ for some $\varphi \in \Phi_X^{+}(\eta)$, with $\eta$ much smaller than $\eta(X)$).
We have a similar description of $\rho'_j(s)$, $\rho'_k(s)$, and maybe $\rho'_l(s)$, with another point
$a(s) \in \S$.

When $a_0 \in V_1$, the three $\cC_j$, $\cC_k$, and $\cC_l$ make $120^\circ$ angles
with each other, and the position of $a(r)$ is determined, within $10C \cE'$, as soon as we know 
the position of the full circle $r^{-1}\wh\rho_j(r)$ and its analogues for $k$ and $l$. 
The same thing holds for the radius $s$, and now \eqref{19.39} implies that
\begin{equation}\label{19.41}
|a(r)-a(s)| \leq C \cE'.
\end{equation}
Once we have this, and by \eqref{19.39} again, we easily deduce that
\begin{equation}\label{19.42}
d_{a_0, 10^{-4} \eta(X)}(r^{-1}\rho'_j(r), s^{-1}\rho'_j(s)) \leq C \cE'
\end{equation}
and similarly for $k$ and $l$. Since we are far from the $\rho_\ell$ and by \eqref{19.36},
we immediately get that
\begin{equation}\label{19.43}
d_{a_0, 10^{-4} \eta(X)}(r^{-1}\rho'(r), s^{-1}\rho'(s)) \leq C \cE'.
\end{equation}
This is good enough for \eqref{19.40}, so we may switch to the case when $a_0 \in V_2$,
and we started from two arcs $\cC_j$ and $\cC_k$ that go in opposite directions.
In this case, we will not control the geodesics separately, but we will be able to control the union.
That is, we may not know so precisely where $a(r)$ and $a(s)$ lie (i.e., \eqref{19.41} may fail),
but nonetheless we claim that \eqref{19.43} still holds, although maybe with a larger constant.
Indeed if the angle of $\wh\rho_j(r)$ and $\wh\rho_k(r)$ at $ra(r)$ is at most $C \cE'$,
the angle of $\wh\rho_j(s)$ and $\wh\rho_k(s)$ at $sa(s)$ is less than $C \cE'$ too,
and in the ball $B(a_0, 10^{-4}\eta(X))$, $r^{-1}\rho'(r)$ is $C\cE'$-close to $r^{-1}\wh\rho_j(r)$
(or to $r^{-1}\wh\rho_k(r)$, since the two are close to each other). 
If now the angle of $\rho_j(r)$ and $\rho_k(r)$ at $ra(r)$ is roughly $\lambda \cE'$,
with $\lambda$ large, then the proof of \eqref{19.41} merely gives
$|a(r)-a(s)| \leq C \lambda^{-1}$, but we still get \eqref{19.43} because the distance
between $\rho_j(r)$ and $\rho_k(r)$ (or similarly $\rho_j(s)$ and $\rho_k(s)$)
varies by at most $C\lambda \cE'$ times this distance. 
Said differently, we look for a Lipschitz graph (for instance $s^{-1}\rho'(s)$)
composed of two arcs of geodesics, knowing these two geodesics with errors of $C\cE'$;
then we can recover the graph within $C\cE'$ (after deciding which way it branches), which we can if
the geodesics make angles $\lambda \cE'$, with $\lambda$ large.

This takes care of small balls $B(a_0, 10^{-4}\eta(X))$ centered on $V_1$ and $V_2$. 
It is even easier to show that
\begin{equation}\label{19.44}
d_{a_0, 10^{-6} \eta(X)}(r^{-1}\rho'(r), s^{-1}\rho'(s)) \leq C \cE'
\end{equation}
when $a_0 \in K$ is such that $\dist(a_0, V_0 \cup V_1 \cup V_2) \geq 10^{-5} \eta(X)$,
because in the ball $B(a_0, 10^{-6} \eta(X))$, $r^{-1}\rho'(r)$ coincides with a single $r^{-1}\wh\rho_j(r)$,
the one for which $a_0 \in \cC_j$.
This comes from the fact that all the other $\rho'_j(r)$ (or $\rho'_\ell(r)$) are far away, by 
\eqref{3.9}-\eqref{3.12} and the fact that we have a good control on the angles that
two arcs $\rho'_j(r)$ make when they have a common endpoint.

This takes care of the part of $\rho'(r)$ and $\rho'(s)$ that lives far from $V_0 = K \cap L$,
and \eqref{19.40} (and hence also the lemma) will follow if we can prove that for $\ell \in V_0$,
\begin{equation}\label{19.45}
d_{\ell, 10^{-4} \eta(X)}(r^{-1}\rho'(r), s^{-1}\rho'(s)) \leq C \cE'.
\end{equation} 
We will need to distinguish cases, depending on the configurations 
that we encounter for $r$ and $s$. A priori, these two configurations may be different.

We start with the case when $K$ has only one branch near $\ell$. 
Since hanging curves never occur when $j(r)$ and $j(s)$ are small (recall \eqref{19.32}), 
there is only one curve $\rho'_j$ near $r\ell$, and this curve ends at $\ell$. The same thing
happens for $s$, and in this case \eqref{19.45} is a simple consequence of \eqref{19.39},
because we know where the curves stop (and on which side they are).

Next assume that $K$ has two branches at $\ell$. Call the corresponding indices $j$ and $k$.
Then (again because there is no hanging curve) we can only be in Configuration 2-
(treated below \eqref{9.19aa}) or Configuration 2+ (treated in Section \ref{S11}).

In the first case, $\rho'(r)$ is composed, near $r \ell$, of the two arcs of geodesic
$\rho'_j(r)$ and $\rho'_k(r)$, and nothing else. They have a common endpoint $r a(r)$,
and even though the position of $\wh\rho_j(r)$ and $\wh\rho_k(r)$ does not necessarily determine $a(r)$
with great precision (because $\rho'_j(r)$ and $\rho'_k(r)$ may make an angle at $ra(r)$ that is 
close to $\pi$), it still determines the union of $\rho'_j(r)$ and $\rho'_k(r)$ with a good enough precision.
That is, if both $r$ and $s$ are subject to Configuration 2-, then we have \eqref{19.45}, by the
same proof as for \eqref{19.43} when $a_0 \in V_2$.

When we have Configuration 2+ for $r$, there are again two cases. We start with the second one
(Case B), because then $\Gamma$ is composed of just two curves that start from $\ell$ (see near
\eqref{11.8}), the geodesics $\rho'_j(r)$ and $\rho'_k(r)$ both start from $\ell$, and their position near
$\ell$ is easy to deduce from the position of $\wh\rho_j(r)$ and $\wh\rho_k(r)$. 
If this happens both for $r$ and $s$, we get \eqref{19.45} right away, and even if 
we have this configuration for $r$ and Configuration 2- for $s$, or the other way around,
we still get \eqref{19.45} for the same reason as in Configuration 2-.

We are left with the case when at least one of the radii, say, $r$, belongs to Case A of Configuration 2+. 
In this case $\rho'(r)$ is composed of three geodesics near $\ell$, 
the usual $\rho'_j(r)$ and $\rho'_k(r)$, that end at a common point 
$ra(r)$, plus the short geodesic $\rho'_\ell(r)$ that goes from $r \ell$ to $ra(r)$.
In this case these three geodesics make large angles at $ra(r)$ (see \eqref{11.2}). 
In fact the proof of \eqref{11.4} (even simplified) shows that then
$\rho'_j(r)$, $\rho'_k(r)$, and $\rho'_\ell(r)$ make angles larger than 
$\frac{2\pi}{3} - \frac{\pi}{9} = \frac{5 \pi}{9}$ at $ra(r)$, and then 
$\rho'_j(r)$ and $\rho'_k(r)$ make an angle smaller than 
$2\pi - 2 \cdot \frac{5 \pi}{9} = \frac{8 \pi}{9} < \pi$ at $ra(r)$. 
In this case, we can recover the position of $a(r)$, within the usual error of $C\cE'$,
from the approximate position of the geodesics $\wh\rho_j(r)$ and $\wh\rho_k(r)$
(known within $C\cE'$).
In addition, in this case the same proof also shows that $\cC_j$ and $\cC_k$
make an angle smaller than $\frac{8 \pi}{9}$ at $\ell$, and we can recover the point
of intersection $sa(s)$ of $\rho'_j(s)$ and $\rho'_k(s)$ with the same sort of precision.
Thus, if $s$ is also coming from case A, we get \eqref{19.45} with the initial proof of \eqref{19.43}.

We are left with the case when $r$ is associated to Case A and $s$ is associated to
Case B or Configuration 2-. Case B is not a problem, because $a(s)$, which is the intersection
near $\ell$ of $s^{-1}\wt\rho_j(s)$ and $s^{-1}\wt\rho_k(s)$, lies very close to $a(r)$ 
(which has a similar definition in terms of $r$), and at the same time is equal 
to $\ell$, so that the additional geodesic $\rho'_\ell(r)$ is very short and we still get \eqref{19.45}.
We are left with the case when $s$ belongs to Configuration 2-. 
But in the present case $\cC_j$ and $\cC_k$ make an angle smaller than $\frac{8 \pi}{9}$ at $\ell$, 
and it is easy to see that our union of curves $\Gamma = \Gamma_1 \cup \Gamma_2$ is not efficient
because we may as well cut its edge near $\ell$. 
We claim that this case (i.e., Configuration 2- with an angle smaller than $\frac{8 \pi}{9}$) 
does not occur for $s$ when $j(s)$ is small enough.
The proof is the same as for Configuration 3 = 2+1, treated below \eqref{18.44}, except that
we don't even need to worry about the extra arc leaving from $\ell$.
This completes our proof of \eqref{19.45} when there are only two arcs $\cC_j$ and $\cC_k$
that leave from $\ell$.

Now may now assume that we have three arcs $\cC_i$, $\cC_j$, and $\cC_k$ that touch $\ell$. 
The three main geodesics $\rho'_j(r)$, $\rho'_j(r)$ and $\rho'_k(r)$ make angles
nearly equal to $\frac{2\pi}{3}$ near $\ell$, so the location of the intersections of the great
circles that contain them is known with good precision. 
In terms of Configurations, recall that there is no hanging curve, and that 
Configuration 3 = 2+1 is also ruled out by the discussion near \eqref{18.44}. We are thus left
with Configuration~3- (where $\rho'(r)$ is composed of the three geodesics $\rho'_i(r)$,
$\rho'_j(r)$, $\rho'_k(r)$ that all leave from a same endpoint that we call $ra(r)$ (see near \eqref{9.31}), 
and Configuration 3+, where again we have have two subcases. In Case A, $\Gamma$
and then $r^{-1}\rho'(r)$ are three-legged spider centered at $\ell$ (see \eqref{10.30}, the comment
that follows it, and then the discussion above \eqref{12.13} that confirms how we cut $\Gamma$
and found geodesics).

In Case B, $\Gamma$ and then and $r^{-1}\rho'(r)$ are authorized to have a fork. That is, 
they are composed of one curve that leaves from $\ell$ (with the notation of \eqref{10.90},
the corresponding piece of $\Gamma$ is called $\wt\Gamma_5$), a short arc $r^{-1}\rho'_\ell(r)$
(corresponding to $\wt\Gamma_4$ in \eqref{10.90}), that goes from $\ell$ to some 
fork point $a(r)$ (corresponding to $x_0$ in Section \ref{S11}), and then two other curves
that leave from $a(r)$ (corresponding to $\Gamma_2$ and $\Gamma_3$ in \eqref{10.90}).

If both $r$ and $s$ both correspond to Configuration 3- or Case A, then we have \eqref{19.45}
because the positions of $a(r)$ and $a(s)$ can be obtained with the desired
precision form the position of the great circles where they cut. The proof is still the same
as for \eqref{19.43}.

So we may assume that for $r$ we have Case B, and (again without loss of generality) that $a(r)$ 
is the common endpoint of $r^{-1}\rho'_j(r)$ and 
$r^{-1}\rho'_k(r)$. First assume that $s$ corresponds to Configuration 3- or Case A. Then
$a(s)$, which is the intersection near $\ell$ of $s^{-1}\wh\rho_j(s)$ and $s^{-1}\wh\rho_k(s)$, 
lies within $C \cE'$ of $a(r)$, which is defined similarly, but with $s$ replaced by $r$
(apply \eqref{19.39} as usual). Also, $\ell$ lies close to $s^{-1}\wh\rho_i(s)$ 
because it lies in $r^{-1}\wh\rho_i(r)$. Moreover, if we assume for the sake of the discussion
that the tangent of $\cC_i$ is horizontal at $\ell$ and leaves from $\ell$ in the direction of the right, 
$a(r)$ is roughly aligned  with the opposite of $r^{-1}\rho'_i(r)$ (see Lemma \ref{t10.5}), i.e.,
lies on the left of $\ell$, and then $a(s)$ also lies on the left (or at least, not far right) of $\ell$;
hence $\ell$ also lies within $C \cE'$ of $s^{-1}\rho'_i(s)$ (and not just $s^{-1}\wh\rho_i(s)$
as we said above). So $r^{-1}(\rho'_i(r) \cup \rho'_\ell(r))$ is $C \cE'$-close to $s^{-1}\rho'_i(s)$ 
and we get \eqref{19.45} by adding the two other geodesics.

We may thus assume that $s$ also corresponds to Case $B$. If $\rho'_i(s)$ is also
the geodesic of $\rho'(s)$ that leaves from $s\ell$, the intersection $a(s)$ of 
$s^{-1}\rho'_j(s)$ and $s^{-1}\rho'_k(s)$ lies close to $a(r)$, as before, and \eqref{19.45} 
holds as usual. So we may assume that $\rho'_j(s)$, say, is the one that starts from $s\ell$, and
$a(s)$ is the common endpoint of $s^{-1}\rho'_i(s)$ and $s^{-1}\rho'_k(s)$.
This is not impossible, but we shall show that then $a(r)$ and $a(s)$ are both close to $\ell$.

By Lemma \ref{t10.5}, $a(r)-\ell$ lies in the direction almost opposite to
the direction of $\cC_i$ at $\ell$; since $ra(r) \in \wh\rho_j(r)$ and $\wh\rho_j(r)$
runs in a quite different direction, this proves that
$\dist(\ell, r^{-1}\wh\rho_j(r)) \geq \frac{1}{10}|a(r)-\ell |$.
On the other hand, $s\ell \in \rho'_j(s) \subset \wh\rho_j(s)$, so
$\ell$ lies $C \cE'$-close to $r^{-1}\wh\rho_j(r)$ (by \eqref{19.37}) 
and altogether $a(r)$ lies $C \cE'$-close to $\ell$.
The same argument (with $r$ and $s$ exchanged) shows that 
$|a(s)-\ell| \leq C\cE'$; then \eqref{19.45} follows as usual: we control the directions
of the geodesics and their origin.

So \eqref{19.45} holds in our last case. We have seen earlier that
\eqref{19.40}, and then \eqref{19.31}, follow. 
This completes the proof of Lemma \ref{t19.3}.
\qed

\section{We finally get a good approximation by cones}
\label{S20} 

In this section we complete the proof of Theorem \ref{t18.1}.
In the previous sections, we took $E$ as in that theorem, selected a 
sliding minimal cone $X$ (see Lemma \ref{t18.2}), 
constructed deformations $Z(r)\in \cZ(X,\eta)$ of $X$, $r\in \cR$
(in Lemma \ref{t18.4}), and proved that they often lie close to each other
(see Lemma \ref{t19.3}).
Now we want to pick one of the cones $Z(r)$ and show that it is close to $E$,
as needed for Theorem \ref{t18.1}. 

So let us choose a radius $r_0 \in \cR$. We simply use Chebyshev to select 
$r_0\in \cR$ such that
\begin{equation}\label{20.1a}
r_0\in (1,2) \ \text{ and }\ 
j(r_0) \leq 2 \int_{1}^2 j(r) dr \leq C \cE,
\end{equation}
where the second inequality comes from \eqref{19.21}. 

Set $Z = Z(r_0)$; want to show that $E$ is close to $Z$ in, say, $B(0,2)$,
but it will be simpler to first take care of the annulus $A_0 = B(0,2) \sm B(0,10^{-1})$;
we will worry later about $B(0,10^{-1})$, with an iteration argument.
First we check that points of $Z \cap A_0$ are close to $E$.

\begin{lem} \label{t20.1}
With $Z$ and $A_0= B(0,2) \sm B(0,10^{-1})$ as above,
\begin{equation} \label{20.3}
\dist(z,E) \leq C \cE^{1/3} \ \text{ for } z\in Z\cap A_0.
\end{equation}
\end{lem}
\ms
Let $z\in Z\cap A_0$ be given, set $r = |z|$, and pick $s \in \cR$ such that $j(s) \leq \cE^{2/3}$.
By Chebyshev, we can find $s$ so that
\begin{equation} \label{20.4}
|s-r| \leq 2 \cE^{-2/3} \int_{0}^{90} j(r) dr  
\leq C \cE^{1/3},
\end{equation}
by \eqref{19.21} again. Set $z_1 = s r^{-1} z$; thus $z_1\in Z \cap \S_s$ and 
$|z_1-z| = |s-r| \leq C \cE^{1/3}$. By Lemma~\ref{t19.3} 
(applied to $r_0$ and $s$), we can find $z_2 \in Z(s)\cap \S_s$ such that
\begin{equation} \label{20.5}
|z_2-z_1| = \dist(z_1,Z(s)) \leq 3d_{0,1}(Z,Z(s)) 
\leq C (j(r_0)+j(s)+\cE)^{1/2}
\leq C \cE^{1/3}.
\end{equation}
Then we use Lemma \ref{t18.4} and find $x\in \gamma^\ast(s)$ such that 
$|x-z_2| \leq C j(s)^{1/2} \leq C \cE^{1/3}$. 
Since $x \in \gamma^\ast(s) \subset E$ (see the first line of Lemma~\ref{t18.4}),
we get that $\dist(z,E) \leq |x-z| \leq C \cE^{1/3}$, as needed.
\qed

\begin{lem} \label{t20.2}
Keep $Z$ and $A_0$ as above; then
\begin{equation} \label{20.6}
\dist(x,Z) \leq C \cE^{1/4} \ \text{ for } x\in E\cap A_0.
\end{equation}
\end{lem}

\ms
We first find some points of $E$ for which \eqref{20.6} holds.
Let $x\in E$ be given, and first assume that $r = |x|$ lies in $\cR$,
with $j(r) \leq \cE^{1/2}$, and that in addition $x \in \gamma^\ast(r)$.
Then by \eqref{18.26}, there is a point $z\in Z(r)$ such that
$\dist(x,z) \leq C j(r)^{1/2}\leq C \cE^{1/4}$. In addition,
Lemma~\ref{t19.3} gives us a point $w\in Z = Z(r_0)$ such that
$|w-z| \leq C [j(r)+j(r_0) + \cE]^{1/2} \leq C \cE^{1/4}$,
by \eqref{19.31} and the definition of $r_0$. That is,
$\dist(x,Z) \leq C_0 \cE^{1/4}$ for some constant $C_0$ that
satisfies the usual requirements. Next we consider
\begin{equation} \label{20.7}
E_0 = \big\{ x\in E \cap B(0,3) \sm B(0,10^{-2}) 
\, ; \,  \dist(x,Z) > C_0 \cE^{1/4} \big\}.
\end{equation}
We want to estimate the measure of $E_0$, and unfortunately we will have 
to single out the ugly set
\begin{equation}\label{21.8n}
E_b = \big\{ x\in E \cap B(0,3) \sm B(0,10^{-2}) 
\, ; \, E \text{ has no tangent plane at $x$ or } 
\cos\alpha(x) = 0 \big\},
\end{equation}
which will be treated separately, after we look at
 \begin{equation} \label{20.8}
E_1 = \big\{ x\in E_0 \sm E_b \, ; \,  |x| \notin \cR \text{ or } j(|x|) > \cE^{1/2} \big\}
\end{equation}
and
\begin{equation} \label{20.9}
E_2 = \big\{ x\in E_0 \sm (E_b \cup E_1) 
\, ; \,  x \in E \cap \S_{|x|} \sm \gamma^\ast(|x|)  \big\}.
\end{equation}
By the discussion above, $E_0 = E_b \cup E_1 \cup E_2$.

We shall now use the coarea formula and Lemma \ref{t19.2} to estimate $\H^2(E_1 \cup E_2)$. 
We write
\begin{eqnarray} \label{20.10}
\H^2(E_0 \sm E_b) &=& \int_{E_0 \sm E_b} \Big\{[1-\cos\alpha(x)]+\cos\alpha(x) \Big\} \, d\H^2(x) 
\leq C\cE + \int_{E_0\sm E_b} \cos\alpha(x) d\H^2(x) 
\nn\\
&=& C\cE + \int_{r=10^{-2}}^3 \H^1((E_0 \sm E_b)\cap \S_r) dr
\nn\\
&\leq& C\cE + \int_{r=10^{-2}}^3 \H^1(E_1 \cap \S_r) dr + \int_{r=10^{-2}}^3 \H^1(E_2 \cap \S_r) dr
\end{eqnarray}
by \eqref{19.15} and \eqref{19.16} with $g = \1_{E_0 \sm E_b}\,$, and where we recall that 
$J(x) = \cos\alpha(x)$. Next 
\begin{equation} \label{20.11}
\int_{r=10^{-2}}^3 \H^1(E_2 \cap \S_r) dr
\leq \int_{r=10^{-2}}^3 \H^1(E \cap \S_r \sm \gamma^\ast(r)) dr
\leq \int_{r=10^{-2}}^3 j(r) dr \leq C \cE
\end{equation}
by \eqref{20.9}, \eqref{18.25} and \eqref{19.21}. 
Now we estimate $\int_{r=10^{-2}}^3 \H^1(E_1 \cap \S_r) dr$. We can drop 
the radii $r \in (10^{-2},3) \sm \cR$, because the corresponding set has vanishing measure 
by \eqref{18.23}. We also restrict to $r$ such that $j(r) > \cE^{1/2}$, by \eqref{20.8},
and notice that since $r\in \cR$,
\begin{eqnarray} \label{20.12}
\H^1(E_1 \cap \S_r) dr &\leq& v'(r) = r^2 \theta'(r) + 2 r^{-1} v(r)
=  r^2 \theta'(r) + 2 r \theta(r) 
\nn\\ &=&  r^2 \theta'(r) + 2 r f(r) + 2r \theta_0
\leq 2r \theta_0 + C r j(r)
\end{eqnarray}
by \eqref{18.22}, \eqref{18.21}, and then \eqref{19.27n}.
Thus
\begin{equation} \label{20.13}
\int_{r=10^{-2}}^3 \H^1(E_1 \cap \S_r) dr
\leq \theta_0 \int_{r=10^{-2}}^3 \, \1_{j(r) > \cE^{1/2}} \, r dr + C \int_{r=10-2}^3 r j(r)
\leq C \cE^{1/2}
\end{equation}
by \eqref{19.21} and Chebyshev. We compare with \eqref{20.11} and \eqref{20.10} and get that
\begin{equation} \label{20.14}
\H^2(E_0 \sm E_b) \leq C \cE^{1/2}.
\end{equation}
Unfortunately, we still have to take care of $E_b$, where the co-area formula does not seem
to work so well. In fact, if we apply \eqref{19.16} with $g = \1_{E_b}$, the left-hand side vanishes
(because $J(x)=0$ almost everywhere on $E_b$), and the right-hand side is 
\begin{equation}\label{21.16n}
0 = \int_{0}^{90} \Big\{ \int_{E \cap \S_r} g(x) d\H^1(x) \Big\} dr
= \int_{0}^{90} \H^1(E_b \cap \S_r) dr.
\end{equation}
Thus $\H^1(E_b \cap \S_r) = 0$ for almost every $r$, and 
the contribution of $E_b$ is not seen when we evaluate the variations of $v(r)$
using the integral of $v'$ and the estimate \eqref{18.22}.
That is, if we set $\wt E = E \sm E_r$, $\wt v(r) = \H^2(\wt E \cap B(0,r))$, 
and $\wt \theta(r) = r^{-2} \wt v(r)$, the proof of near monotonicity for $\theta$ also yields 
the near monotonicity of $\wt \theta$, as in \eqref{18.10}. 
It is a little sad that the author is forcing the reader to trust that the proof
of near monotonicity uses nothing else than \eqref{18.22}; in \cite{C1} the author gave an other proof
that avoids this unpleasant point, but at the same time is more complicated. Anyway, the near monotonicity
for $\wt \theta(r)$ yields
\begin{equation}\label{21.17n}
\wt\theta(90) \geq \theta_0 {\rm exp}\Big(-\alpha_n \int_0^{90} \frac{h(2t) dt}{t} \Big)
\geq \theta_0 - \theta_0 \alpha_n \int_0^{90} \frac{h(2t) dt}{t}
\end{equation}
by \eqref{18.1} (as in \eqref{18.11}), and because $\int_0^{90} \frac{h(2t) dt}{t}$ is small by \eqref{18.8}.
Then $\theta(90) =  \wt\theta(90)+ 90^{-2}\H^2(E_b)
\geq \theta_0 + 90^{-2}\H^2(E_b) - \theta_0 \alpha_n \int_0^{90} \frac{h(2t) dt}{t}$, hence
\begin{equation}\label{21.18n}
f(90) \geq 90^{-2}\H^2(E_b) - \theta_0 \alpha_n \int_0^{90} \frac{h(2t) dt}{t}
\end{equation}
or equivalently
\begin{equation}\label{21.19n}
\H^2(E_b) \leq 90^{2} f(90) + 90^{2} \theta_0 \alpha_n \int_0^{90} \frac{h(2t) dt}{t}
\leq C \cE
\end{equation}
by the definition of $\cE$ in Lemma \ref{t19.2}. With \eqref{20.14}, this shows that
$\H^2(E_0) \leq C \cE^{1/2}$.

Now assume that we can find $x\in E\cap A_0$ such that
$\dist(x,Z) \geq 2C_1 \cE^{1/4}$, where the large constant $C_1 > C_0$ will be chosen soon.
Set $\rho = C_1 \cE^{1/4}$; the set $E_3 = E \cap B(x,\rho)$ stays at distance at least
$C_1 \cE^{1/4}$ from $Z$, and it is also also contained in $B(0,3) \sm B(0,10^{-2})$
(because $x\in A_0 = B(0,2)\sm B(0,10^{-1})$ and
we may assume that $\cE$ is arbitrarily small), so $E_3 \subset E_0$. 
On the other hand, the local Ahlfors regularity of $E$ yields 
$\H^2(E_3) \geq C^{-1} \rho^2 = C^{-1} C_1^2 \cE^{1/2}$,
with a constant of the usual type, and that does not depend on $C_1$; we choose
$C_1$ large enough and get the desired contradiction with our upper bound for $\H^2(E_0)$.
This completes the proof of Lemma \ref{t20.2}.
\qed

\ms
\begin{rem} \label{t20.3}
In \cite{C1} we obtained a better power, namely $1/3$ instead of $1/4$.
We do not try to do this here, and send the reader to \cite{C1} instead in the unlikely
event where something like this would be needed. The general
idea was not hard: because of Lemma \ref{t20.1}, we already know that all points of $Z\cap A$
lie $C\cE^{1/3}$-close to $E$; we also know that $E$ is reasonably close, in any ball
$B_0$ centered on $E \cap A$ and with radius $10^{-2}$, say,
to our initial minimal cone $X$. The point is to use the fact that, in such a ball (and if we want, 
due to the fact that near $B_0$, the cone $Z$ is one step simpler than in the ball centered at the origin), 
we have a good description of $E$ in $10^{-1}B_0$, which we can use to say that it 
cannot look like $Z$, plus a tiny bit that goes away from $Z$. 
In the case of \cite{C1}, we showed that $E$ is locally
H\"older-equivalent to a cone of type $\bY$ or $\bP$; here we would use the results of
\cite{Mono} (or even this paper with a smaller density $\theta_0$) to get a good description of
$E$ near a point, that prevents additional spikes that go away from $Z$. In both case we use 
extra flatness instead of Ahlfors-regularity to get a better control of $E$ at the scale $\cE^{1/3}$
rather than $\cE^{1/4}$.
\end{rem}

\ms
We are now ready to prove Theorem \ref{t18.1}.
Our first observation is that if $E$ is as in the theorem, and we choose a new scale
$\rho \in (0,1/2)$, then the new set $E_\rho = \rho^{-1} E$ satisfies almost the same 
assumptions as $E$ itself. That is, the new gauge function for $E_\rho$ is
$h_\rho(r) = h(\rho r)$, and it satisfies \eqref{18.6} (even with the slightly smaller
constant $C_0 \rho^{\beta_0}$) if $h$ satisfies \eqref{18.6}. As for the 
analogue $f_{\rho}$ of $f$, notice that
\begin{equation} \label{20.15}
\begin{aligned}
f_{\rho}(200) 
&= f(200\rho) = \theta(200 \rho) - \theta_0
\leq  \theta(200) {\rm exp}\Big(\alpha_n \int_0^{200} \frac{h(2t) dt}{t}\Big) - \theta_0 
\cr& \leq [f(200) + \theta_0] \Big(1+2\alpha_n \int_0^{200} \frac{h(2t) dt}{t}\Big) - \theta_0
\leq f(200) + C \int_0^{400} \frac{h(t) dt}{t}
\end{aligned}
\end{equation}
because the density at the origin of $E_\rho$ is still $\theta_0$, and by \eqref{18.10}.
This is essentially as good as $f(200)$, i.e., when we assume \eqref{18.8} for $E$ with a slightly
smaller $\varepsilon_1$, we also get \eqref{18.8} for $E_\rho$ with $\varepsilon_1$.

Let us just consider $\rho_k = 2^{-k}$, with $k \in \bN$. 
For each $k$, we proceed as above, i.e., select a minimal cone $X = X_{\rho_k}$, 
then other cones $Z(r) = Z_{\rho_k}(r)$, $r\in \cR_{\rho_k}$, then 
a radius $r_k \in \cR_{\rho_k}$ that plays the role of $r_0$ above, 
and finally the cone $Z^{(k)} = Z_{\rho_k}(r_k)$
that we used for Lemmas~\ref{t20.1} and \ref{t20.2}. 

Notice that $r_{k+1}$ lies in the set $\cR_k$ that was used for the $k^{th}$ step,
that $j(r_{k+1})$ is actually the same when we think that $r_{k+1} \in \cR_{k}$ or
$r_{k+1} \in \cR_{k+1}$, and that we could have used the same cone 
$Z^{(k+1)} = Z_{\rho_{k+1}}(r_{k+1})$ as the set $Z_{\rho_{k}}(r_{k+1})$.
Then by Lemma \ref{t19.3} (applied with choice of $Z_{\rho_{k}}(r_{k+1})$),
\begin{equation} \label{20.16}
d_{0,1}(Z^{(k)},Z^{(k+1)}) = d_{0,1}(Z_{\rho_{k}}(r_k),Z_{\rho_{k}}(r_{k+1}))
\leq C \big(j(r_k) + j(r_{k+1}) + C\cE_k \big)^{1/2}
\end{equation}
where 
\begin{equation}\label{21.22n}
\cE_k = f(90\rho_k) + \int_0^{180\rho_k} h(t) \frac{dt}{t}
\end{equation}
is the analogue of $\cE$ at stage $k$, (see Lemma \ref{t19.2}). 
But $r_k$ and $r_{k+1}$ were chosen so that $j(r_k) \leq C \cE_k$ and
$j(r_{k+1}) \leq C \cE_{k+1}$ (see \eqref{20.1a}), so
\begin{equation}\label{21.23n}
d_{0,1}(Z^{(k)},Z^{(k+1)}) \leq C\big(\cE_k + \cE_{k+1}\big)^{1/2}.
\end{equation}
Notice that 
\begin{equation}\label{21.24n}
\cE_k \leq C \cE_j \ \text{ for } 0 \leq j < k,
\end{equation}  
by the near monotonicity of $f$ (or $\theta$), with the same proof as for \eqref{20.15}.
We claim that then
\begin{equation}\label{21.25n}
d_{0,1}(Z^{(j)},Z^{(k)}) \leq C (k-j) \cE_j^{1/2} \ \text{ for } 0 \leq j < k.
\end{equation}
For instance, if $z_j\in Z^{(j)} \cap \ol B(0,1)$, \eqref{21.23n} gives a point $z_{j+1} \in Z^{(j+1)}$ 
such that $|z_{j+1}-z_j| \leq C \big( \cE_j + \cE_{j+1} \big)^{1/2} \leq C \cE_j^{1/2}$,
by \eqref{21.24n};  we may assume that $z_{j+1} \in \ol B(0,1)$ (because its projection 
on $\ol B(0,1)$ still lies in the cone $Z^{(j+1)}$ and is not twice further. 
By induction, for every $\ell > j$ we can find $z_\ell \in Z^{(l)} \cap \ol B(0,1)$
such that $|z_{\ell}-z_j| \leq C(\ell-j)\cE_j^{1/2}$; we stop at $\ell = k$ and get 
the first half of \eqref{21.25n}. 
The converse is the same: any point $z_k \in Z^{(k)} \cap \ol B(0,1)$ is within 
$C \big( \cE_k + \cE_{k-1} \big)^{1/2} \leq C \cE_j^{1/2}$
of $Z^{(k-1)} \cap \ol B(0,1)$, and so on until we reach $Z^{(j)} \cap \ol B(0,1)$.

\ms
Now it is easy to see that
\begin{equation} \label{20.18}
d_{0,1}(E, Z) \leq C \cE_1^{1/4}.
\end{equation}
Recall that we chose $Z = Z^{(0)}$.
Let $x\in E \cap B(0,1)$ be given, and choose $k$ so that $2^{-k-1} \leq |x| \leq 2^{-k}$.
By Lemma \ref{t20.2} (applied to $E_{\rho_k}$), we can find $z\in Z^{(k)}$ such that 
$|z-\rho_k^{-1} x| \leq C \cE_k^{1/4} \leq C \cE^{1/4}$ (recall that $\rho_k = 2^{-k}$). 
We may as well take $z\in \ol B(0,1)$, because $|\rho_k^{-1} x| \leq 1$ and so the projection 
of $z$ on $\ol B(0,1)$ cannot be twice further.
Then by \eqref{21.25n} (with $j=0$) we can find $w\in Z$ such that $|w-z| \leq C k \cE^{1/2}$;
thus $|x-\rho_k w| = \rho_k |\rho_k^{-1} x - w| \leq C (1+k)\rho_k \cE^{1/4} \leq C \cE^{1/4}$ 
and we get the first half of \eqref{20.18}. 
The second half is done the same way, using \eqref{21.25n} and then Lemma \ref{t20.1}.

We already noticed in \eqref{18.41} that 
$\alpha(Z) = \alpha(Z(r_0)) \leq C j(r_0) \leq C \cE$, by \eqref{20.1a}. 
So we can use $Z$ to establish \eqref{18.7}; Theorem \ref{t18.1} follows.
\qed

\section{A partial conclusion and the tangent cone is unique}
\label{S21}

In this section we stop and think a little about what we have done so far, and prove the existence
of a tangent cone $X_0$ at the origin (i.e., the uniqueness of blow-up limit) in some circumstances, 
as well as a good approximation result by $X_0$ in small balls $B(0,r)$.

We shall systematically assume that $L$ is a line through the origin, 
\begin{equation} \label{22.1n}
\begin{aligned}
E \ &\text{ is a reduced sliding almost minimal set in } B(0,r_1) \subset \R^n, 
\cr&\hskip1.3cm \text{with a boundary condition coming from } L, 
\end{aligned}
\end{equation}
with a gauge function $h$ such that
\begin{equation} \label{22.2n}
h(r) \leq C_h r^\beta \ \text{ for } 0 < r \leq r_1
\end{equation}
for some constants $C_h \geq 0$, $\beta > 0$, and $r_1 > 0$. We also assume that
\begin{equation} \label{22.3n}
0 \in E \cap L.
\end{equation}
Our simplest result is the following.

\begin{cor} \label{t22.1n}
Let $E$ satisfy \eqref{22.1n}-\eqref{22.3n}, and suppose in addition that 
\begin{equation}\label{21.1}
\text{some blow-up limit of $E$ at $0$ satisfies the full length condition.}
\end{equation}
Then $E$ has a unique blow-up limit $X_0$ at $0$, and we can find $a>0$, 
$r_0 \in (0,r_1)$, and $C_1, C_2 \geq 0$ such that
\begin{equation}\label{22.4n}
f(r) \leq C_1 r^a \ \text{ and }\ d_{0,r}(E,X_0) \leq C_2 r^{a/4}
\ \text{ for } 0 < r < r_0,
\end{equation}
where as usual 
\begin{equation}\label{22.5n}
f(r) = \theta(r) - \theta_0, \text{ with } \theta(r) = r^{-2} \H^2(E \cap B(0,r))
\text{ and } \theta_0 = \lim_{t \to 0} \theta(t).
\end{equation}
Here $a$ depends only on $n$ and the full length constant for $X_0$ (which turns out to be the
unique blow-up limit), while $r_0$, $C_1$, and $C_2$ may depend wildly on $E$.
\end{cor}

The reader should not pay too much attention to the difference between $a$ and $a/4$; this is just how
they come in the proof.

This corollary generalizes Corollary \ref{t17.2}, and will apply automatically when 
$\theta_0 = \lim_{t \to 0} \theta(t) \\ %%% went to the next line to cut better
\leq \frac{3\pi}{2}$ if we ever prove \eqref{17.23}, by the
full length result of Section \ref{S30}.

As usual, we prefer to state a more precise result, with more quantifiers, where we start from
the good approximation of $E$ by a full length cone $X$ in a given ball, and get the existence of a tangent
cone $X_0$ and more precise approximation results for $X_0$ in smaller balls.

\begin{thm}  \label{t22.2n}
Let $X$ be a sliding minimal cone of dimension $2$ in $\R^n$, with sliding boundary condition coming 
from $L$, and assume that $X$ satisfies the full length condition. For each choice of constant
$\beta > 0$, we can find $\varepsilon_0 > 0$, such that if the sliding almost minimal set $E$ 
satisfies \eqref{22.1n}-\eqref{22.3n}, has the same density at $0$ as $X$, i.e., 
\begin{equation} \label{21.6}
\H^2(X\cap B(0,1)) = \lim_{r \to 0} r^{-2}\H^2(E\cap B(0,r)),
\end{equation}
and if in addition we can find $\varepsilon \in (0,\varepsilon_0]$ such that
\begin{equation} \label{21.7}
C_h r_1^{\beta} \leq \varepsilon
\ \text{ and }
d_{0,r_1}(E,X) \leq \varepsilon,
\end{equation}
then $E$ has a unique tangent cone $X_0$ at $0$, 
\begin{equation} \label{21.8}
\ d_{0,r}(E,X_0) \leq c_1(\varepsilon) \Big(\frac{r}{r_1}\Big)^{a/4}
\ \text{ for } 0 < r < r_1,
\end{equation}
and, with $f$ as in \eqref{22.5n}, 
\begin{equation} \label{22.9n}
f(r) \leq \Big(\frac{3r}{r_1}\Big)^a f(r_1/3) + C_3 C_h r_1^\beta \Big(\frac{r}{r_1}\Big)^a
\leq c_2(\varepsilon) (r/r_1)^a
\ \text{ for } 0 < r < r_1/3.
\end{equation}
Here $a > 0$, $\varepsilon_0$, and $C_3$ depend only on $n$, $\beta$, and $X$ through the geometric
constants $\theta_0$, $\eta(X)$, $\eta$, and $c$ associated to $X$ and its full length condition.
The constants $c_1(\varepsilon)$ and $c_2(\varepsilon)$ depend also on $\varepsilon$, and tend to $0$ 
(with $n$, $\beta$, $\eta(X)$, $\eta$, and $c$ fixed) when $\varepsilon$ tends to $0$.
\end{thm}

As we shall see at the end of this section, 
this result is stronger than the combination of Theorems~\ref{t1.3} and \ref{t1.4}, 
but not as trivially as one could think. We cannot apply Theorem~\ref{t22.2n} brutally
because the cones of $\bP$, $\bY$, and $\bT$, for instance, do not really satisfy 
the full length property with uniform constants, since the number $\eta(X)$ 
also measures the distance from $\ell \in L \sm K$ to the closest vertex of $K$, which may 
be arbitrarily small. 
We will finesse the issue by a small covering argument, as we did for Proposition \ref{t17.1}.

Theorem \ref{t22.2n} clearly implies Corollary \ref{t22.1n} (apply it with a full length blow-up limit
$X$ and a small enough radius $r_1$ such that \eqref{21.7} holds). 
In addition to the more precise estimates, it has an advantage over Corollary \ref{t22.1n} 
that we don't need to compute a blow-up limit of $E$; it is enough to approximate $E$ 
well enough by a full length minimal cone. However, we still need to know the density of $E$ at $0$, 
because of \eqref{21.6}. And the small constant $\varepsilon_0$ depends on our choice of $X$,
so we may have to make tough arbitrages between good approximation and large full length constants.

The sets $X$ and $X_0$ are not related a priori, but the proof will show that $X$ and $X_0$, 
in addition to having the same density, are based on the same model. 
That is, $X_0$ is a deformation of $X$ as in Definition \ref{t3.1}.

The general strategy for the proof will be to use Proposition \ref{t16.2}
to get some decay for $f$, and Theorem \ref{t18.1} to deduce from the size of $f(r)$
that $E$ lies close to a nice cone. We will have to do the two things at the same time,
because we also need the good approximation result of Theorem \ref{t18.1} to find a nice 
minimal cone for which \eqref{16.7} holds for smaller radii. That is, we will need to show at the same time 
that \eqref{21.8} and \eqref{22.9n} hold, for smaller and smaller radii $r$.

\ms
Before we turn to the proof, let us say that Theorem \ref{t22.2n} is not enough to 
give a good $C^1$ description of $E$ near $0$, even when the blow-up limits of $E$ at $0$ are simple.
Sure enough, we get a good control on $E$ in every small ball centered at $0$, but what about
small balls contained in $B(0,r_1/10)$, but centered at other points of $E\cap L$, and more importantly
at points of $E \sm L$? If we want to apply something like Reifenberg's topological disk theorem
to describe $E$ near $0$, it seems that we need a uniform control on (the approximation of $E$ by
nice cones in) these ball to get biH\"older descriptions, and even a uniform decay to get a $C^1$,
or slightly better than $C^1$, description.

We managed in \cite{Mono} to get enough uniform control on such balls to get a biH\"older 
description of $E$ near $0$ in some specific situations (when $E$ looks a lot like a half plane 
or a $\bV$-set in $B(0,r_1)$), and in the present paper we want a better control 
(better than $C^1$), and slightly more cases. Both here and in \cite{Mono}, 
we rely on the near monotonicity of a close relative of $\theta$,
the function $F$ of \eqref{1.20}, which is adapted to balls that are centered slightly off $L$.

This is why we will need to redo a large proportion of the arguments of this part in the next one, 
and adapt them to the situation of balls centered on $E \sm L$, but unfortunately with a limited 
list of approximating minimal cones (or truncated cones). In the mean time we prove 
Theorem \ref{t22.2n}, and then the fact that it implies Theorems \ref{t1.3} and \ref{t1.4}.

\msi
{\bf Proof of Theorem \ref{t22.2n}.}
Let $E$ and $X$ be as in the statement, and define $\theta$ and $f$ as in \eqref{22.5n}.
We want to prove that $f$ decays like a power, and for this we want to use Proposition~\ref{t16.2}
and integrate the differential inequality that it gives.

So we want to find a cone $X(r)$ such that \eqref{16.7} holds, and since we don't want a mess
with varying full length minimal cones (that we also would have to find anyway), the simplest will be 
to keep the same cone $X$ and hope that it works for all radii. 
This means that we will have to prove that it stays close to $E$
at small scales, which will be done with the help of Theorem \ref{t18.1}.

Anyway, we we want to apply Proposition \ref{t16.2} with $r_0 = r_1/2$. 
Let us first check the easy assumptions: \eqref{16.1} holds because of \eqref{22.1n}, 
and \eqref{16.2} follows from \eqref{22.2n}. For \eqref{16.8}, 
we work with the fixed cone $X$, so \eqref{16.8} just requires that $C_h r_0^\beta \leq \varepsilon(X)$ 
for some small constant $\varepsilon(X)$, and this follows at once from \eqref{21.7}. 
Finally, \eqref{16.10} holds with $q(r)=0$, by \eqref{21.6}.
We are left with \eqref{16.7}. Again we work with the fixed cone $X$, so \eqref{16.7} demands that 
\begin{equation}\label{22.11n}
d_{0,2r}(E, X) \leq \varepsilon(X);
\end{equation}
maybe we will not be able to prove this directly for all $r\in (0,r_1/2)$, so we define
\begin{equation}\label{22.12n}
r_{00} = \sup \big\{ r \in (0,r_1/2) \, ; \,  \eqref{22.11n} \text{ fails }\big\},
\end{equation}
with the convention that $r_{00} = 0$ if \eqref{22.11n} holds for all $r \in (0,r_1/2)$.
Notice however that since $d_{0,r_1}(E,X) \leq \varepsilon_0$ by \eqref{21.7}, we immediately 
get that \eqref{22.11n} holds for $r > \varepsilon(X)^{-1} \varepsilon_0 r_1$. That is, 
\begin{equation}\label{22.13n}
r_{00} \leq \varepsilon(X)^{-1} \varepsilon_0 r_1,
\end{equation}
which we can make as small as we wish (compared to $r_1$) by taking $\varepsilon_0$ small.

Eventually we shall prove that $r_{00} = 0$; in the mean time, set $I = (r_{00},r_1/2)$.
Our last condition \eqref{22.11n} is only known to hold on $I$, but fortunately it was observed 
a few lines after the proof of Proposition \ref{t16.2} that with our weaker assumptions 
(where \eqref{16.7} only holds for $r \in I$), the conclusion of Proposition \ref{t16.2}, 
i.e., the differential inequality \eqref{16.11}, still holds for all $r\in I$. This means that 
\begin{equation}\label{22.14n}
r f'(r) \geq \frac{4\alpha}{(1-2\alpha)} f(r) - 3 h(r)
\ \text{ for almost every } r \in I,
\end{equation}
for a fixed constant $\alpha = \alpha(X)$ (and because $q(r)=0$).
This differential inequality can be integrated on $I$ as we did in Section \ref{S16},
and we get the inequality \eqref{16.23}, valid for radii in $I$. We change notation
because $r_1$ is already taken, and get that for $r, s \in I$, $r \leq s$, 
\begin{equation} \label{21.10}
\begin{aligned}
f(r) &\leq \Big(\frac{r}{s}\Big)^{a} f(s) + 3 r^a \int_{r}^{s} h(t) \frac{dt}{t^{a+1}}
\leq \Big(\frac{r}{s}\Big)^{a} f(s) + 3 C_hr^a 
\int_{r}^{s} \frac{r^{\beta} dr}{r^{a+1}}
\cr& \leq \Big(\frac{r}{s}\Big)^{a} f(s) + \frac{6C_h}{\beta} \, r^{a} s^{\beta-a},
\end{aligned}
\end{equation}
where the small positive constant $a = \frac{4\alpha}{1-2\alpha}$ from \eqref{16.20}
depends on $n$, $\beta$ and $X$ (as above)
but not on $C_h$, and then by \eqref{17.2} and because we can safely assume that 
$a < \beta/2$.

Let us take $s= r_1/3$. 
We get that for $r \in (r_{00},r_1/3)$, 
\begin{equation}\label{22.17n}
f(r) \leq  \Big(\frac{r}{s}\Big)^{a} f(r_1/3) 
+ \frac{6C_h}{\beta} s^\beta \, \Big(\frac{r}{s}\Big)^a 
\leq \Big(\frac{3r}{r_1}\Big)^{a} f(r_1/3)  + C C_h r_1^{\beta} \Big(\frac{r}{r_1}\Big)^{a}.
\end{equation}
Notice that this is compatible with the first half
of \eqref{22.9n}, which therefore will follow as soon as we prove that $r_{00} = 0$.

Let prove the second inequality of \eqref{22.9n} now. 
Observe that $C_h r_1^{\beta} \leq \varepsilon_0$ by \eqref{21.7}, so we only need to show that
\begin{equation}\label{22.18nn}
f(r_1/3) \leq c(\varepsilon),
\end{equation}
with a constant $c(\varepsilon)$ such that $\lim_{\varepsilon \to +\infty} c(\varepsilon) = 0$.

We deduce this from the fact that $d_{0,r_1}(E,X) \leq \varepsilon$ (by \eqref{21.7}),
with a simple compactness argument, similar to the proof of 
Lemma 16.43 in \cite{Holder}, but based on the limiting arguments of \cite{Sliding}
because of the sliding condition. The point is that if this failed, we could construct a sequence
of sliding almost minimal sets $E_j$, and a sequence of sliding minimal cones $X_j$, 
both associated to the boundary $L$, so that (after a dilation that sends $r_1$ to $1$)
$d_{0,1}(E_j,X_j)$ tends to $0$ but the densities $\theta_j(1/3) = 9 \H^2(E_j \cap B(0,1/3)$ 
and $\H^2(X_j \cap B(0,1))$ stay far from each other. Then we would extract convergent sequences,
use Theorems 10.97 and 22.1 of \cite{Sliding} to control the densities, show that
$f_j(1/3) = 9\H^d(E_j \cap B(0,1/3)) - \H^d(X_j \cap B(0,1))$ tends to $0$, and get the desired
contradiction.

For the moment, we only know the first part for $r\in I$, but we still get that 
\begin{equation}\label{22.18n}
f(r) \leq  c \Big(\frac{r}{r_1}\Big)^{a} \ \text{ for } r \in (r_{00},r_1/3),
\end{equation}
with $c$ as small as we want.
We want to use this, and Proposition~\ref{t18.1}, to control the geometry of $E$, in particular 
in balls that are too small for \eqref{21.7} to give good results.

Set $r_k = 10^{-3k} r_1$, and apply Proposition~\ref{t18.1} to the set $E_k = r_k^{-1} E$.
The assumption \eqref{18.6} (with $\beta_0 = \beta$) follows at once from \eqref{22.2n}
and we can even take $C_0 = C_h r_k^{\beta} \leq 10^{-3k \beta}\varepsilon_0$ (by \eqref{21.7}).
So $C_0$ is small, and Proposition~\ref{t18.1} says that \eqref{18.7} holds, i.e., 
\begin{equation}\label{22.19n}
\beta_{\cX,\eta}(E_k,1) \leq C \Big[ f_k(200) + \int_0^{400} \frac{h_k(t) dt}{t} \Big]^{1/4},
\end{equation}
with $f_k(200) = f(200 r_k)$ and 
\begin{equation}\label{22.20n}
\int_0^{400} \frac{h_k(t) dt}{t} \leq \int_0^{400} \frac{C_h (r_k t)^{\beta}  dt}{t}
\leq C C_h r_k^\beta \leq C \varepsilon_0 10^{-3k \beta}
\end{equation}
by \eqref{21.7} again. We shall restrict our attention to integers $k \geq 2$ such that 
\begin{equation}\label{22.21n}
200 r_k \geq r_{00},
\end{equation}
because this way we can apply \eqref{22.18n} to $r = r_k$ and get that
\begin{equation}\label{22.22n}
\beta_{\cX,\eta}(E_k,1) 
\leq C c \Big(\frac{r}{r_1}\Big)^{a/4} + C \big(\varepsilon_0 10^{-3k \beta}\big)^{1/4}
\leq c_1\Big(\frac{r}{r_1}\Big)^{a/4},
\end{equation}
again with $c_1$ as small as we want, and because $a < \beta/2$ and $\varepsilon_0$ is small. 
By \eqref{18.5} this means that we can find a cone $Z_k \in \cZ(X,\eta)$ such that in particular
\begin{equation}\label{22.23n}
d_{0,r_k}(E,Z_k) = d_{0,1}(E_k,Z_k) \leq c_1 \Big(\frac{r_k}{r_1}\Big)^{a/4} = c_1 10^{-3ka/4}.
\end{equation}
We shall only apply this for $k \geq k_0$, where $k_0$ will be chosen soon. Notice that
for $k \geq k_0+1$, 
\begin{eqnarray}\label{22.24n}
d_{0,1}(Z_k,Z_{k-1}) &=& d_{0,1/2}(Z_k,Z_{k-1}) 
\leq 2d_{0,1}(Z_k,E_k) + 2d_{0,1}(E_k,Z_{k-1})
\nn\\
&=& 2d_{0,1}(Z_k,E_k) + 2 \cdot 10^3 d_{0,10^{-3}}(E_{k},Z_{k-1})
\nn\\
&=& 2d_{0,1}(Z_k,E_k) + 2 \cdot 10^3 d_{0,1}(E_{k-1},Z_{k-1}) \leq c_2 10^{-3ka/4},
\end{eqnarray}
with $c_2$ as small as we want and where for the first line (and similar computations later) 
we actually use the fact that our estimates for the normalized distances (that follow) are small, 
so that we can chase points inside $B(0,1)$.
For $k\leq k_0$, we prefer to use the fact that
\begin{equation}\label{22.25n}
d_{0,r_k}(E,X) \leq 10^{3k} d_{0,r_1}(E,X) \leq 10^{3k}\varepsilon \leq 10^{3k}\varepsilon_0
\end{equation}
by \eqref{21.7}. Let us pick $k_0$ so large that
\begin{equation}\label{22.26n}
\sum_{k \geq k_0} (c_1+c_2) 10^{-3ka/4} \leq 10^{-5}\varepsilon(X),
\end{equation}
where $c_1$ and $c_2$ come from \eqref{22.23n} and \eqref{22.24n} 
(and we do not care yet whether they are small or not) and $\varepsilon(X)$ comes
from \eqref{22.11n}. Also make sure to pick $\varepsilon_0 \leq 10^{-3k_0-5}\varepsilon(X)$, 
so that by \eqref{22.25n}
\begin{equation}\label{22.27n}
d_{0,r_{k_0}}(E,X) \leq 10^{3k} d_{0,r_1}(E,X) \leq 10^{-5} \varepsilon(X)
\end{equation}
and, by the same proof as for \eqref{22.24n},
\begin{eqnarray}\label{22.28n}
d_{0,1}(Z_{k_0},X) &=& d_{0,1/2}(Z_{k_0},X) 
\leq 2d_{0,1}(Z_{k_0},E_{k_0}) + 2d_{0,1}(E_{k_0},X)
\nn\\
&=& 2d_{0,r_{k_0}}(Z_{k_0},E) + 2d_{0,r_{k_0}}(E,X)
\nn\\
&\leq& 2c_1 10^{-3k_0 a/4}+ 2 \cdot 10^{3k_0}\varepsilon \leq 4 \cdot 10^{-5} \varepsilon(X)
\end{eqnarray}
by \eqref{22.23n} and \eqref{22.25n} (and because $\varepsilon \leq \varepsilon_0$). 
We return to \eqref{22.23n} and get that
\begin{eqnarray}\label{22.29n}
d_{0,r_k/2}(E,X) &=& 2d_{0,r_k}(E,Z_k)+ 2d_{0,r_k}(Z_k,X) 
\leq 2c_1 10^{-3ka/4} + 2d_{0,1}(Z_k,X)
\nn\\
&\leq& 2c_1 10^{-3ka/4} + 8 \cdot 10^{-5} \varepsilon(X) \leq 10^{-4} \varepsilon(X).
\end{eqnarray}
Hence
\begin{equation}\label{22.30n}
d_{0,400r_{k+1}}(E,X) = d_{0,4r_{k}/10}(E,X) \leq 2 d_{0,r_k/2}(E,X) \leq 10^{-3} \varepsilon(X).
\end{equation}

We are now ready to prove that $r_{00} = 0$. Let $k_0$ be as above; because of \eqref{22.13n},
we can choose $\varepsilon_0$ so small that $k_0$ still satisfies \eqref{22.21}. Let
$k \geq k_0$ be such that \eqref{22.21} holds. Then \eqref{22.30n} holds too, and says that
$200r_{k+1} > r_{00}$ (compare with \eqref{22.12n} and \eqref{22.11n}). 
That is, we can show by induction that $200r_k > r_{00}$ for all $k$, as needed.

\ms
As was said earlier, \eqref{22.17n}, and hence \eqref{22.9n} are now proved for all $r < r/3$.
Now we go for \eqref{21.8}. Now every $k \geq k_0$ satisfies \eqref{22.21n}, 
and by \eqref{22.24n} the sequence $\{ Z_k \}$ converges to a limit $X_0$.
By \eqref{22.23n},  \eqref{22.24n}, and the same computations as for \eqref{22.29n}
\begin{equation}\label{22.31n}
d_{0,r_k/2}(E,X_0) \leq 2d_{0,r_k}(E,Z_k) +2d_{0,r_k}(Z_k,X_0) \leq c 10^{-3ka/4}
= c \Big(\frac{r_k}{r_1}\Big)^{a/4},
\end{equation}
with $c$ as small as we want. If $0 < r < r_{k_0}/2$, we can apply \eqref{22.31n}
to the smallest $r_k$ such that $B(0,r) \subset B(0,r_k/2)$, and we get that 
\begin{equation}\label{22.32n}
d_{0,r}(E,X_0) \leq 10^3 d_{0,r_k/2}(E,X_0) \leq 10^3 c \Big(\frac{r}{r_1}\Big)^{a/4}.
\end{equation}
Since we can make $c$ arbitrarily small by taking $\varepsilon$ small, this takes care of the small 
radii in \eqref{21.8}. As usual, for the large radii we will try to use \eqref{21.7}. First observe that
\begin{eqnarray}\label{22.33n}
d_{0,1}(X_0,X) &\leq& 2 d_{0,1}(X_0,Z_{k_0}) + 2 d_{0,1}(Z_{k_0},X)
\leq C c_2 10^{-3k_0a/4} + 2 d_{0,1}(Z_{k_0},X)
\nn\\
&\leq& C c_2 10^{-3k_0a/4} +2 c_1 10^{-3k_0a/4} + 4 \cdot 10^{3k_0} \varepsilon
= c_3 \Big(\frac{r_{k_0}}{r_1}\Big)^{a/4} + 4 \cdot 10^{3k_0} \varepsilon
\end{eqnarray}
by \eqref{22.24n} and \eqref{22.28n}, and where $c_3$ is still as small as we want.
Now we deduce from \eqref{21.7} that for $r_{k_0}/2 \leq r \leq r_1/2$,
\begin{eqnarray}\label{22.34n}
d_{0,r}(E,X_0) &\leq& d_{0,r}(E,X) + 2 d_{0,1}(X_0,X)
\leq \frac{r_1}{r} d_{0,r_1}(E,X) + 2 d_{0,1}(X_0,X)
\nn\\
&\leq& \frac{r_1}{r} \varepsilon + 2c_3 \Big(\frac{r_{k_0}}{r_1}\Big)^{a/4} + 8 \cdot 10^{3k_0} \varepsilon
\leq 3 c_3 \Big(\frac{r}{r_1}\Big)^{a/4} + 10^{3k_0+1} \varepsilon
\end{eqnarray}
where for the first inequality we used again that $X$ and $X_0$ are cones.
The first term is all right for \eqref{21.8}, and for the second term, notice that
\begin{equation}\label{22.35n}
\Big(\frac{r}{r_1}\Big)^{-a/4} 10^{3k_0+1} \varepsilon 
\leq 2\Big(\frac{r_{k_0}}{r_1}\Big)^{-a/4} 10^{3k_0+1} \varepsilon
= 20 \Big(\frac{r_{k_0}}{r_1}\Big)^{-1-a/4} \varepsilon
\end{equation}
is still as small as we want, because $k_0$ was chosen in terms of $\varepsilon(X)$, and
$\varepsilon$ is as small as we want.
This completes the proof of Theorem \ref{t22.2n}. 
\qed

\ms
{\bf Proof of Theorems \ref{t1.3} and \ref{t1.4}.}
In the general case, although Theorems \ref{t1.3} and \ref{t1.4} correspond to two different
estimates (decay for $f$ and good approximation by cones when $f$ is small), we prove them at
the same time. Also, the quantifiers in the statement force us to get constants that do not depend
on how close the spine of an initial approximating cone to $E$ can get to $L$, without containing
half of it, so we will use the compactness of the following class of minimal cones.

Denote by $\cX_0$ the class of minimal cones of type $\bP$, $\bY$, $\bT$, $\bH(L)$,
or $\bV(L)$ (the same as in the statement of Theorems \ref{t1.3} and \ref{t1.4}).
Then let $n \geq 3$ and $\beta >0$ be given. For Theorem~\ref{t1.3} we want to find a 
small constant $\varepsilon_0 > 0$ such that the good estimate \eqref{1.17} holds as soon as 
$E$ and $B(0,r_1)$ satisfy the assumptions. For Theorem \ref{t1.4}, we would also 
give ourselves $c > 0$ small, and we would need to get \eqref{1.18} with $c(\varepsilon_0) < c$.

Let $L \subset \R^n$ be fixed (we can always do this by rotation invariance), 
and for each $X \in \cX_0$, observe that $X$ satisfies the full length property
(by Theorem \ref{t30.1}) and denote by $\varepsilon_0(X)$ the small constant 
$\varepsilon$ given by Theorem \ref{t22.2n} (applied so that $c_1(\varepsilon) < c$ and 
$c_2(\varepsilon) < 10^{-10}$). Then cover $\cX_0$, as we did with \eqref{17.13}, 
by the small balls 
\begin{equation} \label{22.36n}
V_X = \big\{ Y \in \cX_0 \, ; \, d_\H^c(X,Y) < 10^{-1} \varepsilon_0(X)) \big\}.
\end{equation}
Since $\cX_0$ is compact, we just need a finite family $\cal{Y} \subset \cX$ to cover,
and we take $\varepsilon_0 = \frac1{10} \min\big\{ \varepsilon_0(Y) \, ; \, Y\in \cal{Y}\big\}$.
Let us check that this work. Let $E$ and $r_1 > 0$ satisfy the assumptions of Theorem \ref{t1.3} 
or \ref{t1.4}; then in particular there is a minimal cone $X \in \cX_0$, with the same density as $E$
(as in \eqref{1.15}), and such that $d_{0,r_1}(E,X) \leq \varepsilon_0$. Then
$X \in V_Y$ for some $Y \in \cal{Y}$, and this implies that
$d_{0,r_1}(E,X) \leq 3 \varepsilon_0 < \varepsilon_0(Y)$. For this, since we find it neater
not to modify the ball where we look, we use the triangle inequality and also the fact that $X$
and $Y$ are cones.

The other assumptions of Theorem \ref{t22.2n} are also satisfied (because $Y$ has the same density
as $X$; we could also have fixed the density of $E$ at $0$ (out of a set of four values), and
restricted to cones that have this density); now the conclusions of Theorem \ref{t22.2n}
implies the conclusion of Theorems \ref{t1.3} and \ref{t1.4}, and this completes the proof of these
theorems. \qed

\vfill\eject
\part{Decay and approximation for balls centered on $E \sm L$}

\ms
In this part we still consider a sliding almost minimal set $E$, with a sliding boundary condition that 
comes from a line $L$, and we generalize some of the results of the previous two parts to balls 
that are centered on $E \sm L$. 

Our starting point is the near monotonicity of the quantity $F$ of \eqref{1.20}, which was proved in
\cite{Mono}. We show that when $E$ is sufficiently close to a half line, a plane, a set of type $\bV$,
or a short truncated cone of type $\bY$, this quantity actually decays like a power. 
This analogue of Part II will be done, with the same sort of method, in Sections \ref{S22}-\ref{S27}. 
For this we will need to prove analogues of the full length condition in specific situations, 
and these computations, done in Sections \ref{S26} and \ref{S27}, will also be used
in Section \ref{S30}, when we complete the full length verification for balls centered on $L$.

In addition, we will show in Section \ref{S28} that in the same circumstances as above, 
$F$ controls the geometry of $E$. This will allow us to get good approximation properties of $E$ 
in balls that are centered slightly off $L$, as needed if we want to apply Reifenberg-type constructions.

Let us just describe a situation where we will obtain something. Suppose that at the unit scale, 
$E$ looks like a set of type $\bY$ truncated by $L$, with a spine $E_Y$ that contains the origin,
runs almost parallel to $L$, and lies very close to $L$. 
At this large scale, $E$ looks essentially like a $\bV$-set, with angle $\2$. 
In very small balls near $0$, $E$ looks like a full $\bY$-set. 
We are interested in what happens at intermediate scales, and in particular in proving some decay
for quantities that show how well $E$ is approximated by truncated $\bY$-sets. 
This will be our way of proving that nothing wild happens between the two extreme scales, 
and even that the approximation at the small scale is better than expected.

We will see this sort of situation in the next part, where we use the decay information from this part 
to start the desired classification of sliding almost minimal sets near the boundary.

\section{Balls centered on $E \sm L$: preliminaries } 
\label{S22}

In this section we set the stage for a study of decay properties of the adapted density
function $F$, for balls that are centered on $E \sm L$. We will proceed like in the
previous sections, except that the functional $F$ has an additional term and the obvious 
competitors for $E$ are no longer cones over $E \cap \S_r$, but slightly larger sets with 
an additional triangular piece that allows retractions on the sets which preserve $L$.

In this section and the next ones, we assume that $L$ is a line, no longer through the origin, 
and that
\begin{equation} \label{23.1n}
\begin{aligned}
E \ &\text{ is a reduced sliding almost minimal set (of dimension $2$)} 
\cr&\hskip0.2cm \text{in $B(0,R)$, with a boundary condition coming from } L, 
\end{aligned}
\end{equation}
with a gauge function $h$ such that
\begin{equation} \label{23.2n}
h(r) \leq C_h r^\beta \ \text{ for } 0 < r \leq R,
\end{equation}
for some constants $C_h \geq 0$, $\beta > 0$, and $R > 0$. Also we assume that
\begin{equation} \label{23.3n}
0 \in E \sm L.
\end{equation}
The results that will be proved here would still hold when $L$ contains $0$, with essentially 
the same proof, but there would be no point because in this case the previous part gives us 
what we need, and also it will be convenient in some places to discuss things in terms of 
\begin{equation} \label{22.1}
d_0 = \dist(0,L) > 0.
\end{equation}
Let us review some of the notation and results of \cite{Mono}.
We shall be interested in balls $B(0,r)$, $0 < r \leq R$.
The shade of $L$ (seen from the origin) is
\begin{equation} \label{22.2}
S = \big\{ z\in \R^n \, ; \,  \lambda z  \in L \text{ for some } \lambda \in [0,1]\big\}.
\end{equation}
We keep $\theta(r) = r^{-2} \H^2(E \cap B(0,r))$ as it was, but now consider
\begin{equation} \label{22.3}
F(r) = \theta(r) + r^{-2} \H^2(S \cap B(0,r)) 
= r^{-2} \big[\H^2(E \cap B(x_0,r)) + \H^2(S \cap B(0,r)) \big]. 
\end{equation}
Notice that we take the sum, and not the measure of the union. 

Let us review the properties of $F$ that we intend to use. First assume more, i.e. that
\begin{equation}\label{22.4}
\begin{aligned}
&\text{$E$ is a sliding reduced minimal set in $B(0,R)$ of dimension $2$, }
\cr&\hskip 1.5cm
\text{with a boundary condition coming from $L$.}
\end{aligned}
\end{equation} 
In this case, $F$ is nondecreasing on the interval $(0,R)$. See Theorem 1.2 in \cite{Mono}.

There are two special cases of sliding minimal sets for which $F$ is constant.
The first one is the half plane $H_0$ bounded by $L$ and that contains the origin; it is
easy to see that for $H_0$, $F$ is constant equal to $\pi$ (the measure of the shade
exactly compensates for the missing half plane).

The second one is the truncated $\bY$-set $Y_0$, which is $Y_0 = \ol{Y_1 \sm S}$, 
where $Y_1$ is the only cone of type $\bY$ that is centered at $0$ and contains $L$ 
(thus its singular set is parallel to $L$ and $S \subset Y_1$).
For this set $Y_0$, $F$ is constant and equal to $\frac{3\pi}{2}$.

We also have local slightly tilted variants of $Y_0$. If $Y_1$ is a cone of type $\bY$
such that $L\cap B(0,R)$ is contained in one of the three faces of $Y_1$
(and hence $S \cap B(0,R)$ is also contained in that same face, by elementary geometry),
$Y = \ol{Y_1 \cap B(0,R) \sm S}$ is also a sliding minimal set in $B(0,R)$ 
(at least, we claim that this is very probable but we won't need to check),
and the function $F$ attached to it is constant and equal to $\frac{3\pi}{2}$ on $(0,R)$.

Theorem 1.3 in \cite{Mono} gives a nice description of $E$ when \eqref{22.4} holds and 
$F$ is constant on an interval, but we shall only need the following two specific cases, 
which give a converse to the examples above.

\begin{lem}\label{t22.1} 
Suppose that \eqref{22.4} holds and $R > d_0 > 0$.
If $F(r) = \pi$ for $r\in (0,R)$, then $E = H_0 \cap B(0,R)$.
If $F(r) = \frac{3\pi}{2}$ for $r\in (0,R)$, then there is a cone $Y_1$ of type $\bY$,
centered at $0$, such that $L\cap B(0,R)$ is contained in one of the three faces of $Y_1$, 
and for which $E = \ol{Y_1 \cap B(0,R) \sm S}$.
\end{lem}

Notice that we already know, from previous work on the situation with no sliding boundary
(probably even  before \cite{Ta}), that since $F(r) = \theta(r)$ is constant and equal to $\pi$ or 
$\frac{3\pi}{2}$ on $(0,d_0)$, $E$ coincides with a plane or a cone of type $\bY$ on
$B(0,d_0)$. But let us apply Theorem~1.3 in \cite{Mono}, with $R_0$ very small and $R_1=R$. 
Recall that ``coral'' is the same as ``reduced'', so the assumptions are satisfied. 
Set $A = B(0,R_1)\sm B(0,R_0)$ as in \cite{Mono}. Let $X$ be the positive cone over
$E \cap A$ (as in (1.13) there). We get that $X$ is a reduced minimal set in $\R^n$ 
(that is, with no boundary condition), and that $A \cap X \sm S \subset E$
(as in (1.14) there), and where $S$ is still the shade of $L$ (see (1.9) there)).
Thus in $B(0,d_0) \sm B(0,R_0)$, $X$ coincides with $E$ (by definition of $X$,
$X \supset E\cap A$), and since $F(r) = \theta(r)$ for $r < d_0$, we get that
the density of $X$ is $\pi$ or $\frac{3\pi}{2}$, hence $X$ is a plane or a cone of type
$\bY$.

It was also observed after the statement of Theorem 1.3 in \cite{Mono} that in $A$,
$E$ and $X \sm S$ coincide modulo a set of vanishing $\H^2$-measure. They also 
coincide in $B(0,R_0)$: either use the fact that $E$ is a plane or a $\bY$ in $B(0,d_0)$,
or observe that $X$ cannot depend on $R_0$ and let $R_0$ tend to $0$).
That is
\begin{equation}\label{22.5}
E \cap B(0,R) = (X \sm S) \cap B(0,R), \text{ modulo a set of vanishing $\H^2$-measure.}
\end{equation}
Then, for $r\in (d_0,R)$,
\begin{eqnarray}\label{22.6}
\H^2(X\cap B(0,r)) &=& r^2 \H^2(X\cap B(0,1)) = r^2 F(r)
= \H^2(E\cap B(0,r)) + \H^2(S\cap B(0,r))
\nn\\
&=& \H^2((X \sm S)\cap B(0,r)) + \H^2(S\cap B(0,r))
\end{eqnarray}
because $X$ is a cone, $F$ is constant, by \eqref{22.3}, and by \eqref{22.5}. 
This forces $X$ to contain almost all of $S\cap B(0,r)$.

If $X$ is a plane, this forces $X$ to contain a bit of $L$, then the whole $L$;
thus $X$ is the plane that contains $H_0$ and the result follows from \eqref{22.5}
(and the fact that $E$ is closed and reduced).

If $X$ is a cone of type $\bY$, $L \cap B(0,r) \subset X$ as above, and
since this is true for all $r\in (d_0,R)$, $X$ contains $L\cap B(0,R)$.
In fact, $L\cap B(0,R)$ is contained in a single face of $X$
(if $L\cap B(0,R)$ crossed the spine of $X$, one piece of it would not lie in $X$),
so we can take $Y = X$ in the description above. Again the fact that
$E = \ol{Y_1 \cap B(0,R) \sm S}$ follows from \eqref{22.5}.
\qed

\ms
We shall also need the simpler version of Lemma \ref{t22.1} where $0 \in L$. 
We start with a description of sliding minimal cones with low density.
Denote by $\bP_0$ the set of planes through the origin.

\begin{lem} \label{t22.2}
There is a small constant $\tau(n) > 0$ such that if $X$ is a sliding minimal cone 
of dimension $2$ in $\R^n$, with a sliding condition coming from a line $L$ that contains 
the origin, and if $\H^2(X\cap B(0,1)) \leq \pi + \tau(n)$, then $X \in \bH(L) \cup \bP_0 \cup \bV(L)$, 
i.e., $X$ is a half plane bounded by $L$, a plane that contains the origin, but not necessarily $L$), 
or a set of type $\bV$ associated to $L$.
\end{lem}

See the beginning of Subsection \ref{S1.2} for the definitions.
Notice that this is a simpler special case of \eqref{17.23}, wich at least we can prove.
We start with the apparently even weaker statement with $\tau(n)=0$.
Let $X$ be as in the statement, with a density $d(X) = \H^2(X\cap B(0,1)) \leq \pi$.
Let us use the description of $K = X \cap \d B(0,1)$ that was given in Proposition~\ref{t2.1}. 
We see that $K$ is a union of great circles and arcs of great circles. 
If $K$ contains a great circle, this eats all the available density, $K$ is a great circle, 
and $X$ is a plane. Otherwise, $K$ is a union of arcs of geodesic.

Suppose two such arcs meet at some point $y\in K \sm L$. Then there are three arcs of $K$
meeting at $y$ (with $120^\circ$ angles, but we don't care), and the density of $X$ at
$y$ is at least $3\pi/2$. This means that $\lim_{r \to 0} F_y(r) = 3\pi/2$, where $F_y$ is
the functional defined as in \eqref{22.3}, but with the set $X$ and an origin at $y$. That is,
$F_y(r) = r^{-2}\big[\H^2(X\cap B(y,r)) + \H^2(S_y\cap B(y,r))\big]$, where $S_y$
denotes the shade of $L$ seen from $y$. It is easy to see that 
$\lim_{r \to +\infty} r^{-2}\H^2(X\cap B(y,r)) = \lim_{r \to +\infty} r^{-2}\H^2(X\cap B(0,r))
= d(X) \leq \pi$, hence $\lim_{r \to +\infty} F_y(r) \leq 3\pi/2$. But $F_y$ is nondecreasing,
so $F_y(r) = 3\pi/2$ for $0 < r < +\infty$. By Lemma \ref{t22.1}, $X$ coincides in large balls
$B(0,R)$ with truncated cones of type $\bY$, but centered at $y$. This contradicts the fact that
$X$ is a cone centered at $0$. 

Thus none of the arcs that compose $K$ ends away from $L$, which means that $K$ is composed
of half circles with endpoints in $L$. There is no more than two arcs, because $d(X) \leq \pi$.
If there is one arc, $X\in \bH$. Otherwise, $X$ is composed of two half planes, and $X\in \bV$
because if these two half plane make an angle smaller than $\2$, it is easy to see (or well known)
that $X$ is not minimal. 

We still need to prove the result with a positive $\tau(n)$. 
Suppose the lemma fails, so that for each large integer $k$ we can find a sliding minimal cone $X_k$, 
such that $d(X_k) \leq \pi + 2^{-k}$ and yet $X_k \notin \bH(L) \cup \bP_0 \cup \bV(L)$.
Take a subsequence (and still call it $\{ X_k \}$) such that $X_k$ converges to a limit cone $X$
(or equivalently here, since we work with cones, $K_k = X_k \cap \d B(0,1)$ 
converges to $K_\infty = X \cap \d B(0,1)$ for the Hausdorf distance on $\d B(0,1)$).

By the various convergence theorems in \cite{Sliding}
(Theorems 21.3, 10.97, and 22.1 there)
$X$ is a minimal cone and $d(X) = \lim_{k \to +\infty} d(X_k) \leq \pi$. 
By the case we already proved, $X \in \bH(L) \cup \bP_0 \cup \bV(L)$.
Let $y_k$ be any endpoint of an arc of $K_k$ that does not lie in $L$.
Such a point exists, because  $X_k \notin \bH(L) \cup \bP_0 \cup \bV(L)$
and by the argument above.

If we could find a subsequence for which $y_k$ converges to a limit
$y\in K \sm L$, then $K$ would have a point of type $\bY$ at $y$,
because $\{ K_k \}$ converges to $K$ and all the arcs of $K_k$ that do not end
on $L$ have lengths at least $\eta_0$ for some constant $\eta_0 = \eta_0(n)$.
Thus the endpoints $y_k$ all tend to $L$. For each $k$ large, 
$K_k$ has at most two short arcs that leave from the points $\ell_\pm$ of
$L \cap \d B(0,1)$ (see \eqref{2.4}), and all the other ones are long, because
they go from a small neighborhood of $\ell_-$ to a small neighborhood of $\ell_+$.
Thus there are at most two long ones (because $\H^1(K_k) = 2d(X_k) \leq 2\pi + 2^{-k+1}$).
If there is a single $y_k$, then $K_k$ is composed of two long arcs (from $\ell_-$, say, to $y_k$)
and a short one (the geodesic from $y_k$ to $\ell_+$). This is impossible, because the long
arcs are geodesics that both leave from $\ell_-$; they can only meet back at $\ell_+$.

We are left with the case when there are two points $y_k$ and $y'_k$, and
$K_k$ is composed of two geodesics from $y_k$ to $y'_k$, plus two short
geodesics from these points to the closest $\ell_\pm$. As before, the long geodesics
can only meet at the antipode, i.e., $y'_k = - y_k$. It is easy to see that the corresponding
set is not minimal. For instance, if the three arcs make the correct angles of $120^\circ$
at $y_k$, then the three arcs at $y'_k$ make acute angles of $60^\circ$.

This completes our contradiction and compactness argument; Lemma \ref{t22.2}
follows.
\qed

\ms
Let us continue our rapid description of the results of \cite{Mono}. 
We return to the more general situation where $E$ is a sliding almost minimal set,
as in \eqref{23.1n}, that $0 \in E \sm L$, and that the associated gauge function $h$ is such that 
\begin{equation}\label{22.7}
A(r) = \int_0^r h(t) \frac{dt}{t} < +\infty \ \text{ for } 0 < r < R,
\end{equation}
and $h(R)$ is small enough (depending on $n$). 
Then Theorem 1.5 in \cite{Mono} says that 
\begin{equation}\label{22.8}
F(r) e^{\alpha A(r)} \ \text{ is a nondecreasing function on } (0,R),
\end{equation}
with a constant $\alpha$ that depends only on $n$.

It will be psychologically useful to know the general idea of the proof, which is essentially
the same as for the (near) monotonicity of $\theta$ when $L$ is a cone centered at the origin. 
We would like to compare $E$ with the cone
\begin{equation} \label{22.9}
\Gamma(E,r) = \big\{ \lambda z\, ; \,  z  \in E \cap \d B(0,r) \text{ and } \lambda \in [0,1] \big\}
\end{equation}
over $E \cap \d B(0,r)$, but since it may no longer be a limit of sliding competitors
(moving a point $z\in E \cap L \cap \d B(0,r)$ in the direction of $0$ may detach it from $L$),
we add to $\Gamma(E,r)$ the set
\begin{equation} \label{22.10}
T(r) = \big\{ \lambda z\, ; \,  z  \in L \cap \ol B(0,r) \text{ and } \lambda \in [0,1] \big\},
\end{equation}
which is the convex hull of the triangle with vertices $0$ and the two points of
$L \cap \d B(0,r)$.  It turns out that $\Gamma(E,r) \cup T(r)$ can be used as a competitor 
(maybe, after taking a limit), just as $\Gamma(E,r)$ before.
Now $\Gamma(E,r) \cup T(r)$ is not as small as $\Gamma(E,r)$, and when 
we do the computation, we find out that we only get the (near) monotonicity of $F(r)$, 
where we added the (sometimes strictly) nondecreasing term $r^{-2} \H^2(S \cap B(x_0,r))$.

We will be more interested in the case when $0$ lies very close to $L$. Then $T(r)$ is quite thin; 
nonetheless it has an effect on the functional $F$ and on our estimates that we cannot neglect.

In the work that we did so far, with balls centered on $L$, the main point was to try to 
construct a competitor for $E$ that was significantly better than the cone $\Gamma(E,r)$,
and then we proved some decay for $\theta$ (i.e., a good differential inequality)
rather than proving that it is nearly monotone. Here we want to do the same thing,
i.e., improve significantly over $\Gamma(E,r) \cup T(r)$, and then we'll get a good differential inequality
involving $F$. As in the previous sections, the main point is the construction of good competitors.
This is what we do in the next two sections, in the two special cases for which we know that the
function $F$ can be constant on some truncated minimal cones.

But before we come to this, let us also show how to use Lemmas \ref{t22.1}
and \ref{t22.2}, and a little bit of compactness, to get similar results for almost minimal
sets. We now assume that $d_0 > 0$ (as in \eqref{23.3n}) and that $h$ satisfies \eqref{23.2n};
this way there exists a density
\begin{equation} \label{22.11}
\theta_0 = \lim_{r \to 0} \theta(r) = \lim_{r \to 0} F(r)
\end{equation}
because $d_0 > 0$, and by \eqref{22.8} or more simply its version in the plain case.
We start with an application of Lemma \ref{t22.1}, where we show that $E$
is some times close to a half plane.

\begin{lem}\label{t22.3} 
For each choice of small constants $\delta>0$ and $\tau > 0$, we can find 
$\varepsilon = \varepsilon(\delta, \tau) > 0$, that depends only on $\delta$, 
$\tau$ and $n$, with the following property. 
Let $E$ satisfy \eqref{23.1n}-\eqref{23.3n},  
and let $r$ be such that
\begin{equation}\label{22.12}
 r \leq \frac{d_0}{\delta}  \ \text{ and }\  \frac{11d_0}{10} \leq r < R.
\end{equation}
Suppose in addition that 
\begin{equation}\label{22.13}
h(r) \leq \varepsilon \ \text{ and } \ \int_0^{r} h(t) \frac{dt}{t} \leq \varepsilon,
\end{equation}
and
\begin{equation}\label{22.14}
F(r) \leq \pi + \varepsilon.
\end{equation}
Let $H_0$ denote the half plane plane bounded by $L$
that contains the origin. Then
\begin{equation}\label{22.15}
d_{0,\frac{20 r}{21}}(E, H_0) \leq \tau,
\end{equation}
and also 
\begin{equation}\label{22.16}
|\H^2(E \cap B(y,t)) - \H^2(H_0 \cap B(y,t))| \leq \tau r^2
\end{equation}
for all $y \in \R^n$ and $t > 0$ such that $B(y,t) \subset B(0,\frac{20 r}{21})$.
\end{lem}

\ms
It is important here to have in mind that when $r$ gets too large compared with $d_0$,
we need to take $\delta$ large (because of \eqref{22.12}), so we may need to take 
$\varepsilon$ very small. This is not shocking, it is just a reminder of the fact that 
limiting arguments (that will be used to prove the lemma) will only lead you so far.
The case when $r >> d_0$ will be discussed later.

We shall deduce this lemma from Theorem 1.6 in \cite{Mono}, whose main point is that
when the function $F$ is nearly constant on an interval, $E$ is quite close to a minimal set 
for which $F$ is constant. We shall apply that theorem with a fixed line $L_0$, which we 
choose so that $\dist(0,L_0) = 1$ (otherwise, the constants would depend on the line, 
and we want to avoid this). Let $f : \R^n \to \R^n$ be a composition of a 
rotation and a dilation, which we choose so that $f(0) = 0$ and $f(L) = L_0$. Thus
the dilation factor is $d_0^{-1}$. We want to apply the theorem to $E' = f(E)$,
so we check the assumptions, with $\tau' = \tau/2$ and the radius $r_1 = d_0^{-1}r$.

But let us first talk about our constant $\delta$. By \eqref{22.13},
$\frac{11}{10} \leq r_1 \leq \delta^{-1}$. On the other hand, Theorem 1.6 in \cite{Mono}
is stated with a single $r_1$, i.e., the small constant $\varepsilon >0$ in that statement depends
also $r_1$, which does not make us happy a priori. It is even noted after the statement that
in the present case, $\varepsilon$ depends on the ratio $\dist(L_0)^{-1} r_1$
(by dilation invariance). A later statement Corollary 9.3 in \cite{Mono}, solves this issue,
and gives a constant $\varepsilon$ that does not depend on $r_1$ as long as $r_1 \leq C$
(or here, $\delta^{-1}$), but the statement is a little more unpleasant because it also allows
more complicated choices of $L_0$ (that is, we are only interested in a line $L$ here,
and the mapping $\xi$ is an isometry), and also because the statement would rather concern
another dilation $\wt f(E)$, with a dilation factor $r^{-1}$, so that now $r$ becomes $1$
and $d_0 = \dist(0,L)$ becomes $r^{-1} d_0 \in [\delta, \frac{10}{11}]$.
The reader should not pay attention to the fact that the statement in \cite{Mono} requires
$\dist(0,L) \leq \frac{9}{10}$; the proof works the same way. We decided to simplify our lives,
and use Theorem 1.6 in \cite{Mono} with the knowledge that $\varepsilon$ does not depend
on $r_1$ as long as $r_1$ stays bounded.

So we check the assumptions. First, $E'$ is sliding minimal in $B(0,d_0^{-1}R)$, 
relative to $L_0$ and with the gauge function $h'(r) = h(d_0r)$.
We need to know that $r_1 \leq d_0^{-1}R$, or equivalently $r \leq R$,
and this is given by \eqref{22.12}. Also, $h'(r_1) = h(d_0 r_1) = h(r) \leq \varepsilon$
and (1.22) in \cite{Mono} holds. For (1.23) there, denote by $F'$ the functional associated to
$L_0$; then 
\begin{eqnarray} \label{22.17}
F'(r_1) &=& F(r) \leq \pi + \varepsilon \leq \theta_0 + \varepsilon
\leq e^{\alpha A(10^{-3}r)} \inf_{0 < \rho < 10^{-3}r} F(\rho)  + \varepsilon
\nn\\
&\leq& e^{\alpha \varepsilon} \inf_{0 < \rho < 10^{-3}r} F(\rho)  + \varepsilon
= e^{\alpha \varepsilon} \inf_{0 < \rho < 10^{-3}r_1} F'(\rho)  + \varepsilon
\end{eqnarray}
by the dilation invariance of densities, \eqref{22.14}, 
the fact that $\theta_0 = \lim_{\rho \to 0} \theta(\rho) = \lim_{\rho \to 0} F(\rho)$ is at least $\pi$,
the near monotonicity estimate \eqref{22.8}, the definition \eqref{22.7}, and \eqref{22.13}.
This gives the desired bound, if $\varepsilon$ is small enough. We can apply the
theorem, and we get a sliding minimal set $E_0$, with all sort of properties.
We want to check that $E_0$ coincides in $B(0,r_1)$ with the half plane bounded by $L_0$
that contains $0$ (or equivalently that $f^{-1}(E_0) = H_0$ in $B(0,r)$), and this way \eqref{22.15}
and \eqref{22.16} will follow from (1.25)-(1.27) in \cite{Mono}.

Now (1.24) in \cite{Mono} says that the analogue of $F'$ for $E_0$ takes a constant value
$D$ on $(0,r_1)$. Notice that $r_1 = d_0^{-1} r \geq \frac{11}{10}$ by \eqref{22.12}.
By (1.27) for $B(y,t) = B(0,1)$, we get that $D$ is as close to $\pi$ as we want. 
Now Theorem~1.3 in \cite{Mono} (about constant density) 
gives the following extra information on $E_0$. 

Set $A = B(0,r_1) \sm \{ 0 \}$, denote by $X$ the cone over $A\cap E_0$,
and by $S$ the shade of $L_0$. We get that $\H^2(A \cap E_0 \cap S) = 0$, that
$A \cap X \sm S \subset E_0$, that $X$ is a minimal cone (no boundary),
and $\H^d(S \cap B(0,r_1) \sm X) = 0$.

Notice that $B(0,1) \sm \{ 0 \} \subset A$ (because $r_1 \geq \frac{11}{10}$),
and that inside $B(0,1) \sm \{ 0 \}$, $E_0 \subset X$ by definition of $X$,
and $X = X \sm S \subset E_0$ because $B(0,1)$ does not meet $S$.
Then $H^{2}(X\cap B(0,1)) = H^{2}(E_0\cap B(0,1)) = D$, which is as close to $\pi$
as we want. Since $X$ is a minimal cone, $X$ is a plane. In addition, 
$\H^d(S \cap B(0,r_1) \sm X) = 0$ and $r_1 \geq \frac{11}{10}$, so 
$X$ contains a nontrivial bit of $S$, hence also the whole $L$. That is,
$X$ is the plane that contains $0$ and $L_0$.

Set $H = f(H_0) = \ol{X \sm S}$; we want to show that $E_0$ coincides with $H$ in
$B(0,r_1)$, or equivalently in $A =  B(0,r_1) \sm \{ 0 \}$ (because $E_0$ is closed).
We know that $A \cap X \sm S \subset E_0$, hence $A \cap H \subset E_0$
(again, $E_0$ is closed). Then $E_0 \cap A \subset X$ (by definition of $X$),
which means that $E_0 \cap A \sm H \subset S$. Since 
$\H^2(A \cap E_0 \cap S) = 0$, and $E_0$ is coral (or more brutally, locally Ahlfors regular),
we get that $E_0 \cap A \subset H$, as needed for Lemma \ref{t22.3}.
\qed

\begin{lem}\label{t22.4} 
Lemma \ref{t22.3} stays valid when instead of \eqref{22.14}, 
we require that the density of $E$ at $0$ is $\theta_0 = \frac{3\pi}{2}$ and that 
$F(r) \leq \frac{3\pi}{2} + \varepsilon$, 
and we get the same conclusion, except that $H_0$ is replaced with the set 
$E_0 = \ol{Y \sm S}$, where $Y$ is a minimal cone of type $\bY$, centered at $0$ 
and such that $L \cap B(0,r)$ is contained in a face of $Y$.
\end{lem}

Set $B = B(0,r)$. We only care about $E_0 \cap B$, because the outside part 
does not interfere with our description of $E$ in \eqref{22.15} and \eqref{22.16}, 
since $B(0,\frac{20 r}{21})$ lies well inside $B$. Inside $B$, $E_0 = \ol{Y \sm S}$ is really
a truncated set of type $\bY$, where we removed from $Y$ the part that lies on the other side of $L$,
of the face of $Y$ that contains $L\cap B$.

For the proof we proceed as for Lemma \ref{t22.3}. We can still apply Theorem 1.6 in \cite{Mono},
after applying the same composition $f$ of a rotation and a dilation by $d_0^{-1}$.
This theorem gives a sliding minimal set, which we now call $E'_0 \subset B(0,r_1)$, where
$r_1 = d_0^{-1} r \geq \frac{11}{10}$, with the additional property that
the analogue of $F$ takes a constant value $D$ on $(0,1)$, and which is 
very close to $f(E)$ in $B(0,r_1)$. In addition, $D$ is still as close as we want to 
the values of $F$ (computed with $E$ and for radii smaller than $r$), which are as 
close to $\theta_0 = \frac{3\pi}{2}$ as we want.

Then we turn to Theorem 1.3 of \cite{Mono} to get a good description of $E'_0$ in 
$A = B(0,r_1) \sm \{ 0 \}$.
We get the same basic properties as above, in terms of some minimal cone $X$, but now 
the density of $X$ is $D$, which is as close to $\frac{3\pi}{2}$ as we want.
Proposition 14.1 of \cite{Holder} gives a description of minimal cones of dimension $2$
that implies that this cannot happen unless $D = \frac{3\pi}{2}$, and hence $X$ is a cone of type $\bY$.

Let us now denote by $S'$ the shade of $L_0$.
We still have that $\H^d(S' \cap B(0,r_1) \sm X) = 0$, so $X$ contains $S' \cap B(0,r_1)$ 
because $X$ is closed. Notice also that $S' \cap B(0,r_1)$ is a nontrivial piece of plane,
because $r_1 \geq \frac{11}{10}$.

Next we check that $E'_0 \cap A = \ol{X \sm S'} \cap A$.
We know that $\H^2(A \cap E'_0 \cap S') = 0$, so each $x\in A \cap E'_0$ is the limit 
of a sequence $\{ x_j \}$ in $E'_0 \sm S'$ (recall that $E'_0$ is coral). 
Clearly $x_j \in A$ for $j$ large, hence $x_j \in X \sm S'$ (because $E'_0 \cap A \subset X$ 
by definition of $X$); thus $x\in \ol{X \sm S'} \cap A$.
Conversely, we know that $A \cap X \sm S' \subset E'_0$, hence 
$A \cap \ol{X \sm S'} \subset E'_0 \cap A$, and our claim follows. 
Both sets contain the origin, so $E'_0 \cap B(0,r_1) = \ol{X \sm S'} \cap B(0,r_1)$.

Set $Y = f^{-1}(X)$ and $E_0 = f^{-1}(E'_0)$. Then $Y$ is also a cone of type $\bY$,
and $E_0 \cap B(0,r) = f^{-1}(E'_0 \cap B(0,r_1)) = f^{-1}(\ol{X \sm S'} \cap B(0,r_1))=
\ol{Y \sm S} \cap B(0,r)$. Thus, inside $B(0,r)$, $E_0$ has the form that was announced in the lemma.
We do not care about what it is outside, because $\R^n \sm B(0,r)$ is far from $B(0, \frac{20r}{21})$
where we want to approximate $E$, as in \eqref{22.15} and \eqref{22.16}.
Finally, the good approximation of $E$ in $B(0, \frac{20r}{21})$ follows from the good approximation
of $E' = f(E)$ in  $B(0,(1-\tau)r_1)$ that is given by (1.24)-(1.27) of \cite{Mono}.
\qed

\ms
For radii $r$ that are much larger than $d_0$, it is easier to use compactness in another
way, and get a good approximation by a plane or a cone of type $\bH$ or $\bV$ centered on $L$.
Here is a statement, whose proof will rely on Lemma \ref{t22.2}.

\begin{lem}\label{t22.5} 
For each choice of small constant $\tau > 0$, we can find constants
$\varepsilon = \varepsilon(\tau) > 0$ and $\delta = \delta(\tau)$, 
that depend only on $\tau$ and $n$, with the following property. 
Let $E$ satisfy \eqref{23.1n} and \eqref{22.7}, and let $r$ be such that
\begin{equation}\label{22.18}
\delta^{-1} d_0 \leq  r  < R.
\end{equation}
Suppose in addition that $0 \in E \sm L$, 
\begin{equation}\label{22.19}
h(r) \leq \varepsilon \ \text{ and } \ \int_0^{r} h(t) \frac{dt}{t} \leq \varepsilon,
\end{equation}
and that there is $\theta_0 \in \{ \pi, \frac{3 \pi}{2} \}$ such that
\begin{equation}\label{22.20}
\lim_{\rho \to 0} F(\rho) = \theta_0  
\ \text{ and } \  
F(r) \leq \theta_0 + \varepsilon.
\end{equation}
Then there is a set $X_0 \in \bH(L) \cup \bV(L) \cup \bP_0$ such that
\begin{equation}\label{22.21}
d_{0,\frac{20 r}{21}}(E, X_0) \leq \tau
\end{equation}
and 
\begin{equation}\label{22.22}
|\H^2(E \cap B(y,t)) - \H^2(X_0 \cap B(y,t))| \leq \tau r^2
\end{equation}
for all $y \in \R^n$ and $t > 0$ such that $B(y,t) \subset B(0,\frac{20 r}{21})$.
If $\theta_0 = \pi$, then $X_0 \in \bH(L)$; if $\theta_0 = \frac{3\pi}{2}$, 
then $X_0 \in \bV(L) \cup \bP_0$.
\end{lem}

\ms
In this statement the planes through the origin (the elements of $\bP_0$) are some sort of a 
stowaway (or party crasher); the proof allows them, but we expect to get rid of them later. 
That is, if we get \eqref{22.21} and \eqref{22.22} for a plane $X_0$
that does not nearly contain $L$ (i.e., the two unit vectors of $L$ are far from $X_0$), then
we shall be able to show that $E$ is smooth near $0$, and $\theta_0 = \pi \neq \frac{3\pi}{2}$,
a contradiction. See Theorem \ref{t28.2} and Remark \ref{r28.3} for another instance of this 
reasoning, where we need to look at different scales to exclude apparently acceptable behaviors.

The proof is a standard compactness argument, similar to what was done for the
proof of Theorem 1.6 in \cite{Mono}. Suppose we can find $\tau > 0$ such that
taking $\varepsilon = \delta = 2^{-k}$ never works. Let $E_k$, $h_k$, $L_k$, $r_k$, etc.
provide a counterexample. By scale and rotation invariance, we may assume that $r_k = 1$ for all $k$,
and that we can find orthogonal unit vectors $e_1$ and $e_2$ such that
$L_k = \big\{d_k e_1+te_2 \, ; \, t\in \R\big\}$, and with positive numbers
$d_k = \dist(0,L_k)$ that tend to $0$ (by \eqref{22.18} and because $\delta_k$
tends to $0$).

We replace $\{ E_k \}$ with a subsequence which has a limit $E_\infty$.
Let $L_\infty$ denote the limit of the $L_k$; this is a line through the origin.
Also consider $E'_k = E_k - d_k e_1$; this is a sliding minimal set, with sliding boundary
$L_k - d_k e_1 = L_{\infty}$, and $E'_k$ also tends to $E_\infty$.
Notice that the gauge functions $h_k$ satisfy \eqref{22.19} uniformly on $(0,1)$, 
and also tend to $0$ uniformly on $(0,1)$. By Theorem 10.8 in \cite{Sliding}, 
 $E_\infty$ is a sliding minimal set in $B(0,1)$, associated to $L_\infty$
 (and the gauge function $h=0$). Next we check that
 \begin{equation} \label{22.23}
\H^{2}(E_\infty \cap B(0,\rho)) = \theta_0 \rho^2
\ \text{ for } 0 < \rho < 1.
\end{equation}
In fact, let $B = B(y,t)$ be given, with $|y| + t < 1$; we first apply the lower semicontinuity 
result in \cite{Sliding} (namely, Theorem 10.97 there) to the same sets $E'_k$, with 
the same assumptions, and get that
\begin{eqnarray} \label{22.24}
\H^{2}(E_\infty \cap B) &\leq& \liminf_{k \to +\infty} \H^{2}(E'_k \cap B)
= \liminf_{k \to +\infty} \H^{2}(E_k \cap B(y+d_k e_1,t))
\nn\\
&\leq& \liminf_{k \to +\infty} \H^{2}(E_k \cap B(y,t+d_k)).
\end{eqnarray}
For the upper semicontinuity, we call Lemma 22.3 in \cite{Sliding}, which we can apply
with $M$ as close as we want to $1$ and $h$ as small as we want, and we get that
for the compact set $\ol B$,
\begin{eqnarray} \label{22.25}
\H^{2}(E_\infty \cap \ol B) &\geq& \limsup_{k \to +\infty} \H^{2}(E'_k \cap \ol B)
= \limsup_{k \to +\infty} \H^{2}(E_k \cap \ol B(y+d_k e_1,t))
\nn\\
&\geq& \limsup_{k \to +\infty} \H^{2}(E_k \cap \ol B(y,t-d_k)).
\end{eqnarray}
Let us apply this with $y=0$; notice that if $S_k$ denotes the shade of $L_k$, then 
\begin{equation} \label{22.26}
\lim_{k \to +\infty} \H^2(S_k \cap B(0,t+d_k)) = \frac{\pi t^2}{2}
\end{equation}
because $d_k$ tends to $0$ and $L_k$ tends to $L_\infty$. Thus \eqref{22.24} implies that
\begin{equation} \label{22.27}
\H^{2}(E_\infty \cap B(0,t)) \leq \liminf_{k \to +\infty} \H^{2}(E_k \cap B(y,t+d_k))
= - \frac{\pi t^2}{2} +  \liminf_{k \to +\infty}\big[(t+d_k)^{2} F_k(t+d_k)\big].
\end{equation}
For $k$ large, $t+d_k < 1$, hence by \eqref{22.8}
\begin{equation} \label{22.28}
F_k(t+d_k) \leq e^{\alpha A_k(1)} F_k(1)\leq e^{\alpha 2^{-k}}F_k(1)
\leq e^{\alpha 2^{-k}}[\theta_0 + 2^{-k}]
\end{equation}
because \eqref{22.19} holds with $\varepsilon = 2^{-k}$, and then by \eqref{22.20}. 
The right-hand side tends to $\theta_0$, hence by \eqref{22.27}
\begin{equation} \label{22.29}
\H^{2}(E_\infty \cap B(0,t)) \leq t^2 \big[\theta_0 - \frac{\pi}{2} \big].
\end{equation}
Conversely, \eqref{22.25} yields
\begin{equation} \label{22.30}
\H^{2}(E_\infty \cap \ol B(0,t)) \geq \limsup_{k \to +\infty} \H^{2}(E_k \cap \ol B(0,t-d_k))
\geq - \frac{\pi t^2}{2} +  \limsup_{k \to +\infty}\big[(t+d_k)^{2} F_k(t-d_k)\big]
\end{equation}
and, since by \eqref{22.8} and \eqref{22.20}
\begin{equation} \label{22.31}
F_k(t-d_k) \geq e^{-\alpha A_k(1)} \lim_{\rho \to 0} F_k(\rho)
\geq e^{-\alpha 2^{-k}}\lim_{\rho \to 0} F_k(\rho)
= e^{-\alpha 2^{-k}}\theta_0,
\end{equation}
which tends to $\theta_0$, we get that
\begin{equation} \label{22.32}
\H^{2}(E_\infty \cap \ol B(0,t)) \geq t^2 \big[\theta_0 - \frac{\pi}{2} \big].
\end{equation}
It follows that 
\begin{equation} \label{22.33}
t^{-2} \H^{2}(E_\infty \cap B(0,t)) = \theta_0 - \frac{\pi}{2}
\ \text{ for } 0 < t < 1,
\end{equation}
i.e., $E_\infty$ has constant density equal to $\theta_0 - \frac{\pi}{2}$ on $(0,1)$.
By the constant density result (Theorem 29.1) in \cite{Sliding}, 
$E_\infty$ coincides with a sliding minimal cone in $B(0,1)$.
Call this cone $X_0$; by Lemma \ref{t22.2}, $X_0 \in \bH(L)$
if $\theta_0 = \pi$ and $X_0 \in \bV(L) \cup \bP_0$ if $\theta_0 = \frac{3\pi}{2}$.

Let us now check that \eqref{22.21} and \eqref{22.22} hold for $k$ large;
this will give the desired contradiction with the definition of $E_k$
and complete the proof of Lemma \ref{t22.5}.
Now \eqref{22.21} holds because $X_0$ is the same as $E_\infty$ in $B(0,1)$,
we normalized things so that $r_k = 1$, and $E_\infty$ is the limit of $E_k$
locally in $B(0,1)$. For a given ball $B = B(y,t)$, notice that for $0 < t_1 < t < t_2$,
with $B(y,t_2) \subset B(0,1)$, \eqref{22.24} and \eqref{22.25} yield
\begin{eqnarray} \label{22.34}
\H^{2}(X_0 \cap B(y,t)) &=& \H^{2}(E_\infty \cap B(y,t))
\leq \liminf_{k \to +\infty} \H^{2}(E_k \cap B(y,t+d_k))
\nn\\
&\leq&
\liminf_{k \to +\infty} \H^{2}(E_k \cap B(y,t_2))
\end{eqnarray}
and similarly
\begin{eqnarray} \label{22.35}
\H^{2}(X_0 \cap B(y,t)) &=& \H^{2}(E_\infty \cap B(y,t))
\geq \limsup_{k \to +\infty} \H^{2}(E_k \cap B(y,t-d_k))
\nn\\
&\geq&
\limsup_{k \to +\infty} \H^{2}(E_k \cap B(y,t_1)).
\end{eqnarray}
From this it is easy so deduce that for a fixed $B(y,t)$, the estimates in \eqref{22.22} 
hold for $k$ large. But we do not want to let $k$ depend on $y$ and $t$, so a little
bit of uniformity is needed to conclude. This is rather easily done, because we control $X_0$
well; we refer to Lemma 9.2 in \cite{Mono} for the proof.
Thus we get the desired contradiction, and Lemma \ref{t22.5} follows.
\qed

\section{Statements of decay for $F$; differential inequalities}
\label{S23}

Recall that we want to generalize the work of Sections \ref{S3}-\ref{S21},
with balls that are no longer centered on $L$, and we decided to replace the usual density $\theta(r)$
with the functional $F(r)$ defined in \eqref{22.3}.
In this section we give the main decay statement for $F$. 
Recall that $F$ is almost nondecreasing; we intend to say that in some
circumstances, it actually decays at some speed, but we shall only be able to do this
when $E$ is close enough to some special minimal sets.

The assumptions for this section and the next ones are the following. 
We still work in $\R^n$, with a line $L$ and a sliding almost minimal set $E$ that satisfies 
\eqref{23.1n} and \eqref{23.2n}; we also assume that 
\begin{equation} \label{23.1}
0 < d_0 := \dist(0,L) < \frac{2R}{3},
\end{equation}
and in the statements additional conditions on the size of $C_h$ in \eqref{23.2n} will arise.

Denote by $\bH = \bH(L)$ the set of half planes bounded by $L$, and by 
$\bV = \bV(L)$ the collection of sets of type $\bV$ bounded by $L$, i.e., unions of
two half planes of $\bH$ that make an angle at least $\frac{2\pi}{3}$ with each other along $L$. 
This includes planes that contain $L$. Still let $\bP_0$ denote the collection of all planes through 
the origin. We will often require $E$ to be close to sets of $\bH \cup \bV \cup \bP_0$, 
and we measure this with the quantities
\begin{equation} \label{23.4}
\beta_H(r) = \inf_{H \in \bH} d_{0,r}(E,H)
\ \text{ and } \ 
\beta_{VP}(r) = \inf_{V \in \bV \cup \bP_0} d_{0,r}(E,V),
\end{equation}
where we will naturally restrict to $r \in (0,R]$. 

Let us give two parallel statements, which will be proved afterwards.
We start with the case when there is a good approximation by a half plane.

\begin{thm} \label{t23.1}
There exist constants $a \in (0,10^{-1})$, $\varepsilon_H > 0$, 
and $C_H \geq 1$, that depend only on $n$ and $\beta$, with the following properties.
Let $L$, $E$, and $h$ satisfy \eqref{23.1}, \eqref{23.1n}, and \eqref{23.2n}, 
with a constant $C_h$ such that
\begin{equation}\label{23.5}
C_h R^{\beta} \leq \varepsilon_H.
\end{equation}
Suppose also that $0 \in E$, and that 
\begin{equation} \label{23.6}
\beta_H(R) \leq \varepsilon_H  \ \text{ or } F(R) -\pi \leq \varepsilon_H.
\end{equation}
Then 
\begin{equation} \label{23.7}
F(r_1) -\pi \leq \Big(\frac{2 r_1}{r_2}\Big)^{a} [F(r_2) -\pi] + C_H C_h r_1^a r_2^{\beta-a}
\end{equation}
for $0 \leq r_1 \leq r_2 \leq 9R/20$. 
\end{thm}

\ms
See \eqref{22.3} for the definition of $F$. 
Within minor modifications, this is the same statement as Theorem \ref{t1.10a} in the introduction.
Notice that because of \eqref{23.2n}, 
the limit density $\theta_0 = \lim_{r\to 0} r^{-2}\H^2(E\cap B(x,r))$
exists (by \eqref{22.11}); we know that for $0 \in E$, this limit cannot be smaller than $\pi$,
and in the present situation, we will see during the proof that $\theta_0 = \pi$. 
Or just notice that you get this when you let $r_1$ tend to $0$ in \eqref{23.7}.
We will also check that the two possible assumptions in \eqref{23.6} imply each other, 
modulo changing the small constant and replacing $R$ with a slightly smaller radius.

It happens that the good decay provided by \eqref{23.7} implies a polynomial control
on $\beta_H(r)$ for $r$ small; see Section \ref{S28}. % OK so far

We have a similar statement for the case when $E$ is well approximated by a
set of $\bV$; this time the relevant value of density is $\theta_0 = \frac{3\pi}{2}$.

\begin{thm} \label{t23.2}
There exist constants $a \in (0,10^{-1})$, $\varepsilon_V > 0$, 
and $C_V \geq 1$, that depend only on $n$ and $\beta$, with the following properties.
Let $L$, $E$, and $h$ satisfy \eqref{23.1}, \eqref{23.1n}, and \eqref{23.2n},  
with a constant $C_h$ such that
\begin{equation}\label{23.8}
C_h R^{\beta} \leq \varepsilon_V.
\end{equation}
Suppose also that $0 \in E_0$, 
\begin{equation}\label{23.9}
\lim_{r\to 0} r^{-2}\H^2(E\cap B(x,r)) = \frac{3 \pi}{2},
\end{equation}
and that
\begin{equation} \label{23.10}
\beta_{VP}(R) \leq \varepsilon_V.  
\end{equation}
Then 
\begin{equation} \label{23.11}
F(r_1) -\frac{3 \pi}{2} \leq \Big(\frac{C_V r_1}{r_2}\Big)^{a} [F(r_2) -\frac{3 \pi}{2}] 
+ C_V C_h r_1^a r_2^{\beta-a}
\end{equation}
for $0 \leq r_1 \leq r_2 \leq R/2$.
\end{thm}

\ms
This time see Theorem \ref{t1.11a} in the introduction.
The same sort of remarks as above apply to this case. Notice the additional constant $C_V^a$
in \eqref{23.11}, which is due to a gap in the set of radii $r$ for which the main differential inequality
described below holds. This could probably be improved, but the additional constant does not 
disturb much.

We did not include the option that $F(R) - \frac{3 \pi}{2} \leq \varepsilon_V$ instead of
\eqref{23.10}, because it does not imply that $E$ is close to a set of type $\bV$ or a plane. 
The difference will not be enormous in the end; we will see in Lemma \ref{t24.2} 
that if $F(R) - \frac{3 \pi}{2} \leq \varepsilon_V$ and $d_0$ is much smaller than $r$,
then $\beta_{VP}(r) \leq \varepsilon$, and we can apply Theorem \ref{t23.2}. 
When instead $d_0$ is not so small compared to $r$ (and $r \leq R/2$, say), 
Lemma \ref{t24.3} will say that we can find a truncated $\bY$-set centered
at $0$ that approximates $E$ well in $B(x,r)$. As hinted above, this set is not 
close to a $\bV$-set because it is centered at $0$ (and $d_0$ is not so small). 
We could try to show that $F(r)$ decays also in this intermediate region, 
but instead we will just use the fact that $F$ is almost nondecreasing there (by \cite{Mono}), 
and this will be fine because the concerned set of radii $r$ is not so large anyway.

See Section \ref{S28} 
for the control of the geometry of $E$ that follows from \eqref{23.11}.

In both statements, the interesting part of the conclusion is when $r_1$ gets much smaller
than $r_2$; otherwise a direct application of \eqref{22.8} gives at least as much.
In both cases the main ingredient in the proof is a differential inequality which we state now.

\begin{pro} \label{t23.3}
There exist constants $a \in (0,10^{-1})$, $\varepsilon_1 > 0$, 
and $C_1 \geq 1$, that depend only on $n$ and $\beta$, with the following properties.
Let $E$ and $h$ satisfy \eqref{23.1}, \eqref{23.1n}, and \eqref{23.2n}, and suppose that $0 \in E$. 
% needed? yes , density needs to be large, right?
For almost every $r$ such that
\begin{equation} \label{23.12}
2d_0 \leq r \leq \frac{R}{2},
\end{equation}
\begin{equation}\label{23.13}
C_h r^{\beta} \leq \varepsilon_1
\end{equation}
and 
\begin{equation} \label{23.14}
\beta_H(2r) \leq \varepsilon_1,
\end{equation}
the function $F$ of \eqref{22.3} is differentiable at $r$, and
\begin{equation} \label{23.15}
r F'(r) \geq a [F(r) - \pi]_+ - C_1 \int_0^{2r} \frac{h(t) dt}{t}.
\end{equation}
\end{pro} % la j'ai mis la partie positive. je mets un commentaire

\ms
This will be proved in Section \ref{S25}. 
We do it on purpose to mention $h$ explicitly in \eqref{23.15}, rather than the estimate
that we could get from \eqref{23.2n}, because we may sometimes get an estimate that is
better than expected. Even though $\varepsilon_1$ needs to be quite small, we think of it as
being roughly constant, while we hope that $F(r) - \pi$, for instance, will become really small.

We took the positive part of $F(r) - \pi$ not to get confused by the case when 
$F(r) - \pi < 0$, in which case \eqref{23.15} is actually better when $a$ is smaller.
This way, at least, our estimate is better when we can take $a$ larger.
However, we will pay a (moderate) price for this simplification, when we prove \eqref{23.15}. 
We could also have used the same sort of computations as in \cite{C1} and Proposition \ref{t16.2}.
This way the reader gets to choose their prefered method.

The next statement is similar, but concerns the approximation with sets of type $\bV \cup \bP_0$ and 
the larger reference density $2\pi/3$. It is a little more complicated for the same reasons 
as for Theorem \ref{t23.2}; it will be proved in Sections \ref{S25}-\ref{S27}.

\begin{pro} \label{t23.4}
There exist constants $a \in (0,10^{-1})$, $N \geq 1$, $\varepsilon_2 > 0$, 
and $C_2 \geq 1$, that depend only on $n$ and $\beta$, with the following properties.
Let $E$ satisfy \eqref{23.1}, \eqref{23.1n}, and \eqref{23.2n}, and suppose that
\begin{equation} \label{23.16}
\lim_{\rho \to 0} \rho^{-2} \H^2(E \cap B(0,\rho)) = \frac{3\pi}{2}.
\end{equation}
For almost every $r$ such that
\begin{equation} \label{23.17}
Nd_0 \leq r \leq \frac{R}{2},
\end{equation}
\begin{equation}\label{23.18}
C_h r^{\beta} \leq \varepsilon_2,
\end{equation}
and 
\begin{equation} \label{23.19}
\beta_{VP}(2r) \leq \varepsilon_2,
\end{equation}
the function $F$ of \eqref{22.3} is differentiable at $r$, and
\begin{equation} \label{23.20}
r F'(r) \geq a \big[F(r) - \frac{3\pi}{2} \big]_+ - C_2 \int_0^{2r} \frac{h(t) dt}{t}.
\end{equation}
\end{pro}

\section{How to deduce decay from differential inequalities}
\label{S24}

In this section we see how to deduce the decay estimates, Theorems \ref{t23.1} and \ref{t23.2},
from the corresponding differential inequalities, Propositions \ref{t23.3} and \ref{t23.4}.
Most of it will consist in checking that the main geometric assumption \eqref{23.14} or \eqref{23.19}
is valid. 

Throughout this section, we assume that the main assumptions of Section \ref{S23} are valid,
i.e., that $L$, $E$, and $h$ satisfy \eqref{23.1}, \eqref{23.1n}, and \eqref{23.2n}.
By \eqref{23.2n} and its consequence \eqref{22.11}, the density
\begin{equation}\label{24.1}
\theta_0 = \lim_{r \to 0} \theta(r) = \lim_{r \to 0} F(r) 
\end{equation}
exists; we shall either assume or prove that $\theta_0 \in \big\{ \pi, \frac{3\pi}{2} \big\}$.
We first check that the conditions of \eqref{23.6} essentially imply each other, and that
\eqref{23.10} implies that $F(r) - \frac{3 \pi}{2}$ is small.

\begin{lem}\label{t24.1} 
For each small $\varepsilon > 0$, there exist $\varepsilon_H > 0$ and $\varepsilon_V > 0$, 
that depend only on $n$ and $\beta$, such that if the assumptions of Theorem \ref{t23.1} 
are satisfied, then
\begin{equation} \label{24.2}
\beta_H(9R/10) \leq \varepsilon  
\end{equation}
and 
\begin{equation} \label{24.3}
F(r) \leq  \pi + \varepsilon  \ \text{ for } \  0 < r \leq 9R/10
\end{equation}
and if the assumptions of Theorem \ref{t23.2} are satisfied, then
\begin{equation} \label{24.4}
F(r) \leq \frac{3 \pi}{2} + \varepsilon \ \text{ for } 0 < r \leq 9R/10.
\end{equation}
\end{lem}

\ms
First assume that $E$ is as in Theorem \ref{t23.1}, with $\beta_H(R) \leq \varepsilon_H$.
Then $\beta_H(9R/10) \leq 10\varepsilon_H/9$, by the definition \eqref{23.4}, so we 
just need to show that \eqref{24.3} holds if $\varepsilon_H$ is small enough. 
Let us first check that
\begin{equation} \label{24.5}
F(9R/10) \leq \pi + \varepsilon/2.
\end{equation}
Let us proceed by compactness. If this fails, then for each large integer $k$, 
we can find $L_k$, $E_k$, $h_k$, $R_k$, as in Theorem \ref{t23.1} with $\varepsilon_H = 2^{-k}$, 
but for which $F_k(9R_k/10) > \pi + \varepsilon/2$, where 
\begin{equation} \label{24.6}
F_k(r) = r^{-2} [\H^2(E_k \cap B(0,r)) + \H^2(S_k \cap B(0,r))],
\end{equation}
and $S_k$ is the shade of $L_k$.
We want to take a limit, but first we use the dilation invariance of our problem to
assume that $R_k = 1$ for all $k$.  Also, choose two unit vectors $e_1$ and $e_2 \perp e_1$;
we can use the rotation invariance to ensure that the following two properties hold for each $k$. 
First let $z_k$ denote the point of $L_k$ that lies closest to $0$. Notice that $z_k \neq 0$ 
(because $0 \notin L_k$ by \eqref{23.1}); we require that $z_k/|z_k| = e_1$. And also that 
$e_2$ is parallel to $L_k$.

Next set $d_k = \dist(0, E_k)$; recall from \eqref{23.1} that $d_k < R_k = 1$, 
so we we may assume (at the price of replacing our sequence by a subsequence, 
to which we automatically give the same name) that $d_k$ has a limit $d_\infty \in [0,1]$.
Then $L_k$ converges to the limit $L_\infty = d_\infty e_1 + \R e_2$.

Since $\beta_H(R) \leq 2^{-k}$ (for the set $E_k$), there is a half plane $H_k$
bounded by $L_k$, such that 
\begin{equation}\label{24.7}
d_{0,1}(E_k, H_k) \leq 2^{-k}.
\end{equation}

We extract a new subsequence, so that after extraction $H_k$ converges 
(say, for the Hausdorff distance in $\ol B(0,2)$) to a half plane $H_\infty$ 
bounded by $L_\infty$. We allow the case when $d_\infty = 0$, 
but notice that $\dist(0,H_k) \leq 2^{-k}$, by \eqref{24.7} and because $0 \in E$.
Thus $H_\infty$ contains the origin.

Extract again a subsequence, so that $\{ E_k \}$ converges, locally in $B(0,1)$, to a closed
set $E_\infty$. In fact, in the present situation this is not even needed, because of \eqref{24.7}, but for
the next lemma it will feel better, and anyway this is costless.
We want to apply a theorem about limits to the sequence $\{ E'_k \}$, where 
$E'_k = E_k + (d_\infty -d_k) e_1$. Since $d_\infty -d_k$ tends to $0$,
$\{ E'_k \}$ also converges to $E_\infty$, but the point of the translation is
that $E'_k$ is sliding minimal, in a domain $B_k = B((d_\infty -d_k) e_1,1)$ that tends to $B(0,1)$,
with a same boundary set $L_k + (d_\infty -d_k) e_1 = L_\infty$. This way we can
apply theorems of convergence with a fixed boundary set.

We put ourselves in $B= B(0,\frac{99}{100})$, which is contained in $B_k$ for $k$ large. 
Thus $E'_k$ is almost minimal in $B$, relative to the boundary $L_\infty$, and with 
a gauge function $h'_k$ such that $h'_k(1) = h_k(R_k) \leq 2^{-k}$, 
by \eqref{23.2n} and \eqref{23.5} with $\varepsilon_H = 2^{-k}$. 

By \eqref{24.7}, both $E_k$ and $E'_k$ converge to $H_\infty$, locally in $B(0,1)$,
i.e., for any localized Hausdorff distance function $d_{0,r}$, $0 < r < 1$.

Let us fix $r \in (\frac{9}{10},\frac{98}{100})$, and apply Theorem 22.1 in \cite{Sliding} to the 
sequence $\{ E'_k \}$ and the compact set $\ol B(0,r)$. Notice in particular that 
the minimizing sequence property (21.14) in \cite{Sliding} is satisfied, with $\delta = 1$, 
and where $k_0$ is simply chosen so that $h_k(1) < \varepsilon$ for $k \geq k_0$. We find that
\begin{equation} \label{24.8}
\H^2(H_\infty \cap B(0,r)) = 
\H^2(H_\infty \cap \ol B(0,r)) 
\geq \limsup_{k \to +\infty} \H^2(E'_k \cap \ol B(0,r)).
\end{equation}
But for $k$ large, 
\begin{equation} \label{24.9}
\begin{aligned}
\H^2(E'_k \cap \ol B(0,r)) &= \H^2(E_k \cap \ol B((d_k-d_\infty) e_1 ,r))
\cr&\geq \H^2(E_k \cap \ol B(0,r-|d_\infty -d_k|)) \geq \H^2(E_k \cap B(0,\frac{9}{10})),
\end{aligned}
\end{equation}
because $r > \frac{9}{10}$. Hence
\begin{equation} \label{24.10}
\begin{aligned}
\Big(\frac{9}{10}\Big)^2 \limsup_{k \to +\infty} F_k(\frac{9}{10}) &= 
\limsup_{k \to +\infty} \Big[\H^2(E_k \cap B(0,\frac{9}{10})) + \H^2(S_k \cap B(0,\frac{9}{10})) \Big]
\cr& \leq \H^2(H_\infty \cap B(0,r)) +  \limsup_{k \to +\infty} \H^2(S_k \cap B(0,\frac{9}{10})).
\end{aligned}
\end{equation}
If $d_\infty = 0$, then $H_\infty$ is a half plane bounded by a line $L_\infty$ that contains
the origin, and the right-hand side of \eqref{24.10} is $\frac{\pi r^2}{2} + \frac{\pi (9/10)^2}{2}
\leq \pi r^2$. Otherwise, we know that $0 \in H_\infty \sm L_\infty$, thus  
$H_\infty$ is the half plane bounded by $L_\infty$ and that contains the origin.
At the same time, $S_k$ converges nicely to the closure of the complement of $H_\infty$
in the plane that contains it. Thus the right-hand side of \eqref{24.10} is also smaller than $\leq \pi r^2$
in this case. We put things together and get that
\begin{equation}\label{24.11}
\limsup_{k \to +\infty} F_k\Big(\frac{9}{10} \Big) \leq \Big(\frac{10 r}{9}\Big)^2 \pi .
\end{equation}
We take $r > \frac{9}{10}$ so close to $\frac{9}{10}$ that the right-hand side is
smaller than $\pi + \varepsilon/2$, and get a contradiction with the fact that $E_k$ was chosen 
so that $F_k(9R_k/10) > \pi + \varepsilon/2$. This concludes our proof of \eqref{24.5}.

We shall now easily deduce \eqref{24.3} from \eqref{24.5} and the near monotonicity 
formula \eqref{22.8}. Let us first recall that if $A$ is as in \eqref{22.7}, then
\begin{equation}\label{24.12}
A(R) = \int_0^R h(t) \frac{dt}{t} \leq C_h \int_0^R t^{\beta -1} dt = \beta^{-1} C_h R^\beta 
\leq \beta^{-1} \varepsilon_H
\end{equation}
by \eqref{23.2n} and because \eqref{23.5} holds. Then for $0 < r \leq 9R/10$, \eqref{22.8}
yields
\begin{equation}\label{24.13}
\begin{aligned}
F(r) &\leq  e^{\alpha A(r)} F(r) \leq  e^{\alpha A(9R/10)} F(9R/10)
\cr&
\leq e^{\alpha\beta^{-1} \varepsilon_H} F(9R/10) 
\leq e^{\alpha\beta^{-1} \varepsilon_H}(\pi + \varepsilon/2)
< \pi + \varepsilon
\end{aligned}
\end{equation}
by \eqref{24.12} and \eqref{24.5} and if $\varepsilon_H$ is small enough. Thus \eqref{24.3} holds.

\ms
Next we assume that $E$ is as in Theorem \ref{t23.1} and that $F(R) \leq \pi + \varepsilon_H$,
and we prove \eqref{24.2} and \eqref{24.3}. We start with \eqref{24.3}. Observe that for
$0 < r < R' < R$, 
\begin{equation}\label{24.14}
F(r) \leq  e^{\alpha A(r)} F(r) \leq  e^{\alpha A(R')} F(R')
\leq e^{\alpha\beta^{-1} \varepsilon_H} F(R') 
\end{equation}
by \eqref{22.8} and \eqref{24.12}. We let $R'$ tend to $R$ in \eqref{24.14} and get that
for $0 < r < R$,
\begin{equation}\label{24.15}
F(r) \leq e^{\alpha\beta^{-1} \varepsilon_H} F(R) 
\leq e^{\alpha\beta^{-1} \varepsilon_H} (\pi + \varepsilon_H) \leq \pi+\varepsilon
\end{equation}
if $\varepsilon_H$ is small enough. This proves \eqref{24.3}, and we are left with \eqref{24.2} to prove.

Let us first try to apply Lemma \ref{t22.5}, to the radius $r = \frac{21}{20} \, \frac{9R}{10}$,
with $\tau = \varepsilon$ and $\theta_0 = \pi$.
If we can do this, \eqref{22.21} says that $d_{0,\frac{20r}{21}}(E,X_0) \leq \varepsilon$
for some $X_0 \in \bH(L)$, and this yields \eqref{24.2}. So we just check the assumptions.
First, \eqref{22.19} follows from \eqref{23.2n} 
and \eqref{23.5} (if $\varepsilon_H$ is small enough).
The second half of \eqref{22.20} (the upper bound for $F$) follows from \eqref{24.15}, which also
implies (when you let the radius in \eqref{24.15} tend to $0$) 
that $\theta_0 < \pi + \varepsilon < \frac{3\pi}{2}$.
This implies that $\theta_0 = \pi$ (there is no other possible value, since $0 \in E$), and so \eqref{22.20}
holds. The second half of \eqref{22.18} is satisfied too, so we can apply the lemma and get the desired conclusion
\eqref{24.2} as soon as $d_0 \leq \delta r$, where $\delta = \delta(\varepsilon)$ is the small constant attached
by the lemma to our choice of $\tau = \varepsilon$. 

Thus we may assume that $d_0 \geq \delta r$,
and now we want to apply Lemma \ref{t22.3} to the same radius 
$r = \frac{21}{20}\,\frac{9R}{10}$ as before, with the constant $\delta$ that we just found,
and with $\tau = \varepsilon$. As before, the assumptions \eqref{22.13} and \eqref{22.14} are satisfied
if $\varepsilon_H$ is small enough, $0 \in E$, and the first part of the remaining assumption \eqref{22.12} 
is satisfied. So we can apply the lemma if $\frac{11 d_0}{10} \leq r = \frac{21}{20}\,\frac{9R}{10}$.
This is the case, because we required in \eqref{23.1} that $d_0 \leq 2R/3$.
So we get \eqref{22.15} for some half plane $H_0 \in \bH$ (in fact the one that contains the origin), 
and this implies \eqref{24.2} as before.

\ms
This completes our verification in the two cases that belong to Theorem \ref{t23.1}.
Now assume that $E$ is as in Theorem \ref{t23.2}. In particular, \eqref{23.10}
says that $\beta_{VP}(R) \leq \varepsilon_V$. If we prove that
\begin{equation}\label{24.16}
F(9R/10) \leq \frac{3\pi}{2} + \frac{\varepsilon}{2},
\end{equation}
then \eqref{24.4} will follow at once, by the proof of \eqref{24.13}.

We prove \eqref{24.16} with the same compactness argument as for \eqref{24.5}.
This time \eqref{23.10} yields the analogue of \eqref{24.7}, but for a set $V_k \in \bV(L_k)$
instead of $H_k$. We may still take a subsequence so that $V_k$ converges to a limit $V_\infty$,
and $V_\infty \in \bV(L_\infty)$. As before, $\dist(0,V_k) \leq 2^{-k}$, hence $0 \in V_\infty$.

We can keep the limiting argument (with the sequence $\{ E'_k \}$) as it was, and we get that for
$r \in (\frac{9}{10},\frac{98}{100})$, 
\begin{equation} \label{24.17}
\H^2(V_\infty \cap B(0,r)) \geq \limsup_{k \to +\infty} \H^2(E'_k \cap \ol B(0,r))
\geq \limsup_{k \to +\infty} \H^2(E_k \cap B(0,\frac{9}{10}))
\end{equation}
as in \eqref{24.8} and \eqref{24.9}, and 
\begin{equation} \label{24.18}
\begin{aligned}
\Big(\frac{9}{10}\Big)^2 \limsup_{k \to +\infty} F_k(\frac{9}{10}) &= 
\limsup_{k \to +\infty} \Big[\H^2(E_k \cap B(0,\frac{9}{10})) + \H^2(S_k \cap B(0,\frac{9}{10})) \Big]
\cr& \leq \H^2(V_\infty \cap B(0,r)) +  \limsup_{k \to +\infty} \H^2(S_k \cap B(0,\frac{9}{10}))
\end{aligned}
\end{equation}
as in \eqref{24.10}.  We start with the simpler case when $L_\infty$
goes through the origin. Then $\H^2(V_\infty \cap B(0,r)) = \pi r^2$,
$\lim_{k \to +\infty} \H^2(S_k \cap B(0,\frac{9}{10})) = \frac{\pi}{2} (9/10)^2 
\leq \frac{\pi r^2}{2}$, and the right-hand side of \eqref{24.18} is less than $\frac{3\pi r^2}{2}$.
If we take $r$ close enough to $9/10$, we get that
$\limsup_{k \to +\infty} F_k(\frac{9}{10}) \leq \frac{3\pi}{2} + \frac{\varepsilon}{3}$,
and for $k$ large this contradicts the fact that $E_k$ was chosen to violate \eqref{24.16}.

\begin{figure}[!h]  
\centering
\hskip0.7cm\includegraphics[width=9cm]{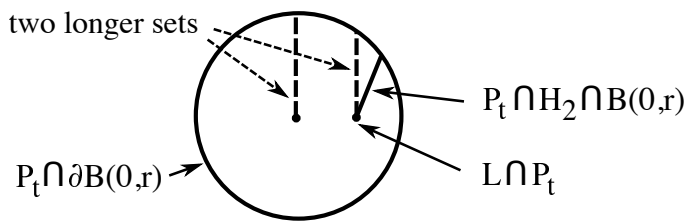}
\caption{A picture in $P_t$ 
\label{f24.1}}
\end{figure}

So we may assume that $d_\infty > 0$. Write $V_\infty = H_1 \cup H_2$, with $H_i \in \bH(L_\infty)$.
Since $0 \in V_\infty$, we may assume that $0 \in H_1$. At the same time,
$S_k$ tends nicely to the shade $S_\infty$ of $L_\infty$, which is just opposite to $H_1$. Thus,
if $A$ denotes the right-hand side of \eqref{24.18}, 
\begin{equation} \label{24.19}
\begin{aligned}
A &= \H^2(V_\infty \cap B(0,r)) + \H^2(S_\infty \cap B(0,\frac{9}{10}))
\cr&\leq \H^2(V_\infty \cap B(0,r)) + \H^2(S_\infty \cap B(0,r))
\cr& = \H^2(H_1 \cap B(0,r)) +  \H^2(H_2 \cap B(0,r)) + \H^2(S_\infty \cap B(0,r))
\cr& = \pi r^2 + \H^2(H_2 \cap B(0,r)).
\end{aligned}
\end{equation}
Now $H_2$ makes an angle at least $2\pi/3$ with $H_1$. One way to compute $\H^2(H_2 \cap B(0,r))$
consists in slicing it by planes. That is, write $L_\infty = d_\infty e_1 + \R e_2$, where $e_1$ and $e_2$
are orthogonal unit vectors, and let $e_3$ be a third unit vector, orthogonal to $e_1$ and $e_2$,  
such that $H_2$ is contained in the $3$-space generated by $e_1$, $e_2$, and $e_3$. Set 
$P_t = \big\{x e_1 + t e_2 + y e_3  \, ; \, (x,y) \in \R^2 \big\}$ for $t \in (-r,r)$.
Then 
\begin{equation}\label{24.20}
\H^2(H_2 \cap B(0,r)) = \int_{-r}^r \H^1(P_t \cap H_2 \cap B(0,r)) dt,
\end{equation}
because $e_2$ is parallel to $H_2$. For each $t$, $\H^1(P_t \cap H_2 \cap B(0,r))$
is less that what it would be if $H$ was the half plane bounded by $L_\infty$ and
that contains $d_\infty e_1 + e_3$, and even less than what it would be (for $H_2$ with
the same direction and) for $d_\infty = 0$.  
See Figure \ref{f24.1}. 
Thus, after integrating, 
$\H^2(H_2 \cap B(0,r)) \leq \frac{\pi r^2}{2}$, $A \leq \frac{3\pi r^2}{2}$, and \eqref{24.18},
with $r$ close enough to $9R/10$, implies that \eqref{24.16} holds for $k$ large. This completes 
our proof of \eqref{24.16} by compactness and contradiction.
As was said earlier, \eqref{24.4} follows from \eqref{24.16}, and this last case ends our proof of
Lemma \ref{t24.1}.
\qed

 % p199

\ms
In the case of Theorem \ref{t23.1}, we were able to replace our assumption that $\beta_H(R)$
is small by a density assumption, but for Theorem \ref{t23.2}, the corresponding assumption that
\begin{equation}\label{24.21}
\theta_0 = \frac{3\pi}{2} \ \text{ and }\ 
F(R) \leq \frac{3\pi}{2} + \varepsilon_V,
\end{equation}
where $\theta_0 = \lim_{r\to 0} r^{-2}\H^2(E\cap B(x,r))$ is still as in \eqref{24.1},
is not enough to give a good approximation by a set of type $\bV \cup P_0$. We can only do this when
$d_0$ is small enough, as in the following lemma that we state for general $r \in (0,R)$.
The initial assumptions are as in Theorem \ref{t23.2}, but we replace \eqref{23.10} with \eqref{24.21}.

\begin{lem} \label{t24.2}
For each $\varepsilon > 0$, there exist $\varepsilon_V > 0$ and $\delta(\varepsilon) > 0$,
that depend only on $n$, $\beta$, and $\varepsilon$, with the following property.
Let $L$, $E$, and $h$ satisfy \eqref{23.1}, \eqref{23.1n}, and \eqref{23.2n} 
with a constant $C_h$ such that \eqref{23.8} holds. Suppose in addition that \eqref{24.21} holds. Then
\begin{equation}\label{24.22}
\beta_{VP}(r) \leq \varepsilon
\end{equation}
for every $r$ such that 
\begin{equation}\label{24.23}
\delta(\varepsilon)^{-1} d_0 \leq r \leq \frac{9R}{10}.
\end{equation}
\end{lem}

\ms
We shall proceed as in the previous lemma.
Let $L$, $E$, and $h$ be as in the statement. By \eqref{24.21} and the proof of \eqref{24.15},
\begin{equation}\label{24.24}
F(r) \leq e^{\alpha A(R)} F(R)
\leq e^{\alpha\beta^{-1} \varepsilon_H} F(R) 
\leq e^{\alpha\beta^{-1} \varepsilon_H} (\frac{3\pi}{2} + \varepsilon_V) 
\ \ \text{ for } 0 < r \leq R.
\end{equation}
Let $r$ satisfy \eqref{24.23}.
We want to apply Lemma~\ref{t22.5}, this time with the density $\theta_0 = \frac{3 \pi}{2}$,
$\tau =  \varepsilon$, and the radius $r_1 = \frac{21r}{20}$. 
The statement gives a small constant $\delta$ that depends on $\varepsilon$, 
and we take $\delta(\varepsilon) = \frac{20 \delta}{21}$.

Let us check the assumptions. We start with \eqref{23.1n}, which is satisfied by assumption,
and \eqref{22.7}, which follows from \eqref{23.2n}.
Next \eqref{22.18} (the size of $r_1$)
follows from \eqref{24.23}, \eqref{22.19} follows from \eqref{23.2n} and  \eqref{23.8}
(if $\varepsilon_V$ is small enough, now depending also on $\varepsilon$), and \eqref{22.20} 
follows from \eqref{24.21} and \eqref{24.24}.
We get a set $X_0 \in \bV(L) \cup \bP_0$ (because $\theta_0 = \frac{3\pi}{2}$) 
such that \eqref{22.21} holds. That is, 
\begin{equation}\label{24.25a}
d_{0,r}(E,X_0) = d_{0,\frac{20r_1}{21}}(E,X_0) \leq \varepsilon,
\end{equation}
where the first part holds because $r_1 = \frac{21r}{20}$.
Lemma \ref{t24.2} follows; we also get an extra estimate \eqref{22.22} on the Hausdorff 
measure, which was not mentioned in the Lemma but holds anyway.
\qed

\ms
The lemma is a little weird (and will be improved seriously later), because we should be able 
to prove that the only sets $X_0$ that can show up in the proof above are sets of type $\bV$ 
whose two faces make an angle close to $\2$. For instance, we expect  that otherwise, 
the good approximation of $E$ by a flat object (say, a plane) at the large scale $r$ implies that $0$ 
is a smooth point of $E$. 
For the moment we have to wait for a more precise statement, because we do not seem to have 
the tools yet to prove this, but in Theorem \ref{t28.2} (also see Remark \ref{r28.3}), 
we will get a better result, that also applies to some intermediate radii, of approximation
by truncated $\bY$-sets. This is of course compatible, because in a ball of radius $r >> d_0$,
a truncated $\bY$-sets looks a lot like a $\bV$-set with an angle close to $\2$.

\ms
In the situation of Lemma \ref{t24.2}, but for radii $r < \delta(\varepsilon)^{-1} d_0$,
we can still get some geometric control on $E$, and show that  it looks like a truncated 
cone of type $\bY$.

\begin{lem} \label{t24.3}
For each choice of $\varepsilon > 0$ and $\delta \in (0,1)$, there exist $\varepsilon_V > 0$,
that depends only on $n$, $\beta$, $\varepsilon$, and $\delta$, with the following property.
Let $L$, $E$, and $h$ satisfy \eqref{23.1n} and \eqref{23.2n} with a constant $C_h$ 
such that \eqref{23.8} holds. Suppose in addition that \eqref{24.21} holds. Then for each $r$ such that
\begin{equation}\label{24.25}
\frac{22 d_0}{21} \leq r < \min\Big(\frac{20R}{21},\delta^{-1} d_0\Big)
\end{equation}
we can find a minimal cone $Y$ of type $\bY$, centered at $0$, such that
$L \cap B(0,\frac{21r}{20})$ is contained in a face of $Y$, and for which
\begin{equation}\label{24.26}
d_{0,r}(E,E_0) \leq \varepsilon, \ \text{ with } E_0 = \ol{Y \sm S}.
\end{equation}
\end{lem}

\ms 
Here $S$ still denotes the shade of $L$, and we may observe that in the ball $B(0,\frac{21r}{20})$
(the only place that counts for \eqref{24.26}), $E_0$ is a nice truncation of $Y$ by $L$.
The proof will also show that $E_0$ approximates $E$ well in measure, in the sense that 
\begin{equation}\label{24.27}
|\H^2(E \cap B(y,t)) - \H^2(E_0 \cap B(y,t))| \leq 2\varepsilon r^2
\end{equation}
for all $y \in \R^n$ and $t > 0$ such that $B(y,t) \subset B(0,r)$.

The proof is easy. We want to apply Lemma \ref{t22.4} to the radius $r_1 = \frac{21r}{20}$,
and with $\tau = \varepsilon$. 
The assumptions \eqref{23.1n}, \eqref{22.7}, and \eqref{22.13} follow, as in Lemma \ref{t24.2},
from \eqref{23.1n}, \eqref{23.2n}, and \eqref{23.8}. The replacement for \eqref{22.14}
follows from \eqref{24.21} and \eqref{24.24}. Finally the requirement \eqref{22.12} on the size of $r_1$
was computed to be the same as \eqref{24.25}. Thus we get $Y$ and $E_0$, as in Lemma \ref{t22.4},
and the properties announced in Lemma \ref{t24.3}, as well as \eqref{24.27}, are the same as
what we get from Lemma \ref{t22.4}.
\qed

\msi{\bf The proof of Theorem \ref{t23.1} modulo Proposition \ref{t23.3}.}
Let $E$ be as in the statement. By Lemma \ref{t24.1}, we get that if $\varepsilon_H$
is small enough, \eqref{24.2} and \eqref{24.3} hold with any small $\varepsilon > 0$ that 
we may have chosen in advance. Let us check that we can also get that
\begin{equation}\label{24.28}
\beta_H(r) \leq \varepsilon \ \text{ for }\ \frac{3d_0}{2} < r \leq 9R/10. 
\end{equation}
For $R/2 < r \leq 9R/10$, this follows directly from \eqref{24.2} (with $\varepsilon/2$).
For $r < R/2$, let $\varepsilon_1$ be the value of $\varepsilon_H$ needed to get $\varepsilon$
in Lemma \ref{t24.1}, and apply first Lemma \ref{t24.1} with $\varepsilon =  \varepsilon_1$
to define $\varepsilon_H$. Then by Lemma \ref{t24.1}, $F(10r/9) - \pi \leq \varepsilon_1$.
Next apply Lemma \ref{t24.1} again, with $\varepsilon$, and to the radius $R' = 10r/9$; 
the initial assumptions of Theorem \ref{t23.1} are still valid for $R'$ with the same constants, 
and \eqref{23.6} (with the constant $\varepsilon_1$) is true because $F(10r/9) - \pi \leq \varepsilon_1$.
We get \eqref{24.2} for $R'$, which is just  \eqref{24.28}.

Next we apply Proposition \ref{t23.3}. If $\varepsilon$ above is chosen smaller than the 
$\varepsilon_1$ from the proposition, the main assumption \eqref{23.14} is satisfied as soon as 
$r \leq 9R/20$. The assumptions \eqref{23.1}, \eqref{23.1n}, \eqref{23.2n},
and \eqref{23.13} are satisfied by assumption
(if $\varepsilon_H < \varepsilon_1$), and so we get the differential inequality \eqref{23.15}, i.e.,
\begin{equation} \label{24.29}
r F'(r) \geq a [F(r) - \pi]_+ - C_1 \int_0^{2r} \frac{h(t) dt}{t}
\end{equation}
for almost every $r$ such that
\begin{equation}\label{24.30}
2 d_0 \leq r \leq \frac{9R}{20}.
\end{equation}
Recall from \eqref{22.3} that $F(r) = \theta(r) + r^{-2} \H^2(S\cap B(0,r))$. The second
term is a smooth function of $r > d_0$, so the differentiability properties of $F$ are the same as
those of $\theta$. Thus, even though we recalled that $F$ is differentiable almost everywhere
in the statement of Proposition \ref{t23.3}, we already knew this from Lemma \ref{16.1}.
The same lemma, and in particular \eqref{16.6}, also says that we can integrate pointwise
inequalities on $\theta'$, and hence on $F'$ too, to get lower bounds on the increase of $\theta$
and $F$. Let us do this with \eqref{24.29}.

We proceed as for \eqref{16.23}, but change some things because we slightly changed
the error term in \eqref{23.15} and \eqref{24.29}.
Set $g(r) = r^{-a} [F(r)-\pi]$ for $r\in I = (2d_0, 9R/20)$; 
then $g$ is also differentiable almost everywhere
(by Lemma \ref{16.1}), with for $r$ in the range of \eqref{24.30}
\begin{equation} \label{24.31}
\begin{aligned}
g'(r) &= - a r^{-a-1} [F(r)-\pi] + r^{-a} F'(r) \geq - C_1 r^{-a-1} \int_0^{2r} h(t) \frac{dt}{t}
\cr&\geq - C_1 C_h r^{-a-1} \int_0^{2r} t^{\beta-1}dt \geq -2C_1 \beta^{-1} C_h r^{\beta-a-1}
=: - C_3 C_h r^{\beta-a-1}
\end{aligned}
\end{equation}
by \eqref{24.29} (and we don't need the positive part), and then \eqref{23.2n}; 
the last identity is a definition of $C_3$.

We may as well assume that $a \leq \beta/2$ (notice that our form of \eqref{23.15},
with the positive part inside, allows us to make $a$ smaller), then $\beta-a-1 \in (-1,0)$. 
Lemma \ref{16.1} allows us to integrate this
(recall that $F-\theta$ is continuously differentiable), and we get that
for $2d_0 \leq r_1 \leq r_2 \leq 9R/20$, 
\begin{equation} \label{24.32}
\begin{aligned}
g(r_1) &\leq g(r_2) - \int_{r_1}^{r_2} g'(r) dr
\leq g(r_2) + C_3 C_h \int_{r_1}^{r_2} r^{\beta-a-1} dr
\cr&\leq g(r_2) + C_3  C_h (\beta -a)^{-1} r_2^{\beta-a}
= g(r_2) + C_4 C_h  r_2^{\beta-a},
\end{aligned}
\end{equation}
with $C_4 = (\beta -a)^{-1} C_3 \leq 2C_3/\beta$, and now
\begin{equation}\label{24.33}
F(r_1)-\pi = r_1^{a} g(r_1) 
\leq \Big(\frac{r_1}{r_2}\Big)^a [F(r_2) - \pi] + C_4 C_h  r_1^a r_2^{\beta-a}.
\end{equation}
This holds for $2d_0 \leq r_1 \leq r_2 \leq 9R/20$, and in this region it is better than
\eqref{23.7} because we don't use the extra $2^a$.

Now we need to consider radii $r_1 < 2d_0$. Let us first check what happens on the
interval $I_1 = [d_0,2d_0]$. In this range, we simply use the fact that by \eqref{22.8},
\begin{equation}\label{24.34}
F(r_1) \leq F(r_1) e^{\alpha A(r_1)} \leq F(r_2) e^{\alpha A(r_2)}
\end{equation}
for $d_0 \leq r_1 \leq r_2 \leq 2d_0$ and, since 
\begin{equation}\label{24.35}
A(r_2) = \int_0^{r_2} h(t) \frac{dt}{t} \leq C_h \beta^{-1} r_2^{\beta}
\leq C_h \beta^{-1} \min(2d_0,R)^\beta \leq C_h \beta^{-1} \varepsilon_H
\end{equation}
by \eqref{23.2n} and  \eqref{23.5}, we see that
$e^{\alpha A(r_2)} \leq 1+3\alpha \beta^{-1} C_h d_0^\beta$.
Also, $F(r_2) \leq 2\pi$ by \eqref{24.3}, and \eqref{24.34} yields
\begin{equation}\label{24.36}
F(r_1) \leq F(r_2) + 20\alpha \beta^{-1} C_h d_0^\beta 
= F(r_2) + C_5 C_h d_0^\beta.
\end{equation}
This is better than \eqref{23.7} because $\frac{2r_1}{r_2} \geq 1$. This is the place where 
we lost the extra $2^a$.

If $d_0 \leq r_1 \leq 2d_0 < r_2 \leq 9R/10$, we combine \eqref{24.35} with \eqref{24.33}
and get that
\begin{equation}\label{24.37}
\begin{aligned}
F(r_1) - \pi &\leq \big[F(2d_0) -\pi \big] + C_5 C_h d_0^\beta
\cr&
\leq \Big(\frac{2d_0}{r_2}\Big)^a[F(r_2) - \pi] + C_4 C_h  (2d_0)^a r_2^{\beta-a}
+ C_5 C_h d_0^\beta
\cr&
\leq \Big(\frac{2r_1}{r_2}\Big)^a[F(r_2) - \pi] + C_6 C_h r_1^a r_2^{\beta-a},
\end{aligned}
\end{equation}
which implies \eqref{23.7}. When $0 < r_1 \leq r_2 \leq d_0$, there is no visible 
sliding boundary condition, and we can use the decay estimates from \cite{Ta}, as
one may find them in \cite{C1}, and which take the same form as in the previous sections,
or even just above with $d_0 = 0$. That is, we get that for some choice of $a > 0$ and $C_7 \geq 1$,
\begin{equation}\label{24.38}
\begin{aligned}
F(r_1) - \pi &= \theta(r)-\pi \leq \Big(\frac{r_1}{r_2}\Big)^a [\theta(r_2) - \pi] 
+ C_7 C_h  r_1^a r_2^{\beta-a}
\cr&= \Big(\frac{r_1}{r_2}\Big)^a [F(r_2) - \pi] + C_7 C_h  r_1^a r_2^{\beta-a}.
\end{aligned}
\end{equation}
For the remaining case when $r_1 < d_0 \leq r_2 \leq 9R/20$, we 
glue this estimate to \eqref{24.37} and get that
\begin{equation}\label{24.39}
\begin{aligned}
F(r_1) - \pi &\leq \Big(\frac{r_1}{d_0}\Big)^a [F(d_0) - \pi] + C_7 C_h  r_1^a d_0^{\beta-a}
\cr&
\leq \Big(\frac{r_1}{d_0}\Big)^a 
\Big\{\Big(\frac{2d_0}{r_2}\Big)^a[F(r_2) - \pi] + C_6 C_h d_0^a r_2^{\beta-a} \Big\}
+ C_7 C_h  r_1^a d_0^{\beta-a}
\cr&
\leq \Big(\frac{2r_1}{r_2}\Big)^a + C_6 C_h r_1^a r_2^{\beta-a} + C_7 C_h  r_1^a d_0^{\beta-a},
\end{aligned}
\end{equation}
which is also as good as \eqref{23.7}.
Theorem \ref{t23.1} follows, modulo Proposition \ref{t23.3} which will be proved 
in the next section. \qed

\msi{\bf The proof of Theorem \ref{t23.2} modulo Proposition \ref{t23.4}.}

This will work essentially as for Theorem \ref{t23.1}.
Let $E$ be as in the statement. 
Let $\varepsilon$ be small, to be chosen soon. 
By Lemma \ref{t24.1}, \eqref{24.4} holds (for $0<r \leq 9R/10$), while by Lemma \ref{t24.2},
\begin{equation} \label{24.40}
\beta_{VP}(r) \leq \varepsilon \ \text{ for }\ \delta(\varepsilon)^{-1} d_0 \leq r \leq \frac{9R}{10}.
\end{equation}

Next we apply Proposition \ref{t23.4}. 
The assumptions \eqref{23.1} (if $r \geq 2d_0$), \eqref{23.1n}, \eqref{23.2n}, \eqref{23.6},
\eqref{23.16}, and \eqref{23.18} come directly from the assumptions of Theorem \ref{t23.2}, 
\eqref{23.19} follows from \eqref{24.40} if $\varepsilon$ is chosen smaller than $\varepsilon_2$, 
and we are left with \eqref{23.17}, which requires that
\begin{equation} \label{24.41}
Nd_0 \leq r \leq \frac{R}{2}.
\end{equation}
Here $N$ is a constant that depends only on $n$ and $\beta$, and since we may now 
choose $\varepsilon$, $\delta(\varepsilon)$ also becomes such a constant. 
Set $N_1 = \max(N,\delta(\varepsilon)^{-1})$.
We get that for almost all $r \in (N_1 d_0,\frac{R}{2})$, $F$ is differentiable at $r$ 
and \eqref{23.20} holds. The same argument as for \eqref{24.33} shows that
\begin{equation}\label{24.42}
F(r_1)-\frac{3\pi}{2} \leq \Big(\frac{r_1}{r_2}\Big)^a [F(r_2) - \frac{3\pi}{2}] + C'_4 C_h  r_1^a r_2^{\beta-a}
\end{equation}
for $N_1d_0 \leq r_1 \leq r_2 \leq 9R/20$.
For $0 < r_1 \leq r_2 \leq d_0$, the proof of \eqref{24.38} yields
\begin{equation}\label{24.43}
F(r_1) - \frac{3\pi}{2} \leq \Big(\frac{r_1}{r_2}\Big)^a [F(r_2) - \frac{3\pi}{2}] + C'_7 C_h  r_1^a r_2^{\beta-a}.
\end{equation}
In the intermediate regions where $d_0 \leq r_1 \leq r_2 \leq N_1 d_0$, we simply use 
the near monotonicity of $F$, as in \eqref{24.36}. Finally, we glue all these estimates as above, 
and get that
\begin{equation}\label{24.44}
F(r_1) - \pi \leq \Big(\frac{N_1 r_1}{r_2}\Big)^a + C_8 C_h r_1^a r_2^{\beta-a} 
\end{equation}
in the full range of $0,R/2$. 
This proves \eqref{23.11} and Theorem \ref{t23.2}, modulo Proposition \ref{t23.4} 
which will be proved in Sections \ref{S25}-\ref{S27}. \qed

\section{Construction of competitors, with the triangle $T(r)$}
\label{S25}

In this section we adapt the main part of the construction of competitors that was done in
Sections \ref{S4}-\ref{S15}, to the case of balls centered away from $L$. The goal is to
get the differential inequalities of Section \ref{S23}, so we will restrict our attention
to the case when $E$ looks a lot like a set $X \in \bH \cup \bV \cup \bP_0$. 

It would be logical to deal also with the case when $E$ looks like a truncated set of type 
$\bY$, but this would only give a differential inequality that holds for a relatively small range of 
radii, and we decided that in this range we will just use the near monotonicity of $F$, and not lose so much.

We will concentrate our attention more on the case when $E$ looks like a set of type $\bV \cup \bP_0$,
because it is a little more complicated, and also seems more useful. That is,
we could probably manage without the case of a half plane. But this case is easier anyway.

We start the construction with assumptions relative to a fixed radius $r < R$, where 
$R$ is as in \eqref{23.1n} or \eqref{23.2n}, 
and we assume that the assumptions of Proposition \ref{t23.3}
or \ref{t23.4} are satisfied. In particular, we assume that $0 \in E$,
\begin{equation}\label{25.1}
\theta_0 = \lim_{\rho \to 0} \rho^{-2} \H^2(E \cap B(0,\rho)) \in \big\{ \pi, \frac{3 \pi}{2} \big\}
\end{equation}
and that there is a set $X \in \bH \cup \bV \cup \bP_0$, such that
\begin{equation}\label{25.2}
d_{0,2r}(E,X) \leq \varepsilon.
\end{equation}
Here $\varepsilon$ is a shortcut for $\varepsilon_H$ or $\varepsilon_V$, 
$X$ is a half plane if $\theta_0 = \pi$ and a set of type $\bV$ or $\bP_0$ if $\theta_0 = \frac{3 \pi}{2}$,
and \eqref{25.2} comes from \eqref{23.14} or \eqref{23.19}. But we also have the extra information 
that 
\begin{equation}\label{25.3}
2d_0 \leq r \leq \frac{R}{2} \ \text{ if } \theta_0 = \pi
\end{equation}
and
\begin{equation}\label{25.4}
Nd_0 \leq r \leq \frac{R}{2} \ \text{ if } \theta_0 = \frac{3 \pi}{2},
\end{equation}
where we can choose $N$ as large as we want, and that
\begin{equation}\label{25.5}
F(\rho) \leq \theta_0 + \varepsilon \ \text{ for } 0 < \rho \leq r
\end{equation}
by \eqref{24.3} or \eqref{24.4}. 

Our proof of differential inequalities will follow the same route as for Proposition \ref{t16.2}; 
fortunately, we do not need to repeat everything, and the geometric situation will be simpler. 
We explain how it goes here, but the truth is that no real difference with what was 
done before, except for some occasional simplification, happens before the description of 
Section~\ref{S13}, where we build a competitor, and where the triangle $T(r)$ will show up.

We start as in Section \ref{S4}; our assumptions \eqref{23.2n} and \eqref{23.5} replace 
\eqref{4.1} and \eqref{4.2}, and \eqref{25.2} replaces \eqref{4.3}. 
We also assume that $r$ satisfies the extra properties \eqref{4.4}, \eqref{4.7}, and \eqref{4.8};
this is all right, because we noticed in Section \ref{S4} that they hold almost everywhere.
These extra assumptions were used to take some limits, for instance when we proved
Lemma~\ref{t13.2}, and we will apply the same arguments here. 

So we fix $r$ with all these properties. For the moment, let us not normalize $r$ away 
(i.e., take $r=1$) as we did earlier. We want to construct nice competitors for $E$ in 
$\ol B(0,r)$, that probably beat the natural one. Earlier, the natural one was just the cone 
$\Gamma(E,r)$ over $E \cap \S_r$ (where $\S_r = \d B(0,r)$), but now it is 
\begin{equation}\label{25.6}
\wt\Gamma(E,r) = \Gamma(E,r) \cup T(r),
\end{equation}
where $T(r)$ is the triangle with vertices $0$ and the two points 
of $L \cap \S_r$, which we denote by $\ell_\pm = \ell_\pm(r)$. 
This is the same set that was used to prove the near monotonicity formula \eqref{22.8}, 
and hopefully if we do even better than $\wt\Gamma(E,r)$, we will get the desired
differential inequality. Set 
\begin{equation}\label{25.7}
K_r = X \cap \S_r = X \cap \d B(0,r).
\end{equation}

When $X \in \bH$, $K_r$ is composed of one nice curve $\rho_1$. It is the intersection of
$\S_r$ with a half plane bounded by $L$, which by \eqref{25.3} passes rather near $0$; it is not a 
piece of geodesic unless $0$ is exactly in front of $X$, but it is still an arc of circle with a not too 
large curvature. And it has two endpoints in $L \cap \S_r$.

When $X \in \bV$, $K_r$ is composed of two nice curves $\rho_1$ and $\rho_2$, both ending at the
two points of $L \cap \S_r$. They are not geodesics either in general, but since \eqref{25.4} says that
$L$ passes very near $0$ (as near as we want), they are very close to being two arcs of great circle with
length $\pi r$.

Finally, when $X \in \bP_0$, $K_r$ is a full great circle, that may or may not contain points of $L$.
This case is slightly different from the others, but we keep it along for some time. We cut $K_r$
in two roughly equal parts $\rho_1$ and $\rho_2$ of lengths nearly equal to $\pi r$; we may be more 
specific later on where we cut, to make some estimates easier to prove.

For $i = 1$ and maybe $2$, let $w_i$ denote the point of $\rho_i$ that lies at equal 
distance from its two endpoints. 
We cut $\rho_i$ at $w_i$, into two sub-arcs $\rho_{i,\pm}$ that go from $w_i$ to $\ell_\pm$
in the first two cases; in the last case, $\rho_{i,\pm}$ goes from $w_i$ to an endpoint $\ell'_\pm$
of the two $\rho_i$, which we choose close to $\ell_\pm$ if $K_r$ gets close to $L$.
Since we want uniform notation, let us also set $\ell'_\pm = \ell_\pm$ in the first two cases.
At this point we have two or four nice arcs $\rho_{i,\pm}$ from the $w_i$ to the $\ell'_\pm$.

When $X \in \bV \cup \bP_0$, we have a constant $N \geq 1$, as in \eqref{25.4}, 
which we can make larger if we want, so that some geometric estimates are satisfied; 
when $X\in \bH$, let us also introduce a large constant $N$ too,
which this time is not related to the constraint \eqref{25.3}. 
We will pick $N$ so large (in both cases), that some geometric properties are satisfied, 
and then $\varepsilon$ will be allowed to be small enough, depending on $N$. 
For instance, we claim that for $N$ large, $\rho_{i,\pm}$ is close to the geodesic $\rho(w_i,\ell_\pm)$ 
with the same endpoints, and more precisely
\begin{equation}\label{25.8}
d_{0,2r}(\rho_{i\pm}, \rho(w_i,\ell'_\pm)) \leq 10 N^{-1}.
\end{equation}

When $X\in \bP_0$, this is trivial because $\rho_{i\pm}$ is actually an arc of geodesic; 
in our remaining case, recall that $\ell'_\pm = \ell_\pm$.

When $r \geq Nd_0$ (which is automatically the case when $X\in \bV$), 
this is because $\S_r$ is almost centered on $L$ (and we put the large constant $10$ to
be sure that neither author nor reader has to think about it). Otherwise, $X$ is a half plane
bounded by $L$, $\dist(0,X) \leq 2r\varepsilon$ by \eqref{25.2}, and we just assumed that $d_0 \geq N^{-1} r$. 
Then $X$ makes a very small angle with the half plane $H_0$ bounded by $L$ and that contains $0$, 
and $K_r$ lies very close to an arc of $H_0$ through the $\ell_\pm$, which happens to be a geodesic. 
So \eqref{25.8} holds in this last case too, if $\varepsilon$ is small enough.

In the discussion below, we shall some times say things as if $X$ were of type $\bV$ or $\bP_0$
and we had two curves $\rho_i$, but the case of a half plane will be included, and easier.

We give ourselves a small constant $\tau > 0$, that depends on the geometry; probably
we can take $\tau = 10^{-5}$. We set $D_\pm = D_\pm(r) = \S_r \cap B(\ell'_\pm, \tau r)$.
The discussion of Section \ref{S5}, with the local regularity of $E$ far from $L$,
gives a nice description of $E \cap \S_r \sm (D_+ \cup D_-)$, as a union of one or two
nice $C^1$ curves $\cL_i$, that are also small Lipschitz graphs over $\rho_i \sm (D_+ \cup D_-)$.
We cut the curve $\cL_i$ into two pieces $\cL_{i,\pm}$, at a point which we call $m_i$ and
which we choose very close to $w_i$.
The curve $\cL_{i,\pm}$ leaves from $m_i$, and ends at a point $c_{i,\pm} \in E \cap \d D_\pm$,
where it actually reaches $\d D_\pm$ transversally. In addition, near each $\d D_\pm$, the intersection
$E \cap \S_r$ is just composed of two nice $C^1$ curves, that are extensions of the 
$\cL_{i,\pm}$, and which cross a thin annulus near $\d D_\pm$ transversally.

The behavior of $E \cap \S_r$ in each $D_\pm$ can be classified into what we shall call 
configurations. Their description is essentially the same as 
in Section \ref{S16}, except with only one or two points in $\d D_\pm$, but we use 
the opportunity to change the vocabulary slightly.

We start with \ub{Configuration $1$} (which we may also call Configuration $1_\pm$ if
we want to specify near which $\ell_\pm$ we work), where
both points $c_{i,\pm}$ lie in the same connected component of $E \cap D_\pm$
as $\ell_\pm$ (which therefore lies in $E$). This is the most likely situation, hence the name.

We also have \ub{Configuration $2$} (or $2_\pm$) where we have what we call a hanging curve,
i.e., when at least one of the $c_{i,\pm}$ (and say that it is $c_{1,\pm}$) does not lie in
the component of $c_{2,\pm}$ in $E \cap D_\pm$ (if $X$ is a $\bV$-set), 
nor in the component of $\ell_\pm$ (if $\ell_\pm$ happens to lie in $E$). 
We will like this case because it is easy to construct a better competitor.

When $X \subset \bH$, these are the only two options, since there is only one $c_{i,\pm}$.
Otherwise, we still have one possibility, \ub{Configuration $3$} (or $3_\pm$), where 
$c_{1,\pm}$ and $c_{2,\pm}$ lie in the same component of $E \cap D_\pm$, but not $\ell_\pm$.
We call this a free attachement; we expect this thing to happen, but only when $X$ is 
very close to a plane, and then $E$ may leave $L$. 
Recall that every plane is a sliding minimal set, independently of its position relative to $L$,
and that $X$ may also be a plane that does not contain $L$. We treat this case like the other
ones for the moment, except that maybe $\ell'_\pm \neq \ell_\pm$, and in the case when
$|\ell'_\pm -\ell_\pm| \geq \tau/10$, say, and $\varepsilon$ is small, we are sure to be 
in configuration $3$, and near $\ell'_\pm$ we may also have kept the curve $\cL_1 \cup  \cL_2$, 
which is perfectly nice (and does not get close to $L$).

\ms
Next we construct a net of curves, as in Section \ref{S6}.

When $X$ is a half plane and we are in Configuration 1, we find a simple curve 
$\gamma'_{1,\pm} \subset E \cap D_\pm$ that goes from $c_{1,\pm}$ to $\ell_\pm$,
and we add it to $\cL_{1,\pm}$ to get a curve $\gamma_{1,\pm} \subset E$
that goes from $m_1$ to $\ell_\pm$. Also set $\gamma_{\pm} =\gamma_{1,\pm}$ for unified notation.

When $X$ is of type $\bV$ or $\bP_0$ and we are in Configuration 1, 
we find a point $z_\pm\in E \cap D_\pm$,
and two simple curves $\gamma'_{i,\pm} \subset E \cap D_\pm$ that go from $c_{i,\pm}$ to $z_\pm$,
$i=1,2$, and a last one, $\gamma_{\ell_\pm}\subset E \cap D_\pm$ that goes from $z_\pm$
to $\ell_\pm$. We add to $\gamma'_{i,\pm}$ the corresponding $\cL_{i,\pm}$ to get a curve 
$\gamma_{i,\pm} \subset E$ that goes all the way to $m_i$, and call
$\gamma_\pm = \gamma_{1,\pm} \cup \gamma_{2,\pm} \cup \gamma_{\ell_\pm}$. 
We allow the degenerate case when $z_\pm = \ell_\pm$.

In configuration 2, when $\cL_{1,\pm}$ is a hanging curve (say), we decide to essentially remove 
$\cL_{1,\pm}$ and the component of $c_{1,\pm}$ in $D_\pm$ from the game,
and we set $\gamma_{1,\pm} = \{ m_1 \}$. If $X$ is a half plane, we are finished with $D_\pm$.
Otherwise, if $\cL_{2,\pm}$ is also a hanging curve, we set $\gamma_{2,\pm} = \{ m_2 \}$.
If not, we construct $\gamma_{2,\pm}$ as in Configuration 1, when $X$ is a half plane.
We also set $\gamma_{\pm} =\gamma_{1,\pm} \cup \gamma_{2,\pm}$.

We are left with configuration 3. In this case we select a simple curve 
$\gamma'_{\pm} \subset E \cap D_\pm$ that goes from $c_{1,\pm}$ to $c_{2,\pm}$,
pick a point $z_\pm \in \gamma'_{\pm}$ close to $\ell'_\pm$
(a point of $\gamma'_{\pm}$ closest to $\ell_\pm$ seems to be the simplest), 
cut $\gamma'_{\pm}$ at the point $z_\pm$, into two pieces
$\gamma'_{i,\pm}$ that go from $c_{i,\pm}$ to $z_\pm$, and add $\cL_{i,\pm}$ to 
$\gamma'_{i,\pm}$ to get a longer curve $\gamma_{i,\pm}$ that starts from $m_i$ and ends 
at $z_\pm$. Finally set $\gamma_{\pm} =\gamma_{1,\pm} \cup \gamma_{2,\pm}$ as usual.

This gives a collection of simple curves. We may call $\gamma^\ast$ the union of these curves,
and we see it as a first net. The curves don't intersect, because we kept $D_+$ and $D_-$ essentially
disjoint from the rest.

Next we proceed as in Section \ref{S7}. Each of the simple curves (call it $\gamma$)
that was constructed above is replaced with a small Lipschitz graph $\Gamma$ with the same endpoints.
When we deal with a hanging curve, of course, we don't see the difference, because both curves
$\gamma$ and $\Gamma$ are reduced to a point. In Configuration 1 when $X \in \bV \cup \bP_0$,
it could be that the three Lipschitz graphs that we constructed do not make nice angles or, 
even worse, intersect; then we will modify it later appropriately, but let us not worry for the moment 
and continue as if this did not happen. 
We take the union of all these curves $\Gamma$ and get a net that we call $\Gamma^\ast$.

Recall that in \eqref{7.1} we required the endpoints of the curve $\gamma$ that we transform 
into $\Gamma$ not to be too far from each other. In the present case, if $X$ is of type 
$\bV$, then $d_0/r$ is quite small, the curves $\cL_i$ stay quite close to diameters of $\S_r$,
the $w_i$ are at distance about $\pi r/2$ from the $\ell_\pm$, so do the $m_i$,
and we get \eqref{7.1} because we can choose $\tau$ small and hence the point $z_{\pm}$ 
(when it exists) lies quite close to $\ell_\pm$. 
If $X \in \bP_0$ and the $\ell'_\pm$ lie close to the $\ell_\pm$, we can argue as when 
$X\in \bV$, while otherwise, when for instance $\ell'_+$ lies at distance at least $\tau r$ from $L$,
we may as well have taken it opposite to $\ell'_-$ and get what we want. Finally, 
assume that $X$ is a half plane; if $d_0$ is small compared to $r$, we can still say the 
same thing. And even when $d_0$ is large, $X$ is quite close to the half plane $H_0$ 
that is bounded by $L$ and goes through the origin, $w_1$ and $m_1$ are quite close 
to the middle point of the long geodesic arc $\S_r \cap H_0$. 
But we required in \eqref{25.3} that $2d_0 \leq r$, and it can be 
checked that this forces the length of $\S_r \cap H_0$ to be significantly smaller than 
$\frac{3 \pi r}{2}$, so that  $\H^1(\rho(m_1,\ell_\pm)) \leq \frac{3 \pi r}{4}$ as needed for \eqref{7.1}.
The reader is invited not to do the precise computation, since $\frac{3 \pi r }{4}$ could have been 
any number smaller than $\pi r$, and $\pi r$ corresponds to $d_0 = r$. 
So we did not cheat with the assumption \eqref{7.1} here.

We do not need to modify what we did in Section \ref{S8}. That is, for each of the curves
$\Gamma$ that compose $\Gamma^\ast$, we construct a graph $\Sigma_G(\Gamma)$, 
which is bounded by $\Gamma$ and the two line segments from $0$ to the endpoints of $\Gamma$,
and which in general has a smaller area than the cone $\Sigma_F(\Gamma)$ over 
$\Gamma$. See for instance \eqref{8.18}.

The reader may be worried about the fact that in the present situation,
$\Sigma_G(\Gamma^\ast) = \bigcup_\Gamma \Sigma_G(\Gamma)$ does not give a competitor 
for $E$ (even after we do the small modification needed to glue things near $\S_r$), 
because probably $\Sigma_G(\Gamma^\ast)$ detaches itself from $L$ when it leaves 
the two points $\ell_\pm$. We already had this problem in \cite{Mono}, because
in the proof of \eqref{22.8} we could not use 
$\Sigma_F(\Gamma^\ast) = \bigcup_\Gamma \Sigma_F(\Gamma)$;
this is why we added the triangle $T(r)$ to $\Sigma_F(\Gamma^\ast)$, 
and here again we will need to add $T(r)$ to $\Sigma_G(\Gamma^\ast)$ when needed; 
we will take care of this later.

But let us first continue with the flow of the previous sections and discuss 
analogue of Sections \ref{S9}-\ref{S11}. We said above how to construct curves $\Gamma$
by taking the same endpoints as for our initial curves $\gamma$ and applying Section \ref{S7}.
This is what we do in most cases, but there is one case when we apply the construction of 
Section \ref{S7} to slightly different curves. This  is when we are in Configuration 1, 
with a set $X \in \bV \cup \bP_0$, and in addition the three endpoints 
$\ell_\pm$, $m_1$, and $m_2$, seen from our vertex $z_\pm$, make wrong angles. 
That is, if we are so lucky that 
\begin{equation}\label{25.9}
\Angle_{z_\pm}(\ell, m_i) \geq \frac{\pi}{2}+10^{-1} \ \text{ for } \ i=1,2,
\end{equation} 
we proceed exactly as we said above, and construct three curves $\Gamma$ with 
the same endpoint $z_\pm$, and glue them together. This is all right, because then the 
curves $\Gamma$ of our net $\Gamma^\ast$ still make large angles at $z_\pm$, 
and this will allow us to produce nice retractions on $\Gamma^\ast$.
Notice that we do not need to say that $\Angle_{z_\pm}(m_1, m_2) \geq \frac{\pi}{2}+10^{-1}$,
because it is automatic, as either $X\in \bV$ and the two faces of $X$ make an angle at least $\2$
at the $\ell_\pm$ (and $z_\pm$ lies in the small disk $D_\pm$ centered at $\ell_\pm$),
or $X \in \bP_0$ and this is even more clear.

The worse picture we have when \eqref{25.9} holds is that $\Gamma^\ast$ is composed of 
six arcs (with two short ones) that make $\Gamma^\ast$ look like like two long arcs of circle, plus two
little branches that connect the ends to the $\ell_\pm$. But of course we could also have a free attachment on one side, or simpler pictures.

In the bad case when \eqref{25.9} fails, it seems that we have no other choice than proceed as in 
Configuration 2+ of the earlier discussion, which is treated in Section \ref{S11}. 
More precisely, as in Subcase B where \eqref{11.2} fails (just as \eqref{25.9} fails here).
In this case we decide that $z_\pm$ is not a nice enough center, and we use only two curves 
$\wt\Gamma_1$ and $\wt\Gamma_2$, that go directly from $\ell_\pm$ to $m_1$ and $m_2$.
The two curves $\wt\Gamma_i$ are constructed as in \eqref{11.7}, as unions of 
the geodesic $\rho(\ell_\pm,z_i)$ and the part of our old Lipschitz curve $\Gamma_{\pm,i}$ 
(from $z_\pm$ to $m_i$) that lies between $z_i$ and $m_i$, 
where $z_i$ is a point of $\Gamma_{\pm,i}$ that lies reasonably far from $z_\pm$
(as defined below \eqref{11.6}, but beware that $r$ there has a different, local, meaning).

Fortunately, the computations of Section \ref{S11} are still valid in this case, 
and we do not repeat them. Their result is twofold. First, what we get when we add 
$\wt\Gamma_1$ and $\wt\Gamma_2$ to the other curves that we construct is still a nice net, 
composed of at most five Lipschitz curves (four long ones whose union look likes the union 
of two half circles with common endpoints, and a short one that connects the other
$\ell_\mp$ to its $z_\mp$) that are disjoint except for their endpoints, 
and make large angles with each other at these points. And we have good estimates 
like \eqref{11.33}-\eqref{11.38} that say that the measure of the symmetric difference
between our initial $\gamma_\pm$ and $\wt\Gamma_1 \cup \wt\Gamma_2$ is controlled 
by what we will win in the estimates, as in \eqref{9.6} and \eqref{9.7}.

We return to the general case. At this point, we have a nice net $\Gamma^\ast$ composed of 
at most six Lipschitz graphs, which we now decide to call $\Gamma_j$ 
(hence, $1 \leq j \leq 6$, maybe less), which are glued together at their endpoints and make 
reasonably large angles there. For each $j$ there is a cone $\Sigma_{F}(\Gamma_j)$ over $\Gamma_j$ 
(and these cones are glued nicely along segments that go from $0$
to the endpoints of the $\Gamma_j$), and a nicer graph $\Sigma_{G}(\Gamma_j)$
(and these graphs are also nicely glued along the same line segments, with reasonably large angles). 
We set $\Sigma_{F}(\Gamma^\ast) = \cup_j \Sigma_{F}(\Gamma_j)$ and 
$\Sigma_{G}(\Gamma^\ast) = \cup_j \Sigma_{G}(\Gamma^\ast)$.

Denote by $a_j$ and $b_j$ the endpoints of $\Gamma_j$; recall that $\Sigma_{G}(\Gamma_j)$ 
is a small Lipschitz graph over its projection on the plane $P_j$ that contains the geodesic 
$\rho(a_j,b_j)$. 
Also, $\Sigma_G(\Gamma_j)$ is bounded by $\Gamma_j$ (on the sphere) 
and the two segments $[0,a_j]$ and $0,b_j]$, where it is glued to the rest 
of $\Sigma_{G}(\Gamma^\ast)$. 
That is, if $\Gamma_j$ and $\Gamma_k$ share an endpoint, which we call $a = a_j = a_k$,
then $\Sigma_G(\Gamma_j)$ and $\Sigma_G(\Gamma_k)$ also share the segment $[0,a]$, 
and they make an angle at of least $\pi/2$ along that segment. 

Indeed, by the small Lipschitz graph description of $\Sigma_G(\Gamma_j)$ and 
$\Sigma_G(\Gamma_k)$, we just need to control the angle of $\Gamma_j$ and $\Gamma_k$. 
This is rather easy when $a$ is one of the $m_i$, because then the other endpoints $b_j$ 
and $b_k$ essentially lie in opposite directions seen from $a$.
This is also easy when $a = \ell'_\pm$ and $b_j$, $b_k$ are the two points $m_i$, 
because either $X \in \bV$ and its two faces make a large angle at $\ell_\pm = \ell'_\pm$, 
or $X\in \bP_0$ and the $w_i$ lie in a geodesic $X \cap \S_r$, in different directions.

We are left with the case when $a = z_\pm$.
When $b_j$ and $b_k$ are the two points $m_i$, we can apply the same reasoning as above,
since $z_\pm$ lies in a small disk $D_\pm$ centered at $\ell'_\pm$. So we may assume
that we are in Configuration $1_\pm$, and for instance $b_j = m_1$ and $b_k = \ell_\pm$.
But in this case we have \eqref{25.9} (because otherwise we decided to start the curves from
$\ell_\pm$), which is exactly what we need. 

So $\Sigma_G(\Gamma_j)$ and $\Sigma_G(\Gamma_k)$ 
make an angle at of least $\pi/2$ along $[0,a]$, and the union $\wt\Sigma_G(\Gamma^\ast)$
is a nice object.  Seen from far (and if $X\in \bV$), it still looks like a set of type $\bV$, 
but maybe pinched twice along two thin triangular surfaces. 
Notice also that we do not say that $\wt\Sigma_G(\Gamma^\ast)$ lives in a $3$-dimensional space,
but it stays quite close to the $3$-space that contains $X$ (but maybe not $0$).
We shall also use later the fact that inside $B(0,\kappa r)$, it coincides with 
the cone over a net of geodesics $\rho^\ast$. 

Our next task is to project on $\wt\Sigma_G(\Gamma^\ast) = \Sigma_G(\Gamma^\ast) \cup T(r)$ 
or maybe, in the case of Configuration 2 or 3, on $\wt\Sigma_G(\Gamma^\ast)$ alone. 
This will be useful because we want to find a deformation that starts as the identity outside of $B(0,r)$, 
and maps roughly on $\wt\Sigma_G(\Gamma^\ast)$ inside.

Because of the hanging and free cases, it may be that $\ell_+$ or $\ell_-$ does not lie
in $\Gamma^\ast$, so we add them, which means that now $\Gamma^\ast$ may also have
one or two isolated components (which we call $\Gamma$ again) composed of just a point $\ell_\pm$;
in this case $\Sigma_G(\Gamma)$ is just the line segment $[0,\ell_\pm]$. This is the manipulation 
described at the beginning of Section \ref{S12}.

In the computations that follow, and in order to simplify the notation, we shall return 
momentarily to the convention of Section \ref{S13}, where we had decided to normalize things 
so that $r=1$. So let us suppose that $r=1$, and forget $r$ from some of the notation. We will
return to the correct scaling afterwards.

The next stage, as in Section \ref{S12}, is to construct a Lipschitz projection $p$ that maps points
from a neighborhood of $E\cap \S$ to the net of curves $\Gamma^\ast$. See Proposition \ref{t12.1}
Let us recall its main properties. Set $\wh E = E \cap \S \cup \{ \ell_-, \ell_+ \}$
(again we add the $\ell_\pm$ if they are not there already, because of the free case).
Then there is a very small number $\tau_3$, such that $p$ is defined on a set $R_+ \subset \S$ 
that contains a $\tau_3$-neighborhood of $\wh E$. It is Lipschitz (but maybe with a huge norm),
and it is also locally $30$-Lipschitz, in the sense that for each $x\in R_+$,
the restriction of $p$ to $R_+ \cap B(x,\tau_3)$ is $30$-Lipschitz.
Here $\tau_3$ is allowed to depend on $r$ in a wild way; nonetheless, the local $30$-Lipschitz 
property is useful, because it is enough to give good bounds on the measure of images of sets by $p$.
The reason for this strange local Lipschitz property is that it is rather easy to construct a Lipschitz 
mapping $p_c$ near each component of $\Gamma^\ast$ (see Lemma~\ref{t12.2} for a local version), 
but we need to split $R_+$ into regions where we use different maps $p_c$ 
(i.e., that belong to a given component), but are far from each other (so that $p$ 
is Lipschitz, but maybe with a bad constant).
The main property of $p$ is that
\begin{equation}\label{25.10}
p(R_+) \subset \Gamma^\ast 
\ \text{ and } \ 
p(x) = x \text{ for } x \in \Gamma^\ast.
\end{equation}
Recall that $\ell_\pm \in \Gamma^\ast$ now, so in particular $p(\ell_\pm) = \ell_\pm$.
Also, the local construction component by component is such that, as in \eqref{12.44}
\begin{equation}\label{25.11}
|p(z)-z| \leq 10 \dist(z,\Gamma^\ast) \ \text{ for } z\in R_+.
\end{equation}
Let us now extend $p$ to $E \cap A$, where $A$ is a small annulus around $\S^1$;
we do not take a radial extension as before (see \eqref{13.3}), because we want to preserve $L$,
so we prefer ``radial with respect to $x_0$'', where $x_0$ is the projection of $0$ on $L$.
That is, for $z \in B(0,2) \sm B(0,1/2)$, denote $\xi(z)$ the point of $\S$ such that
$\xi(z)-x_0 = t (z-x_0)$ for some $t > 0$; we take
\begin{equation}\label{25.12}
p(z) = p(\xi(x))  \ \text{ when } \xi(z) \in R_+.
\end{equation}
We refer the reader to Section \ref{S12} for a more precise description of $p$, 
and now turn  to Section~\ref{S13} where we start the description of a new competitor. 
We shall use the set
\begin{equation}\label{25.13}
\wt\Sigma_G(\Gamma^\ast) = \Sigma_G(\Gamma^\ast) \cup T(r),
\, \text{ with as 	above } \Sigma_G(\Gamma^\ast) = \bigcup_j \Sigma_G(\Gamma_j)
\end{equation}
and where $r=1$ here, a set which is a little larger than $\Sigma_G(\Gamma^\ast)$, as the basis 
for our first competitor. Recall that the triangle $T(r)$ is the convex hull of $0$, $\ell_+$, 
and $\ell_-$; it is nicely glued to the sets $\Sigma_G(\Gamma_j)$ for which $\ell_\pm$ is 
an endpoint of $\Gamma_j$. With the recent addition of $\ell_\pm$ to $\Gamma^\ast$, 
these $\Gamma$ exist, but they may be reduced to one point. 
But we do not say that $T(r)$ makes a large angle with the faces $\Sigma_G(\Gamma)$ 
in question; it could even be that $\Gamma$ is almost opposite to $0$ and 
$\Sigma_G(\Gamma)$ has a big intersection with $T(r)$.
We shall see later that this is not a problem.

We want to construct competitors for $E$, in $\ol B(0,1)$ (recall that we often write 
estimates with $r=1$ now), 
and for this we construct two deformations $\varphi^i$. We start with a first one, 
$\varphi^0$, and we define it in rings, starting with the outside. 
We first set $\varphi^0(x) = x$ for $x\in E \sm B(0,1)$, then take a very small
number $\sigma > 0$, that will depend on $r$, and try to do the interesting modifications on
the very thin annulus $A(2\sigma)$.

By the proof of Lemma \ref{t13.1}, if we take $\sigma$ small enough, then 
\begin{equation} \label{25.14a}
\xi(x) \in R_+ \ \text{ for } x \in E \cap A(2\sigma),
\end{equation}
which implies that $p(x)$ is defined there. We first set
\begin{equation}\label{25.14}
\varphi^0(x) = \frac{|x| + \sigma - 1 }{ \sigma}\, x + \frac{1 - |x| }{ \sigma}\, p(x)
\ \text{ for } x \in E \cap A(\sigma),
\end{equation}
just as in \eqref{13.6}. 

Here $A(\sigma)$ is a gluing region; on the exterior boundary $\S_1$,
we just took $\varphi^0(x)=x$, and on the inside boundary, we now have
\begin{equation}\label{25.14bis}
\varphi^0(x) = p(x)  \in \Gamma^\ast
\ \text{ for } x \in E \cap \S_{1-\sigma}.
\end{equation}
The same proof as before yields that if we set $F(\sigma) = \varphi^0(E\cap A(\sigma))$ 
and $M(\zeta) = \H^2(F(\sigma))$, then $M(\sigma)$ is small, as in Lemma \ref{t13.2}.

Next we want all the variation of $\varphi^0$ to occur on the next small ring
$A_2 = A(2\sigma) \sm A(\sigma) = \ol B(0,1-\sigma) \sm B(0,1-2\sigma)$, 
and we shall make sure that 
\begin{equation}\label{25.15}
\varphi^0(x) \in \wt\Sigma_G(\Gamma^\ast) \ \text{ for } x\in E \cap A_2.
\end{equation}
We also want to make sure that $\varphi^0(x)=x_0$ on $\d B(0,1-2\sigma)$,
where $x_0$ is the orthogonal projection of $0$ on $L$, because we will take
\begin{equation}\label{25.16}
\varphi^0(x) = x_0 \ \text{ for } x\in E \cap \ol B(0,1-2\sigma).
\end{equation}
Then we will use $\varphi^0$ to build our first competitor for $E$. 

The construction of $\varphi_0$ on $A_2$ will take some time, because we prefer 
to be explicit. Yet it will probably comfort the reader to know that the specific construction 
that we adopt does not matter much. What counts is the measure of the set 
$\varphi^0(A_2) \subset \wt\Sigma_G(\Gamma^\ast)$, and things like the respect of our 
boundary conditions.

First we want to construct is a deformation of $\Gamma^\ast$
to the the point $x_0$, through $\wt\Sigma_G =  \wt\Sigma_G(\Gamma^\ast)$.
We will define this mapping independently for each $\Gamma = \Gamma_j$, but so that the
different pieces will glue well. 

So let $\Gamma$ be one of our Lipschitz curves, and first assume that 
it is a nontrivial curve that starts at $\ell \in L \cap \d \B$, and ends at a point that we call $a$.
For each $z\in \Gamma$, we want to find a path $t \to w(z,t)$, $t\in [0,2]$,
that goes from $z$ to the final target $x_0$.
We cut the path in two, and assign to $z$ an intermediate target $w(z,1)$ that lies on
the line $[0,\ell]$. We choose $w(z,1)$ ``linearly'', as follows. Parameterize $\Gamma$ by
$v= [0,1] \to \Gamma$  at constant speed, in such a way that $v(0)=\ell$ and
$v(1) = a$. Then set $w(z,1) = (1-s) \ell$ when $z=v(s)$.

The second part of the trip (going from $w(z,1) = (1-s) \ell$ to $x_0$) is easy to organize, 
because $T(r)$ is just a triangle and $[0,\ell]$ one of its sides.
We simply set $w(z,t) = (t-1)x_0 + (2-t) w(z,1) = (t-1)x_0 + (2-t)(1-s) \ell$ for
$z\in \Gamma$ and $1 \leq t \leq 2$.

\begin{figure}[!h]  
\centering
\includegraphics[width=7cm]{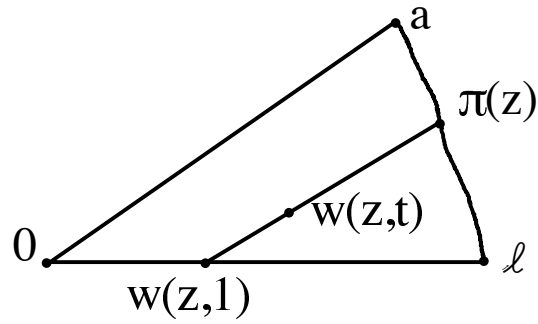}
\caption{The projection $\pi(w(z,t))$ in $S_\Gamma \subset P_\Gamma$
\label{f25.1}}
\end{figure}

For the first part of the trip, we recall that $\Sigma_G(\Gamma)$ is a small Lipschitz
graph over some vaguely triangular sector, which we call $S_\Gamma$, in the plane $P_\Gamma$ 
that contains $a$, $\ell$, and $0$. The two segments $[0,a]$ and $[0,\ell]$ that bound the sector 
are contained in $\Sigma_G(\Gamma)$, and there is a third curvy part of the boundary,
such that $\Gamma$ (the third part of the boundary of $\Sigma_G(\Gamma)$) is a small
Lipschitz graph over that curve.
The triangular sector has a third boundary, that goes from $\ell$ to $a$, and
which is the projection of $\Gamma$ (and is a small Lipschitz curve too). 
We want to make sure that $w(a,t) = (1-t) a$ for $t\in [0,1]$ (we go linearly from
$a$ to its intermediate target $0$), and on the other hand $w(\ell,t) = \ell$ for $t\in [0,1]$.
For the intermediate points $z$, it turns out that we can take
\begin{equation}\label{25.17a}
w(z,t) = \wt G[(1-t)\pi(z) + t w(z,1)] \ \text{ for } t\in [0,1]
\end{equation}
where $\pi$ denotes the orthogonal projection on $P_\Gamma$,
$G : S_\Gamma \to P_\Gamma^\perp$ is the function whose graph is $\Sigma_G(\Gamma)$,
and $\wt G(u) = u + G(u)$ is the parameterization of the graph by $u\in S_\Gamma$.
The point is that although $S_\Gamma$ is not necessarily convex, the small Lipschitz
property of $\Gamma$ implies that all the segments $[\pi(z),w(z,1)] = [\pi(z),(1-s) \ell]$
are contained in $S_\Gamma$ (Figure \ref{f25.1} shows such a segment). That is, 
\begin{equation}\label{25.18a}
w(z,t) \in \Sigma_G(\Gamma) \ \text{ for $z\in \Gamma$ and } t\in [0,1].
\end{equation}

 % p213

This completes the definition of our path function $w(z,t)$ on $\Gamma$, when $\Gamma$
ends at $\ell \in L$. In the other case when $\Gamma$ goes from $a$ to $b$, with 
$a, b \notin L$, we first send points to $0$, and then move them to $x_0$.
That is, write $\Sigma_G(\Gamma)$  as above as the graph over the sector $S_\Gamma$
of the small Lipschitz $G : S_\Gamma \to P_\Gamma^\perp$, and this time take the 
direct path to the origin defined by 
\begin{equation}\label{25.19a}
w(z,t) = \wt G[(1-t)z] \ \text{ for } t\in [0,1];
\end{equation}
as before, \eqref{25.18a} holds (this time because the rays from the curved boundary of 
$S_\Gamma$ to the origin are contained in $S_\Gamma$), and we end with $w(z,1)= 0$.
We complete the path by taking $w(z,t) = (t-1) x_0$ for $t\in [1,2]$, which just
moves from $0$ to $x_0$.
Notice that when two curves $\Gamma$ end at the same vertex $a$, the corresponding
functions $w(a,t)$, $t\in [0,2]$, coincide. This is true both when $a\in L$ and when 
$a \notin L$.

There is a third, trivial case when $\Gamma = \{\ell \}$ for some $\ell \in \S \cap L$;
then we take the same definitions as before (when $\ell$ was an endpoint of $\Gamma$):
we set $w(\ell,t) = \ell$ for $0 \leq t \leq 1$, and $w(\ell,t) = (t-1)x_0 + (2-t) \ell$
for $t\in [1,2]$.

Now we are ready to define $\varphi^0$ on $E \cap A_2$. Recall that 
$\varphi^0(x) = p(x) \in \Gamma^\ast$ for $x\in \S_{1-\sigma}$. For $x\in A_2$,
we still start from $z= p(x) \in \Gamma^\ast$ (which is defined, as before, 
by \eqref{25.12} and \eqref{25.14a}). We also set 
$t(x) = 2 \sigma^{-1} (1 - |x| - \sigma)\in [0,2]$, and finally take
\begin{equation} \label{25.21a}
\varphi^0(x) = w(z,t(x)) = w(p(x),t(x)) \ \text{ for } x\in E \cap A_2.
\end{equation}
Thus on $\S_{1-2\sigma}$, we have $\varphi^0(x) = w(z,2) = x_0$, as needed.

This completes our definition of $\varphi^0$.
Notice also that $\varphi^0(x) \in \wt\Sigma_G$ for $x\in E \cap \ol B(0,1-\sigma)$, as promised, 
and it is easy to check that $\varphi^0$ is Lipschitz (although possibly with a huge norm).
We should also mention that it is easy to find a one-parameter family of Lipschitz mappings 
that preserve $L$ and go from the identity to $\varphi^0$; 
we just need to make sure that points of $L$ stay in $L$, and we don't need to control 
where the intermediate images lie, so we can interpolate linearly and the convexity of $L$
does the rest. With all these remarks in mind, we we just need to check that $\varphi^0 \in L$
for $x\in A(2\sigma) \cap E \cap L$. We made sure to project radially from $x_0$
so that $[x,p(x)] \subset L$, and then $p(x) = \ell_\pm$, and we made sure when we 
retracted $\Gamma^\ast$ to $x_0$ that $w(p(x), t(x)) \in L$ too, as needed.

Now we need to control the measure of $F^0 = \varphi^0(E)$.
We don't need to worry about $E \sm B(0,1)$, because we did not change anything there;
see above \eqref{25.14a}.
Next consider the image of $F(\sigma) = \varphi^0(E \cap A(\sigma))$ (the exterior annulus). 
Fortunately, we took for $\varphi^0$, and in particular in the annulus $A(\sigma)$,
the same sort of formula as in Section \ref{S13} and we can do the same estimates, 
which lead to \eqref{13.26} for $M(\sigma) = \H^2(F(\sigma))$;
recall that the idea is to choose correctly arbitrarily small values of $\sigma$, and define $\varphi^0$
with such values (and later take a limit to get a sharp estimate).
We are left with $\varphi^0(B(0,1-\sigma) \subset \wt\Sigma_G$.
We can keep the same estimates as before, on $\Sigma_G(\Gamma^\ast)$, and then we 
just need to add $\H^2(T(1))$, the additional triangular piece that we decided to add to take care 
of the boundary constraint. We essentially copy \eqref{14.46}, and we get that 
\begin{equation}\label{25.20}
\H^2(E \cap \ol B(0,1)) \leq \frac{1}{2} \H^1(E \cap \S_1) 
- 10^{-5} [\H^1(E\cap \S_1) - \H^1(\rho^\ast)]
+ \H^2(T(1)) +h(1),
\end{equation}
where $\rho^\ast$ is defined as in \eqref{14.43}; it is the union of geodesics that we get 
when we join the endpoints of each $\Gamma$ with a geodesic. Thus there is a maximum of 
six geodesics $\rho_\Gamma$ (they were called $\rho_j$ in Section \ref{S14}), 
four long ones and two short ones. Finally observe that we systematically added 
$\H^2(T(1))$ in our estimate, as if it were disjoint from the other parts. That is,
if by luck $T(1)$ intersects some other piece of $F^0$, we could perhaps have obtained a better
estimate on the total $\H^2$-measure of $F^0$, but we did not try to do this, and this way,
if we later modify $F^0$ by modifying $\Sigma_G(\Gamma^\ast)$ (including a piece that may have intersected $T(1)$), when we do the further estimates, we will be able to compare the measure of the replaced piece of $\Sigma_G(\Gamma^\ast)$ with what it becomes, without having $T(1)$ interfere
in the computation. That is, we shall not actually compare the sets, but the estimates that we use for the sets. Hopefully this comment will become clear when we do this.

This was our main estimate, modulo the full length story below. 
We now forget our normalization $r=1$, and rewrite \eqref{25.20} as 
\begin{equation}\label{25.21}
\H^2(E \cap \ol B(0,r)) \leq \frac{r}{2} \H^1(E \cap \S_r) 
- 10^{-5} r \, [\H^1(E\cap \S_r) - \H^1(\rho^\ast_r)]
+ \H^2(T(r)) + r^2h(r),
\end{equation}
where we now write $\rho^\ast_r$ instead of $\rho^\ast$, and 
which is valid almost everywhere under the assumptions of this section.

In order to know whether we need the full length trick and the construction of 
an additional competitor, we introduce the following set $X_0$, which we will see
as the reference minimal set in the given situation. 
If $\theta_0 = \pi$, $X_0$ is the half plane bounded by $L$ that contains $0$. 
If $\theta_0 = \frac{3\pi}{2}$, $X_0$ is the truncated set of type $\bY$, centered at $0$, 
and with a spine parallel to $L$. That is, we take the cone $Y$ of type $\bY$, centered at $0$, 
and with a face that contains $L$, and we take $X_0 = \ol{Y \sm S}$, where $S$ is still the shade of $L$.
The reader should not get confused (as the author has been a few times); our choice of $X_0$
is just a way to encode some numbers (such as $\H^1(X_0 \cap \S_r)$ below), but we will not compare
$X_0$ with $E$ directly. It is just pleasant to compute things in terms of $X_0$, because we know
that for $X_0$, the functional $F$ is constant, so we know in advance that some simplifications will
occur. Also recall that we know from \cite{Mono}, and with a simpler competitor, that $F$ is almost
nondecreasing, so whatever small improvement that we have should lead to a good differential inequality.
When $X \in \bP_0$, $X_0$ does not look like $X$, but this is all right.

The computations will be simpler when
\begin{equation} \label{25.22}
\H^1(\rho_r^\ast) \leq \H^1(X_0 \cap \S_r),  
\end{equation} 
which will be our analogue of \eqref{15.1}, because in this case we will not need the full
length trick.

So let us assume for the moment that \eqref{25.22} holds, and see how to deduce 
from \eqref{25.21} the differential inequalities \eqref{23.15} and \eqref{23.20}.
In fact, let us check that if we have an inequality like
\begin{equation} \label{25.23}
\H^2(E \cap \ol B(0,r)) \leq \frac{r}{2} \H^1(E\cap \S_r) 
+ \H^2(T(r)) - q r \big[ \H^1(E \cap \S_r) - \H^1(\rho_r^\ast)\big] + r^2 h(r),
\end{equation}
for some number $q \in (0, 10^{-1})$, then we have 
\begin{equation} \label{25.24}
r F'(r) \geq a [F(r) - \theta_0]_+  - C_2 \int_0^{2r} \frac{h(t) dt}{t}.
\end{equation}
where $\theta_0$ is as in \eqref{25.1}, as long as 
such that
\begin{equation} \label{25.25}
a \leq 3q  \ \text{ and } C_2 \geq \max\big(\alpha,\frac{5}{\ln(2)}\big),
\end{equation}
where $\alpha$ is the almost monotonicity constant from \eqref{22.8}. 
Of course \eqref{25.21} is included; it corresponds to $q = 10^{-5}$.

The proof will be similar to what we did in Proposition \ref{t16.2}, but we need to check the algebra.
Otherwise, for the differentiability of $F$, for instance, the justifications are the same as before.

Set $v(r) = \H^2(E \cap B(0,r))$ and $x(r) = (2r)^{-1}\H^1(E \cap \S_r)$
as in Section \ref{S16}; then 
\begin{equation} \label{25.26}
v'(r) \geq \H^1(E \cap \S_r) = 2r x(r)
\end{equation}
as in \eqref{16.13}.  Next 
\begin{eqnarray} \label{25.27}
r^2\theta(r) = v(r) &\leq& \frac{r}{2} \H^1(E\cap \S_r) + \H^2(T(r))
- q r \big[ \H^1(E \cap \S_r) - \H^1(\rho_r^\ast)\big]
+ r^2 h(r)
\nn\\
&\leq& r^2 x(r) + \H^2(T(r)) - q r [\H^1(E \cap \S_r)  - \H^1(X_0 \cap \S_r)] + r^2 h(r)
\end{eqnarray} 
by \eqref{25.23}, the definition of $x(r)$, and \eqref{25.22}.

We want to add $\H^2(S\cap B(0,r))$ to both terms. 
Let $\ol\rho$ denote the arc of great circle that is contained in $P_0$
(the plane that contains $0$ and $L$, goes from $\ell_-$ to $\ell_+$, and
lies on the opposite side of $0$. Thus $\ol\rho$ is the geodesic $\rho(\ell_-,\ell_+)$.
Notice that the positive cone over $\ol\rho$, i.e., 
$\Xi = \big\{ t\xi \, ; \, t\in [0,1] \text{ and } x\in \ol\rho\big\}$ 
is the almost disjoint union of the triangle $T(r)$ and $S \cap B(0,r)$ (a piece of the shade);
thus 
\begin{equation} \label{25.28}
\H^2(T(r)) + \H^2(S \cap B(0,r)) = \H^2(\Xi)
= \frac{r}{2} \H^1(\ol\rho). 
\end{equation}
Also, the union $X_0 \cup S$ is essentially disjoint, and is either the plane $P_0$
(if $\theta_0= \pi$) or a full cone of type $Y$ (if $\theta_0= \frac{3\pi}{2}$); hence
\begin{equation} \label{25.29} 
\H^2(X_0 \cap B(0,r)) + \H^2(S \cap B(0,r)) = \theta_0 r^2. 
\end{equation}
 Next, by \eqref{25.27} and \eqref{25.28}
\begin{eqnarray} \label{25.30} 
r^2 F(r) &=& \H^2(S \cap B(0,r)) + r^2\theta(r) 
\nn\\
&\leq& \H^2(S \cap B(0,r)) 
+ r^2 x(r) + \H^2(T(r)) - q r [\H^1(E \cap \S_r)  - \H^1(X_0 \cap \S_r)] + r^2 h(r)
\nn\\
&=& \frac{r}{2} \H^1(\ol\rho) + r^2 x(r) 
- q r [\H^1(E \cap \S_r)  - \H^1(X_0 \cap \S_r)] + r^2 h(r). 
\end{eqnarray}  
Set $K_0 = X_0 \cap \S_r$; 
then $K_0 \cup \ol\rho$ is the intersection with
$\S_r$ of the full cone $P_0$ or $Y$, and since the union is almost disjoint, we get that
\begin{equation} \label{25.31}
\H^1(K_0) + \H^1(\ol\rho) = 2r \theta_0. 
\end{equation}
Thus \eqref{25.30} becomes
\begin{equation} \label{25.32}
r^2 F(r) \leq  r^2 x(r) + \big[\theta_0 r^2 - \frac{r}{2} \H^1(K_0)\big]
- q r [\H^1(E \cap \S_r)  - \H^1(K_0)] + r^2 h(r). 
\end{equation}
We multiply this by $2r^{-2}$ and get that
\begin{equation} \label{25.33}
2F(r) \leq 2 x(r) + 2\theta_0 - r^{-1}\H^1(K_0) - 2 q r^{-1} [\H^1(E \cap \S_r)  - \H^1(K_0)] 
+ 2 h(r). 
\end{equation}
Next we compute $F'(r)$. The derivative of $\H^2(S \cap B(0,r))$ is $\H^1(\ol\rho)$, so
\begin{equation} \label{25.34}
\begin{aligned}
r F'(r) &= - 2 F(r) + r^{-1} v'(r) + r^{-1} \H^1(\ol\rho)
\geq - 2 F(r) + 2 x(r) + r^{-1} \H^1(\ol\rho)
\cr& =  - 2 F(r) + 2 x(r) + 2\theta_0  -  r^{-1}\H^1(K_0)
\cr& \geq  2 q r^{-1} [\H^1(E \cap \S_r)  - \H^1(K_0)] - 2h(r)
\end{aligned} 
\end{equation}
by \eqref{22.3}, \eqref{25.26}, \eqref{25.31}
and (for the last line) \eqref{25.33}. Next by \eqref{25.33} and the definition of $x(r)$,
\begin{equation} \label{25.35}
\begin{aligned}
F(r) - \theta_0 &\leq x(r) - \frac{1}{2} \, r ^{-1} \H^1(K_0) 
- q r ^{-1} [\H^1(E \cap \S_r)  - \H^1(K_0)] + h(r)
\cr& = \big(\frac{1}{2} - q\big) r ^{-1} [\H^1(E \cap \S_r)  - \H^1(K_0)] + h(r),
\end{aligned}
\end{equation}
hence
\begin{equation} \label{25.36}
r ^{-1} \big[\H^1(E \cap \S_r)  - \H^1(K_0) \big]
\geq \big(\frac{1}{2} - q\big)^{-1} [F(r) - \theta_0] - \big(\frac{1}{2} - q\big)^{-1} h(r).
\end{equation}
We plug this back in \eqref{25.34} and get that
\begin{equation} \label{25.37}
r F'(r) \geq a_0 [F(r) - \theta_0] - b_0 h(r),
\end{equation} 
with $a_0 = 2 q \big(\frac{1}{2} - q\big)^{-1} \geq 3q$
and $b_0 = 2 + 2q\big(\frac{1}{2} - q\big)^{-1} \leq 5$. 

First assume that $F(r) - \theta_0 \geq 0$, so as not to get in trouble with 
the positive part in \eqref{25.24}. 
Since $h(r) \leq \frac{1}{\ln(2)} \int_0^{2r} \frac{h(t) dt}{t}$, \eqref{25.37} is better than
\eqref{25.24} for all the values of $a \leq a_0$ and $C_2 \geq \frac{5}{\ln(2)}$, 
which is a little better than announced in \eqref{25.25}.

Now suppose that $F(r) < \theta_0$. We may also assume that $F'(r)$ exists,
since this is the case almost everywhere. Then, when we differentiate the
monotonicity formula \eqref{22.8}, we get that $F'(r) \geq - \alpha A(r)$.
Since the positive part in \eqref{25.24} vanishes, this establishes \eqref{25.24} in this case,
with any value of $a$ and as soon as $C_2 \geq \alpha$. 

So we finally proved that the desired differential inequality \eqref{25.24} holds, 
with $a$ and $C_2$ as in \eqref{25.25}, as soon as \eqref{25.22} holds. 
We are thus left with the complementary case, when
\begin{equation} \label{25.38}
\H^1(\rho_r^\ast) > \H^1(X_0 \cap \S_r).  
\end{equation}
and we will need the help of a full length condition that we state soon.
But for the moment let us exclude a few cases to make our life simpler later.

We start with the case when $\theta_0 = \pi$, and in addition we have a hanging curve. 
In this case what is left of $\rho_r^\ast$ is just a single geodesic $\rho$, from $m_1$ 
(a point of $E\cap \S_r$ near the middle of $X \cap \S_r$), to one of the two points
of $L \cap \S_r$, say, the point $\ell_+$. Plus a degenerate curve reduced to $\{ \ell_-\}$,
that counts for nothing in the length computations. If we had taken $m_1 = w_1$,
the midpoint of the arc of $X \cap \S_r$, we would have exactly the length
\begin{equation} \label{25.39}
\H^1(\rho(\ell_+,w_1)) = \frac{1}{2} \H^1(X \cap \S_r) \leq \frac{1}{2} \H^1(X_0 \cap \S_r),
\end{equation}
where the second part comes from the fact that $X_0$ is the half plane bounded by $L$
for which $\H^1(X_0 \cap \S_r)$ is the largest. Now we replace $w_1$ by $m_1$,
this hardly changes the lengths, and we still get a contradiction with \eqref{25.38}.
If there are two hanging curves, $\rho^\ast_r$ is reduced to three points, and \eqref{25.38}
is even more impossible.

Let us also exclude the case when $\theta_0 = \frac{3\pi}{2}$ and we have a hanging curve.
This time \eqref{25.4} allows us to take $d_0/r$ as small as we want. Then 
$\H^1(X_0 \cap \S_r)$ is as close as we want to $2\pi r$, while we are still missing
at least one large curve in $\rho_r^\ast$, out of the four, and we get that
$\H^1(\rho_r^\ast)$ is quite close to $\frac{3 \pi r}{2}$, in contradiction with \eqref{25.38}.

Next we use a trick to exclude the remaining case when $\theta_0 = \pi$ and there is
no hanging curve, with only a small computation. In this case $\rho^\ast_r$ is just composed
of two geodesics $\rho_\pm = \rho(m_1, \ell_\pm)$, that connect the $\ell_\pm$ to the point 
$m_1 \in E \cap \S_r$ that we chose below \eqref{25.8}. Recall that $m_1$ is a point
of the curve $\cL_1$, which is the part of $E \cap \S_r$ that lies far from $L$,
that we need to choose near the point $w_1$ in the middle of (the unique arc of) $X\cap \S_r$.
In that region, $E \cap \S_r$ is a nice $C^1$ curve that stays very near $X\cap \S_r$, and by the
intermediate value theorem, we can choose $m_1$ at equal distance from $\ell_+$ and $\ell_-$.
We claim that for such an $m_1$,
\begin{equation} \label{25.40}
\H^1(\rho_r^\ast) = 
\H^1(\rho(\ell_+,m_1))+\H^1(\rho(\ell_-,m_1)) \leq \H^1(X_0 \cap \S_r);
\end{equation}
as soon as we prove this, we will get the desired contradiction with \eqref{25.38}.
So we consider points $m \in \S_r$, at equal distance from the $\ell_\pm$, and
show that $f(m) = \H^1(\rho(\ell_+,m))$ is maximal when $m$ is the point $m_0$
of $P_0 \cap \S_r$ that lies just opposite to $L$ (seen from $0$). 
For this we may assume that $r=1$, and work
in the $3$ space that contains $m$ and $P_0$. Equivalently, we work in $\R^3$,
and we study $f$ on the great circle $\S_1 \cap P$, where $P$ is the vector plane
perpendicular to $L$.  The derivative of $f$ in the direction $v$ is the scalar product of
$v$ with the direction of the geodesic $\rho(\ell_+,m)$ when it arrives at $m$, 
and it is easy to see that this is nonnegative when $v$ points in the direction of $m_0$.
So $f(m)$ is maximal when $m=m_0$, and \eqref{25.40} follows because 
$f(m_0) = \H^1(X_0 \cap \S_r)/2$.

\ms
Return to the proof of the differential inequality \eqref{25.20}. We are left with only 
the case when $\theta_0 = \frac{3 \pi}{2}$, and on each side we have Configuration $1$ or $3$
(also called free attachment). This is where we need a full length estimate.
We have constructed a network $\rho^\ast = \rho^\ast_r$
(we shall often drop the index $r$ again), and it is of the following type.
In all cases, we have selected two points $m_1$ and $m_2$
(near the middle points $w_1$ and $w_2$ of the two arcs of $X \cap \S_r$), 
and two points $z_\pm$, close to the $\ell'_\pm$ (themselves often equal to the $\ell_\pm$). 
By taking $\tau$, and then $\varepsilon$, very small, we can assume
that these four distances are as small as we want compared to $r$. In addition, by taking
$N$ large we may also assume that $d_0/r$ is as small as we want, by \eqref{25.4}.

When we have Configuration 1 on both sides, we take
\begin{equation} \label{25.41}
\rho^\ast = \rho(z_-,m_1) \cup \rho(z_+,m_1) \cup \rho(z_-,m_2) \cup \rho(z_+,m_2)
\cup \rho(\ell_-,z_-) \cup \rho(\ell_+,z_+).
\end{equation}
When we have Configuration 3 on both sides, we only take
\begin{equation} \label{25.42}
\rho^\ast = \rho(z_-,m_1) \cup \rho(z_+,m_1) \cup \rho(z_-,m_2) \cup \rho(z_+,m_2),
\end{equation}
in fact plus the two additional single points $\ell_\pm$ that we added half a page above
\eqref{25.10}, and when we have Configurations $1_-$ and $3_+$, say, we take
\begin{equation} \label{25.43}
\rho^\ast = \rho(z_-,m_1) \cup \rho(z_+,m_1) \cup \rho(z_-,m_2) \cup \rho(z_+,m_2)
\cup \rho(\ell_-,z_-),
\end{equation}
plus the single point $\ell_+$.

To each $\rho^\ast = \rho_r^\ast$ as above, we associate the truncated cones
\begin{equation} \label{25.44}
X'(\rho^\ast) = \big\{ t\xi \, ; \, \xi \in \rho^\ast \text{ and } t\in [0,1] \big\}
\ \text{ and } \ 
X(\rho^\ast) =  T(r) \cup X'(\rho^\ast).
\end{equation}
Notice that even if we had not added the single points $\ell_\pm$ in the free case,
we would add the corresponding segment $[0,\ell_\pm]$ now, with $T(r)$.
There was an additional constraint in the definition of the graphs $\Sigma_G(\Gamma)$ 
associated to our various Lipschitz curves $\Gamma$, which is that 
\begin{equation} \label{25.48a}
\Sigma_G(\Gamma) \cap \ol B(0,\kappa r) = X'(\rho) \cap \ol B(0,\kappa r),
\end{equation}
where $X'(\rho)$ is the cone over the geodesic $\rho$ with the same endpoints as $\Gamma$.
See \eqref{8.13}, which forces the graph of $G$ to be contained in the plane of $\rho$
near the origin. When we take the union, we get that for $\Sigma_G = \Sigma_G(\Gamma^\ast)
= \bigcup_\Gamma \Sigma_G(\Gamma)$,
\begin{equation} \label{25.45}
\Sigma_G \cap \ol B(0,\kappa r) = X'(\rho^\ast) \cap \ol B(0,\kappa r)
\end{equation}
and then, adding $T(r)$,
\begin{equation} \label{25.50a}
\wt \Sigma_G \cap \ol B(0,\kappa r) = X(\rho^\ast) \cap \ol B(0,\kappa r).
\end{equation}
If we have a good competitor for $X(\rho^\ast)$, we can glue it at the tip of $\wt \Sigma_G$,
get a better competitor than $E^0 = \varphi^0(E)$, and improve our main estimate. 
We can even try to do this for $X'(\rho^\ast)$ and $\Sigma_G$, but let us explain what we mean 
by good competitors and how we operate the substitution.

We start with the simpler substitution of a sliding competitor for $X(\rho^\ast)$.
Suppose $Z$ is a sliding competitor for $X(\rho^\ast)$ in $B = \ol B(0,\kappa r/2)$. 
This means that we have a deformation $(x,t) \to f_t(x) = f(x,t)$, defined and continuous
on $X(\rho^\ast) \times [0,1]$, with the usual constrains and in particular
$f_t(x) \in L$ when $x\in L$ and $f_t(x) = x$ when $x\in X(\rho^\ast)\sm B$,
and then we set $Z = f_1(X(\rho^\ast))$. We talk about the whole one-parameter
family $\{ f_t\}$ because it comes with the definition, but (as in the next case), giving $f_1$
alone would be enough as a linear interpolation would complete well (since $L$ is convex).
Extend $f_1$ by setting $f_1(x)=x$ for $x\in \R^n\sm 2B$.
It is easy to see that $f_1$ is still Lipschitz.

We use this to construct a competitor $E^1 = f_1 \circ \varphi^0(E) = f_1(E^0)$.
It is easy to check that this is a sliding competitor for $E$, and the difference between 
$E^0$ and $E^1$ comes from the replacement of $\wt \Sigma_G \cap 2B = X(\rho^\ast) \cap 2B$ 
by $Z\cap 2B$. We are only interested in the replacement if
\begin{equation} \label{25.51a}
\Delta S = \H^2(X(\rho^\ast) \cap 2B) - \H^2(Z \cap 2B) 
= \H^2(X(\rho^\ast) \cap B) - \H^2(Z \cap B) > 0
\end{equation}
(recall that $X(\rho^\ast) = Z$ outside of $B$ anyway),
but when this is the case, we can replace $\varphi^0$ with $\varphi^1$ in the computations
above, find out that we win $\Delta S$ in the intermediate estimate \eqref{25.20}, and proceed
from there on.

\ms
Now let us try to see how we can try to modify a piece of $X'(\rho^\ast)$. 
We try to leave $T(r)$ alone and modify $X'(\rho^\ast)$, but there will be a constraint, 
because we do not want to move the contact region between the two. Set
\begin{equation} \label{25.52a}
L'_\pm = \big\{ t\ell_\pm \, ; \, t \in [0,1] \big\} \ \text{ and } \ L' = L'_+ \cup L'_- .
\end{equation}
A \ub{good competitor} for $X'(\rho^\ast)$ in the same ball $B = \ol B(0,\kappa/2)$ 
as above is a set $Z' = f(X'(\rho^\ast))$, where $f$ is a Lipschitz mapping defined on 
$X'(\rho^\ast) \cup L'$, such that 
\begin{equation} \label{25.53a}
f(x)=x \ \text{ when $x \in L'$ and when $x\in  X'(\rho^\ast)\sm B$},  
\end{equation}
and such that $f(B \cap X'(\rho^\ast)) \subset B$. Notice that we are overkilling something here:
since we added the points $\ell_\pm$ in the free boundary case above, we already have that
$X'(\rho^\ast)$ contains $L'$. But let us keep things like this, because we sometimes tend 
to forget about the one or two extra points.

With this definition we give ourselves a little bit more freedom, because even if $X'(\rho^\ast)$
casually intersects $T(r)$ in an unexpected place, we can pretend not to notice and proceed with
our modification. But we need to be slightly careful when we define our next competitor
$E^1 = \varphi^1(E)$.

So let us define $\varphi^1$. As before, extend $f$ by setting $f(x) = x$ on $\R^n \sm 2B$.
We keep $\varphi^1(x) = \varphi^0(x)$ unless all the following conditions are satisfied:
\begin{equation} \label{25.54a}
x\in A_2, t(x) \in [0,1], p(x) \in \Gamma^\ast \sm \{\ell_+,\ell_-\}, \text{ and }
\varphi^0(x) \in X'(\rho^\ast) \cap B.
\end{equation}
If these conditions are satisfied, we take $\varphi^1(x) = f(\varphi^0(x))$. 
Notice that when $x\in A_2$ and $t(x) \in [0,1]$, the construction gives
$\varphi^0(x) = w(p(x),t(x)) \in \Sigma_G(\Gamma^\ast)$, and
if in addition $\varphi^0(x) \in 2B$, then 
$\varphi^0(x) \in X'(\rho^\ast) \cap 2B$ (by \eqref{25.45}). Then
$f(\varphi^0(x))$ is well defined, and we still have that $\varphi^1(x) = f(\varphi^0(x))$
even if $p(x) \in \{\ell_+,\ell_-\}$ (because then $\varphi^0(x) = w(\ell_\pm, t(x)) = \ell_\pm \in L'$
by the line above \eqref{25.17a}.

We first need to check that that $\varphi^1$ is Lipschitz. Since $\varphi^0$ and $f$
are Lipschitz, we just need to check that
\begin{equation} \label{25.55a}
|\varphi^1(x)-\varphi^1(x')| \leq C |x-x'| 
\end{equation}
when $x\in E$ satisfies \eqref{25.54a} and $x'\in E$ does not. Given the fact that
$|\varphi^0(x)-\varphi^1(x')| = |\varphi^0(x)-\varphi^0(x')| \leq C |x'-x|$
because $\varphi^0$ is Lipschitz, it is enough to show that
\begin{equation} \label{25.56a}
|\varphi^1(x)-\varphi^0(x)|  = |f(\varphi^0(x))-\varphi^0(x)| \leq C |x-x'| 
\end{equation}
under the same conditions, and since $f(x)=x$ on $L'$, this will be proved as soon as
\begin{equation} \label{25.57a}
\dist(\varphi^0(x), L') \leq C |x-x'|.
\end{equation}
We first check this when $x'\notin A_2$. 
Recall from the line above \eqref{25.21a} that $t(x) = 2\sigma^{-1}(1-|x|-\sigma)$;
since $t(x) \leq 1$, $x$ lies at distance at least $\sigma/2$ from $B(0,1-2\sigma)$. 
On the other hand, $t(x) \geq C^{-1}$ because $w$ is Lipschitz, $p(x) \in \S$, and yet 
$\varphi^0(x) = w(p(x),t(x)) \in B$. Then $x$ is also far from $\d B(0,1)$; so
\eqref{25.57a} holds when $x'\notin A_2$.

If $x'\in A_2$ and $t(x') \leq 1$, then either $\varphi^0(x') \notin 2B$, and
\eqref{25.57a} holds because $\varphi^0$ is Lipschitz and 
$|\varphi^0(x)-\varphi^0(x')| \geq \kappa/2$, or else 
$\varphi^1(x') = f(\varphi^0(x'))$ by the remark below \eqref{25.54a}, and we can 
prove \eqref{25.55a} directly without \eqref{25.57a}.

We are left with the case when $x'\in A_2$ and $t(x') \geq 1$.
Recall from the discussion below \eqref{25.16} that $w(p(x),1) = (1-s) \ell \in L'$
(for some $s$), so
\begin{eqnarray} \label{25.58a}
|f(\varphi^0(x))-\varphi^0(x)| 
&\leq& C \dist(\varphi^0(x),L') \leq C |\varphi^0(x)-w(p(x),1)| 
\nn\\
&=& C  |w(p(x),t(x)) -w(p(x),1)| 
\nn\\
&\leq& C |t(x)-1| \leq C |t(x)-t(x')| \leq C |x'-x|
\end{eqnarray}
because $f$ is Lipschitz, $f(x) = x$ on $L'$, and $w$ is Lipschitz; then \eqref{25.56a} 
holds and $\varphi^1$ is Lipschitz.

Let us check that $\varphi^1$ preserves $L$. Let $x \in E \cap L$ be given; 
we want to show that $\varphi^1(x) \in L$, and we already know this when 
$\varphi^1(x) = \varphi^0(x)$, so we may assume that \eqref{25.54a} holds.
But the construction above yields $p(x) = \ell_\pm$ when $x \in E \cap A_2\cap L$,
so \eqref{25.54a} fails and we don't need to prove anything new.

Finally, we should construct a one parameter family $\{ \varphi^1_t \}$ 
that ends with $\varphi^1$, and this is easy; the linear interpolation 
$\varphi^1_t(x) = (1-t)x + t\varphi^1(x)$ does the trick, because $L$ is convex.

We may now use $\varphi^1$ instead of $\varphi^0$ in the computations above.
We compare what we get for the intermediate estimate \eqref{25.20}. Here we replaced
a piece of $\wt \Sigma_G$, namely $X'(\rho^\ast) \cap B$, with its image by $f$,
namely $Z' = f(X'(\rho^\ast)) \cap B = f(X'(\rho^\ast) \cap B)$. These two pieces are
disjoint from the rest of $\wt \Sigma_G$, maybe not from the triangular piece $T(1)$,
but this does not matter, because on all our estimates for $\H^2(E^1)$, where
$E^1 = \varphi^1(E)$ is our competitor for $E$, we sum $\H^2(\wt \Sigma_G)$
and $\H^2(T(1))$. Thus in \eqref{25.20} we can save
\begin{equation} \label{25.59a}
\Delta S = \H^2(X'(\rho^\ast) \cap B) - \H^2(Z' \cap B). 
\end{equation}
Of course we only use $\varphi^1$ instead of $\varphi^0$ when $\Delta S > 0$
(for the radius $r$ under consideration). 

\ms
The next lemma says that in the present situation, when \eqref{25.38} holds,
we can always do one of the two replacements above, and save at least
$\Delta S \geq c \big[\H^1(\rho_r^\ast) - \H^1(X_0 \cap \S)\big]$.

\begin{lem} \label{t25.1}
There is a small constant $c > 0$ such that, for $\rho^\ast = \rho^\ast_r$ as above, 
and keeping the convention that $r=1$ to simplify the statement,
either there is a sliding competitor $Z$ for $X(\rho^\ast)$ in $B = \ol B(0,\kappa/2)$, 
such that
\begin{equation} \label{25.46}
\H^2(Z \cap B) \leq \H^2(X(\rho^\ast) \cap B) 
- c \big[\H^1(\rho^\ast) - \H^1(X_0 \cap \S)\big],
\end{equation}
or there is a good competitor $Z'$ for $X'(\rho^\ast)$ in $B$ such that
\begin{equation} \label{25.46bis}
\H^2(Z' \cap B) \leq \H^2(X'(\rho^\ast) \cap B) 
- c \big[\H^1(\rho^\ast) - \H^1(X_0 \cap \S)\big].
\end{equation}
\end{lem}

\ms 
We postpone the proof of this lemma to the next sections, 
and in the mean time see why it is easy to deduce our differential inequality \eqref{23.20} 
from the lemma. We proceed as explained above, and save 
$\Delta S \geq c \Delta_L$ in  the intermediate estimate \eqref{25.20}, where 
\begin{equation} \label{25.47}
\Delta_L = \H^1(\rho_r^\ast) - \H^1(X_0 \cap \S_r) > 0
\end{equation}
(the inequality comes from \eqref{25.38}, and otherwise we don't do the last step
and don't win anything). 
Thus instead of \eqref{25.21} we now have
\begin{equation}\label{25.48}
\H^2(E \cap \ol B(0,r)) \leq \frac{r}{2} \H^1(E \cap \S_r) 
- 10^{-5} r \Delta_E - c r \Delta_L
+ \H^2(T(r)) + r^2h(r),
\end{equation}
where we decided to set 
\begin{equation} \label{25.49}
\Delta_E = \H^1(E\cap \S_r) - \H^1(\rho^\ast_r);
\end{equation}
notice that even in Configuration 1 when \eqref{25.9} failed and we tampered a little
with the curves, we always made sure to take $\Gamma^\ast$, and a fortiori $\rho^\ast$,
shorter than $E \cap \S_r$, so $\Delta_E \geq 0$; see in particular \eqref{7.16}, and recall that
$\gamma \subset E \cap \S_r$.
We may assume that $c \leq 10^{-5}$, so \eqref{25.48} (and the fact that
$\Delta_E$ and $\Delta_L$ are nonnegative) yield
\begin{eqnarray}\label{bis}
\H^2(E \cap \ol B(0,r)) &\leq& \frac{r}{2} \H^1(E \cap \S_r)  - c r (\Delta_E + \Delta_L)
+ \H^2(T(r)) + r^2h(r)
\nn\\
&\leq& \frac{r}{2} \H^1(E \cap \S_r)  - c r \big[\H^1(E\cap \S_r) -  \H^1(X_0 \cap \S_r)\big]
+ \H^2(T(r)) + r^2h(r).
\end{eqnarray}
This gives directly the second line of \eqref{25.27}, with $q = c$.
Then the same computations as below \eqref{25.27} lead to \eqref{25.37}
(again with $q = c$). We still have the two cases, but as before we obtain \eqref{25.25},
hence also \eqref{23.20} (we no longer care about \eqref{23.15}, because the case when
$\theta_0 = \pi$ was settled before the lemma).

This completes our proof of differential inequalities; hence now Proposition \ref{t23.3}
is established, Theorem \ref{t23.1} follows because of the previous section,
and Proposition \ref{t23.4} and Theorem \ref{t23.2} will follow from Lemma \ref{t25.1}. 
\qed

\section{Basic gain estimates and full length for flat $\bV$ sets}
\label{S26}

% make sure we can let $d = \dist(0,L)$ tend to $0$, or even be $0$ in this section.
% in fact no, we do otherwise
In this section we prove Lemma \ref{t25.1} in most cases. 
The author's initial plan was to use the estimates of this section also for the full length
verifications corresponding to Section \ref{S3b} (with balls centered on $L$), but finally 
decided that this may be confusing; instead we'll do a special argument in Section \ref{S30},
and only import some estimates from this section and the next one.

We are given a net $\rho^\ast_r$ as near \eqref{25.41}-\eqref{25.43}, 
we assume that \eqref{25.38} holds, and we want to find a competitor $Z$ 
for the truncated cone $X(\rho^\ast_r)$, or rather $Z'$ for the truncated cone $X'(\rho^\ast_r)$,
such that \eqref{25.46} (or rather \eqref{25.46bis}) holds.
We shall fulfill this program in this section for most cases, and will be left with a
last, more complicated case, to study in the next one.

We shall try to systematically use the $Z'$ approach, and reserve the approach with 
the sliding competitor $Z$ for a more subtle estimate that may come up later.
The $Z'$ approach required more work to start with, but is more pleasant now because 
we can forget about $T(r)$ and its intersections with the rest of the sets.

We may as well assume again that $r=1$, and we set $\rho^\ast = \rho^\ast_r$ again,
and $X' = X'(\rho^\ast_r)$, the cone over $\rho^\ast$, for simplicity.
The idea of the proof, as for the property of ``full length because of angles'' in \cite{C1}, 
is to show that when  $\Delta_L = \H^1(\rho^\ast) - \H^1(X_0 \cap \S)$ is positive, 
then something in the geometry of $X'$, for instance an angle, allows us to find a better 
competitor.

We keep the same notation as before for $X_0$ (a truncated cone of type $\bY$
with a spine parallel to $L$) and $K_0 = X_0 \cap \S$. Notice incidentally that
$\H^1(X_0 \cap \S) = \H^1(K_0) = 3 \pi - \H^1(\S \cap S)$ 
would stay the same if $X_0$ were replaced by
another truncated cone $X'_0$ of type $\bY$, with a face that contains $L \cap B(0,1)$,
but with a spine that is not parallel to $L$ (but crosses it outside of $B(0,1)$.
This means in particular that if $\rho^\ast = X'_0 \cap \S_1$, we have
$\Delta_L = 0$ and this is fortunate because we could not find a better competitor
$Z$ or $Z'$, since in this case $X(\rho^\ast_r)$ is probably minimal in $B(0,1)$.

Let us give a name to the maximal amount of area that we can save with a competitor
for $X'$, i.e., 
\begin{equation} \label{26.1}
\begin{aligned}
\sigma &= \sup\big\{ \, \H^2(X' \cap B(0,1)) - \H^2(Z' \cap B(0,1)) ; \, 
\cr&\hskip 2cm
Z' \ \text{ is a good competitor for $X'$ in $\ol B(0,1)$} \big\}.
\end{aligned}
\end{equation}
Normally, if we want to relate to what we may win in Lemma \ref{t25.1},
we should consider a competitor in $\ol B(0,\kappa/2)$, but here $X'$ is a cone,
our boundary condition \eqref{25.53a} in the definition of a good competitor
also concerns a truncated line $L'$ through the origin, so it is easy to see that 
the number $\sigma_\kappa$ that would come from replacing $\ol B(0,1)$ 
with $\ol B(0,\kappa/2)$ is simply $(\kappa/2)^2 \sigma$. 
Thus Lemma \ref{t25.1} will follow if we can prove that
\begin{equation} \label{26.2a}
\sigma \geq C^{-1} \Delta_L = C^{-1}[\H^1(\rho^\ast) - \H^1(X_0 \cap \S)]
\end{equation}
when \eqref{25.38} holds, i.e., when $\Delta_L > 0$.

We shall try various sets $Z'$ and get some lower bounds for $\sigma$; later on
we may proceed by contradiction, assume $\sigma$ is small, and contradict something
in the geometry of $X'$. 

We first study the angle of the two branches of $\rho^\ast$ that leave from some
$m_i$, where $i=1,2$.
Denote by $e_{i,\pm}$ the unit vector tangent at $m_i$ to the branch
$\rho(m_i,z_\pm)$ (or $\rho(m_i,\ell_\pm)$, depending on the situation),
pointing in the direction of the other endpoint of the branch. Then
set $\alpha_i = |e_{i,+} + e_{i,-}|$. It is a good measure of the complement to $\pi$ 
of the angle of $e_{i,+}$ and $e_{i,-}$. We claim that
\begin{equation} \label{26.2}
\sigma \geq C^{-1}  \alpha_i^2.
\end{equation}

This is proved in Lemma 10.23 in \cite{C1}, but let us say how it goes because we shall
use similar proofs soon. On the ball $B= B(m_i/2,1/10)$, the set $X'$ is 
just composed of two half planes, that make the same angle with each other as 
$e_{i,+}$ and $e_{i,-}$. 
Also, $B$ is far from $L'$, so we are not worried by the boundary condition \eqref{25.53a}.
We find a competitor in $B$ that smoothes the angle, where near the middle of $B$ we
essentially move the common boundary of the two half planes by a small fixed vector; on
the rest of $B$ there is a gluing piece, but altogether we save some area. 
Computations are done with the help of the area formula.

\ms
When we have Configuration $3$ near our point $\ell'_\pm$, we claim that we can proceed 
the same way with the two branches of $\rho^\ast$ that leave from $z_\pm$. 
The point is that we do not need to worry about the boundary condition in this case. 

If none of the  two branches of $\rho^\ast$ that leave from $z_\pm$ contain
$\ell_\pm$, we can simply use the statement: any Lipschitz mapping defined on $X'$ 
and that only moves points in $B(z_\pm/2, 1/4)$, say, will satisfy condition \eqref{25.53a} 
because in this ball $X'$ stays far from $L'$. 

But even otherwise, the definition of Configuration 3 makes that we do not need to worry about 
sliding conditions for $X(\rho^\ast)$ or $X'$ near $\ell_\pm$. That is, even though we may detach
$X'$ from $L'$, this does not prevent the competitor $E^1$ that we build with the present construction 
from being a sliding competitor for $E$, because we had no sliding constraint near $\ell_\pm$ 
by definition of Configuration 3.
That is, we should modify the definition of ``good competitor'' to suit Configuration 3, 
but yet we don't need to worry about the estimate.
Now the proof of \eqref{26.2} also yields that
\begin{equation} \label{26.3}
\sigma \geq C^{-1}  \alpha_{\pm,3}^2 \, ,
\end{equation}
where we put the index $3$ to remind the reader of Configuration 3, and where 
$\alpha_{\pm,3} = |v_{\pm,1} + v_{\pm,2}|$, where $v_{\pm,i}$ is the unit vector
that points in the direction of $\rho(z_\pm,m_i)$ at the point $z_\pm$.

\ms
Next we want some control when we have Configuration 1 near $\ell_\pm$. 
Still denote by $v_{\pm,i}$ the unit direction of $\rho(z_\pm,m_i)$ at $z_\pm$, and also let
$v_{\pm,0}$ denote the unit direction of $\rho(z_\pm,\ell_\pm)$ at $z_\pm$.
For some time we will forget the subscript $\pm$ in our notation. Set 
\begin{equation} \label{26.4}
s = v_0 + v_1 + v_2 \, ;
\end{equation}
typically, we want to build competitors for $X'$ by moving the point $z_\pm$ in the 
general direction of $s$, but at the same time we will need to be careful because 
of the boundary constraint along $L'$.

\ms
Let us explain what is our basic competitor. We choose a small 
multiple $v$ of $v_0$ (positive or negative), and we push the points of $X'$ in 
the direction of $v$ (using a cut-off function). For this we repeat the construction of 
Lemma 10.23 in \cite{C1}. 

Let us choose coordinates so that $z = (1,0) \in \R \times \R^{n-1}$, and decide to 
work in the region $A_0 = [1/5,3/4] \times B(0,2a)$, where $a$ is a small geometric constant,
for instance $a = 10^{-2}$. Nothing will happen outside of $A_0$.

Notice that in $A_0$, $X'$ is a truncated set of rough type $\bY$, in the sense that
it is composed of three faces $F_0$, $F_1$ and $F_2$, which are the positive cone over
the three geodesics from $z$ to $\ell$ and to the $m_i$. Only $F_0$ is truncated in $A_0$
(the geodesics $\rho(z,m_i)$ go too far), and we shall consider the half plane $F'_0$
that contain $F_0$, and $X'_1= F'_0 \cup F_1 \cup F_2$, which in $A_0$
coincides with a cone. This cone is not exactly of type $\bY$ because the angles may be 
wrong. Notice however that these angles are not too small either, by the construction of
our nets of curves. What we will do is construct a competitor $Z'_1$ for $X'_1$ in $A_0$,
and later on we will see that we can use the construction to restrict to $X'$ and get a 
competitor $Z'$ for $X'$.

Our competitor for $X'_1$ will be $Z'_1 = f(X'_1)$, where
\begin{equation} \label{26.5}
f(x) = x + \psi(x) v, 
\end{equation}
for some appropriate bump function $\psi$ and a small vector $v$ collinear with $v_0$. 
There is an interest in taking the vector $v$ parallel to $v_0$, which is that with this move, 
the restriction of $f$ to the face $F'_0$ is simpler: if $v$ goes in the direction of $-v_0$,
the face only gets larger (that is, we only add a piece to $F'_0$ in the plane that contains it), and if
$v$ goes in the direction of $v_0$, we just remove a piece of $F'_0$.
For the other faces $F_1$ and $F_2$, they are moved sideways, as in \cite{C1}.

Let us say a little more about $\psi$. 
We take $\psi$ supported in $A_1 = [1/4,1/2] \times B(0,a) \subset \R \times \R^{n-1}$, 
with the same coordinates as above, and as in Lemma 10.23 of \cite{C1}.
We shall mention the other (natural) properties of $\psi$ as we need them.
For the moment, let us not worry too much about the boundary condition, 
and compute the area of $Z' _1\cap A_0$. 
If we choose $v$ small enough (depending on $a$ and our choice of $\psi$,
that we may consider fixed), $f$ is a smooth diffeomorphism (see \cite{C1}), and
we van use the area formula to compute the area of the images $f(F'_0)$, $f(F_1)$, and 
$f(F_2)$ that compose $Z'_1$ in the region $A_0$. 

We proceed as in \cite{C1}, to which we shall refer for some computations.
Let a face, for instance $F_2$, be given. The plane $P_2$ that contains $F_2$
is spanned by $e_1 = (1,0,0) = z$ (where now the third coordinate lies in $\R^{n-2}$),
and, by choice of a suitable basis of $\R^n$, $e_2 = (0,1,0)$. Also write $v = (0,\beta,v')$, with
$v' \in \R^{n-2}$ (or, with a slight twist of notation, $v'$ is orthogonal to $e_1$ and $e_2$).

We need to compute the differential of $f$ on $P_2$, which means 
$Df(e_1) = e_1 + \d_1\psi \, v$ and $Df(e_2) = e_2 + \d_2\psi \, v$.
Here we did not yet write the variables $(x_1,x_2) \in P$, and the notation $\d_1\psi$
is rather clear. But in fact we take $\psi$ to be a function of the first variable $x_1 \in \R$
and $r = (x_2^2 + \cdots x_n^2)^{1/2}$ (i.e., radial in all the other variables), and this way
$\d_2\psi$ will be the same function (of $x_1$ and the other variable) for all the planes $P$ 
that contain the line through $e_1$. Thus
\begin{eqnarray} \label{26.6}
Df(e_1) \wedge Df(e_2) &=& [e_1 + \partial_1\psi \, v] \wedge [e_2 + \partial_2\psi \, v]
\nn\\
&=& e_1 \wedge e_2 + [\partial_2\psi \, e_1 - \partial_1\psi \, e_2] \wedge v
\nn\\
&=& [1+\beta \partial_2\psi] \, e_1 \wedge e_2 + \partial_2\psi\, e_1 \wedge v'
- \partial_1\psi \, e_2 \wedge v'
\end{eqnarray}
and the jacobian determinant of the restriction of $f$ to $F_2$ is
\begin{eqnarray} \label{26.7}
J_2(x) &=& |Df(e_1) \wedge Df(e_2)| 
= \big\{ [1+\beta \partial_2\psi]^2 + (\partial_2\psi)^2 |v'|^2 
+ (\partial_1\psi)^2 |v'|^2 \big\}^{1/2}
\nn\\
&\leq& 1 + \beta \partial_2\psi + C |v|^2
\end{eqnarray}
because $|v|^2 = \beta^2 + |v'|^2$.
Notice also that $\beta = \langle v,e_2\rangle$ is the size of the projection of $v$
on $P_2$; hence, when we apply the area formula to compare ${\H}^2(f(P_2\cap A_0))$ to
$\H^2(P_2\cap A_0)$, we get that
\begin{eqnarray} \label{26.8}
\H^2(f(F_2\cap A_0)) - \H^2(F_2\cap A_0)
&=& \int_{F_2\cap A_1} [J_2(z)-1] \, d\H^2(x)
\nn\\&\,& \hskip-3.5cm
\leq  \int_{F_2\cap A_1} \big[\beta \partial_2\psi(x) 
+ C |v|^2 \big] \, d\H^2(x)
\nn\\&\,& \hskip-3.5cm\leq 
 \langle v,e_2 \rangle \int_{F_2\cap A_1} \partial_2\psi
+ C |v|^2 {\H}^2(F_2\cap A_1)
\leq \langle v,e_2 \rangle \int_{F_2\cap A_1} \partial_2\psi + C |v|^2.
\end{eqnarray}
where we use the fact that $f(P_2) \cap A_0 = f(P_2 \cap A_0)$, and also
that $f(x)=x$ on $\R^n \sm A_1$.

What we computed for $F_2$ is also valid for $F_1$ and $F'_0$, except that we need to
replace the unit vector $e_2$ by a unit vector of $P_1$ or $P_0$ that is perpendicular to 
$e_1 = z$. As was explained before, the derivative $\d_2\varphi$ in that direction is the same,
because we took $\psi$ radial in the directions orthogonal to $e_1$.
Notice also that by rotation invariance, we can use the same coordinates (say, $(x_1,x_2) \in F_2$) 
to write the three integral. We get that
\begin{eqnarray} \label{26.9}
\H^2(f(Z'_1\cap A_0)) - {\H}^2(X'_1\cap A_0)
&=& \sum_{i=0}^2 \H^2(f(F_i\cap A_0)) - \H^2(F_i\cap A_0) 
\nn\\&\,& \hskip-3.5cm \leq C |v|^2 \H^2(F_2\cap A_1)
+ \Big(\sum_{i=0}^2 \langle v,e_i \rangle \Big) \int_{F_2\cap A_1} \partial_2\psi.
\end{eqnarray}
The integral $\int_{F_2\cap A_1} \partial_2\psi = 1$ is a constant, which is even computed
in (10.33) of \cite{C1} to be equal to $-1/5$ (what matters is that it is strictly negative). 
Since $\sum_{i=0}^2 e_i = s$ by \eqref{26.4}, we see that
\begin{equation} \label{26.10}
\H^2(f(Z'_1\cap A_0)) - {\H}^2(X'_1\cap A_0) \leq - \frac{\langle v,s \rangle}{5} + C |v|^2.
\end{equation}
We take 
\begin{equation} \label{26.12a}
v = (10C)^{-1}\langle v_0 ,s \rangle v_0
\end{equation}
with the same constant $C$; then 
$\frac{\langle v,s \rangle}{5} = (50C)^{-1} \langle v_0 ,s \rangle^2$, while
$C |v|^2 = 10^{-2} C^{-1} \langle v_0 ,s\rangle^2$ is twice smaller. Thus
\begin{equation} \label{26.11}
\H^2(f(Z'_1\cap A_0)) - {\H}^2(X'_1\cap A_0) \leq -10^{-2} C^{-1} \langle v_0 ,s\rangle^2.
\end{equation}

This was our version of Lemma 10.23 of \cite{C1}. Notice that the estimate
is not so good when $s$ is almost orthogonal to $v_0$, but let us keep the option to use
this open. Now we need to worry about the set $X'$ and the boundary condition.

\begin{figure}[!h]  
\centering
\hskip0.7cm\includegraphics[width=5cm]{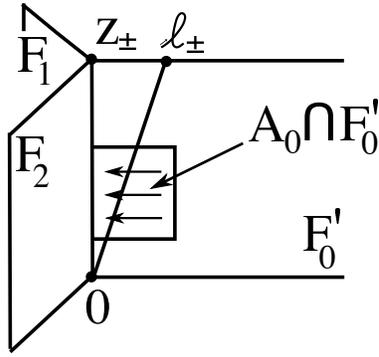}
\caption{A picture in $P_t$ 
\label{f26.1}}
\end{figure}

Let us review how the mapping $f$ works; see Figure \ref{f26.1} already. 
On the two faces $F_1$ and $F_2$, the mapping pushes points in the direction of $v$, 
and the only case when the boundary condition \eqref{25.53a} may be violated is if 
$F_1$ or $F_2$ touches $L'$. Due to the fact that
$\rho_0$ makes a large angle with $\rho_1$ and $\rho_2$, this can only happen if
$\rho_0$ is reduced to one point and $\rho_1$ and $\rho_2$ start from $\ell_\pm$. 
Let us assume that this does not happen for the moment.

Then on $F_0$, the mapping $f$ slides points in the direction of $v$, which is parallel
to the plane $P_0$ that contains $F_0$. Let us start the discussion with the case 
when $v$ is a negative multiple of $v_0$. The effect of $f$ is to extend the faces $F_0$
and $F'_0$, by adding to them a piece that lies further than the boundary $[0,z_\pm]$
(where $z_\pm$ is the common point of the three $\rho_i$). With the way we wrote
$f$, we probably moved points of $F'_0$ that lie on $[0,\ell_\pm] \subset L'$ and beyond,
and \eqref{25.53a} forbids us to do this. 
But this is easy to fix: we replace $f$ on $F_0'$ by a mapping that coincides with $f$ on 
$[0,z_\pm]$ (so that we can still glue with $f_{\vert (F_1 \cup F_2)}$, 
is the identity on $F'_0 \sm F_0$ (and in particular on $[0,\ell_\pm]$), 
and just moves the points faster in $A_0 \cap F_0$ if needed.
This shows that $Z'= F'_1 \cup F'_2 \cup F'_3$ is a good competitor for $X'$,
and since $Z'_1 \sm Z' = X'_1 \sm X' = F'_0 \sm F_0$, we deduce from \eqref{26.11}
that 
\begin{eqnarray} \label{26.14a}
\sigma &\geq&
{\H}^2(X'\cap A_0) - \H^2(f(Z'\cap A_0)) 
\nn\\
&=& {\H}^2(X'_1\cap A_0) - \H^2(f(Z'_1\cap A_0)) 
\geq 10^{-2} C^{-1} \langle v_0 ,s\rangle^2.
\end{eqnarray}
This was our estimate when $v$ is a negative multiple of $v_0$, 
which by \eqref{26.12a} means that $\langle v_0 ,s \rangle < 0$.

When $\langle v_0 ,s \rangle \geq 0$, the mapping $f$ tends to move points of $F'_0$
in the direction of $v_0$, i.e., make $F'_0$ smaller. We can still argue as before, but under
the condition that the points of the common boundary $[0,z_\pm]$ of $F_0$ and $F'_0$
do not go all the way to $[0,\ell_\pm]$. This means that when we choose $v$, we are safe
if $|v| \leq 10^{-3} |z_\pm - \ell_\pm| = 10^{-3}|z-\ell|$, for instance. 
We may as well assume that the constant $C$ in \eqref{26.10} is larger than $100$, and this way, 
if $0 \leq \langle v_0 ,s \rangle \leq |z-\ell|$, we can keep the same choice of $v$ as in 
\eqref{26.12a}, and get the same estimate as in \eqref{26.14a}. Altogether,
\begin{equation} \label{26.13}
\sigma \geq C^{-1}  \langle v_0 ,s\rangle^2 \ \text{ when } \langle v_0 ,s \rangle \leq |z-\ell|.
\end{equation}
When $\langle v_0 ,s \rangle \geq |z-\ell|$, we take the smaller 
$v = (10C)^{-1} |z-\ell | \, v_0$, and the same computation as for \eqref{26.11}
and \eqref{26.14a} yields the less good result
\begin{equation} \label{26.14}
\sigma \geq C^{-1}  |z-\ell| \langle v_0 ,s\rangle \ \text{ when } \langle v_0 ,s \rangle \geq |z-\ell|.
\end{equation} 

These will be our main estimates, but there are some cases when 
$\langle v_0 ,s \rangle$ is really too small for us, or (almost equivalently) the author
did not manage to prove easily that it is large, and we want to try a different
competitor. We shall try this when $v_1$ and $v_2$ make a small angle, and more precisely
$|v_1+v_2| > 1$. We could of course try to control the scalar product above when this happens,
but the author did not manage to do this, and instead we shall try a slightly different competitor,
where we move the points of the faces $F_1$ and $F_2$ above in the direction of $v_1+v_2$
(instead of $-v_0$ above). We need to be more specific, because we want to use the same 
computations as above, but not the same mapping. Suppose we keep $F_1$ and $F_2$ as they are, 
but complete them with a third face $F_3$, starting from their common boundary $[0,z]$, and
going in the opposite direction $v'_0 = - (v_1+v_2)/|v_1+v_2|$. Then use the same algorithm
as before, where $f$ is given by \eqref{26.5}, with for $v$ a positive multiple of $v_1+v_2$.
We want to do the same computations as above, with $v_0$ replaced by $v'_0$, and hence
$s$ replaced by $s' = v'_0 + v_1 + v_2 = - v'_0 (|v_1+v_2|-1)$. In particular, we take
$v= (10C)^{-1}\langle v'_0 ,s' \rangle v'_0 = - (10C)^{-1} (|v_1+v_2|-1) v'_0$ as in \eqref{26.12a}.

This gives a competitor $\wt Z$ for $F_1 \cup F_2 \cup F_3$, for which the estimate leading
to \eqref{26.11} are true. We remove the whole face 
$F_3$, both from $F_1 \cup F_2 \cup F_3$ and from $\wt Z$, and we get a new set
$Z''$ such that 
\begin{equation} \label{26.17a}
\H^2(Z'' \cap A_0)) - \H^2((F_1 \cup F_2) \cap A_0) \leq - 10^{-2}C^{-1} \langle v'_0 ,s'\rangle^2
= - 10^{-2}C^{-1} (|v_1+v_2|-1)_+^2,
\end{equation}
where we added the positive part to remember that we do this only when $|v_1+v_2| > 1$.
Now in $A_0$, $Z''$ is composed of slightly distorted faces $F'_1 = f(F_1)$ and $F'_2 = f(F_2)$,
plus a vaguely triangular piece of $f(F_3)$, which is bounded by a piece of the common boundary
$f([0,z])$ of $F'_1$ and $F'_2$ on one side, and the corresponding piece of $[0,z]$
on the other side. We add to $Z''$ (and in $A_0$ only; since we did not change anything outside 
of $A_0$) the old face $F_0$ (bounded by $[0,z]$, $[0,\ell]$, and the arc $\rho_0$)
 and get a set $Z'$, which is a competitor for $X'$ (which in $A_0$ coincides with the cone over 
 $\rho_0 \cup \rho_1 \cup \rho_2$). It follows from \eqref{26.17a} that
\begin{equation} \label{26.18a}
\H^2(Z' \cap A_0)) - \H^2(X' \cap A_0) \leq  - 10^{-2} C^{-1} (|v_1+v_2|-1)_+^2.
\end{equation}
Now we claim that $Z'$ is a good competitor for $X'$. We do not want to use the same mapping
$f$ as for proving the estimates, but instead observe that $Z' = \wt f(X')$, for some mapping
$\wt f$ that pinches partially the two faces $F_1$ and $F_2$ (in the direction orthogonal
to the plane of $F_3$), is Lipschitz, but will not be written explicitly here. We need to make sure
that $\wt f$ can be extended by setting $\wt f(x) = x$ on $L'$, because of our boundary
constraint, and this follows from the angle condition \eqref{9.2}, which says that
$\rho_1$ and $\rho_2$ make large angles with $\rho_0$ at the point $z$,
so that in the present situation where $\rho_0$ is a nontrivial arc, we only move points
that are far from $L'$. Thus $Z'$ is a good competitor and \eqref{26.18a} can be used
to prove that
\begin{equation} \label{26.19a}
\sigma \geq  C^{-1} (|v_1+v_2|-1)_+^2.
\end{equation}
In fact, we claim that the present estimate also works when $\rho_0$ is reduced to $\{ \ell \}$.
In this case, the two vectors $v_1$ and $v_2$ (the direction of $\rho_1$ and $\rho_2$ when
they leave $z = \ell$) are still well defined, we can define $Z'$ as above, and the fact that it
is a good competitor for $X'$ comes from the fact that the pinching mapping $f$ can be taken to be
the identity on $[0,\ell]$. The claim follows.

There is a last estimate on $\sigma$ that we may use, essentially when all the other ones fail,
which says that
\begin{equation} \label{26.15}
\sigma \geq C^{-1}|z-\ell| |s|^2.
\end{equation}

This estimate seems less good, because the right-hand side is of order $3$.
Its proof also relies on (the proof of) Lemma 10.23 of \cite{C1}.
We use the fact that we can find a tube of width roughly equal to $|z-\ell|$, 
centered on the segment $[z/3,2z/3]$, that does not meet $L$, and where $X'$ 
coincides with a cone which is roughly of type $\bY$, except that at least one of its angles 
is off by roughly $|s|$. We can apply the proof of Lemma 10.23 in \cite{C1} to get an estimate.
With the same value of $s$, and in the unit ball, we would save $C^{-1}|s|^2$;
in a ball of size roughly $|z-\ell|$, and by homogeneity, it would be $C^{-1} |z-\ell|^2|s|^2$. 
But here we are in a thin tube of roughly unit length, and the proof of \cite{C1} allows
us to save $C^{-1} |z-\ell| |s|^2$. This gives the quite general, but not so good estimate 
\eqref{26.15}.

\ms
We shall now start distinguishing between cases. To make our life easier (at least,
in the cases that will be settled in this section), let us decide that the two points $m_i$ 
(that were selected on the curves $\cC_i$, near the $w_i$) 
are chosen at equal distance from $\ell_+$ and $\ell_-$. This is easy to arrange, 
as in the case of a half plane, by the intermediate value theorem and because each 
$\cC_i$ is nicely transverse to the plane equidistant from the $\ell_\pm$.

\msi
{\bf Case 1.}
We start with the case when we have Configuration $1$ near each point of $L\cap \S$
(recall that $\S = \d B(0,1)$), and in addition $X$ is flat, by which we mean that 
either $X \in \bP_0$ or else $X \in \bV$ and for some $\delta > 0$,
\begin{equation} \label{26.16}
\text{ the two half planes that compose $X$ make an angle larger than } \frac{2\pi}{3}+\delta.
\end{equation}

Let us explain why this case is simpler. Assume first that $X \in \bV$ and \eqref{26.16} holds; 
if $\varepsilon$ is small enough, then $|z-\ell|$ also is as small as we want, then the two directions 
$v_1$ and $v_2$ are very close to the directions of the two half planes of \eqref{26.16} 
(understand,  the unit vectors perpendicular to the direction of $L$ that are tangent to these planes 
and go away from $L$). Then $|v_1 + v_2| \leq 1 - \delta/3$, say, by \eqref{26.16}, and 
\begin{equation} \label{26.17}
\langle s, v_0 \rangle = \langle v_0 +v_1+v_2, v_0 \rangle  =
1 + \langle v_1+v_2, v_0 \rangle  \geq \delta/3.
\end{equation}
When $X \in \bP_0$, the situation is even better: the two geodesics of $X$
that start from $\ell'_\pm$ go in opposite directions, and since $z = z_\pm$ lies close to
$\ell'_\pm$ and $w_i$ lies close to $m_i$, we get that $|v_1 + v_2| < 1/2$
and \eqref{26.17} holds as well.

Most probably, $|z-\ell| \leq \langle s, v_0 \rangle$, and then \eqref{26.14} says that
\begin{equation} \label{26.18}
\sigma \geq C^{-1}  |z-\ell| \langle v_0 ,s\rangle \geq C(\delta)^{-1} |z-\ell|.
\end{equation}
But even if $|z-\ell| \geq \langle s, v_0 \rangle$, we may apply \eqref{26.13}
instead and get that $\sigma \geq C^{-1} \delta^2$, which is better than \eqref{26.18}.
This estimate holds near both points $\ell_\pm$; we use this to majorize
\begin{equation} \label{26.19}
\begin{aligned}
\H^1(\rho^\ast) &= \ddist(z_+,\ell_+) + \ddist(z_-,\ell_-)
+ \sum_{i=1}^2 (\ddist(z_+,m_i)+\ddist(z_-,m_i))
\cr&
\leq C(\delta) \sigma + \sum_{i=1}^2 (\ddist(\ell_+,m_i)+\ddist(\ell_-,m_i)),
\end{aligned}
\end{equation}
where the first part is just the definition of $\rho^\ast$ as a concatenation of geodesics.

\begin{lem} \label{t26.1} Denote by $H$ the set of points that lies at equal distance from 
$\ell_+$ and $\ell_-$. For all choices of $m_1, m_2 \in \S \cap H$,
\begin{equation} \label{26.20}
 \sum_{i=1}^2 (\ddist(\ell_+,m_i)+\ddist(\ell_-,m_i)) \leq \H^1(X_0 \cap \S).
 \end{equation}
\end{lem}

We shall prove the lemma soon, but let us see how it implies Lemma \ref{t25.1}
in the present case. We deduce from \eqref{26.19} and \eqref{26.20}
that $\H^1(\rho^\ast) \leq C(\delta) \sigma + \H^1(X_0 \cap \S)$ and, if 
$\Delta_L = \H^1(\rho^\ast)-\H^1(X_0 \cap \S)$ (the same as in \eqref{25.47}) 
is positive, this just means that $\sigma \geq C(\delta)^{-1} \Delta_L$. We use the definition
\eqref{26.1} of $\sigma$ and get the conclusion of Lemma \ref{t25.1}.

\ms
So let us prove Lemma \ref{t26.1}, and our Case 1 will follow. Set 
\begin{equation}\label{26.21}
D_\pm = \ddist(\ell_\pm,m_1)+\ddist(\ell_\pm,m_2) ;
\end{equation}
then \eqref{26.20} will follow at once if we prove that for each sign,
\begin{equation}\label{26.22}
D_\pm \leq \frac{1}{2} \H^1(X_0 \cap \S).
\end{equation}
Let us prove this. Fix a sign $\pm$, and drop it from the notation. 
That is, we set $\ell = \ell_\pm$ and $D = D_\pm$.

We need to compute a few things. We start with the relation between the geodesic and 
Euclidean distances on the sphere. We claim that for $x,y \in \S$,
\begin{equation}\label{26.23}
2-2\cos(\ddist(x,y)) = |x-y|^2
\end{equation}
For this computation we may assume that $y,z \in \R^2$, and even that $x = (1,0)$ and
$y = (\cos\theta, \sin\theta)$ for some $\theta \in [0,\pi]$; in this case $\ddist(x,y) = \theta$
and $|x-y|^2 = (1-\cos\theta)^2 + \sin^2\theta = 2(1-\cos\theta)$; the claim follows.
Notice also that $\ddist(x,y) \in [0,\pi]$ and $|x-y|^2 \in [0,4]$, so $1- \frac{1}{2} |x-y|^2 \in [-1,1]$, 
and \eqref{26.23} is the same as
\begin{equation}\label{26.24}
\ddist(x,y) = \arccos\big(1- \frac{1}{2} |x-y|^2 \big).
\end{equation}
Next we compute numbers like $\ddist(\ell,m)$, where $\ell \in L \cap \S$ and 
$m$ lies in the median hyperplane $H$.
Without loss of generality, we may assume that there are three orthogonal unit vectors
$e_1, e_2, e_3$ such that 
\begin{equation}\label{26.25}
L = \big\{ -d_0 e_2 + t e_1 \, ; \, t\in \R \big\} \ \text{ and }\ 
m= \sin\alpha \, e_3 - \cos\alpha \, e_2 \text{ for some } \alpha \in [0,\pi].
\end{equation}
Thus $\alpha > 0$ small corresponds to a point $m$ just above the shade of $L$ (or if you prefer $-e_2$), 
$\alpha=\pi$ corresponds to an $m$ just opposite to the shade; we decided that we did not need 
the case when $\alpha \in (\pi,2\pi)$ by symmetry.
We may also assume that 
\begin{equation}\label{26.26}
\ell = -d_0 e_2 + \sqrt{1-d_0^2}e_1
\end{equation}
(the other choice $\ell = -d_0 e_2 - \sqrt{1-d_0^2}e_1$ would be equivalent), and then
\begin{equation}\label{26.27}
|m-\ell|^2 = (1-d_0^2) + (d_0-\cos\alpha)^2+ \sin^2\alpha 
= 2 - 2 d_0 \cos\alpha;
\end{equation}
thus by \eqref{26.23} or \eqref{26.24},
\begin{equation}\label{26.28}
\ddist(\ell,m) = \arccos\big(d_0 \cos\alpha \big).
\end{equation}
Since $d_0$ is small, we see that $\ddist(\ell,m)$ is close to $\pi/2$.
Notice that $\ddist(\ell,m)$ is a nondecreasing function of $\alpha \in [0,\pi]$.

Let us return to the two points $m_i \in \S \cap H$, which we write
$m_i= \sin\alpha_i e_{3,i} - \cos\alpha_i e_2$ as in \eqref{26.25}, with possibly different
vectors $e_3 = e_{3,i}$ if we work in $\R^n$, $n > 3$. We want to estimate
\begin{equation}\label{26.29}
D = \ddist(\ell,m_1)+\ddist(\ell,m_2)  = \arccos\big(d_0 \cos\alpha_1 \big)
+ \arccos\big(d_0 \cos\alpha_2 \big).
\end{equation}
This is again a nondecreasing function of $\alpha_1$ and $\alpha_2$. 
We also need to evaluate the angles $\alpha_1$ and $\alpha_2$ in terms of 
the geometry of $X$. Start when $X \in \bV$, and denote by 
$\Angle(X)$ the angle that the two half planes that compose $X$ make along $L$;
thus $\Angle(X) > \frac{2 \pi}{3} +\delta$ by \eqref{26.7}. Call
$\Angle_{x_0}(m_1,m_2)$ the angle of $m_1$ and $m_2$, seen from $x_0 = -d_0 e_2$
(the midpoint of $[\ell_+,\ell_-]$).
Since both $m_i$ lie within $2\varepsilon$ from $X$ (by \eqref{25.2}), we get that
$|\Angle_{x_0}(m_1,m_2)-\Angle(X)| \leq 5\varepsilon$ hence (by \eqref{26.16} and if 
$\varepsilon$ is small enough, depending on $\delta$)
\begin{equation}\label{26.30}
\Angle_{x_0}(m_1,m_2) > \frac{2 \pi}{3} + \frac{\delta}{2}.
\end{equation}
When $X \in \bP_0$, $\Angle_{x_0}(m_1,m_2)$ is almost $\pi$, because $x_0$ is not so far
from $0$, and the two points $w_i$ almost lie in opposite directions.
In both cases, \eqref{26.30} holds, and so
\begin{equation}\label{26.31}
\frac{\langle m_1-x_0, m_2 -x_0\rangle }{ |m_1-x_0| |m_2-x_0|}
= \cos(\Angle_{x_0}(m_1,m_2))
< - \frac{1}{2} - \frac{\delta}{4}.
\end{equation}
Notice that
\begin{eqnarray}\label{26.32}
\langle m_1-x_0, m_2 -x_0 \rangle 
&=& \langle\sin\alpha_1 e_{3,1} - \cos\alpha_1 e_2+d_0 e_2, \sin\alpha_2 e_{3,2} 
- \cos\alpha_2 e_2+d_0 e_2\rangle
\nn\\
&=&\sin\alpha_1 \sin\alpha_2 \langle e_{3,1},e_{3,2} \rangle
+ \cos\alpha_1\cos\alpha_2
- d_0 (\cos\alpha_1 +\cos\alpha_2) + d_0^2
\nn\\
&\geq& - \sin\alpha_1 \sin\alpha_2 + \cos\alpha_1\cos\alpha_2 - 3d_0
= \cos(\alpha_1+\alpha_2) - 3d_0
\end{eqnarray}
because $\sin\alpha_i \geq 0$. Let us take $N \geq 100/\delta$ in \eqref{25.4}, so that
$d_0 \leq N^{-1} \leq \delta/100$. Notice that $\big| |m_i-x_0|-1\big| \leq d_0 \leq \delta/100$ for $i=1,2$, so by \eqref{26.32} and \eqref{26.31}
\begin{equation}\label{26.33}
\cos(\alpha_1+\alpha_2) 
\leq \langle m_1-x_0, m_2 -x_0\rangle + 3d_0 
\leq \Big[- \frac{1}{2} - \frac{\delta}{4}\Big] |m_1-x_0| |m_2-x_0| + \frac{3\delta}{100}
\leq - \frac{1}{2} - \frac{\delta}{100},
\end{equation}
hence (since $0 \leq \alpha_1 + \alpha_2 \leq 2\pi$)
\begin{equation}\label{26.34}
\frac{2 \pi }{ 3} + \frac{\delta}{200} \leq \alpha_1 + \alpha_2 
\leq \frac{4 \pi}{3} - \frac{\delta}{200}.
\end{equation}

Return to $D$ in \eqref{26.21} and \eqref{26.29}. A Taylor expansion of order $2$ near $0$ yields
$\arccos(d_0\cos\alpha_i) = \frac{\pi}{2} - d_0 \cos\alpha_i + O_i$, with $|O_i| \leq d_0^2/2$, then
\begin{equation}\label{26.35}
\begin{aligned}
D &= \pi - d_0(\cos\alpha_1 + \cos\alpha_2) + O_1+O_2 
\leq \pi - d_0(\cos\alpha_1 + \cos\alpha_2) +d_0^2 
\cr&=\pi - 2d_0(\cos((\alpha_1 + \alpha_2)/2) \cos((\alpha_1 - \alpha_2)/2)+d_0^2.
\end{aligned}
\end{equation}
If both cosines have the same sign, this yields $D \leq \pi + d_0^2$ and we'll be happier
than in the next case. Otherwise, since both $0 \leq \alpha_i \leq \pi$ and hence 
$\cos((\alpha_1 - \alpha_2)/2) \geq 0$, we get that $\cos((\alpha_1 + \alpha_2)/2)< 0$,
hence $(\alpha_1 + \alpha_2)/2 \geq \pi/2$ and by \eqref{26.34}
\begin{equation}\label{26.36}
\frac{\pi}{2} \leq \frac{\alpha_1 + \alpha_2 }{2} \leq \frac{2 \pi}{3} - \frac{\delta}{400},
\end{equation}
\begin{equation}\label{26.37}
0 \geq \cos((\alpha_1 + \alpha_2)/2) \geq \cos\Big(\frac{2 \pi}{3} - \frac{\delta}{400}\Big)
\geq - \frac{1}{2} + \frac{\delta}{800},
\end{equation}
and 
\begin{equation} \label{26.38}
D \leq \pi + 2d_0 |\cos((\alpha_1 + \alpha_2)/2)| + d_0^2
\leq \pi + d_0 - \frac{d_0 \delta}{400} + d_0^2 \leq \pi + d_0 - \frac{d_0 \delta}{800}
\end{equation}
if $N$ is large enough.
We also get this in the other case when the cosines have the same sign.
We need to compare this with the right-hand side of \eqref{26.22}, so we 
compute $\H^1(X_0 \cap \S)$.
Recall that $X_0$ is composed of two half great circles, that end at two antipodal 
points $y_\pm$, plus the two short arcs of geodesics $\rho(\ell_\pm,y_\pm)$. 
The two half circles account for $2\pi$,
and with the same choice of basis as for \eqref{26.25}, $y_\pm = \pm e_1$
(because the spine of $X_0$ is parallel to $L$).
Recall from \eqref{26.26} that $\ell_\pm = -d_0 e_2 \pm \sqrt{1-d_0^2}e_1$, hence 
$\ddist(\ell_\pm,y_\pm) = \arcsin d_0 \geq d_0$, and $\H^1(X_0 \cap \S) \geq 2\pi + 2d_0$. 
This completes our proof of \eqref{26.22}, Lemma \ref{t26.1} follows, and we get the desired 
estimate for Lemma \ref{t25.1} in our Case 1.
\qed

\msi
{\bf Case 2.}
Next we assume that we have Configuration $3$ near the two points of $L\cap \S$,
regardless of whether \eqref{26.16} holds or not.
In this case $\rho^\ast$ is merely composed of four curves from the two $m_i$ to the
two $z_\pm$. The complement to $\pi$ of the angles at the $w_i$ are less than
$C \sqrt\sigma$, by \eqref{26.2}, and the the angles at the $z_\pm$ are also less than
$C \sqrt\sigma$, by \eqref{26.3}. Notice that we may assume that $\sigma$ is small, 
because otherwise the conclusion of Lemma \ref{t25.1} is obvious. Then the four vertices of
$\rho^\ast$ lie at distance at most $C \sqrt\sigma$ from some great circle
(we can follow the curve $\rho^\ast$ from $z_-$ back to $z_-$, without turning away from
the geodesic more than $C\sqrt\sigma$), and by standard 
computations (that can be found in \cite{C1}, for instance), $\H^1(\rho^\ast) \leq 2\pi + C \sigma$.
This is better than what we need for Lemma \ref{t25.1}, because $\H^1(X_0 \cap \S) \geq 2\pi$.

\msi
{\bf Case 3.}
Suppose now that we have Configuration $1$ near $\ell_-$ and Configuration $3$ near $\ell_+$.
We still have that $\alpha^2_{3,+} \leq C \sqrt\sigma$, by \eqref{26.3}, 
and $\alpha_i \leq C \sqrt\sigma$, by \eqref{26.2}. With the same reasoning as above, 
all the vertices $m_i$ and $z_\pm$ lie within $C\sqrt\sigma$ of a great circle, and then
\begin{equation}\label{26.39}
\sum_{i=1}^2 (\ddist(m_i,z_+) + \ddist(m_i,z_-)) \leq 2\pi + C \sigma.
\end{equation}
But this time we also have to account for the additional small piece $\rho(\ell_-,z_-)$. 
Since $\alpha^2_{3,+} \leq C \sqrt\sigma$ and we can assume that $\sigma$ 
is very small (because otherwise the thesis is trivial), the two half planes that 
compose $X$ when $X \in \bV$ make an angle $\Angle(X) \geq 9\pi/10$;
when $X \in \bP_0$, they are even in front of each other.
The sum $s = v_0+v_1+v_2$ of \eqref{26.4} (and for the point $z_-$) 
is then quite close to $v_0$ (because $v_1+v_2$ is small), 
so $\langle v_0,s \rangle \geq 1/2$, which is better than \eqref{26.17}. 
As for \eqref{26.18}, we also get that $\sigma \geq C^{-1} |z_--\ell_-|$,
which gives a good enough estimate for $\H^1(\rho(\ell_-,z_-)) = \ddist(\ell_-,z_-)$.
We add this to \eqref{26.39} and get that $\H^1\rho^\ast) \leq 2 \pi + C \sigma
\leq \H^1(X_0 \cap \S) + C \sigma$, as needed.

Recall that we excluded the case of Configuration 2 earlier. At this stage, we have only
one case left, which is when we have Configuration $1$ near both $\ell_\pm$, 
and in addition $X \notin \bP_0$ and \eqref{26.16} fails, i.e., $X \in \bV(L)$ and
\begin{equation}\label{26.40}
\frac{2\pi}{3} \leq \Angle(X) \leq \frac{2\pi}{3} + \delta,
\end{equation}
where the first part comes from our assumption that $X\in \bV(L)$.
Recall also that for this remaining case we are allowed to take $\delta > 0$ as small 
as we want.

\section{Full length for sharp $\bV$ sets}
\label{S27}

In this section we study the last left case for Lemma \ref{t25.1}, when we have 
Configuration $1$ near both $\ell_\pm$, and in addition $X$ satisfies \eqref{26.40}.
We talk about sharp $\bV$ sets because we could even argue that in this remaining
situation, since we have good approximation by a set $X \in \bV$ such that \eqref{26.40}
holds, and in addition we can take $\delta$ small, we are left with the case where we have
a reasonably good approximation by a set $X \in \bV$ with dihedral angle exactly $2\pi/3$.
We shall not try proceed like this, because it would not really help simplify the proof, 
and also we would at least need to be quite careful with the quantifiers. 
Our last case is somewhat more complicated than the other ones, which is why we left it for the end.

We shall keep some of the notation of the previous cases, concerning the two points
$z = z_\pm$ near the vertices $\ell = \ell_\pm$, and two intermediate points $m_1$ and $m_2$,
except that we may not choose the two $m_i$ exactly as we did in the previous section
(that is, at equal distance from $\ell_+$ and $\ell_-$).
We shall first try to estimate the length
\begin{equation} \label{27.1}
L_{12} = \sum_{i=1}^2 \ddist(z_+,m_i) + \ddist(z_-,m_i),
\end{equation}
but for this some additional notation will be useful. We shall think of $z=z_-$ as
the lowest point of $\S$, and will project various things along the unit vector
\begin{equation} \label{27.2}
e_0 = (v_1(z)+v_2(z))/|v_1(z)+v_2(z)|,
\end{equation}
where for $i =1, 2$, $v_i(z)$ is the direction of $\rho(z,m_i)$ at the point $z$.
Let us also write
\begin{equation} \label{27.3}
\Angle_{z}(v_1(z),v_2(z)) = \frac{2\pi}{3} + 2\alpha,
\end{equation}
where we know that $\alpha$ is small because $z=z_-$ lies close to $\ell_-$,
the $m_i$ lie close to $X$, and we have \eqref{26.40}, but it will be useful later
to have a more precise estimate.

When we continue the two geodesics $\rho(z,m_i)$ past the points $m_i$, they eventually
meet back at the point $-z$, with the same angle $\frac{2\pi}{3} + \alpha$. But at the 
point $m_i$, we turned a little and used the geodesic $\rho(m_i,z_+)$ instead. Notice however that
we turned by at most $\alpha_i \leq C \sqrt\sigma$, by \eqref{26.2}, and because of this
the new meeting point $z_+$ (we know it exists!) does not move by more than $C \sqrt\sigma$. 
That is,
\begin{equation} \label{27.4}
|z_+ + z| \leq C \sqrt\sigma.
\end{equation}
Set $f_i(w) = \ddist(w,m_i)$; notice that $f_i(z)+f_i(-z) = \pi$ because the union of the two
corresponding geodesics is a half great circle. We want to evaluate $f(z_+)$ by estimating
the derivative of $f_i$ near $-z$. It is easy to see that at a given point $x \in \S$, $x\neq w_i$, 
the derivative of $f_i$ at $x$ in the direction $e$ is
\begin{equation} \label{27.5}
Df_i(x)(e) = \d_e f_i(x) = - \langle e,v_i(x) \rangle,
\end{equation}
where again $v_i(x)$ is the direction of $\rho(x,w_i)$ at $x$. Moreover, if $x$ is any
point of $\rho(-z,z_+)$, $|v_i(x)-v_i(-z)| \leq 10|z_++z| \leq C \sqrt\sigma$. 
In addition, $v_i(-z) = v_i(z)$ because these are the endpoint directions of the half circle from
$z$ to $-z$ through $w_i$, so $|v_i(x)-v_i(z)| \leq C \sqrt\sigma$. 
Now write $v(x,z_+)$ the direction of $\rho(-z,z_+)$ at $x$; then
$$
f_i(z_+) - f_i(-z) = \int_{\rho(-z,z_+)} Df_i(x)(v(x,z_+)) 
= -  \int_{\rho(-z,z_+)} \langle v(x,z_+),v_i(x) \rangle.
$$
The length of the geodesic is at most $C \sqrt\sigma$, by \eqref{27.4}.
When we replace $v_i(x)$ by $v_i(z)$, we make an error of at most $C\sqrt\sigma$.
When we replace $v(x,z_+)$ by $(z_+ +z)/|z_+ +z|$, we also make an error of at most $C\sqrt\sigma$
(the geodesic does not turn much). We integrate the error and get at most $C\sigma$.
Replacing the length of $\rho(-z,z_+)$ with $|z_+ +z|$ also generates an error of at most $C\sigma$;
altogether 
\begin{equation} \label{27.6}
\big|f_i(z_+) - f_i(-z) + \langle z_+ +z, v_i(z)\rangle \big| \leq C \sigma.
\end{equation}
We sum over $i$, observe that if $z_+$ were equal to $-z$ we would have $L_{12}=2\pi$,
and get that 
\begin{eqnarray} \label{27.7} 
\big| L_{12} - 2\pi + \langle z_+ +z, v_1(z)+v_2(z) \rangle\big|
&=& \nn\\ &\,& \hskip-5cm
\big| f_1(z_+) + f_2(z_+) - f_1(-z) - f_2(-z) - 2\pi + \langle z_+ +z, v_1(z)+v_2(z) \rangle\big|
\leq C \sigma.
\end{eqnarray}
Observe that by \eqref{27.2} and \eqref{27.3}
\begin{equation} \label{27.8}
v_1(z)+v_2(z) - e_0 = e_0 |v_1(z)+v_2(z)| - e_0
= 2e_0 \cos\big(\frac{\pi}{3}+\alpha\big) - e_0
= -\wt \alpha e_0,
\end{equation}
where $\wt \alpha = 2 \big(1-\cos\big(\frac{\pi}{3}+\alpha\big)\big)$ is of the order
of $\sqrt 2 \alpha$, but the precise value will not be so important. Simply notice that by 
\eqref{27.7} and then \eqref{27.4},
\begin{eqnarray} \label{27.9}
L_{12} - 2\pi + \langle z_+ +z, e_0 \rangle
&=& L_{12} - 2\pi + \langle z_+ +z, v_1(z)+v_2(z) \rangle + \langle z_+ +z, e_0 - v_1(z)-v_2(z) \rangle
\nn\\
&\leq&  C \sigma + \wt\alpha \langle z_+ +z, e_0 \rangle
\leq  C \sigma + C \alpha \sqrt\sigma.
\end{eqnarray}
Call $L_\pm = \ddist(\ell_\pm,z_\pm)$ the lengths of our two remaining short arcs
$\rho_\pm = \rho(\ell_\pm,z_\pm)$; then the decomposition of $\rho^\ast$ yields
\begin{equation} \label{27.10}
\H^1(\rho^\ast) = L_{12}+ L_- + L_+.
\end{equation}
 Observe also that
\begin{equation} \label{27.10b}
\H^1(X_0 \cap \S) = 3\pi - \ddist(\ell_-,\ell_+) = 2\pi + \ddist(\ell_-,-\ell_+),
\end{equation}
because when we add $\rho(\ell_-,\ell_+) = S \cap \S$ to $X_0 \cap \S$, we get a union of three half 
great circles. Recall that we want to show that $\Delta_L \leq C \sigma$, as in \eqref{26.2a},
where by \eqref{27.9}
\begin{eqnarray} \label{27.12a}
\Delta_L &=& \H^1(\rho^\ast) - \H^1(X_0 \cap \S)
= L_{12}+ L_- + L_+ 
-2\pi - \ddist(\ell_-,-\ell_+)
\nn\\
&\leq& L_- + L_+ - \langle z_+ +z, e_0 \rangle 
- \ddist(\ell_-,-\ell_+) 
+ C \sigma + \wt\alpha \langle z_+ +z, e_0 \rangle 
\nn\\
&\leq& L_- + L_+ - \langle z_+ +z, e_0 \rangle - \ddist(\ell_-,-\ell_+) 
+ C \sigma + C \alpha \sqrt\sigma.
\end{eqnarray}

The estimate that we want to do now looks like the following. 
Imagine that there is no curvature in the sphere and that the three geodesics $\rho_- =\rho(\ell_-,z_-)$, 
$\rho(z_-,-z_+)$ and $\rho(-z_+,-\ell_+) = -\rho_+$ are all contained in a line parallel to $e_0$.
Then $\ddist(\ell_-,-\ell_+) = L_- - \langle z+z_+, e_0 \rangle + L_+$, where
the middle term may be positive or negative, but in all cases we get that
$\Delta_L \leq C \sigma + C \alpha \sqrt\sigma$. 
We would still need a good estimate on $\alpha$, but would get close to the desired goal.
In the mean time we will try to deal with the curvature of $\S$ and the
alignment of our three geodesics.

Let $\rho_0$ denote the geodesic that contains $z$ and points in the direction $\pm e_0$
at $z$; we want to project all sort of points on $\rho_0$, and then try to follow the sketchy 
argument above. Denote by $\ell'_-$ the point of $\rho_0$ that lie closest to $\ell_-$; 
we earlier used this notation for a point like $\ell_-$ when $X\in \bP_0$, but there is no relation. 
We want to locate $\ell'_-$ more precisely. Denote by $\beta$ the angle of $v_0(z)$ 
(the direction of $\rho(z,\ell_-)$) with $-e_0$. 
Simple estimates (that we do not do because we will do more precise computations below) 
show that since $L_-$ is quite small,
\begin{equation} \label{27.11}
|\ell'_- -\ell_-| \leq 2 \beta L_- \ \text{ and } \ |\ddist(\ell'_-,z) - L_-| \leq C \beta^2 L_-.
\end{equation}
In addition, $\beta \leq C|s(z)|$, where $s(z) = v_0(z) + v_1(z) + v_2(z)$, since 
the projection of $s(z)$ on the line orthogonal to $e_0$ is the same as for $v_0(z)$,
so its length is $|\sin\beta|$. Hence by \eqref{26.15}
\begin{equation} \label{27.12}
\beta^2 L_- \leq 2 |z-\ell_-| \beta^2
\leq C |s(z)|^2 |z-\ell_-| \leq C \sigma.
\end{equation}
We should also mention that by \eqref{9.2}, the angle of $v_0(z)$ with $v_1(z)$
or $v_2(z)$ is at least $\pi/2$, and these two vectors make an angle roughly equal to
$2\pi/3$ with each other; this forces $v_0$ and $e_0$ to make an angle larger than
$5\pi/6 - 10^{-2}$, say. At any rate, seen from $z$, both $\ell_-$ and $\ell'_-$ lie in 
a direction almost opposite to $e_0$. Let us restate this and the second part of \eqref{27.11}
in terms of the coordinates $h(\ell'_-)$ and $h(z)$ of the points $\ell'_-$ and $z$, along $\rho_0$, 
and which we orient in the direction of $e_0$; we find that
\begin{equation} \label{27.12b}
h(z) = h(\ell'_-) + L_- + \cE_1, \ \text{ with } |\cE_1| \leq C \sigma.
\end{equation}

Next consider the closest point projection $z'$ of $-z_+$ on $\rho_0$; its position on $\rho_0$ is
$-\langle z+z_+, e_0 \rangle$ from $z$ in the direction of $e_0$, modulo an error
of at most $C\sigma$ (because $|z+z_+| \leq C \sqrt\sigma$, so the geodesic does not have 
much time to turn). In terms of coordinates $h(z')$ and $h(z)$ along $\rho_0$, still
oriented in the direction of $e_0$, we find that
\begin{equation} \label{27.12bis}
h(z') = h(z) -\langle z+z_+, e_0 \rangle +\cE_2, \ \text{ with } |\cE_2| \leq C \sigma.
\end{equation}

Finally denote by $\beta_+$ the angle of $-v_0(z_+)$ with 
$e_+ = (v_1(z_+)+v_2(z_+))/|v_1(z_+)+v_2(z_+)|$; the proof of 
\eqref{27.12} also implies that
\begin{equation} \label{27.13}
\beta_+^2 L_+ \leq C |s(z_+)|^2 |z_+-\ell_+| \leq C \sigma.
\end{equation}
Now $|v_i(z_+) - v(-z,w_i)| = |v(z_+,w_i) - v(-z,w_i)| \leq C |z_+ + z| \leq C \sqrt\sigma$ 
by \eqref{27.4}, and $v(-z,w_i) = v(z,w_i)$ (we look at the other tip of the half circle), so
$|e_+-e_0| \leq C \sqrt\sigma$ and 
\begin{equation} \label{27.14}
|v_0(z_+) + e_0| \leq C \beta_+ + C \sqrt\sigma.
\end{equation}
Now $v_0(z_+)$ is the direction of $-\rho_+$ when it leaves from $-z_+$.
Let us compute some more. A parameterization of $-\rho_+$ is given by
\begin{equation} \label{27.15}
z(t) = -z_+ \cos t + v_0(z_+) \sin t , \qquad	 t\in [0,L_+]
\end{equation}
(because a parameterization of a great circle can always be written as 
$z(t) = v_1 \cos t + v_2 \sin t$, for two orthogonal unit vectors $v_1$ and $v_2$,
and then we just need to identify). 

Let $P_0$ denote the $2$-plane that contains $\rho_0$
and let $\pi$ be the orthogonal projection on $P_0$. Define $w= -z_+ - z'$, where $z'$ is
the projection of $-z_+$ on $\rho_0$; we know that 
$|w| \leq |-z_+ - z|\leq C \sqrt\sigma$ because $z\in \rho_0$ and by \eqref{27.4}. 
Also, $w$ is orthogonal to $e'_0$, 
the direction of $\rho_0$ at $z'$, and its orthogonal projection
on the direction of $z'$ is of norm at most $C\sigma$. In fact, when $a$ and $b$ are
two unit vectors (such as $-z_+$ and $z'$), then the projection of $w=b-a$ on the line through
$a$ (or $b$, this is the same) has norm at most $|w|^2/2$, because 
$1 = ||b||^2 = ||a+w||^2 = 1 + 2\langle a, w \rangle + ||w||^2$.
Altogether $|\pi(w)| \leq C \sigma$. 

Similarly write $v_0(z_+) = e'_0 + \xi$; then $|\xi| \leq C \beta_+ + C \sqrt\sigma$ 
by \eqref{27.14} and because $|e'_0-e_0| \leq \ddist(z,z') \leq C \sqrt\sigma$.
Next the projection of $\xi$ on the direction of $v_0(z_+)$ has a norm at most 
$C \beta_+^2 + C \sigma$ by the same argument as above (take $a=v_0(z_+)$
and $b=e'_0$). In addition, 
$$
|\langle \xi, z_+\rangle| = |\langle e'_0, z_+\rangle| = |\langle e'_0, z_+ +z'\rangle|
= |\langle e'_0, w \rangle| \leq C |w|^2 \leq C \sigma
$$
because $v_0(z_+)$ is orthogonal to $z_+$ and $z'$ is orthogonal to $e'_0$, and then, as before,
because $z'$ is the ``orthogonal'' projection of $z$ on $\rho_0$.
If $\pi'$ denotes the orthogonal projection on the plane that contains $- z_+$ and $v_0(z_+)$,
we see that $|\pi'(\xi)| \leq C \beta_+^2 + C \sigma$, but since $||\pi-\pi'|| \leq C \sqrt\sigma$
and $|\xi| \leq  C \beta_+ + C \sqrt\sigma$, we get that
\begin{equation} \label{27.16}
|\pi(\xi)| \leq C \beta_+^2 + C \sigma + C \beta_+  \sqrt\sigma \leq C \beta_+^2 + C \sigma.
\end{equation}
Set $\wt\ell = z' \cos L_+ + e'_0 \sin L_+$; this is the point of $\rho_0$ that lies at 
distance $L_+$ from $z'$ (in the direction of $e'_0$); in terms of coordinates along $\rho_0$, 
this means
\begin{equation} \label{27.16a}
h(\wt\ell) = h(z') + L_+.
\end{equation}
Notice that since $-\ell_+ = z(L_+)$ (the final point of $-\rho_+$), \eqref{27.15} yields
\begin{eqnarray} \label{27.17}
\wt\ell + \ell_+ &=&  z' \cos L_+ + e'_0 \sin L_+ - z(L_+)
\nn\\
&=& (z'+z_+) \cos L_+ + (e'_0 - v_0(z_+))\sin L_+
= - w  \cos L_+ - \xi \sin L_+
\end{eqnarray}
which implies that
\begin{equation} \label{27.18}
|\wt\ell + \ell_+| \leq |w| + L_+ |\xi|
\leq C \sqrt\sigma + C L_+ \beta_+
\ \text{ and } 
|\pi(\wt\ell + \ell_+)| \leq C \sigma + C L_+ \beta_+^2
\end{equation}
Let $\ell'_+$ denote the projection of $-\ell_+$ on $\rho_0$; then by \eqref{27.18}
$$
|\wt\ell - \ell'_+ | \leq 2 |\pi(\wt\ell + \ell_+)| \leq C \sigma + C L_+ \beta_+^2
$$
and 
$$
|\ell'_+ + \ell_+| \leq C \sqrt\sigma + C L_+ \beta_+.
$$
The first estimate yields
\begin{equation} \label{27.18a}
|h(\ell'_+)-h(\wt\ell)| \leq 2|\ell'_+ - \wt\ell |  \leq C \sigma + C L_+ \beta_+^2
\end{equation}
and when we combine with \eqref{27.12b}, \eqref{27.12bis}, and \eqref{27.16a}, we get that
\begin{equation} \label{27.19a}
\ddist(\ell'_-, \ell'_+) = |h(\ell'_+)-h(\ell'_-)| = L_- + L_+ - \langle z+z_+,e_0 \rangle + \cC_3,
\end{equation}
with $|\cC_3| \leq C\sigma + C L_+ \beta_+^2 \leq C \sigma$ by \eqref{27.13}.

We now add the orthogonal complement, which may remove some distance
because $\ell_-$ and $-\ell_+$ may turn out to be on the same side of $\rho_0$, and thus be
closer to each other than their projections are. But the estimates above yield
\begin{equation} \label{27.20}
d : = |\ell_- -\ell'_-| + |-\ell_+ - \ell'_+| \leq C \sqrt\sigma + C \beta L_- + C \beta_+ L_+
\end{equation}
and we claim that when $a, b \in \S$ lie within $d$ of a geodesic $\rho_0$, $d$ is small enough,
and $\ol a$ and $\ol b$ denote their respective projections on $\rho_0$, then
$\ddist(a,b) \geq \ddist(\ol a, \ol b) - C d^2$. Indeed, let $\pi$ be the projection on the plane 
that contains $\rho_0$, observe that $|\pi(a)-\ol a| \leq C d^2$ and similarly for $b$, and that 
$|a-b| \geq |\pi(a)-\pi(b)| \geq |\ol a- \ol b| - C d^2$, from which we deduce the result
because $\ddist(a,b) = 2\arcsin(|a-b|/2)$. 

From the claim, \eqref{27.20}, \eqref{27.19a}, and the fact that $d^2 \leq C \sigma$ we deduce that 
\begin{equation} \label{27.21}
\ddist(\ell_-, -\ell_+) \geq  L_- + L_+ - \langle z+z_+,e_0 \rangle - C\sigma.
\end{equation}
We combine this with \eqref{27.12a} and get that
\begin{equation} \label{27.24}
\Delta_L \leq \wt\alpha \langle z_+ +z, e_0 \rangle + C \sigma
\leq C \sigma + C \alpha \sqrt\sigma.
\end{equation}
Recall that Lemma \ref{t25.1}, our goal for this section, will follow as soon as we prove
that $\Delta_L \leq C \sigma$; see near \eqref{26.2a}. So we may assume that
$|\alpha| \geq C_1 \sqrt\sigma$, with $C_1$ quite large. 
Recall from \eqref{27.8} that
\begin{equation} \label{27.25}
v_1(z_-)+v_2(z_-) = 2e_0\cos\big(\frac{\pi}{3}+\alpha\big) =: (1-\wt\alpha)e_0.
\end{equation}
When $\alpha < 0$, $\wt\alpha$ is negative too, and $|\wt\alpha| > C^{-1} |\alpha|$.
In this case $|v_1(z_-)+v_2(z_-)|-1 = |\wt\alpha|$, and 
\eqref{26.19a} implies that $\sigma \geq C^{-1} |\wt\alpha|^2 \geq C^{-1} |\alpha|^2$.
We choose $C_1$ above large enough and exclude this case. So we assume that $\alpha > 0$,
and now \eqref{27.25} implies that $|v_1(z_-)+v_2(z_-)| \leq 1-\wt\alpha \leq 1 - C^{-1}\alpha$.
Recall that $s(z_-) = v_0(z_-)+v_1(z_-)+v_2(z_-)$; then
\begin{equation} \label{27.26}
\langle s(z_-), v_0(z_-) \rangle = 1 + \langle v_1(z_-)+v_2(z_-), v_0(z_-) \rangle
\geq C^{-1}\alpha.
\end{equation}
If we could apply \eqref{26.13}, we would get that 
$\alpha^2 \leq C \langle s(z_-), v_0(z_-) \rangle^2 \leq C \sigma$, and we excluded this case.
Then we can apply \eqref{26.14} and we get that 
\begin{equation} \label{27.27}
\sigma \geq C^{-1} |z_--\ell_-| \langle s(z_-), v_0(z_-) \rangle
\geq C^{-1} \alpha L_- .
\end{equation}
Let us also try the same estimate near $\ell_+$. Define $\alpha_+$ by
$\Angle(v_1(z_+),v_2(z_+)) = \frac{2\pi}{3}+2\alpha_+$.
Recall that $|v_i(z_+) - v_i(z_-)| =|v_i(z_+) - v_i(-z_-)| \leq C |z_+ +z| \leq C \sqrt\sigma$
by \eqref{27.4}, so $\alpha_+ \geq \alpha - C \sqrt\sigma \geq \alpha/2$ if $C_1$ is large enough,
and the proof of \eqref{27.27} also yields
\begin{equation} \label{27.28}
\sigma \geq C^{-1} \alpha_+ L_+ .
\end{equation}
We complete this with a lower bound on $L_-+L_+$. If $\H^1(\rho^\ast) \leq 2\pi$,
then $\Delta_L \leq 0$ simply because $\H^1(X_0 \cap \S) \geq 2\pi$, so we may assume
that $L_{12}+ L_- + L_+ = \H^1(\rho^\ast) \geq 2\pi$ (by \eqref{27.10}). We combine with
\eqref{27.9} and get that 
\begin{equation} \label{27.29}
L_- + L_+ \geq 2\pi - L_{12} \geq  \langle z_+ +z, e_0 \rangle 
-C \sigma - \wt\alpha \langle z_+ +z, e_0 \rangle,
\end{equation}
hence, since $|\wt\alpha| \leq 1/2 $, $\langle z_+ +z, e_0 \rangle \leq 2(L_- + L_+) + 2C \sigma$.
We may now return to \eqref{27.24}, which yields
\begin{equation} \label{27.30}
\Delta_L \leq \wt\alpha \langle z_+ +z, e_0 \rangle + C \sigma
\leq 2\wt\alpha (L_- + L_+) + C\sigma
\leq C \alpha (L_- + L_+) + C \sigma \leq C \sigma
\end{equation}
by \eqref{27.29}, \eqref{27.27}, and \eqref{27.28}. 

This finally completes our proof of \eqref{26.2a} and Lemma \ref{t25.1} in our last case.
As was mentioned at the end of Section \ref{S25}, this also completes our proof of 
Proposition \ref{t23.3}, Theorem~\ref{t23.1} (which in fact were
finished before), Proposition~\ref{t23.4} and Theorem~\ref{t23.2}. 

\section{More cases where the free attachment is allowed}
\label{Sfree}

We interrupt the study of $E$ in balls centered on $E \sm L$ with some comments on the free attachment.
In the construction of competitors, both in Sections \ref{S13}-\ref{S15} (with balls centered on $L$)
and Sections \ref{S25}-\ref{S27} (with balls centered on $E \sm L$), there are situations where we can
use what we call the ``free attachment'', near one or two of the points $\ell_\pm$ of $\S_r \cap L$.
Recall that the main part of the construction of curves in $E \cap \S_r$ happens in two small disks
$D_\pm$ near the $\ell_\pm$, and we used the free attachment in the following situations.

An extreme case of free attachment is what we called a hanging curve, when 
$E \cap \d D_\pm$ has a point that is not in the same connected component $E \cap D_\pm$
as any other point of $\{ \ell_\pm \} \cup E \cap \d D_\pm$. We like this situation a lot, because
we can contract the hanging curve, use this to find a competitor which is much better than the curve,
and at the end of the estimate show that $r \theta'(r) \geq C^{-1}$ or $r F'(r) \geq C^{-1}$. 
In the phase where we try to deduce geometric properties from the small size of $f$, 
as in Section \ref{S18} and the upcoming Section \ref{S28}, we can forget about this case, 
because this never happens for the good radii $r$ that we select, 
by \eqref{18.27} and \eqref{19.27n} in the centered
case, and similar upcoming estimates in the non centered case.

Next assume that $E \cap \d D_\pm$ has exactly two points; then we talk about free attachment
when these two points lie in the same connected component of $E \cap D_\pm$
and in addition $\ell_\pm$ does not lie in the same connected component of $E \cap D_\pm$ 
as these two points (or just $\ell_\pm \in L \sm E$).
Except for hanging curves, this is the only case of free attachment that we have in the context of 
Sections \ref{S25}-\ref{S27} (and we called this Configuration $3$).

Another case, that shows up in Sections \ref{S13}-\ref{S15},
is when $E \cap \d D_\pm$ has three points, that all lie in the same component of
$E \cap D_\pm$, but this component does not contain $\ell$ (either because $\ell_\pm \notin E$ or because
it lies in some other component); we called this Configuration $3-$. And the last case is when
two of the three points of $E \cap \d D_\pm$ lie in a same component of $E \cap \d D_\pm$
but this component does not contain $\ell_\pm$; we call this Configuration $2+1$ (when $\ell_\pm$
is connected to the third point of $E \cap \d D_\pm$).

When we have a free attachment near $\ell_\pm$, we are happier because when we construct
competitors, we don't need to worry about the sliding condition near $\ell_\pm$.
Typically, we select a point $z_\pm \in E \cap D_\pm$,
the net $\gamma$ of curves of $E \cap \S_r$ that we construct consists near $D_\pm$ in two 
curves $\gamma_j$ that start from $z_\pm$, plus maybe (in Configuration $2+1$) a curve that
leaves from $\ell_\pm$ and does not get near the $\gamma_j$.
The same thing happens with the Lipschitz curves $\Gamma_j$ that we construct starting from 
$\gamma_j$. It is often very convenient to have a free attachment, because for instance if the two 
curves $\Gamma_j$ that end at $z_\pm$ make an angle at $z_\pm$ that is far from $\pi$, 
we can modify our first main competitor (the set $F^1$ built in Section \ref{S14} or the set 
$F^0 = \varphi^0(E)$ that shows up above \eqref{25.20}), using the same method as when we 
use the full length property.
That is, we use the fact that the tip of the current competitor coincides with the cone 
over the union $\rho^\ast = \rho^\ast_r$ of the geodesics with the same endpoints as the $\gamma_j$ 
and the $\Gamma_j$, to save some area near the tip if the angle $\alpha_\pm$ of the two 
geodesics $\rho_j$ that end at $z_\pm$ if far from flat. With this manipulation, we save about
$C^{-1} r^2 (\pi-\alpha_\pm)^2$ in area. If $\pi - \alpha_\pm \geq 10^{-2}$, say, this leads
to a very good estimate like the one that we get in \eqref{25.27} or \eqref{25.30}, which itself leads 
to a good lower bound on $\theta'(r)$ or $F'(r)$ and later on, when we try to get a geometric control,
excludes $r$ of the list of good radii, again by \eqref{19.27n}.
This is typically what happens in the situation of Theorem \ref{t23.2} and Proposition \ref{t23.4}.
In principe it means that when $E$ is well approximated by a non-flat set of type $\bV$, the free
attachment situation will not occur. 

In the non centered case of Sections \ref{S25}-\ref{S27}, we also have to think about the triangular
face $T(r)$. For the moment, when we have a free attachment near $\ell_+$, we are simply allowed
to detach $z_+$ from $L$ (or $T(r)$), but we shall see soon that we may also consider that there
is a free attachment near $\ell_-$, and even we'll be able to drop $T(r)$ because we can get away with 
the sliding condition. 

The goal of this short section is to observe that when $E$ does not contain $L \cap \ol B(0,\rho)$,
then we can use the estimates that come with the free attachment 
for all the radii $r$ near $\rho$, even if for some of them, $\ell_\pm(r) $ lies in the same component 
as the other points of $E \cap \d D_\pm(r)$. We first give a statement for the case when $0 \in E \cap L$,
prove the statement, and then discuss a variant for the non centered case and how this could be applied.

\begin{lem}\label{tf.1}
Suppose that $0 \in E \cap L$ and for some $\rho > 0$, $E$ does not contain $L \cap \ol B(0,\rho)$.
Then, for $C^{-1}\rho < r \leq 2\rho$, we can do the estimates 
that lead to differential inequalities of Sections \ref{S16}-\ref{S21} as if we always had free 
configurations in the description of Section \ref{S9}. Yet we need to replace $r^2 h(r)$ with 
$9\rho^2 h(3\rho)$ in the estimates.
\end{lem}

The estimate that we have in mind are \eqref{14.46}, \eqref{15.4}, and their variants that were
used in Sections \ref{S18}-\ref{S20}.
These estimates in turn imply some differential inequalities, which we don't mention here.

As we will see in the proof, the reason for the replacement of $r^2 h(r)$ is that we have to use
competitors of $E$ where we modify $E$ near $B(0,\rho)$, hence the error terms get that large. 
Here $C$ is any given positive constant given in advance, and it should be noted that the only
price that we pay for taking $C$ large is the fact that the error term $9\rho^2 h(3\rho)$ is not
necessarily that small compared to $r^2$.

The main point of the proof will be that when $L \cap \ol B(0,\rho) \sm E \neq \emptyset$,
we can prepare the work by finding a first (sliding) competitor $F_0$ of $E$, 
in the ball $B(0,3\rho)$, which is almost as good as $E$ itself, but for which 
$F_0 \cap L \cap \ol B(0,2\rho) = \emptyset$.
Then we replace $E$ with $F_0$ in all the proofs above, and get almost the same results, 
except for the following details. First we lose a small quantity $\eta > 0$ when we replace $E$ with $F_0$, but this does not matter because $\eta$ will be as small as we want. 
But also, and this is the reason for the replacement discussed above, the competitors that we construct now are only competitors for $E$ in the ball $B(0,3\rho)$, 
so the error terms get a little larger. Oh course when we use the fact that $h(r) \leq C_h r^{\beta}$,
this amounts to multiplying $C_h$ by $9C^{2+\beta}$, which is not too bad.

Let us now prove the main estimate for the lemma. Let $\rho$ be as in the statement, and also find
coordinates of $\R^n$ so that $R^n \simeq L \times \R^{n-1} \simeq \R \times \R^{n-1}$.
By assumption, one of the points of $L \cap \ol B(0,\rho)$ does not lie in $E$; 
let us write this point $y = (t\rho,0)$, with $t\in [-1,1]$. 
Let $\eta > 0$ be given, as small as we want, and let us construct our competitor $F_0$ so that
\begin{equation}\label{f.14}
\H^2(F_0 \sm E) \leq \eta.
\end{equation}
We start with the choice of a very thin tube $T$, where most of the construction will happen.
For reasons that will be clear soon, we prefer $T$ to be composed of cubes.
Let $\tau > 0$ be small, to be chosen later (depending on $\eta$), but certainly so small that 
$B(y,3\tau)$ does not meet $E$. Identify $L$ with $\R$ and $y$ with $t\rho \in \R$,
and denote by $K$ the set of integers $k\in \bZ$ such that 
$I_k := [\rho + k\tau,\rho+(k+1)\tau]$ meets $[-2\rho,2\rho]$. Then
set $I = \cup_{k\in K} I_k$; thus $[-2\rho,2\rho] \subset I \subset (-3\rho,3\rho)$.
Also write $I = [a,b]$, denote by $Q$ the cube in $\R^{n-1}$ of side length $\tau$ and 
centered at $0$, set $Q_k = I_k \times Q \subset L \times R^{n-1} \simeq \R^n$, 
and finally set $T = I \times Q = \cup_{k\in K} Q_k$.

We start with a Lipschitz mapping $f_0$ such that $f_0(x)=x$ on $\R^n \sm T$, that
maps $T$ to its boundary $\d T$, the interval $[t\rho+\tau,b] \subset L$ to the point 
$b \in L$, and similarly $[a,t\rho-\tau]$ to $a \in L$. This is because we want to respect the 
sliding boundary condition.

When $n=3$, we can take $F_0 = f(E)$, notice that $F_0$ is a (sliding) competitor for $E$
in $B(0,3\rho)$ (because the linear interpolation between the identity and $f_0$ gives a one parameter family of mappings with the desired properties, and that \eqref{f.14} holds. More precisely, if we set 
$W_0 = \big\{ x\in E \, ; \, f_0(x) \neq x \big\}$, then
\begin{equation}\label{f.15}
\H^2(F_0 \sm E) \leq \H^2(f_0(W_0)) \leq \H^2(\d T) \leq C \tau \rho < \eta
\end{equation}
if $\tau$ is small enough. When $n > 3$, we cannot estimate like this because $\H^2(\d T) = +\infty$,
and even though $\H^2(f(E \cap T))$ is finite because $f_0$ is Lipschitz, it may be much too large
for our taste. So we shall compose $f_0$ with a Federer-Fleming projection.
Write each $Q_k$, $k\in K$, as a union of $2^{n-1}$ cubes $Q'_j$ of side length $\tau/2$,
and thus write $T$ as a union of smaller cubes $Q'_j$, $j\in J$. 
We do this because we want $0$ to be a vertex and $L$ to be contained in the $1$-skeleton
of $T$ (seen as the union of the $Q'_j$). We add to the $Q'_j$ the cubes of the same ``dyadic net''
(and the same side length $\tau/2$) that touch the $Q'_j$; we then get a new tube $T' \supset T$,
twice thicker and a tiny bit longer, which is a union of cubes $Q'_j$, $j\in J$. 

The Federer-Fleming projection will occur in $T'$, which means that we shall use the composition 
$f_1 = \varphi \circ f$, where $\varphi$ is a new Lipschitz mapping such that $\varphi(x) = x$ 
for $x\in \R^n \sm T'$, $\varphi(T') \subset T'$, and even $\varphi(Q'_j) \subset Q'_j$
for $j \in J'$. This mapping is constructed with the same standard scheme as in Chapter 3 of \cite{DSQM},
so we only recall how the construction goes and the properties of $\varphi$ that will be helpful.
We start with the observation that $T' \cap f(T)$ has a finite (although possibly large) $\H^2$ measure.
Our mapping $\varphi$ is itself a composition of elementary Federer-Fleming projections that act on faces
of various dimensions. Each elementary Federer-Fleming projection consists in choosing ``centers'' $x_F$ 
inside the faces $F$ of cubes that compose $T'$, so that they are not contained in the current image 
(we start with $f(E)$, but as the construction goes, we consider the images of that set by the 
previously constructed mappings), and we compose with a Lipschitz mapping that coincides 
on the current image with the radial projection on $F$, centered
at $x_F$, that maps $F \sm \{ x_F \}$ to $\d F$ and is the identity on $\d F$. We proceed independently
on all the faces of the same dimension, but thanks to the fact that we always take the identity on $\d F$,
we get a global Lipschitz map. We first do this on the faces of dimension $n$, then $n-1$, and so on,
and end with a projection of the $3$-faces on their $2$-dimensional boundaries. Each time, we use the fact
that the $\H^2$ measure of the image of $f(E)$ by the previous mapping is finite to choose
$x_F$ outside of that image, and in fact sufficiently far from that image in average, so that the projection 
will never multiply the measure by more than $C$.

In fact, we only do this on some of the faces of the $Q'_j$. On the $n$-faces (i.e., the interiors) 
of the cubes that compose $T$, we don't really need to do this, because we have no piece of $f(E)$ 
left there anyway, but it does not hurt either. In the faces that are not contained in $\d T'$, 
we do the construction as described above, so as to get a $2$-dimensional set.
But on the faces that are contained in $\d T'$, we do not do anything, i.e., we keep the identity.
This is important because we take $\varphi(x)=x$ on $\R^n \sm T'$.

Notice that $\varphi$ preserves the cubes, but also the faces. Because of this, it preserves $L$
and so does $f_1$; thus $f_1(E)$ is a sliding competitor for $E$. We need to estimate 
$H^2(f_1(E) \cap T') = H^2(f_1(E \cap T'))$. One piece is $f_1(E \cap T)$, and for this piece we know that
we followed the construction down to $2$-faces. That is, this set is contained in the $2$-skeleton of
$T'$, which has a $H^2$-measure smaller than $C \tau^2 (\sharp K) \leq C \tau \rho$.
For $f_1(E \cap T' \sm T) = \varphi(E \cap T' \sm T)$, we observe that if we choose the centers $c_F$
correctly, its measure is multiplied by at most $C$, so that
\begin{equation}\label{f.16}
\H^2(\varphi(E \cap T' \sm T)) \leq C \sum_{j\in J'} \H^2(E \cap Q'_j)
\leq C (\sharp K) \tau^2 \leq C \tau \rho
\end{equation}
by the local Ahlfors regularity of $E$. We may now choose $\tau$ so small that
\begin{equation}\label{f.17}
\H^2(F_0 \sm E) \leq \H^2(f_1(W_0)) \leq C \tau \rho \leq \eta,
\end{equation}
where $W_0 = \big\{ x\in E \, ; \, f_0(x) \neq x \big\}$ as above. This proves \eqref{f.14}.
Notice also that $W_0 \subset T'$, and hence 
$\H^2(W_0) \leq C \tau^2 (\sharp K) \leq C \tau \rho$ by the same argument as above, 
using the the local Ahlfors regularity of $E$.

The reader may be worried, because the set $F_0$ that we just constructed is no longer almost
minimal. So we don't want to use estimates that would rely on the almost minimality of $F_0$. 
The natural solution would be to adapt the construction to $F_0$, but this is not what we will do. 
Instead, we just compute brutally with our initial set $E$, construct ``competitors'' 
$F^i = \varphi^i(E)$ with the free attachment if needed, and estimate $\H^2(F^i)$. 
Now the $F^i$ are perhaps not competitors, because using the free attachment may violate 
the boundary condition that $\varphi^i(E \cap L) \subset L$, so we are not allowed to compare 
$F^i $ with $E$ directly. 
There is no such problem with $F_0$, because $F_0 \cap L = \emptyset$ on $B(0,3\rho)$
where $\varphi^i$ moves points, and so $\varphi^i(F_0)$ is really a competitor for $E$
(but in the larger ball $B(0,3\rho)$). Now we use the fact that $\varphi^i$ is Lipschitz,
and let $\tau$ and $\eta$ tend to $0$ in the estimate above. Observe that then 
$H^2(F_0 \cap B(0,3\rho))$ and $\H^2(\varphi^i(F_0 \cap B(0,3\rho)))$
 tend to $\H^2(E \cap B(0,3\rho))$ and $\H^2(\varphi^i(E \cap B(0,3\rho)))$, so that we 
 get the desired estimates on $E$ by applying the almost minimality of $E$ to the competitor
 $\varphi^i(F_0)$, and then taking a limit.
Lemma \ref{tf.1} follows.
\qed

\ms
Let us now state the variant of Lemma \ref{tf.1} for balls centered on $E \sm L$. 

\begin{lem}\label{tf.2}
Suppose that $0 \in E \sm L$ and for some $\rho > 0$, 
$E$ does not contain $L \cap \ol B(0,\rho)$.
Then, for $C^{-1}\rho < r \leq 2\rho$, we can do the estimates 
that lead to differential inequalities of Sections \ref{S25}-\ref{S27} as if we always had free 
configurations in the description of Section \ref{S9}. 
In particular, we don't need $T(r)$ and we may drop $\H^2(T(r))$ from the estimates.
Yet we need to replace $r^2 h(r)$ with 
$9\rho^2 h(3\rho)$ in the estimates.
\end{lem}

This sounds a little bit like winning the jackpot, but of course what this means is that in most 
situations, $E$ contains $L \cap \ol B(0,\rho)$. The proof is the same. First we construct a 
competitor $F_0$ for $E$ in $\ol B(0,3\rho)$, such that \eqref{f.14} holds, and which no longer
meets $L \cap \ol B(0,2\rho)$. The proof goes as before (we never used the fact that $L$ contains $0$),
and then we can end the argument as above. As was suggested earlier, not only we can use the free
attachment for the estimates, but since we no longer have to enforce the sliding condition for our
competitors, we don't need to add the triangular piece $T(r)$ either. The lemma follows.
\qed

\ms
Let us just give an example of how we may use the Lemmas.
Suppose that $0 \in E$, $h(R)$ is small enough, and that in addition $E$ is quite close to
a generic set $X \in \bV$, such that the half planes that compose $X$ make an angle smaller 
than $\pi - 10^{-2}$, say. We may either assume that $0 \in L$ as in the early sections, 
or that $0 \in E \sm L$ and $R^{-1}\dist(0,L)$ is small enough.
Then $L \cap B(0,R/2) \sm B(0,10^{-2}R)$ is contained in $E$.
Indeed otherwise we may apply Lemma \ref{tf.1} or Lemma \ref{tf.1}, find that we can apply 
the free attachment construction for all the nearby radii $r$, get a very good estimate 
for such $r$ that imples that $\theta'(r) \geq C^{-1} r^{-1}$ or $F'(r) \geq C^{-1}r^{-1}$, 
and get a contradiction with the fact that, when $E$ lies close enough to a $\bV$ set, 
$\theta$ or $F$ is nearly constant in the range under consideration. 
In fact, we can also iterate this argument (apply it to $R/2$, $R/4$, and so on)
and get that $L \cap B(0,R/2) \subset E$.
We will detail the argument during the proof of Lemma \ref{t29b.4}, 
mostly as an example of  how it may go and to give a flavor of why we get estimates 
like $\theta'(r) \geq C^{-1} r^{-1}$.

We may even apply the same argument to the case when $E$ is very close to a half plane
in $B(0,R)$, and get the same conclusion that $L \cap B(0,R/2) \subset E$. This time,
when we apply Lemma \ref{tf.1} or \ref{tf.2}, instead of a standard free attachment, 
we immediately get a hanging curve near $\ell_\pm$, which also gives a bound on 
$\theta'(r)$ or $F'(r)$ that is incompatible with the fact that $\theta$ or $F$ is nearly constant. 
We shall also sketch a more direct argument, when we discuss the proof of \eqref{29a.9},
and we will find the proofs of Lemmas \ref{tf.1} and \ref{tf.2} convenient in Section \ref{S30},
when we check the full length property in some special cases.

\section{Geometric estimates follow from the decay of $F$}
\label{S28}

The decay of $F$ that we got in Sections \ref{S26} and \ref{S27} is not so much good in itself,
but it will allow us to control the geometry of $E$. 
In this section we prove two main statements to this effect, corresponding to the densities
$\theta_0 = \pi$ and $\theta_0 = \frac{3\pi}{2}$ of Theorems \ref{t23.1} and \ref{t23.2}.

We start with a statement in the simpler case of Theorem \ref{t23.1}, with an approximation by half planes, 
where we will see that under the assumptions of Theorem \ref{t23.1}, we also have a good control 
(with decay) on the approximation numbers $\beta_H(r)$, in the interesting region where $r \geq d_0$. 
We give the statement first, and then comment.

\begin{thm}\label{t23.5} 
There exist constants $\varepsilon_3 > 0$ and $C_6 \geq 1$, that depend only on $n$
and $\beta \in (0,1]$, such that the following holds.
Let $E$ be a reduced sliding almost minimal set in $B(0,400R)$, with a boundary condition 
coming from the line $L$, and a gauge function $h$ such that 
\begin{equation}\label{23.48}
h(r) \leq C_h r^{\beta} \ \text{ for } 0 < r \leq 400R,
\end{equation}
for some $C_h$ such that $C_h R^{\beta} \leq \varepsilon_3$.
Suppose that $0 \in E$ and $0 < d_0 = \dist(0,L) \leq R/2$.
Then 
\begin{equation} \label{23.49}
d_{0,R}(E,H_0) \leq C_6 \Big[[F(200R) -\pi] + C_h R^\beta\Big]^{1/4} 
\end{equation}
where $F$ is defined by \eqref{22.3} and $H_0$ denotes the half plane bounded by $L$ 
that contains the origin.
\end{thm}

Notice the analogy with Theorem \ref{t18.1}, but here the center is off $L$.
Of course this is only useful when the right-hand side of \eqref{23.49} is small, so that in particular 
the density excess $F(200R) -\pi$ is small. Here $\pi$ is the smallest value that 
$\lim_{t \to 0} F(t)$ could possibly take (because $0 \in E \sm L$); this is also why we do 
not need to put in the assumption that $0$ is a point of density $\pi$.

We required that $d_0 \leq R/2$, but we do not feel bad about this; for $R < d_0$,
there is no sliding condition in $B(0,R)$, so we may still show $E$ is very well approximated
by planes in $B(0,R)$, using the regularity theorems for plain almost minimizers. 
This is just a different story. Notice however that if $200R < d_0$, the other assumptions 
of the theorem allow $E$ to coincide with any plane in $B(0,R)$,  not just the ones that
nearly contain $L$.

The point of this estimate is not to give some rough control on $d_{0,R}(E,H_0)$
(we will see something like this as soon as \eqref{23.53}), but to use this rough control to 
get much better estimates that depend only on the density excess and $h$.
Since we proved earlier that this excess tends to decay like a power, this will give a good decay 
for geometric quantities as well. 

\begin{rem}\label{r30a2} 
We can prove an even better control when the gauge function is even smaller than $C_h r^{\beta}$.
Set
\begin{equation} \label{30a3}
J(R) = \int_0^{2R} \frac{h(t) dt}{t} \ \text{ and } 
J_+(R) = \sum_{k \geq 0 \, ; \, 10^{-k} R \geq d_0} J(10^{-k} R)^{1/2}.
\end{equation}
We shall also prove that, under the assumptions of Theorem \ref{t23.5}, we have the estimate
\begin{equation} \label{30a4}
d_{0,R}(E,H_0) \leq C_6 [F(200R) -\pi]^{1/4} + C_6 J(200R)^{1/4} + C_6 J_+(200R)^{1/2}.
\end{equation}
\end{rem}

Notice that this is better than \eqref{23.49}, because $J(R) \leq C C_h R^{\beta}$
and $J_+(R) \leq C C_h^{1/2} R^{\beta/2}$ when \eqref{23.48} holds.
The strange definition of $J_+(R)$ reflects some of the trouble we will have with the proof,
where we will need to fetch information at the scale $d_0$ (to get the relative position of
$H_0$, $L$, and $0$) and return to the possibly much larger scale $R$.

We can use Theorem \ref{t23.5} to prove the regularity of $E$ when it satisfies the 
assumptions of Theorem \ref{t23.1}. Indeed, that theorem gives us good estimates on
the density excess $F(200R) -\pi$, even with some decay, and Theorem \ref{t23.5} then says
that $E$ is close to $H_0$ in all the balls $B(0,R)$, $R \geq 2d_0$. 
We can even get a good control in smaller balls $B(0,R)$, $R < d_0$, by first applying
the result to $R = 2d_0$ to show that $E$ is close to a plane (the plane that contains $H_0$)
in $B(0,d_0/2)$, and then applying the regularity results for plain almost minimal sets (with no
sliding boundary) in smaller balls; we get additional decay there. The consequence is that we get a 
very good $C^1$ description of $E$ near $0$. See Section \ref{S29a} for more details.

Yet Theorem \ref{t23.5} and the proof of regularity sketched above are not really 
needed to control of $E$ in balls that are not centered on $E$ (first via the decay the 
functional $F$, and then through the geometric control that follows), because we may 
get the desired regularity result otherwise.
When $E$ satisfies the assumptions of Theorem \ref{t23.1} in the large ball $B(0,R)$, 
$R > 10^3 d_0$, say, it turns out that every point of $L \cap B(0,d_0)$ lies in $E$ 
(and has density $\pi/2$). This is proved in \cite{Mono}. 
Then we may also apply the simpler decay results for balls $B(x, r)$ centered on $E \cap L$
(see for instance Corollary \ref{t22.1n}), and get the same geometric information 
in these balls $B(x,r)$, $r \geq 2d_0$, as given by Theorem \ref{t23.5}. 
This is fortunate, because this proof of regularity will help us simplify our proof of
Theorem \ref{t23.5} itself. We will return to this in due time. 

Yet the fact that we can find enough points in $E \cap L$ with the
right density is quite lucky, it seems, and if we could not find these points in $E \cap L \cap B(0,2d_0)$,
we would not be able to apply Corollary \ref{t22.1n} to them!

We will have a second statement (Theorem \ref{t28.2}) similar to Theorem \ref{t23.5}, 
but with points of density $\frac{3\pi}{2}$ and where we approximate $E$ by 
truncated $\bY$-sets. There the story will be different: it seems that we cannot easily
get the regularity results of Sections \ref{S29b}-\ref{S29d} 
without actually applying Theorem \ref{t28.2} to some points of type $\bY$ in $E \sm L$.

The proof of Theorem \ref{t23.5} will be rather long and complicated, and to save some energy
we will group it with the proof of the upcoming Theorem \ref{t28.2}.

\ms
We shall use the following notation concerning \ub{truncated sets of type $\bY$}. 
First denote by $\bY(L, r)$ the set of cones $Y$ of type $\bY$ that are centered at the origin,
and for which $L \cap B(0,r)$ is contained in a face of $Y$. For $Y \in \bY(L, r)$,
we set $Y^t = \ol{Y \sm S}$, where $S$ still denotes the shade of $L$ seen from $0$,
but in fact we are only interested in $Y^t \cap \ol B(0,r)$, where $Y^t$ truly looks like a truncated
cone of type $\bY$, but not necessarily with a straight truncation parallel to the spine of $Y$.
Notice that $Y^t \cap \S_r$ is a net of geodesics like the ones that we studied in Section~\ref{S27},
with two large arcs of great circles (in fact, half circles) and two small tips that connect to 
the points of $L \cap \S_r$ (and may be reduced to one point $\ell_\pm$).

\begin{thm}\label{t28.2} 
There exist constants $\varepsilon_3 > 0$ and $C_6 \geq 1$, 
that depend only on $n$ and $\beta \in (0,1]$, such that the following holds.
Let $E$ be a reduced sliding almost minimal set in $B(0,400R)$, with a boundary condition 
coming from $L$, and a gauge function $h$ such that 
\begin{equation}\label{28.7}
h(r) \leq C_h r^{\beta} \ \text{ for } 0 < r \leq 400R,
\end{equation}
for some $C_h$ such that $C_h R^{\beta} \leq \varepsilon_3$.
Suppose that $0 \in E$, with $F(0) = \frac{3\pi}{2}$, and $0 < d_0 := \dist(0,L) \leq R/2$.  
Then we can find a cone $Y \in \bY(L,R)$ such that
\begin{equation} \label{28.8}
d_{0,R}(E,Y^t) \leq C_6 \Big[ [F(200R) -\frac{3\pi}{2}] + C_h R^{\beta} \Big]^{1/4},
\end{equation}
where $F$ is defined by \eqref{22.3} and $Y^t$ is as above the statement. 
\end{thm}

In fact, under the assumptions of the theorem, we also get that
\begin{equation} \label{30a7}
d_{0,R}(E,Y^t) \leq C_6 [F(200R) - \frac{3\pi}{2}]^{1/4} + C_6 J(200R)^{1/4} + C_6 J_+(200R)^{1/2},
\end{equation}
with $J$ and $J_+$ as in \eqref{30a3}.

As before, we restrict to $R \geq 2 d_0$ because for $r << d_0$ we would get a set of type
$\bY$, but unless we can apply Theorem \ref{t28.2} to a radius $R > 2d_0$, we cannot really
say how it is oriented. Notice however that the approximation in \eqref{28.8} or \eqref{30a7}
is valid on the whole ball $B(0,R)$. The proof will even give some
uniform approximation in all the smaller balls, even leading to the existence to a tangent $Y$-set
that lies close to $Y^t$. See Remark \ref{R28.4} 
and Sections \ref{S29b}-\ref{S29d}. 

\begin{rem}\label{r28.3}
There is more in this statement that one may have expected.
The most important assumption is that the modified density excess $F(200R) -\frac{3\pi}{2}$
is very small, which implies that $F'(r)$ is often small for $r < 200R$. 
Yet, for instance, it could a priori happen that $F'(r)$ is very small for some $r >> d_0$, 
but $E$ looks a lot like a plane, or a flat set of type $\bV$, in $B(0,r)$. 
So we will need to exclude these cases from the discussion, by
comparing all the different scales between $d_0$ and $r$, and then using the fact that the 
density at $0$ is $F(0) = \frac{3\pi}{2}$.
\end{rem}

\ms
We intend to prove Theorems \ref{t23.5} and \ref{t28.2} together, because there are
many common points. The proof will be quite long, even though we shall rely on some of the
computations and estimates that we did for Theorem \ref{t18.1}, so we'll try to cut the proof
into steps, often coming with their own tiny introduction.

One of the features of the proof is that we'll have to go up and down between scales,
and most of our estimates will be obtained by constructing a competitor for $E$ at some 
intermediate scale $d_0 \leq r \leq R$, typically as in the proof of the decay estimate for $F$.
This time the point of the computation is that if the geometry is not almost perfect, then 
we can find a better competitor, which implies that the derivative of $F$ for the corresponding radii
is not too small, and in principle this does not happen much when $F$ almost has the minimal value.

\msi
{\bf Step 1. We make sure that we can apply the construction of Sections \ref{S25}-\ref{S27}.}

We start the proof with a small reduction, that will allow us to apply the construction and estimates 
of Sections \ref{S25}-\ref{S27}, to all radii roughly between $2d_0$ and $180R$. 
For this we apply the near monotonicity of $F$ and our implicit assumption that $F$ essentially 
keeps its minimal value, to get a rough control of the geometry.
Since we want to unify some estimates, it will be convenient to set
\begin{equation} \label{30a8}
\left\{ \  \begin{aligned}
&\theta_0= \pi \ \text{ when $E$ and $R$ are as in Theorem \ref{t23.5},} 
\cr& \theta_0 = \frac{3\pi}{2} \ \text{ when they are as in Theorem \ref{t28.2},} 
\end{aligned}
\right .
\end{equation}
and then
\begin{equation} \label{30a9}
f(r) = F(r) - \theta_0 \ \text{ for } 0 < r < 400R.
\end{equation}

Let $\varepsilon_4 > 0$ be very small, to be chosen later. We may assume that
\begin{equation} \label{23.50}
f(200R) + \int_0^{400R} \frac{h(t) dt}{t} = f(200R) + J(200R) \leq \varepsilon_4,
\end{equation} 
because otherwise \eqref{30a4} or \eqref{30a7} holds trivially.
This is the same justification as for \eqref{18.8}.
Then by \eqref{22.8}, we also get that for $0 < r \leq 200R$,
\begin{equation}\label{23.51}
F(r) \leq \exp\Big(\alpha \int_0^{400R} \frac{h(t) dt}{t}\Big) F(200R) 
\leq e^{\alpha \varepsilon_4} F(200R) \leq e^{\alpha \varepsilon_4} (\theta_0 + \varepsilon_4)
\leq \theta_0 + C \varepsilon_4
\end{equation}
by \eqref{22.7} and \eqref{23.50}. 
We claim that
\begin{equation} \label{30a}
\theta_0 = \lim_{r \to 0} F(r) = \lim_{r \to 0} \theta(r).
\end{equation}
In the case of Theorem \ref{t28.2}, this is our assumption that $F(0)= \frac{3\pi}{2}$,
plus the fact that $F(r) = \theta(r)$ for $r < d_0$. In the case of Theorem \ref{t23.5},
we know that $\lim_{r \to 0} \theta(r)$ exists because $\theta$ is almost monotone, 
and is the density of any blow-up limit of $E$ at $0$. Recall also that $E$ is reduced and 
contains $0$, so these blow-up limits are nontrivial minimal cones. But the only minimal
cones of density smaller than $3\pi/2$ are the planes; now \eqref{30a} follows from \eqref{23.51}
if $\varepsilon_4$ is small enough.

Next we show that $E$ is as close as we want to a set of constant density.
Let $\tau > 0$ be small. We start with the case of Theorem \ref{t23.5}, and show that
\begin{equation}\label{23.53}
d_{0,r}(E,H_0) \leq \tau \ \text{ for } 2d_0 \leq r \leq 180R.
\end{equation}
For this we apply Lemma \ref{t22.3} to $E$, the radius $r_1 = \frac{21r }{ 20}$,
and the large radius $200 R > r_1$ (to play the role of $R$ in the lemma).
The initial assumptions \eqref{22.4} and \eqref{22.7} are satisfied (by \eqref{23.48} in particular),
the constraint \eqref{22.12} too, because $r \geq 2d_0$, we just checked that $\theta_0 = \pi$,
and \eqref{22.14} holds by \eqref{23.51}, because $r_1 < 200 R$, and if $\varepsilon_4$ is small enough. 
Here we are only interested in the first conclusion, which is that 
$d_{0,\frac{20r_1}{21}}(E,H_0) \leq \tau$. This is precisely \eqref{23.53}.
Let us set $X(r) = H_0$ in the present case, so as to unify the notation with the next one.

In the case of Theorem \ref{t28.2}, we claim that there is a constant $\delta > 0$,
that depends on $\tau$, such that if $\varepsilon_4$ is small enough, then for 
\begin{equation} \label{28.11}
\delta^{-1} d_0 \leq r \leq 180R,
\end{equation}
we can find a cone 
$X(r) \in \bV \cup \bP_0$ such that
\begin{equation} \label{28.10}
d_{0,r}(E,X(r)) \leq \tau. 
\end{equation}
This is the same argument, but we replace Lemma \ref{t22.3} with Lemma \ref{t24.2}.
This forces us to restrict to radii $r$ such that \eqref{28.11} holds
(as in \eqref{24.23}), and we need to take $\delta \leq \delta(\tau)$; the rest is the same.

These approximation properties will be useful (see below), but they are not what we want
eventually. First, they come from compactness arguments and are far from being precise enough.
That is, $\tau$ is fixed and we are interested in the cases when the right-hand sides of
\eqref{23.49} and \eqref{28.8} (or their even smaller variants) tend to $0$.
Also, in the case of \eqref{28.10}, we want to prove that $X(r)$ is a nearly sharp
set of type $\bV$, or a truncated cone of type $\bY$, which is more precise than our 
description of $X(r)$. Of course, when $r$ is much larger than $d_0$, a truncated cone of type $\bY$
(centered at $0$) looks a lot like a sharp set of type $\bV$ at the scale $r$.

When $r$ is not much larger than $d_0$, we can deduce the existence of an approximating
truncated set from Lemma \ref{t24.3}.
That is, for any $\delta_1 \in (0,1)$, we claim that if $\varepsilon_4$ is small enough
(depending in particular on $\delta_1$ and $\tau$), then for
\begin{equation} \label{30a16}
2d_0 \leq r \leq \min\Big(\frac{d_0}{\delta_1},180R\Big)
\end{equation}
we can find a minimal cone $Y \in \bY(0,\frac{21r}{20})$ such that 
\begin{equation} \label{30a17}
d_{0,r}(E, Y^t) \leq \tau,
\end{equation}
where $Y^t = \ol{Y\sm S}$ is the corresponding truncated cone.
Just apply Lemma \ref{t24.3}, with the same $r$, the large radius $180R$,
and $\varepsilon = \tau$. 

Again this will be useful. It is closer in spirit to our goal, but we'll have to extend it to
larger radii, and also get much smaller bounds than $\tau$.

Return to the cone $X(r)$ of \eqref{23.53} or \eqref{28.10}. We now claim that 
we can perform all the construction of Sections \ref{S25}-\ref{S27}, which is good enough
to prove the two differential inequalities \eqref{23.15} and \eqref{23.20}, and then 
Theorems \ref{t23.1} and \ref{t23.2}. For all these estimates to hold for (almost every) given $r$,
we need to be able to find $R'$ such that \eqref{23.1}, \eqref{23.1n}, \eqref{23.2n} hold
(as usual, but for that $R'$, which we could for instance take equal to $240R$), but also
$2d_0 \leq r \leq R'/2$ when $\theta_0 = \pi$ and $N d_0 \leq r \leq R'/2$ otherwise, 
as in \eqref{23.12} and \eqref{23.17} or \eqref{25.3} or \eqref{25.4}. 
The main assumptions, though, are that $C_h r^\beta$
and $f(r)$ be small enough, which follow from our assumption that 
$C_h R^{\beta} \leq \varepsilon_3$ and \eqref{23.51}, and that $E$ be close enough to
a minimal cone of type $\bH, \bP_0$, or $\bV$ 
(see \eqref{23.6}, \eqref{23.10}, \eqref{23.14}, \eqref{23.19}, or \eqref{25.2}). 
This last follows from \eqref{28.10}, and the reader should
not worry about the constants depending on $X(r)$, as we can always choose it from a 
fixed finite family. So we'll remember that the construction of Sections \ref{S25}-\ref{S27}
works well, provided that 
\begin{equation} \label{30a18}
2d_0 \leq r \leq 180R \text{ when $\theta_0 = \pi$, and } 
Nd_0 \leq r \leq 180R \text{ otherwise.} 
\end{equation}
Also recall, if you are worried about $\delta$, that we can take $N$ somewhat larger than $\delta$.
This completes this first step of preparation. Next we follow for some time
the argument given in Sections \ref{S18}-\ref{S20}.

\msi
{\bf Step 2. We approximate $E$ in spheres $\S_r$, by some nets of geodesics.}
We try to estimate $E$ on the annulus 
\begin{equation} \label{30a19}
A_0 = B(0,90R) \sm \ol B(0,10^{-1}R),
\end{equation}
and we first proceed independently on most spheres $\S_r$. We assume that
\begin{equation} \label{30a20}
90R \geq N d_0 \, ;  
\end{equation}
otherwise, some parts of the construction will be slightly different but simpler, 
but we shall discuss this in Steps 8 and 9.

We use our first step and \eqref{30a20} to select, as in \eqref{18.21}-\eqref{18.23}, 
a set $\cR$ of full measure in $(10^{-1}R,90R)$ such that we can apply the construction 
of competitors of Sections \ref{S25}-\ref{S27} to any $r \in \cR$, based on 
the approximation by the set $X(R) \in  \bH \cup \bV \cup \bP_0$ that we got in
\eqref{23.53} or \eqref{28.10}.
This yields different nets of curves on the sphere, and in particular the initial net 
$\gamma^\ast = \gamma^\ast_r \subset E \cap \S_r$, and a net of geodesics 
$\rho^\ast = \rho^\ast_r$. 

We also introduce a function $j$, defined as in \eqref{18.24} by
\begin{eqnarray}\label{30a21}
j(r) &=& rf'(r) + f(r) + (1+2 \theta_0\alpha_n) h(2r) 
+ (1+\theta_0\alpha_n) \int_0^r \frac{h(2t) dt}{t}
\nn\\
&=& r F'(r) + f(r) + (1+2 \theta_0\alpha_n) h(2r) 
+ (1+\theta_0\alpha_n) \int_0^r \frac{h(2t) dt}{t},
\end{eqnarray}
where the density excess $f$ is defined by \eqref{30a8} and \eqref{30a9}.
As before, the cosmetic addition of the terms with $\alpha_n$ is done so that
\begin{equation}\label{30.14nn}
j(r) \geq (r \theta'(r))_+ + f(r)_+ + h(2r) + \int_0^r \frac{h(2t) dt}{t} \ \text{ for } r\in \cR,
\end{equation}
as in \eqref{19.27n}. We will prefer to work with the radii $r \in \cR$ such that $j(r)$ is rather
small, and $j(r)$ will control various geometric quantities.

We start with the estimate \eqref{18.25} in Lemma \ref{t18.4}, which says that 
\begin{equation} \label{30.15nn}
\H^1(E \cap \S_r \sm \gamma^\ast_r) \leq C j(r),
\end{equation}
where $\gamma^\ast_r \subset E \cap \S_r$ is our first net of curves; see below \eqref{18.32}. 
The proof can be repeated here; it consists in checking that the various differences of 
lengths $\Delta_j(r)$ that show up in the estimates are dominated by $j(r)$ 
(or else we are in one of the exceptional cases and then $j(r) \geq f'(r)$ was large in the first place).

Next we check that $\rho^\ast_r$ approximates $\gamma^\ast_r$ well, in the sense that
\begin{equation} \label{30.16nn}
d_{0,2r}(\rho^\ast_r , \gamma^\ast_r) \leq C j(r).
\end{equation}
The proof is the same as for Lemma \ref{t18.4} in Section \ref{S18}; 
we prefer not to define a cone $Z(r)$ yet, because the fact that $0 \notin L$ complicates the geometry,
but Lemma \ref{t18.4} concerns only $Z(r) \cap \S_r = \rho^\ast_r$ anyway.
In fact, there is a small lie here: in the special case where we have a free attachment near $\ell_\pm$,
$\rho^\ast_r$ has an additional, isolated point $\ell_\pm$, which we remove from 
$\rho^\ast_r$ before we check \eqref{30.16nn}. That is, $\rho^\ast=\rho^\ast_r$ should be 
replaced in \eqref{30.16nn} with $\rho'$, obtained from $\rho^\ast$ by removing 
the points $\ell_\pm$ with a free attachment.
See the discussion that leads to \eqref{18.45}. Also, it will turn out (later in the argument below,
and independently) that there is no free attachment when $j(r)$ is small, so the issue does not arise after all.

During the proof of Lemma \ref{t18.4}, one also shows that $j(r)$ controls various geometric quantities
that show up in the construction of competitors, such as $\Delta_0(r), \Delta_1(r), \Delta_2(r)$ in
\eqref{18.32}-\eqref{18.34}. In particular, \eqref{18.37} (still valid here with the same proof) says
that
\begin{equation} \label{30.17nn}
\Delta_0(r) + \Delta_1(r) + \Delta_2(r) \leq 10^6 r j(r) \ \text{ for } r \in \cR.
\end{equation}
It is also proved that some configurations, such as hanging curves of free attachments when
$E$ is close to a half plane or a non-flat set of type $\bV$, are impossible when $j(r)$ is small; 
we will return to this issue, but see the discussion below \eqref{18.41}.

\msi
{\bf Step 3. We control the variations of the main part of $\rho^\ast_r$.}
Recall that our nets of curves, and in particular the $\rho^\ast_r$, are initially constructed with
the model $X(R)$. When $\theta_0 = \pi$, $X(R) \in \bH(L)$ and $\rho^\ast_r$
is composed of two main geodesics, which we shall call $\rho_{1,+}$ and $\rho_{1,-}$,
and which go from a midpoint $w_1$ near $X(R)$ to $\ell_+$ and $\ell_-$. We can exclude free
attachments here, at least if we restrict to $r$ such that $j(r)$ is small, because they correspond 
to hanging curves. 

When $\theta_0 = \frac{3\pi}{2}$, $X(R) \in \bV \cup \bP_0$ and $\rho = \rho^\ast_r$
is composed of four main geodesics, which we shall call $\rho_{j,\pm}$, plus maybe some
additional short geodesics $\rho_{\pm}$, depending on which type of attachment.
For the moment, let us not discuss attachment, and concentrate on the large 
$\rho_{j,\pm} = \rho_{j,\pm,r}$.

We want to show that the $\rho_{j,\pm,r}$ vary slowly with $r$ (both when $\theta_0 = \pi$
and when $\theta_0 = \frac{3\pi}{2}$). We proceed as in Section \ref{S19},
isolate any of the two or four $\rho_{j,\pm,r}$, construct vertical curves on $E$,
near the middle of the corresponding interval $I = I_{j,\pm}$ of $X(R)$ (where $E$
is actually a nice $C^1$ graph), and use the co-area formula to control the variation
of angles along these curves.
This starts with Lemmas \ref{t19.1} and \ref{t19.2}, which we can keep as they are. 
In the mean time, we prove the inequality \eqref{19.21}, which
says that (since we no longer normalize $R$ away any more)
\begin{equation}\label{30.18nn}
R^{-1} \int_{r\in (10^{-1} R,90 R)} j(r) \leq C\cE, 
\ \text{ with }  \cE=  f(90R) + \int_{0}^{180R} h(r) \frac{dr}{r};
\end{equation}
this will be useful, because then there are lots of $r\in \cR$ such that $j(r)$ is small.
We use the Lemmas to prove an easier version of Lemma \ref{t19.3}, i.e., the fact that for 
$r, s \in \cR$,
\begin{equation}\label{30.19nn}
d_\H(r^{-1}\wh\rho_{j,\pm,r},s^{-1}\wh\rho_{j,\pm,s}) 
\leq C j(r)^{1/2} + C j(s)^{1/2} + C \cE^{1/2},
\end{equation}
where $d_\H$ denotes the standard Hausdorff distance on the unit sphere (we could also
have used $d_{0,2}$), $\wh\rho_{j,\pm,r}$ denotes the full great circle that contains 
$\rho_{j,\pm,r}$, and similarly for $\wh\rho_{j,\pm,s}$.
That is, for the moment we do not want to control the place where these geodesics stop,
but just their position near $I$; this way we can use \eqref{19.35}, and skip the slightly unpleasant 
discussion below \eqref{19.35}, about guessing where the geodesics meet, and what happens near
$\ell_\pm$.

Let us also observe, as in \eqref{18.41}, that $j(r)$ also controls some geometric information
on $\rho^\ast_r$ relative to its near minimality. We claim that if $v_{j,\pm,r}$ denotes the
unit tangent vector at $m_i$ of $\rho_{j,\pm,r}$ (going in the direction of $z_{\pm}$ or $\ell_\pm$),
then
\begin{equation} \label{30.20nn}
|v_{j,+,r} + v_{j,-,r}| \leq C j(r)^{1/2}.
\end{equation}
In other words, the two main geodesics that leave from $m_j$ are almost aligned.
The reason is the same as for \eqref{18.41}: if not, we can modify our construction
of competitors a little near its tip (where the competitor is a cone near the direction of $m_j$), 
to make the angle flatter and the measure a tiny bit smaller.
Because of this, \eqref{19.31} is really an information on the the pairs of geodesics ending at a 
same point $m_i$, or the sets $\wh\rho_{j,+,r} \cup \wh\rho_{j,-,r}$. Anyway, we shall some times find it 
more convenient to forget some information and just remember that by \eqref{30.19nn}
\begin{equation}\label{30.21nn}
d_\H(r^{-1}\wh\rho_{r},s^{-1}\wh\rho_{s}) \leq C j(r)^{1/2} + C j(s)^{1/2} + C \cE^{1/2},
\end{equation}
where $\wh\rho_{r}$ is the union of the (two or four) pieces $\wh\rho_{j,\pm,r}$,
and similarly for $\wh\rho_{s}$.

\msi
{\bf Step 4. We fetch information from the scale $d_0$.}
In Section \ref{S19}, we used a sharper version of \eqref{30.19nn} directly to control $E$ 
near the $\S_r$; let us not try to do this yet, and consider the variations of the $\rho_{j,\pm,r}$ 
across smaller annuli
\begin{equation} \label{30.22nn}
A_k = B(0,90R_k) \sm \ol B(0,10^{-1} R_k), \ \text{ with } R_k = 10^{-k} R. 
\end{equation}

Recall from \eqref{23.51} that $f(r) \leq C \varepsilon_4$ for $0 < r \leq 200R$.
So $f(R_k)$ is as small as we want, but we shall restrict to $k$ such that $90R_k \geq Nd_0$ 
(or just to $90R_k \geq 2d_0$ when $\theta_0 = \pi$), as in \eqref{30a20}, because this way 
we can find a nice approximating set $X_k = X(R_k)$ as in \eqref{23.53} or \eqref{28.10} and do the 
same construction as above for $R = R_0$. 
So let $e$ denote the largest value of $k$ for which $90R_k > N d_0$
(think, $e$ like ``end'', but the truth is that not so many letters were left); thus
\begin{equation} \label{30.23n}
N d_0 \leq 90 R_e \leq 10 N d_0.
\end{equation}
When $90R < 10 N d_0$, let us still take $e = 0$ and not worry if some of our statements 
below are slightly wrong. We shall return to this case in the last steps
and only small adaptations will be needed, because $\cE$ includes a control on radii 
$r \in [10R,100R]$ which is more than enough. 

For $0 \leq k \leq e$, we define a set $\cR_k$ of full measure in $(10^{-1}R_k,90R_k)$
with the same properties as before (namely, we can do the construction of competitors as in
Sections \ref{S25}-\ref{S27}) and, for $r \in \cR_k$, define the number 
$j(r)$ as in \eqref{30a21}. Then construct the nets $\gamma^\ast_r$ and $\rho^\ast_r$. 
We also get the same estimates as above, but
the small quantity $\cE$ needs to be replaced by 
\begin{equation} \label{30.24n}
\cE_k=  f(90R_k) + \int_{0}^{180R_k} h(r) \frac{dr}{r}.
\end{equation}
Then we select for each $k$ a radius $r_k \in \cR_k \cap [R_k, 2R_k]$, 
so that
\begin{equation} \label{30.25n}
j(r_k) \leq 10R_k^{-1} \int_{\cR_k} j(r) dr \leq C \cE_k,
\end{equation}
where the second inequality comes from \eqref{19.21}.

When we choose $r_{k+1}$, the reference cone $X_{k+1} = X(R_{k+1})$ may be a little different
than $X_k$; yet we claim that the proof of \eqref{30.21nn} also yields
\begin{equation}\label{30.26n}
d_\H(r_k^{-1}\wh\rho_{r_k},r_{k+1}^{-1}\wh\rho_{r_{k+1}}) 
\leq C j(r_k)^{1/2} + C j(r_{k+1})^{1/2} + C \cE_k^{1/2}
\leq C (\cE_k + \cE_{k+1})^{1/2},
\end{equation}
where $\wh\rho_{r_{k+1}}$ is defined in terms of $X_{k+1}$,
or else 
\begin{equation} \label{30.27n}
\theta_0 = \frac{3\pi}{2} \ \text{ and $\gamma^\ast_{r_{k+1}}$ has a free attachment} 
\end{equation}
(that is, for at least one of the points $\ell_\pm$ and one of the choices of $X_{k}$ or $X_{k+1}$).
Indeed, when our curves are attached to the points $\ell_\pm$ in the usual (non free) way,
the algorithm for choosing our nets of curves is the same, i.e., does not depend
on our choice of $X_{k}$ or $X_{k+1}$, and the variation of $\wh\rho_{r_{k+1}}$
is just the same as when we pick a different net $\gamma^\ast$ to start with; 
this matters no more than it did above. And we have seen earlier that there is no free attachment
when $\theta_0 = \pi$ and $j(r)$ is small enough, because this would mean a hanging curve.
Hence the claim.

Notice that when \eqref{30.27n} happens, say, with a free attachment at the point $z_+$, 
the proof of \eqref{30.20nn} also shows that ($\theta_0 = \frac{3\pi}{2}$ and)
\begin{equation} \label{30.28n}
|v_{1} + v_{2}| \leq C j(r_{k+1})^{1/2},
\end{equation}
where $v_j$ is the direction at $z_+$ of the geodesic $\rho(z_+,m_j)$. Since we also have
\eqref{30.20nn} at the two vertices $m_j$, we see that 
\begin{equation} \label{30.29n}
\text{the whole $r_{k+1}^{-1}\wh\rho_{r_{k+1}}$
is $C j(r_{k+1})^{1/2}$-close to a great circle.} 
\end{equation}

Anyway, let us return to \eqref{30.26n}; observe that if $0 \leq k_1 < k_2 \leq e$,
and if \eqref{30.27n} fails for $k_1 \leq k \leq k_2$, then we may sum \eqref{30.26n} and get that
\begin{equation}\label{30.30n}
d_\H(r_{k_1}^{-1}\wh\rho_{r_{k_1}},r_{k_2}^{-1}\wh\rho_{r_{k_2}}) 
\leq C \sum_{k_1 \leq k \leq k_2} \cE_k^{1/2} \leq C \cF_{k_1},
\end{equation}
where we set
\begin{equation}\label{30.31n}
\cF_{k_1} = \sum_{k_1 \leq k \leq e} \cE_k^{1/2}.
\end{equation}
We will estimate the $\cF_k$ more precisely in the next step, but let us start with some basic decay.
Recall that we chose our constants (such as $\varepsilon_4$) so that we can apply
the differential inequality \eqref{23.15} and \eqref{23.20}, and then even Theorem \ref{t23.1} 
or \ref{t23.2} (depending on $\theta_0$). When $\theta_0 = \frac{3\pi}{2}$, we get that
\begin{eqnarray} \label{30.32n}
f(90 R_k) &=& F(90 R_k) - \frac{2\pi}{3} 
\leq (C_V 10^{-k})^a f(90R) + C_V C_h 10^{-ka} R^\beta 
\nn\\
&\leq& C 10^{-ka} \varepsilon_4 + C C_h 10^{-ka} R^\beta 
\leq C 10^{-ka} \varepsilon_4 + C 10^{-ka} \varepsilon_3,
\end{eqnarray}
by Theorem \ref{t23.2}, \eqref{23.51}, and our fortunate assumption in 
Theorems \ref{t23.5} and \ref{t28.2} that $C_h R^{\beta} \leq \varepsilon_3$.
When $\theta_0 = \pi$, we get an even better result. Since
\begin{equation} \label{30.33n}
\int_{0}^{180R_k} h(r) \frac{dr}{r} \leq C C_h 10^{-\beta} R^{\beta}
\leq C 10^{-k\beta} \varepsilon_3
\end{equation}
 by the same assumption, we see that $\cE_k \leq C 10^{-ka} (\varepsilon_4 + \varepsilon_3)$
 and, summing over $l \geq k$, 
 \begin{equation} \label{30.34n}
\cF_{k} \leq C 10^{-ka/2} (\varepsilon_3 + \varepsilon_4)^{1/2}.
\end{equation}
Thus, under the assumptions of \eqref{30.30n}, 
\begin{equation}\label{30.35n}
d_\H(r_{k_1}^{-1}\wh\rho_{r_{k_1}}, r_{k_2}^{-1}\wh\rho_{r_{k_2}}) 
\leq C  10^{-k_1 a/2} (\varepsilon_3 + \varepsilon_4)^{1/2}.
\end{equation}
This is as small as we want, even for $k_1$ small. 

We are finally ready to use the small scale and prove that, in the case when 
$\theta_0 = \frac{3\pi}{2}$, our approximating cone $X_k = X(r_k)$ is never flat. 
For this, we shall first use Lemma \ref{t24.3} 
to show that $E$ is also close to a truncated cone of type $\bY$ in $B(0,R_e)$.

Let $\tau_3 > 0$ be small, and apply Lemma \ref{t24.3}, with $r = R_e$, $R = 2R_e$, 
and $\varepsilon = \tau_3$. If $\varepsilon_3$ and $\varepsilon_4$
are small enough, the assumptions \eqref{23.1n} and \eqref{23.2n}
with \eqref{23.8} hold by \eqref{23.50}, and \eqref{24.21} follows from 
\eqref{23.51} and our assumption that $C_h R^{\beta} \leq \varepsilon_3$.
Then we need to check \eqref{24.25}, but since we have \eqref{30.23n}, 
this is true as soon as we take $N$ large enough, depending on $\delta$, so that 
$\delta(\tau_3) R_e \leq 20d_0/21$. 
So Lemma \ref{t24.3} applies, and gives $E_0$, a truncated cone of type $\bY$ 
centered at $0$, such that 
\begin{equation}\label{30.36n}
d_{0, R_e}(E,E_0) \leq \tau_3.
\end{equation}

Recall that the approximating cone $X_e = X(r_e)$ was also such that 
$d_{0, R_e}(E,R_e) \leq \tau$, by \eqref{28.10}. Since $E_0$ has its two big
faces that make a $\frac{2\pi}{3}$ angle, and $\tau+\tau_3$ is as small as we want,
we deduce from \eqref{30.36n} that $X_e$ is of type $\bV$ (not $\bP_0$),
and that its two faces make an angle $\alpha(X_e)$ which is at most $\frac{2\pi}{3} + 10^{-3}$.

By \eqref{30.23n} and if $N$ is large enough, $R_e$ is quite large compared to $d_0$,
and then $r_e^{-1}\wh\rho_{r_e}$ is very close to the union of two great circles 
(we knew this already, because of \eqref{30.20nn}), and that make an angle 
$\alpha_e \in [\frac{2\pi}{3}-2 \cdot 10^{-3}, \frac{2\pi}{3} + 2 \cdot10^{-3}]$.
Now we prove by backwards induction that for $k \leq e$, $r_k^{-1}\wh\rho_{r_k}$
is very close to a union of two great circles (as usual), and that make an angle 
\begin{equation}\label{30.37n}
\alpha_k \in [\frac{2\pi}{3}-10^{-2}, \frac{2\pi}{3} + 10^{-2}].
\end{equation}
Indeed, as long as the free attachment event \eqref{30.27n} does not happen for $k+1$,
we have \eqref{30.35n} with $k_1 = k$ and $k_2 = e$, and then $\alpha_k$ is in 
the right range. But \eqref{30.27n} never happens unless $r_{k+1}^{-1}\wh\rho_{r_{k+1}}$
gets flat as in \eqref{30.29n}; this completes our induction. Hence \eqref{30.37n}
holds for all $k$, \eqref{30.35n} holds for $0 \leq k_1 < k_2 \leq e$, 
\eqref{30.29n} never happens, and (by the proof of \eqref{30.29n}), there is no free attachment 
associated to any of the $R_k$. This will simplify the discussion a little.

\msi
{\bf Step 5. We estimate the $\cE_k$.} We use two auxiliary sequences $\{ a_j \}$ and $\{ b_j\}$
to estimate the $\cE_k$. Let $T$ be a large integer, to be chosen soon, and set
\begin{equation} \label{30.38n}
a_j = f(90 R_{jT}) \ \text{ and } \  
b_j = \int_{0}^{180R_{jT}} h(r) \frac{dr}{r} = J(90R_{jT}),
\end{equation}
where we use the notation of \eqref{30a3}, and which we define only 
when $0 \leq jT \leq e$. We want to estimate $\sum_{j} a_j^{1/2}$
in terms of $a_0^{1/2}$ and $\sum_{j} b_j^{1/2}$. So we apply the proof of 
Theorem \ref{t23.2} or \ref{t23.1} (depending on the value of $\theta_0$), where we integrate the
differential inequality \eqref{23.15} or \eqref{23.20}, between the radii $r_1 = 90 R_{jT}$ and 
$r_2 = 90 R_{(j-1)T}$. We get that
\begin{eqnarray} \label{30.39n}
a_j &=& f(90 R_{jT}) \leq C_1 10^{-aT} f(90 R_{(j-1)T}) + C_2 \int_{0}^{180R_{(j-1)T}} h(r) \frac{dr}{r}
\nn\\ &=& C_1 10^{-aT} a_{j-1} + C_2 b_{j-1}.
\end{eqnarray}
where we don't care if $C_2$ depends on $T$. We choose $T$, depending on $C_1$, so large 
that $ C_1 10^{-aT} \leq 1/4$; then \eqref{30.39n} yields $a_j \leq a_{j-1}/4 + C_2 b_{j-1}$.
Then we take square roots, iterate, and get that
\begin{eqnarray} \label{30.40n}
a_j^{1/2} &\leq& \frac12 a_{j-1}^{1/2} + C b_{j-1}^{1/2}
\nn\\
&\leq& \frac14 a_{j-2}^{1/2} + C b_{j-1}^{1/2} + \frac{C}{2}\, b_{j-2}^{1/2}
\\
&\leq& \frac18 a_{j-3}^{1/2} + \frac{2C}{2}\, b_{j-1}^{1/2} + \frac{2C}{4}\, b_{j-2}^{1/2} 
+ \frac{2C}{8}\, b_{j-3}^{1/2} \, ,
\nn
\end{eqnarray}
and so on. Eventually
\begin{equation} \label{30.41n}
a_j^{1/2} \leq 2^{-j} a_0^{1/2} + 2C \sum_{0 \leq l \leq j} 2^{j-l} b_{l}^{1/2}.
\end{equation}
Then for $k \leq e$, we let $j$ be the integer such that $jT \leq k < (j+1)T$, and
\begin{eqnarray} \label{30.42n}
\cE_k^{1/2} &\leq& \Big(f(90 R_k) + \int_{0}^{180R_{k}} h(r) \frac{dr}{r}\Big)^{1/2} 
\leq (2a_j + b_j)^{1/2} \leq 2a_j ^{1/2}+ b_j^{1/2}
\nn\\
&\leq& 2^{1-j} a_0^{1/2} + b_j^{1/2} + 4C \sum_{0 \leq l \leq j} 2^{l-j} b_{l}^{1/2}
\end{eqnarray}
by the near monotonicity of $f$ and \eqref{30.41n}. Now we sum over $k \geq m$ 
to get an estimate for $\cF_m = \sum_{k \leq m} \cE_k^{1/2}$. 
Each estimate \eqref{30.42n} is used less than $T$ times, and becomes a sum over
indices $j \geq m/T -1$. The first term yields $C 2^{-m/T}a_0^{1/2}$.
The second term yields $C \sum_{j \geq m/T -1} b_j^{1/2}$. For the last term,
$b_l^{1/2}$ is multiplied by $C \sum 2^{l-j}$, where we sum over $j \geq l$ such that $j \geq m/T -1$.
When $l \leq m/T$, we get a geometric series that starts near $j \geq m/T -1$, with a sum
less than $C 2^{l -m/T}$. For $l \geq m/T$, we get a full sum bounded by $C$, but anyway the
contribution is similar to that of the second term. That is,
\begin{equation} \label{30.43n}
\cF_m \leq C 2^{-m/T} a_0^{1/2} + C \sum_{0 \leq l \leq m/T} 2^{l -m/T} b_{l}^{1/2}
+ C \sum_{ l > m/T} b_{l}^{1/2}
\leq C a_0^{1/2} + C \sum_{ l \geq 0} b_{l}^{1/2}.
\end{equation}
We translate back into integrals, and get that
\begin{eqnarray} \label{30.52n}
\cF_m &\leq& \cF_0 \leq C f(90 R)^{1/2} + C \sum_{k \geq 0} 
\Big(\int_{0}^{2^{-k}\cdot180R} h(r) \frac{dr}{r}\Big)^{1/2}
\nn\\
&=& C f(90 R)^{1/2} + C \sum_{k \geq 0} J(2^{-k}\cdot90R)^{1/2}
\leq C f(90 R)^{1/2} + C J_+(90R)
\end{eqnarray}
by \eqref{30a3}, and where for the last inequality we observe that we only
(defined and) used the $b_j$ when $jT \leq e$, hence
when $90R_{jT} \geq 90R_{e} \geq Nd_0$, by \eqref{30.23n}; so the 
restriction in the definition of $J_+$ in \eqref{30a3} is respected.

\msi
{\bf Step 6. We show that most $\rho_r^\ast$ lie close to a single truncated cone.}
We are now ready to check that when $\theta_0 = \pi$, most of our nets $\rho^\ast_r$
lie close to the half plane $H_0$ that contains $0$ and $L$, and when $\theta_0 = \frac{3\pi}{2}$
they lie close to some truncated set of type $\bY$ that we will choose.
We start with the easier first case. 

\begin{lem}\label{t30.3n} 
When $\theta_0 = \pi$,
\begin{equation} \label{30.53n}
\dist(z,H_0) \leq  Cr_e \cE_e^{1/2}  \ \text{ for } z \in \rho^\ast_{r_e},
\end{equation}
and then, for $0 \leq k \leq e$, 
\begin{equation} \label{30.54n}
\dist(z,H_0) \leq  C r_k \cF_0  \ \text{ for } z \in \rho^\ast_{r_k}.
\end{equation}
\end{lem}

First consider $r = r_e$, and recall that there is no free attachment. 
Thus the two arcs of geodesics that compose $\rho^\ast_r$ are the $\rho_\pm = \rho(\ell_\pm, m_1)$, 
where the $\ell_\pm$ are the points of $L \cap \S_r$.
We are going to use the flatness property \eqref{30.20nn} to estimate how far $\rho_\pm$ goes
from the plane $P_0$ that contains $0$ and $L$. Let us compute in the $3$-space that contains
$P_0$ and $m_1$ (and hence also the $\ell_\pm$), assume for the sake of the computation
that $r=1$, and choose coordinates where $P_0 = \big\{ z = 0 \big\}$, 
$L = \big\{ (t, -d_0, 0) \, ; \,  t\in \R\big\}$, $\ell_\pm = (\pm\sqrt{1-d_0^2}, -d_0, 0)$,
and $m_1 = (0,a,b)$, with $b \geq 0$, and where we used the possibility to take $w$ at equal 
distance from $\ell_-$ and $\ell_+$ to simplify the computation. The advantage of working with
$r_e$ is that with the present normalization, \eqref{30.23n} says that
$d_0 \geq (10N)^{-1}$ is bounded from below;
then it is also easy to see (because $E$ is close to a half plane that contains $0$) that $b$ is
small and $a > 0$.

We are interested in the unit tangent vector $v_\pm$ to $\rho_\pm$ at $m_1$, 
that points away from $m_1$. It must lie in the vector plane that contains $m_1$ and $\ell_\pm$,
and be orthogonal to $m_1$, hence be proportional to
\begin{equation} \label{30.45n}
\begin{aligned}
\xi_\pm &= \ell_\pm - \langle \ell_\pm, m_1 \rangle m_1= \ell_\pm + d_0 a m_1
\cr&
= (\pm\sqrt{1-d_0^2}, -d_0+d_0 a^2, d_0ab) = (\pm\sqrt{1-d_0^2}, -d_0 b^2, d_0ab).
\end{aligned}
\end{equation}
The square of the norm is $n^2 = 1-d_0^2+d_0^2b^2(b^2+a^2) = 1 - d_0^2(1-b^2) = 1 - d_0^2 a^2$,
which is the same for both signs, and close to $1$ (because $d_0 \leq N^{-1}$). Then
\begin{equation} \label{30.46n}
v_+ + v_- = n^{-1}(\xi_+ +\xi_-) = 2 n^{-1} d_0 b (0,b,a),
\end{equation}
whose norm is $2 n^{-1} d_0 b \geq d_0 b$. Whence, by \eqref{30.20nn} and forgetting about
our last normalization $r_e = 1$,
\begin{equation} \label{30.47n}
\dist(z,P_0) \leq b r_e \leq C r_e j(r_e)^{1/2} \leq C r_e \cE_{e}^{1/2}
\ \text{ for } z \in \rho^\ast_{r_e},
\end{equation}
by \eqref{30.25n}. And then, by \eqref{30.30n}, 
\begin{equation} \label{30.48n}
\dist(z,P_0) \leq b r_k \leq C r_k \cF_{0}
\ \text{ for $0 \leq k \leq e$ and } z \in \rho^\ast_{r_k}.
\end{equation}
We announced distances to $H_0$ instead of $P_0$, but this is the same because $m_1$
lies near $H_0$, not on the other side of $L$. 
\qed

\ms
Let us now consider $\theta_0 = \frac{3\pi}{2}$, and try to define a cone $Y$
that will work for \eqref{28.8}. In fact, we'll try to find it at the scale $r_{e}$,
and then use it a larger scales too, but let us discuss general radii $r_k$ for the moment.

Recall from the discussion that leads to \eqref{30.30n} that there is no free attachment for the $r_k$ 
(because the $j(r_k)$ are small). Then $\rho^\ast = \rho_{r_k}^\ast$ is composed of four main arcs 
$\rho_{j,\pm} = \rho(m_j, z_\pm)$, plus two short ones $\rho_{\pm} = \rho(z_\pm,\ell_\pm)$, 
where the last two may be reduced to a single point $z_\pm = \ell_\pm$.

We already observed for \eqref{30.20nn} that $j(r)$ controls the near minimality of the 
$\rho^\ast_r$, and in particular the angle that the two main arcs that compose $\rho^\ast_r$
make at their common endpoint $m_1$, but we also want to show that it control the angles at the 
points $z_\pm$. Suppose that $z_\pm \neq \ell_\pm$ (otherwise, let us not worry). 
Denote by $v_{\pm,j}$ the direction at $z_\pm$ of $\rho_{j,\pm}$, and by $v_{\pm}$
the direction at $z_\pm$ of $\rho_\pm$. We claim that
\begin{equation} \label{30.59n}
(r_k^{-1} \dist(z_\pm,L))^{1/2} \big|v_{\pm,1} + v_{\pm,2}+ v_{\pm}\big| \leq C j(r_k)^{1/2}.
\end{equation}
That is, we have the same estimates as for the central angles (in \eqref{30.20nn}),
but less precise because the proof only allows us to modify the tip of our competitors
by moving the point $z_\pm$ in a ball of size $C^{-1}\dist(z_\pm,L)$. 
The proof is as before, except that we replace \eqref{26.2} by the less performant \eqref{26.15}. 
Here $j(r)$ plays the role of $\sigma$, $\big|v_{\pm,1} + v_{\pm,2}+ v_{\pm}\big|$ 
is like $s$ in \eqref{26.4}, and $r_k^{-1} \dist(z_\pm,L)$ (correctly normalized) plays the role of $|z-\ell|$. 
The power $1/2$ on the left-hand side is unexpected (it comes from the fact that 
we can do a replacement in a tube rather than a ball), but plays for us. 
And for the worried reader, let us observe that we'll only use this estimate when $r_k^{-1} \dist(z_\pm,L)$
is reasonably large, where the subtle difference does not exist and the proof of \eqref{30.20nn} works too.

To complete the geometric information that we have, recall that by construction,
\begin{equation}\label{30.60n}
\text{the angles at $z_\pm$ of  $v_{\pm,1}, v_{\pm,2}, v_{\pm}$ are all at least $\pi/2$;}
\end{equation}
see the description below \eqref{25.9}. 

Now let us check that when $k = e$, the unpleasant factor
$r_k^{-1} \dist(z_\pm,L)$ is bounded from below, i.e., $z_\pm \neq \ell_\pm$ and 
\begin{equation}\label{30.61n}
\dist(z_\pm(r_e), L) \geq C^{-1} N^{-1} r_e,
\end{equation}
where we now mention explicitly the fact that $z_\pm$ comes from $r_{e}$.
Recall that $r_e \leq 90 R_e$ (see below \eqref{30.23n}), so we can apply 
the proof of \eqref{30a17} to the radius $r = 2r_e$; we find a minimal cone 
$Y_r \in \bY(0,\frac{21r}{20})$ such that $d_{0,r}(E,Y_r^t) \leq \tau$, 
with $\tau$ as small as we want. Now \eqref{30.16nn} says that the points of 
$\rho_{r_e}^\ast$ lie within $C j(r_e) r_e \leq C \cE_e r_e$
of $\gamma_{r_e}^\ast$, which itself lies in $E$, hence close to $Y_r$. Because
of this (and because they are geodesics), the four arcs $\rho_{j,\pm}$ lie very close to
the two main arcs of $Y_r^t \cap \S_{r_e}$. Which means that the two points 
$z_\pm$ are very close to the two points of $Spine(Y_r) \cap \S_{r_e}$,
where $Spine(Y_r)$ is the intersection of the three faces of $Y_r$.
But in turn, the two points of $Spine(Y_r) \cap \S_{r_e}$ cannot be too close to $L$,
because $Y_r \in \bY(0,\frac{21r}{20})$, with $r = 2r_e$; our claim \eqref{30.61n}
follows.

So \eqref{30.59n} says that for $r_e$, all the angles at the $z_\pm(r_e)$ 
are $C j(r_k)^{1/2}$-close to $\2$; we also control the angles at the $m_j$; altogether,
there is a cone $Y_0 \in \bY(0,r_e)$ such that 
\begin{equation} \label{30.62n}
\dist_\H (r_e^{-1} \rho_{r_e}^\ast, r_e^{-1}(Y_0^t \cap \S_{r_e}))
\leq C j(r_k)^{1/2},
\end{equation}
where as usual $Y_0^t$ denotes the truncated cone associated to $Y_0$.
We would like to use $Y = Y_0$, but maybe $Y_0$ does not lie in $\bY(L,R)$,
because $L \cap B(0,R)$ is not contained in a face of $Y_0$. So we need to discuss
a little more. Call $F_1$ and $F_2$ the two faces of  $Y_0$ that pass near $m_1(r_e)$ and 
$m_2(r_e)$ respectively, and $Y_0$ the remaining face, which contains 
$L \cap B(0,\frac{21r_e}{20})$ since $Y_0 \in \bY(0,\frac{21r_e}{20})$.
Thus $F_0$ is contained in the $2$-plane $P_0$ that contains $0$ and $L$, and one of our
concern is the angle of its boundary with $Spine(Y_0)$. 

But let us first assume that $Y_0 \in \bY(L,R)$, take $Y = Y_0$, and check that for $0 \leq k \leq e$,
\begin{equation}\label{30.63n}
d_{0,2r_k}(\rho_{r_k}^\ast, Y^t \cap \S_{r_k})
\leq C \cF_0.
\end{equation}
When $k = e$, this is just \eqref{30.62n}; and since all the angles are rather large,
the approximation is valid separately for the two long arcs and the two short ones. 
Then for a general $k$, \eqref{30.62n} shows that the the two long arcs of $\rho_{r_k}^\ast$ 
are $C \cF_{0} r_k$-close to the two arcs $F_j \cap \S_{r_k}$, $j=1,2$. We are left with
the two short geodesic arcs $\rho_{\pm}$, from $z_\pm = z_\pm(r_k)$ to $\ell_\pm$, which we want
to approximate by the corresponding arcs of $F_0 \cap \S_{r_k}$. Call these two arcs $\rho'_\pm$;
thus $\rho'_\pm = \rho(\ell_\pm, z'_\pm)$, where $z_\pm'$ is the point of 
$Spine(Y_0) \cap \S_{r_k}$ that lies close to $\ell_\pm$. 
Since $z_\pm'$ lies at the intersection of the two long arcs $F_1 \cap \S_{r_k}$ and 
$F_2 \cap \S_{r_k}$ and these arcs make a large angle at $z'_\pm$,
our estimate on the long arcs shows that $|z_\pm - z'_\pm| \leq C \cF_0 r_k$. Hence 
$\rho_{\pm} = \rho(\ell_\pm,z_\pm)$ and $\rho'_{\pm} = \rho(\ell_\pm,z'_\pm)$ are close to each other, 
and \eqref{30.63n} follows.

We are left with the case when $Y_0 \notin \bY(L,R)$, and in this case we want to replace
$Y_0$ with $Y = R(Y_0)$, for some rotation $R$ with a small angle. 
Consider the largest radius $r_0$; when we chose $r_0$ (just below \eqref{30.24n}), we made sure that
$r_0 \geq R$, so that if $Y_0 \notin \bY(L,R)$ as here, then $Y_0 \notin \bY(L,r_0)$ either.
Consider $r_0$, and still denote by $z_\pm'$ the point of 
$Spine(Y_0) \cap \S_{r_0}$ that lies close to $\ell_\pm$. 
Recall that both points $z_\pm'$ lie in the plane $P_0$ that contains $0$, $L$, and $F_0$,
but in the present case at least one of the points $z_\pm'$ lies on the other side of $Spine(Y_0)$
(compared to the projection on $P_0$ of most of $F_1$ and $F_2$, or the points $m_1(r_0)$ and $m_2(r_0)$,
for instance). Notice that this happens only at one of the two points, say, $z'_+$, because 
$Y_0 \in \bY(L,r_e)$ and hence $0$ lies on the same side of $P_0 \sm Spine(Y_0)$ as the projections
above. The other point $z'_-$ even lies further from $L$ as $0$.

Another way to say that $z_+'$ lies on the other side of $Spine(Y_0)$ is to say that seen from $\ell_+$,
it lies in the direction almost opposite to $v_{\pm,1}+v_{\pm,2}$. On the other hand, 
\eqref{30.60n} says that  $z_+$ (still seen from $\ell_+$), lies roughly in the direction of $v_{\pm,1}+v_{\pm,2}$.
But $|z_+ - z'_+| \leq C \cF_0 r_0$ for the same reason as before, and we deduce from this that 
$|z'_+- \ell_+| \leq C \cF_0 r_0$. Now let $R$ be the rotation that preserves $P_0$, is the identity on
$P_0^\perp$, and maps $z'_+$ to $\ell_+$. Notice that it moves the points very little, i.e., 
$|R(\xi)-\xi| \leq C \cF |\xi|$ for all $\xi$, and $Y = R(Y_0)$ lies in $\bY(L,r_0) \subset \bY(L,R)$
by construction. 

Return to a general $k\in [0,e]$, and set $z''_\pm(r_k) = R(z'_\pm(r_k)) = (r_k/r_0) R(z'_\pm)$; 
we see that $|z''_\pm(r_k)-z'_\pm(r_k)| \leq C  \cF_{0} r_k$, and also $Y \cap \S_{r_k}$
is $C  \cF_{0} r_k$-close to $Y_0$ in $\S_{r_k}$, so \eqref{30.63n} now holds for the same reason as in our first
case when we did not need to move $Y_0$.

This concludes our proof of \eqref{30.63n}, which we see as the correct analogue of Lemma \ref{t30.3n}.

\msi
{\bf Step 7. Our truncated cone approximates $E$ well in the exterior annulus.}
When $\theta_0 = \frac{3\pi}{2}$, we just constructed a truncated cone $Y^t$, and 
in order to unify the notation, let us also denote by $Y^t$ the half plane $H_0$ when $\theta_0 = \pi$.
This way, we also have \eqref{30.63n} when $\theta_0 = \pi$. Indeed, by \eqref{30.54n} every point of 
$\rho_k^\ast$ lies close to $Y^t \cap S_{r_k}$, but then conversely every point of 
$Y^t \cap S_{r_k}= H_0 \cap S_{r_k}$ lies close to $\rho_k^\ast$, because $\rho_k^\ast$ is just
the concatenation of two geodesics that end on $L$ and $S_{r_k}= H_0 \cap S_{r_k}$ is simple too.
We are ready to see that $E$ stays quite close to $Y^t$ in the region $B(0,R) \sm B(0,Nd_0)$.

\begin{lem}\label{t30.4n} Let $Y^t$ be as above (thus $Y^t = H_0$ when $\theta_0 = \pi$). 
Set $A_{00} = B(0,2R) \sm B(0,2Nd_0)$; then
\begin{equation} \label{30.64}
\dist(z,E) \leq C (\cF_0 + \cE_0^{1/4}) |z|
\ \text{ for } z\in Y^t \cap A_{00}
\end{equation}
and 
\begin{equation} \label{30.65}
\dist(z,Y^t) \leq C (\cF_0 + \cE_0^{1/4}) |z|
\ \text{ for } z\in E \cap A_{00}.
\end{equation}
\end{lem}

We will follow the argument of Section \ref{S20}. First we claim that when $r\in \cR_k$ for some
$k \in [0,e]$, 
\begin{equation} \label{30.66}
\dist_{0,2r}(\gamma^\ast_r, Y^t \cap S_{r}) \leq C \cF_0 + C j(r),
\end{equation}
where $\gamma^\ast_r$ is our initial curve in $E \cap \S_r$.
Indeed, \eqref{30.16nn} says that $d_{0,2r}(\rho^\ast_r , \gamma^\ast_r) \leq C j(r)$,
then \eqref{30.19nn} says that $\wh \rho^\ast_r$ is quite close to $\frac{r}{r_k} \rho^\ast_{r_k}$.
We now know that for $j(r)$ small there is no free attachment and we control the angle between the 
various pieces of $\rho^\ast_{r_k}$, so in fact the proof of \eqref{30.63n} also yields
\begin{equation}\label{30.67}
d_{0,2r}(\rho_{r}^\ast, Y^t \cap \S_{r})
\leq C \cF_0 + C j(r)^{1/2}
\end{equation}
(see \eqref{30.62n} in particular), and \eqref{30.66} follows.
As we did for Lemma \ref{t20.1}, we first restrict our attention to the set $\cR_k^\sharp$ of radii 
$r \in \cR_k$ such that $j(r) \leq C \cE_k^{1/3}$. Thus by Chebyshev
\begin{equation} \label{30.68}
|(10^{-1}R_k,90R_k) \sm \cR_k^\sharp| 
= |\cR_k \sm \cR_k^\sharp| \leq C \cE_k^{-1/3} \int_{\cR_k} j(r) dr \leq C \cE_k^{2/3} R_k
\end{equation}
by \eqref{30.25n} (or by the proof of \eqref{19.21}). Thus every radius $\rho \in (10^{-1}R_k,90R_k)$
lies within $C \cE_k^{2/3} R_k$ of a radius $r\in \cR_k^\sharp$. For each point 
$z\in Y^t \cap B(0,90R) \sm \ol B(0, 10^{-1}R_e)$, we can find $k \in [0,e]$
such that $\rho = |z|$ lies in $(10^{-1}R_k,90R_k)$, then $r \in \cR_k^\sharp$ such that 
$|r-\rho| \leq C \cE_k^{2/3} R_k$, then $z_1 \in Y^t \cap \S_r$ such that $|z_1-z| \leq C \cE_k^{2/3} R_k$, 
and finally by \eqref{30.66} a point
$z_2 \in \gamma^\ast_r \subset E$ such that $|z_2-z_1| \leq C (C \cF_0 + C j(r)) r \leq C (\cF_0 + \cE_k^{1/3}) r$.
Thus 
\begin{equation} \label{30.69}
\dist(z,E) \leq C (\cF_0 + \cE_k^{1/3}) |z|
\end{equation}
for $z\in Y^t \cap B(0,90R) \sm \ol B(0, 10^{-1}R_e)$.
We also need to evaluate $\dist(z, Y^t)$ when $z\in E \cap B(0,90R) \sm \ol B(0, 10^{-1}R_e)$.
We proceed as in Lemma \ref{t20.2}. We start with the points $z$ that lie in some $r \in \cR_k^\sharp$, 
(so that $j(r) \leq C \cE_k^{1/3}$), and in addition lie in corresponding set $\gamma_r^\ast$. For those 
we can use \eqref{30.66} and get that $\dist(z,Y^t) \leq C (\cF_0 + \cE_k^{1/3}) r$. Then we evaluate 
the measure of the piece of $E \cap B(0,3R_k) \sm 10^{-1}R_k)$ for which \eqref{30.69} fails
(because $z\in E \cap \S_r \sm \gamma_r^\ast$, or because $j(r)$ is too large, or because the co-area
formula does not cooperate), and get a set of measure at most $C R_k^2 \cE_k^{1/4}$ (see below \eqref{21.19n}).
Finally, we use the Ahlfors regularity of $E$ and get that
\begin{equation} \label{30.70}
\dist(z,Y^t) \leq C (\cF_0 + \cE_k^{1/4}) |z|
\end{equation}
for $z\in E \cap B(0,2R) \sm \ol B(0, 10^{-1}R_e)$. See Lemma \ref{t20.1}.

The lemma follows because $10^{-1}R_e \geq 2Nd_0$ by \eqref{30.23n}, and 
$\cE_k \leq C\cE_0$ by its definition \eqref{30.24n}.

\msi
{\bf Step 8. Decay in the intermediate and small ranges, in the flat case.}

So far we assumed that $90R \geq N d_0$ (for some large $N$), as in \eqref{30a20},
and we proved in Lemma \ref{t30.4n} that $E$ is close enough to a truncated $\bY$-set $Y^t$
in the annulus $A_{00} = B(0,2R) \sm B(0,2Nd_0)$.

We need to complete this description, both to allow the radii $R \in [2d_0, N d_0]$
and to get the same control (with the same set $Y^t$) in the interior ball $B(0,3R)$, say.

We shall first do this in the flat case (when $\theta_0 = \pi$). If $E$ were minimal, we could deduce
a control on $B(0,3R)$ from what happens in $A_{00}$, or even $\d B(0,3R)$, by some form of
maximum principle (compare $E \cap B(0,3R)$ with its projection on a small convex neighborhood
of $H_0 \cap B(0,3R))$. Probably there is a way to do a similar argument for almost minimal sets,
but it does not seem so pleasant, and the author fears that we would need the co-area formula, in much the
same way as above, to control angles of tangents and integrate; the case when 
$\theta_0 = \frac{3\pi}{2}$ would probably be even worse. Instead we shall return to our main decay estimate
for $f(r)$, extend it to the present situation where \eqref{25.3} or \eqref{25.4} fails, and
prove a similar decay anyway. And then we shall proceed as in the first step, say that a slightly bad
geometric configuration, even for these intermediate radii, implies a corresponding lower bound for $f'(r)$,
which cannot happen too often under the present assumptions.

Why didn't we do this earlier? A first (bad) excuse is that the author did not think that this would be needed, 
and to be fair, what we are going to do now does not change the main $C^1$ results of the present paper, but only
gives more precise estimates (including a quantitative $C^{1+\beta}$ estimate). 
The second excuse is that this slightly simplifies the apparent structure of the proof. In general, the geometric
situation for the intermediate radii looks a little bad for our earlier construction of competitors to run smoothly; 
here we shall use some extra information on $E$ to make things easier.

Let us now assume that $\theta_0 = \pi$, and merely assume that $R \geq 2d_0$
(instead of \eqref{30a20}). 
We shall use the following fact to control the intersection $E \cap \S_r$ for intermediate radii.

\begin{lem}\label{t30.5n}
The set $E$ coincides, in $B(0,10R)$, with the graph over $H_0$ of some $C^1$ and $10^{-1}$-Lipschitz
function $\varphi : H_0 \to H_0^\perp$ such that $\varphi(x)=0$ for $x\in L$.
\end{lem}

Here $H_0^\perp$ is the $(n-2)$-plane orthogonal to $P_0$, the $2$-plane that contains $H_0$.
The graph of $\varphi$ is $G_\varphi = \big\{ x+\varphi(x) \, ; \, x\in H_0 \big\}$.
The proof of the lemma will use the $C^1$ and Lipschitz part of Theorem \ref{t29a.1}, which
will be proved later but independently. In fact, the proof of Theorem~\ref{t29a.1} will never
involve decay estimates for balls that are not centered on $L$. Let us explain how we deduce 
Lemma \ref{t30.5n} from Theorem \ref{t29a.1}, and rapidly sketch the part of Theorem \ref{t29a.1}
that we need.

We start from \eqref{23.53} in the preparation Step 1, which says that 
$d_{0,r}(E,H_0) \leq \tau$ for $2d_0 \leq r \leq 180R$. Here $\tau$ is as small as we want,
and with our usual assumption that the gauge function $h$ is small enough, this is enough to
apply Theorem \ref{t29a.1} and get the conclusion of the lemma. The arguments for Theorem \ref{t29a.1}
are not that complicated. Since $\tau$ is as small as we want, a compactness argument shows
that in $B(0,160R)$, $E$ is also close to $H_0$ in measure; this means in particular that for 
$\xi \in L\cap B(0,60R)$, $(60R)^{-2}\H^2(E \cap B(\xi,60R))$ is as close to $\pi$ as we want.

There is an additional argument to show that $L \cap E \cap B(0,10R)$ is not empty
(otherwise we could easily cut out a big part of $E$), and a variation of the same argument
shows that $L \cap B(0,60R) \subset E$. All this allows us to to apply the results of Part III
(or Theorems~\ref{t1.3} and \ref{t1.4}) and get a good flatness control on $E$ in balls centered on 
$L \cap B(0,60R)$. This also mechanically gives a good control on balls $B(x,r)$ centered
on $E$ such that $10^{-2} r \leq \dist(x,L) \leq 100 r$, and for the other ones, we 
start from $B(x,\dist(x,L))$ and then use the usual regularity theorem with no boundary $L$.
\qed

Let us now consider radii $r$ such that
\begin{equation} \label{30.71}
d_0 < r \leq \min(10R, 10N d_0)
\end{equation}
(we could use $N=2$ here, because we are in the flat case, but let us not bother).
By Lemma~\ref{t30.5n} and the implicit function theorem, and since $\S_r$ is perpendicular to 
$H_0$ where it meets it, $E \cap \S_r $ is a nice $C^1$ curve that starts from one point
$\ell_-(r) \in L$, runs very close to $H_0 \cap \S_r$, and ends at the other point $\ell_+(r)$
of $L \cap \S_r$. It is also a small Lipschitz graph over $H_0 \cap \S_r$ (its tangent lies close to
the direction of $H_0$, and belongs to the tangent hyperplane to $\S_r$). 
That is, even when $r$ is barely larger than $d_0$ and $\S_r$ is nearly tangent to $L$,
it stays nicely transverse to $E$ and nothing bad happens.

Set $\gamma_r^\ast = E \cap \S_r$; this is consistant with the notation of Section \ref{S25},
but the situation is simpler now. When $r$ is barely more than $d_0$, say, when
\begin{equation} \label{30.72}
d_0 < r  < \frac{5d_0}{4},
\end{equation}
let us modify our construction and now cut $\gamma_r^\ast$ in roughly three equal parts (instead of two), 
with two intermediate points $m_1$ and $m_2$; we do this because in this case the length of 
$\gamma_r^\ast$ may be close to $2\pi r$, and we want arcs of lengths significantly less than $\pi r$ 
to apply the construction and estimates of Section \ref{S7}. 
When $r \geq \frac{5d_0}{4}$, let us just proceed as above and cut $\gamma_r^\ast$ only once, near the middle.

Then we do the construction of Section \ref{S25}, i.e., build a competitor where,
instead of taking the union of the cone over $\gamma^\ast_r$ and the triangle $T(r)$
(as we would do in order to prove the near monotonicity of $f$), we replace the cones over the two or three 
pieces $\gamma_{i,r}$ of $\gamma_r^\ast$ with harmonic graphs that end with a small flat plate.

Denote by $\rho_{i,r}$ the geodesic with the same endpoints as $\gamma_{i,r}$, and set
$\rho_r^\ast = \cup_{i=1}^3\rho_{i,r}$ (when \eqref{30.72} holds; otherwise we stop at $i=2$). 
We get the same control as before on the distance 
between $\gamma_{i,r}$ and $\rho_{i,r}$, in terms of $f'(r)$.
That is, we still get that \eqref{30.16nn} holds, with $j(r)$ as in \eqref{30a21} or \eqref{18.24}.
And in addition the same argument as usual, where we modify the tips of our competitor to get an even 
better one when the geodesics are not aligned at some $m_i$, also gives the same control as in 
\eqref{30.20nn}. That is, if we still denote by $v_{j,\pm,r}$ the two unit 
tangent vectors of $\rho^\ast_{r}$ at the point $m_i$, we still have that
$|v_{j,+,r} + v_{j,-,r}| \leq C j(r)^{1/2}$. 

Let us first use this information to control the geometry when
\begin{equation} \label{30.73}
\frac{5d_0}{4} \leq r \leq \min(10R, 10 N d_0).
\end{equation}
Then the proof of Lemma~\ref{t30.3n} yields 
\begin{equation} \label{30.74}
\dist_{0,2r}(E \cap \S_r, H_0 \cap \S_r) \leq C j(r)^{1/2}.
\end{equation}
When instead $d_0 < r < \frac{5d_0}{4}$, the proof does not work 
(because we cannot get a good control on the angle with $P_0$ of a great circle that almost contains 
$\rho^\ast_{r}$ when $r$ is too close to $d_0$), but fortunately we can use the same trick as before, i.e., 
use transverse curves in $E$ and the co-area formula to control the variations of the 
great circles $\wh \rho_{i,r}$ that contain the $\rho_{i,r}$. 
That is, the analogue of \eqref{30.30n} holds in the present case too, and 
allows us to deduce from \eqref{30.74} (for $r$ close to $\frac{5d_0}{4}$) that for our $r$
\begin{equation} \label{30.75}
\dist_{0,2r}(E \cap \S_r, H_0 \cap \S_r) \leq C j(r)^{1/2} + C\cE^{1/2},
\end{equation}
where $\cE$ is still as in \eqref{30.18nn}, and could even be replaced by 
$f(90d_0) + \int_{0}^{180d_0} h(r) \frac{dr}{r}$.

We can also do like this for radii $r \in (\frac{d_0}{10}, d_0)$.
For these we proceed as usual, but since $L$ does not meet $B(0,r)$,
we just do as for the interior regularity result, cut $E \cap \d B(0,r)$ into
three roughly equal parts, construct the $\rho_{i,r}$, and get good estimates on their
angles. As when \eqref{30.72} holds, we only know a priori that the three $\rho_{i,r}$
lie $C j(r)^{1/2} r$-close to some great circle, and in order to show that this circle lies close to
the plane $P_0$, we use the transverse curve and the co-area estimate that relates their variations
to the decay of $f$, to get that \eqref{30.75} holds also for these $r$. 

At this point we obtained a stronger analogue of Lemma \ref{t30.3n}, where \eqref{30.53n}
still holds for $\frac{d_0}{10} < r < \min(10R, 10 N d_0)$. This time we managed to include
unpleasant case when \eqref{30a20} fails, which we had left alone before.
Now we also have the analogue of Lemma \ref{t30.4n} for the same radii, 
again with the same proof as above.

We still need to control $E \cap B(0,d_0/2)$, say. For this, the analogue of 
Theorems \ref{t18.1} and \ref{t23.5} is valid, with a simpler proof (see \cite{C1}), 
and we get a plane $P$ through the origin such that 
\begin{equation} \label{30.76}
\dist_{0,d_0/2}(E,P) \leq C \Big[f(100d_0) + \int_{0}^{200d_0} \frac{h(t)dt}{t} \Big]^{1/4}.
\end{equation}
In addition, we know that on the outside rim $B(0,d_0) \sm B(0,d_0/2)$, $E$ is quite close to 
$H_0$, by the analogue of Lemma \ref{t30.4n}; this provides the desired extension of Lemma \ref{t30.4n}
to the ball, and concludes the proof of Theorem \ref{t23.5} and Remark \ref{r30a2}.

\msi
{\bf Step 9. Decay in the remaining ranges, when $\theta_0 = \frac{3\pi}{2}$.}

We shall now complete the proof in the remaining case when $\theta_0 = \frac{3\pi}{2}$.
We start with the fact that for $0 < r \leq 200R$,
\begin{equation} \label{30.77}
f(r) + \int_0^{2r} \frac{h(t)dt}{t} = f(r) + J(r) \leq f(r) + J(200R) \leq C \varepsilon_4
\end{equation}
by \eqref{30a3} (the definition of $J(r)$), \eqref{23.50}, and \eqref{23.51}, 
and where $\varepsilon_4$ is as small as we want. We shall first restrict to $r$ such that 
\begin{equation} \label{30.78}
\frac{5d_0}{4} \leq r \leq 190\min(R, 10N d_0),
\end{equation}
and then \eqref{30.77} allows us to apply Lemma \ref{t24.3}, with $\delta$ small enough (depending on $N$),
and then $\varepsilon_4$ small (depending on $\delta$); we get that there is a set $Y \in \bY(L,r)$ such that 
\begin{equation} \label{30.79}
d_{0,r}(E,Y^t) \leq \varepsilon_5
\end{equation}
for the corresponding truncated cone $Y^t$, where $\varepsilon_5$ is also as small as we want.
In the definition of $\bY(L,r)$, there is the fact that the $L \cap B(0,r)$ is contained in a single face of
$Y$, and in particular does not meet the spine $L_Y$ of $Y$ (the singular set of $Y$). Since
$Y$ is centered at $0$ and $d_0 = \dist(0,L) \geq (1900N)^{-1} r$ by \eqref{30.78},
we see (for instance by drawing the two lines $L$ and $L_Y$ in the plane that contains them) that 
\begin{equation} \label{30.80}
\dist(L \cap B(0,r/2), L_Y) \geq c r,
\end{equation}
where $c > 0$ depends on $N$ but this does not matter. If $\varepsilon_5$ is small enough,
we can apply the usual regularity results to prove that in $B(0,r/3)$ and when we stay at distance 
at least $cr/3$ from $L$, $E$ is a $C^1$ version of the $\bY$-set $Y^t$.
And at distance less  than $2cr/3$, we can apply Theorem \ref{t29a.1}, as we did for Lemma\ref{t30.5n}, 
but on smaller balls, to show that $E$ coincides with the graph over $H_0$ of some $C^1$ and 
$10^{-1}$-Lipschitz function $\varphi : H_0 \to H_0^\perp$, with as above $\varphi(x)=0$ for $x\in L$.

This is good, because this shows that for $d_0 < s < r/3$, $E \cap \S_s$ is also a $C^1$ version
of $Y^t \cap \S_r$. That is, $E \cap \S_s$ is composed of two long arcs that leave from nearly
opposite points $m_\pm \in E \cap \S_r$ (near the points of $L_Y \cap \S_r$) with angles close to 
$\frac{2\pi}{3}$, plus two short ones that go from those points to the two points $\ell_\pm$ of $L \cap \S_s$.
We stop at $d_0$, because when $s < d_0$ the sphere no longer meets $L$, but when $s$ approaches
$d_0$ from above, nothing bad happens, the two short arcs simply get longer and go to points $\ell_\pm$
that lie very close to each other, as we already saw in the flat case.

This allows us to implement the construction that we did the previous sections, to construct competitors
for $E$ in $\ol B(0,s)$. We thus get for each such $s$ a net $\gamma_s^\ast$, composed of two ``short''
curves $\gamma_{0,\pm,s}$ and four ``long'' ones (we cut the two initial long arcs in two parts, 
near the middle, as usual) $\gamma_{i,\pm,s}$, $i = 1, 2$; 
then we construct $6$ Lipschitz graphs with flat tips, add the triangle 
$T(s)$ as usual, and also try to improve it near the tips, depending on the angles that the six geodesics 
with the same ends as the six curves $\gamma$ (we call these geodesics the $\rho_{i,\pm,s}$, $0 \leq i \leq 2$)
make with each other. The initial Lipschitz graphs are enough to get the near monotonicity of $f$,
so lower bounds on the defect of angles yield lower bounds on $f'(s)$. 

More precisely, denote by $\alpha = \alpha(s)$ the largest angle defect, i.e., the largest of the differences between
between the angle between two curves $\rho$ of $\rho^\ast_s$ and the expected value at that point
(i.e., $\pi$ or $\frac{2\pi}{3}$). Then, due to the fact that in spite of their code names,
the ``short'' geodesics are never too short, by \eqref{30.80}, we get that $\alpha \leq C j(s)^{1/2}$, 
as in \eqref{30.20nn} for instance. 

Now we can recover some of the geometry of $\rho^\ast_s$ in terms of $\alpha$.
Start from the three geodesics of $\rho^\ast_s$ that leave from $z_-$ (we called them
$\rho_{i,-,s}$), and denote by $v_{i}$ the tangent vector to $\rho_{i,-,s}$ at $z_-$;
the three $v_i$ make angles that are $\alpha$-close to $\frac{2\pi}{3}$. If we followed
the two geodesics $\rho_{i,-,s}$, $i=1,2$, they would meet again, with the same angles,
at $-z_-$. But we allow them to turn near their middle (when they become the $\rho_{i,+,s}$), 
but by less than $\alpha$. They still meet at $z_+$ (this was a property of $E \cap \S_s$), 
transversally,  and in addition $|z_+ + z_-| \leq C \alpha s$, and they make an angle
that is still $C\alpha$-close to $\frac{2\pi}{3}$ at $z_+$. Here we skip some of the details,
but we made similar computations in Sections \ref{S26} and \ref{S27}, in situation that were a little
more complicated because we did not have \eqref{30.78}.
Said differently, there is a $\bY$-set $Y_\rho$ centered at $0$, whose spine contains $z_-$,
and such that the four $\rho_{i,\pm,s}$, with $i=1,2$, stay $C \alpha s$-close to two of the three 
arcs that compose $Y_\rho \cap \S_s$. 

Let us now look at the two remaining geodesics $\rho_{0,\pm,s}$ as they leave from $z_\pm$
(with angles nearly $\frac{2\pi}{3}$ with the other ones); recall that they meet $L$
at the points $\ell_\pm$. By rotating $Y_\rho$ by less than $\alpha$ along its spine, we may
assume that it contains the arc $\rho_{0,-,s}$ in one of its faces $F$. Then on the opposite side, since 
$\rho_{0,-,s}$ also lies close to $Y_\rho$, we see that $F$ also passes within $C\alpha s$ of $\ell_+$.

Let us assume that $s \geq \frac{5d_0}{4}$, so that $|\ell_+ - \ell_-| \geq C^{-1} s$.
This way, the distance from $\ell_-$ to $F$ also controls the angle between $P_0$
(the plane that contains $0$ and $L$) and the plane $P_\rho$ that contains $F$
(recall that $P_\rho$ contains $0$ and $\ell_-$). Thus this angle is small, and we can find a rotation,
which is $C\alpha$-close to the identity, and that fixes $0$ and $\ell_-$ and sends $F$ to $P_0$. 
The image $Y_s$ of $Y_\rho$ by this rotation has the extra advantage that it lies in $\bY(L,r)$
(because $F$ is mapped into the plane $P_0$, and the spine of $Y_\rho$ was far from 
$L$ in $B(0,r/2)$, by \eqref{30.80}). Notice that we do not say that
$Y_s$ is the cone $Y_0$ of type $\bY$ that contains $L$, because the image of $Y_0$ by a small 
rotation that fixes $P_0$ is also allowed. Since we only moved $Y_\rho$ a little, it is still true that 
every point of $\rho^\ast_s$ lies within $C \alpha s$ from $Y_s$. In fact, given the structure
of the two sets, we even get that
\begin{equation}\label{30.81}
d_{0,2s}(\rho^\ast_s, Y_s^t \cap \S_s) \leq C \alpha,
\end{equation}
where $Y_s^t$ is the truncated cone associated to $Y_s$.

By the same argument a for \eqref{30.75} (and simplified from the same one in earlier proofs,
because here we know that $E \cap \S_s$ is equal to $\gamma_s$), we also deduce from 
\eqref{30.81} and the fact that $\alpha \leq C j(s)^{1/2}$ that
\begin{equation}\label{30.82}
d_{0,2s}(E \cap  \S_s, Y_s^t \cap \S_s) \leq C j(s)^{1/2}.
\end{equation}
This takes care of most $s$ such that $s \geq \frac{5d_0}{4}$. For
$d_0 < s < \frac{5d_0}{4}$, we still have a good description of $\rho^\ast_s$ in terms of
$Y_\rho$, and we can use transverse curves and the co-area theorem to control the 
variations of $Y_\rho$, and show that it lies within $C j(s)^{1/2} + C\cE^{1/2}$
of some fixed $Y_{s_0}$ (chosen with $s_0 \in (\frac{5d_0}{4}, 2d_0)$, by a Chebyshev argument),
as in \eqref{30.75}. Thus \eqref{30.81} even holds for $d_0 < s < 60\min(R, 10N d_0)$
(we divide the bound of \eqref{30.78} by $3$ because $s \leq r/3$), and even with a fixed set
$Y_{s_0} \in \bY(L, r_0)$, with $r_0 = 60\min(R, 10N d_0)$ if we choose $s_0$ large enough or use 
\eqref{30.79} and the fact that all $Y_s$ have a face contained in $P_0$ and lie close to $Y$.

As before, this controls most of $E \cap B(0,r_0) \sm B(0,d_0)$ (because $j(s)$ may be too large for 
some $s$), but then we can use the local Ahlfors regularity of $E$ to control the rest.
Finally, for $E \cap B(0,d_0)$, we use the same argument, but without the boundary $L$, to control
the variations of an approximating set of type $\bY$ that approximates $E \cap \S_s$
for $s < d_0$. The argument is the same as in the flat case; we just use the standard regularity result 
near a cone of type $\bY$, far from the boundary. At the end of the game, we find that 
\begin{equation}\label{30.83}
d_{0,r_0}(E , Y_{s_0}^t)  \leq C \Big[f(3r_0) + \int_{0}^{6r_0} \frac{h(t)dt}{t} \Big]^{1/4},
\end{equation}
as in \eqref{30.76}; this is enough because $f(3r_0)$ and $\int_{0}^{6r_0} \frac{h(t)dt}{t}$
are controlled by the same quantity for $R$, by \eqref{30.77}.

Recall that we had also established a similar control on the exterior annulus $A_{00}$;
with this last estimate, we end the proof of Theorem \ref{t23.5}, Remark \ref{r30a2}, 
and Theorem \ref{t28.2}.
\qed

\begin{rem} \label{R28.4}
We proved a little more than what we said: we proved that $E$ can be approximated well
by truncated $\bY$-sets in every ball $B(0,R')$, $R' \leq R$, with uniform, and even improving estimates,
leading to the existence of a tangent $\bY$-cone at $0$ (which we already knew), but which is also
quite close to the the $\bY$-cone that contains $Y^t$.
\end{rem}

At this point we have almost all the information needed to have a good description of $E$ near a
point of type $\bV$. We shall summarize this and similar local descriptions in the next part.

\vfill\eject
\part{Geometric descriptions of $E$ near some cones}

\ms
In this part we start the local description of sliding almost minimal sets of dimension $2$ near
a one-dimensional, smooth sliding boundary $L$. We will concentrate on the case when $L$ is
a line, and only explain in Section \ref{S31} how to extend our results to the case when 
$L$ is a $C^2$ curve, say.

In the next few sections, we give ourselves a (reduced) sliding almost minimal set $E$ of dimension $2$,
associated to boundary $L$ which is a straight line through the origin, assume that in the ball $B(0,R)$, 
$E$ is close enough to a given sliding minimal cone $X$, that the gauge function $h$ is small enough,
and we want a good description of $E$ in, say $B(0,R/2)$. In a few good cases, we will see that $E$
is just $C^1$-equivalent to $X$ in the smaller ball, but in more interesting cases (in particular when $X$ is a sharp
set of type $\bV$), we will get a good description, but where $E$ may have a different topology.

What we can get depends on the cone $X$; for simple cones we get a good result, and unfortunately
for some cones, such as $\bY$-sets with the spine $L$, we have reasonable conjectures but no proof.

In Sections \ref{S29a}-\ref{S36} we take the possible cones one after the other, and say what local
regularity result we have (or not) near these cones.
Then in Section \ref{S30} we will complete the verification of the full length properties that were
announced throughout this paper, and in Section \ref{S31} we say why $L$ can be replaced with a smooth 
curve in all our regularity results. We also decided to add a Section \ref{S32} where we check that 
sets of type $\H$ or $\bV$ are sliding minimal (the verification was not done yet).
% no story about convex domains yet

\section{Local regularity of $E$ near a half plane}
\label{S29a}

We start our list of local regularity results with the description of $E$ when it lies close to a
half plane. We state our main assumptions so that we can use them in later sections, with different
cones $X$. Let $L$ be a line in $\R^n$, that contains the origin, and suppose that
\begin{equation} \label{29a.1}
\begin{aligned}
&\text{$E$ is a reduced almost minimal set in $B(0,R)$, with sliding boundary $L$}
\cr&\hskip3cm \text{and with the gauge function $h$,}
\end{aligned}
\end{equation}
where we shall also assume that 
\begin{equation} \label{29a.2}
h(r) \leq C_h r^{\beta} \ \text{ for } 0 < r < R
\end{equation}
for some choice of power $\beta \in (0,1]$ and $C_h \geq 0$, and that 
\begin{equation} \label{29a.3}
C_h R^{\beta} \leq \varepsilon_0
\end{equation}
for some $\varepsilon_0 > 0$ that we can chose very small, depending on $n$ and $\beta$.
The main geometric assumption that complements this is that there is a sliding minimal cone $X$,
centered at $0$, such that 
\begin{equation} \label{29a.4}
d_{0,R}(E,X) \leq \varepsilon_0,
\end{equation}
where $d_{0,R}$ is the local Hausdorff distance of \eqref{1.13}. In this section we are interested
in the special case when $X \in \bH$ is a half plane bounded by $L$; this case was partially treated in 
\cite{Mono}, but we add the $C^1$ nature of the estimate here.
 
\begin{thm}\label{t29a.1} 
There is a constant $a>0$ that depends only on $n$ and $\beta$ and, for each small $\tau > 0$, 
a constant $\varepsilon_0 > 0$, that depends only on $n$, $\beta$, and $\tau$,
with the following properties.
Let $E$, $h$, $R$, satisfy \eqref{29a.1}-\eqref{29a.3}, and assume that \eqref{29a.4}
holds for some half plane $X \in \bH$ bounded by $L$. Then $E$ coincides in $B(0,R/10)$
with the graph of some $C^{1}$ function $\varphi : X \to X^\perp$. In addition,
$\varphi(x) =x$ for $x\in L$, $\varphi$ is $\tau$-Lipschitz, and 
\begin{equation} \label{29a.5}
\Angle(T_x E, T_y E) \leq \tau |x-y|^a R^{-a} \ \text{ for } x,y \in E \cap B(0,R/10),
\end{equation}
where $T_x$ denotes the tangent plane to $E$ at $x\in E \cap B(0,R/10)$.
\end{thm}

Here $X^\perp$ is the $(n-2)$-space perpendicular to the plane that contains $X$.
When $x\in L$, $T_x$ is only a half-tangent plane. And the simplest definition of 
$\Angle(T_x E, T_y E)$ is probably $||\pi_x-\pi_y||$, where $\pi_x$ is the orthogonal 
projection on $T_x E$, and similarly for $\pi_y$.

There is nothing special with the constant $1/10$; any constant smaller than $1$ would work, at the 
price of taking $\varepsilon_0$ smaller and complicating the argument.

We may also replace the assumption \eqref{29a.4} in Theorem \ref{t29a.1} with the density assumption
\begin{equation} \label{29a.6}
\dist(0,E) \leq \varepsilon_0 R \ \text{ and } R^{-2} \H^2(E \cap B(0,R)) \leq \frac{\pi}{2} + \varepsilon_0 ;
\end{equation}
as we shall see, under the other assumptions, \eqref{29a.4} with $X \in \bH$ and \eqref{29a.6}
are essentially equivalent to each other (modulo taking a different $\varepsilon_0$ and a slightly smaller $R$).

Most of the proof of the theorem goes as in the previous papers \cite{Mono} and \cite{C1}; we shall
not repeat the arguments when they are the same. 

Let $E$ be as in the theorem. First observe that $E \cap B(0,10^{-2}R) \cap L \neq \emptyset$,
because otherwise, $E$ is a plain minimal set in $B(0,10^{-2}R)$ (with no sliding boundary)
that looks a lot like a half plane. If $\varepsilon_0$ is small enough, this is impossible: either
take a limit and find that $X$ is plain minimal in $(0,10^{-2})$, or (alas, by another limiting argument)
say that the density of $E$ at some point near $0$ is $<\pi$.

So choose $x_0 \in E \cap L \cap B(0,10^{-2}R)$. It follows from our assumption 
\eqref{29a.4}, and the usual upper semicontinuity lemma for limits (see Lemma 22.3 in \cite{Sliding},
which we may apply with $M$ arbitrarily close to $1$ and $h$ arbitrarily small, or if you prefer 
Theorem 22.1 in \cite{Sliding}) that if $\varepsilon_0$ is small enough and if $R_1 = 9R/10$,
\begin{equation} \label{29a.7}
R_1^{-2} \H^2(E \cap B(x_0,R_1)) \leq \frac{\pi}{2} + \varepsilon_1,
\end{equation}
with $\varepsilon_1$ as small as we want (provided that $\varepsilon_0$ is small enough).
Notice even that $\varepsilon_0$ does not depend on $X$ or our choice of $x_0$; the point is that if
this failed, we could find a sequence of almost minimizers, with $\varepsilon_0$ tending to $0$, 
with $R = 1$, and so that the points $x_0$ are all translated back to $0$; even that way the corresponding
sets $E$ converge to a set like $X$ and the upper semicontinuity lemma gives the desired contradiction.

By the almost monotonicity of density (and again if $\varepsilon_0$ is small enough), we deduce from
\eqref{29a.7} that
\begin{equation} \label{29a.8}
r^{-2} \H^2(E \cap B(x_0,r)) \leq \frac{\pi}{2} + 2\varepsilon_1
\ \text{ for } 0 < r \leq R_1.
\end{equation}
So we also control the density.
Next we claim that
\begin{equation} \label{29a.9}
L \cap B(x_0, \frac{8R}{10}) \subset E.
\end{equation}
This is also proved in \cite{Mono}, but let us sketch an argument that should convince the reader.
Suppose that for some $r < \frac{8R}{10}$, $E$ does not contain the two points of 
$L \cap \d B(0,r)$. Then the construction of Section \ref{Sfree} gives a competitor $F_0$
of $E$, in $B(0,\frac{9R}{10})$, that looks a lot like $E$ on $B(0,\frac{9R}{10})$ but does not meet
$L \cap B(0,r)$. Rather than using the whole story about free attachments, let us just observe that it
is now rather easy to contract most of $E \cap B(0,r)$ onto a piece of $\d B(0,r)$ that is very close to
$X$. When $n=3$, this is enough to save substantial $\H^2$-measure (because $\d B(0,r)$ has 
a finite measure); in higher co-dimensions, we also need to do an additional Federer-Fleming projection,
as we did in the proof of Lemma \ref{tf.1}. Even in this case, we get a competitor with substantially less
area, and the ensuing contradiction proves \eqref{29a.9}.

Next the proof of \eqref{29a.7} also gives that for $x \in L \cap B(0, \frac{8R}{10})$,
\begin{equation} \label{29a.10}
(10^{-1} R)^{-2} \H^2(E \cap B(x,10^{-1}R)) \leq \frac{\pi}{2} + \varepsilon_1,
\end{equation}
and by near monotonicity as above
\begin{equation} \label{29a.11}
r^{-2} \H^2(E \cap B(x,r)) \leq \frac{\pi}{2} + 2\varepsilon_1
\ \text{ for } 0 < r \leq 10^{-1} R.
\end{equation}
But $x\in E$ by \eqref{29a.9}, and there is no possible density smaller than $\frac{\pi}{2} + 2\varepsilon_1$
other than $\frac{\pi}{2}$, so this means that 
\begin{equation} \label{29a.12}
\lim_{r \to 0} r^{-2} \H^2(E \cap B(x,r))= \frac{\pi}{2}.
\end{equation}
Now we apply Theorem \ref{t22.2n}, with for $X$ the half plane provided by \eqref{29a.4};
the excess density assumption comes from \eqref{29a.11} and \eqref{29a.12}; we get that 
(if $\varepsilon_1$ is small enough) $E$ has a tangent half plane $T(x)$ at $x$, 
\begin{equation} \label{29a.13}
d_{x,r}(E, T(x)) \leq c_1(\varepsilon_1) (10r/R)^{a/4} \ \text{ for } 0 < r < \frac{R}{10}
\end{equation}
(as in \eqref{21.8} and with $c(\varepsilon_1)$ as small as we want), and (as in \eqref{22.9n}
\begin{equation} \label{29a.14}
f(r) \leq c_1(\varepsilon_1) (r/R)^{a/4} \ \text{ for } 0 < r \leq \frac{R}{10}.
\end{equation}
This is where there is a difference with \cite{Mono}, because we get some decay.
Observe also that when we apply \eqref{29a.13} with $r=R/10$ and compare with \eqref{29a.4},
we get that $T(x)$ is close to $X$, i.e.,
\begin{equation} \label{29a14aa}
d_{0,1}(T(x),X) \leq 10 c_1(\varepsilon_1).
\end{equation}
Then we also have \eqref{29a.13} for $R/10 \leq r \leq R/2$, even though with a larger
constant $20 c_1(\varepsilon_1)$ (compare with  \eqref{29a.4} again).

The estimate \eqref{29a.13} is good enough to control the approximation of $E$
by half planes in small balls $B(x,r)$ centered on $L \cap B(0, \frac{8R}{10})$,
but we also want to consider balls $B(y,t)$ for which $y\in E \cap B(0,R/5) \sm L$ and, 
say, $t < R/10$.

So let $y\in E \cap B(0,R/5) \sm L$ and $t < R/10$ be given. 
Set $d = \dist(y,L)$ and pick $x\in L$ such that $|x-y| = d$. 
We start when $d\leq 4t$. Then $T(x)$ is still close to $E$ in $B(y, t)$, because
\begin{equation} \label{29a.14a}
d_{y,t}(E,T(x)) \leq \frac{d+t}{t}\, d_{x,d+t}(E,T(x)) 
\leq 20 c_1(\varepsilon_1) (10(d+t)/R)^{a/4}
\leq 60c_1(\varepsilon_1) (t/R)^{a/4}
\end{equation}
because $B(y, t) \subset B(x, d+t)$, then by our extension of \eqref{29a.13}, applied 
to $d+t \leq 5t \leq R/2$, and because we may take  $a\leq 1$. 
Notice that $T(x)$ coincides with a plane $P(x)$ in $B(x,d(x))$, so 
\eqref{29a.14a} says that $d_{y,t}(E,P(x)) \leq 60c_1(\varepsilon_1) (t/R)^{a/4}$.

When $t \geq d/3$, the simplest is to first use \eqref{29a.14a} with $t=d/2$,
to show that $E$ is very close to a plane in $B(y,d/2)$. Then we can apply the regularity
result for plain almost minimal sets, i.e., the analogue of Theorem \ref{t22.2n} near planes 
and with no sliding boundary. We find that $E$ has a tangent plane $T(y)$ at $y$, and even
\begin{equation} \label{29a.15}
d_{y,t}(E, T(y)) \leq c (t/d)^{a/4},
\end{equation}
with $c$ as small as we want.
This would be enough to prove the more precise estimate 
\begin{equation} \label{29a.15a}
d_{y,t}(E, T(y)) \leq c (t/R)^{b/4}
\end{equation}
for some other $b>0$, by distinguishing the cases in terms of the relative position of 
$t<d<R$; we shall do this with Theorem \ref{t29b.1} below, for instance, but in the present
we can obtain \eqref{29a.15a} in a slightly more direct way, which we present here.

Instead of applying Theorem \ref{t22.2n} to get \eqref{29a.14} for balls centered on 
$E \cap L$, we can use decay estimates for the functional $F$ defined as in \eqref{22.3},
but associated to the center $y \in E \sm L$. The density $F(0)$ is $\pi$
because $E$ has a tangent at $y$; we could also use the approximation of $E$ by a plane in 
$B(y,d/2)$ to show that it cannot be larger than this. On the other hand,
the proof of \eqref{29a.7} also shows that $\H^2(E\cap B(y,2R/3)) \leq \H^2(X\cap B(y,2R/3))
+ \varepsilon_1 R^2$, with $\varepsilon_1$ as small as we want. Then, denoting by $S$
the shade of $L$ lit by $y$,
\begin{eqnarray} \label{29a.15b}
F(2R/3) &=& (2R/3)^{-2} [\H^2(E\cap B(y,2R/3))+\H^2(S\cap B(y,2R/3))]
\nn\\
&\leq& (2R/3)^{-2} [\H^2(X\cap B(y,2R/3))+\H^2(S\cap B(y,2R/3))] + 3\varepsilon_1
= \pi + 3\varepsilon_1,
\end{eqnarray}
because $X$ is a half plane bounded by $L$, thus $\H^2(X\cap B(y,2R/3))$ is largest when
$y\in X$ and then $X\cup S$ is a plane through $y$, whose density is precisely $\pi$.

Then Theorem \ref{t23.1} (applied with the second assumption in \eqref{23.6}, and for instance
with $r_2 = R/2$, gives a good decay estimate \eqref{23.7} for $F(t)$. Then for $0<t<d$,
\begin{equation} \label{29a.15c}
\theta(t) = t^{-2} \H^2(E\cap B(y,t)) = F(t) \leq \pi + (4t/R)^a + CC_h t^a R^{\beta-a}
\leq \pi + (4t/R)^a (1+ C\varepsilon_0) 
\end{equation}
by \eqref{29a.2} and \eqref{29a.3}, and now Theorem \ref{t23.5} gives \eqref{29a.15a}.

Once we have \eqref{29a.13}, it is easy to obtain the $C^1$ description of Theorem \ref{t29a.1},
and we don't even need to use a form of Reifenberg topological theorem. 
For instance, we need to define a function $\varphi$ whose graph contains $E \cap B(0,R/10)$.
Let $\pi$ denote the orthogonal projection on the plane that contains $X$ and set
$\pi^\perp = I - \pi$; we need to know that
\begin{equation} \label{29a.15d}
|\pi^\perp(x)-\pi^\perp(y)| \leq \tau |\pi(x)-\pi(y)| \ \text{ for } x,y\in E \cap B(0,R/10).
\end{equation}
By symmetry, we just need to check this when $\dist(y,L) \leq \dist(x,L)$. If 
$|x-y| \geq \dist(x,L)/10$, we select $x'\in L$ such that $|x-x'| = \dist(x,L)$,
notice that $x'\in E \cap B(0,R/2)$ by \eqref{29a.9}, use \eqref{29a.13}
with $r= 3|x-y|$ to find out that $x$ and $y$ are both much closer to $T(x')$ than they are
to each other. In addition, \eqref{29a14aa} says that $T(x')$ is almost parallel to $X$,
and \eqref{29a.15d} in this case follows.

When $|x-y| \leq \dist(x,L)/10$, we first observe that by \eqref{29a.15a} 
(applied to $t=2|x-y|$) says that $x$ and $y$ are much closer to $T(x)$
as they are to each other; in addition, $T(x)$ is as close to $T(x')$ as we want
because on $B(x,|x'-x|/4)$, \eqref{29a.15a} says that $E$ is quite close to $T(x)$,
while \eqref{29a.13} or \eqref{29a.14a} says that $E$ is close to $T(x')$. So by 
\eqref{29a14aa} $T(x)$ is also as close to $X$ as we want, and \eqref{29a.15d} follows.

So $E \cap B(0,R/10)$ is the graph of a Lipschitz function (defined on $\pi(E \cap B(0,R/10))$).
The estimate \eqref{29a.1}, which can also be seen as an estimate on the H\"older norm of
$D\varphi$, follows easily from \eqref{29a.13} and \eqref{29a.15a}, as in the proof of
\eqref{29a14aa}. Then $\varphi = 0$ on $L \cap B(0,R/10)$ because this set is contained in $E$.
The fact that $\pi((E \cap B(0,R/10))$ lies on one side of $L$ also follows from \eqref{29a.13}
(recall that $T(x)$ is a half plane that lies close to $X$ when $x\in L$). 

Finally, we did not mention in the statement that $\pi(E \cap B(0,R/10))$ contains 
$X \cap B(0,R/11)$, for instance. This would follow from the information that we have
from the local regularity of $E\sm L$ and a fairly simple degree argument.
Theorem \ref{t29a.1} follows.
 \qed

\ms 
We also promised a version of Theorem \ref{t29a.1} with \eqref{29a.4} replaced by \eqref{29a.6}.
We claim that with \eqref{29a.6} (and the other assumptions of the theorem), it is easy to find a
half plane $X$ such that
\begin{equation} \label{29a.16}
d_{0,98R/100}(E,X) \leq \varepsilon_1,
\end{equation}
with $\varepsilon_1$ as small as we want. We consider the functional $F$ associated to a
center $x_0 \in E \sm L$ chosen such that $|x_0-0| \leq 2\varepsilon_0 R$, but otherwise
computed as in \eqref{1.20} or \eqref{22.3}. We know that $F(0) \geq \pi$ because there is no smaller
density at a point of $E$. And 
\begin{eqnarray} \label{29a.17}
F((1-2\varepsilon_0) R) 
&\leq& (1-2\varepsilon_0)^{-2} R^{-2} \H^2(E \cap B(x_0,(1-2\varepsilon_0) R)) + \frac{\pi}{2}
\nn\\
&\leq& (1-2\varepsilon_0)^{-2} R^{-2} \H^2(E \cap B(0,R)) + \frac{\pi}{2}
\nn\\
&\leq& (1-2\varepsilon_0)^{-2} \big[\frac{\pi}{2}+\varepsilon_0\big] + \frac{\pi}{2}
\leq \pi + 10\varepsilon_0
\end{eqnarray}
by \eqref{22.3} and \eqref{29a.6}. So $F$ stays approximately constant on $(0, (1-2\varepsilon_0) R)$
(because $F$ is almost monotone), and Theorem 1.6 in \cite{Mono} says that in $B(x_0,99R/100)$,
$E$ is as close as we want to a minimal sliding set $E_0$ for which $F$ is constant and very close to $\pi$.
Then by Theorem 1.3 in \cite{Mono}, and the discussion of Lemma \ref{t22.1}, $E_0$ coincides with a half plane
in $B(x_0,99R/100)$, and we get \eqref{29a.16}.
Once we have  \eqref{29a.16}, we can end the proof of Theorem~\ref{t29a.1} exactly as if we had 
\eqref{29a.4}.

A last comment is in order before we go to other cases: if $0 \in E \cap L$ is a point of density
$\pi/2$, then the assumptions of Theorem~\ref{t29a.1} are satisfied for $R$ small
(let us say, with the alternate assumption \eqref{29a.6} to save some time); hence we get a 
good description of $E$ near $0$. If all the cases (depending on the blow-up limits of $E$ at
$0$) were as friendly as this one, we would get a nearly perfect description of the singularities 
of $E$ near $L$. We want to continue in this direction, but some cases will not be as friendly.

\section{When $E$ is close to a generic $\bV$ set} 
\label{S29b}

The second case when we have no surprise is when the cone $X$ in \eqref{29a.4} is a generic set
of type $\bV$. That is, $X$ is the union of two half planes $H_1$ and $H_2$ bounded by $L$,
and the angle of $H_1$ and $H_2$ along $L$ is such that
\begin{equation} \label{S29b.1}
\2 < \Angle(H_1,H_2) < \pi.
\end{equation}
Let us add some notation to simplify our description. Denote by $v_0$ a unit vector parallel
to $L$, and for $i=1,2$, let $v_i$ be the unit vector of $H_i$ that is orthogonal to $v_0$.
That is, $v_i$ points directly in the direction of $H_i$. We may also need to use the plane $P_i$
that contains $H_i$, the $(n-2)$-dimensional vector space $H_i^\perp$ orthogonal to
$P_i$, and the orthogonal projection $\pi$ from $\R^n$ to $P_i$. Another way to state \eqref{S29b.1}
is to say that
\begin{equation} \label{S29b.2}
-1 < \langle v_1,v_2\rangle < -\frac12 \, .
\end{equation}

\begin{thm}\label{t29b.1} 
There is a constant $a>0$ that depends only on $n$ and $\beta$ and, for each 
value of $\Angle(H_1,H_2) \in (\2,\pi)$ and $\tau > 0$, 
a constant $\varepsilon_0 > 0$, that depends only on $n$, $\beta$, $\Angle(H_1,H_2)$, and $\tau$,
with the following properties.
Let $E$, $h$, $R$, satisfy \eqref{29a.1}-\eqref{29a.3}, and assume that \eqref{29a.4}
holds for some set $X \in \bV$ such that \eqref{S29b.1} holds. Then $E$ coincides in $B(0,R/10)$
with the union of two graphs of $C^{1}$ functions $\varphi_i : H_i \to H_i^\perp$. In addition,
$\varphi_i(x) =x$ for $x\in L$, $\varphi_i$ is $\tau$-Lipschitz, and for $i=1, 2$,
\begin{equation} \label{29b.3}
|D\varphi_i(x)-D\varphi_i(y)| \leq \tau |x-y|^a R^{-a} \ \text{ for } x,y \in H_i \cap B(0,R/10).
\end{equation}
\end{thm}

\ms
That is, $E$ coincide with the union of two $C^{1+a}$ faces that meet along the common edge $L$,
which is thus locally contained in $E$;
the proof will show that, since the the tangent cone at $x\in L$ is in fact very close to $X$, these
faces make an angle at $x$ which is still generic. But this angle is expected to depend (slowly) on $x$.

We decided to write the H\"older-continuity of the tangent direction of $E$ at $z\in E$,
on each of the two faces of $E$ near the origin, in terms of the derivative of $\varphi_i$; but
this is the same sort of estimate as \eqref{29a.5}.

When the angle of $H_1$ and $H_2$ tends to $\2$ or $\pi$, we make $\varepsilon_0$ tend to $0$,
because we don't want to allow the angle of the two half tangents to $E$ along $L$ to take the limit values.
This is because when this happens the topology of $E$ and $E\cap L$ may change, 
as we will see in later sections.

Finally observe that if $0 \in E \cap L$ and one of the blow-up limits of $E$ at $0$ is a generic
$\bV$-set $X$, then we can apply Theorem \ref{t29b.1} for some small radii $r$, 
and we get a local description of $E$ near $0$. And in particular $X$ is the only blow-up limit
of $E$ at $X$.

\ms
The proof of Theorem \ref{t29b.1} will follow a similar route as when $X$ was a half plane. 
Some parts of the argument will stay valid in the next cases when $X$ is not generic, 
and we shall mention that along the way.
For instance, the following lemma is still valid when $X$ is a plane that contains $L$ 
(but not a sharp $\bV$ set).

\begin{lem}\label{t29b.2}
The density of $E$ at every point of $E \cap B(0,2R/3) \sm L$ is $\pi$.
\end{lem}

Suppose not, and let $x_0 \in E \cap B(0,2R/3) \sm L$ have a density larger than $\pi$. 
Set $d_0 = \dist(x_0,L)$, and first assume that $d_0 \geq 10C \varepsilon_0 R$, where
$C$ is a large constant that will be chosen soon. 
Set $B_0 = B(x_0,d_0/10)$; observe that $B_0 \subset B(0,R)$, so 
$d_{x_0,d_0/10}(E,X) \leq \frac{10R}{d_0} d_{0,R}(E,X) \leq \frac{10R \varepsilon_0}{d_0} \leq C^{-1}$. This means that $X$ meets $B(x_0,d_0/100)$ (because $x_0 \in E$),
and also, since $X$ coincides with a plane $P$ in $2B_0$ (the two faces of $X$ 
make a large angle along $L$, and $10B_0$ does not meet $L$), we get that 
$d_{x_0,d_0/10}(E,P) \leq C^{-1}$.
Now if $C$ is large enough (maybe depending on $n$), and since we have 
\eqref{29a.2} and \eqref{29a.3}, the standard regularity theorem of \cite{Ta} for 
plain almost minimal sets implies that $x_0$ is a point of density $\pi$. 

In fact, we don't even need the result of \cite{Ta} to prove this; it is enough to observe that 
by the usual use of Theorem 22.1 in \cite{Sliding} and a comparison with a plane,
the density $(d_0/20)^{-2} \H^2(E \cap B(x_0,d_0/20))$ is smaller than $\pi+\eta$, 
with $\eta$ as small as we want, hence by the near monotonicity of density, the density of $E$ 
at $x_0$ is at most $\pi + 2\eta$, and now we observe that there is no possible density between 
$\pi$ and $\pi + 2\eta$.
So $x_0$ is a point of density $\pi$, and this contradiction with its definition implies that
\begin{equation} \label{29b.4}
d_0 \leq 10C \varepsilon_0 R.
\end{equation}
We are now ready to apply Theorem \ref{t28.2}, to the set $E-x_0$ so that $x_0$
becomes the origin, and the radius $R' = 10^{-4}R$ so that $B(x_0,400R') \subset B(0,R)$. 
We obtain a cone $Y \in \bY(L,R')$ centered at $x_0$, such that if $Y^t$ is the same cone 
truncated by $L$, then
\begin{equation} \label{29b.5}
d_{x_0,R'}(E,Y^t) \leq C_6 \Big[ [F(200R') -\frac{3\pi}{2}] + \varepsilon_2,
\end{equation}
as in \eqref{28.8}, and with $\varepsilon_2$ as small as we want. We also need to estimate
\begin{equation} \label{29b.6}
F(200R') = (200R')^{-2} \H^2(E \cap B(x_0,200R')) + (200R')^{-2} \H^2(S \cap B(x_0,200R')),
\end{equation}
where $S$ is the shade of $L$ seen from $x_0$. Set $\rho = 200R' + d_0$, and call
$y_0$ the point of $L$ that lies closest to $x_0$; then
\begin{eqnarray} \label{29b.7}
\H^2(E \cap B(x_0,200R')) + \H^2(S \cap B(x_0,200R'))
&\leq& \H^2(E \cap B(y_0,\rho)) + \H^2(S \cap B(y_0,\rho))
\nn\\
&\leq& \H^2(E \cap B(y_0,\rho)) + \frac{\pi \rho^2}{2}.
\end{eqnarray}
We apply the upper semicontinuity estimate again, as for \eqref{29a.7} and find that if $\varepsilon_0$
is small enough, 
\begin{equation} \label{29b.8}
\H^2(E \cap B(y_0,\rho)) \leq \H^2(X \cap \ol B(y_0,\rho))+ \varepsilon_2 \rho^2
\leq \pi \rho^2 + \varepsilon_2 \rho^2;
\end{equation}
again $\varepsilon_0$ does not depend on $y_0$ or $X$ (provided that is is a minimal cone of density $\pi$).
When we combine everything we find out that $F(200R) -\frac{3\pi}{2}$ is as small as we want, 
and so we get that $d_{x_0,R'}(E,Y^t) \leq \varepsilon_3$
with $\varepsilon_3$ as small as we want.
But $R' = 10^{-4}R$, and in the larger ball $B(0,R)$ our set $E$ is very close to $X$, which is a plane or a generic set of type $\bV$; this yields the desired contradiction if $\varepsilon_3$ and $\varepsilon_0$ 
are small enough. Not surprisingly, we need to take $\varepsilon_0$ even smaller when $X$ is close to
being sharp.
\qed

\begin{lem}\label{t29b.3}
The density of $E$ at every point of $E\cap L \cap B(0,R/2)$ is $\pi$.
\end{lem}

\ms
This lemma is still true when $X$ is a plane but we expect it to fail when $X$ is a sharp 
$\bV$ set (and also its proof uses Lemma \ref{t29b.2}).

Let $z \in E\cap L \cap B(0,R/2)$ be given. By the same upper semicontinuity argument 
as for \eqref{29a.7} and \eqref{29b.8}, we find that if $\varepsilon_0$ is small enough,
\begin{equation} \label{29b.9}
\H^2(E \cap B(z,R/3)) \leq \H^2(X \cap \ol B(z,R/3))+ \varepsilon_4 R^2
\leq (R/3)^2 (\pi+ 9\varepsilon_4),
\end{equation}
with an $\varepsilon_4$ which is as small as we want. Then by the near monotonicity of the density,
\begin{equation} \label{29b.10}
r^{-2}\H^2(E \cap B(z,r)) \leq \pi + 10\varepsilon_2 \ \text{ for } 0 < r < R/2.
\end{equation}
So the density of $E$ at $z$ is at most $\pi + 10\varepsilon_2$.
Recall from Lemma \ref{t22.2} that (if $\varepsilon_2$ is small enough) all the cones with such
a density are half planes, planes, and sets of type $\bV$. Hence, the density of $E$ at $z$ is either
$\pi/2$ or $\pi$, and we just need to exclude the case when it is $\pi/2$. 

So we assume that the density is $\pi/2$; then some blow-up limit of $E$ at $z$
is a half plane $H$, and Theorem \ref{t29a.1} says that there is a small ball centered at $z$
where $E$ is $C^1$-equivalent to $H$.

The topological argument that follows is almost the same as in Section 17 of \cite{Holder}, 
starting at (17.9), so we will only give an outline, and send to \cite{Holder} for details. 
Notice however that our life is a little simpler here, because we are willing to use a 
$C^1$ description of $E$ near every point of $E \cap B(0,R/2) \sm L$, 
whereas in \cite{Holder} we wanted to merely use a H\" older description.

Let us first assume that $n = 3$, because the topology is simpler then, and consider 
a small circle $C_0$ of radius $\rho$ centered on $z$ and contained in the plane orthogonal 
to $L$ at $z$. 
If $\rho$ is small enough, $C_0$ meets $E$ exactly once, and transversally. 
On the other hand, let $C_1$ be the circle centered at $0$, contained in the plane 
orthonormal to $L$ and with radius $R/3$. Recall that \eqref{29a.4} says that $E$ is very close to 
$X$ in $B(0,R)$; denote by $x_1$ and $x_2$ the two points of $C_1 \cap X$, and observe that
$X$ coincides with a plane in both $B(x_i,R/10)$. Then we can apply the local regularity theorem
in $B(x_i,R/20)$ and find that in each of these balls, $E$ meets $C_1$ transversally exactly once.
Of course $E$ does not meet the rest of $C_1$ (too far!), so $C_1$ has just two transverse 
intersections with $E$. 

Now there is a homotopy $\{ h_t \}$, $0 \leq t \leq 1$, that goes from $C_0$ to $C_1$, 
and whose image lies in $B(0,2R/3) \sm L$. But $E$ is locally $C^1$ there, 
because Lemma \ref{t29b.2} says that each point of $E \cap B(0,2R/3) \sm L$ is a point
of density $\pi$, and we know that $E$ is a $C^1$ surface near such points.
We claim that along the homotopy, the number of intersections of $E$ with $C_t = h_t(C_0)$ 
stays the same modulo $2$, which leads to the desired contradiction. 
The proof of the claim would consist in transforming the homotopy slightly, so that for each $t$ the loop
$h_t(C_0)$ meets $E$ a finite number of times, and each time transversally. But let us describe 
the general case first, and anyway refer to \cite{Holder} for details.

When $n \geq 3$, we want to proceed as above, but replace the circles $C_j$ with
$(n-2)$-spheres. That is, $C_0$ is now a small $(n-2)$-sphere centered at $z$ and contained
in the hyperplane orthogonal to $L$ at $z$, and $C_1$ is the $(n-2)$-sphere in the hyperplane 
orthogonal to $L$ at $0$, with radius $R/3$. As before, we can find a homotopy from $C_0$ to $C_1$,
among spheres that are contained in $B(0,2R/3)\sm L$, and prove that (after a suitable modification
to put things in general position) the number of intersections with $E$ stays the same modulo $2$.

Now we follow a suggestion of Christopher Collins: rather than trying to use degree theory too
soon, we simply write the two equations of $C_t$, and thus get a one parameter family of
functions $h_t$, defined on $\R^n$ and with values in $\R^2$. Thus 
$C_t = \big\{ x\in \R^n \, ; \, h_t(x) = 0\big\}$, and we want to compute the number
of solutions in $E$ of the system of two equations $h_t(x)=0$.

Then we discretize, and at the same time modify the functions $h_t$ to put $C_t$ in the
general position with respect to $E$, so that it intersects $E$ transversally. Also, we cut the
elementary move from $h_{t_i}$ to $h_{t_{i+1}}$ into smaller modifications, where each time 
the function $h_t$ is only modified in a small ball where we have a good description of $E$
as a $C^1$ surface. This gives a new collection of mappings, that we shall call
$\wt h_j$, such that $\wt h_j^{-1}(0)$ is always transverse to $E$ and meets it a finite
number of times, and $\wt h_{j+1}-\wt h_j$ is small (in supremum norm) and supported in a small ball.

Then we have to show that $E \cap \wt h_{j+1}^{-1}(0)$ and $E \cap \wt h_j^{-1}(0)$ have the same
number of points modulo $2$, and for this we finally use some degree theory. We replace the equation
$\wt h_j(y)=0$ (in $E$) with the equivalent equation $\wh h_j(y)=0$, where $\wh h_j$ is just
a version of $\wt h_j$ which takes values in the $2$-sphere. Because of our transversality
condition and modulo $2$, the number of solutions is the degree of $\wh h_j$, and this degree
is the same for $\wh h_{j+1}$ and $\wh h_j$, because $\wt h_{j+1}-\wt h_j$ is small.

This completes our rapid proof of the fact that $\theta(x) = \pi/2$ never occurs, and
 Lemma~\ref{t29b.3} follows.
\qed

\ms
\begin{lem}\label{t29b.4}
The set $E$ contains $L \cap B(0,2R/3)$.
\end{lem}

This lemma stays valid when $X$ is a sharp $\bV$-set. 
That is, we just need to exclude flat $\bV$-sets here, i.e, planes that contain $L$. 

First we check that $E$ meets $L \cap B(0,10^{-2}R)$. Otherwise, $E$ is also a plain
almost minimal set in $B(0,10^{-2}R)$, with the same gauge function $h$, but no sliding boundary.
Yet \eqref{29a.4} says that $X$ is $\varepsilon_0 R$-close in $E\cap B(0,R)$ to the non-flat set $X$ 
of type $\bV$. If this can happen for arbitrarily small $\varepsilon_0$, a small limiting argument 
shows that $X$ (or another non-flat $\bV$ set) is minimal in $B(0,10^{-2}R)$. 
This is false, so $E$ meets $L \cap B(0,10^{-2}R)$.

Pick $x_0 \in E \cap L \cap B(0,10^{-2}R)$; if \eqref{29b.11} fails, then there is a radius 
$t \in (0,3R/4)$ such that at least one of the two points of $L \cap \d B(x_0,t)$ lies outside of $E$.
So we suppose so and get a contradiction.

Let us first describe a rather brutal argument that uses Lemma \ref{tf.1}.
This will also be an opportunity to describe what we meant there; after this we will sketch a more direct argument with the same ideas.

Let us run the argument of Sections \ref{S16}-\ref{S21}, with the set $E-x_0$
(because we want to use $x_0$ as an origin), the same approximating
cone $X$ that we have, and radii $r\in (\frac{R}{10},\frac{R}{2})$.
The various smallness assumptions that we need to do are satisfied. In addition,
because of our assumption on the existence of $t$, Lemma \ref{tf.1} allows us to use the estimates of 
Sections~\ref{S16}-\ref{S21} with the free attachment. 

Recall that in this case we build four Lipschitz curves $\Gamma_{\pm,i}$, where the notation
makes sense because $\Gamma_{\pm,i}$ goes from a point $m_i \in E$ that lies close to
the midpoint $w_i$ of $H_i \cap \S_r$, to a point $z_\pm$ that lies close to a point 
$\ell_\pm \in L \cap \S_r$. We also construct a net $\rho^\ast = \rho_r^\ast$ of geodesics,
composed of the four geodesics $\rho_{\pm, i} = \rho(w_i, z_\pm)$.
Now these endpoints lie on $E$, very close to points of $X=H_1 \cup H_2$, and
if $\varepsilon_0$ is small enough we get that the angle at $z_\pm$ of $\rho_{\pm,1}$ 
and $\rho_{\pm,2}$ is at most $\pi - \alpha/2$, where we set $\alpha = \pi - \Angle(H_1,H_2) > 0$.

The simplest for us now is to catch the argument during the construction of our third competitor, 
where we modify the tip of the second one, which coincides near the origin with the cone 
$X^\ast$ over $\rho^\ast$.
We know that we can modify this tip, and we do so by softening the angles near the $z_\pm$.
Recall that we don't need to worry about what happens along $L$; this is the advantage of having 
the free attachment at both points of $\S_r \cap L$.
It turns out that we did this sort of modification and the associated computation in Section \ref{S26},
for \eqref{26.2}. It does not matter that there the point where we soften the angle was $m_i$
rather than $z_\pm$, the computation is the same and we can save an area comparable to 
$C^{-1} \alpha^2 r^2$. The estimate that we get, instead of \eqref{15.4} is now
\begin{eqnarray}\label{29b.11}
\H^2(E \cap \ol B(0,r)) &\leq&  \frac{r}{2} \, H^1(E \cap \S_r)
- 10^{-5} [\H^1(E \cap \S_r) - \H^1(\rho^\ast_r)]  - C^{-1} \alpha^2 r^2 + R^2 h(R)
\nn\\
&\leq& \frac{r}{2} \, H^1(E \cap \S_r) - C^{-1} \alpha^2 r^2 + R^2 h(R),
\end{eqnarray}
where $R^2 h(R)$ comes from the almost monotonicity property (and the fact that we use a competitor
in $B(0,R)$), and we can drop the middle term with $10^{-5}$, which is nonnegative by \eqref{14.45},
\eqref{14.33}, and the comment below \eqref{14.30}. 

Now Lemma \ref{t16.1} says that $v(r) = \H^2(E \cap \ol B(0,r))$ is differentiable 
almost-everywhere, that its distribution derivative is at least as large as its almost everywhere 
derivative $v'(r)$, and that the same thing holds for 
$\theta(r) = r^{-2} \H^2(E \cap B(0,r)) = r^{-2} v(r)$, with 
\begin{equation}\label{29b.12}
\theta'(r) = - 2r^{-3} v(r) + r^{-2} v'(r) \geq - 2r^{-3} v(r) + r^{-2} \H^1(E \cap \S_r),
\end{equation}
where the second inequality comes from \eqref{16.13} and is valid for almost every $r \in (0,1)$. 
Notice that $\H^2(E \cap \S_r) = 0$ for almost every $r$ too, so \eqref{29b.11} and \eqref{29b.12}
yield
\begin{eqnarray}\label{29b.13}
\theta'(r) &\geq& - 2r^{-3} v(r) + r^{-2} \H^1(E \cap \S_r) 
= - 2 r^{-3} \H^2(E \cap B(0,r)) + r^{-2} \H^1(E \cap \S_r)
\nn\\
&\geq&  - r^{-2} \, H^1(E \cap \S_r) + 2C^{-1} r^{-1}\alpha^2 - 2 r^{-3} R^2 h(R) 
+ r^{-2} \H^1(E \cap \S_r) 
\nn\\
&=& 2 C^{-1}r^{-1} \alpha^2 - 2 r^{-3}R^2 h(R) \geq (3C)^{-1} r^{-1} \alpha^2 
\geq (30C)^{-1} R^{-1} \alpha^2
\end{eqnarray}
when $\frac{R}{10} < r < \frac{R}{2}$ and if $\varepsilon_0$ (and hence $h(R)$) is small enough. 
As we said above, Lemma~\ref{t16.1} allows us to integrate this, and we get that
$\theta(R/2) \geq \theta(R/10 ) + (100C)^{-1} \alpha^2$. 

This contradicts either the fact that $\theta(0) = \pi$ by Lemma \ref{t29b.3}, or that
$\theta(R/2) \leq \pi + 10\varepsilon_4$ by the proof of \eqref{29b.10}, 
or the near monotonicity of $\theta$. This completes the first proof of Lemma \ref{t29b.4}.

For a more natural and direct proof, start by observing that 
\begin{equation} \label{29b.14}
\H^2(E \cap B(x_0, 2R/3)) \geq (2R/3)^2 (\pi - \varepsilon_5),
\end{equation}
with $\varepsilon_5$ as small as we want; this time, the simplest is to use the lower semicontinuity
of measure along minimizing sequences (Theorem 10.97 in \cite{Sliding}), plus the usual limiting 
argument with a sequence of counterexamples with $\varepsilon_0$ tending to $0$, and that tend 
to a $\bV$-set with density $\pi$.

Then, since we still assume that $L \cap \d B(x_0,r)$ is not contained in $L$, we can 
use the proof of Lemma \ref{tf.1} to construct a sliding competitor $F_0$ for $E$ in $B(0,R)$, 
which does not meet $L \cap B(x_0, 3R/4)$,
and which is as close to $E$ in measure as we want. That is, in particular, 
$\H^2(F_0 \sm E) \leq \eta$, with $\eta$ as small as we want, as in \eqref{f.17}. 
And now, since the shape of $E$ in $B(x_0, 2R/3)$ is well known ($E$ is close to the $\bV$-set $X$), 
we can construct by hand a competitor for $F_0$ and $E$ in $B(0,R)$, 
essentially without changing anything in $B(x_0, 2R/3)$ (except for the part of $F_0\sm E$
that is already there, in the thin tube near $L$), and that does better than the quantity that 
we get from \eqref{29b.14}.
As usual, the construction is rather easy when $n=3$ because we can do the gluing in a very small
portion of $\d B(x_0, 2R/3)$ near $X$, and when $n>3$ we would need to do a piece of Federer-Fleming
projection as in the proof of Lemma \ref{tf.1}. We leave the details to the reader, as a 
punishment for not trusting the brutal but complicated proof above. But either way Lemma \ref{t29b.4}
is proved.
\qed

\ms
We are now ready to use decay estimates for the density and the distance to a $\bV$-set.
Let $x_0 \in L \cap B(0,R/2)$ be given. We know that $x_0 \in E$ (by Lemma \ref{t29b.4})
and its density is $\pi$ (by Lemma \ref{t29b.3}). We may thus apply Theorem \ref{t22.2n} to $E-x_0$. 
The cone $X$ has full length (see Section \ref{S30}), 
we just checked the density condition \eqref{21.6}, and \eqref{21.7} with $r_1 = R/2$
follows from \eqref{29a.2}-\eqref{29a.4}. We get the existence of a tangent cone 
$X(x_0)$ to $E$ at $0$, with 
\begin{equation} \label{29b.15}
d_{x_0,r}(E,X(x_0)) \leq c_1(\varepsilon_0) \Big(\frac{2r}{R}\Big)^{a/4}
\ \text{ for } 0 < r < R/2,
\end{equation}
where $c_1(\varepsilon_0)$ is as small as we want.

This cone has density $\pi$ (just like $x_0$), and by \eqref{29b.15} and \eqref{29a.4}
$X(x_0)$ is as close to $X$ as we want. We go to the list of Lemma \ref{t22.2}
and find that $X(x_0)$ is a generic set of type $\bV$.

In addition to \eqref{29b.15}, Theorem \ref{t22.2n} also says that
\begin{equation} \label{29b.16}
f(r) := r^{-2} \H^2(E\cap B(x_0,r)) - \pi 
\leq c_2(\varepsilon_0) (2r/R)^a
\ \text{ for } 0 < r < R/4,
\end{equation}
with $c_2(\varepsilon_0)$ as small as we want.

With this we get a good control of $E$ on all balls centered on $L \cap B(0,R/2)$
and radii $r < R/2$. Then we proceed as for the end of Theorem \ref{t29a.1}.
That is, we also need to control $E$ in other small balls $B(x,t)$; the main case is 
when $t < \dist(x,L)/10$, and for those we start from \eqref{29b.15} 
(applied to $r = 2 \dist(x,L)$ and the point $x_0 \in L$ that lies closest to $x$), 
which that $E$ is very close to a plane in $B(x, \dist(x,L)/2$. Then we can apply the usual regularity
result for plain almost minimizers to prove that $E$ is $C^1$ in $B(x,t)$.

As for the precise control \eqref{29b.3}, which we of course prove with a smaller constant $a$, 
the same argument with the different functional $F$ no longer works, because $x$ has density $\pi$
and $F(R)$ is more like $\frac{3\pi}{2}$. So let us cheat instead. Set $d = \dist(x,L)$
and call $x_0$ the point of $L$ that lies closest to $x$.
If $\frac{d}{t} \leq \big(\frac{R}{d}\big)^{a/8}$, then \eqref{29b.15} yields
\begin{eqnarray} \label{29b.17}
d_{x,t}(E,X(x_0)) &\leq& \frac{d}{t} \, d_{x,d}(E,X(x_0)) 
\leq 2 \Big(\frac{R}{d}\Big)^{a/8} d_{x_0,2d}(E,X(x_0))
\leq 2 \Big(\frac{R}{d}\Big)^{a/8} c_1(\varepsilon_0) \Big(\frac{4d}{R}\Big)^{a/4}
\nn\\
&\leq& 4  c_1(\varepsilon_0) \Big(\frac{d}{R}\Big)^{a/8}
= 4  c_1(\varepsilon_0) \Big(\frac{d}{R}\Big)^{a/16} \Big(\frac{d}{R}\Big)^{a/16}
\leq 4  c_1(\varepsilon_0) \Big(\frac{t}{R}\Big)^{a/16},
\end{eqnarray}
where for the last line we used the fact that $\frac{d}{R} \leq \frac{t}{d}$ because 
$\frac{d}{t} \leq \big(\frac{R}{d}\big)^{a/8} \leq \frac{R}{d}$ (since $a < 1$).
If $P$ is the plane that coincides with $X(x_0)$ near $B(x,t)$, we also get that
\begin{equation} \label{29b.18}
d_{x,t}(E,P) \leq 4  c_1(\varepsilon_0) \Big(\frac{t}{R}\Big)^{a/16}.
\end{equation}
If instead $\frac{d}{t} \geq \big(\frac{R}{d}\big)^{a/8}$, we just use the fact that
$d_{x_0,2d}(E,X)$ is as small as we want to start anew from a good flat approximation 
in $B(x,d/2)$. The analogue of Theorem \ref{t22.2n} for plain almost minimal sets gives a
plane $P$ such that
\begin{equation} \label{29b.19}
d_{x,t}(E,P) \leq c (t/d)^{b}  
\end{equation}
for some $b > 0$, and a constant $c > 0$ that we can take as small as we want.
But $(t/d)^{b} \leq (t/d)^{b/2}(d/R)^{ab/16} \leq (t/R)^{ab/16}$, so we have an analogue 
of \eqref{29b.18} with a different power.

Once we have \eqref{29b.18} and its analogue, we get a good control on the tangent planes and their
variations. We still need to write $E \cap B(0,R/10)$ as a union of two Lipschitz graphs on
(pieces of) the two half planes $H_i$ that compose $X$. Set 
\begin{equation} \label{29b.20}
F_i = \big\{ \, x\in E \cap B(0,R/10) ; \, \dist(x,H_i) \leq 10^{-2} \dist(x,L) \big\}
\end{equation}
for $i=1,2$. It is clear that $F_1 \cap F_2 \subset L$. Also, we observed below \eqref{29b.15} 
that the $\bV$-sets $T(x_0)$ are as close to $X$ as we want; then if $x\in E \cap B(0,R/10)$,
\eqref{29b.15} (applied to the projection $x_0$ of $x$ on $L$ and $r = 2\dist(x,L)$) shows that
$x\in F_1$ or $F_2$, depending on which piece of $T(x_0)$ lies closer to $x$. 
Now we want to show that $F_i$ is a nice Lipschitz graph on $H_i$.
For this we can follow quietly the final argument given for Theorem \ref{t29a.1}, applied
to each $F_i$ separately. Theorem \ref{t29b.1} follows.
\qed

\section{When $E$ is close to a plane that contains $L$}
\label{S29c}

We now want a variant of Theorems \ref{t29a.1} and \ref{t29b.1} for the case when the approximating 
cone $X$ is a plane that contains $L$. We waited this long because this is the first time 
where we may have a slightly complicated singular set. 

\begin{figure}[!h]  
\centering
\includegraphics[width=12.cm]{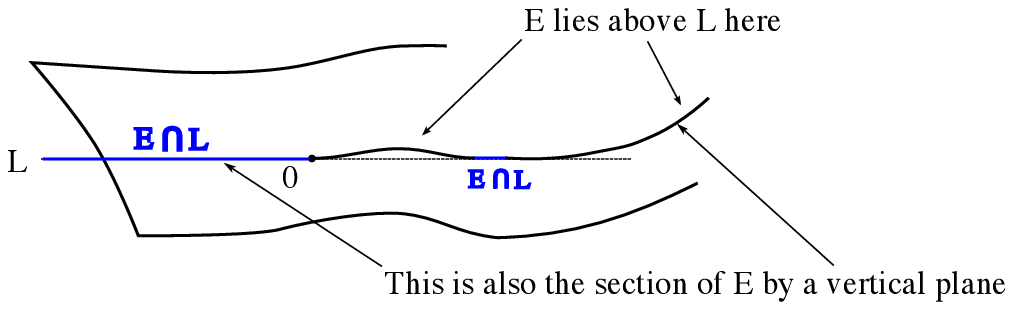}
\caption{ Behavior of $E$ near a plane through $L$; creases may exist along pieces of $E \cap L$\label{f29c1}}
\end{figure}

Let us assume that $E$ satisfies \eqref{29a.1}-\eqref{29a.4}, and that $X$ is a plane that contains
$L$. We intend to prove that in $B(0,R/10)$, $E$ looks like a nice $C^{1+a}$ surface which is also
a small Lipschitz graph over $X$, except that along some part of $E \cap L$, $E$ may have a crease
where all the points admit a tangent cone which is a generic (and in fact almost flat) set
of type $\bV$. See Figure \ref{f29c1}. %

% mis page 278

Let us give a little more notation to prepare the statement. If $\varphi : X \to X^\perp$
is a function, the graph of $\varphi$ is the set 
\begin{equation} \label{29c.1}
{\rm Graph}(\varphi) = \big\{ x+\varphi(x) \, ; \, x\in X \big\}.
\end{equation}
Our crease set will be (the intersection of $B(0,R/10)$ with) 
\begin{equation}\label{29c.2}
L_g = \big\{ y\in E \cap L \, ; \, E \text{ has a blow-up limit at $y$ which is a generic $\bV$ cone}\big\}.
\end{equation}
We know from Theorem \ref{t29b.1} that $L_g$ is an open subset of $L$, and we even have a nice
description of $E$ near each point of $L_g$.

\begin{thm}\label{t29c.1} 
There is a constant $a>0$ that depends only on $n$ and $\beta$ and, for each $\tau > 0$, 
a constant $\varepsilon_0 > 0$, that depends only on $n$, $\beta$, and $\tau$,
with the following properties.
Let $E$, $h$, $R$, satisfy \eqref{29a.1}-\eqref{29a.3}, and assume that \eqref{29a.4}
holds for some plane $X$ that contains $L$. Then there is a $\tau$-Lipschitz function 
$\varphi : X \to X^\perp$ such that 
\begin{equation} \label{29c.3}
E \cap B(0,R/10) = {\rm Graph}(\varphi) \cap B(0,R/10).
\end{equation}
In addition, $\varphi(x) = x$ for $x\in L_g$, $\varphi$ is $C^1$ on $X \sm L_g\cap B(0,R/10)$, 
and
\begin{equation} \label{29c.4}
|D\varphi(x)-D\varphi(y)| \leq \tau |x-y|^a R^{-a} \ \text{ for } x,y \in B(0,R/10) 
\text{ such that }  [x,y] \subset X \sm L_g.
\end{equation}
\end{thm}

\ms
We tried to make the statement short, but there are a few implicit things there that we prefer to explain now.
We said that $\varphi$ is defined on the whole $X$, but of course only the values on $X\cap B(0,R/10)$
matter. This is not too costly, because we can always extend. At points of $L_g$, we know that $E$ has 
a tangent cone $T(x) \in \bV$, which is nearly flat, and we can compute the two half tangents as the limits,
when $y \in X \sm L_g$, of (half of) the tangent plane $T(y)$ given by the derivative $D\varphi(y)$.
We get two different directions, because we have two ways to approach $x$.
Yet $T(x)$ varies in a H\" older-continuous way along $L_g$, because \eqref{29c.4} controls 
the variations of the two half planes. That is,
\begin{equation} \label{29c.5}
d_{0,1}(T(x),T(y)) \leq 2\tau |x-y|^a R^{-a} \ \text{ for } x,y \in L_g \cap B(0,R/10).
\end{equation}
When $x \in B(0,R/10)$ lies in the boundary in $L$ of $L_g$, then $E$ has a tangent plane
at $x$; this comes from \eqref{29c.4} as well (and so it made sense to include $x$
in the domain $X \sm L_g$).

The author expects that near $0$, $E \cap L$ can be almost any closed set that contains $0$.
It is less clear to him whether $L_g$ can be very complicated or not, and whether such pathologies
are also possible for locally minimal sets. 
See Remark \ref{r29b.1} below for a short discussion about this.

As usual, Theorem \ref{t29c.1} can be applied to give a nice description of $E$ near
any point $x\in E$ where some a bow-up limit of $E$ at $x$ is a plane that contains $L$.

\ms
Let us now prove Theorem \ref{t29c.1}. Let $y\in E \cap L \cap B(0,R/2)$ be given.
Lemma \ref{t29b.3}, which is still valid here when $X$ is a plane, says that $y$ is a 
point of density $\pi$. This allows us to apply Theorem \ref{t22.2n} with the cone $X$
(certainly a full length cone), after a translation by $-y$, and to the radius $r_1 = R/2$.
We find that $E$ has a tangent cone $T(y)$ at $y$, and 
\begin{equation} \label{29c.6}
\ d_{y,r}(E,T(y)) \leq c_1(\varepsilon_0) \Big(\frac{2r}{R}\Big)^{a/4}
\ \text{ for } 0 < r < R/2,
\end{equation}
with $c_1(\varepsilon_0)$ as small as we want, as in \eqref{21.8}. 
We know that $T(y)$ has density $\pi$ as well, and (by \eqref{29c.6}
and \eqref{29a.4}) that it is fairly close to $X$. That is, 
\begin{equation}\label{29c.7}
d_{0,1}(T(y),X) \leq 2 c_1(\varepsilon_0).
\end{equation}
But $T(y)$ could be a plane (that contains $L$ or not) or a flat generic $\bV$ set. 
Moreover, a comparison between the various estimates \eqref{29c.6}, 
with $y,z\in E \cap L \cap B(0,R/2)$, yields
\begin{equation}\label{29c.8}
d_{0,1}(T(y),T(z)) \leq 100 c_1(\varepsilon_0) \Big(\frac{|y-z|}{R}\Big)^{a/4}
\ \text{ for } y, z \in E \cap L \cap B(0,R/2).
\end{equation}
This, and in particular \eqref{29c.6}, gives a good description of $E$ in all the balls $B(y,r)$, with 
$y \in E \cap L \cap B(0,R/2)$ and  $0 < r < R/2$.

Notice that all the points $y\in E \cap B(0,R)$ are points of density $\pi$, either
by Lemma \ref{t29b.3} when $y\in L$, or by  Lemma \ref{t29b.2}, which is also valid
when $X$ is a plane, when $x\in E \sm L$. So, except for the points of $L_g$ where
$E$ has a tangent $\bV$-set, $E$ has a tangent plane $T(x)$ at $x$ (recall that a sharp
$\bV$-set would not satisfy \eqref{29c.7}).

At this point we know that $E$ is $C^1$ everywhere on $B(0,R/2) \sm L_g$,
and Theorem \ref{t29b.1} gives a nice description of $E$ near the points of 
$E \cap B(0,R/2) \cap L_g$.  Yet we want more precise and uniform estimates on 
the variations of the direction of $T(x)$, $x\in E \cap B(0,R/2)$, or (essentially equivalently) 
of the numbers $d_{x,r}(E,T(x))$ that control the good approximation of $E$.

For $x\in E \cap L \cap B(0,R/2)$ and $r < R/2$, $d_{x,r}(E,T(x))$
is directly controlled by \eqref{29c.6}. For $x\in E \cap B(0,R/2)\sm L$
and $10^{-1}R \leq r \leq R/2$, we can use \eqref{29a.4} to show that
$d_{x,r}(E,X) \leq 10\varepsilon_0$.
Let us now assume that $x\in E \cap B(0,R/3)\sm L$ and $r < 10^{-1}R$.
Set $d = \dist(x,E\cap L) = \dist(x,E\cap L\cap B(0,R/2))$ and pick 
$y\in E\cap B(0,R/2)\sm L$ such that $|y-x|=d$. If $r \geq d/5$, 
we get a good approximation by the nearly flat $\bV$-set $T(y)$,
since by \eqref{29c.6}
\begin{equation}\label{29c.9}
d_{x,r}(E,T(y)) \leq \frac{r+d}{d} \, d_{y,r+d}(E,T(y)) 
\leq 6c_1(\varepsilon_0) \Big(\frac{12r}{R}\Big)^{a/4}.
\end{equation}
Finally for $r<d/4$, we first observe that $T(y)$ coincides with a plane
$P(y,x)$ near $B(x,2d/3)$ (because $T(y)$ is a rather flat $\bV$-set),
then $d_{x,d/2}(E,P(y,x))$ is as small as we want (by \eqref{29c.6} or \eqref{29c.9}),
then we can apply the analogue of Theorem \ref{t22.2n}
for approximations by a plane $X$ in the plain case (no sliding boundaries); we
get that for $0 < r < d/4$,
\begin{equation} \label{29c.10}
d_{x,r}(E,P(x)) \leq c \Big(\frac{r}{d}\Big)^{a/4},
\end{equation}
where $P(x)$ denotes the tangent plane to $E$ at $x$ (we already knew its existence)
and $c>0$ is as small as we want (provided, as usual, that we take $\varepsilon_0$ small enough). 
The constant $a>0$ may be different (it depends on the full length constants for a plane).
And also, we may replace \eqref{29c.10} by the apparently better
\begin{equation} \label{29c.11}
d_{x,r}(E,P(x)) \leq c' \Big(\frac{r}{R}\Big)^{b/4}
\end{equation}
by the same small trick where we distinguish between cases depending on $r/d$
and use \eqref{29c.6} or \eqref{29c.9} when $\frac{d}{r} << \frac{R}{d}$ as for the end of the
proof of Theorem \ref{t29b.1} (near \eqref{29b.17}).

We may now compare \eqref{29c.10} and \eqref{29c.9} (use $r=d/2$);
we get that $P(x)$ is indeed quite close to $P(y,x)$. Because of this, and also \eqref{29c.7},
all the directions of the planes $P(x)$, $x\in E \cap B(0,R/3)\sm L$ and the half planes
that compose the $T(y)$, $y\in E\cap L\cap B(0,R/2)$, are as close to the direction of $X$
as we want. Then $E \cap B(0,R/3)$ is the graph of a $\tau$-Lipschitz function $\varphi$ 
defined on a subset of the plane $X$, by \eqref{29c.6}, \eqref{29c.9}, \eqref{29c.10}, 
and the proof of \eqref{29a.15d}. Then,  the estimate \eqref{29c.4} on the derivative of
$\varphi$ (or equivalently the direction of $P(x)$ or $T(y)$) follow from these same estimate
(compare the $P(x)$ or $T(y)$ to $E$ on intersecting balls). The fact that $\varphi(y) = y$
for $y\in L \cap E \cap B(0,R/10)$ (and in particular on $L_g \cap B(0,R/10)$) comes from the
graph description. As usual, we would deduce additional information, such as the fact that 
$\pi(E \cap B(0,R/10))$ contains $X \cap B(0,R/11)$, with a little bit of topology.
This completes our proof of Theorem \ref{t29c.1}.
\qed

\ms\begin{rem} \label{r29b.1} 
In the description above, $L \cap E \cap B(0,R/10)$ can probably be just any closed subset of 
$L\cap B(0,R/10)$ that contains the origin, and we could even take $L_g = \emptyset$ to see this. 
That is, we can probably make $E$ leave $L$ and return to $L$ as we wish, 
provided that we keep it extremely close to $L$.
We will not prove this here, but the following argument may convince the reader. 
It is easy to see that a plane $P_0$ is minimal, no matter which choice of boundary constraint
we take (because it is already minimal without constraint); in particular we could take a boundary
$L$ that is a smooth curve, that comes and leaves $P_0$ in a very tangential way, but 
along any given closed subset of $L$. The part of the argument
that we will not do is to show that we can go from this situation to the situation where $L$ is
a straight line, with a change of variable in $\R^n$ that maps $L$ to a line, and which is so close to
the identity that the image of $P_0$ is still almost minimal.
See Section \ref{S31} for a discussion of changes of variables though.

Similarly, the author tends to expect that the set $L_g$ can be, locally, just about any 
open subset of $L$ that does not contain $0$, but constructing examples may be much 
more delicate than for the assertion above, because there will be a balance between 
pulling $E$ enough in one direction (so as to create a crease),
but not too much (to control the almost minimality and also not make the creases too long).

Of course the situation of actual minimal sets could be very different, 
because minimality probably forces some rigidity conditions that the author does not understand
but that may prevent many pathologies, in the same way as analytic functions do not always 
do what we want.
\end{rem}

\section{When $E$ is close to a sharp $\bV$-set}
\label{S29d}
 
 In this section we give a local description of $E$ near $0$, under the usual assumptions
 of Theorems~\ref{t29a.1}, \ref{t29b.1}, and \ref{t29c.1}, except that this time the 
 approximating cone $X$ is a sharp $\bV$-set, i.e., such the two
 half planes $H_1$ and $H_2$ that compose $X$ make $\2$ angles along $L$.
 This case looks more interesting than the previous one, in particular because 
 the topology of $E$ near $0$ will possibly be different from the topology of $X$.
 
\begin{figure}[!h]  
\centering
\includegraphics[width=12.cm]{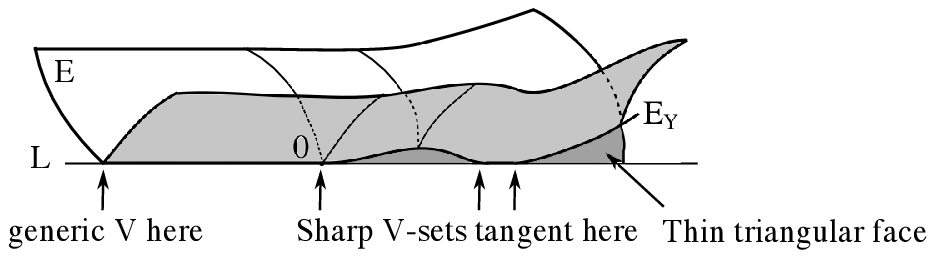}
\caption{ The set $E$ near a sharp $\bV$ set; $\gamma$ is composed of $E_Y$ and a bit of $L$
\label{f29d1}}
\end{figure}

 Another new thing in this case is the possible presence near $0$ of a curve, 
 contained in the set 
\begin{equation} \label{29d.1}
E_Y = \big\{ x\in E \sm L \, ; \, \theta_x(0) = \frac{3\pi}{2} \big\},
\end{equation}
where $\theta_x(0) = \lim_{r \to 0} r^{-2}\H^2(E\cap B(x,r))$ is the density of $E$
at $x$, and along which the $\bV$-part of $E$ detaches itself from $L$, leaving a thin
triangular extra face of $E$ between the curve and a corresponding piece of $L$.
We try to give a more precise description below; in the mean time see 
Figure \ref{f29d1} %
for a good idea of what $E$ may look like.

% mis page 282

But first recall that since the $\bY$-cones are the only plain minimal cones in $\R^n$ with a density 
$\frac{3\pi}{2}$, $E_Y$ is also the set of points of $E \sm L$ where all the blow-up limits of
$E$ are $\bY$-cones. We also know, by the regularity theorem of \cite{Ta} or \cite{C1}, 
that near these points $E$ is $C^1$-equivalent to a $\bY$-cone.

Let us put the description of $E$ near $0$ before the statement of Theorem \ref{t29d.1}
because it is a little longer than the previous ones. 
Set $B_0 = B(0,R/10)$, call $H_1$ and $H_2$ the two half planes that compose $X$, then let
$e_i$ denote the unit vector of $H_i$ that is perpendicular to $L$, and let
$H_0$ be the half plane bounded by $L$ and that contains $e_0 = e_1+e_2$;
we may call it the vertical half plane.
First of all (we will show that)
\begin{equation} \label{29d.2}
L \cap B_0 \subset E.
\end{equation}
Next, there is a curve $\gamma$, which is the graph over $L$ of some 
$\tau$-Lipschitz function $\psi : L \to L^\perp$, with $\tau$ as small as we want, 
and which is also of class $C^{1+a}$, (for some small constant $a>0$ which 
may depend on $n$ and $\beta$), with
\begin{equation} \label{29d.3}
|D\psi(x)-D\psi(y)| \leq \tau |x-y|^{a} \ \text{ for } x, y \in L \cap B_0.
\end{equation}
The curve $\gamma$ is contained in $E$, meets $B(0,10^{-10}R)$, and lies in a small sector near $H_0$. That is,
\begin{equation} \label{29d.4}
\gamma \subset H_0^\sharp := \big\{ x\in \R^n \, ; \, \dist(x,H_0) \leq \tau \dist(x,L) \big\}.
\end{equation}
Also
\begin{equation} \label{29d.5}
\gamma \cap B_0 \sm L \subset E_Y \, .
\end{equation}
Then $E$ has a tangent cone at every point of $E \cap B_0$, as follows.
\begin{equation} \label{29d.6}
\text{For $x\in \gamma \cap L \cap B_0$, $E$ has a unique tangent cone $V(x)$ at $x$}
\end{equation} 
such that 
\begin{equation} \label{29d.7}
V(x) \in \bV(L) \ \text{ and } d_{0,1}(V(x),X) \leq \tau.
\end{equation}
In addition, if $x\in \gamma \cap L \cap B_0$ is a boundary point (in $L$) 
of $\gamma \cap L$, then $V(x)$ is sharp.
\begin{equation} \label{29d.8}
\text{For $x\in E \cap L \cap B_0 \sm \gamma$, $E$ has a unique tangent cone $H(x)$ at $x$,
$H(x) \in \bH(L)$}
\end{equation}
(that is, $H(x)$ is a half plane bounded by $L$), and 
\begin{equation} \label{29d.9}
\text{for $x\in E \cap B_0 \sm (L \cup \gamma)$, $E$ has a unique tangent cone $P(x)$ at $x$,
which is a plane.}
\end{equation}
In fact, when $x\in E \cap B_0 \sm (L \cup \gamma)$, there is a small neighborhood
of $x$ where $E$ is a $C^{1+a}$ submanifold of dimension $2$ of $\R^n$.

\begin{figure}[!h]  
\centering
\includegraphics[width=6.5cm]{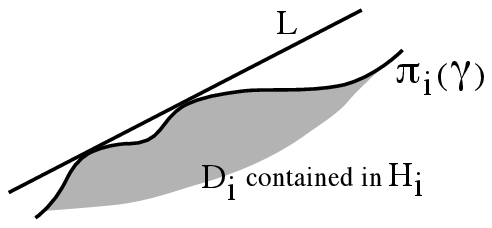}
\hskip1.5cm 
\includegraphics[width=6.cm]{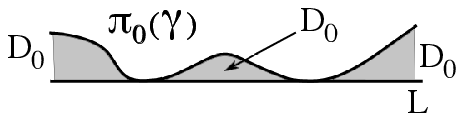}
\caption{ On the left, $D_i$ and the projection of $\gamma$ on $H_i$. On the right, the projection $D_0$
of $F_0$ on the vertical half plane $H_0$
\label{f29d2}}
\end{figure}

Next $E \cap B_0$ is composed of three main (closed) pieces that meet along $\gamma$. 
The first two, which will be called the faces $F_1$ and $F_2$, are bounded by $\gamma$ 
and correspond to the two half planes $H_i$ that compose $X$. 
More precisely, inside $B_0$, $F_i$ is the Lipschitz graph of 
some $\tau$-Lipschitz function $\varphi_i : D_i \to H_i^\perp$, where $D_i$ is 
(inside $B_0$) the subdomain of $H_i$ bounded by the orthogonal projection $\pi_i(\gamma)$ 
and that contains the largest part of $H_i$ near $0$ (See the left part of Figure \ref{f29d2} %%
and notice that $\pi_i(\gamma) \subset H_i$ because $\gamma \subset H_0^\sharp$). In addition, 
\begin{equation} \label{29d.10}
|D\varphi_i(x)-D\varphi_i(y)| \leq \tau |x-y|^{a} \ \text{ for } x,y \in D_i \cap B_0.
\end{equation}
These are the two largest pieces. The last one, which will be called $F_0$,
is bounded by $\gamma \sm L$ on one side and $L\sm\gamma$ on the other side;
it may be empty (if $\gamma \subset L$), or on the opposite composed of infinitely many pieces
(if $\gamma \sm L$ has infinitely many connected component), and looks like a succession of
thin nearly vertical surfaces that connects $L \sm \gamma$ to $F_1$ and $F_2$ along 
$\gamma \sm L$.
Inside of $B_0$, this piece is the graph of a $\tau$-Lipschitz function $\varphi_0$, 
defined on a domain $D_0 \subset H_0$ bounded by $L$ and the orthogonal projection 
$\pi_0(\gamma)$ on $H_0$, and with values in $H_0^\perp$. 
See the right part of Figure \ref{f29d2}. %
And as usual
\begin{equation} \label{29d.11}
|D\varphi_0(x)-D\varphi_0(y)| \leq \tau |x-y|^{a} \ \text{ for } x,y \in D_0 \cap B_0.
\end{equation}
This completes our description of $E$; hopefully we did not forget anything important.
Now we give the corresponding statement.

% mis page 283

\begin{thm}\label{t29d.1} 
There is a constant $a>0$ that depends only on $n$ and $\beta$ and, 
for each small $\tau > 0$, a constant $\varepsilon_0 > 0$, that depends only on $n$, $\beta$, 
and $\tau$, with the following properties.
Let $E$, $h$, $R$, satisfy \eqref{29a.1}-\eqref{29a.3}, and assume that \eqref{29a.4}
holds for some sharp $\bV$-set $X$ bounded by $L$. Then $E$ has the description in 
$B_0 = B(0,R/10)$ that was given just before this statement.
\end{thm}

\ms
The theorem looks complicated but hopefully the pictures give a good idea of what
we claim is going on. Thus $E$ does not leave $L$ entirely, but two of its branches may 
(hesitate and then) go away, keeping their $2\pi/3$ angles, and a thin triangular piece 
stays and connects this main piece $F_1 \cup F_2$ to $L$. 
We really expect this sort of behavior to happen, typically when the film is pulled up by some
force that lifts the two wings $F_1$ and $F_2$.
But for minimal sets it is not clear to the author that we can produce curves $\gamma$ 
that leave from $L$ and return to it as many times as we want.
Also, we proved in the previous two sections that this sort of leaving behavior only happens for
sharp tangent $\bV$-sets, not when they are flatter. The creases of Section \ref{S29c}
that are produced when $E$ looks like a plane and hesitates to go away are different, a plane
may even cross $L$ transversally without even noticing (see the next section for confirmation).

\ms
We start the proof as in the previous case, and first observe that all the bad points are 
close to $L$. That is, if we call
\begin{equation} \label{29d.12}
E_P = \big\{ x\in E \sm L \, ; \, \theta_x(0) = \pi \big\},
\end{equation}
the set of points $x\in E \sm L$ where all the blow-up limits of $E$ at $x$ are planes
(and hence $E$ is a smooth surface near $x$ and has a unique blow-up limit at $x$), 
then there is a constant $C \geq 0$, that depends only on $n$ and $\beta$, such that
\begin{equation} \label{29d.13}
\dist(x,L) \leq C \varepsilon_0 R \ \text{ for } E \cap B(0,R/2) \sm E_P.
\end{equation}
The point is that if $x \in E \cap B(0,R/2)$ and if $d : = \dist(x,L) \geq C \varepsilon_0 R$, 
then \eqref{29a.4} implies that $d_{x,d/2}(E,P) \leq 2C^{-1}$ for the plane $P$
that coincides with $X$ near $B(x,d/2)$; if $C$ is large enough and $\varepsilon_0$
is small enough (to control the gauge function $h$), the standard regularity theorem
implies that $E$ is of class $C^1$ near $x$. The details are the same as for the first lines of 
Lemma~\ref{t29b.2}.

Next we control some densities. By the upper semicontinuity lemma 
(see for instance Lemma~22.3 in \cite{Sliding}), and as for \eqref{29a.7}, 
we find that (if $\varepsilon_0$ is small enough) for $x\in E \cap B(0,R/2)$
and $R/10 \leq r \leq R/3$, 
\begin{equation} \label{29d.14}
\H^2(E \cap B(x,r)) \leq \H^2(X \cap B(x,r))) + \varepsilon_1 R^2,
\end{equation}
where $\varepsilon_1 > 0$ is as small as we want. This will be useful to control
our density functionals. We start our study with the points $x\in E_V \cap B(0,R/2)$, where
\begin{equation} \label{29d.15}
E_V = \big\{ x\in E \cap L \, ; \, \theta_x(0) = \pi \big\}.
\end{equation}
We do not assume a priori that their blow-up limits are $\bV$ sets, but this will come soon.
Pick $x\in E_V \cap B(0,R/2)$. Notice that for $r = R/3$, \eqref{29d.14} says that
\begin{equation} \label{29d.16}
\theta_x(r) := r^{-2} \H^2(E \cap B(x,r)) \leq r^{-2}\H^2(X \cap B(x,r))) + 9\varepsilon_1
= \pi + 9\varepsilon_1.
\end{equation}
Then, by the almost monotonicity of density and if $\varepsilon_0$ is small enough,
\begin{equation} \label{29d.17}
\theta_x(r) \leq \pi + 10\varepsilon_1 \ \text{ for } 0 < r \leq R/3.
\end{equation}
In particular, the blow-up limits of $E$ at $x$ are planes or $\bV$-sets (we excluded
half planes because $x\in E_V$). But what matters to us at this point is that we 
can apply Theorem \ref{t22.2n}; we find that $E$ has a tangent cone $V(x)$ at $x$,
and that
\begin{equation} \label{29d.18}
d_{x,r}(E,V(x)) \leq c_1(\varepsilon_0) (r/R)^{a/4}
\ \text{ for } 0 < r \leq R/3.
\end{equation}
In addition, since both $V(x)$ and $X$ are close to $E$ in $B(x,R/3)$, we also get that
\begin{equation} \label{29d.19}
d_{0,1}(V(x),X) \leq 2c_1(\varepsilon_0).
\end{equation}
Thus $V(x)$ is a set of type $\bV$, and even sharp or almost sharp, and 
\eqref{29d.18} will give us enough control on balls centered on $E_V \cap B(0,R/2)$.
Notice that \eqref{29d.18} also holds (with a slightly worse constant)
when $R/3 < r < R/2$, by \eqref{29a.4} and 
\eqref{29d.19}.

\ms
Next we consider points of $E \cap  B(0,R/2) \sm (L \cup E_P)$, and for these points 
we prefer to use the functional $F_x(r)$ defined by \eqref{22.3}, but for the set $E-x$ 
because of our choice of origin. We claim that for $r = R/3$ and as for \eqref{29d.16},
\begin{eqnarray} \label{29d.20}
F_x(r) &=&  r^{-2} \H^2(E \cap B(x,r)) + r^{-2} \H^2(S_x \cap B(x,r))
\\ \nn
&\leq& r^{-2} \H^2(X \cap B(x,r)) + 9 \varepsilon_1 + \frac{\pi}{2}
\leq \frac{3\pi}{2} + 10\varepsilon_1,
\end{eqnarray}
where $S_x$ is the shade of $L$ seen from $x$, $r^{-2} \H^2(S_x \cap B(x,r)) \leq \frac{\pi}{2}$
is always true because $S_x$ is at most a half plane, and the last line comes from the fact
that since $\dist(x,L) \leq C \varepsilon_0 R$ by \eqref{29d.13},
$r^{-2} \H^2(X \cap B(x,r))$ is as close as we want to $\pi$. As before, the near monotonicity of
$F_x$ then yields
\begin{equation} \label{29d.21}
F_x(r) \leq \frac{3\pi}{2} + 11\varepsilon_1 
\ \text{ for } 0 < r \leq R/3,
\end{equation}
and since $\theta_x(0) := \lim_{r \to 0} \theta_x(r) = \lim_{r \to 0} F_x(r)$ cannot take
values strictly between $\frac{3\pi}{2}$ and $\frac{3\pi}{2} + 11\varepsilon_1$,
we see that $\theta_x(0) \leq \frac{3\pi}{2}$. In fact,  $\theta_x(0) = \frac{3\pi}{2}$
because we assumed that $x \in E \sm (L \cup E_P)$, so $x \in E_Y$. In other words, we checked
along the way that
\begin{equation} \label{29d.22}
E \cap  B(0,R/2) \sm L \subset E_P \cup E_Y.
\end{equation}
We continue as when $x\in E_V$, but with the decay statement associated to $F_x$.
We first apply Theorem \ref{t23.2} to $E-x$ in $B(0,R/2)$; the constraint \eqref{23.8} 
on $h$ is satisfied if $\varepsilon_0$ is small enough, the density of $E$ at $x$ 
is $\frac{3\pi}{2}$ as required in \eqref{23.9},
and $d_{x,R/2}(E,X)$ is as small as we want, by \eqref{29a.4}. Recall that the numbers $\beta_{VP}$
allow approximation by $\bV$ sets or planes (see \eqref{23.4}), but we don't need planes here. 
We get that \eqref{23.11} holds, we take $r_1 = r$ and $r_2 = R/3$, and obtain that
\begin{equation} \label{29d.23}
F_x(r) - \frac{3 \pi}{2} \leq \Big(\frac{C_V 3r}{R}\Big)^{a} [F_x(R/3) -\frac{3 \pi}{2}] 
+ C_V C_h r^a R^{\beta-a}
\leq C \Big(\frac{r}{R}\Big)^{a} (\varepsilon_1 + \varepsilon_0)
\end{equation}
for $0 \leq r \leq R/4$, by \eqref{29d.20}, \eqref{29a.2}, and \eqref{29a.3}.
Now we use Theorem \ref{t28.2}, with a radius $r \leq R/1200$
(so that $400r \leq R/3$). Because of our assumption that $\dist(x,L) \leq r/2$, we restrict to $r$ such that
\begin{equation} \label{29d.24}
2\dist(x,L) \leq r \leq R/1200.
\end{equation}
Notice that there are still lots of radii $r$ available, since \eqref{29d.13} says that 
$\dist(x,L) \leq C \varepsilon_0 R$. 
The theorem provides a set $Y = Y(x,r) \in \bY(L,r)$ such that 
\begin{equation} \label{29d.25}
d_{x,r}(E,Y^t(x,r)) \leq C_6 \Big[ [F_x(200r) -\frac{3\pi}{2}] + C_h R^{\beta} \Big]^{1/4}
\leq C \Big(\frac{r}{R}\Big)^{a/4} \varepsilon_1^{1/4}
\end{equation}
as in \eqref{28.8}, where $Y^t(x,r)$ is corresponding truncated cone, and by \eqref{29d.24}.
It may be a little awkward to compare the sets $Y^t(x,r)$, because it may be that the 
truncation of $Y^t(x,r)$ is no longer clean in larger balls, but at least we can control their
directions.  Denote by $Y'(x,r) = Y(x,r) - x$ the translation of $Y(x,r)$ that is centered at $0$; 
we claim that 
\begin{equation} \label{29d.26}
d_{0,1}(Y'(x,r),Y'(x,s)) \leq C \big(\frac{r}{R}\big)^{a/4} \varepsilon_1^{1/4}
\ \text{ for } 2\dist(x,L) \leq s \leq r \leq R/1200.
\end{equation}
When $r \leq 2s$, this is because we can easily evaluate the position of the
two half planes (among the three that compose $Y(x,r)$) that do not contain 
$L \cap B(x,r)$, by knowing \eqref{29d.25}; of course the third half plane comes for free
when we have the two main ones. The general case follows as usual: we just compare successive 
cones and sum a geometric series.

When we use \eqref{29d.25} with $r = R/1200$ and compare with \eqref{29a.4}, we see that
the two half planes of $Y(x,R/1200)$ that compose $Y(x,R/1200)$ and do not meet $L \cap B(x,R/1200)$
are very close to the two pieces $H_1$ and $H_2$ of $X$. This is not shocking, because the
third part of $Y^t(x,R/1200)$ is extremely thin, because $\dist(x,L) \leq C \varepsilon_0 R$.
Because of this closeness, and then by \eqref{29d.26}, we can identify in 
$Y(x,r)$ two half planes $H_1(x,t)$ and $H_2(x,t)$, whose directions $H'_1(x,t) = H_1(x,t)-x$ 
and $H'_2(x,t)= H_2(x,t)-x$ are very close to $H_1$ and $H_2$ respectively. That is,
\begin{equation} \label{29d.27}
d_{0,1}(H'_i(x,t),H_i) \leq 10 \varepsilon_1 \ \text{ for } i=1,2.
\end{equation}
And when we use \eqref{29d.26} and check the labelling, we find out that for $i=1,2$,
\begin{equation} \label{29d.28}
d_{0,1}(H'_i(x,s),H'_i(x,t)) \leq C \big(\frac{r}{R}\big)^{a/4} \varepsilon_1^{1/4}
\ \text{ for } 2\dist(x,L) \leq s \leq r \leq R/1200.
\end{equation}
Now we keep our point $x \in E_Y \cap  B(0,R/2)$ and worry about smaller balls. 
Set $d(x) = \dist(x,L)$ and $Y_d(x) = Y(x,2d(x))$, and observe that by \eqref{29d.25}
\begin{equation} \label{29d.29}
d_{x,d(x)/2}(E,Y_d(x)) = d_{x,d(x)/2}(E,Y^t_d(x))
\leq 4 d_{x, 2d(x)}(E,Y^t_d(x))
\leq C \big(\frac{d(x)}{R}\big)^{a/4} \varepsilon_1^{1/4},
\end{equation}
so $E$ lies pretty close to a $\bY$-set in $B(x,d(x)/2)$.
In fact, if $\varepsilon_0$ and $\varepsilon_1$ are small enough, this information is
sufficient to use the analogue of Theorem \ref{t22.2n} for plain almost minimal sets near
a $\bY$-set. This result dates from \cite{Ta}, we can also refer to \cite{C1}, but since the exact
same statement does not seem to be explicitly written there, let us also observe that the 
proof of Theorem \ref{t22.2n} also goes through (with simplifications). 
Anyway, we get that $E$ has a tangent cone $Y(x)$
(a $\bY$-set centered at $x$), and the analogue of \eqref{21.8} is that
\begin{equation} \label{29d.30}
d_{x,r}(E,Y(x)) \leq c \Big(\frac{r}{d(x)}\Big)^{a/4}
\ \text{ for } 0 < r \leq d(x)/2,
\end{equation}
maybe with a different constant $a > 0$, and with a constant $c > 0$ that can be taken as small as
we want (by making $\varepsilon_0$ and $\varepsilon_1$ smaller).
In fact, when we combine this and \eqref{29d.29} (to take care of the radii that are too close to $d$),
the same trick as for \eqref{29b.18} shows that 
\begin{equation} \label{29d.31}
d_{x,r}(E,Y(x)) \leq c \Big(\frac{r}{R}\Big)^{b/4}
\ \text{ for } 0 < r \leq d(x)/2,
\end{equation}
for some (other) $b < a$. When we compare this to \eqref{29d.29}, we find that
\begin{equation} \label{29d.32}
d_{0,1}(Y'(x),Y'(x, 2d(x))) \leq c' \Big(\frac{d(x)}{R}\Big)^{b/4},
\end{equation}
where $Y'(x) = Y(x) - x$ is the parallel cone through the origin, and 
$c'$ is another positive constant that we can make as small as we want. And for even larger 
radii, we combine this with \eqref{29d.26} and get that 
\begin{equation} \label{29d.33}
d_{x,1}(Y(x),Y(x,r)) = d_{0,1}(Y'(x),Y'(x,r)) \leq c' \Big(\frac{r}{R}\Big)^{b/4}
\ \text{ for } 2d(x)  < r \leq R/1200.
\end{equation}
With this, \eqref{29d.31}, and \eqref{29d.25}, we get good enough estimates 
on the approximations of $E$ in the small balls $B(x,r)$ centered on $E_Y \cap  B(0,R/2)$.

\ms
We want to know a bit more about the set $E_Y$ (the set of  \eqref{29d.1}) itself. Observe that
\begin{equation} \label{29d.34}
 E_Y \cap  B(0,R/2) \subset H_0^\sharp,
\end{equation}
the small sector $H_0^\sharp$ of \eqref{29d.4}. 
Indeed, \eqref{29d.27} says that the two main branches of 
$Y_d(x) = Y(x,\delta^{-1}d(x))$ go in directions very close to $H_1$ and $H_2$,
hence the third one goes in a direction almost opposite to $H_0$. 
And it is precisely in this direction that we can find the point $y\in L$ that lies closest to $x$.
So the direction of $x-y$ lies very close to $H_0$ (compared to $d(x)$), which proves
\eqref{29d.34}.

Next, the regularity theorem for plain almost minimal sets that $E_Y$
is a $C^1$ curve (locally in $B(0,R/2) \sm L$), and its tangent line at $x \in E_Y$ 
is the spine $\ell(x)$ of $Y(x)$. Also denote by $\ell'(x)$ the direction of $\ell(x)$; that is,
$\ell'(x)=\ell(x)-x$.
It follows from \eqref{29d.32} that $Y(x)$ and $Y_d(x) = Y(x,\delta^{-1}d(x))$
are as close to each other as we want, and hence also their spines are as close to each other
as we want. Then by \eqref{29d.27}, the spine of $Y_d(x)$ is also as close to the direction 
of the spine of $X$ as we want. The spine of $X$ is $L$, and its direction is $L'=L-x$;
then $\ell'(x)$ is as close to $L'$ as we want.
That is, in $B(0,R) \sm L$, $E_Y$ is locally the graph of a $\tau$-Lipschitz function over 
its projection on $L$. But we can say a bit more on the variations of $\ell$. We claim that
\begin{equation} \label{29d.35}
d_{0,1}(\ell'(x),\ell'(y)) \leq  c \Big(\frac{|x-y|}{R}\Big)^{b/4}
\ \text{ for } x,y\in E_Y \cap B(0,R/2),
\end{equation}
where as usual $c$ can be made as small as we want.
Set $r = 3 |x-y|$. For $r \geq R/1200$, we use \eqref{29d.27}, which as we saw recently 
implies that $\ell'(x)$ and $\ell'(y)$ are as close to $L'$ as we want. 
Otherwise, we compare the two approximation estimates that we have for $E$ in $B(x,r)$ 
and $B(y,r)$, which both contain $B(x,2|x-y|)$. 
For $x$, we use $Y(x)$ and \eqref{29d.31} if $r < d(x)$, and $Y^t(x,r)$ and \eqref{29d.25}
otherwise. For $y$ we proceed similarly, and when we compare we get that $Y^t(x,r)$,
for instance, is close to $Y(y)$ in $B(x,2|x-y|)$. We use \eqref{29d.33} if needed 
to return to $Y(x)$ and $Y(y)$, observe that for \eqref{29d.35} it is enough to control 
the two main faces, and conclude. Let us also check that
\begin{equation} \label{29d.36}
d_{0,1}(\ell'(x),L') \leq  c \Big(\frac{\dist(x,E_V \cap B(0,R/2))}{R}\Big)^{b/4}
\ \text{ for } x\in E_Y \cap B(0,R/2).
\end{equation}
If $\dist(x,E_V) \geq R/2500$ (which also happens if $E_Y \cap B(0,R/2)$ is empty),
we use \eqref{29d.27} again and get that $d_{0,1}(\ell'(x),L') \leq c$. 
Otherwise, pick $y \in E_Y \cap B(0,R/2)$ whose distance to
$x$ is almost minimal, and apply the same argument as above with $r = 2|x-y|$, 
but where we use $V(y)$ and \eqref{29d.18} to approximate $E$ in $B(y, r)$. This time
the analogue of $\ell(y)$ (the spine of $V(y)$) is $L$, and \eqref{29d.36} follows as 
for \eqref{29d.35}.

Before we continue with $E_Y$, let us observe that
\begin{equation} \label{29d.37}
L \cap B(0,2R/3) \subset E,
\end{equation}
by Lemma \ref{t29b.4}, which is still valid for sharp $\bV$-sets $X$. In addition, we claim that
\begin{equation} \label{29d.38}
L \cap B(0,2R/3) \sm E_V \subset E_H,
\end{equation}
where 
\begin{equation} \label{29d.39}
E_H = \big\{ x\in E \cap L \, ; \, \text{ the density of $E$ at $x$ is } 
\theta_x(0) = \frac{\pi}{2}\big\}.
\end{equation}
By Theorem \ref{t29a.1}, if $x\in E_H$, $E$ has a tangent half plane $H(x)$ at $x$, and 
$E$ is even $C^1$-equivalent to $H(x)$ near $x$; whence the name $E_H$.

Let us check \eqref{29d.38}. Let $x\in L \cap B(0,R/2)\sm E_V$ be given.
By \eqref{29d.14} and the proof of \eqref{29d.17} (still valid),
we get that $\theta_x(0) \leq \pi + 10\varepsilon_1$.
Since $\theta_x(0)$ is also the density of the blow-up limits of $E$ at $x$,
Lemma~\ref{t22.2} says that $\theta_x(0) \in \{ \frac{\pi}{2}, \pi \}$; 
in addition, since $x \notin E_V$ and by the definition \eqref{29d.15}, we are left with 
$\theta_x(0) = \frac{\pi}{2}$; \eqref{29d.38} follows.

We return to $E_Y$ and $E_V$. Set $\gamma = E_Y \cup E_V$. Notice that
$\gamma$ is closed in $B(0,R/2)$, because \eqref{29d.22} and \eqref{29d.38} say that
\begin{equation} \label{29d.40}
E \cap B(0,R/2)) = B(0,R/2) \cap [E_V \cup E_H \cup E_Y \cup E_P],
\end{equation}
and Theorem \ref{t29a.1} implies that $E_H$ is open in $L$, while 
$E_P$ is closed in $E$ by the regularity theorem. Let us check that
\begin{equation} \label{29d.41}
\pi(\gamma \cap B(0,R/2)) \supset L \cap B(0,R/3),
\end{equation}
where $\pi$ still denotes the orthogonal projection on $L$.
Suppose this fails, and pick $y\in L\cap B(0,R/3) \sm \pi(\gamma \cap B(0,R/2))$.
We start with the case when $\dist(y,\gamma) \leq R/20$. Choose $x\in \gamma$
such that $|\pi(x)-\pi(y)| + |\pi^\perp(x)-\pi^\perp(y)|$ is minimal, where $\pi^\perp = I - \pi$;
the existence and the reason for this strange choice of distance will appear very soon.
The infimum is at most $R/10$ because $\dist(y,\gamma) \leq R/20$, so the 
good competitors lie well inside of $B(0,R/2)$, and $x$ exists because $\gamma$ is closed there.
Of course $x\in B(0,R/2)$. Suppose for a minute that $x\in E_Y$. Recall from the discussion
above that there is a small neighborhood of $x$ where $E_Y$ is a $C^1$ curve that runs
nearly parallel to $L$. This is because $\ell'(x)$ is so close to $L'$.
But when we follow that curve so that the projection gets closer to $y$, we see that 
our strange distance to $y$ gets strictly smaller (the variations of $\pi^\perp(x)$ are too
small to compensate). This contradicts the definition of $x$, and so $x\in E_V$. 

We now run the same topological argument as for Lemma \ref{t29b.3}.
In the ball $B(x,2|x-y|)$, \eqref{29d.18} says that $E$ looks a lot like $V(x)$.
Yet in small balls centered at $y$, $E$ is $C^1$-equivalent to $H(y) \in \bH$.
We can find a small $(n-d)$ sphere $C_0$ centered at $y$ and contained in the hyperplane
orthogonal to $L$, that meets $E$ exactly once, transversally. We can also find another
$(n-d)$-sphere centered $C_1$, in the same plane, of larger radius $r/100$, and this one meets 
$E$ twice transversally (once near each wing of $V(x)$). But the obvious homotopy from $C_0$
to $C_1$ stays far from $E_Y$, because 
$|\pi(z)-\pi(y)| + |\pi^\perp(z)-\pi^\perp(y)| \geq |x-y|$ by minimality of $x$.
Then we can proceed as in Lemma \ref{t29b.3} (and in fact \cite{Holder}) to find a nicer
homotopy, show that the number of intersections stays constant modulo $1$, and 
obtain the desired contradiction. 

We are left with the case when $\dist(y,\gamma) \geq R/20$; let us check that
this never happens, and in fact that 
\begin{equation} \label{29d.42}
\dist(y,\gamma) \leq C \varepsilon_0 R
\ \text{ for } y \in L \cap B(0,R/3).
\end{equation}
Indeed suppose that $y \in L \cap B(0,R/3)$ and $\dist(y,\gamma) > C \varepsilon_0 R$
Obviously $x\notin E_V$, so $x\in E_H$, and as before we can find  
a small $(n-d)$-sphere $C_0$ centered at $y$, contained in the hyperplane
orthogonal to $L$, and that meets $E$ exactly once and transversally.
Also consider the $(n-d)$-sphere  $C_1$ of radius $C \varepsilon_0 R/2$ centered at 
$y$ and contained in the same hyperplane. Recall from \eqref{29a.4} that in $B(0,R)$, 
$E$ stays $\varepsilon_0 R$-close to the $\bV$-set $X$. If $C$ is large enough, we can apply the
standard regularity theorem, as we did for \eqref{29d.13} and the first lines of 
Lemma~\ref{t29b.2}, to find that $E$ is a $C^1$ surface near the two points of $C_1 \cap X$.
Moreover (because $E$ is also a small Lipschitz graph over the plane that contains the corresponding
$H_i$), $C_1$ cuts $E$ transversally, once near each point of $C_1 \cap X$. 
The same topological argument as above, using the fact that by definition $B(y, C \varepsilon_0 R)$
does not meet $E_Y$, gives the contradiction that proves \eqref{29d.42}.
This also concludes our proof of \eqref{29d.41}, and of the fact that $\gamma$ meets 
$B(0,10^{-10}R)$.
\ms
We are now able to prove that $\gamma \cap B(0,R/3)$ satisfies all the requirements
mentioned before Theorem~\ref{t29d.1}. First observe that when $x\in E_V$ and $r < R/3$, 
it follows from \eqref{29d.18} and the proof of \eqref{29d.13} that
\begin{equation} \label{29d.43}
\dist(y,L) \leq C c_1(\varepsilon_0) (r/R)^{a/4} r 
\ \text{ for } y\in E_Y \cap B(x,r/2).
\end{equation}
Since $C c_1(\varepsilon_0) (r/R)^{a/4}$ tends to $0$ with $r$, it follows that
$\gamma$ has a tangent at $x$, equal to $L$.

Next we check that $\gamma$ coincides with a Lipschitz graph in $B(0,R/3)$.
Let $x, y \in \gamma \cap B(0,R/3)$ be given; we want to show that 
\begin{equation} \label{29d.44}
|\pi^\perp(x)-\pi^\perp(y)| \leq \tau |\pi(x)-\pi(y)|,
\end{equation}
with $\tau >0$ as small as we want. When $x, y \in E_V$, this is trivial. 
When $x \in E_V$ and $y\in E_Y$, this follows from \eqref{29d.43} with $r=3|x-y|$ 
(or directly from \eqref{29d.13} if $r$ is large). When $x, y \in E_Y$, 
either there is a curve in $E_Y$ that goes from $x$ to $y$, and then we use the fact 
that the direction $\ell'(x)$ of the tangent stays close to $L'$, or else we can find $z\in L$, 
between $\pi(x)$ and $\pi(y)$, which does not lie in $\pi(E_Y \cap B(0,R/2))$; 
this point lies in $E_V$ by \eqref{29d.41}, and now we can apply
\eqref{29d.43} to the pairs $(z,x)$ and $(z,y)$ to get \eqref{29d.44}.

So $\gamma \cap B(0,R/3)$ is the graph of a $\tau$-Lipschitz function $\psi$; 
the estimate \eqref{29d.3} on $D\psi$ comes from \eqref{29d.35}, \eqref{29d.36}, and
\eqref{29d.43} (because $E_V \subset L$ anyway). We already checked 
\eqref{29d.2} (see \eqref{29d.38}), \eqref{29d.4} 
(see \eqref{29d.34} and notice that $E_Y \subset L \leq H_0^\sharp$), 
\eqref{29d.5} is part of the definition of $\gamma$, the tangent cone 
properties \eqref{29d.6}-\eqref{29d.9} were also checked before, and the fact that
$V(x)$ is sharp when $x$ is a boundary point of $E_V$ in $L$ follows from Theorem \ref{t29b.1}.

\ms
We are left with the task of establishing the description of $E$ in terms of faces.
As usual, this mostly mean proving estimates on the approximation of $E$ in balls $B(x, r)$.
After $x\in E_V$ (see \eqref{29d.18}) and $x\in E_Y$ (see \eqref{29d.25} and \eqref{29d.31}),
the next piece of $E$ in the hierarchy is $E_H$. So we take $x\in E_H \cap B(0, R/4)$ 
and try to approximate $E$ in $B(x, r)$. 
Set $d(x) = \dist(x,\gamma)$ and pick $y\in \gamma$ such that $|y-x| = \dist(x,\gamma)$. 
Recall from \eqref{29d.42} that $d(x) \leq C \varepsilon_0 R$, so $y$ lies well inside $B(0,R/3)$. 
Also, we easily deduce from the small Lipschitz description of $\gamma$ that $y\in E_Y$. 
Notice also that $\dist(y,L) \leq d(x) \leq 2 \dist(y,L)$, again by the Lipschitz description of $\gamma$.

For $r > d(x)/2$, set $D = \max(2\dist(y,L), d(x)+r)$
(so that $B(y,D)$ contains $B(x, r)$ and \eqref{29d.24} holds for $D$); notice that $D \leq 4r$ and
apply \eqref{29d.25} in $B(y,D)$; we get that
\begin{equation} \label{29d.45}
d_{x, r}(E,Y^t(y,D)) \leq \frac{D}{r} \, d_{y,D}(E,Y^t(y,D)) 
\leq C \frac{D}{r} \Big(\frac{D}{R}\Big)^{a/4} \varepsilon_1^{1/4}
\leq  C \big(\frac{r}{R}\big)^{a/4} \varepsilon_1^{1/4},
\end{equation}
which is as good as usual. For $r < d(x)/2$, we first apply \eqref{29d.24}
with the radius $D=2\dist(y,L)$, and notice that near $B(x, d(x)/2)$,
the truncated set $Y^t(y,D)$ coincides with a half plane $H_0(x) \in \bH$
(compare with the definition of $\bY(y,D)$ above Theorem \ref{t28.2}); hence
\begin{equation} \label{29d.46}
d_{x, d(x)/2}(E,H_0(x)) = 
d_{x, d(x)/2}(E,Y^t(y,D)) \leq \frac{2D}{d(x)} \, d_{y,D}(E,Y^t(y,D)) 
\leq  C \Big(\frac{d(x)}{R}\Big)^{a/4} \varepsilon_1^{1/4}.
\end{equation}
The right-hand side is still as small as we want, which allows us to apply Theorem \ref{t22.2n}.
We get that $E$ has a tangent half plane $H(x)$ at $x$ (which we knew already), and that
\begin{equation} \label{29d.47}
d_{x,r}(E,H(x)) \leq c  \Big(\frac{r}{d(x)}\Big)^{a/4}
\ \text{ for } 0 < r \leq d(x)/2.
\end{equation}
With the usual manipulation (use \eqref{29d.46} instead of \eqref{29d.47} when
$r$ is rather close to $d(x)$, and proceeding as for \eqref{29b.18}), we also have 
\begin{equation} \label{29d.48}
d_{x, r}(E,H(x)) \leq c  \Big(\frac{r}{R}\Big)^{b/4}
\ \text{ for } 0 < r \leq d(x)/2,
\end{equation}
with $c$ as small as we want and an even smaller $b > 0$. In the proof we first get it with
$H_0(x)$, but we can compare $H_0(x)$ and $H(x)$ using an intermediate radius
(the radius where we switch from \eqref{29d.46} to \eqref{29d.47}).

\ms
We now go to the last level of the hierarchy, and take $x\in E_P \cap B(0,R/9)$.
As usual, set $d(x) = \dist(x, L \cup E_Y) \leq R/8$ (because $E$ contains points near the origin), 
and pick $y\in L \cup E_Y$ such that $|y-x| = \dist(x,L \cup E_Y)$. 
For $r \geq d(x)/2$, we can use the good description of $E$ that we already got in
$B(y, d(x)+r)$, and get a good approximation of $E$ by a set $Z(x, r)$, of the form 
$V(y)$ (if $y\in E_V$), $Y^t(y, d(x)+r)$ or $Y(y)$ (if $y \in E_Y$ and depending on 
whether $y$ lies far from $L$ or not), or $Y^t(z, d(x)+r)$ or $H(y)$ (if $y\in E_H$ and 
depending on its distance to $E_Y$). 

When $r \leq d(x)/2$, we first use the estimate for $d(x)$, notice that $Z(x, d(x))$ coincides 
with a plane in $B(x,2d(x)/3)$ and so $E$ is well approximated by a plane in $B(x, d(x)/2)$, 
apply the regularity theorem, get an analogue of \eqref{29d.47}, and then by the now usual trick
an analogue of \eqref{29d.48}. We skip the details because they are quite similar to what we did before.

\smallskip
At this point we have good H\"older estimates on the variations of the direction of the tangent 
to $E$ at $x\in E_P$, which are thus valid as long as we stay in a face, i.e., when we do not cross
$L$ or $\gamma$. We still need to check that the organization of the faces is as in our initial
description of $E$ near $0$. In particular we need to be able to recognize faces.

Let $x\in E_P$ be given, and let $y\in L \cup E_Y$ such that $|y-x| \leq 10 d(x) :=\dist(x,L \cup E_Y)$
(we will use the factor $10$ to choose points preferably in $E_V$ or $E_Y$). 
The simplest case is when we can take 
$y\in E_V$; then by \eqref{29d.18} (applied with $r=20d(x)$) 
$x$ lies close to one of the two branches $H_i(y)$ that compose $V(y)$, and more precisely
$\dist(x,H_i(y)) \leq \tau d(x)$, with $\tau$ as small as we want. In this case, $x\in F_i$,
and it is even easy to find a nice path from $x$ to $y$, in $E_P$, by applying \eqref{29d.18}
and the regularity theorem in the successive balls $B(y,2^{-k+1} d(x))$. We also get the Lipschitz
description with $\varphi_i$ and \eqref{29d.10} near $x$, because we control both the position of
points (directly by \eqref{29d.18}), and because variations of the direction of the tangent planes 
are controlled by our approximation estimates.

The second simplest case is when $x$ is closer to $E_Y$ than to $E_V$, and we can pick $y\in E_Y$; 
then we apply \eqref{29d.25} with $r=20d(x)$. 
Again $x$ lies close to one of the three faces of $Y^t(y,20d(x))$ or $Y(y)$,
depending on whether $20d(x) \leq \dist(y,L)$ or not. If in addition this face is
$H'_1(y,20d(x))$, then we can proceed as before and prove that $x\in F_1$, there is nice path 
in $E_P$ that connects $x$ to $y$, and we even have the desired Lipschitz description of $E = E_P$
near this path. The same thing happens, with $F_2$, if this face is $H'_2(y,20d(x))$.

Otherwise, $\dist(x,H_0(y),20d(x)) \leq \tau d(x)$ for the third face $H_0(y)$ of 
$Y^t(y,20d(x))$ or $Y(x)$. In this case, $x$ lies in the third piece of $E \cap B_0$, 
and we can connect $x$ to $y$ as above. Notice that if we see the vector $e_0 = e_1+e_2$
of $H_0$ as pointing upwards; then $x$ lies almost right below $y$, or in other words
$y-x \in H_0^\sharp$. 

The third and last case is when we have to take $y \in E_H$. 
That is, we may now assume that $|y-x| = d(x)$ but $\dist(x, E_V \cup E_Y) \geq 10 d(x)$.
Again we have a nice approximation of $E$ in $B(y,5d(x))$, which is given given by a half plane $H(y)$
(as in the discussion above, for balls centered on $E_H$). In this case, we get a good description
of $E$ between $x$ and $y$, coming from a mixture of the standard regularity result for plain almost
minimizers that look look a plane (near $x$ and not too close to $y$), and Theorem \ref{t29a.1}
(close to $y$). By looking at radii larger than $d(x)$, up to $d_1(x) = 10\dist(x, E_V \cup E_Y)$,
we find out that $E_Y$ lies much closer to $x$ and $y$ than $E_V$, and $x$ lies in the vertical 
piece $F_0$, essentially below $E_Y$. Hence we can recognize the faces and their regularity comes from
the slow variations of the tangent planes.

Here and for the Lipschitz description of the faces $F_1$ and $F_2$, we are skipping a small part of the 
argument, where we check that the orthogonal projection, for instance from $F_1$ to $H_1$, is
surjective on the domain $D_1$ bounded by $\pi_1(\gamma)$ and away from $L$. Due to the fact
that we have a small bound on the the angle between the tangent plane and, here, the direction of $H_1$,
we don't need a complicated degree argument, and instead we can proceed as for \eqref{29d.41},
and say that if $z\in D_1 \cap B_0$ does not lie in the projection of the face $F_1$, we can consider the 
point $x\in F_1$ that minimizes $|\pi(x)-\pi(z)|+|\pi^\perp(x)-\pi^\perp(z)|$, and then observe
that $x$ does not lie on the boundary of the face (here $\gamma$) and, since near $z$, 
$E$ is a $C^1$ surface with a tangent plane almost parallel to $H_1$, we can find
points near $x$ that do strictly better. The same argument works for $F_2$, and even the third piece
$F_0$, even though the boundary may be more complicated.

This completes our description of $E \cap B_0$ and our proof of Theorem~\ref{t29d.1}.
\qed

\section{The missing case of parallel $\bY$-sets}
\label{S35}

The next interesting case in our study should be when, in the same setting as for 
Theorems~\ref{t29a.1}-\ref{t29d.1}, $E$ is approximated by a minimal cone $X \in \bY(L)$,
thus composed of three half plane bounded by $L$ and that make $\2$ angles along $L$.

Unfortunately, the methods of this paper do not seem to be enough to treat this case directly.
Yet we try to explain why and what could be expected in this situation.

The main problem that we have is that, when $E$ is a $\bY$-set centered at $0$ and that
contains $L$ in one of its face, the monotonicity formula that we have (for the function $F$ of 
\eqref{22.3}) is not adapted. For sure $F$ is monotone, but in this case $F(r) = \frac{3\pi}{2}$
for $r$ small, $F$ starts to increase when $B(0,r)$ meets $L$, and it tends to $2\pi$ at $\infty$.
By contrast, for a truncated set of type $\bY$, $F$ is constant equal to $\frac{3\pi}{2}$, while for
the full $\bY$-set the measure of the shade is computed twice. This is not good, because our typical
proof is based on the observation that in many situations $F$ is nearly constant and then we can control
the geometric situation. In the case of our $\bY$-set, we would be very happy with a different 
monotonicity formula which is adapted to it (either because it counts the shade in a different way, 
or for some other reason), but the author is not sure that such a formula exists.
Less ambitiously one could hope that, in the specific situation where $X \in \bY(L)$
and $E$ is close enough to $X$, there is a quantity that can be controlled and proved to be nearly
monotone. That is, the near monotonicity would specifically use the description of $E \cap \S_r$ 
that we have, and possibly also the net of geodesics that was constructed in Section \ref{S25},
even though different competitors may have to be constructed. That is, it does not seem optimal to
add the triangle $T(r)$ as we did. This last scheme looks more plausible but at this time the author 
was not able to make it work.

Notice that $X \in \bY(L)$ satisfies the full length condition; this is not the point.
So if  we assume that $0$ is a $\bY$-point of $L$, i.e., a point of $E \cap L$ where
the density of $E$ is $\frac{3\pi}{2}$, Theorem~\ref{t22.2n} applies, shows that $E$
has a tangent cone $Y \in \bY(L)$ at $0$, and says that $E$ is very close to $Y$
in small balls $B(0,r)$.
But in the situation of Theorems~\ref{t29a.1}-\ref{t29d.1} with $X\in \bY(L)$,
it could be that $E \cap L$ contains no point of density $\frac{3\pi}{2}$.
And even if we assume that $0$ is such a point, it could be that $E \cap L$ contains 
no other point of density $\frac{3\pi}{2}$; see the expected description below.
Then we do not know how to control directly $E$ in other small balls centered on $L$.
We could try to control $E$ in small balls centered on $E_Y$, as we have done when $X$ is a $\bV$-set,
but for this we used the monotonicity formula, which no longer works here.

\ms
We try now to describe what $E$ should look like in the situation of 
Theorems~\ref{t29a.1}-\ref{t29d.1} but with $X\in \bY(L)$. 
Let us even assume that $0$ is a $\bY$-point of $L$, with a blow-up limit $Y \in \bY(L)$ that contains $L$,
so that the description is interesting starting from $0$.

Just as when $X$ is a plane, one option is that $E$ leaves $L$ right away. 
After all, $\bY$-sets are minimal even without a sliding condition, so $E$ could have been any such set, 
even transverse to $L$, and in the present case it may start parallel to $L$ and just leave it. 
It is also possible that two of the foils of $Y$ leave from $E$, and the last one becomes composed of a 
large face bounded by $L$ and a triangular face, bounded by $L$ on one side, and on the other side
by the set $E_Y$ where it connects to the two foils above.
Finally, we expect a complicated combination of these behaviors to be possible, at least for 
almost minimal sets. As usual, minimal sets may behave differently because have more rigidity.

\begin{figure}[!h]  
\centering
\includegraphics[width=16.cm]{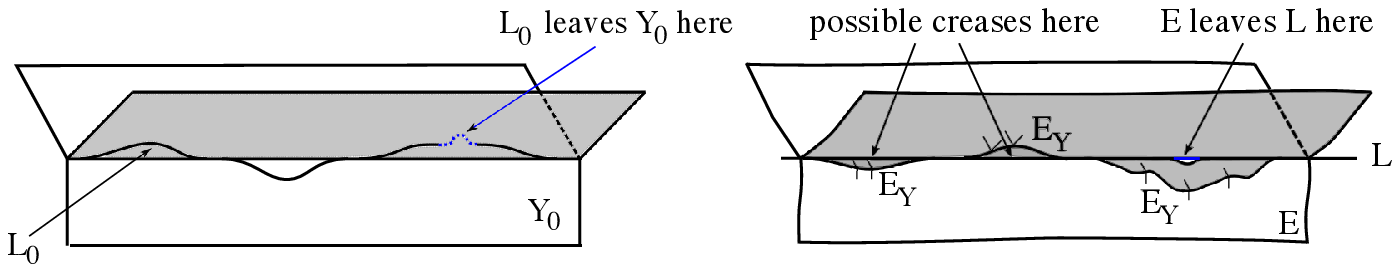}
\caption{Left: A minimal set $Y_0$ and a boundary curve.
Right: The probably almost minimal set $E = G(Y_0)$. \label{f35.1}}
\end{figure}

Before we try to describe this with more detail, let us say how apparently complicated examples
may arise (but we will never actually check the almost minimality of the suggested example below). 
Start from a $\bY$-set $Y_0$ in $\R^3$, say through the origin, but consider a boundary set $L_0$ 
which is not a straight line, but is smooth and is (very!) tangent to the spine $L$ of $Y_0$ at $0$; 
see the left part of the somewhat exaggerated Figure \ref{f35.1}
(the same as Figure \ref{f1.7}).
Obviously $Y_0$ is also sliding minimal with respect to the boundary $L_0$, 
and if $L_0$ is nice enough, we can quite probably (but will not check here) find
a $C^1$ diffeomorphism $G : \R^3 \to \R^3$ that maps $L_0$ to $L$ and $Y_0$ to an almost minimal 
set $E = G(Y_0)$ with sliding boundary $L$. The nice property of $G$ for this to work is that
its derivative should be close to conformal in many places, and in particular along $L_0$, but nothing so
special that it would be impossible to get. Yet $E$ can have a slightly complicated behavior, in particular
concerning its position with respect to $L$. This is what is suggested by the right part of Figure \ref{f35.1},
which will inspire the tentative description below. But first let us remark that
$Y_0$ itself is already an interesting example, because we try in these notes to give descriptions of 
almost minimal sets that stay valid also when the boundary $L$ is a smooth curve, as will be mentioned
in Section \ref{S31}, and the left part of Figure \ref{f35.1} is a perfectly valid example in this context.

 % 294

Return to our tentative description of almost minimal sets $E$ that look like $Y \in \bY(L)$.
First there should be a curve $\gamma$, which is the graph of a $\tau$-Lipschitz function 
$\psi : L \to L^\perp$ and is also $C^{1+a}$-continuous near $0$, as in \eqref{29d.3}.
This curve may leave and return to $L$ many times, but for $x\in \gamma$, $E$ has a unique 
tangent cone which is a $\bY$-set. In the example above, $\gamma = G(L)$.
Part of it lies on $L$ (for instance, the origin), and part of it may be away from $L$,
creating a set $E_Y$ as in the previous section.

There may also be a set $E_V \subset L$ of points of $E \cap L$ where $E$ has a tangent 
$\bV$-set, flat or almost flat. In the example above, this corresponds to parts of $L = G(L_0)$
where $L_0$ has left the spine, but still lies in $Y_0$. At such points, $E$ may have a crease
(and on one side we have a face that leads to $E_Y$, while the other side looks like an infinite
face of $E = G(Y_0)$). At the boundary of $E_V$ in $L$, we may find points of $\gamma$
where $\gamma$ leaves $L$, or points where $E$ has a tangent plane (and $E$ may in fact leave
$L$), but otherwise, if $x\in E_V$ is a point where $E$ has a real crease (i.e., the tangent cone
of $E$ at $x$ is a $\bV$ set which is not a plane), we know from Section \ref{S29b} that
there is a small neighborhood of $x$ where $E$ has a crease and $L$ is composed of points
where $E$ has a tangent cone which is a generic $\bV$-set, almost flat.

Finally, in our attempted description, it is important to say that all the other points of 
$E$ are points of type $\bP$, near which $E$ is a $C^1$ surface. These include
points $x\in E \cap L$ where $E$ casually crosses $L$ (see the discussion of fully transverse cones
in the next section). Also notice that $E$ is allowed to leave $L$ (at a point of $L \cap \gamma$, or at a
point of $E_V$ where the tangent is a plane). But there is no point of type $\bH$ in this discussion,
where $E$ would have a tangent half plane. And there is no point of $E_Y \sm\gamma$.

The description should involve a $C^{1+a}$-behavior, where the direction of the tangent plane (or cone) 
to $x$ varies in a H\" older-continuous way along the faces (and with natural limits on $L$ and $\gamma$),
as in the previous section, but we do not write this down.

\ms
In order to  establish the description above, we could also try to use the information that we have so far 
(essentially from the previous sections) and prove things with topological arguments. 
Obviously the author did not succeed so far. One of the difficulties seems to be that in addition to
the expected nice curves that compose $\gamma \sm L$, the set $E_Y$ of $\bY$-points
of $E \sm L$ may also, a priori, contain some other curves. 
We know that away from $L$, $E_Y$ is locally a simple $C^1$ curve with no endpoints, and 
we may hope to control entrance and exit points in some planes perpendicular to $L$, but there seems
to be no topological reason why $E_Y$ cannot also contain some small loops near $0$.
Even worse, the monotonicity of $F$ only allows us to exclude points of $E \sm L$ with a density 
at least $2\pi$, which means that we could also imagine that different curves $E_Y$ let at some
points of type $\bT$, and then topological counting arguments immediately become more complicated.
The point of excluding extra curves in $E_Y$ is that then we could start arguments as in 
Lemma \ref{t29b.3}, where we move curves (or $(n-2)$-spheres) and say that since they never cross
$E_Y$, they keep the same number of intersections with $E$ modulo~$2$. 

\begin{figure}[h]  
\centering
\includegraphics[width=7.cm]{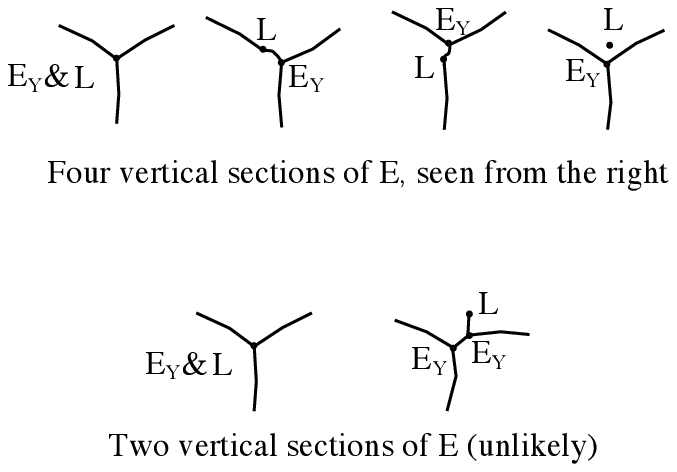}
\caption{Four sections of $G(Y_0)$ and then two less probable sections
of a set $E$  \label{f35.2}}
\end{figure}

But anyway the author expects a mixture of topological and metric information to be needed, to avoid 
various strange cases. A topological possibility is evoked in the second part of Figure \ref{f1.8}, 
which we repeat here as Figure \ref{f35.2}. The top  part  is a warm-up picture of four sections of
the example of $G(Y_0)$ above, and the strange case is depicted in the lower two sections.
On the left, the section that contains the origin does not show any strange behavior, just a set
that looks like a $Y$. And on the right, we see two nodes, that correspond to two  
curves of $E_Y$ that both leave from $L$ at the origin and stay close together, allowing a thin
triangular face (seen as a curve on the picture) that connects one of the curves to $L$,
as it happens near sharp $\bV$-sets, another triangular face that connects the two curves of $E_Y$
(also seen as a small curve on the picture), and as expected three large faces
that make the usual $\2$ angles at the large scale. Some of the faces would have
to turn fast, because of the various $\2$ angles constraints, and again the author does not trust that
this is a viable model.

 %295

Another natural argument does not seem to work easily here. 
In some configurations, when we have a curve $E_Y$ that leaves $L$ and is attached to it by a thin 
triangular piece, it seems tempting to remove the third piece of $E$ (the face in front of the thin face, 
directly across from $L$), and then say that what remains is still an almost 
minimal set, to which we can apply the results of the previous sections. 
Of course it is not certain that the truncated set is still minimal, but more importantly the truncation 
only works as long as $E_Y$ does not do something crazy, like returning to $L$ so that we may need 
a different combination of faces. That is, in strange situations like the one depicted in 
Figure \ref{f35.2}, we may have a hard time defining the face that needs to be removed.
Or it could be that what looks like two different faces near $L$ is actually two pieces of a same face, with a connection somewhere near $L$, and then the truncation is impossible as well.

As we shall see when we deal with other cones $X$, the same issue as for $x\in \bY(L)$ arises when 
$X$ is a minimal cone whose spine contains a half line of $L$. For instance, when $X$ is a cone of type
$\bT$ centered at $0$ and that comes from a tetrahedron with a vertex in $L$. Then the description of
$E$ near the corresponding side of $L$ is essentially as hard to get as when $X \in \bY(L)$.

\section{Approximation by other minimal cones}
\label{S36}

There are many other possible cones $X$ in the setting of Theorems~\ref{t29a.1}-\ref{t29d.1}, 
and the purpose of this section is to say that the local study of $E$ when it is close to one of these cones 
is not very different from what we did (or unfortunately did not manage) in the previous sections.
The general theme is that we can often study $E$ in annuli centered at the origin, and then glue the 
pieces together.

The situation will be simpler if we assume also that 
\begin{equation}\label{36.1}
\text{$0 \in E$ and $X$ is a blow-up limit of $E$ at $0$.}
\end{equation}
We did not always assume this, because in the special cases where we were able to prove
something, we were able to find points $x\in E$ near the origin, with an acceptable density.
For instance, in the situation of Theorem~\ref{t29d.1} (sharp $\bV$-sets), in the simplest case 
$0$ would be a $\bV$-point of $L$, but even if this were not the case, we could still find points
of $\gamma$ very near $0$, and these points have the maximal $F$-density.

Notice also that ``$X$ is a blow-up limit of $E$ at $0$'' has always been the main case where we intended
to apply Theorems~\ref{t29a.1}-\ref{t29d.1}, and we observed earlier that the assumptions for these 
theorems are satisfied as soon as \eqref{36.1} holds, for arbitrarily small radii $R$.

It is often quite convenient to know that
\begin{equation}\label{36.2}
\text{ $X$ satisfies the full length condition,}
\end{equation}
because in this case if $\varepsilon_0$ is small enough and if \eqref{36.1} holds, or even only
\begin{equation}\label{36.3}
\theta(0) := \lim_{r \to 0} r^{-2}\H^2(E\cap B(0,r)) = \H^2(X\cap B(0,1)),
\end{equation}
Theorem \ref{t22.2n} says that $E$ has a unique tangent cone $X_0$ at $0$ 
(if \eqref{36.1} holds, $X = X_0$) and 
\begin{equation}\label{36.4}
d_{0,r}(E,X_0) \leq c_1(\varepsilon_0) \Big(\frac{r}{R}\Big)^{a/4}\ \text{ for } 0 < r \leq R
\end{equation}
as in \eqref{21.8}. This is the typical tool which will allow us to cut $B(0,R/10)$ into annuli where 
hopefully we can control $E$.

Recall that we do not know whether all the minimal cones satisfy the full length property.
At least the planes, the cones of type $\bY$, $\bT$, or $\bV$ do. 
For the products of two orthogonal $Y$-sets of dimension $1$, V. Feuvrier showed that they 
do not satisfy the ``full length property because of angles'',
but it is not known whether they satisfy the full length property itself.

Even if $0 \in E$ and we do not know that $X$ satisfies the full length condition
for one (or equivalently all) of its blow-up limits at $0$, we may still be able to do a decomposition
of $B(0,R/10)$ into annuli where $E$ is well approximated by sliding minimal cones. Indeed,
if $\cX(0)$ denotes the set of blow-up limits of $E$ at $0$, then every element $X$ of $\cX(0)$
is a sliding minimal cone, and (since $\cX(0)$ is actually a compact set)
for each $\varepsilon_0 > 0$ we can find $R_0 > 0$ such that for $0 < R \leq R_0$, we can find
$X \in \cX(0)$ such that $d_{0,R}(E,X) \leq \varepsilon_0$. Thus we are in the situation of
Theorems~\ref{t29a.1}-\ref{t29d.1} for $R$ small, even though in this case we do not know
whether $X$ depends on $R$, and at which speed $d_{0,R}(E,X)$ tends to $0$.

\subsection{Cones $X$  which are fully transverse to $L$}
\label{S36.1}

We shall say that $X$  is fully transverse to $L$ when $K = X \cap \S$ does not meet $L$.
That is, $L \cap X = \{ 0\}$. In this case, $X$ is a plain minimal cone (no sliding condition),
and it has full length as a sliding minimal cone (as in Definition \ref{t3.1} and \eqref{36.2})
if and only if it has full length as a plain minimal cone (as in \cite{C1}).
Set 
\begin{equation}\label{36.5}
A = B(0,R) \sm B(0, 10^{-2}R) \ \text{ and } \, A_1 = B(0,R/2) \sm B(0, 2 \cdot 10^{-2}R).
\end{equation}
Notice that $X \cap A$ is composed of a finite collection of planar faces, bounded by
(radial line segments and) arcs of great circles, that may meet by sets of three and with $\2$ angles,
and it does not get close to $L$.
Under the assumptions of Theorems~\ref{t29a.1}-\ref{t29d.1}, and as usual if $\varepsilon_0$
is small enough, $E \cap A_1$ can be described as a nice $C^1$ version of $X \cap A$, i.e.,
composed of $C^1$ faces that may meet by sets of $3$ with $\2$ angles along $C^1$ curves.
This merely uses the regularity results in the plain case, from \cite{Ta} or \cite{C1}.

Once we have the same assumption in the annuli $2^{-k}A$, and hence the same description in 
the annuli $2^{-k}A_1$, we can glue them. When we do not know whether $X$ satisfies the 
full length property, we still get a reasonably nice description of $E$ near $0$, as a sort of
spiral, not exactly $C^1$, but not far.

If in addition $X$ has full length as in \eqref{36.2}, we can use the extra decay in \eqref{36.4} to apply the 
regularity results from \cite{Ta} and \cite{C1} with constants $\varepsilon_k$ in geometric decay,
get regularity estimates for $E \cap 2^{-k}A_1$ with some geometric decay too, and when we glue we
get $C^{1+a}$ faces (including near $0$) and thus a good $C^{1+a}$ description of $E \cap B(0,R/2)$,
i.e., that near the origin $E$ is the image of $X_0$ by a $C^{1+a}$ diffeomorphism, with uniform estimates.

This is enough for us; the main challenge is to get a full list of minimal cones and establish the full length.
But there is nothing special in this subsection concerning the sliding boundary $L$. We like to say that
$E$ crosses $L$ casually, without even noticing that it exists.

\subsection{Cones $X$ that contain half of $L$ in the interior of a face}
\label{S36.2}

The simplest example of the situation that we describe now is obviously when $X$ is a plane that contains $L$.
But $X$ could also be a $\bY$-set or a $\bT$ set, for which $X \cap \S$ meets $L$, but not at a a vertex.
Or, in higher dimensions, the union of a plane that contains $L$ and some additional transverse stuff
(the simplest being another, nearly orthogonal plane, as in \cite{PUP} or \cite{PUPn}).

We shall first try to describe $E$ under some additional properties;
we shall discuss the other cases later, and in the mean time we will get a simpler discussion.
First assume that $X$ has the full length property (as in \eqref{36.2}); this way we can use \eqref{36.4}, 
get a description in concentric annuli, and deduce from this a global description.
In fact, if the tangent cone $X_0$ at the origin were different from $X$, what we should really study is
the behavior of $E$ near $X_0$, so let us just assume that $X_0 = X$, i.e., $X$ happens to be the
unique tangent cone to $E$ at $x$. This way we get the following simpler form of \eqref{36.4}, i.e.,
\begin{equation}\label{36.6}
d_{0,2r}(E,X) \leq c_1 \big(\frac{r}{R}\big)^{a/4}
\ \text{ for } 0 < r \leq R/2,
\end{equation}
for some $c_1 > 0$ which is as small as we want. 

So assume that $\ell \in L \cap \S$ lies in $K = X \cap \S$, but not as a vertex of $K$.
By the standard description of $K$ near $\ell$ (See Section \ref{S3}), there is 
a small $\eta > 0$, which depends only on the dimension, such that in $B(\ell,40\eta)$, 
$K$ coincides with a geodesic through $\ell$.
Call this geodesic $\cC$, and denote by $P$ the plane that contains $\cC$. Thus 
$L \subset P$, and \eqref{36.6} also says that 
\begin{equation}\label{36.7}
d_{r\ell, 20\eta r}(E,P) \leq (10\eta)^{-1} c_1 \big(\frac{r}{R}\big)^{a/4}.
\end{equation}
This is good, because if $\varepsilon_0$ and hence also $c_1$ are small enough,
$E$ satisfies the assumptions of Theorem \ref{t29c.1} in $B(r\ell,20\eta r)$.
Hence $E$ coincides in $B(r\ell, 2\eta r)$ with a small Lipschitz graph over $P$, with maybe
a crease along part of $L$. 

We have a description like this for every $r \in (0,R/2)$, and now we glue the various descriptions,
and see that in the cone 
\begin{equation} \label{36.8}
T_+(\ell) = \big\{ x\in \R^n \, ; \, \dist(x, L_+(\ell)) \leq \eta |x| \big\},
\end{equation}
where $L_+(\ell) = \big\{ t\ell \, ; \, t \geq 0\big\}$ is the half line that contains $\ell$,
$E$ still coincides with the graph over $P$ of a Lipschitz function $\varphi$ with small Lipschitz constant. 
Again we may have a crease along part of $L$, i.e., a discontinuity of the derivative $D\varphi$ along 
a part $L_g$ of $L$, and we also have the analogue of \eqref{29c.4} on $P \cap T_+(\ell)$, namely that
\begin{equation} \label{36.9}
|D\varphi(x)-D\varphi(y)| \leq \tau |x-y|^b R^{-b} \ \text{ for } x,y \in 
B(0,R/2) \cap T_+(\ell) \cap P \sm L_g 
\end{equation}
with a constant $\tau$ which is as small as we want.
Here we observe that because we have a power decay in \eqref{36.7}, 
and if we follow the proof of of \eqref{29c.4}, we also get \eqref{36.9} with an estimate that
is also valid when $x$ or $y$ approaches $0$ (while staying in $T_+(\ell) \cap P \sm L_g$).
Since $X$ is tangent to $E$ at $0$, observe that $D\varphi(x)$ tends to $0$
(like a power, by \eqref{36.9}) when $x$ tends to $0$.

This was the main point of our description. Let us first assume in addition that 
$X \cap \S$ meets $L$ only once (not at a vertex), as in the examples of a $\bY$-set or a $\bT$ set 
above. We also get a nice description of $E$ in $B(0,R/2)$, obtained by gluing what we have in 
$T_+(\ell)$ and what we can easily deduce from the regularity theorem for plain almost minimal sets.
Let us  say a little more about this. Recall the standard description of $K$; in the present case,
it is composed of a finite collection arcs of circles $\cC_j$, all with length at least $\eta$,
and which connect vertices $\xi\in V$. Here because of our assumption, we only have the special vertex
$\ell$, plus some natural vertices $\xi\in V_1$ where $K$ has a $Y$-shape. 
The points of $V_1$ are at distance at least $10\eta$ from $\ell$ and from each other, and also
(maybe, at the price of redefining $\eta$ and requiring $\varepsilon_0$ to be smaller
if $X$ passes near $-\ell$), the whole $K$ stays at distance at least $10\eta$ from $-\ell$.
For $\xi \in V_1$ and $0 < r \leq R/2$, the regularity theorem for plain almost minimal sets
gives us a description of $E \cap B(r \xi, \eta r)$ as a 
$C^1$ version of a $\bY$-set centered at $\xi r$ and with a nearly radial spine. We can glue all these
results and get a nice description of $E$ in the cone $T_+(\xi)$, again as three faces of class $C^{1+a}$
(including all the way to $0$).

Finally, we can deal with $\xi \in K$ such that $\dist(\xi, \{\ell \} \cup V_1) \geq \eta/10$.
For each such point $\xi$ and $0 < r < R/2$, $X$ coincides with a plane near $B(r \xi, \eta r/20)$, so
we may apply the flat case of the regularity theorem to show that in $B(r \xi, \eta r/100)$,
$E$ is a $C^{1+a}$-version of an almost radial plane, with good H\"older estimates on the 
direction of the tangent plane. We can glue together the local descriptions in the 
$B(r \xi, \eta r/100)$, taking into account that near $0$ we have an extra decay coming from \eqref{36.6},
complete the information with what we already have in $T_+(\ell)$ and the $T_+(\xi)$, $\xi\in V_1$,
and get a description of $E \cap B(0,R/2)$ as a finite union of nice $C^{1+b}$ faces,
with the $\2$ angle condition along $C^{1+b}$ curves $\gamma_\xi$, $\xi \in V_1$
(the set $E_Y$, where the various faces meet), and the creases that were described in $T_+(\ell)$.

We consider ourselves happy with the description above. In the similar case when $-\ell$ also
lies in $K$, and is not a vertex either, we have the same description of $E$ in the opposite cone 
$T_+(-\ell)$ (as a flat surface with creases), and we can glue it to the same other local descriptions
as above, to get a full description of $E$ in $B(0,R/2)$. For the moment the only examples that we have
are when $X$ is a plane (already treated in Section \ref{S29c}) and transverse unions of planes,
where we just get transverse unions of flat surfaces with maybe some creases near $L$.

\ms
For the description above, we left out the case when $E$ is close to $X$ in $B(0,R)$, but maybe not
in smaller balls. If we know that all the blow-up limits of $E$ at $0$ are cones $X$ such 
that no point of $L \cap K$ is a vertex of $K$, we can still proceed as above, except that 
instead of \eqref{3.6} we only get thad $d_{0,r}(E, X(r))$ tends to $0$ (but we don't know the speed)
and $X(r)$ is a minimal cone as above, but that may depend on $r$. In this case, we get nice
descriptions of $E$ in annuli $B(0,r) \sm B(0, r/10)$ that we can glue, but we don't know that
the faces and $E_Y$ become better at $0$, and for instance the faces of $E$ may spiral
near the origin. That is, we get a reasonable description of $E \cap B(0,R/2)$, but in the 
small bilipschitz category rather than $C^{1+b}$.
Notice that this hybrid case can only happen when $X$ does not satisfy the full length
condition. Such cones may exist, but we have no known example.

In the worse case when we do not know about the blow-up limits of $E$ at $0$, we still get a description of
$E$ in the annulus $A_1$, for instance. Again without something like the full length, and especially before
we have a concrete list of minimal cones, it seems a little too abstract to ask for a concrete description.
Notice that when $-\ell \notin K$, the full length should be easier to prove because Section \ref{Sfree} allows
us to use the free attachment in the proofs. At any rate, the global conclusion of this section is that we do 
not really fear the case when $X$ meets $L$, but not at a vertex.

\subsection{Other behaviors of $X$ near $L$, exotic sliding minimal sets}
\label{S36.3}

We continue with the description of $E$ in $B(0,R)$, depending on the behavior of $K=X \cap \S$ near $L$.
We will proceed as before, assume that $X$ satisfies the full length property, so that \eqref{36.4} holds, 
and in addition the blow-up limit of $E$ at $0$ is $X_0 = X$, so that the approximation at smaller 
scales still comes from $X$, as in \eqref{36.6}. Then we take $\ell \in K \cap L$ and try to get a 
description of $E$ in the cone $T_+(\ell)$, depending on the type of $\ell$ in $K$. When $\ell \notin K$,
there is nothing to study, because for $\varepsilon_0$ small enough $T_+(\ell)$ does not meet 
$E \cap B(0,R/2)$.
We studied the case when $\ell \in K$ is not a true vertex, and got the description of the previous subsection.
And there are three more cases that we can treat essentially as before. This is when only one geodesic of $K$
leaves from $\ell$, or when exactly two geodesics leave from $\ell$, either with a generic angle or with a sharp 
$\2$ angle. 

Again we use the fact that thanks to Section \ref{S3}, the length of the geodesics that leave from
$\ell$ is at least $40\eta$, and there is no other geodesic that meets $B(\ell, 40\eta)$. Because of this,
and for every radius $r \in (0,R/2)$, \eqref{36.6} implies an analogue of \eqref{36.7}, but where $P$
is now a half plane bounded by $L$, a generic $\bV$ cone, or a sharp $\bV$ cone.
So, if $\varepsilon_0$ is small enough, $c_1$ is very small too, and we can apply Theorem \ref{t29a.1}, \ref{t29b.1}, or \ref{t29d.1} in $B(r \ell, 20\eta r)$. This gives a good description of $E \cap B(r \ell, 2\eta r)$. Then we glue the various pieces and get a good description of $E$ in $T_+(\ell)$, which actually becomes better
when we approach $0$, because of the extra decay in \eqref{36.6} and \eqref{36.7}.

In the slightly more interesting case when $K$ makes a sharp $\2$ angle at $\ell$, we get that 
there is a small Lipschitz graph $\gamma$ over $L_+(\ell)$, with a part in $L_+(\ell)$ and a 
part in $E_Y$, and then $E \cap B(0,R/2) \cap T_+(\ell)$ consists in two main folds with roughly 
the same direction as the two branches of $X$ near $\ell$, plus maybe a thin vertical piece that 
attaches $E_Y$ and the two folds to $\gamma \cap L_+(\ell)$. 
As before, $E$ has a tangent cone $X$ at the origin, so the description
becomes flat there. The precise description is the same as at the beginning of Section \ref{S29d},
except that we restrict to $T_+(\ell)$. The half plane and generic cases are the same, 
except that there is no curve $\gamma$ and the foils are directly attached to $L_+(\ell)$.

If $-\ell$ also lies in $K$, either as an edge point or a vertex with valence at most two, we have a similar
description of $E$ in $B(0,R/2) \cap T_+(-\ell)$. On the complement of a thinner cone around $L$,
we also have a good description of $E$, as a union of $C^{1+b}$ faces that (maybe) meet along 
$C^{1+b}$ curves of $E_Y$. We glue all  this and get a beautiful description of $E \cap B(0,R)$,
with $C^{1+b}$ faces bounded by $C^{1+b}$ curves of $E_Y$, and maybe thin ``vertical'' faces
and creases along parts of $L$. This was the pleasant part of the subsection.

\ms
So we are left with two issues. The first one is that we do not have an explicit list of sliding minimal cones.
That is, in addition to the planes, unions of strongly transverse planes, and sets of type $\bY$, $\bH$, 
and$\bV$, there are probably other sliding minimal cones, that we'll call exotic. The first candidate for
an exotic minimal cone is the cone over the union of the faces of a cube centered at $0$ and such that $L$
contains a great diagonal of the cube, but there may be many other ones, in particular in higher dimensions.
Another one (in dimensions $4$ and above) is the product $Y \times Y$ of \cite{YxY}. 
At least we know that it is minimal. We don't know whether all these cubes satisfy the full length property, 
and if not we may have situations where we don't know whether $E$ has a unique blow-up limit at $0$ 
(in which case we can still try to give a biLipschitz description of $E$ near $0$, as suggested at the end 
of the previous subsection), or not even (and then we have to distinguish cases depending on the list 
of blow-up limits, but even so we can glue good descriptions in concentric annuli if we can get them). 
This does not seem so bad to the author.

But our main problem arises again, full length property or not, when $\ell \in K \cap L$
is a point of type $\bY(L)$, where three geodesics of $K$ leave from $\ell$ with $\2$ angles.
In this case we don't have an analogue of Theorems \ref{t29a.1}-\ref{t29d.1}, so $E$ may
have an erratic behavior in any of the balls $B(r\ell, \eta r\ell)$. This seems to be the only important
missing piece in our puzzle. See Section \ref{S35} concerning our difficulties in this case.

\section{Verification of full length for the basic minimal cones}
\label{S30}

In this section, we check the full length property for the usual minimal cones.
There would be other cones to study, but we will not do this in the present paper.

\begin{thm}\label{t30.1}
Let $X$ be a cone centered at the origin, and also assume that the line $L$ contains the origin.
The full length property of Section \ref{S3} is satisfied by $X$ when $X$ is a plane, 
a half plane bounded by $L$, a cone of type $\bV$ (two half planes bounded by $L$ and making 
an angle of at most $2\pi/3$ along $L$), a cone of type $\bY$ (regardless of its position 
with respect to $L$), or a cone of type $\bT$ (the cone over the edges of a regular tetrahedron centered at $0$, again regardless of its position with respect to $L$).
\end{thm}

There seems to be lots of cases here, but fortunately the computations were often done earlier.
Recall from Definition \ref{t3.1} that what we need to do is the following. We start from a standard
decomposition of $K = X \cap \S$ (and we can even choose it if we want), and we consider various 
deformations $\varphi_\ast(K)$ of $K$ associated, through some simple rules that may vary a little,
to a mapping $\varphi$ defined on the set of edges of $K$. When the total length of $\varphi_\ast(K)$ 
is larger than $\H^1(K)$, we need to find a sliding competitor for the cone over $\varphi_\ast(K)$ that 
does significantly better than that cone. We will be more specific soon.

\ssi\ub{\bf Case~$0$}.
There is a first case that was already treated in Section 14 of \cite{C1}, which is when $K$ does 
not meet $L$ (and $X$ is a plane or a cone of type $\bY$ or $\bT$). In this case
(we shall call it Case~$0$), the deformation 
of $K$ simply consists, when $K$ is the union of the geodesics $\rho(a_i,b_i)$, in taking 
\begin{equation}\label{30.1}
\varphi_\ast(K) = \bigcup_i \rho(\varphi(a_i),\varphi(b_i)).
\end{equation}
When we do this, we do not need to check anything related to the sliding condition along $L$,
and we can simply import the computations from \cite{C1}. So we will not need to worry about this case,
even though we shall partially redo the case of a cone of type $\bY$. 

When $K$ contains at least one point $\ell \in L$, we cannot do this, because we have to do the
computation alluded to above also for other constructions of $\varphi_\ast(K)$, 
which we call attached, where typically we add the short geodesic $\rho(\ell, \varphi(\ell))$ 
and connect the rest of the geodesics in a way that depends on the number of geodesics of $K$ 
that end at $\ell$. Again, we shall be more specific soon.

\ssi\ub{\bf Case~$1$}.
Nonetheless, there is a second case where we can still rely on the computations of \cite{C1}.
This is when one point of $L \cap \S$ (call it $\ell$) lies in $K$, but not the other one.
We shall call this Case~$1$. Due to the short list of cones $X$ under scrutiny, this happens only
when $X$ is a $\bY$-set and $\ell$ is in the interior of one of the faces of $X$ (Case 1a), 
$X$ is a $\bT$-set and $\ell$ is in the interior of one of the faces of $X$ (Case 1b),
and $X$ is a $\bT$-set and $\ell$ is one of the vertices of $K$ (Case 1c).

The simplest way to get rid of the computation entirely would be to observe that,
thanks to Lemma~\ref{tf.1}, we never need to use the attached option near $\ell$
when $-\ell \notin K$; then we may as well change the definition of full length in this case,
remove the attached option, and in the other case rely on the computations of \cite{C1}
because $L$ no longer plays a role in the computations.

But we announced full length, so let us take a slightly less lazy attitude and yet try to avoid
complicated computations. We only need to do the full length verification with the attached
deformations $\varphi_\ast(K)$ (which are the new ones), but we observe that, 
by the proof of Lemma~\ref{tf.1}, it is enough to find competitors for $\varphi_\ast(X)$, 
not necessarily satisfying the sliding condition, that satisfy \eqref{3.25}. 
Recall that we can do this because, given such a competitor, we can
always modify it as in Section \ref{Sfree} (i.e., by projecting along thin tubes, starting from
the hole near $-\ell$) so that it satisfies the sliding condition and is nearly as good.

\ssi\ub{\bf Subcases~$1a$ and $1b$}.
Let us start with Cases 1a and b, where $\ell$ is connected to two vertices $a_1$ and
$a_2$ of $K$, and near $\ell$ the attached deformation $\varphi_\ast(K)$ coincides
with $\rho_\ast = \rho(\varphi(a_1),\varphi(\ell)) \cup \rho(\varphi(a_2),\varphi(\ell)) 
\cup \rho(\ell, \varphi(\ell))$; see \eqref{3.18}. 
With the free attachement, we would have used the simpler set 
$\rho_\sharp = \rho(\varphi(a_1),\varphi(\ell)) \cup \rho(\varphi(a_2),\varphi(\ell))$ instead.
Call $\varphi_\sharp(K)$ the corresponding deformation of $K$, and $\varphi_\sharp(X)$ the 
corresponding cone. One possibility is that $\varphi(\ell) = \ell$. Then, as sets, 
$\varphi_\ast(K) = \varphi_\sharp(K)$ and $\varphi_\ast(X) = \varphi_\sharp(X)$.
In this case, due to the fact that we are allowed to forget the attachment condition,
\eqref{3.25} for $\varphi_\ast(X)$ is the same as \eqref{3.25} for $\varphi_\sharp(X)$,
which was checked in \cite{C1}. So we may assume that $\varphi(\ell) \neq \ell$.
But then, and again because we no longer need to worry about the sliding condition, it is very easy to 
deform $\varphi_\ast(X)$ in $\ol B(0,1)$ into a subset of $\varphi_\ast(X)$ which coincides with
$\varphi_\sharp(X)$ in $\ol B(0,1/2)$; just contract the  additional geodesic $\rho(\ell, \varphi(\ell))$
along itself, and follow the contraction partially on the cone over $\rho(\ell, \varphi(\ell))$.
When we do this, we cut off a substantial part of the cone over $\rho(\ell, \varphi(\ell))$, and
we save at least $\Delta_1 = \frac18 \H^1(\rho(\ell, \varphi(\ell))) = \frac18 \ddist(\ell,\varphi(\ell))$
in measure.

But since this may not be enough (for instance if the geodesic is very short), 
can even compose this deformation with a deformation of $\varphi_\sharp(X)$,
done in the smaller ball $\ol B(0,1/2)$, into a competitor $\wt X$. We take $\wt X$ from
our verification of \eqref{3.25} for $\varphi_\sharp(X)$, except that we divide the scale by $2$
to allow an easy composition. Of course we only do this when 
$\Delta_\sharp(\varphi) := \H^1(\varphi_\sharp(K)) - \H^1(K) > 0$ (as in \eqref{3.24}),
and then we save an additional measure of
$\Delta_2 = \H^2(\varphi_\sharp(X)\cap B(0,1/2)) - \H^2(\wt X \cap B(0,1/2)) 
\geq \frac{c}{4}\, \Delta_\sharp(\varphi)$, by \eqref{3.25} for $\varphi^\sharp(X)$
(which was checked in \cite{C1}). Altogether, since we may assume that $c < 1/2$,
\begin{eqnarray}\label{30.2}
\H^2(\varphi_\ast(X)\cap B(0,1)) - \H^2(\wt X \cap B(0,1)) 
&\geq& \Delta_1 + \Delta_2
\nn\\
&\geq& \frac18 \ddist(\ell,\varphi(\ell)) + \frac{c}{4} \, \big[\H^1(\varphi_\sharp(K)) - \H^1(K) \big]_+
\nn\\
&\geq& \frac{c}{4} \, \big[\ddist(\ell,\varphi(\ell)) + \H^1(\varphi_\sharp(K)) - \H^1(K) \big]_+
\\
&\geq& \frac{c}{4} \, \big[\H^1(\varphi_\ast(K)) - \H^1(K) \big]_+
\nn
\end{eqnarray}
because $\H^1(\varphi_\ast(K)) \leq \H^1(\varphi_\sharp(K)) + \ddist(\ell,\varphi(\ell))$.
This  proves \eqref{3.25} in our first two cases.

\ssi\ub{\bf Subcase~$1c$}.
In Case 1c, $\ell$ is connected in $K$ to three vertices $a_i$, $1 \leq i \leq 3$, the free deformation 
$\varphi_\sharp(K)$ of $K$ coincides near $\ell$ with 
$\rho_\sharp = \bigcup_{i} \rho(\varphi(\ell),\varphi(a_i))$, while we are interested
in the attached deformation where we select an index, say, $i=1$, and replace
$\rho(\varphi(\ell),\varphi(a_1))$ with $\rho(\varphi(\ell),\ell) \cup \rho(\ell, \varphi(a_i))$.
That is, we force one of the branches $\rho(\varphi(\ell),\varphi(a_1))$ to make a detour through $\ell$.
See \eqref{3.20} and \eqref{3.21}. When $\ell \in \rho(\varphi(\ell),\varphi(a_1))$, 
we did not change anything to the final sets, $\varphi_\ast(K) = \varphi_\sharp(K)$ and 
$\varphi_\ast(X) = \varphi_\sharp(X)$, and \eqref{3.25} for $\varphi_\ast(X)$ is the same as 
\eqref{3.25} for $\varphi_\sharp(X)$, which was checked in \cite{C1}. Otherwise, and since 
Definition \ref{t3.1} allows us to restrict to so-called injective mappings $\varphi$, 
$\rho(\varphi(\ell),\ell) \cup \rho(\ell, \varphi(a_i))$ is disjoint from $\rho_\sharp$, and  
\begin{equation}\label{30.3}
\H^1(\varphi_\ast(K)) - \H^1(\varphi_\sharp(K)) = \ddist(\ell,\varphi(\ell))
+ \ddist(\ell,\varphi(a_1))
- \ddist(\varphi(\ell),\varphi(a_1)).
\end{equation}
As before, we take $\varphi_\ast(X)$ and deform it first, in $\ol B(0,1)$, into a set $\wt X_1$
that coincides with $\varphi_\sharp(X)$ in $B(0,\kappa)$, where $\kappa$ will be chosen soon. 
Recall that we do not need to worry about the sliding condition; we claim that the same sort of 
computations as in Lemma 10.23 in \cite{C1} (see below \eqref{26.2} for a description of the proof) 
yield that
\begin{eqnarray}\label{30.4}
\Delta_1 &:=&\H^2(\varphi_\ast(X) \cap B(0,1)) - \H^2(\wt X_1 \cap B(0,1))
\nn\\
&\geq& C^{-1} \alpha^2 |\varphi(\ell) - \ell| 
\geq C^{-1} \big[\H^1(\varphi_\ast(K)) - \H^1(\varphi_\sharp(K))\big],
\end{eqnarray}
where $\pi - \alpha$  denotes the angle at $\ell$ of the two geodesics $\rho(\varphi(\ell), \ell)$
and $\rho(\ell,\varphi(a_1))$. If we take $\kappa = 0$, the first inequality is a direct consequence of 
Lemma 10.23 in \cite{C1}. Recall that there we use the angle at $\ell$ to push the two faces 
bounded by $\rho(\varphi(\ell), \ell)$ and $\rho(\ell,\varphi(a_1))$ a little bit
in the direction of the cone $H$ over $\rho(\varphi(\ell),\varphi(a_1))$. 
With a slightly more complicated discussion, we could even make sure that in $B(0,\kappa)$,
for some small enough $\kappa$, we push $X$ all the way to that cone $H$. 
Let us proceed another way. First let the first competitor (with $\kappa = 0$, coming from \cite{C1}) 
be as it is (i.e., lying between $H$ and the cone over 
$\rho(\varphi(\ell),\ell) \cup \rho(\ell, \varphi(a_1))$, and push it again inside the smaller 
ball $B(0,C\kappa)$, so that it coincides with $H$ in $B(0,\kappa)$. We then need to show that
we do not increase the measure by more than $C \kappa^2\alpha^2 |\varphi(\ell) - \ell|$
when we do this, because then we take $\kappa$ small and get the first half of \eqref{30.4}.
We proceed as in Lemma 10.23 of \cite{C1}, and in particular follow the error terms, use the
fact that at some point the Jacobian of our soft transformation has an $\alpha^2$ in it (by Pythagorus),
and get the desired estimate. We skip the details because we believe they would just make the 
reader feel sick uselessly.

For the second part of \eqref{30.4}, the strong reader will use elementary geometry in the sphere,
while the author would rely on the size of the derivative of 
$f(z) = \ddist(z,\varphi(a_0)) + \ddist(z,\varphi(\ell))$ along a path from $\ell$
to $\rho(\varphi(\ell),\varphi(a_1))$, and use the computations below \eqref{27.5}.
Again we skip the (now easy) details.

\ms
Once this first modification is done, we use \eqref{3.25} for $\varphi_\sharp(X)$, which was checked 
in \cite{C1}, to construct a second competitor $\wt X$, where we modify $\wt X_1 = \varphi_\sharp(X)$
in $B(0,\kappa)$. This time we save only
$\Delta_2 = \H^2(\varphi_\sharp(X)\cap B(0,1/2)) - \H^2(\wt X \cap B(0,1/2)) 
\geq c \kappa^2 \Delta_\sharp(\varphi)$. We add the two gains as in \eqref{30.2} and get that
\begin{eqnarray}\label{30.5}
&\,& \H^2(\varphi_\ast(X)\cap B(0,1)) - \H^2(\wt X \cap B(0,1)) 
 \geq \Delta_1 + \Delta_2
\nn\\
&\,& \hskip3cm
\geq C^{-1} \big[\H^1(\varphi_\ast(K)) - \H^1(\varphi_\sharp(K))\big] + 
c \kappa^2\,  \big[\H^1(\varphi_\sharp(K)) - \H^1(K) \big]_+
\nn\\
&\,& \hskip3cm \geq \tau \big[\H^1(\varphi_\ast(K)) - \H^1(K) \big]_+,
\end{eqnarray}
with $\tau = \min(C^{-1},c \kappa^2)$ and because $\H^1(\varphi_\ast(K)) \geq \H^1(\varphi_\sharp(K))$.
Again this is \eqref{3.25} for $\varphi_\ast(X)$.

\msi\ub{\bf Case~$2$}.
We may now turn to Case 2, when both points $\ell_\pm$ of $L \cap \S$ lie in $K$.
With the present list of minimal cones $X$, this happens only when $X$ is a $\bV$ set, 
including a plane through $L$, when $X \in \bY(L)$ (a $\bY$-set for which $\ell_+$ lies
on the interior of a face would not do, because then $-\ell \notin K$), and also when $x\in \bT$ and
$\ell_\pm$ lie in the middle of two opposite edges of $K$.

\ssi\ub{\bf Subcase~$2a$}.
We start with the apparently most interesting case when $X \in \bY(L)$. 
First we recall how we decompose $K = X \cap \S$ and define the 
possible deformations $\varphi_\ast(K)$. The set $K$ is composed of three half circles $\cC_i$,
$1 \leq i \leq 3$, we choose a point $w_i \in \cC_i$,  in the middle of 
of $\cC_i$ (because this is allowed and may simplify some computations), denote by $\ell_+$ and $\ell_-$ 
the two points of $L \cap \S$, and also denote by $\cC_{i,\pm}$ the geodesic arc $\rho(w_i,\ell_\pm)$. 
Thus we use a set $V$ of five vertices, and $K$ is the union of the six arcs $\cC_{i,\pm}$ of length $\pi/2$.
Recall that for each verification that we have to do, we are given a mapping $\varphi: V \to \S$,
such that $\sup_{v\in V} |\varphi(v)-v|$ is as small as we want, then we define a set
$\varphi^\ast(K)$ by some set of rules that will be explained soon, denote by $\varphi_\ast(X)$
the cone over $\varphi_\ast(K)$, and, if 
\begin{equation}\label{30.6}
\Delta_L = \H^1(\varphi_\ast(K)) - \H^1(K) > 0
\end{equation}
(as in \eqref{3.24} but with a different name), we need to prove that 
\begin{equation}\label{30.7}
\sigma \geq C^{-1} \Delta_L
\end{equation}
where $\sigma$ is, a little bit as in \eqref{26.1}, the supremum of what we can save in terms 
of $\H^2$-measure when we replace $\varphi_\ast(X)$ with one of its sliding competitors in the unit ball. 
As earlier, we shall use lower bounds on $\sigma$ that come from simple geometric information, 
like the angles between the geodesic arcs that compose $\varphi_\ast(K)$. 
And the constant $C$ in \eqref{30.7} is not allowed to depend on $\varphi$.

We start with the most interesting subcase when $\varphi_\ast(K)$ is attached at both points
$\ell_\pm$, and $\varphi_\ast(K)$ is described near \eqref{3.20}. 
We select for each choice of sign $\pm$ an index $i_\pm \in \{ 1, 2, 3 \}$ and set
\begin{equation}\label{30.8}
\varphi_\ast(K) = K_+^{i_+} \cup K_-^{i_-},
\end{equation}
where for each sign
\begin{equation}\label{30.9}
K_\pm^{i_\pm} = \rho(\varphi(w_{i_\pm}),\ell_\pm) \cup \rho(\ell_\pm, \varphi(\ell_\pm))
\bigcup \cup_{j \neq i_\pm} \rho(\varphi(w_j), \varphi(\ell_\pm)).
\end{equation}
Thus $\varphi_\ast(K)$ is composed of six long geodesics starting from the three $w_i$ and that end
at or near the $\ell_\pm$, plus two short arcs $\rho(\ell_\pm, \varphi(\ell_\pm))$ to connect them.
We allow the degenerate case when $\varphi(\ell_\pm) = \ell_\pm$.

The way we choose the three points $m_i = \varphi(w_i)$ in the arguments, that is, when we choose
the competitors, allows us to take them in the hyperplane $H$ at equal distance from $\ell_+$ and $\ell_-$.
That is, even though in principle the definition would force us to study the case when $m_i \notin H$,
we know that we do not need this case and so we will not study it. Even though this would be possible,
at the price of an additional comparison between such a choice of $m_i$ and the closest choice where
$m_i \in H$. Now set $f_i(z) = \ddist(z, m_i) = \ddist(z, \varphi(w_i))$ for $z\in \S$, and notice that 
\begin{equation}\label{30.10}
\H^1(K_\pm^{i_\pm}) - \frac{3\pi}{2} 
= \ddist(\varphi(\ell_\pm),\ell_\pm) - \pi + \sum_{i \neq i_\pm} f_i(\varphi(\ell_\pm))
\end{equation}
because $\H^1(\rho(m_{i_\pm},\ell_\pm)) = \frac{\pi}{2}$ and, since we may assume that $\varphi$ 
is ``injective'' as in Definition~\ref{t3.1}, the arcs that compose $K_\pm^{i_\pm}$ are disjoint.
We add up these two identities and get that 
\begin{equation}\label{30.11}
\Delta_L = f_+(\varphi(\ell_+)) + f_-(\varphi(\ell_-)) -2\pi, \ \text{ where }
f_\pm(z) := \ddist(z,\ell_\pm) + \sum_{i \neq i_\pm} f_i(z).
\end{equation}
What we want next is an estimate for $\Delta_L$ in terms of various angles and distances, that
we shall then estimate in terms of $\sigma$. For the computation that follows, we fix a sign $\pm$,
drop it from the notation, set $m_0 = \ell_\pm$, and assume, without loss of generality, 
that $i_{\pm} = 3$, so that
\begin{equation}\label{30.12}
f_\pm(z) = f(z) = \sum_{i=0}^2 f_i(z),
\end{equation}
with $f_0(z) = \ddist(z,\ell_\pm) =  \ddist(z,\ell) = \ddist(z,m_0)$. Denote by 
$v_i(z)$ the unit vector pointing in the direction of $\rho(z,m_i)$ (as it leaves from $z$,
and assuming that $z \neq m_0, m_1, m_2$). Then set $s(z) = \sum_{i=0}^2 v_i(z)$,
and for $z = \varphi(\ell)$, set $v_i = v_i(\varphi(\ell))$ and $s = s(\varphi(\ell))$. Finally
define $\delta = \ddist(\ell,\varphi(\ell)) = \H^1(\rho_0)$. We want to estimate 
$f(\varphi(\ell))$ in terms of all these numbers.

Let $\xi$ lie on the geodesic $\rho_\ell = \rho(\ell, \varphi(\ell))$ and let
$w = w(\xi)$ denote the tangent vector to $\rho_\ell$ at $\xi$ pointing in the direction 
of $\varphi(\ell)$; thus $w(\varphi(\ell)) = -v_0$.
The derivative of $f_i$ in the direction $w$ is $- \langle v_i(\xi), w \rangle$ and so the derivative
of $f$ is $- \langle s(\xi), w \rangle$. Thus 
\begin{equation}\label{30.13}
f(\varphi(\ell)) - \pi = f(\varphi(\ell)) - f(\ell)
= - \int_{\rho_\ell} \langle s(z), w(\xi) \rangle \, dl(\xi)
\leq \delta \langle s, v_0 \rangle + 10\delta^2
\end{equation}
because $f(\ell) = \pi$, and with an easy estimate on the variations of $w(\xi)$ and 
the $v_i(\xi)$ along the short geodesic. Thus this tends to be larger when $v_1$ and $v_2$ 
make an angle which is larger than $\2$.
Our next goal is to prove that 
\begin{equation}\label{30.14}
\langle s, v_0 \rangle \delta +10\delta^2 \leq C \sigma
\end{equation}
(precisely the large angle situation alluded to above)
and for this we shall use some lower bounds that were obtained in Sections \ref{S26} and \ref{S27}. 
The point $z = \varphi(\ell)$, and the vectors $v_i$ and $s$ play the same role as there; see near \eqref{26.5}.
Here we are thinking about estimates like \eqref{26.2}, \eqref{26.3}, \eqref{26.14}, and \eqref{26.15}.
The reader may be worried that for these estimates we used a notion of ``good competitors''
(see near \eqref{25.52a}) which is different from our definition of sliding competitor. But here the origin
lies in $L$, the triangle $T(r)$ reduces to the interval $[\ell_-,\ell_+]$, and in this context the new condition \eqref{25.53a} is actually stronger than the usual sliding condition, which means that the good competitors
that we build are also sliding competitors and so we can use the estimates \eqref{26.2}-\eqref{26.15}.
Thus, with the same proof as above, \eqref{26.14} yields
\begin{equation}\label{30.15}
\sigma \geq C^{-1} \delta \langle s, v_0 \rangle \ \text{ when } \langle s, v_0 \rangle \geq |z-\ell|,
\end{equation}
 and by \eqref{26.15}
\begin{equation}\label{30.16}
\sigma \geq C^{-1}  \delta |s|^2.
\end{equation}
We start with the case when $|s| \geq 10^{-10}$, say. Then $\sigma \geq C^{-1} \delta$
by \eqref{30.16}, which is better than \eqref{30.14}. So we may assume that $|s| \leq 10^{-10}$.

Next suppose that $\sigma \leq C_1^{-1} \delta^2$, where $C_1$ will be computed soon. 
For $i = 1, 2$, let $\alpha_i$ denote as near \eqref{26.2} the norm of the sum of the unit directions at 
$m_i = \varphi(w_i)$ of the two geodesics of $\varphi_\ast(K)$ that leave from $m_i$, 
and recall from \eqref{26.2} that $\alpha_i^2 \leq C \sigma \leq C C_1^{-1} \delta^2$. 
Now follow $\varphi_\ast(K)$ when it leaves $z = \varphi(\ell)$, starting
in the direction of $m_i$. When it reaches $m_i$, it turns by a small angle, of size at most 
$2\alpha_i \leq 2 \sqrt{C C_1^{-1}} \delta << \delta$. 
After nearly half a turn, it has deviated by at most $C \alpha_i <<\delta$ from the continuation 
of $\rho(z,m_i)$. 
This last arc (call it $\rho_i$) goes through $-z$, where it arrives with the direction $-v_i$. 
Because $|s| \leq 10^{-10}$, the three vectors $v_i$ make angles with each other that are roughly 
equal to $2\pi/3$, which means that $\rho_i$, when it leaves from $z$, 
meets $\rho(z,\ell)$ transversally, and does not get close to $\ell$. Or equivalently that 
$\rho_i$, when it arrives at $-z$, meets $\rho(-z,-\ell)$ transversally, and does not get close to $-\ell$.
Hence, if $C_1$ is large enough, the continuation of $\varphi_\ast(K)$, which stays so close to $\rho_i$,
does not end at $-\ell$. That is, the two arcs of $\varphi_\ast(K)$ that continue the geodesics 
$\rho(z,m_i)$ meet back at some point $z'$, the next vertex of $\varphi_\ast(K)$,
and which is also the other vertex $z_{\mp}$, because we made sure that we avoid $-\ell$.
Another way to say this is that the index $i_\mp$, associated to $\ell_\mp = - \ell$, is also equal to $3$. 
Now the directions of the three geodesics of $\varphi_\ast(K)$ that leave from $z'$ are 
$C(\alpha_1 + \alpha_2)$-close to $v_1$ and $v_2$, and for the third one the direction of 
$\rho(z',-\ell)$, which is close to $-v_0$ (and it is important that this is just the wrong sign!).
Thus, near the point $z'$, the situation is the following: $s' = s(z')$ is large, 
$\delta' = \ddist(z',-\ell) \geq \delta/2$, and then $\sigma \geq C^{-1} \delta$ by our first
case above. Again this is better than \eqref{30.14}. 

Hence we may assume that $\sigma \geq C_1^{-1} \delta^2$.
Thus, in order to prove \eqref{30.14}, we may assume that $\langle s, v_0 \rangle \geq \delta$,
and now \eqref{30.14} follows from \eqref{30.15}.

We also have the analogue of \eqref{30.14} neat the other point $\ell_\mp$, we sum, we use \eqref{30.13}
and \eqref{30.11}, and we get that $\Delta_L \leq C \sigma$, as needed for \eqref{30.7}.
This completes our proof of full length in the case when $\varphi_\ast(K)$ is attached at both points 
$\ell_\pm$.

\ms
Now consider the case when $\varphi_\ast(K)$ is free near both points $\ell_\pm$.
The set $\varphi_\ast(K)$ is now composed of the six geodesics $\rho_{i,\pm} = \rho(m_i, z_\pm)$,
where $m_i = \varphi(w_i)$ and $z_\pm = \varphi(\ell_\pm)$. For the same reason as before, 
we may assume that the $m_i$ all lie in the hyperplane $H'$ at equal distance between $z_-$ and $-z_-$, 
and this will simplify the computation below a little bit. 
In addition to the three numbers $\alpha_i \geq 0$, that measure how
flat $\varphi_\ast(K)$ is near each $m_i$ (as above), we also define $s_\pm$, the sum of the three unit
directions of the geodesics of $\varphi_\ast(K)$ that leave from $z_\pm$, and the proof of
\eqref{30.16} (or \eqref{26.15}) and \eqref{26.3} now yield
\begin{equation}\label{30.17}
\alpha_i^2 + |s_\pm|^2 \leq C \sigma.
\end{equation}
Consider the geodesic $\rho_{i,-}$ that leaves from $z_-$; at the point $m_i$, it turns by at most
$2\alpha_i$, and then it becomes $\rho_{i,+}$ and ends at $z_+$. We may assume that $\sigma$
is as small as we want, so all the $\alpha_i$ are as small as we want, the three geodesics $\rho_{i,+}$
meet with large angles, and with a little bit of geometry we get that 
$|z_++z_-| \leq C \sum_i \alpha_i \leq C \sqrt\sigma$. We now observe that
since we assumed that $\ddist(m_i,z_-) = \pi/2 = \ddist(m_i,-z_-)$,
\begin{equation}\label{30.18}
\Delta_L = \H^1(\varphi_\ast(K)) - 3 \pi = \sum_{i=1}^3 \ddist(m_i,z_+) - \frac{3\pi}{2}
= f(z_+) - f(-z_-),
\end{equation}
where we set $f(z) = \sum_{i=1}^3 \ddist(m_i,z)$.
Then we estimate the derivative of $f$ along the geodesic $\rho$ from $z_+$ to $-z_-$,
which is bounded by $|\sum_i v_i(z)| \leq |s_\pm| + 10 |z_++z_-| \leq C \sqrt\sigma$.
We integrate along this geodesic and find that $|f(z_+) - f(-z_-)| \leq C \sigma$, which yields
$\Delta_L \leq C \sigma$, as needed.

The case when $\varphi_\ast(K)$ is free near $\ell_-$ and attached near $\ell_+$ is not needed (see
the discussion above where we use Section \ref{Sfree}), but would not really be harder than the two 
previous ones; we would start near $\ell_+$ with the hyperplane $H$ perpendicular to $L$, 
compute the position of the opposite free vertex $z'$ as in the attached case, and end with the two 
computations of $f$ from the two arguments above.
This completes our verification of full length when $X \in \bY(L)$.

\msi\ub{\bf Subcase~$2b$}.
Our next case is when $X \in \bV(L)$ or $X$ is a plane that contains $L$.
One possibility to prove the full length in this case would have been to follow the proof
of Sections~\ref{S26} and \ref{S27} and notice that we can let $d_0$ tend to $0$, but we can
instead follow the proof that was given when $X \in \bY(L)$, and simply remove some branches
from the computation. That is, suppose first that $\varphi_\ast(K)$ is attached at both $\ell_\pm$.
Set $m_i = \varphi(w_i)$ for $i=1,2$ and $z_\pm = \varphi(\ell_\pm)$,
and recall that in the present case $\varphi_\ast(K)$ is composed of the four arcs 
$\rho_{i,\pm} = \rho(m_i, z_\pm)$, plus two short connections $\rho_\pm = \rho(z_\pm,\ell_\pm)$.
Then $\varphi_\ast(K) = K_+ \cup K_-$, with 
$K_\pm = \rho_{1,\pm} \cup \rho_{2,\pm} \cup \rho_\pm$ (as in \eqref{30.8} and \eqref{30.9},
with one less piece each time), $\H^1(K_\pm) = f_\pm(z_\pm)$ with
$f_\pm(z) = \ddist(z,\ell_\pm) + \ddist(z,m_1)+\ddist(z,m_2)$, and then
\begin{equation}\label{30.19}
\Delta_L = \H^1(\varphi_\ast(K)) - \H^1(K) = f_+(z_+)+f_-(z_-) - \pi
\end{equation}
by \eqref{30.6} and as in \eqref{30.11} and \eqref{30.12}, but where we don't even have to 
worry about the two additional arcs leaving directly from $\ell_+$ and $\ell_-$.

As before, we may assume that the two points $m_i =\varphi(w_i)$ lie on the hyperplane $H$
that lies at equal distance from $\ell_+$ and $\ell_-$, prove that 
$f_+(z_+) - \pi \leq \delta \langle s_+, v_0 \rangle + 10\delta^2$ as in \eqref{30.13}
and $\langle s_+, v_0 \rangle \delta +10\delta^2 \leq C \sigma$ as in \eqref{30.14}
(and similarly for $f_-(z_-)$), and conclude from there.

We now stay with the same $X \in \bV(L) \cup \bP(L)$, and assume that we have a free attachment
at both $\ell_\pm$. We still have that $\alpha_1^2 + \alpha_2^2 \leq C \sigma$ as in \eqref{30.17}
or \eqref{26.3}, but we also have numbers $\alpha_\pm = |v_1(z_\pm) + v_2(z_\pm)$
(coming from the angles of the two geodesics $\rho_{i,\pm}$ at $z_\pm$), and the same proof also yields
$\alpha_\pm^2 \leq C \sigma$. Notice that even if by bad luck some $\rho_{i,\pm}$ contains $\ell_\pm$,
we still do not need to check the sliding condition at that point, by definition of a free attachment. 
Here we are in the situation when $\varphi_\ast(K)$ is very close to a great circle, and we may
appeal to computations that were done in \cite{C1}, which yield
\begin{equation}\label{30.20}
\Delta_L = -2\pi + \sum_{i=1}^2 \sum_{\pm} \ddist(m_i,z_\pm)
\leq C (\alpha_1^2 + \alpha_2^2 + \alpha_-^2) \leq C \sigma,
\end{equation}
where we do not even need the last angle $\alpha_+^2$ because it is controlled by the other ones. 
We refer to \cite{C1} for the computation. 

If we have a free attachment at $\ell_-$ but not at $\ell_+$, we still have that 
$\alpha_1^2 + \alpha_2^2 + \alpha_-^2 \leq C \sigma$ as before. In our estimate for $\Delta_L$,
we also have to add the length of $\rho_+ = \rho(z_+,\ell_+)$, so the proof of \eqref{30.20} yields
\begin{equation}\label{30.21}
\Delta_L \leq C (\alpha_1^2 + \alpha_2^2 + \alpha_-^2) + \ddist(z_+,\ell_+)
\leq C \sigma + \ddist(z_+,\ell_+).
\end{equation}
We may assume that the three angles $\alpha_i$ and $\alpha_-$ are small, because otherwise
$\Delta_L \leq 1 \leq C (\alpha_1^2 + \alpha_2^2 + \alpha_-^2) \leq C\sigma$ directly.
Then the two main geodesics $\rho_{1,+}$ and $\rho_{2,+}$ make an angle close to $\pi$
at $z_+$, hence $s_+ = v_0 +v_1+v_2$ (the sum of the three directions of the geodesics that
leave from $z_+$) is large because $v_1+v_2$ is small. That is, $|s_+| \geq 1/2$. 
Then $\sigma \geq C^{-1} \delta = C^{-1}\ddist(z_+,\ell_+)$
by the analogue of \eqref{30.16} or \eqref{26.15}, and $\Delta_L \leq C \sigma$ as needed.

\begin{figure}[!h]  
\centering
\hskip0.2cm\includegraphics[width=5.8cm]{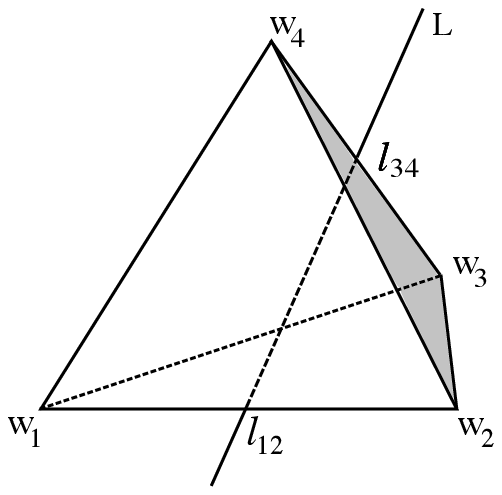}
\hskip2.6cm
\includegraphics[width=7cm]{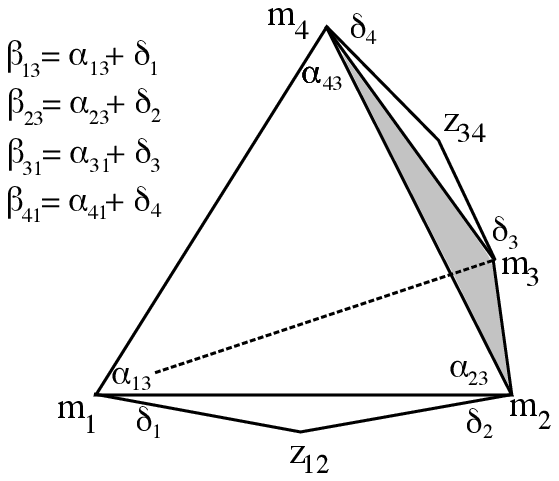}
\caption{The situation of Subcase $2c$ ($T$ on the left, $T_\varphi$ on the right)
\label{F30.A}}
\end{figure}

\msi\ub{\bf Subcase~$2c$}.
We are left with only one possibility in Case 2, when $X \in \bT$ and $L$ goes through
opposite points at the middle of two opposite edges of $X$. See the right part of 
Figure \ref{F30.A}. 
The general plan is, as we did in \cite{C1}, to use the angles of the deformed tetrahedron
to control the lengths and then $\Delta_L$.

Some general notation will be useful. Denote by $w_1, w_2, w_3, w_4 \in \S$ the four edges of the 
unit tetrahedron $T$ that defines $X$. Set $m_i = \varphi(w_i) \in \S$. We are also
interested in the tetrahedron $T_\varphi$ with vertices $m_i$.
We may label the points so that the two points of $L \cap \S$ are
$\ell_{12} = (w_1+w_2)/2$ and $\ell_{34} = (w_3+w_4)/2$. Then set
$z_{12} = \varphi(\ell_{12})$ and $z_{34} = \varphi(\ell_{34})$. See the left part of 
Figure \ref{F30.A}. 
Some times it will not matter where they are relative to $L$; the relative position
of these points with respect to $T_\varphi$ will be more important.

What will be controlled best is angles at the vertices $m_i$. First consider the angles
of the tetrahedron $T_\varphi$; denote by $\alpha_{ij}$ the angle at $m_i$
of the face of $T_\varphi$ that does not contain $m_j$. With the notation of \eqref{9.12},
$\alpha_{12} = \Angle_{m_1}(m_3,m_4)$, and similarly for the other ones. That is, 
\begin{equation}\label{30.22}
\alpha_{ij} = \Angle_{m_i}(m_k,m_l), \ \text{ where $i,j,k,l$ are different.}
\end{equation}
Again see Figure \ref{F30.A}, on the right. 
All these angles are close to $\2$.

We are also interested in the angle $\beta_{ij}$, which is close to $\alpha_{ij}$ and obtained 
from $\alpha_{ij}$ as follows. if the geodesic $\rho(m_i,m_k)$ is one of the two special geodesics 
$\rho(m_1,m_2)$ and $\rho(m_3,m_4)$, we replace $m_k$ with the vertex 
$z_{ik} = \varphi(\ell_{ik})$ that is about half way; for this notation to work fine we also set 
$z_{kj}=z_{jk}$ when $j<k$. Thus for instance 
$\beta_{13} = \Angle_{m_1}(z_{12},m_4)$, $\beta_{23} = \Angle_{m_2}(z_{12},m_4)$,
and $\beta_{43}= \alpha_{43} = \Angle_{m_4}(m_1, m_2)$, to name the three angles relative 
to the front face of Figure \ref{F30.A}. 

Finally, for $1 \leq i \leq 4$ we define the sum $s_i$ of the three unit vectors that are used 
to define the three $\beta_{ij}$. That is,
$s_1 = v(m_1,z_{12}) + v(m_1,m_3)+v(m_1,m_4)$,
$s_2 = v(m_2,z_{12}) + v(m_2,m_3)+v(m_2,m_4)$,
$s_3 = v(m_3,z_{34}) + v(m_3,m_1)+v(m_3,m_2)$, and
$s_4 = v(m_4,z_{34}) + v(m_4,m_1)+v(m_4,m_2)$.

Next we explain about $\varphi_\ast(K)$ and how we control the angles.
Except for the two exceptional geodesics $\rho(w_1,w_2)$ and $\rho(w_3,w_4)$,
we just replace $\rho(w_i,w_j)$ in $K$ with $\rho(m_i,m_j)$ and get the corresponding piece
of $\varphi_\ast(K)$. Then let us explain for $\rho(w_1,w_2)$; the case of $\rho(w_3,w_4)$
is similar. If we have a free attachment, we just replace $\rho(w_1,w_2)$ with 
$\rho(m_1,z_{12}) \cup \rho(z_{12},m_2)$. If instead we have an attached configuration, we
also add the short geodesic $\rho(z_{12},\ell_{12})$. We do this for both points $\ell$,
take the union, and get $\varphi_\ast(K)$. 

Notice that near the half line through $m_i$, $\varphi_\ast(X)$ (the cone over $\varphi_\ast(K)$)
is a $\bY$-set with angles $\beta_{ij}$, and $|s_i|$ is the same number as $\alpha_\varphi(w_i)$
in (10.20) of \cite{C1}; thus $\sigma \geq C^{-1} |s_i|^2$ by Lemma 10.23 there.
See the discussion below \eqref{26.4}, too. Let us record this, i.e.,
\begin{equation} \label{30.23}
\sum_{i=1}^4 |s_i|^2 \leq C \sigma,
\end{equation}
and now try to control the geometry with this information. Later on, we will take care of the 
short geodesics of attached configuration. 

\ssi\ub{\bf Subcase~$2c1$}.
We start with a case that is easier to understand, when $z_{12}$ and $z_{34}$ both lie in 
the $3$-space that contains the four $m_i$. Then we can rely a little more safely on 
Figure~\ref{F30.A}, 
which we rather see as a picture on the unit sphere, with straight lines replaced with geodesics.
Let $\delta_1$ denote the angle $\Angle_{m_1}(m_2,z_{12})$, counted positive if
$z_{12}$ lies outside of the face $(m_1,m_2,m_4)$ as in the picture. 
Define $\delta_2 = \Angle_{m_2}(m_1,z_{12})$, with the same sign convention,
and then $\delta_3 = \Angle_{m_3}(m_4,z_{34})$ and $\delta_4 = \Angle_{m_4}(m_3,z_{34})$,
counted positively when $z_{34}$ lies out of the face $(m_2,m_3,m_4)$, as in the picture.
Since $\varphi$ moves the points very little, $\delta_1/\delta_2$ and $\delta_3/\delta_4$
are as close to $1$ as we want.

A second advantage of the fact that $z_{12}$ lies in the $3$-space that 
contains the $m_i$ is that the three unit vectors whose sum is $s_1$ lie in a same plane
(the tangent plane to $\S$ at $m_1$ in that $3$-space). 
Then the fact that $|s_1|^2 \leq C \sigma$ implies that $|\beta_{1,j}-\2|^2 \leq C \sigma$
for $j=2,3,4$. In fact, the same argument works at every vertex, and we get that for all $i,j$,
\begin{equation} \label{30.24}
\big|\beta_{i,j}-\2 \big|^2 \leq C \sigma.
\end{equation}
But we prefer to have a similar control the simpler angles $\alpha_{i,j}$, and for this we want
to show that
\begin{equation} \label{30.25}
\delta_j \leq C \sqrt\sigma \ \text{ for } 1 \leq j \leq 4.
\end{equation}
We need a bit of spherical geometry. Set $L_{ij} = \sin\ddist(m_i,m_j)$
and concentrate on the spherical triangle $(m_1,m_2,m_4)$ in front of the picture.
By 18.6.13.4 in \cite{Berger}, % and 18.6.13.7 
\begin{equation} \label{30.26}
\frac{L_{14}}{\sin\alpha_{23}} = \frac{L_{12}}{\sin\alpha_{43}} = \frac{L_{42}}{\sin\alpha_{13}}.
\end{equation}
In terms of $\beta_{i,j}$, 
\begin{equation} \label{30.27}
\frac{L_{14}}{\sin(\beta_{23}-\delta_2)} = \frac{L_{12}}{\sin\beta_{43}} 
= \frac{L_{42}}{\sin(\beta_{13}-\delta_1)} \, .
\end{equation}
Then we do the same thing with the bottom face $(m_1,m_2,m_3)$, and we get that
\begin{equation} \label{30.28}
\frac{L_{13}}{\sin(\beta_{24}+\delta_2)} = \frac{L_{12}}{\sin\beta_{34}} 
= \frac{L_{32}}{\sin(\beta_{14}+\delta_1)} \, ,
\end{equation}
where we exchanged $3$ and $4$, and changed signs because now $z_{12}$ lies in the face.
The same computation in the face $(m_3,m_2,m_4)$ yields
\begin{equation} \label{30.29}
\frac{L_{24}}{\sin(\beta_{31}-\delta_3)} = \frac{L_{34}}{\sin\beta_{21}} 
= \frac{L_{23}}{\sin(\beta_{41}-\delta_4)} 
\end{equation}
and in the face $(m_3,m_1,m_4)$ (back there; exchange $m_1$ and $m_2$)
\begin{equation} \label{30.30}
\frac{L_{14}}{\sin(\beta_{32}+\delta_3)} = \frac{L_{34}}{\sin\beta_{12}} 
= \frac{L_{13}}{\sin(\beta_{42}+\delta_4)}.
\end{equation}
Let us assume that $\delta_1$ and $\delta_2$ are positive (otherwise, 
we could just exchange the names of two faces), and also that $\delta_3$ and $\delta_24$ 
are positive (we'll do the other case later).

We use \eqref{30.24} to estimate the $\sin(\beta_{ij})$.
A very simple computation that gives a hint of what will happen next would be to 
assume that all the $\beta_{ij}$ are equal to $\2$, and use \eqref{30.27}-\eqref{30.30} 
to find that if both $\delta_i$ are positive, $L_{14} < L_{12}$ by \eqref{30.27}, 
$L_{12} < L_{32}$ by \eqref{30.28}, $L_{32}= L_{23} < L_{34}$ by \eqref{30.29}, and
$L_{34} < L_{14}$ by \eqref{30.30}, a contradiction.

Here we take into account small errors that come from \eqref{30.24}.
We deduce from the first part of \eqref{30.27} that
\begin{equation} \label{30.31}
L_{14} = \frac{\sin(\beta_{23}-\delta_2)}{\sin\beta_{43}} L_{12}
\leq (1-\delta_2/10)(1+C\sqrt\sigma) L_{12}.
\end{equation}
Then by the second part of \eqref{30.28}
\begin{equation} \label{30.32}
L_{32} = \frac{\sin(\beta_{14}+\delta_1))}{\sin\beta_{34}} L_{12}
\geq (1+\delta_1/10)(1-C\sqrt\sigma) L_{12}.
\end{equation}
We compare, use the fact that $\delta_1/\delta_2$ is close to $1$, and get that
\begin{equation} \label{30.33}
L_{14} \leq (1-\delta_1/10)(1+C\sqrt\sigma) L_{32}.
\end{equation}
Now \eqref{30.29} yields
\begin{equation} \label{30.34}
L_{32} = L_{23}= \frac{\sin(\beta_{41}-\delta_4)}{\sin\beta_{21}} L_{34}
\leq (1-\delta_4/10)(1+C\sqrt\sigma) L_{34}
\end{equation}
while by \eqref{30.30}
\begin{equation} \label{30.35}
L_{34} = \frac{\sin\beta_{12}}{\sin(\beta_{32}+\delta_3)} L_{14}
\leq (1-\delta_3/10)(1+C\sqrt\sigma) L_{14}
\end{equation}
and, when we compare and use the fact that $\delta_3/\delta_4$ is close to $1$, 
\begin{equation} \label{30.36}
L_{32} \leq (1-\delta_3/10)(1+C\sqrt\sigma) L_{14}.
\end{equation}
This is only compatible with \eqref{30.33} when $\delta_1 + \delta_3 \leq C \sqrt\sigma$;
\eqref{30.25} follows.

We still need to check \eqref{30.25} when $\delta_3$ and $\delta_4$ are negative.
The simple sketchy computation modulo errors from \eqref{30.24} is now that 
$L_{42} \leq L_{12}$ by \eqref{30.27}, $L_{12} \leq L_{13}$ by \eqref{30.28}, 
$L_{13} \leq L_{34}$ by \eqref{30.30}, and $L_{34} \leq L_{42}= L_{24}$ by \eqref{30.29}. 
The details with the errors are the same as for \eqref{30.33} and \eqref{30.36}.

At this point we know that \eqref{30.25} holds and all the angles $\alpha_{ij}$ and 
$\beta_{ij}$ are $C\sqrt\sigma$-close to $\2$. We now try to reconstruct $T_\varphi$
from this information. By \eqref{30.26}-\eqref{30.30}, there is a number 
$L$ such that $|L_{ij}-L| \leq C \sqrt\sigma$ for $i \neq j$. 
Recall that $L_{ij} = \sin\ddist(m_i,m_j)$; let us set $l_{ij} = \ddist(m_i,m_j)$
(so $L_{ij} = \sin l_{ij}$) and $x_{ij} = \cos l_{ij}$. 
Since $\varphi$ does not move the points too much) the $L_{ij}$, and hence also $L$,
are close to $\sin l_0$, where $l_0$ is the common value of the distances $\ddist(w_i,w_j)$ in $T$.
Then, inverting the sine function locally near $l_0$, $|l_{ij}-l| \leq C \sqrt\sigma$, for the number $l$ such that $\sin l = L$ which lies
close to $l_0$, and also $|x_{ij}-x| \leq C \sqrt\sigma$, where $x = \cos l$.

We want to evaluate $L, l$, and $x$ more precisely. We use 18.6.13.7 in \cite{Berger}, 
which says that in a spherical triangle with angles $\alpha_1, \alpha_2, \alpha_3$ 
and opposite lengths $l_1, l_3, l_3$, 
\begin{equation}\label{30.37}
\cos\alpha_1 = \frac{\cos l_1-\cos l_2 \cos l_3}{\sin l_2 \sin l_3}.
\end{equation}
We can use this to compute $l_0$, because for the equilateral triangles that compose $T$,
$\cos \alpha_1 = \cos\2 = -1/2$, and the right-hand side is 
$\frac{\cos l_0-\cos^2 l_0}{\sin^2 l_0}$, so that $x_0 = \cos l_0$ is a solution of
$x-x^2 = -\frac12 (1-x^2)$, or $3x^2 -2x -1 = 0$. The solutions are $x=-\frac13$
and $x=1$ (which we exclude), and so $x_0 = \cos l_0 = -1/3$. 
A confirmation is that if $w_1 = (1,0,0)$, the common first coordinate of the other 
vertices $w_i$ is $-1/3$, because $\sum w_i = 0$.

Return to \eqref{30.37}, which we now apply to one of the triangles that compose $T_\varphi$.
Set $x_i = \cos l_i$; then the right-hand side is $(x_1-x_2 x_3)(1-x_2^2)^{-1/2} (1-x_3^2)^{-1/2}$.
Since $|x_{i}-x| \leq C \sqrt\sigma$ for $i=1,2,3$, the partial derivatives of the expression are less
than $100$ in the region of interest near $x_0 = -1/3$, and the left-hand side of \eqref{30.37} is 
$C \sqrt\sigma$-close to $\cos\2 = -1/2$,
we see that $\big|(x-x^2)(1-x^2)^{-1} + \frac12\big| \leq C \sqrt\sigma$.
Or (multiplying by $1-x^2$), we can find $\tau \in [-C \sqrt\sigma, C\sqrt\sigma]$ such that
$x-x^2 + \frac{1+\tau}{2}(1-x^2) = 0$. We expand, solve by radicals, keep the root that lies close to $x_0$,
and get that $|x-x_0| \leq C \sqrt\sigma$. Then we take the cosine and get that 
$|l - l_0| \leq C \sqrt\sigma$.
Returning to the triangles of $T_\varphi$, we see that $|\ddist(m_i,m_j)-l_0| = |l_{ij}-l_0|
\leq C\sqrt\sigma$ for all the distances.

We may conclude our initial length computation as we did in the last pages of \cite{C1}.
We may reconstruct $T_\varphi$ (modulo an isometry), from the length $l_{ij}$ and the angles 
$\alpha_{ij}$, with errors less than $C\sqrt\sigma$.
Then we use the fact that $T$ is a critical point of the sum of lengths to prove that
\begin{equation} \label{30.38}
\sum_{i < j} \ddist(m_i,m_j) \leq \sum_{i < j} \ddist(w_i,w_j) + C \sigma
= 6 l_0 + C \sigma.
\end{equation}
For the total length of the geodesics of $\varphi_\ast(K)$, we also need to add
the extra length  
\begin{eqnarray}\label{30.39}
&\,& \hskip -1cm [\ddist(m_1,z_{12})+\ddist(z_{12},m_2) - \ddist(m_1,m_2)]
+ [\ddist(m_3,z_{34})+\ddist(z_{34},m_4) - \ddist(m_3,m_4)]
\nn\\
&\,& \hskip 4cm\leq C \delta_1^2 + C \delta_3^2 \leq C \sigma.
\end{eqnarray}
In the simpler case where we have the free attachment at both points $\ell_{12}$ and $\ell_{34}$,
there is no other term and we get that 
\begin{equation} \label{30.40}
\Delta_L = \H^1(\varphi_\ast(K)) - 6 l_0 \leq C \sigma, 
\end{equation}
as needed for the full length property.

We are left with the case where $\varphi_\ast(K)$ is attached, and we have to add the length
$\ddist(\ell_{12},z_{12})$ or $\ddist(\ell_{34},z_{34})$ (or both) to get $\H^1(\varphi_\ast(K))$.
Suppose for instance that $\varphi_\ast(K)$ is attached near $\ell_{12}$ and we have to add 
$\ddist(\ell_{12},z_{12}) > 0$. Near the segment $[0,z_{12}]$, $\varphi_\ast(X)$ is composed
of two large faces that leave from $[0,z_{12}]$ in almost opposite directions, plus a thin face 
(the cone over $\rho(z_{12},\ell_{12})$) which is disjoint from the rest because we assumed
that $\varphi$ is ``injective,'' so $\rho(z_{12},\ell_{12})$ does not meet the other geodesics
$\rho(z_{12},m_1)$ and $\rho(z_{12},m_2)$. But then, by a minor variant of 
\eqref{26.14} (notice that in this case, since $v_1+v_2$ in \eqref{26.4} is small, $s$ is close to $v_0$),
we get that $\sigma \geq C^{-1} |z_{12}-\ell_{12}|$, so the additional term $\ddist(\ell_{12},z_{12})$
is controlled. The same estimate would hold near $\ell_{34}$, and we still get the conclusion of 
\eqref{30.40}. This  completes our proof of full length in our Subcase~$2c1$ where  
$z_{12}$ and $z_{34}$ both lie in the $3$-space that contains the four $m_i$.

\ssi\ub{\bf Subcase~$2c2$}.
Suppose this is not the case; we want to use the previous computation, so we denote by
$V$ the $3$-space that contains the $m_i$, set $\S_V = \S \cap V$, and denote by
$\wt z_{12}$ and $\wt z_{34}$ the closest point projection of $z_{12}$ and $z_{34}$
on $\S_V$. We start with the case when
\begin{equation}\label{30.41}
|\wt z_{12}-z_{12}| + |\wt z_{34}-z_{34}| \leq C_0\sqrt\sigma,
\end{equation}
where $C_0$ will be chosen later. Even though $C_0$ may be large, the estimates that follow
hold because $\varphi$ does not move the points much (hence $|\wt z_{12}-z_{12}|$ is very small,
for instance). Since $|\wt z_{12}-z_{12}|$ is minimal,
the direction $v(\wt z_{12},z_{12})$ is orthogonal to $v(\wt z_{12},m_1)$ and 
$v(\wt z_{12},m_1)$ (that lie in the tangent space of $\S_V$), 
so the derivative of $\ddist(z, m_i)$ in the direction of $v(\wt z_{12},z_{12})$
vanishes at $\wt z_{12}$ (for $i=1,2$),
and integrating this derivative on $\rho(\wt z_{12},z_{12})$ yields
\begin{equation}\label{30.42}
\ddist(z_{12}, m_1)+\ddist(z_{12}, m_2) 
\leq \ddist(\wt z_{12}, m_1)+\ddist(\wt z_{12}, m_2) + C \sigma.
\end{equation}
Of course we have the same estimate for $z_{34}$ and the index $34$.
Our estimates \eqref{30.24} for the angles remains valid when we replace 
$z_{12}$ and $z_{34}$ by $\wt z_{12}$ and $\wt z_{34}$, by \eqref{30.41}.
Then we can follow the computations that we did above, 
applied with $\wt z_{12}$ and $\wt z_{34}$, and we get that \eqref{30.38} holds, and 
also \eqref{30.39} with $z_{12}$ and $z_{34}$ replaced by $\wt z_{12}$ and $\wt z_{34}$.
We add this with \eqref{30.42} and its analogue for $34$, and get \eqref{30.40} for the case when
we do not have additional pieces $\ddist(\ell_{12},z_{12})$ or $\ddist(\ell_{34},z_{34})$
to worry about.
Finally the estimate for adding these pieces is the same as in Subcase~$2c1$.

So we may now assume that \eqref{30.41} fails, and for instance 
$|\wt z_{12}-z_{12}| \geq \frac{C_0}{2} \sqrt\sigma$.
Notice that if $\ol z$ denotes the projection of $z_{12}$ on $V$, then
$|\wt z_{12}-z_{12}| \leq |\ol z-z_{12}| + \dist(\ol z, \S_V) \leq |\ol z-z_{12}|+||\ol z|-1|
\leq 2 |\ol z-z_{12}|$, so $|\ol z-z_{12}| \geq \frac{C_0}{4} \sqrt\sigma$.
Now let $s$ denote the sum of the unit vectors $v(m_1,m_3)$, $v(m_1,m_4)$ and 
$v(m_1,z_{12})$ that describe the three faces of $\varphi_\ast(X)$ near $m_1$.
We just keep the coordinate $s^\perp$ of $s$ along $V^\perp$; we get that
\begin{equation}\label{30.43}
|s| \geq |s^\perp| =  v(m_1,z_{12})^\perp \geq \frac12 \dist(z_{12},V)
= \frac12 |\ol z-z_{12}| \geq \frac{C_0}{8} \sqrt\sigma.
\end{equation}
This is good, because Lemma 10.23 in \cite{C1} (or if you prefer the discussion near \eqref{26.4})
yields $|s|^2 \leq C\sigma$, a contradiction if $C_0$ is chosen large enough. 
So \eqref{30.41} was our only case, this completes the verification of full length in our last
Subcase~$2c2$, and Theorem \ref{t30.1} follow.
\qed

\section{Extension to curvy boundaries}
\label{S31}

The main theorems of this paper are still valid when the boundary $L$ is a curve of class $C^{1+b}$,
for some $b > 0$, rather than a line as in the previous sections. The proof consists 
in checking that all our arguments still work, with rather minor modifications, but let us be a little more specific here.

One of the main engines of our proofs is the use of near monotonicity formulae, provided by
\cite{Sliding} and \cite{Mono}. In the first cases, for balls centered on $L$, the relevant extension
is presented in Remark 28.11 and Theorem 28.15 of \cite{Sliding}. Our present assumption that
$L$ is of class $C^{1+b}$ is stronger than the sufficient condition given there; it would be enough to
assume that $L$ is a Lipschitz curve, say, and that, if we want to prove the near monotonicity
of $\theta(r) = r^{-2} B(0,r)$) for $0< r \leq r_0$, we can find a bilipschitz mapping 
$\xi : B(0,2r_0) \to \R^n$, that maps $0$ to itself and $L \cap B(0,2r_0)$ to a line,
and such that the restriction of $\xi$ to $B(0,r)$ is $(1+\rho(r))$-bilipschitz, with for instance
$\rho(r) \leq c r^{b}$ for some $b > 0$ and $c > 0$ small enough.

We also need the near monotonicity of (a minor variant of) the function $F$ of \eqref{22.3},
for balls that are centered slightly outside of $L$. Here again, there is a statement when $L$
is a curve of class $C^{1+b}$, which is given in Theorem 7.1 on page 380 of \cite{Mono}; notice also
that for $r$ small enough (so that $L$ is flat enough in the given balls), Remark 7.3 on 
page 383 of \cite{Mono} says that the added term in the functional is still the normalized 
measure of the shade of $L$.

Associated to the near monotonicity of $\theta$ or $F$ is the fact that we control the geometry
of $E$ in balls where it is almost constant. See Proposition 30.19 in \cite{Sliding} and Theorem 9.1 in \cite{Mono}.

Another fact that is used a lot in our constructions is the local regularity theory of $E$ far from the
boundary, that we import from \cite{C1}, and that we use to control $E$ far from $L$ and in particular
construct the curves in $E \cap \S_r$ that lead to competitors. For this, the precise shape of $L$ does 
not matter.

And finally, there is the main construction of competitors, where one starts with a fixed origin,
almost any small radius $r$, and one constructs curves, and finally competitors that we compare with
$E$ to get differential inequalities that eventually lead to a decay of $\theta$ or $F$, and also the geometric
control on the approximation of $E$ (as in Sections \ref{S18}-\ref{S20} and \ref{S28}).
For this, the simplest way seems to proceed as in the proof of Theorem 28.15 of \cite{Sliding},
which consists in using a bilipschitz change of variable $\xi$ as above to transform the pair $(E,L)$
into another pair $(\wt E, \wt L)$ for which $\wt L$ is a line. We may then construct the same competitor
for $\wt E$ as above, using in particular the local regularity of $E$ (or equivalently $\wt E$) far from
$L$ (or $\wt L$). The estimate on the competitors $\wt F$ for $\wt E$ that we construct then yield
the desired estimates for corresponding competitors $F$ for $E$; the main point is that in $B(0,r)$,
if the mapping $\xi$ is $(1+\rho(r))$-bilipschitz, the errors on the $\H^2$-measure of competitors
are less than $C \rho(r) r^2$, which is of the same order as the other error terms that we had already.
We skip the details of the computations, which are very similar to what was done in \cite{Sliding}.

For this part of the argument, we could formalize what we are doing, by defining the notion, for a set $E$
which is already known to be quasiminimal with sliding boundary $L$, of being almost-minimal, at the point
$0$ (say, and $0$ does not need to be in $L$), with the same boundary $L$ and a given function $h$.
This just means that when $F = \varphi_1(E)$ is a sliding competitor for $E$ in $B(0,r)$, as in
Definition \ref{t1.1}, then we have \eqref{1.11}. Thus the difference is that we only consider balls
centered at the origin.

Then we may observe that if $E$ is almost minimal at the point $0$ (with gauge function $h$),
and $\xi$ is a bilipschitz, and asymptotically optimally bilipschitz as above, then $\wt E = \xi(E)$
is also an almost minimal set at the origin, with a gauge function $\wt h$ such that
$\wt h(r) \leq 2 h(2r) + C \rho(2r)$, where $C$ also depends on the local Ahlfors regularity constant
for $E$ near $0$ (which exists because we assumed $E$ to be quasiminimal), and it is enough to take
$2r$ in the argument because we assume $\xi$ to be bilipschitz with a constant that is close to $1$.

This part is easy to check, just using the definitions and the local Ahlfors regularity of $E$ and $\wt E$.
Then the main point of the proof that follows is that the main decay and approximation results of this paper
are still true if we only assume that $E$ is almost minimal at the point $0$, with $h(r) \leq C r^b$
as in the previous sections, and in addition $E$ is quasiminimal near $0$ (so that it is locally Ahlfors-regular), 
and satisfies the regularity estimates of \cite{C1} far from $L$. For this part, we would just need to read
our proof again, checking where each estimate comes from.

The author would have preferred to say that it is enough to use the pointwise almost minimality of $E$
and its quasiminimality (but not the regularity results far from $L$), but this does not seem to 
be the case, or at least the proof above does not say this (because we often use 
$C^{1+a}$-regularity, possibly often for convenience). This is why we decided not to insist 
so much on the notion of almost minimality at a point, even though it would be very convenient
during the proof.

So the conclusion of this section is that the generalization of our theorems to smooth
boundaries $L$ is rather straightforward, but tedious and boring enough for us to skip the proof.

\section{Sets of $\bH(L) \cup \bV(L)$ are sliding minimal cones}
\label{S32}

In this section we prove two sliding minimality results that were apparently 
not written down yet. They are not needed for this paper, but of course the results 
above make more sense because they are true. In both cases  the simplest version
is when $L$ is a line through the origin, but we also included some larger boundaries because it is 
not much harder, and the dimension of the ambient space $\R^n$ does not matter.

\begin{lem}\label{t32.1}
Let $H \in \bH$, denote by $\d H$ its boundary (a line), suppose that $0 \in \d H$,
and choose an orthonormal basis $(e_1, \ldots, e_n)$ of $\R^n$ so that $\d H = \R e_1$ and $e_2 \in H$.
Let a boundary set $L$ be given, such that $\d H \subset L \subset e_2^\perp$.
Then $H$ is a sliding minimal set in $\R^n$, associated to the sliding boundary $L$.
\end{lem}

Let $E = \varphi_1(H)$ be a sliding competitor for $H$ in some ball $B$,
coming from a one parameter family $\{ \varphi_t \}$ as in Definition \ref{t1.1}; 
we want to show that 
\begin{equation} \label{32.1}
\H^2(H \cap B) \leq \H^2(E \cap B).
\end{equation}
In fact, we shall not need to know that $\varphi_1$ is Lipschitz, as in \eqref{1.8},
so the lemma also gives the \ub{slightly stronger minimality} property where \eqref{1.8}
is not required. 
Also, we may assume that $B$ is centered on $0$ (otherwise we replace it with a larger ball). 
 
We use the orthonormal basis above to write coordinates; 
vectors of $\R^n$ will be denoted by $w = (x,y,z)$, with $x, y\in \R$
and $z \in \R^{n-2}$. Denote by $\pi$ the closest distance projection on $H$, 
defined by $\pi(x,y,z)= (x, y_+,0)$, where $y_+ = \max(y,0)$. 
We define a new deformation $\{ \varphi_t^\ast \}$ by $\varphi_t^\ast(w) = \pi(\varphi_t(w))$. 
The main constraints that we need to check are \eqref{1.5} and \eqref{1.7};
\eqref{1.5}, which demands in particular that $\varphi_t^\ast(x) = x$ for $x\in H \sm B$,
holds because $\pi(x) =x$ on $H$. As for \eqref{1.7}, let $x\in H \cap L = \d H$ be given;
we know that $\varphi_t(x) \in L$, which means that its second coordinate vanishes (because 
$L \subset e_2^\perp$). Then $\varphi_t^\ast(x) = \pi(\varphi_t(x)) \in \d H \subset L$, as needed.
The other constraints hold easily, and in particular $\varphi_1^\ast$ is Lipschitz if $\varphi_1$ is Lipschitz;
hence $E^\ast = \varphi_1^\ast(H) = \pi(E)$ is another sliding competitor for $H$ in $B$.
Since $\pi$ is $1$-Lipschitz, we see that $\H^2(E^\ast \cap B) \leq \H^2(E \cap B)$
(the points from $\R^n \sm B$ do not contribute because $\varphi_1(w) = \varphi_1^\ast(w) = w$
for $w \in \R^n \sm B$). Thus it is enough to show that 
$\H^2(H \cap B) \leq \H^2(E^\ast \cap B)$, or that $E^\ast$ contains $H \cap B$, or also
(since $\varphi_1^\ast(w) = w$ on $\R^n \sm B$) that $\varphi_1^\ast(H) \supset H$.

For the topological argument that follows, it is easier to work with the plane $P$ 
that contains $H$. Let $\sigma$ denote the symmetry with respect to $L$. 
Then extend $\varphi^\ast$ to $P$ by setting
$\varphi_t^\ast(w) = \sigma(\varphi_t(\sigma(w)))$ for $w\in \sigma(H)$. 
When $x\in H$ tends to a point $x_0$ of $\d H$, $\varphi_t(x)$ tends to
$\varphi_t(x_0) \in L$, and $\varphi_t^\ast(x)$ tends to $\varphi_t^\ast(x_0) 
= \pi(\varphi_t(x_0)) \in \d H$. Then $\sigma(\varphi_t^\ast(x_0)) = \varphi_t^\ast(x_0)$.
Because of this, our extension $\varphi_t^\ast$ is continuous across $\d H$. 
In addition, it takes values in $P$, and it is the identity on $P \sm B$. So it is surjective, 
and for $\xi \in H \sm L$, we can find $w\in P$ such that 
$\varphi_t^\ast(w) = \xi$. But $\varphi_t^\ast(w) \in \sigma(H)$ when $w\in \sigma(H)$,
so $w \in P \sm \sigma(H) \subset H$. We are left with $H \cap L = \d H$, which is also contained
in $\varphi_t^\ast(H)$ because $H$ is closed and $\varphi_t^\ast$ is continuous.
This completes our proof of \eqref{32.1}; Lemma \ref{t32.1} follows.
\qed

\ms
The next lemma is a similar result that concerns $\bV$-sets. Since we also want to include
larger boundary sets $L$, we give some of the notation before the statement.

Let $V$ be a $\bV$ set, thus composed of two half spaces $H_1$ and $H_2$ bounded by 
a same line $\ell$, and that make make an angle at least $\frac{2\pi}{3}$ along $\ell$. 
Let us choose an orthonormal basis of $\R^n$ such that if 
$v_i$ denotes the unit vector in $H_i$ that is orthogonal to $\ell$, then
\begin{equation} \label{39.2}
\text{$\ell = e_1 \R$, $v_1 = (0,\cos\alpha,\sin\alpha,0) \in H_1$, and 
$v_2 = (0,-\cos\alpha,\sin\alpha,0) \in H_2$,}
\end{equation}
with $0 \leq \alpha \leq \frac{\pi}{6}$.
We will also accept boundaries $L$ larger than $\ell$, but require that $\ell \subset L$ and, when 
$\alpha > 0$, that
\begin{equation} \label{39.3}
 L \subset  (v_1+v_2)^\perp = e_3^\perp .
\end{equation}

\ms
\begin{lem}\label{t32.2}
Let $V$ be a $\bV$-set and $L \supset \ell$ be as above, and in particular satisfy \eqref{39.3}.
Then $V$ is a sliding minimal set, associated to the sliding boundary $L$.
\end{lem}

We first prove this lemma with the official Definition \ref{t1.1}, because this makes 
our life much simpler. Yet, we claim that, as for Lemma \ref{t32.1},
Lemma \ref{t32.2} stays true when we even allow  mappings $\varphi_1$ that are not 
Lipschitz in the definition of competitors for $V$, which makes the notion of minimal set 
and the lemma a little stronger. We'll prove this, and the minimality of $V$ for slightly 
different problems (related to separation conditions, and where $L$ is a plane) in the 
next lemma, but the proof will be less pleasant.

Let us also mention that E. Cavallotto \cite{Ca} already has a shorter proof of this 
with a slicing argument (where this time slices are defined as for currents, the proof
is done with polyhedral chains with coefficients in $\bZ_2$,
and boundary computations for chains replace the topological separation argument).
This proof also seems to require the official Definition \ref{t1.1} where $\varphi_1$ is 
required to be Lipschitz.

\ms
We start the proof with some reductions.
We know that the planes are minimal regardless of the sliding boundary $L$;
for this we may proceed as in Lemma \ref{t32.1}, but directly with the orthogonal projection
$\pi$ on the plane $V$. So we may assume that $\alpha > 0$.

Let $E = \varphi_1(V)$ be a sliding competitor for $V$ in some ball $B$,
coming from $\{ \varphi_t \}$ as in Definition \ref{t1.1}; we can assume that $B$
is centered at $0$ and want to show the analogue of \eqref{32.1} for $V$.
We can further assume that $B = \ol B(0,1)$, by scale invariance; this will just simplify
the notation a little.

Notice that $E = \varphi_1(V)$ is also a sliding competitor for $H$ in $B$, associated to
the largest possible boundary $L'  = e_3^\perp$, because $L \cap V = L' \cap V$, and
so $\{ \varphi_t \}$ also satisfies \eqref{1.7} with respect to  $L'$. That is, it is enough to prove
the lemma when $L = L'  = e_3^\perp$, which we assume now.

Denote by $W$ the $3$-dimensional space that contains $V$, and by 
$W_+$ the upper half space in $W$ defined by 
\begin{equation} \label{39.4}
W_+ = \big\{ (x,y,z,0) \, ; \, (x,y,z) \in \R^3 \text{ and } z \geq 0 \big\} \subset W.
\end{equation}
Denote by $\pi$ the closest distance projection on $W_+$, defined by $\pi(x,y,z,t) = (x,y,z_+,0)$.
Notice that $\pi$ is $1$-Lipschitz, is the identity on $W_+ \supset V$, and maps $L$ to $L \cap W$.

Set $\varphi_t^\ast = \pi \circ \varphi_t$ and $E^\ast = \pi(E) = \varphi_1^\ast(V)$.
Since $\varphi_t^\ast(x) \in L$ when $x\in V \cap L = \ell$,  $E^\ast$ is a sliding competitor
for $V$ in $B$. As before, $\H^2(E^\ast \cap B) = H^2(\varphi_1^\ast(V\cap B)) 
\leq H^2(\varphi_1(V\cap B))$, so is enough to prove that
\begin{equation} \label{32.2}
\H^2(V \cap B) \leq \H^2(E^\ast \cap B).
\end{equation}

We will do this with a slicing and separation argument. 
Observe that all the values taken by the $\varphi_t^\ast$ lie in $W_+$
(we composed by $\pi$). We will forget the last coordinates and 
denote by $(x,y,z)$ the coordinates of points of $W$. Set, for $x\in \R$,
\begin{equation} \label{32.3}
E^\ast_x = \big\{ (y, z) \in \R^2 \, ; \, (x,y,z) \in E^\ast \cap B \big\} \subset \R^2 
\end{equation}
and 
\begin{equation} \label{32.4}
V_x = \big\{ (y, z) \in \R^2 \, ; \, (x,y,z) \in V \cap B \big\} \subset \R^2;
\end{equation}
we want to show that 
\begin{equation} \label{32.5}
\H^1(V_x) \leq \H^1(E^\ast_x) \ \text{ for } x\in \R,
\end{equation}
but let us first explain why \eqref{32.2} follows from this, when we use the full
Definition \ref{t1.1} and therefore assume that $\varphi_1$ is Lipschitz.

Indeed, on the one hand, it follows from Fubini that
$\H^2(V \cap B) = \int_x \H^1(V_x) dx$, where here and below 
all the integrals are in fact taken on compact sets (projections of $B$). 
On the other hand, we claim that
\begin{equation} \label{32.6}
\int_x \H^1(E^\ast_x) dx \leq \H^2(E^\ast \cap B).
\end{equation}
As soon as we check this, we can integrate on $x$, use \eqref{32.4}, and get that
$\H^2(V \cap B) = \int_x \H^1(V_x) dx \leq \int_x \H^1(E^\ast_x) dx \leq \H^2(E^\ast \cap B)$,
as needed for \eqref{32.2}. Now \eqref{32.6} holds because $E^\ast$ is rectifiable,
which is true because $V$ is rectifiable and $E^\ast =  \varphi_1^\ast(V) = \pi(\varphi_1(V))$,
with mappings $\varphi_1$ and hence $\varphi^\ast_1$ that are Lipschitz.

To prove that \eqref{32.6} holds when $E^\ast$ is rectifiable, the simplest is to apply
the co-area formula (Theorem~3.2.22 in \cite{Federer}), to the restriction to $E^\ast \cap B$
of the orthogonal projection on $\ell$, after noticing that the appropriate Jacobian is at most $1$
and the corresponding level sets are the $E^\ast_x$.
But in fact, we are only using the easy part of the co-area formula, and it would not be 
too difficult to check \eqref{32.6} directly. 
We would first prove it for measurable subsets of $C^1$ surfaces (essentially by Fubini), 
and then take countable disjoint unions to go to the general case. The author wrote a little more
about this and the next lines near (4), in Section 76.b of \cite{MSBook}, in a slightly
similar context (but for Mumford-Shah minimizers).

When $E^\ast$ is not rectifiable, the author does not know whether \eqref{32.6}
necessarily holds. It would hold for some other variants of the Hausdorff measure, essentially 
as good as $\H^d$, but this is not an excuse. Since we do not necessarily want to assume that 
our competitor $E$ comes from a Lipschitz deformation $\varphi_1$, we will give a different
proof of minimality in a next lemma.

But here we assumed that $\varphi_1$ is Lipschitz, and then $E^\ast =\varphi_1^\ast(V)$ 
is rectifiable and \eqref{32.6} holds. Now we just need to show \eqref{32.5}
(for $x$ in the projection of $B$) because \eqref{32.2} will ensue, and for this we want to 
show that $V_x$ is minimal for some problem. We will use separation properties of $V$ and 
$\varphi_1^\ast(V)$, and again it is easier to work on $W$ rather than $\R^n$, which is why 
we composed with $\pi$ in the first place.
We have mappings $\varphi_t^\ast : V \to W_+$, that we would like to extend to $W$.

Let $\sigma$ denote the reflection across 
$P =  \big\{ (x,y,0) \, ; \, (x,y) \in \R^2 \big\} = L \cap W$; 
thus $\sigma(x,y,z) = (x,y,-z)$ for $(x,y,z) \in W$. We define $\varphi_t^\ast$ on $\sigma(V)$
by $\varphi_t^\ast(w) = \sigma(\varphi_t(\sigma(w)))$ for $w\in \sigma(V)$.
Notice that $\sigma(\varphi_t^\ast(w)) = \varphi_t^\ast(w)$ for $w\in V \cap \sigma(V) = \ell$,
because  $\varphi_t(w) \in L$ and hence $\varphi_t^\ast(w) \in L\cap W = P$. 
Because of this, we now have a continuous mapping from $[V \cup \sigma(V)] \times [0,1]$ 
to $W$. 

Consider the three points $A_1 = (0,2,0)$, $A_2 = (0,-2,0)$, and
$A_3 = (0, 0, 2)$ (above $V$). We want to show that 
\begin{equation} \label{32.7}
\varphi_t^\ast(V \cup \sigma(V)) \ \text{ separates the three points $A_i$ from
each other in } W
\end{equation}
and for this we apply 4.3 in Chapter~XVII of \cite{Dugu}, on page 360.
The point is that for $t=0$, $\varphi_t^\ast(V \cup \sigma(V)) = V \cup \sigma(V)$,
which separates these three points in $W$, and then when we deform this set continuously, 
it never crosses (or get close to) our three points $A_i$. 
To be fair, Dugundji only gives the result when the initial set
is compact, so we should modify things a little bit. That is, assume to the contrary
that $A_1$ and $A_2$ are connected in $W \sm \varphi_t^\ast(V \cup \sigma(V))$ 
by a path $\gamma$, and choose $R > 2$ such that $\gamma \subset B(0,R)$. 
Then consider the sets $Z_t = [\varphi_t^\ast(V \cup \sigma(V))] \cap B(0,R) \cup \d B(0,R)$.
They are compact, still represent a continuous deformation of $Z_0$, and never get close 
to the $A_i$, hence $Z_t$ separates $A_1$ from $A_2$ because $Z_0$ does,
and this contradicts the existence of $\gamma$. So \eqref{32.7} holds.

The same deformation argument, without the reflection, also shows that 
\begin{equation} \label{39.11}
\varphi_t^\ast(V) \ \text{ separates $A_3$ from $A_1$ and $A_2$ in  } W,
\end{equation}
because this is true for $t=0$.
We return to \eqref{32.7} and claim that
\begin{equation} \label{32.8}
\varphi_t^\ast(V) \ \text{ separates the three points $A_i$ in  } W_+,
\end{equation}
where $W_+$ as in \eqref{39.4}.
Otherwise, there is a path $\gamma \subset W_+ \sm \varphi_t^\ast(V)$ that goes from
some $A_i$ to some other $A_j$. Since $\varphi_t^\ast(V \cup \sigma(V))$ separates,
$\gamma$ meets $\varphi_t^\ast(\sigma(V))$. Most of this last set lies in the complement of
$W_+$, with the exception of $P \cap \varphi_t^\ast(\sigma(V)) = P \cap \varphi_t^\ast(V)$.
But we assumed that $\gamma$ does not meet $\varphi_t^\ast(V)$; this contradiction proves 
\eqref{32.8}.

\ms
Now let $x \in \R$ be given, and consider slices by the corresponding vertical plane.
Call $W'_+ = \big\{ (y, z) \in \R^2 \, ; \, z \geq 0 \big\}$ the (constant) slice of $W_+$, 
and denote by $A'_i$ the ``projection'' of $A_i$, where we just forget the first coordinate.
Finally let 
\begin{equation} \label{39.13}
F_x = \big\{ (y, z) \in \R^2 \, ; \, (x,y,z) \in \varphi_1^\ast(V)  \big\}
\end{equation}
be the slice of $\varphi_1^\ast(V)$.
Notice that $A_{i,x} = A_i + (x,0,0)$ lies in the same component of $W_+ \sm \varphi_t^\ast(V)$
as $A_i$, so \eqref{32.8} (with $t=1$) implies that 
\begin{equation} \label{39.14}
F_x \ \text{ separates the three points $A'_i$ inside of } W'_+.
\end{equation}
Similarly, by \eqref{39.11},
\begin{equation} \label{39.15}
F_x \ \text{ separates $A'_3$ from $A'_1$ and $A'_2$ in  } \R^2.
\end{equation}

\ms
At this point we would probably have enough information for a calibration argument,
but we decided to use connectedness, so we need some plane topology. We first use 
\eqref{39.15} in $B' = \ol B(0,2) \subset \R^2$. Since $\varphi_1^\ast(V)$ only differs
from $V$ in $\ol B(0,1)$, the set $F_x$ meets $\d B'$ only twice, at the points 
$B_\pm = (\pm\cos\alpha,\sin\alpha)$. By Theorem~14.3 on p.~123 of [Ne], there is a 
connected component of $H_0$ of $F_x \cap B'$ that separates $A'_1$ from $A'_3$;
this component contains both points $B_\pm$, because otherwise one of the two arcs of $\d B'$
that go from $A'_1$ to $A'_3$ does not meet $F_x$. 

We now use \eqref{39.14}, in $B' \cap W'_+$ (a closed half disk). This time we get a
connected component $H_1$ of $F_x \cap B' \cap W'_+$ that separates $A'_1$ from
$A'_2$. This set meets $\d B' \cap W'_+ = \{ B_+, B_- \}$, because otherwise
we could use an arc of $\d B'$ to go from $A'_1$ to $A'_2$ without meeting $F_x$.
So $H_1$ meets $H_0$.
But we can also try to connect $A'_1$ to $A'_2$ directly with $[A'_1,A'_2] \subset \d W'_+$,
so $H_1$ meets this segment too.

The set $H_0 \cup H_1$ is connected, contained in $F_x$, and it contains
$B_+$, $B_-$, and some point of $[A'_1,A'_2] \subset \d W'_+$. Said in other words,
$F_x \cap B'$ connects $B_+$, $B_-$, and $[A'_1,A'_2]$. 

We claim that among all connected sets (like $H_0 \cup H_1$) that do this, 
the slice $V'_x$ of $V$ is the shortest.
For the standard proof, we would first replace any connected set that contains the two $B_\pm$ 
and a point $B_0 \in [A'_1,A'_2]$ with a shorter one composed of at most $3$ line segments emanating 
from a point of connection between two arcs from $B_0$ to the $B_\pm$, and then optimize 
the position of that point. Our assumption that $\alpha \leq \frac{\pi}{6}$ is of course used there.
So $\H^1(F_x \cap W_0) \geq \H^1(H_0 \cup H_1) \geq \H^1(V'_x)$, and when we remove the 
identical contribution from $\ol B(0,2) \sm \ol B(0,1)$, we get \eqref{32.5}.

This completes our proof of Lemma \ref{t32.2} when we use the official Definition \ref{t1.1}
to define competitors and minimal sets.
\qed

\ms
We shall now discuss the minimality of $V \in \bV(L)$ when we insist on removing \eqref{1.8}
from the definition of competitors. The author did not find any other way than going 
to the separation properties suggested by the proof above, finding a set that 
satisfies these properties and minimizes the Hausdorff measure, and finally proving 
that this set is rectifiable, hence has a larger measure than $V$ by the proof of Lemma \ref{t32.2}. 
Since this set is a minimizer, we will also get a control on the other competitors, even the ones 
that are not rectifiable.

It is a little strange that we should need to do this, but the point is that given a competitor $E$ to the 
separation problem below, it is not so easy to find directly, by modifying $E$, a rectifiable 
competitor $E'$ that does better than $E$, even though it is much easier to show that minimizers
are rectifiable. The difficulty will then be to find an equivalent weak problem, with the same solutions,
but for which we can prove that minimizers exist. The compactness properties of $BV$ and 
Caccioppoli sets will be useful for this. 

We shall now state the generalization of Lemma \ref{t32.2}, set the strong and weak separation 
problems, solve the weak one, show that it is equivalent to the strong one, and then conclude.

\ms
\begin{lem}\label{t32.3}
Lemma \ref{t32.2} (about the minimality of $\bV$ sets) is still valid when we forget the Lipschitz
property \eqref{1.8} in Definition \ref{t1.1}.
\end{lem}

\ms
Before we start the proof, let us add some last comments. So far the prevailing definition
of minimality is with the Lipschitz condition \eqref{1.8}, both because Almgren asked for it,
it accommodates currents and varifolds better, and at the same time adding it makes
the regularity theorems a little better and essentially costs nothing. But the author feels that 
there may be situations soon where we can only prove existence theorems for our sliding Plateau
problem or variants, where the minimizer that we find is a competitor of our initial data $E_0$,
but perhaps not with the Lipschitz condition \eqref{1.8};
then we will need to make tough decisions, and possibly leave \eqref{1.8} behind. 
In this context Lemma \ref{t32.3} makes more sense.

Another related situation that arises some times is that we initially set a problem concerning all
closed sets that satisfy some topological condition, and only have a proof of existence in the category
of rectifiable sets. A simple solution is to pretend that only the rectifiable sets are interesting
(which makes some sense because if there is a minimizer, it is rectifiable), and replace our initial problem
with its variant for rectifiable sets only. The author does not think that it is entirely satisfactory, 
but yet does not have a general solution to this issue. In the present case, it turns out that thanks
to the existence of an equivalent weak problem, we can prove an existence theorem first, and then
prove after the fact that rectifiable sets are enough. Since the author does not know, 
apart from \cite{Pacific} which is in a slightly different context, of cases where this was done, 
he decided to put the argument here, for possible later reference.
It also turned out that setting the weak problem right was more complicated than he had expected,
partially due to the sliding condition, and this is one more reason for not addressing the issue upfront.

\ms
Our proof will rely on ideas and results from \cite{Pacific}.
We will state separation problems that are specific to the problem at hand,
but one of the reason why the author thought it would be a good idea to prove the
present lemma is that the techniques should probably be useful for similar problems.
He also found out during the writing process that although separation problems should always
end up with a BV statement, the precise statement is not always easy to find.

Let us first state a \ub{strong separation problem}. 
Consider the balls $B = \ol B(0,1) \subset \R^3$ and $2B = \ol B(0,2)$, 
the upper half ball $2B_+ = \big\{ (x,y,z) \in 2B \, ; \, z \geq 0\big\}$, and let $V$
denote the same $\bV$-set as above \eqref{39.2}, with $\alpha > 0$. We consider the class $\cF_s$
of compact sets $E \subset 2B_+$ such that $\H^2(E) < +\infty$,
$E$ coincides with $V$ in $2B \sm B$, $E$
separates the two points $A_1 = (0,2,0)$ and $A_2 = (0,-2,0)$ from $A_3 = (0, 0, 2)$ 
inside of $2B$, and also separates $A_1$ from $A_2$ in $2B_+$. Thus 
the sets $E^\ast \cap B = \varphi_1^\ast(V) \cap B$ of Lemma~\ref{t32.2} lie in $\cF_s$,
by \eqref{39.11} and \eqref{32.8}, as soon as they have a finite measure
(and otherwise \eqref{32.2} is clear). A minimizer for the strong separation problem (SSP)
will just be a set $E_0 \in \cF_s$ such that 
\begin{equation} \label{39.16}
\H^2(E_0) = m_s, \ \text{ with } m_s = \inf_{E \in \cF_s} \H^2(F).
\end{equation}
As soon as we know that there exists a minimizer for the SSP and prove that it is rectifiable, 
the proof of Lemma \ref{t32.2}, and in particular of \eqref{32.2}, will say that 
$\H^2(E_0) \geq \H^2(V \cap B)$, so in fact
$\H^2(V \cap B) = m_s$ (because $V \cap B \in \cF_s$)
and then \eqref{32.2} also holds for sets $E^\ast$ that are not rectifiable.
This is what we want, even though in the present case the notion of minimizer for the SSP is not so 
interesting because we'll find out a posteriori that $E_0$ was in fact equivalent to $V$.
We could of course imagine other situations where it is more interesting.

\ms
Unfortunately, getting minimizers for the SSP directly seems unpleasant, so we introduce a
\ub{weak separation problem}. Since the author is always a bit afraid of trace conditions in the
set $BV$ of functions of bounded variation, we use a security ring and set our problem
in $D = B(0,2) \subset \R^3$. Set 
$D_+ = \big\{ (x,y,z) \in 2B \, ; \, z \geq 0\big\}$ and 
$P = \big\{ (x,y,z) \in D \, ; \, z = 0\big\}$, and denote by $\sigma$ the reflection across $P$.

Denote by $\cC(D)$ the set of Caccioppoli subsets of $D$, i.e. measurable subsets
of $D$ such that $\1_{F} \in BV$. For such a set, we will denote by $|D\1_{F}|$
the total variation of $D\1_{F}$ (in $\R^3$, not just in $D$); it is a finite positive 
Borel measure, and its total mass is called the perimeter of $F$. 

Our set of competitor will be the set $\cF$ of quadruples 
$\ub F = (F_1,F_2,F_3,F_4) \in \cC(D)^4$ with the following properties.
First, $D$ is the essentially disjoint union of the $F_i$, i.e., 
\begin{equation} \label{39.17}
\sum_{i=1}^3 \1_{F_i} = \1_D.
\end{equation}
Incidentally, we work modulo a set of vanishing Lebesgue measure, as always with
Caccioppoli sets. Next denote by $G_i$, $1 \leq i \leq 3$, the connected component of 
$\ol D \sm (V \cup \sigma(V))$ that contains $A_i$. We require that 
\begin{equation} \label{39.18}
\sigma(F_i) = F_i  \ \text{ and } \ F_i \supset G_i
\ \text{ for }    i = 1, 2,
\end{equation}
and that
\begin{equation} \label{39.19}
G_3 \subset F_3 \subset D_+  \ \text{ and }  F_4 = \sigma (F_3).
\end{equation}

\begin{figure}[!h]  
\centering
\includegraphics[width=12.cm]{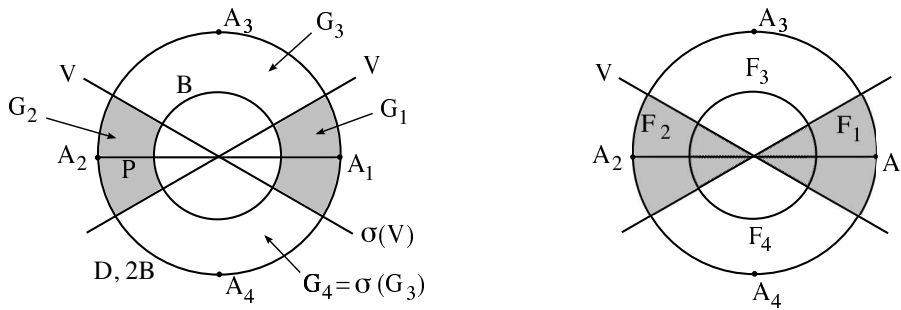}
\caption{The central section of $D$ and the sets $G_i$, and on the right 
the competitor associated with $V$
\label{f39.1a}
}
\end{figure}

\begin{figure}[!h]  
\centering
\includegraphics[width=9.cm]{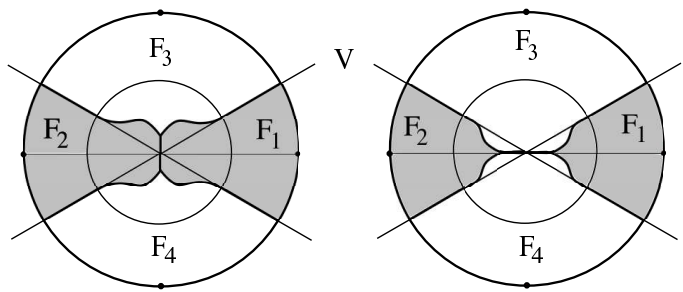}
\caption{The same section of $D$, with two examples of slightly different competitors
\label{f39.1b}
}
\end{figure}

See the left part of Figure \ref{f39.1a} for a picture of the $A_i$ and $G_j$, 
the right part for the most obvious quadruple $\ub F$, and Figure \ref{f39.1b} for 
two slightly different examples of quadruples $\ub F$; 
we'll explain soon why we decided to double everything.

We also need to define a functional on $\cF$. For $\ub F = (F_1,F_2,F_3,F_4) \in \cF$, set
\begin{equation} \label{39.20}
f_i = \1_{F_i} \in BV(\R^3) \ \text{ and call } \mu_i = |Df_i| = |D\1_{F_i}|
\end{equation}
the total variation of $D\1_{F_i}$. 
So $\mu_i(A)$ is the perimeter of $F_i$ in $A$ when $A$ is an open set. 
Here and below, we will refer to \cite{Gi} for the various properties of $BV$ functions and 
Caccioppoli sets that will be used. 
We set
\begin{equation} \label{39.21}
J(\ub F) = \mu_3(P)+\mu_4(P) + \sum_{i=1}^4 \mu_i(D).
\end{equation}

% mis page 323

% mis page 323

Let us explain why we do this. 
We treat $F_3$ and $F_4$ differently because in simple situations like the one
suggested by Figure \ref{f39.1a} (right), where $F_1 \cap D_+$, $F_2 \cap D_+$, and $F_3$ are 
the components of $D_+ \sm E$ for a nice competitor $E \in \cF_s$ with the same
topology as for $V$, we want to obtain $J(\ub F) = 4\H^2(E)$; the additional term 
$\mu_3(P)+\mu_4(P)$ is designed for this. 
The author's first attempt was to work on the interior of $D_+$
and add $2\mu_3(P)$; this should give the same result in the nice cases, but defining
$J$ as in \eqref{39.21} is a way to prevent ugly mixtures of $F_1$ and $F_2$ on
$P$, that would not have been counted in $\mu_1+\mu_2$ of the interior of $D_+$,
but would be counted here. The reader is invited to check that \eqref{39.21} does what
we want in the simple examples suggested by Figure \ref{f39.1b}. %

The reader may be shocked because we integrate $\mu_3$ and $\mu_4$
against a function which is not lower semicontinuous (its value on $P$ is larger than
its limit), which is usually not a good idea when we want to prove existence for the minimizers,
and indeed we will need to be careful when we show that
\begin{equation} \label{39.22}
\text{we can find $\ub F \in \cF$ such that $J(\ub F) = m$, where } 
m = \inf_{\ub F' \in \cF} J(\ub F').
\end{equation}
So let us look for $\ub F$. Let $\{ \ub F_k \}$ be a minimizing sequence in $\cF$, i.e., 
such that $\lim_{k \to +\infty} J(\ub F_k) = m$. 
Write $\ub F_k = (F_{1,k},F_{2,k},F_{3,k}, F_{4,k})$, and set as above
\begin{equation} \label{39.23}
f_{i,k} = \1_{F_{i,k}} \in BV(\R^3) \ \text{ and } \mu_{i,k} = |Df_{i,k}| = |D\1_{F_{i,k}}|.
\end{equation}
Notice that for $1 \leq i \leq 4$, the total perimeter of $F_{i,k}$ stays bounded, 
because it is less than $m + \H^2(\d D)$.
So $\{ f_{i,k} \}$ is a bounded sequence in $BV(\R^3)$, and since $f_{i,k}$ is supported 
in $B_2$, $\{ f_{i,k} \}$ is a relatively compact subset of $L^1(B_2)$. 
We extract a subsequence (still denoted the same way) that converges in $L^1$ to a limit $f_i$. 
Let us even extract a subsequence which works for the four indices at the same time, 
and so that we have pointwise convergence almost everywhere as well. 
The fact that $f_i$ is the characteristic function of a set $F_i$ (defined modulo
a set of vanishing measure), and the properties \eqref{39.17}, \eqref{39.18}, and \eqref{39.19},
follow from the almost everywhere convergence. Also, $f_i \in BV(\R^3)$ by the lower 
semicontinuity of the $BV$ norm; so $\ub F = (F_1,F_2,F_3,F_4)$ lies in $\cF$.
The lower semicontinuity of the $BV$ norm also says that for $1 \leq i \leq 4$,
\begin{equation} \label{39.24}
\mu_i(U) \leq \liminf_{k \to +\infty} m_{i,k}(U)
\ \text{ for every open set } U \subset D.
\end{equation}
If we did not have the terms $\mu_3(P)$ and $\mu_4(P)$ in the definition of $J$, 
we would apply this with $U=D$, get that $J(\ub F) \leq m$, and deduce at once that $\ub F$ 
is the desired minimizer. We now have to show that the additional terms $\mu_3(P)$ and $\mu_4(P)$ 
can be estimated by contributions of the measures $m_{i,k} = \mu_{i,k} = |D\1_{F_{i,k}}|$ to sets 
that lie close to $P$ (and may disappear in the limit of measures).

We will need a little more information on the $\mu_i$, which the reader may find in \cite{Gi}.
For $F_i \in \cC(D)$, there is a measurable and rectifiable set $\d^\ast F_i \subset \d F_i$,
called the reduced boundary of $F_i$, such that 
\begin{equation} \label{39.25}
\mu_i = \H^2_{\vert \d^\ast F_i}
\end{equation}
(so that in particular $\mu_i(P) = \H^2(P \cap \d^\ast F_i)$)
and with some additional geometric properties that we recall now.

For $x \in \d^\ast F_i$, there is an approximate tangent plane $P(x) = P_i(x)$ 
to $\d^\ast F_i$ at $x$ and a unit normal $n(x)=n_i(x)$ to $P(x)$, 
with the following properties. For each $\varepsilon > 0$,
there exists $r(x) > 0$ such that for $0 < r < r(x)$,
\begin{equation}\label{39.26}
\big| \big\{ y\in B(x,2r) \cap F_i \, ; \, \langle y-x, n(x)\rangle \leq 0 \big\} \big| \leq \varepsilon r^3
\end{equation}
\begin{equation}\label{39.27}
\big| \big\{ y\in B(x,2r) \sm F_i \, ; \, \langle y-x, n(x)\rangle \geq 0 \big\} \big| \leq \varepsilon r^3,
\end{equation}
and the density of $\d^\ast F_i$ is close to $\pi$, i.e.,
\begin{equation}\label{39.28}
\big| \H^2(\d^\ast F_i \cap B(x,r)) - \pi r^2 \big| \leq \varepsilon r^2.
\end{equation}

Notice that when $i=3$ and $x\in \d^\ast F_3 \cap P$, \eqref{39.27} leave us no choice because 
$F_3$ is contained in the upper half of $B(x,r)$ : $P(x)$ must be the plane that contains $P$,
and $n(x) = (0,0,1)$. Notice that $x$ is then an interior point of $D$, again by \eqref{39.27}.

Let $\varepsilon > 0$ be small, and let $\tau > 0$ be small too, to be chosen later 
(depending on $\varepsilon$). Since for $x\in P \cap \d^\ast F_3$ we can choose arbitrarily small 
radii $r < \tau$ such that \eqref{39.26}-\eqref{39.28} hold and $B(x,2r) \subset D$, 
Vitali's covering lemma (see for instance \cite{Mattila}) allows us to pick disjoint balls
$B_j = B(x_j, r_j)$ centered on $P \cap \d^\ast F_3$, with the properties above, and such that
\begin{equation}\label{39.29}
\mu_3\big(P \sm \cup_{j} B_j \big) = 
\H^2\big(P \cap \d^\ast F_3 \sm \cup_{j} B_j \big) \leq \varepsilon
\end{equation}
by \eqref{39.25}.
Then we can even find a finite subcollection $\{ B_j \}$, $j\in J$, with the same property
(but with $2\varepsilon$ in the right-hand side). By \eqref{39.25} and \eqref{39.28}, we get that
\begin{equation}\label{39.30}
\mu_3(P) = \H^2\big(P \cap \d^\ast F_3)
\leq (\pi+\varepsilon) \sum_{j\in J} r_j^2 + 2\varepsilon
\leq \pi \sum_{j\in J} r_j^2 + C \varepsilon
\end{equation}
because $\sum_{j\in J} r_j^2 \leq \frac{2}{\pi} \sum_{j\in J} \H^2(\d^\ast F_3 \cap B_j)
\leq \frac{2}{\pi} \H^2(\d^\ast F_3) = \frac{2}{\pi} \mu_3(\R^3) \leq m$ by \eqref{39.28},
because the $B_j$ are disjoint, and by \eqref{39.21}.

We want to  compare this with what we get for $\ub F_k$, $k$ large. 
Define small cylinders $T_j$ and $T_j^+$ by 
\begin{equation} \label{39.31}
\begin{aligned}
&T_j = \big\{(x,y,z)  \, ; \, (x,y,0) \in B_i \text{ and } - r_j  < z < r_j \big\}, 
\cr& T_j^+ = \big\{(x,y,z)  \, ; \, (x,y,0) \in B_i \text{ and } 0  < z < r_j \big\}.
\end{aligned}
\end{equation}
Since the $f_{i,k}$ converge in $L^1$ to  $f_i$, we deduce from 
\eqref{39.27} that for $k$ large, 
\begin{equation} \label{39.32}
\big|T_j^+ \sm F_{3,k}\big| \leq 2\varepsilon r_j^3.
\end{equation}
For the computation that follows, we may fix $k$ and $j$. We see $T_j$ as the product 
$B' \times (-r_j,r_j)$, where $B' = \big\{(x,y,z)  \, ; \, (x,y,0) \in B_j \big\} \simeq B_j \cap P$.
By \eqref{39.32}, Fubini and Chebyshev, we can choose $\rho \in (r_j/2,r_j)$ such that 
if we define a good set $A \subset B'$ by
\begin{equation} \label{39.33}
A =  \big\{(x,y) \in B'  \, ; \,  (x,y,\rho) \in F_{3,k} \big\},
\end{equation}
then 
\begin{equation} \label{39.34}
\H^2(B' \sm A) \leq C \varepsilon r_j^2.
\end{equation}
We are now going to use the fact that $f_{i,k} \in BV(T_j)$, seen as a function of
the last variable $z$, lies in $BV(-r_j,r_j)$, with the logical estimate for the norm.
That is, Theorem 3.103 on page 195 of \cite{AFP} says that
for almost every $(x,y) \in B'$, the function $z \to f_{i,k}(x,y,z)$
lies in $BV(-r_j,r_j)$, and its derivative is a measure $\nu_{x,y,i}$. In addition,
\begin{eqnarray} \label{39.35}
\int_{B'} |\nu_{x,y,i}|(-r_j,r_j) dxdy 
&=& \int_{B'} ||f_{i,k}(x,y,\cdot)||_{BV(-r_j,r_j)} dxdy 
\nn\\
&=& \int_{(x,y,z) \in T_j} \Big|\frac{\d f_{i,k}}{\d z}\Big|
\leq \int_{(x,y,z) \in T_j} \big|Df_{i,k}\big|  = \mu_{i,k}(T_j),
\end{eqnarray}
where we use a homogeneous norm on $BV$,
systematically denote the total variations of measures like absolute values,
and the only step which is not a definition is the second identity, which holds by a combination
of the weak definition of $BV$, Fubini, and a correct choice of product test functions.

The result above also says that $\nu_{x,y,i}$ is a finite measure for almost every $(x,y) \in B'$.
Since $f_{i,k}$ only takes the values $0$ and $1$, this means that (maybe after changing it
on a set of measure $0$) it has a finite number
$N_{i,k}(x,y)$ of jumps, and a simpler way to write \eqref{39.35} is 
\begin{equation} \label{39.36}
\int_{B'} N_{i,k}(x,y) dx dy \leq \mu_{i,k}(T_j).
\end{equation}
Consider points $(x,y) \in A$. Thus $(x,y,\rho) \in F_{3,k}$, and by symmetry 
$(x,y,-\rho) \in F_{4,k}$. 

Start with the set $A_0$ of points of $A$ such that $N_{1,k}(x,y)=N_{2,k}(x,y) =0$.
For such a point, $(x,y,z)$ stays in $F_{3,k} \cup F_{4,k}$, and in fact (by \eqref{39.19})
in $F_{3,k}$ for $z>0$ and in $F_{4,k}$ for $z<0$.
We claim that for almost every point $(x,y)\in A_0$, the point $(x,y,0)$ lies in both sets
$\d^\ast F_{3,k}$ and $\d^\ast F_{4,k}$. Indeed, if $(x,y)$ is a Lebesgue density point of $A_0$,
 then the density of $F_{3,k}$ at $(x,y)$ is $1/2$ (because all the vertical half lines that are
 contained in $F_{3,k}$ locally), and similarly for $F_{4,k}$ (this time, we look down). Then 
 (by Poincar\'e) $(x,y)$ is a point of positive lower density for $\mu_{3,k}$ and $\mu_{4,k}$,
 and almost all such points lie in $\d^\ast F_{3,k}$ and $\d^\ast F_{4,k}$. Said otherwise,
 $(x,y)$ lies on what is called the geometric measure boundaries of $F_{3,k}$ and $F_{4,k}$, with the
 same conclusion. So we like $A_0$ because 
 \begin{equation} \label{39.37}
2\H^2(A_0) \leq \mu_{3,k}(P\cap B_j)+\mu_{4,k}(P\cap B_j).
\end{equation}
Now let $A_1$ be the set of points of $A$ such that $N_{1,k}(x,y)>0$. This means that
$f_{1,k}(x,y,z)$, which starts at $0$ for $z=\rho$, becomes $1$ for some $z<\rho$. By symmetry
of $F_{1,k}$, this happens first for $z \geq 0$.
Since for almost every $(x,y)$, the total variation of $f_{1,k}(x,y,\cdot)$ on $(-r_j,r_j)$ is finite,
and also we only count jumps when there is a positive measure of each set near the jump
point $z$, this means that we can take $z > 0$. Then, by symmetry, $f_{1,k}(x,y,\cdot)$
also has a jump at $-z$, so $N_{1,k}(x,y) \geq 2$. These jumps also count for $N_{3,k}(x,y)$ 
when $z > 0$ (because $z$ was chosen to be the first point below $\rho$) and hence $N_{4,k}(x,y)$ 
when $z < 0$, so $N_{3,k}(x,y) \geq 1$ and $N_{4,k}(x,y) \geq 1$.

Similarly if $A_2$ is the set of points of $A$ such that $N_{2,k}(x,y)>0$, then 
for almost every $(x,y)$, $N_{2,k}(x,y) \geq 2$, $N_{3,k}(x,y) \geq 1$, and $N_{4,k}(x,y) \geq 1$.
Thus 
\begin{equation} \label{39.38}
4\H^2(A_1 \cup A_2) \leq \int_{A_1 \cup A_2} 
\sum_{i=1}^4 N_{i,k}(x,y) dx dy \leq \sum_{i=1}^4 \mu_{i,k}(T_j \sm P).
\end{equation}
by the proof of \eqref{39.35}. 
We add this to (twice) \eqref{39.37}, use \eqref{39.34}, notice that 
$A \subset A_0 \cup A_1 \cup A_2$, and get that
\begin{eqnarray} \label{39.39}
4 \pi r_j^2 =  4 \H^2(B') &\leq& 2\mu_{3,k}(P\cap B_j)+2\mu_{4,k}(P\cap B_j) 
+ \sum_{i=1}^4 \mu_{i,k}(T_j \sm P) + C \varepsilon r_j^2
\nn\\
&\leq& \mu_{3,k}(P\cap B_j)+\mu_{4,k}(P\cap B_j) 
+ \sum_{i=1}^4 \mu_{i,k}(T_j) + C \varepsilon r_j^2,
\end{eqnarray}
which is the contribution of $T_j$ to $J(\ub F_k)$ in \eqref{39.21}. We return to the full collection
of balls $B_j = B(x_j,r_j)$, sum \eqref{39.39} over $j$, compare to \eqref{39.30}, and get that
\begin{equation} \label{39.40}
4\mu_3(P) \leq 4\pi \sum_{j\in J} r_j^2 + C \varepsilon
\leq C \varepsilon + \sum_{j\in J} \Big[\mu_{3,k}(P\cap B_j)+\mu_{4,k}(P\cap B_j) 
+ \sum_{i=1}^4 \mu_{i,k}(T_j) 
\Big],  
\end{equation}
where we also used the line below \eqref{39.30} to sum the $r_j^2$.
Recall that all the $r_i$ were chosen smaller than the small $\tau > 0$, to be chosen small
soon. Thus, if we set
\begin{equation} \label{39.41}
H(\tau) = \big\{ (x,y,z) \in D \, ; \,  |z| \leq \tau \big\}
\ \text{ and  } H'(\tau) = H(\tau)\sm P,
\end{equation}
then \eqref{39.40} says that
\begin{equation} \label{39.42}
4\mu_3(P) \leq C \varepsilon + \mu_{3,k}(P) + \mu_{4,k}(P) + \sum_{i=1}^4 \mu_{i,k}(H(\tau))
\end{equation}
because the $B_j$ and the $T_j$ are disjoint. We will return to \eqref{39.42}, but we
also need to control $\mu_i(P)$ for $i=1,2$. We claim that
\begin{equation} \label{39.44}
\mu_i(P) = \H^2(P \cap \d^\ast F_i) = 0 \ \text{ for } i=1,2, 
\end{equation}
because $P \cap \d^\ast F_i$ cannot have a Lebesgue density point in $P$.
Indeed, for such a point, the approximate tangent plane would need to be the plane
that contains $P$, and then \eqref{39.26} or \eqref{39.27} would contradict the symmetry
of $A_i$. Finally, for $D \sm H(\tau)$, \eqref{39.24}
says that $\mu_i(D \sm H(\tau)) \leq \mu_{i,k}(D \sm H(\tau)) + \varepsilon$
for $k$ large. We add everything and get that
\begin{eqnarray} \label{39.45}
J(\ub F) &=& \mu_3(P) + \mu_4(P) + \sum_{i=1}^4 \mu_i(D)
= 4 \mu_3(P) + \sum_{i=1}^4 \mu_i(D \sm P)
\nn\\
&\leq& C \varepsilon + \mu_{3,k}(P) + \mu_{4,k}(P) + \sum_{i=1}^4 \mu_{i,k}(H(\tau))
+ \sum_{i=1}^4 \mu_i(D \sm P)
\nn\\
&\leq& C \varepsilon + 2\mu_{3,k}(P) + 2\mu_{4,k}(P) + \sum_{i=1}^4 \mu_{i,k}(H'(\tau))
+ \sum_{i=1}^4 \mu_i(D \sm P)
\end{eqnarray}
by \eqref{39.21}, \eqref{39.19} (for the symmetry), \eqref{39.44}, \eqref{39.42},
and the definition \eqref{39.41}. Next
\begin{equation} \label{39.46}
\sum_{i=1}^4 \mu_i(D \sm P) \leq \varepsilon + \sum_{i=1}^4 \mu_i(D \sm H(\tau))
\leq 2\varepsilon  + \sum_{i=1}^4 \mu_{i,k}(D \sm H(\tau))
\end{equation}
by \eqref{39.43} and \eqref{39.24}. 
We now chose $\tau$ so small that
\begin{equation} \label{39.43}
\sum_{i=1}^4 \mu_i(H'(\tau)) < \varepsilon
\end{equation}
(recall that the intersection of the sets $H'(\tau)$ is empty). We sum and get that
\begin{equation} \label{39.47}
J(\ub F) \leq C \varepsilon + 2\mu_{3,k}(P) + 2\mu_{4,k}(P) + \mu_{i,k}(D \sm P)
\leq C\varepsilon + J(\ub F_k)
\end{equation}
by \eqref{39.21} for $\ub F_k$. For each $\varepsilon > 0$ we found that this holds for $k$ large.
Since $\{ \ub F_k \}$ was chosen to be a minimizing sequence, we get that $J(\ub F) \leq m$;
\eqref{39.22} follows, because we knew already that $\ub F \in \cF$.

\ms
Next we take the minimizer $\ub F$ given by \eqref{39.22}, and study its regularity
to show that it also gives a minimizer for the strong problem.
For this part we will follow the first steps of \cite{Pacific}, and often refer to it for details.
Denote by $Z$ the closed support of
$\mu = \sum_{i=1}^3 \mu_i$; our next task is to show that $\mu$ is locally Ahlfors regular
in $B = \ol B(0,1)$, which means that
\begin{equation} \label{39.48}
C^{-1} r^2 \leq \mu(B(x,r)) \leq C r^2 \ \text{ for $x \in Z \cap B$ and $0 < r < 10^{-1}$.} 
\end{equation}
Compared to \cite{Pacific}, there is a small difference, because we have 
four sets $F_i$ rather than two (a set and its complement), but this will not matter 
for what we have to do. 

Let $B(x,r)$ be as in \eqref{39.48}, and set $B^\ast = B(x,r) \cup \sigma(B(x,r))$.
We start with the upper bound, and even the special case when
\begin{equation} \label{39.49}
B(x,2r) \subset B.
\end{equation}
If $\mu(B(x,r))$ is too large, we replace 
$\ub F = (F_1,F_2,F_3,F_4)$ with $\ub F' = (F'_1,F'_2,F'_3,F'_4)$,
obtained by taking $F'_1 = F_1 \cup B^\ast$, and 
$F'_i = F_i \sm B^\ast$ for $i \geq 2$. 
This way we replace the total contribution $\mu(\ol B(x,r))$ with at most 
$\H^2(\d B^\ast) \leq C r^2$. This gives a quadruple $\ub F' \in \cF$, because 
$B(x,2r) \subset B$ does not see the constraints about the $G_i$.
Since $\ub F$ is a minimizer, we get the second part of \eqref{39.48} in this case.
When \eqref{39.49} fails, we just proceed with more caution, and only modify things inside
of $B$; that is, we take $F'_1 = F_1 \cup (B^\ast \cap B)$, and 
$F'_i = F_i \sm (B^\ast \cap B)$ for $i \geq 2$; the added measure is still at most $Cr^2$,
and we did not need to worry about $D \sm B$ anyway, because $\mu(\ol B(x,r) \sm B) \leq C r^2$
by \eqref{39.18} and \eqref{39.19} anyway.

So let us now worry about the lower bound. This will be a little more complicated, 
and will rely on the isoperimetric inequality. Again we start with the simplest case
when \eqref{39.49} holds, and in addition 
\begin{equation} \label{39.50}
B(x,2r) \subset W_+
\end{equation}
(the upper half space). Of course the case when $B(x,2r) \subset B \sm W_+$ would be
similar, and then we will need to worry about balls centered on (or near) $P$.
 
Suppose that $\mu(B(x,r)) \leq c r^2$, with $c$ very small. 
First observe that, by the isoperimetric inequality in $B(x,r)$, three of the sets $F_j$ have 
very small measures in $B(x,r)$, and the last one is most of $B(x,r)$; more precisely, 
since we know that $F_4$ does not contribute, either
\begin{equation} \label{32.17}
|F_1 \cup F_2 \cap B(x,r)| \leq C \mu(B(x,r))^{3/2} \leq C c^{3/2} r^3,
\end{equation}
or we have the same estimate for some other combination of indices ($1$ and $3$ or $2$ and $3$). 

Let us assume that we have \eqref{32.17} (the two other cases will be similar),
or even more generally that 
\begin{equation} \label{32.18}
r^{-3} |F_1 \cup F_2 \cap B(x,r)| \leq c',
\end{equation}
with some very small $c'$. We want to construct a competitor $\ub F'$ for $\ub F$ in 
$B(x,r) \cup \sigma(B(x,r))$, and we first select, by a Fubini and Chebyshev argument, 
a radius $r_1 \in (9r/10, r)$ such that 
\begin{equation} \label{32.19}
H^2((F_1 \cup F_2) \cap \d B(x,r_1)) \leq 10 r^{-1} |F_1 \cup F_2 \cap B(x,r)| \leq 10 c' r^2.
\end{equation}
Then we set $B_1 = B(x,r_1)$, $B_1^\ast = B_1 \cup \sigma(B_1)$, and
define $\ub F' = (F'_1,F'_2,F'_3,F'_4)$ by $F'_1 = F_1 \sm B_1^\ast$,
$F'_2 = F_2 \sm B_1^\ast$, $F'_3 = F_3 \cup B_1$, and $F'_4 = F_4 \cup \sigma(B_1)$. 
Again $\ub F' \in \cF$ because we are far from $D \sm B$ and we respected the symmetry
constraints. We save $\mu(B_1)$ for the total perimeter in $B_1$ because we no longer 
have boundaries inside of $B_1$, and similarly on $\sigma(B_1)$, but we may
need to spend an extra $2 H^2((F_1 \cup F_2) \cap \d B_1))$ because of the new
discontinuity along $\d B_1 \cup \sigma(\d B_1)$. 
Since $\ub F$ is minimal, we get that
\begin{equation} \label{32.20}
\mu(B_1) \leq 2H^2((F_1 \cup F_2) \cap \d B_1).
\end{equation}
In turn the proof of \eqref{32.17} yields for this smaller ball
\begin{equation} \label{32.21}
\begin{aligned}
|F_1 \cup F_2 \cap B_1| &\leq C \mu(B_1)^{3/2}
\leq C H^2((F_1 \cup F_2) \cap \d B_1)^{3/2}
\cr& \leq C^2 r^{-3/2} |F_1 \cup F_2 \cap B(x,r)|^{3/2},
\end{aligned}
\end{equation}
by \eqref{32.19}. We keep a power $1/2$ for decay, and record that
\begin{equation} \label{32.22}
r_1^{-3} |F_1 \cup F_2 \cap B_1| \leq a r^{-3} |F_1 \cup F_2 \cap B(x,r)|,
\end{equation}
with $a = 2 C^2 r^{-3/2} |F_1 \cup F_2 \cap B(x,r)|^{1/2} \leq 2 C_2 \sqrt{c'}$
because $r^{-3} |F_1 \cup F_2 \cap B(x,r)| \leq c'$. That is, the density of $F_1 \cup F_2$
really decreased, and we can apply the same argument again and again, on a sequence of balls
$B(x,r_k)$. Even more: as the density decreases, $a$ gets smaller and this leaves some 
room to take the $k^{th}$ ratio $r_k/r_{k-1}$ closer and closer to $1$, to the extent 
that actually $r_k \geq r/2$ for all $k$. Yet $r_k^{-3} |F_1 \cup F_2 \cap B(x,r_k)|$ tends to
$0$ and at the end 
\begin{equation} \label{32.23}
 |F_1 \cup F_2 \cap B(x,r/2)| = 0.
\end{equation}
We refer to \cite{Pacific} for the organization of the sequence $\{ r_k \}$, the computations
and also details on the argument that follows.

When we have $F_1 \cup F_3$ in \eqref{32.17}, we proceed as above, except that we take
$F'_1 = F_1 \sm B_1^\ast$, $F'_2 = F_2 \cup B^\ast$, $F'_3 = F_3 \sm B(x,r)$, 
and $F'_4 = F_4 \sm \sigma(B(x,r))$. And the last case of $F_2 \cup F_3$ is similar.

Let us rephrase what we proved: if $r^{-2}\mu(B(x,r))$ is small enough, and for instance
$F_1$ and $F_2$ are the two small sets in $B(x,r)$, we get that $|F_1 \cup F_2 \cap B(x,r/2)| = 0$,
which implies that $\mu(B(x,r/2)) = 0$ and the closed support $Z$ does not meet $B(x,r/2)$.
We would get the same conclusion for other choices of small sets, and of course this excludes
the possibility that $x\in Z$; the first part of \eqref{39.48} follows, in the special case when
$B(x,2r) \subset B \cap W_+$ (as in \eqref{39.49} and \eqref{39.50}).

We are left with the case when $B(x,2r)$ crosses $P$ or $\d B$. 
First consider the case when $x\in P$ and $B(x,2r) \subset B$. 
Again suppose that $\mu(B(x,r)) \leq c r^2$, with $c$ very small. This means that 
all the $F_i$, except one, have a small measure in $B(x,r)$. The unique large one cannot be 
$F_3$ or $F_4$, because by symmetry they would both be large, so it is $F_1$ or $F_2$.
Let us assume that it is $F_1$; the other case would be similar. Instead of \eqref{32.17},
we now have 
\begin{equation} \label{32.58}
|F_2 \cup F_3 \cup F_4 \cap B(x,r)| \leq C \mu(B(x,r))^{3/2} \leq C c^{3/2} r^3,
\end{equation}
where we write the information about $F_4$, but it is the same as for $F_3$.
Notice that here $\sigma(B(x,r)) = B(x,r)$; 
we now use the competitor $\ub F'$ defined by
$F'_1 = F_1 \cup B(x,r)$ and $F'_j = F_i \sm B(x,r)$ for $i \geq 2$. Then we follow
the same proof as above when $F_2 \cup F_3$ is small in $B(x,r)$, just replacing
$F_3$ by $F_3 \cup F_4$, and working with balls centered at $x\in P$.
This takes care of the case when $x\in P$ and $B(x,2r) \subset B$. 

When $B(x,2r) \subset B$ but $B(x,2r)$ meets $P$, we distinguish between two cases.
If $B(x,2r/3)$ meets $Z \cap P$ at some point $y$
(recall that $Z$ is the closed support of $\mu$), we can apply the case when $x\in P$ 
to the ball $B(y,r/4)$, and get the desired lower bound. Otherwise, $P \cap B(x,2r/3)$
is contained in a single component $F_i$, and clearly $i \neq 3,4$ because the corresponding
components lie in half spaces. So we may assume that $P \cap B(x,2r/3) \subset F_1$
(the other case would be similar).
In this case, we still take $F'_1 = F_1 \cup B^\ast$, with $B^\ast : = B(x,r) \cup \sigma(B(x,r))$,
and $F'_j = F_i \sm B^\ast$ for $i \geq 2$;
this does not introduce new discontinuities along $P$, and we can continue with the same 
estimates as in the general case. 

We may now assume that $B(x,2r)$ meets $\d B$. Again the case when $x\in \d B$ 
is a little easier. Suppose first that $x\in F_1 \cap \d B$; then since at least one third of $B(x,r)$ 
lies in $F_1$, the large component has to be $F_1$. In particular, $B(x,r/2)$ does not meet
$V \cup \sigma(V)$ (or, this amounts to the same thing, $G_3$ or $G_4$).
We can then proceed (almost) as in the general case, i.e., systematically replace $F_1$
with $F_1 \cup B(x,r/2) \cup \sigma(B(x,r/2))$ and the other $F_i$ by 
$F_i \sm [B(x,r/2) \cup \sigma(B(x,r/2))]$; the estimates stay the same. 
In the situation when $x\in F_3$ ($F_4$ would be similar), we first observe that
$B(x,r/2)$ does not meet $V \cup \sigma(V)$ (otherwise there would be too much of 
$F_1 \cup F_2$ in $B(x,r)$), so we can proceed as in the general case, with
$F'_3 = F_3 \cup B(x,r)$, $F'_4 = F_4 \cup \sigma(B(x,r))$, and 
$F'_i = F_i \sm [B(x,r/2) \cup \sigma(B(x,r/2))]$. 

In the slightly more general case where $B(x,2r)$ meets $\d B$, either $\d B$ meets
$Z \cap B(x,2r/3)$ at some point $y$, and we can use the previous case on $B(y,r/4)$,
or else the whole $\d B \cap B(x,2r/3)$ is contained in a single $F_i$.
In this last case, $F_i$ is the large component, and we can proceed as in the general case,
pouring all the other components inside $F_i$ (or $F_3 \cup F_4$) if $i = 3, 4$).
The computations are the same as usual, and with this this last case we end the long
proof of our local Ahlfors-regularity estimate \eqref{39.48}.

\ms
Incidentally, we claim that the proof of \eqref{39.48} also shows that in the same circumstances, 
the largest of the $|F_i \cap B(x,r)|$ cannot be too large. That is, there is a small constant 
$c > 0$ such that
\begin{equation} \label{39.59}
|B(x,r) \sm F_i| \geq c r^3
\ \text{ for $i = 1, 2, 3, 4$, $x \in Z \cap B$, and } 0 < r < 10^{-1}.
\end{equation}
Indeed, if this fails, then $|B(x,r) \cap F_j| \leq c r^3$ for all $j \neq i$,
and we can follow the proof of \eqref{39.48}, except that we start from the analogue of 
\eqref{32.19} and \eqref{32.58} (now with three small components) to find a first ball 
$B(x,r_1)$ for which $r_1^{-2} \mu(B(x,r_1))$ is small, proceed as above, and eventually 
get that $|F_j \cap B(x,r/2)| = 0$ for these $j$. As before, if this happens, $Z$ does not meet 
$B(x,r/2)$, and this contradiction proves \eqref{39.59}.

By \eqref{39.48} and a small covering argument, $\mu$ is locally equivalent to 
$\H^d_{\vert Z}$ on $B$. See for instance Lemma 18.11 and Exercise 18.25 in \cite{MSBook}. 
Hence $Z$ is a locally Ahlfors regular of dimension $2$, and of course $|Z| = 0$ 
(we control things well on $B_2 \sm B$).

Notice that on the open set $D \sm Z$, the four functions $\1_{F_i}$ are locally constant 
(almost everywhere), because $\mu = 0$ there (and by Poincar\'e). 
This means that we can assume that the $F_i$ are, after modification on a set of vanishing measure
that we'll assume done, open subsets of $D$, with boundaries that are contained in $Z$. 
In particular, 
\begin{equation} \label{32.25}
\H^2(\d F_i) \leq \H^2(Z) \leq C \mu(D) < +\infty
\end{equation}
where $C$ comes from \eqref{39.48} and the comparison between $\H^2_{\vert Z}$ and $\mu$. 

Here we could quietly continue as in \cite{Pacific} and prove a form of ``Condition B'',
and then the local uniform rectifiability of $Z$, but for the present purposes we don't need as much.
We still need to investigate the relations between $\H^2(Z)$ and our functional $J$.

Return to the reduced boundaries $\d^\ast F_i$ that were used before. 
Recall from \eqref{39.25} that $\mu_i = \H^2_{\vert \d^\ast F_i}$.
We claim that, due to the local Ahlfors regularity of $\mu$,
\begin{equation} \label{32.26}
\H^d(Z \sm \bigcup_i \d^\ast F_i) = 0.
\end{equation}
Indeed, otherwise (since $\H^2(Z) < +\infty$) we can find a point
$x\in Z$, which is a Lebesgue density point for $Z$, but where $\cup_i \d^\ast F_i$ has vanishing
density (see for instance the density properties of sets in \cite{Mattila}). 
This contradicts \eqref{39.48} or the definition of $Z$ as the closed support of $\mu = \sum \mu_i$. 

We are ready to compare $J(\ub F)$ to Hausdorff measures that comes from the 
strong separation problem. Recall from \eqref{39.21} that
\begin{equation} \label{39.62}
J(\ub F) = \mu_3(P)+\mu_4(P) + \sum_{i=1}^4 \mu_i(D)
= \H^2(P \cap \d^\ast F_3) + \H^2(P \cap \d^\ast F_4) + \sum_{i=1}^4 \H^2(D \cap \d^\ast F_i);
\end{equation}
we need to know how the different reduced boundaries $\d^\ast F_i$ match up.
For instance, if $x\in P \cap \d^\ast F_3$, then it also lies in $\d^\ast F_4$,
by symmetry. In general, a point of some $\d^\ast F_i$ lies in at most one other
$\d^\ast F_j$, because \eqref{39.27} say that each $F_i$ takes at least about 
half of the measure in small balls. We claim that conversely, $\H^2$-almost every
point of $Z$ lies in at least two sets $\d^\ast F_j$.

So let $x\in Z$ be given; \eqref{32.26} says that (almost surely) $x\in \d^\ast F_i$ for some $i$.
For each small enough $r$, \eqref{39.26} says that nearly half of $B(x,r)$ lies 
outside of $F_i$, so we can find $j \neq i$ such that $|F_j \cap B(x,r)| \geq 10^{-1} |B(x,r)|$.
By Poincar\'e's inequality inside of $B(x,r)$ (and since $F_i$ also is large in that ball),
$\mu_j(B(x,r)) \geq C^{-1} r^2$. Taking a sequence of radii $r$ that tends to $0$ and for which
this happens for the same $j$, we see that $x$ is a point of positive density for $\mu_j$.
But, since $\mu_j = \H^2_{\vert \d^\ast F_j}$ and $\H^2(\d^\ast F_j) < +\infty$,
this $\H^2$-almost never happens when $x \notin \d^\ast F_j$ (see \cite{Mattila} again). 
Our claim follows. That is, almost every point of $Z$ lies in exactly two sets $\d^\ast F_i$,
and \eqref{39.62} yields
\begin{equation} \label{39.63}
J(\ub F) =  \H^2(P \cap \d^\ast F_3) + \H^2(P \cap \d^\ast F_4) + 2 \H^2(D \cap Z).
\end{equation}
We have seen that $P \cap \d^\ast F_3$ is the same as $P \cap \d^\ast F_4$.
In addition, $\H^2(P \cap \d^\ast F_i) = 0$ for $i=1,2$, because $\H^2$-almost everywhere
on $P$, the approximate tangent plane to $\d^\ast F_i$ is $P$ itself (by the uniqueness almost everywhere of the approximate tangent, or because almost everywhere on $P \cap \d^\ast F_i$,
the density of the difference of these two sets vanishes), and then \eqref{39.26} and \eqref{39.27}
are incompatible with the symmetry of $F_i$. Thus
$\H^2(P \cap Z) = \H^2(P \cap \d^\ast F_3) = \H^2(P \cap \d^\ast F_4)$, and \eqref{39.63}
now says that
\begin{equation} \label{39.64}
J(\ub F) =  2 \H^2(D \cap Z) + 2 \H^2(P \cap Z) = 4 \H^2(D_+ \cap Z),
\end{equation}
with $D_+ = D \cap W_+$ as above. 

Recall that we want to show that $\ub F$ also yields a minimizer $E_0$ for the strong 
separation problem. We take $E_0 = Z \cap 2B_+ = Z \cap 2B \cap W_+$;
it lies in the class $\cF_s$ of strong competitors defined above \eqref{39.16}, 
because $Z$ separates the various $F_i$ also in the smaller $B_+$, and 
\begin{equation} \label{39.65}
J(\ub F) = 4 \H^2(E_0)
\end{equation}
by \eqref{39.64} and because $\H^2(E_0 \sm [D_+ \cap Z])  = \H^2(E_0 \cap \d D) = 0$. 

The last step of our long proof consists in taking any other competitor $E \in \cF_s$
and showing that $\H^2(E) \geq \H^2(E_0)$, and for this we first associate to $E$
a competitor $\ub F' = (F'_1, F'_2, F'_3, F'_4)$. 

First double the set $E$, i.e. set $E^\ast = E \cup \sigma(E) \subset 2B$, and denote by 
$F'_i$, $1 \leq i \leq 3$, the connected component of $G_i$ of in $2B \sm E^\ast$.
Also let $F'_4$ denote the component of $G_4 = \sigma(G_3)$, or equivalently (by symmetry) 
set $F_4 = \sigma(F_3)$. We want to show that all these sets are disjoint, i.e., that
$E^\ast$ separates the four $G_i$ in $2B$.

Our assumption that $E$ separates $A_1$ and $A_2$ from $A_3$ in $2B$ implies that
$F_1$ and $F_2$ do not meet $F_3$. They do not meet $F_4$ either, by symmetry.
Also, we know that $E$ separates $A_1$ from $A_2$ in $2B_+$, and then 
$E^\ast$ still separates them in $2B$: suppose instead that there exists a path 
$\gamma \subset 2B \sm E^\ast$ that joins them, replace $\gamma$ with a path 
$\gamma_+ \subset 2B_+$ by replacing  the last coordinate $z$ with $|z|$ along $\gamma$, 
and notice that $\gamma_+$ does not meet $E$, a contradiction. So $F_2$ is disjoint 
from $F_1$. 

We still need to see that $F_3$ does not meet $F_4$, and for this we simply use the fact that
$E$ separates $A_1$ and $A_2$, hence also $A_4 = \sigma(A_3)$, from $A_3$ in $2B$.
So all these connected components are different. It could be that there are other components
in $D \sm E^\ast$, and if this is the case just add them to $F'_1$, say. If we do this, the symmetry
is preserved, because we know that our initial sets $F'_i$ have the required symmetries.
Since $\H^2(E) < +\infty$, we know that all the $F'_i$ are Caccioppoli sets, so 
$\ub F' = (F'_1, F'_2, F'_3, F'_4) \in \cF$. Consequently,
\begin{equation} \label{39.66}
J(\ub F') \geq  J(\ub F) =  4 \H^2(E_0)
\end{equation}
by \eqref{39.65}. We still need prove that 
\begin{equation} \label{39.67}
J(\ub F') \leq  4 \H^2(E), 
\end{equation}
and for this we just follow the argument between \eqref{39.62} and \eqref{39.64},
which we apply to $\ub F'$, but simply replace $Z$ with the union $Z'$ of reduced boundaries
$\d^\ast F'_i$. This yields $J(\ub F') =  4 \H^2(D_+ \cap Z')$, and \eqref{39.67} holds because 
$Z' \subset E$ and $E \subset D_+$.
So we proved that $\H^2(E) \geq \H^2(E_0)$, $E_0$ is a strong minimizer, and it was enough to consider
rectifiable competitors in the proof of Lemma \ref{t32.2}. This completes our proof of
Lemma \ref{t32.3}.
\qed

\bigskip
\vfill \vfill \vfill\vfill
\noindent Guy David,  
\par\noindent 
Laboratoire de Math\'{e}matiques d'Orsay, 
\par\noindent 
Univ. Paris-Sud, CNRS, Universit\'e Paris-Saclay
\par\noindent 
91405 Orsay, France
\par\noindent 
guy.david@math.u-psud.fr
\par\noindent 
http://www.math.u-psud.fr/$\sim$gdavid/

\end{document}